\documentclass[11pt]{article}
\usepackage[a4paper,margin=25mm]{geometry}
\usepackage[T1]{fontenc}
\usepackage{lmodern}
\usepackage{amsmath,amssymb,amsthm,mathtools}
\usepackage{microtype}
\usepackage{needspace}
\usepackage{xcolor}
\usepackage{tikz}
\usepackage[unicode,colorlinks=true,allcolors=black]{hyperref}
\hypersetup{pdfauthor={Yuanyang Hu},pdftitle={Periodic Fisher--KPP equations on a regular tree},
 pdfsubject={Autonomous spatial heterogeneity: spectral threshold, spreading speeds, and logarithmic estimates}}
\allowdisplaybreaks[1]
\numberwithin{equation}{section}

\newtheorem{lemma}{Lemma}[section]
\newtheorem{theorem}[lemma]{Theorem}
\newtheorem{corollary}[lemma]{Corollary}
\newtheorem{proposition}[lemma]{Proposition}
\theoremstyle{remark}
\newtheorem{remark}[lemma]{Remark}
\newtheorem*{remark*}{Remark}
\newcommand{\Tq}{T_q}
\newcommand{\Z}{\mathbb Z}
\newcommand{\R}{\mathbb R}
\newcommand{\N}{\mathbb N}
\newcommand{\C}{\mathbb C}
\newcommand{\ind}{\mathbf1}
\newcommand{\norm}[1]{\lVert#1\rVert}
\newcommand{\supp}{\operatorname{supp}}
\newcommand{\spec}{\operatorname{spec}}
\newcommand{\diag}{\operatorname{diag}}
\newcommand{\re}{\operatorname{Re}}
\newcommand{\im}{\operatorname{Im}}
\newcommand{\dd}{\mathrm d}
\title{Periodic Fisher-KPP equations on a regular tree}
\author{Yuanyang Hu\thanks{Corresponding author.}\\[3pt]
School of Mathematics and Statistics, Henan University\\
475004 Kaifeng, China}
\date{}
\begin{document}
\maketitle
\begin{abstract}
We study Fisher-KPP equations on regular trees with reactions periodic in a Busemann coordinate or in distance from a vertex or a bi-infinite geodesic. For nonzero finitely supported data, we characterize the spectral threshold between uniform extinction and persistence, including critical extinction. We establish Bramson corrections for Busemann step data at positive, zero and negative speeds and, above the periodic extinction threshold, for the largest Busemann coordinate and radial level positions of finitely supported solutions. When the periodic spectral minimum is nonpositive but the geodesic spectral threshold is positive, we construct a smallest positive stationary solution with exponential transverse decay, prove uniqueness under an exponential bound, and determine longitudinal spreading with its logarithmic correction. We also classify positive stationary solutions depending only on distance from the geodesic and prove existence and stability results for moving and stationary fronts.
\end{abstract}


\clearpage
\setcounter{tocdepth}{2}
\makeatletter
\IfFileExists{\jobname.toc}{%
  \tableofcontents
}{%
  \section*{\contentsname}
\contentsline {section}{\numberline {1}Introduction}{3}{section.1}%
\contentsline {subsection}{\numberline {1.1}The tree and coefficients}{3}{subsection.1.1}%
\contentsline {subsection}{\numberline {1.2}Background and the questions addressed here}{6}{subsection.1.2}%
\contentsline {subsection}{\numberline {1.3}Main results}{7}{subsection.1.3}%
\contentsline {paragraph}{Extinction thresholds.}{7}{section*.2}%
\contentsline {paragraph}{Logarithmic level positions.}{9}{section*.3}%
\contentsline {paragraph}{A common derivative at zero.}{11}{section*.4}%
\contentsline {paragraph}{Moving and stationary fronts.}{11}{section*.5}%
\contentsline {paragraph}{Positive stationary solutions near a geodesic.}{12}{section*.6}%
\contentsline {paragraph}{Notation.}{15}{section*.7}%
\contentsline {section}{\numberline {2}Comparison and periodic kernel estimates}{16}{section.2}%
\contentsline {subsection}{\numberline {2.1}The periodic matrix and the associated chain}{23}{subsection.2.1}%
\contentsline {subsection}{\numberline {2.2}The continuous-time chain}{25}{subsection.2.2}%
\contentsline {subsection}{\numberline {2.3}Absorption at a straight line}{29}{subsection.2.3}%
\contentsline {subsection}{\numberline {2.4}A lower bound away from the boundary}{32}{subsection.2.4}%
\contentsline {subsection}{\numberline {2.5}Changing the exponential parameter}{36}{subsection.2.5}%
\contentsline {subsection}{\numberline {2.6}A decreasing logarithmic boundary}{40}{subsection.2.6}%
\contentsline {paragraph}{The backward logarithmic boundary.}{42}{section*.8}%
\contentsline {subsection}{\numberline {2.7}The linearized equation}{45}{subsection.2.7}%
\contentsline {section}{\numberline {3}Subsolutions and moving boundaries}{46}{section.3}%
\contentsline {subsection}{\numberline {3.1}Comparison at a moving boundary}{49}{subsection.3.1}%
\contentsline {subsection}{\numberline {3.2}A half-line equation}{57}{subsection.3.2}%
\contentsline {section}{\numberline {4}Tree spectra and fronts}{59}{section.4}%
\contentsline {subsection}{\numberline {4.1}Comparison reactions with a common derivative at zero}{64}{subsection.4.1}%
\contentsline {subsection}{\numberline {4.2}Front equations}{69}{subsection.4.2}%
\contentsline {subsection}{\numberline {4.3}Moving fronts}{69}{subsection.4.3}%
\contentsline {subsection}{\numberline {4.4}Stationary fronts}{76}{subsection.4.4}%
\contentsline {subsection}{\numberline {4.5}Weighted difference estimates}{79}{subsection.4.5}%
\contentsline {section}{\numberline {5}Distance from a geodesic}{83}{section.5}%
\contentsline {subsection}{\numberline {5.1}The spectral upper endpoint}{83}{subsection.5.1}%
\contentsline {subsection}{\numberline {5.2}A uniqueness identity}{89}{subsection.5.2}%
\contentsline {subsection}{\numberline {5.3}Uniqueness under an exponential bound}{91}{subsection.5.3}%
\contentsline {subsection}{\numberline {5.4}Propagation and an eight-periodic example}{94}{subsection.5.4}%
\contentsline {subsection}{\numberline {5.5}Heat estimates}{98}{subsection.5.5}%
\contentsline {section}{\numberline {6}Proofs of the main results}{101}{section.6}%
\contentsline {subsection}{\numberline {6.1}Extinction and spreading}{101}{subsection.6.1}%
\contentsline {subsection}{\numberline {6.2}Periodic linear growth at zero}{104}{subsection.6.2}%
\contentsline {subsection}{\numberline {6.3}Common linear growth at zero}{107}{subsection.6.3}%
\contentsline {subsection}{\numberline {6.4}Moving and stationary fronts}{109}{subsection.6.4}%
\contentsline {subsection}{\numberline {6.5}Stability}{109}{subsection.6.5}%
\contentsline {subsection}{\numberline {6.6}Distance from a geodesic}{111}{subsection.6.6}%
\contentsline {subsection}{\numberline {6.7}Tangential spreading along the geodesic}{115}{subsection.6.7}%
\contentsline {subsection}{\numberline {6.8}Tangential Bramson correction}{123}{subsection.6.8}%
\contentsline {paragraph}{Two-coordinate moving-domain formula.}{147}{section*.9}%
\contentsline {subsection}{\numberline {6.9}Data with infinite support}{157}{subsection.6.9}%
\contentsline {subsection}{\numberline {6.10}Two-sided kernel bounds at a straight line}{158}{subsection.6.10}%
\contentsline {subsection}{\numberline {6.11}A logarithmic boundary}{159}{subsection.6.11}%
\contentsline {section}{\numberline {7}Examples and further questions}{162}{section.7}%
\contentsline {subsection}{\numberline {7.1}A two-periodic reaction}{162}{subsection.7.1}%
\contentsline {subsection}{\numberline {7.2}Dependence on the branching parameter}{163}{subsection.7.2}%
\contentsline {subsection}{\numberline {7.3}Level positions and profile convergence}{164}{subsection.7.3}%
\contentsline {subsection}{\numberline {7.4}Stationary solutions with transverse decay}{166}{subsection.7.4}%
\contentsline {subsection}{\numberline {7.5}Further questions}{167}{subsection.7.5}%
\@starttoc{toc}%
}
\makeatother
\clearpage

\section{Introduction}

How does a periodic arrangement of reactions on a regular tree
determine extinction, persistence and the positions of level sets?
We answer this question for reactions that depend periodically on
three integer coordinates of the tree.

Let $\Tq$ be the infinite tree in which every vertex has $q+1$
neighbors, where $q\ge2$ is an integer. Write $x\sim y$ when $x$
and $y$ are adjacent. For a function $v:\Tq\to\R$, define
\begin{equation}\label{eq:tree-laplacian}
 \Delta_{\Tq}v(x)=\sum_{y\sim x}(v(y)-v(x))
 \quad\text{ for }x\in\Tq.
\end{equation}
For $D>0$, $u_0:\Tq\to[0,1]$ and a reaction
$f:\Tq\times[0,1]\to\R$, we study
\begin{equation}\label{eq:tree-cauchy}
 \left\{\begin{aligned}
 u_t(t,x)&=D\Delta_{\Tq}u(t,x)+f(x,u(t,x))
 &&\text{for }t>0\text{ and }x\in\Tq,\\
 u(0,x)&=u_0(x)&&\text{for }x\in\Tq.
 \end{aligned}\right.
\end{equation}
Even if $\partial_s f(x,0)>0$ at every vertex, a solution with
$u_0\not\equiv0$ and $\#\{x\in\Tq:u_0(x)>0\}<\infty$ may satisfy
$\sup_{x\in\Tq}u(t,x)\to0$ as $t\to\infty$.
The thresholds in the main results quantify the balance between
$D\Delta_{\Tq}$ and $f$. We first define the three integer coordinates used to assign $f$.

\subsection{The tree and coefficients}

Let $E=\{\{x,y\}:x\sim y\}$ be the edge set of $\Tq$.
Write $d(x,y)$ for the number of edges in the unique simple path
from $x$ to $y$ in $\Tq$; see \cite[Chapter~1]{Diestel} for graph
terminology. Fix $o\in\Tq$, and define $r(x)=d(o,x)$ and
$B(o,R)=\{x\in\Tq:r(x)\le R\}$ for $R\in\R$.
Then $B(o,R)=\varnothing$ for $R<0$, while
\[
 \#\{x\in\Tq:r(x)=r\}=(q+1)q^{r-1}
 \quad\text{ for every integer }r\ge1.
\]
Thus the number of vertices at distance $r$ grows exponentially in $r$. We use $\N=\{1,2,\ldots\}$ and
$\N_0=\{0,1,2,\ldots\}$. Unless a weight is explicitly indicated,
sequence spaces use counting measure. Their index set $S$ is one of
$\Tq$, $\Z$, $\N_0$, $E$, an explicitly specified subset of one
of these sets, the vertex set of an explicitly defined rooted tree,
or one of the products $\Z\times\N_0$ and $\N_0\times\N_0$.
For each such $S$, the norms are
\[
 \|v\|_p=\left(\sum_{x\in S}|v(x)|^p\right)^{1/p}\quad\text{ for } 1\le p<\infty,
 \qquad \|v\|_\infty=\sup_{x\in S}|v(x)|.
\]
For $1\le p\le\infty$, the space $\ell^p(S)$ consists of the sequences
on $S$ with finite norm $\|\cdot\|_p$. We omit $S$ when the underlying
set is clear. For $v,w\in\ell^2(S)$, put
$\langle v,w\rangle=\sum_{x\in S}\overline{v(x)}w(x)$.
For $A\subseteq S$, let $\mathbf1_A:S\to\{0,1\}$ equal one on
$A$ and zero on $S\setminus A$. For $x\in S$, set
$\delta_x=\mathbf1_{\{x\}}$. For $v:S\to\C$, write
$\supp v=\{x\in S:v(x)\ne0\}$.
For real $z$, set $z_+=\max\{z,0\}$; for complex $z$,
$\re z$ and $\im z$ denote its real and imaginary parts. Constants denoted by $C$ may change from line to line.
A constant described as uniform is independent of the parameters over
which uniformity is asserted. For real-valued functions $A,B$ on $[0,\infty)$ with $B(t)>0$,
the notation $A(t)=O(B(t))$ as $t\to\infty$ means that there are
$C>0$ and $T\ge0$ such that $|A(t)|\le CB(t)$ for $t\ge T$;
$A(t)=o(B(t))$ as $t\to\infty$ means $A(t)/B(t)\to0$ as $t\to\infty$. A subscript on $O$ lists
additional quantities on which $C$ and $T$ may depend. For sequences,
the same definitions apply as the integer index tends to infinity.
For limits at finite points or at negative infinity, the variable and
direction are stated explicitly.
If $A(t,\xi)$ and $B(t,\xi)>0$ are indexed by $\xi$ in a
nonempty set $K$, uniform $O(B)$ means that the same $C,T$ work
for every $\xi\in K$. Uniform $o(B)$ means
\[
 \sup_{\xi\in K}\frac{|A(t,\xi)|}{B(t,\xi)}\to0
 \quad\text{ as }t\to\infty.
\]

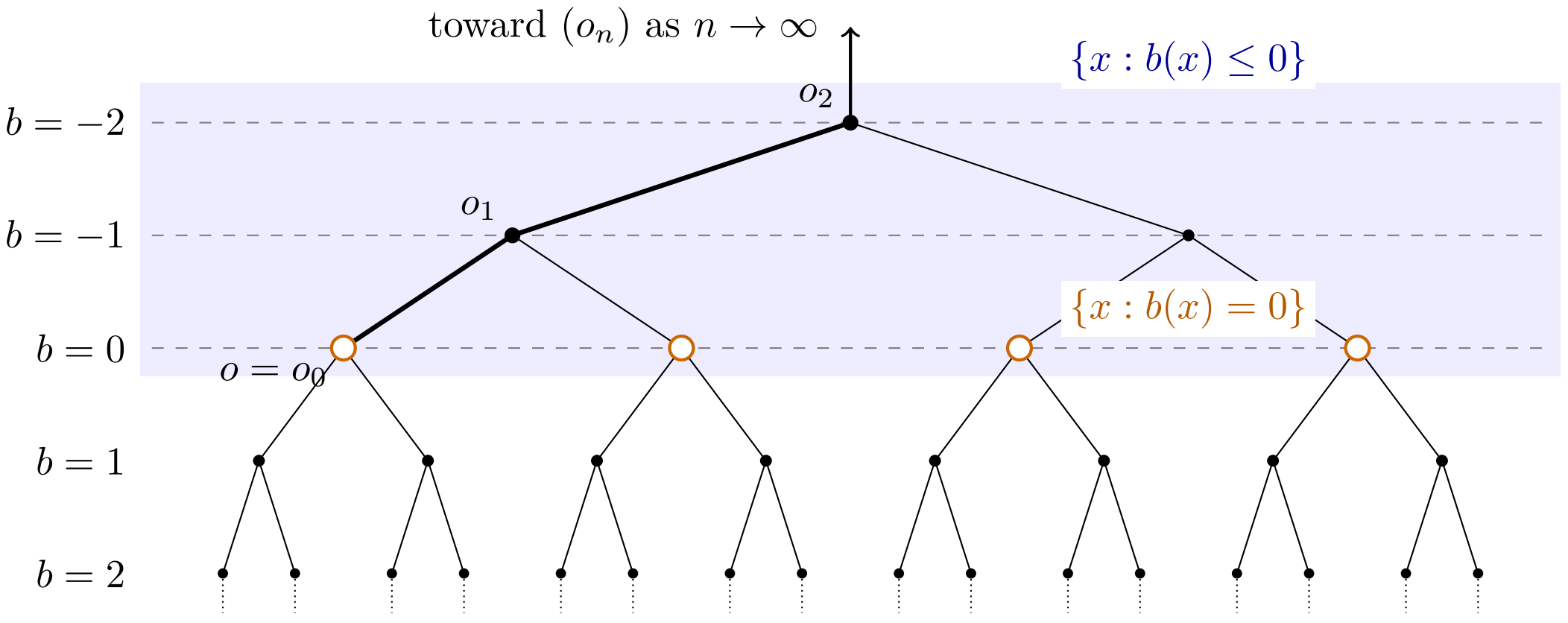
\begin{figure}[htbp]
\centering
\begin{tikzpicture}[x=1cm,y=1cm, every node/.style={font=\small}]
 \fill[blue!7] (-6.3,-.25) rectangle (6.3,2.35);
 \draw[gray,dashed] (-6.2,0)--(6.2,0);
 \draw[gray,dashed] (-6.2,1)--(6.2,1);
 \draw[gray,dashed] (-6.2,2)--(6.2,2);
 \node[left] at (-6.3,2) {$b=-2$};
 \node[left] at (-6.3,1) {$b=-1$};
 \node[left] at (-6.3,0) {$b=0$};
 \node[left] at (-6.3,-1) {$b=1$};
 \node[left] at (-6.3,-2) {$b=2$};
 \draw[->,thick] (0,2)--(0,2.85) node[left,xshift=-4pt] {toward $(o_n)$ as $n\to\infty$};
 \foreach \x in {-3,3} {\draw (0,2)--(\x,1);}
 \foreach \x/\parent in {-4.5/-3,-1.5/-3,1.5/3,4.5/3} {\draw (\parent,1)--(\x,0);}
 \foreach \x/\parent in {-5.25/-4.5,-3.75/-4.5,-2.25/-1.5,-.75/-1.5,.75/1.5,2.25/1.5,3.75/4.5,5.25/4.5} {
  \draw (\parent,0)--(\x,-1);
  \draw (\x,-1)--(\x-.32,-2);
  \draw (\x,-1)--(\x+.32,-2);
  \foreach \dx in {-.32,.32} {\fill (\x+\dx,-2) circle (1.3pt);\draw[densely dotted] (\x+\dx,-2)--(\x+\dx,-2.35);}
  \fill (\x,-1) circle (1.5pt);
 }
 \draw[very thick] (-4.5,0)--(-3,1)--(0,2);
 \fill (0,2) circle (2pt) node[above left] {$o_2$};
 \fill (-3,1) circle (2pt) node[above left] {$o_1$};
 \fill (3,1) circle (1.5pt);
 \foreach \x in {-4.5,-1.5,1.5,4.5} {\draw[orange!80!black,thick,fill=white] (\x,0) circle (3pt);}
 \node[below left] at (-4.5,0) {$o=o_0$};
 \node[blue!60!black,fill=white,inner sep=2pt] at (3,2.55) {$\{x:b(x)\le0\}$};
 \node[orange!70!black,fill=white,inner sep=2pt] at (3,.35) {$\{x:b(x)=0\}$};
\end{tikzpicture}
\caption{A finite portion of $T_2$. The thick path follows the chosen ray
from $o=o_0$ toward $o_1,o_2,\ldots$; along this ray, $b(o_n)=-n$.
The outlined vertices lie in the horocycle $\{x:b(x)=0\}$.
The shaded region represents the horoball $\{x:b(x)\le0\}$.
Horizontal dashed lines join vertices with equal Busemann coordinate;
solid sloping lines represent tree edges. The tree and each indicated level continue beyond the drawing.}
\label{fig:busemann}
\end{figure}

Fix vertices $(o_n)_{n\ge0}$ with $o_0=o$ and
$d(o_m,o_n)=|m-n|$ for $m,n\in\N_0$. This sequence is a geodesic ray.
Define its Busemann coordinate $b:\Tq\to\Z$ by
\[
 b(x)=\lim_{n\to\infty}\bigl(d(x,o_n)-n\bigr)\quad\text{for }x\in\Tq.
\]
For $x\in\Tq$, let $\ell\in\{0,\ldots,r(x)\}$ be the number
of initial edges shared by the path from $o$ to $x$ and the ray.
Then $d(x,o_n)=r(x)+n-2\ell$ for every integer $n\ge\ell$. Thus the
limit stabilizes, $b(x)=r(x)-2\ell$, and $b(o)=0$. Every vertex has
one neighbor $x^-$ with $b(x^-)=b(x)-1$, lying on the path from $x$
to $o_n$ for all sufficiently large $n$, and $q$ neighbors with
Busemann coordinate $b(x)+1$. We call $x^-$ the parent of $x$ and
the other $q$ neighbors its children. For $j\in\Z$, the sets
$\{x\in\Tq:b(x)=j\}$ and $\{x\in\Tq:b(x)\le j\}$
are a horocycle and a horoball, respectively; see
\cite{Woess,BartholdiNeuhauserWoess}. Figure~\ref{fig:busemann} illustrates
these sets with the sign convention used here.

\begin{figure}[htbp]
\centering
\begin{minipage}[t]{.38\textwidth}
\centering
\begin{tikzpicture}[scale=.78,every node/.style={font=\small}]
 \foreach \a in {90,210,330} {
  \draw (0,0)--(\a:1.2);
  \fill (\a:1.2) circle (2pt);
  \foreach \b in {-20,20} {
   \draw (\a:1.2)--(\a+\b:2.4);
   \fill (\a+\b:2.4) circle (2pt);
   \foreach \c in {-8,8} {
    \draw (\a+\b:2.4)--(\a+\b+\c:3.5);
    \fill (\a+\b+\c:3.5) circle (1.5pt);
    \draw[dotted] (\a+\b+\c:3.5)--(\a+\b+\c:3.8);
   }
  }
 }
 \draw[blue!70!black,very thick] (0,0)--(90:1.2)--(70:2.4)--(62:3.5);
 \fill (0,0) circle (2.7pt) node[left] {$o$};
 \fill[blue!70!black] (62:3.5) circle (2.7pt) node[above] {$x$};
 \node[right,xshift=3pt,blue!70!black,fill=white] at (62:3.5) {$r(x)=3$};
 \node at (0,-3.6) {(a) Distance from $o$};
\end{tikzpicture}
\end{minipage}\hfill
\begin{minipage}[t]{.60\textwidth}
\centering
\begin{tikzpicture}[x=1.25cm,y=.88cm,every node/.style={font=\small}]
 \draw[very thick] (-3,0)--(3,0);
 \draw[very thick,dotted] (-3.5,0)--(-3,0);
 \draw[very thick,dotted] (3,0)--(3.5,0);
 \foreach \a in {-2,-1,0,1,2} {\fill (\a,0) circle (2pt);}
 \foreach \a in {-2,-1,1,2} {
  \draw (\a,0)--(\a,-.65);
  \draw[dotted] (\a,-.65)--(\a-.25,-.95);
  \draw[dotted] (\a,-.65)--(\a+.25,-.95);
 }
 \draw (0,.8)--(-.65,1.4);
 \draw[dotted] (-.65,1.4)--(-1,1.7);
 \draw[dotted] (-.65,1.4)--(-.45,1.8);
 \draw (.8,1.6)--(1.5,2.05);
 \draw[dotted] (1.5,2.05)--(1.8,2.3);
 \draw[dotted] (1.5,2.05)--(1.9,2);
 \draw[orange!85!black,very thick] (0,0)--(0,.8)--(.8,1.6)--(.8,2.4);
 \foreach \a/\b in {0/.8,.8/1.6,.8/2.4} {\fill (\a,\b) circle (2pt);}
 \draw[dotted] (.8,2.4)--(.45,2.9);
 \draw[dotted] (.8,2.4)--(1.1,2.9);
 \node[above] at (.8,2.4) {$x$};
 \draw[orange!85!black,thick,fill=white] (0,0) circle (3pt);
 \node[above,yshift=2pt] at (-2,0) {$\gamma(-3)$};
 \node[above,yshift=2pt] at (-1,0) {$\gamma(-2)$};
 \node[below,yshift=-4pt] at (0,0) {$\gamma(-1)$};
 \node[above,yshift=2pt] at (1,0) {$\gamma(0)$};
 \node[above,yshift=2pt] at (2,0) {$\gamma(1)$};
 \node[right] at (3.4,0) {$\mathcal G$};
 \node[right,xshift=3pt,orange!85!black,fill=white] at (.8,2.4) {$s(x)=3$};
 \node at (0,-3.6) {(b) Distance from the geodesic $\mathcal G$};
\end{tikzpicture}
\end{minipage}
\caption{Distance from a fixed vertex and from a bi-infinite geodesic in $T_2$.
The two panels show separate configurations.
In (a), the blue path is the unique shortest path from $o$ to $x$.
It contains three edges, so $r(x)=d(o,x)=3$.
In (b), $\mathcal G=\gamma(\mathbb Z)$ is the bi-infinite geodesic,
and $\gamma(n)$ denotes its vertex with index $n$.
The hollow circle marks the unique nearest vertex
$\gamma(\ell_{\mathcal G}(x))$ of $\mathcal G$ to $x$;
$\ell_{\mathcal G}(x)$ is its integer index.
In (b), $\ell_{\mathcal G}(x)=-1$.
The orange path is the unique shortest path from $x$ to $\mathcal G$.
It contains three edges, so
$s(x)=\min_{n\in\mathbb Z}d(x,\gamma(n))=3$.
Dotted segments indicate continuation of the tree;
$\mathcal G$ extends infinitely in both directions beyond the displayed
vertices $\gamma(-3)$ and $\gamma(1)$.}
\label{fig:tree-distances}
\end{figure}
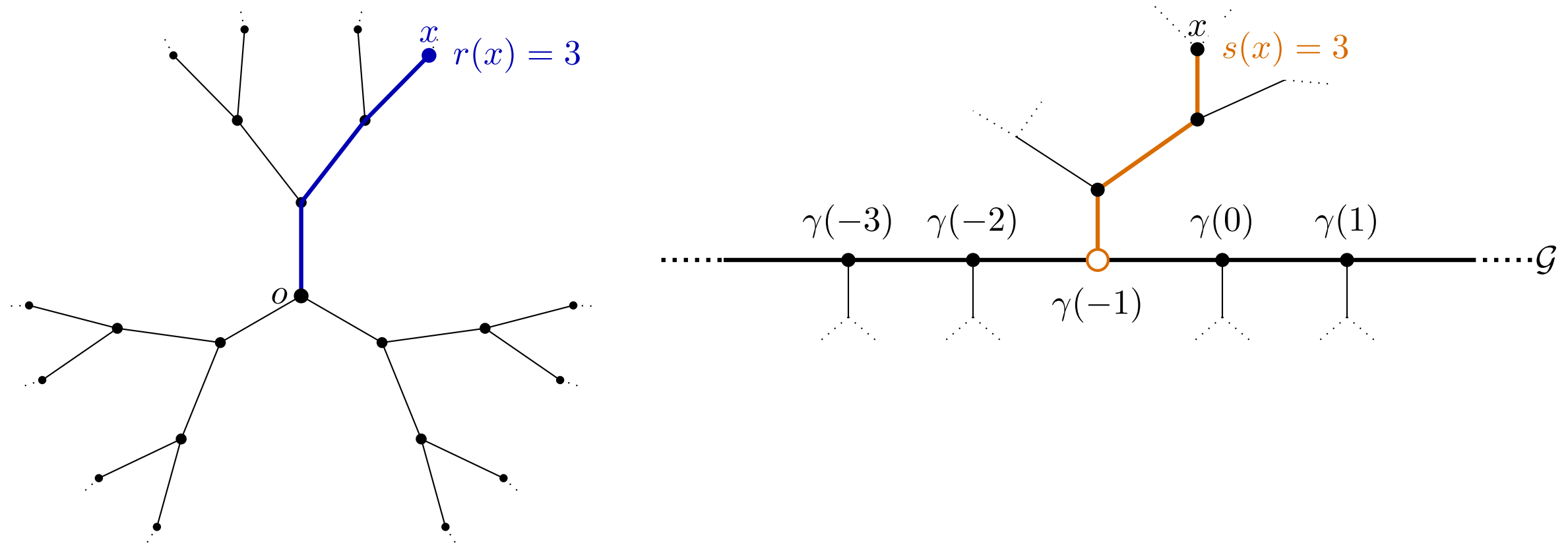

Fix also a bi-infinite geodesic $\gamma:\Z\to\Tq$, with
$d(\gamma(m),\gamma(n))=|m-n|$ for $m,n\in\Z$, and put
\[
 \mathcal G=\gamma(\Z),\qquad
 s(x)=\min_{n\in\Z}d(x,\gamma(n))\quad\text{for }x\in\Tq.
\]
Since $d(x,\gamma(n))\ge |n|-d(x,\gamma(0))$, the minimum defining
$s(x)$ is attained. The unique-path property of the tree gives a
unique minimizing vertex. Define $\ell_{\mathcal G}(x)\in\Z$ by
$d(x,\gamma(\ell_{\mathcal G}(x)))=s(x)$. The neighborhoods
$\mathcal N_R(\mathcal G)=\{x\in\Tq:s(x)\le R\}$, $R\in\R$,
are empty for $R<0$ and infinite for $R\ge0$. If $s(x)\ge1$, then $x$ has one neighbor
$y$ with $s(y)=s(x)-1$ and $q$ neighbors $y$ with
$s(y)=s(x)+1$. A vertex on $\mathcal G$ has two neighbors on
$\mathcal G$ and $q-1$ neighbors $y$ with $s(y)=1$. Figure~\ref{fig:tree-distances}
illustrates $r(x)$ and $s(x)$.

Fix an integer $P\ge1$ and functions $(f_j)_{j\in\Z}$
from $[0,1]$ to $\R$ satisfying $f_{j+P}=f_j$. Thus $f_j=f_i$
for the unique $i\in\{0,\ldots,P-1\}$ with $j-i\in P\Z$.
In particular, a subscript $b(x)$, $r(x)$ or $s(x)$ retains its
coordinate value, and the periodic extension determines the reaction.
We consider the following three assignments separately:
\begin{equation}\label{eq:reaction-assignments}
 f(x,u)=f_{b(x)}(u),\qquad
 f(x,u)=f_{r(x)}(u),\qquad
 f(x,u)=f_{s(x)}(u)\quad\text{for }x\in\Tq,\ u\in[0,1].
\end{equation}
The initial datum $u_0$ is prescribed separately in each theorem.
Since $r(x)-b(x)$ is even, the assignments
$f(x,u)=f_{b(x)}(u)$ and $f(x,u)=f_{r(x)}(u)$ agree for $P=2$. For $P>2$, they can differ already at
$o_1$, where $f_{b(o_1)}=f_{P-1}$ and $f_{r(o_1)}=f_1$.

Let $w=(w_j)_{j\in\Z}$, $U=(U_r)_{r\in\N_0}$ and
$V=(V_r)_{r\in\N_0}$ be families of $C^1$ functions
from $[0,\infty)$ to $[0,1]$. A function $u(t,x)=w_{b(x)}(t)$
satisfies the differential equation
in \eqref{eq:tree-cauchy} with $f(x,u)=f_{b(x)}(u)$ if and only if
\begin{equation}
\label{eq:lattice}
 \dot w_j=D(w_{j-1}+qw_{j+1}-(q+1)w_j)+f_j(w_j),
 \quad\text{ for }t\ge0,\ j\in\Z.
\end{equation}
For $f(x,u)=f_{r(x)}(u)$, a radial solution
$u(t,x)=U_{r(x)}(t)$ satisfies
\begin{equation}
\label{eq:radial}
 \begin{aligned}
 \dot U_r&=D(U_{r-1}+qU_{r+1}-(q+1)U_r)+f_r(U_r),&&\text{ for }t\ge0,\ r\ge1,\\
 \dot U_0&=D(q+1)(U_1-U_0)+f_0(U_0),&&\text{for }t\ge0.
 \end{aligned}
\end{equation}
For $u(t,x)=V_{s(x)}(t)$ and $f(x,u)=f_{s(x)}(u)$, the differential equation
in \eqref{eq:tree-cauchy} reduces to
$\dot V_r=D(V_{r-1}+qV_{r+1}-(q+1)V_r)+f_r(V_r)$ for $r\ge1$,
together with
\[
 \dot V_0=D(q-1)(V_1-V_0)+f_0(V_0)\quad\text{for }t\ge0.
\]
The factors $q+1$ and $q-1$ in the boundary equations count the
neighbors of $o$ at distance one from $o$ and the neighbors of a
vertex of $\mathcal G$ outside $\mathcal G$, respectively.

For $i\in\{0,\ldots,P-1\}$, assume
\begin{equation}
\label{eq:hyp}
 \begin{gathered}
 f_i\in C^2([0,1]),\qquad f_i(0)=f_i(1)=0,\\
 a_i:=f_i'(0)>0,\qquad f_i'(1)<0,\\
 0<f_i(s)\le a_i s\quad\text{ for } 0<s<1.
 \end{gathered}
\end{equation}
Extend $a_i$ periodically to $\Z$. By \eqref{eq:hyp}, Taylor's theorem
and the finiteness of $\{0,\ldots,P-1\}$, we may choose
\[
 L=\max_{\substack{0\le i<P\\0\le s\le1}}|f_i'(s)|,
 \qquad
 C_2=\frac12\max_{\substack{0\le i<P\\0\le s\le1}}|f_i''(s)|.
\]
Then, for $0\le i<P$ and $s\in[0,1]$,
\begin{equation}
\label{eq:taylor}
 |f_i'(s)|\le L,\qquad a_i s-C_2s^2\le f_i(s)\le a_i s.
\end{equation}
Write $a_{\min}=\min_i a_i$, $a_{\max}=\max_i a_i$, and
$\bar a=P^{-1}\sum_{i=0}^{P-1}a_i$, with both extrema taken over
$i\in\{0,\ldots,P-1\}$.

For $\lambda\in\R$, define the $P\times P$ matrix $M_\lambda$
by its action on $v=(v_i)_{i=0}^{P-1}\in\R^P$:
\begin{equation}
\label{eq:matrix}
 (M_\lambda v)_i
 =De^\lambda v_{i-1}+Dqe^{-\lambda}v_{i+1}
   +(a_i-D(q+1))v_i.
\end{equation}
The component index satisfies $0\le i<P$, and all indices are taken
modulo $P$. If $P=1$ or $P=2$, terms with the same
index are added. Let $k(\lambda)$ be the principal eigenvalue of
$M_\lambda$, characterized by a strictly positive eigenvector
$\phi(\lambda)$ normalized by $\sum_{i=0}^{P-1}\phi_i(\lambda)=1$.
Extend the entries periodically by
$\phi_{j+P}(\lambda)=\phi_j(\lambda)$ for $j\in\Z$.
The proof of Lemma~\ref{lem:matrix} establishes the existence and
simplicity of $k(\lambda)$ and the uniqueness of this normalized
positive eigenvector. Put
\begin{equation}
\label{eq:parameters}
 \lambda_0=\frac12\log q,\qquad \mu=k(\lambda_0),\qquad
 c_+=\min_{\lambda>0}\frac{k(\lambda)}{\lambda},\qquad
 c_-=\max_{\lambda<0}\frac{k(\lambda)}{\lambda}.
\end{equation}
By Lemma~\ref{lem:matrix}, the minimum defining $c_+$ is attained
at a unique $\lambda_+>0$, and the maximum defining $c_-$ is attained
at a unique $\lambda_-<0$. Writing $\operatorname{sgn}$ for the sign
function, the lemma also gives
\[
 k(\lambda)=k(\log q-\lambda)\quad(\lambda\in\R),\qquad
 \mu=\min_{\lambda\in\R}k(\lambda),\qquad
 \operatorname{sgn}c_+=\operatorname{sgn}\mu.
\]
If $a_i=a>0$ for $0\le i<P$, write
\begin{equation}
\label{eq:homogeneous-parameters}
 \begin{aligned}
 \Psi_a(\lambda)&=a+D(e^\lambda+qe^{-\lambda}-(q+1))\quad(\lambda\in\R),\\
 c_*(a)&=\min_{\lambda>0}\frac{\Psi_a(\lambda)}{\lambda},\qquad
 \lambda_*(a)=\operatorname*{argmin}_{\lambda>0}
                   \frac{\Psi_a(\lambda)}{\lambda},\\
 \rho_q&=D(q+1-2\sqrt q).
 \end{aligned}
\end{equation}
For $a_i=a$, the constant vector is a positive eigenvector of
$M_\lambda$, so $k(\lambda)=\Psi_a(\lambda)$ and
$\mu=a-\rho_q$. Throughout, $0\le u_0\le1$. A solution of
\eqref{eq:tree-cauchy} is a function
$u\in C^1([0,\infty);\ell^\infty(\Tq))$ satisfying the initial
condition and the differential equation at each vertex.
Lemma~\ref{lem:comparison} proves existence, uniqueness and
preservation of the interval $[0,1]$. We use the spatial norm
$\|u(t,\cdot)\|_\infty=\sup_{x\in\Tq}|u(t,x)|$.

\subsection{Background and the questions addressed here}

Fisher and Kolmogorov, Petrovskii and Piskunov introduced the
homogeneous reaction--diffusion equation in 1937
\cite{Fisher,KPP}. Aronson and Weinberger established spreading
results in one and several space dimensions
\cite{AronsonWeinberger,AronsonWeinberger1978}.
McKean's representation connects the logistic equation with branching
Brownian motion \cite{McKean}. Bramson identified the logarithmic
displacement and proved convergence to traveling waves
\cite{Bramson1978,Bramson}.
Hamel, Nolen, Roquejoffre and Ryzhik gave a PDE proof of the
logarithmic correction in the homogeneous case in 2013 and treated
periodic media in 2016 \cite{HNRRshort,HNRR}.
For compact perturbations of homogeneous step data, Nolen,
Roquejoffre and Ryzhik proved convergence with a constant phase
correction \cite{NRRsingle}.

The formula $c_+=\min_{\lambda>0}k(\lambda)/\lambda$ in
\eqref{eq:parameters} is related to spreading and traveling waves
in heterogeneous media. Weinberger studied spreading speeds and
traveling waves for order-preserving evolution maps
\cite{Weinberger2002}.
Berestycki, Hamel and Nadirashvili characterized KPP speeds by
periodic elliptic eigenvalue problems \cite{BHN2005}.
Liang and Zhao studied spreading speeds and traveling waves for
monotone semiflows \cite{LiangZhao2007,LiangZhao2010}, and Fang, Yu and
Zhao treated time--space periodic monotone systems \cite{FangYuZhao2017}.
For lattice equations, Zinner, Harris and Hudson proved existence
of traveling fronts \cite{ZinnerHarrisHudson}; Chen and Guo,
and Chen, Fu and Guo, studied their existence, uniqueness and
asymptotic properties \cite{ChenGuo2002,ChenGuo2003,ChenFuGuo}.
Guo and Hamel considered periodic diffusion and reaction
\cite{GuoHamel}. Ma and Zhao studied stability of minimal fronts
\cite{MaZhao2008}; Guo and Wu established uniqueness under prescribed
exponential asymptotics and stability above the minimal speed
\cite{GuoWu}. Huang and Shen, Shen, and Cao and Shen treated
time-dependent coefficients \cite{HuangShen,Shen2010,CaoShen}.
Liang and Zhou developed discrete Harnack estimates and
generalized-principal-eigenvalue methods for spatially heterogeneous
lattice equations \cite{LiangZhou2020}.
For bistable lattice equations with obstacles, entire solutions were
constructed by Hoffman, Hupkes and Van Vleck
\cite{HoffmanHupkesVanVleck2017}.

For the logarithmic displacement in discrete propagation,
Addario-Berry and Reed obtained estimates for minima of branching
random walks \cite{AddarioBerryReed2009}.
Lubetzky, Thornett and Zeitouni proved results for centered maxima
in periodic branching Brownian motion and periodic nearest-neighbor
branching random walks \cite{LubetzkyThornettZeitouni}.
Graham's correction for one-dimensional nonlocal Fisher--KPP
equations includes positive jump rates at $1$ and $-1$
\cite{Graham}. For the homogeneous symmetric lattice, Besse, Faye,
Roquejoffre and Zhang proved the correction and approximation by a
minimal front with a bounded time-dependent phase
\cite{BFRZ}. Boutillon derived a logarithmic delay for a linearized
nonlocal equation and a corresponding nonlinear upper estimate
\cite{Boutillon2024}. Roquejoffre gives a systematic treatment of
one-dimensional nonlocal front propagation \cite{Roquejoffre2024}.
For smooth compactly supported even dispersal kernels and smooth
nonzero compactly supported initial data, Roquejoffre and Zhang
obtained directional level-position expansions through
a time-independent directional shift and an $o(1)$ remainder
\cite{RoquejoffreZhang2026}. The assumptions on their homogeneous reaction $f$ in
\cite[(1.2)--(1.3)]{RoquejoffreZhang2026} include
$f'(s)\le f'(0)$ for $0<s<1$, and
$f(s)\ge f'(0)s-Ms^2$ for $0<s\le s_0$, with
$M>0$ and $s_0\in(0,1)$.

Geometry enters the persistence question through the spectrum of
the diffusion operator. Matano, Punzo and Tesei studied propagation
and extinction for nonlinear diffusion on hyperbolic space
\cite{MatanoPunzoTesei}. For homogeneous trees and an initial
population at one vertex, Hoffman and Holzer identified the
diffusion threshold and the linearly determined spreading speed
\cite{HoffmanHolzer}. Gonz\'alez, Gonz\'alvez and Quir\'os treated
critical extinction and traveling-wave asymptotics in hyperbolic
space for three classes of symmetry \cite{GGQ}.
For the constant logistic reaction $f(u)=au(1-u)$, $a>0$, on a homogeneous tree,
\'Alvarez, Murillo-Arcila, Ponce and Pozo proved the logarithmic
correction for nonzero finitely supported radial data when
$a>D(q+1-2\sqrt q)$ \cite{AlvarezMurilloPoncePozo2026}.
The spectral theory of periodic Jacobi operators on trees is
developed in \cite{AvniBreuerSimon,BanksBreuerGarzaSeeligSimon}.

The reductions for $f(x,u)=f_{r(x)}(u)$ and
$f(x,u)=f_{s(x)}(u)$ differ at the index zero.
For the former, \eqref{eq:radial} gives
$\dot U_0=D(q+1)(U_1-U_0)+f_0(U_0)$; for the latter,
$\dot V_0=D(q-1)(V_1-V_0)+f_0(V_0)$.
How does this difference affect extinction and persistence?
Can a solution remain small at large $s(x)$ while positive levels
propagate in $\ell_{\mathcal G}(x)$?

\subsection{Main results}\label{sec:intro-main-results}

Unless a more specific dependence is stated, constants in the
asymptotic estimates may depend on $D$, $q$, $P$,
$(f_i)_{i=0}^{P-1}$, the initial datum, and the level.

\paragraph{Extinction thresholds.}
For nonzero finitely supported data, both assignments
$f(x,u)=f_{b(x)}(u)$ and $f(x,u)=f_{r(x)}(u)$ have the extinction
threshold $\mu=0$.

\begin{theorem}\label{thm:threshold}
Assume \eqref{eq:hyp} and $f(x,s)=f_{b(x)}(s)$.
Let $u_0:\Tq\to[0,1]$ be nonzero and finitely supported,
and let $u$ solve \eqref{eq:tree-cauchy} with initial datum $u_0$.
Set $\mu=k(\tfrac12\log q)$.
\begin{enumerate}
\item[\textup{(i)}] \textbf{Uniform extinction: $\mu\le0$.}
We have
\[
 \lim_{t\to\infty}\norm{u(t,\cdot)}_\infty=0.
\]
If $\mu<0$, then
\begin{equation}
\label{eq:exp-extinction}
 \norm{u(t,\cdot)}_\infty\le e^{\mu t}\norm{u_0}_2
 \quad\text{ for }t\ge0.
\end{equation}
\item[\textup{(ii)}] \textbf{Persistence: $\mu>0$.}
For every fixed vertex $x\in\Tq$,
\[
 \lim_{t\to\infty}u(t,x)=1.
\]
\end{enumerate}
\end{theorem}

\begin{remark*}
Under the assumptions on $f$ in Theorem~\ref{thm:threshold},
let $u$ solve \eqref{eq:tree-cauchy} with initial datum
$u_0\in\ell^2(\Tq)$ satisfying $0\le u_0\le1$.
Then $\|u(t,\cdot)\|_\infty\to0$ as $t\to\infty$ when $\mu\le0$.
For these initial data, estimate \eqref{eq:exp-extinction} also holds
when $\mu<0$.
\end{remark*}

\Needspace{16\baselineskip}
\begin{theorem}\label{thm:radial-threshold}
Assume \eqref{eq:hyp} and $f(x,s)=f_{r(x)}(s)$.
Let $u_0:\Tq\to[0,1]$ be nonzero and finitely supported,
and let $u$ solve \eqref{eq:tree-cauchy} with initial datum $u_0$.
\begin{enumerate}
\item[\textup{(i)}] \textbf{Uniform extinction: $\mu\le0$.}
We have
\[
 \lim_{t\to\infty}\|u(t,\cdot)\|_\infty=0.
\]
If $\mu<0$, then
\[
 \|u(t,\cdot)\|_\infty\le e^{\mu t}\|u_0\|_2
 \quad\text{ for }t\ge0.
\]
\item[\textup{(ii)}] \textbf{Persistence: $\mu>0$.}
For every fixed vertex $x\in\Tq$,
\[
 \lim_{t\to\infty}u(t,x)=1.
\]
\end{enumerate}
\end{theorem}

\begin{remark*}
Under the assumptions on $f$ in Theorem~\ref{thm:radial-threshold},
let $u$ solve \eqref{eq:tree-cauchy} with initial datum
$u_0\in\ell^2(\Tq)$ satisfying $0\le u_0\le1$.
Then $\|u(t,\cdot)\|_\infty\to0$ as $t\to\infty$ when $\mu\le0$.
For these initial data and $\mu<0$, we also have
\[
 \|u(t,\cdot)\|_\infty\le e^{\mu t}\|u_0\|_2
 \quad\text{ for }t\ge0.
\]
\end{remark*}

For the reaction assignment $f(x,u)=f_{s(x)}(u)$, define the bounded self-adjoint
tridiagonal operator $\mathcal J_{\mathcal G}$ on $\ell^2(\N_0)$ by
its action on $v=(v_r)_{r\in\N_0}\in\ell^2(\N_0)$:
\[
 \begin{aligned}
 (\mathcal J_{\mathcal G}v)_0
   &=(a_0-D(q-1))v_0+D\sqrt{q-1}\,v_1,\\
 (\mathcal J_{\mathcal G}v)_1
   &=D\sqrt{q-1}\,v_0+(a_1-D(q+1))v_1+D\sqrt q\,v_2,\\
 (\mathcal J_{\mathcal G}v)_r
   &=D\sqrt q\,(v_{r-1}+v_{r+1})+(a_r-D(q+1))v_r,
   \quad\text{ for }r\ge2.
 \end{aligned}
\]
For an operator $\mathcal A$ on a specified sequence space,
$\spec(\mathcal A)$ denotes its spectrum. Set
$\sigma_{\mathcal G}=\sup\spec(\mathcal J_{\mathcal G})$.
By Lemma~\ref{lem:geodesic-spectrum}, $\sigma_{\mathcal G}$ is also
the upper spectral endpoint of $D\Delta_{\Tq}+a_{s(x)}$ on
$\ell^2(\Tq)$, and $\sigma_{\mathcal G}\ge\mu$.
Theorem~\ref{thm:geodesic-persistence} gives an example with
$\mu<0<\sigma_{\mathcal G}$.

\begin{theorem}\label{thm:geodesic-threshold}
Assume \eqref{eq:hyp} and $f(x,u)=f_{s(x)}(u)$.
For each initial datum $u_0:\Tq\to[0,1]$, let $u$ solve
\eqref{eq:tree-cauchy} with initial datum $u_0$.
\begin{enumerate}
\item[\textup{(i)}] \textbf{Uniform extinction: $\sigma_{\mathcal G}\le0$.}
For every initial datum $u_0:\Tq\to[0,1]$ in $\ell^2(\Tq)$,
\[
 \lim_{t\to\infty}\|u(t,\cdot)\|_\infty=0.
\]
If $\sigma_{\mathcal G}<0$, then
\[
 \|u(t,\cdot)\|_\infty
 \le e^{\sigma_{\mathcal G}t}\|u_0\|_2
 \quad\text{ for }t\ge0.
\]
\item[\textup{(ii)}] \textbf{Persistence: $\sigma_{\mathcal G}>0$.}
Among the stationary solutions $p:\Tq\to(0,1]$ of the differential
equation in \eqref{eq:tree-cauchy}, there is a pointwise smallest
solution $p_*$. For every $u_0\not\equiv0$,
\begin{equation}
\label{eq:geodesic-minimal-persistence}
 \liminf_{t\to\infty}u(t,x)\ge p_*(x)>0
 \quad\text{for every fixed }x\in\Tq.
\end{equation}
For every $\varepsilon>0$, there is a finitely supported initial datum
$u_0:\Tq\to[0,1]$ with $0<\|u_0\|_\infty\le\varepsilon$ whose
solution increases pointwise to $p_*$ as $t\to\infty$.
\item[\textup{(iii)}] \textbf{Stationary solution: $\sigma_{\mathcal G}>0$.}
The function $p_*$ depends only on $s(x)$.
If $\mu>0$, then $p_*\equiv1$.
If $\mu\le0<\sigma_{\mathcal G}$, then
\begin{equation}
\label{eq:geodesic-minimal-tail}
 0<p_*(x)<1,\qquad
 p_*(x)\le q^{-s(x)/2}
 \frac{\phi_{s(x)}(\lambda_0)}{\phi_0(\lambda_0)}
 \quad\text{ for }x\in\Tq.
\end{equation}
\item[\textup{(iv)}] \textbf{Transverse bound: $\mu\le0$.}
For every finitely supported $u_0$, there exists $C<\infty$ such that
\begin{equation}
\label{eq:geodesic-transverse-confinement}
 u(t,x)
 \le\min\{1,Cq^{-s(x)/2}\phi_{s(x)}(\lambda_0)\}
 \quad\text{ for }t\ge0,\ x\in\Tq.
\end{equation}
\end{enumerate}
\end{theorem}

The following logistic example satisfies $\mu<0<\sigma_{\mathcal G}$.

\begin{theorem}\label{thm:geodesic-persistence}
Let $q=2$, $D=1$, $P=8$, and
\[
 a_0=\frac34,\qquad
 a_i=\frac1{100}\quad\text{ for }1\le i\le7,
 \qquad
 f(x,u)=a_{s(x)}u(1-u)\quad\text{for }x\in\Tq,\ u\in[0,1].
\]
Then the following conclusions hold.
\begin{enumerate}
\item[\textup{(i)}] 
We have $\mu<0<\sigma_{\mathcal G}$.

\item[\textup{(ii)}] 
For every $\varepsilon>0$, there is a finitely supported initial datum
$u_0:\Tq\to[0,1]$ with $0<\|u_0\|_\infty\le\varepsilon$ such that
the solution of \eqref{eq:tree-cauchy} with initial datum $u_0$
increases pointwise as $t\to\infty$ to the smallest positive
stationary solution $p=p_*$ of the differential equation in
\eqref{eq:tree-cauchy}. The function $p$ satisfies
\begin{equation}
\label{eq:localized-steady}
 0<p(x)<1\quad\text{ for }x\in\Tq,
 \qquad
 \lim_{R\to\infty}\sup_{s(x)\ge R}p(x)=0.
\end{equation}

\item[\textup{(iii)}] 
For every nonzero finitely supported initial datum
$u_0\colon\Tq\to[0,1]$, the solution $u$ of
\eqref{eq:tree-cauchy} with initial datum $u_0$ satisfies
\[
 \liminf_{t\to\infty}u(t,x)>0
 \quad\text{for every fixed }x\in\Tq.
\]
\end{enumerate}
\end{theorem}

\paragraph{Logarithmic level positions.}
For $t\ge0$, define
\[
 m_*(t)
 =
 c_+t-\frac{3}{2\lambda_+}\log(1+t).
\]
The logarithmic term in $m_*(t)$ is the Bramson correction.

Let $\theta\in(0,1)$. For a solution $w$ of \eqref{eq:lattice}, set
$H_\theta(t)=\max\{j\in\Z:w_j(t)\ge\theta\}$ whenever the set
is nonempty and bounded above. For a solution $u$ of
\eqref{eq:tree-cauchy}, define
\[
\begin{gathered}
 H_\theta^b(t)=\max\{b(x):x\in\Tq,\ u(t,x)\ge\theta\},\\
 R_\theta^o(t)=\max\{r(x):x\in\Tq,\ u(t,x)\ge\theta\},\\
 R_\theta^{\mathcal G}(t)=\max\{s(x):x\in\Tq,\ u(t,x)\ge\theta\}
 \quad\text{ for }t\ge0,
\end{gathered}
\]
whenever the corresponding set of coordinate values is nonempty and
bounded above. Under the respective hypotheses of
Theorems~\ref{thm:periodic-bramson-step},
\ref{thm:horo-compact-bramson}, and~\ref{thm:periodic-bramson-distance},
each maximum exists for all sufficiently large $t$.
Indeed, the coordinate values are integers, so a nonempty set bounded
above contains its maximum.

\begin{theorem}\label{thm:periodic-bramson-step}
Assume \eqref{eq:hyp}, and fix integers $A\le B$.
The following conclusions hold for every $\theta\in(0,1)$.
\begin{enumerate}
\item[\textup{(i)}] \textbf{The lattice equation.}
Let $w$ solve \eqref{eq:lattice} with
\[
 \mathbf1_{\{j\le A\}}\le w_j(0)\le\mathbf1_{\{j\le B\}}
 \quad\text{ for }j\in\Z.
\]
There are $C_\theta,T_\theta<\infty$ such that
\begin{equation}
\label{eq:periodic-bramson-step}
 \{j:j\le m_*(t)-C_\theta\}
 \subseteq
 \{j:w_j(t)\ge\theta\}
 \subseteq
 \{j:j\le m_*(t)+C_\theta\}
 \quad\text{ for }t\ge T_\theta.
\end{equation}
Consequently,
\[
 H_\theta(t)
 =
 c_+t-\frac{3}{2\lambda_+}\log(1+t)
 +O_{\theta,A,B}(1)
 \quad\text{ as }t\to\infty.
\]

\item[\textup{(ii)}] \textbf{The tree equation.}
Let $u$ solve \eqref{eq:tree-cauchy} with
$f(x,u)=f_{b(x)}(u)$ and
\[
 \mathbf1_{\{b\le A\}}
 \le u_0
 \le\mathbf1_{\{b\le B\}}\quad\text{on }\Tq.
\]
There are $C_\theta,T_\theta<\infty$ such that
\[
 \{x:b(x)\le m_*(t)-C_\theta\}
 \subseteq\{x:u(t,x)\ge\theta\}
 \subseteq\{x:b(x)\le m_*(t)+C_\theta\}
 \quad\text{ for }t\ge T_\theta.
\]
\end{enumerate}
Both conclusions hold for $c_+>0$, $c_+=0$ and $c_+<0$ under \eqref{eq:hyp}.
\end{theorem}

\begin{theorem}\label{thm:horo-compact-bramson}
Assume \eqref{eq:hyp} and $\mu>0$, and let
$f(x,s)=f_{b(x)}(s)$ for $x\in\Tq$ and $s\in[0,1]$.
Let $u_0:\Tq\to[0,1]$ be nonzero and finitely supported,
and let $u$ solve \eqref{eq:tree-cauchy} with initial datum $u_0$.
The following conclusions hold for every $\theta\in(0,1)$.
\begin{enumerate}
\item[\textup{(i)}] 
The set $\{x\in\Tq:u(t,x)\ge\theta\}$ is finite for every
$t\ge0$ and nonempty for all sufficiently large $t$.

\item[\textup{(ii)}] 
There are $0\le C_\theta,T_\theta<\infty$, depending only on
$\theta$, $u_0$, $D$, $q$, $P$, and $(f_i)_{i=0}^{P-1}$, such that
\begin{equation}\label{eq:horo-compact-bramson}
 \left|H_\theta^b(t)-c_+t+
 \frac{3}{2\lambda_+}\log(1+t)\right|\le C_\theta
 \quad\text{ for }t\ge T_\theta.
\end{equation}
\end{enumerate}
\end{theorem}

\begin{theorem}\label{thm:periodic-bramson-distance}
Assume \eqref{eq:hyp} and $\mu>0$, equivalently $c_+>0$.
Let $u$ solve \eqref{eq:tree-cauchy} with
$u_0:\Tq\to[0,1]$.
\begin{enumerate}
\item[\textup{(i)}] \textbf{Distance from $o$.}
Assume $f(x,u)=f_{r(x)}(u)$, and let $u_0$ be nonzero and
finitely supported. For every $\theta\in(0,1)$, there are
$C_\theta,T_\theta<\infty$ such that
\begin{equation}
\label{eq:periodic-bramson-balls}
 B(o,m_*(t)-C_\theta)
 \subseteq
 \{x:u(t,x)\ge\theta\}
 \subseteq
 B(o,m_*(t)+C_\theta)
 \quad\text{ for }t\ge T_\theta.
\end{equation}
Consequently,
\[
 R_\theta^o(t)=m_*(t)+O(1)
 \quad\text{ as }t\to\infty.
\]

\item[\textup{(ii)}] \textbf{Distance from $\mathcal G$.}
Assume $f(x,u)=f_{s(x)}(u)$ and
\[
 \eta\mathbf1_{\mathcal N_{R_1}(\mathcal G)}
 \le u_0
 \le\mathbf1_{\mathcal N_{R_2}(\mathcal G)},
 \qquad
 0<\eta\le1,\quad 0\le R_1\le R_2<\infty.
\]
For every $\theta\in(0,1)$, there are
$C_\theta,T_\theta<\infty$ such that
\begin{equation}
\label{eq:periodic-bramson-tubes}
 \mathcal N_{m_*(t)-C_\theta}(\mathcal G)
 \subseteq
 \{x:u(t,x)\ge\theta\}
 \subseteq
 \mathcal N_{m_*(t)+C_\theta}(\mathcal G)
 \quad\text{ for }t\ge T_\theta.
\end{equation}
Consequently,
\[
 R_\theta^{\mathcal G}(t)=m_*(t)+O(1)
 \quad\text{ as }t\to\infty.
\]
\end{enumerate}
\end{theorem}

In the geodesic assertion of
Theorem~\ref{thm:periodic-bramson-distance}, the condition $R_1\ge0$
gives $\mathcal G\subseteq\mathcal N_{R_1}(\mathcal G)$ and hence
$u_0(x)\ge\eta$ for every $x\in\mathcal G$. For
$f(x,u)=f_{s(x)}(u)$, $\mu\le0<\sigma_{\mathcal G}$ and nonzero
finitely supported $u_0:\Tq\to[0,1]$,
Theorems~\ref{thm:tangential-spreading}
and~\ref{thm:tangential-bramson} give the propagation speed and
logarithmic level positions along $\mathcal G$.

\paragraph{A common derivative at zero.}
When $f_i'(0)=a$ for $0\le i<P$, the reactions may be assigned
by any map from $\Tq$ to $\{0,\ldots,P-1\}$.

\begin{theorem}\label{thm:common}
Assume \eqref{eq:hyp} and $f_i'(0)=a$ for $0\le i<P$.
Let $\iota:\Tq\to\{0,\ldots,P-1\}$ be arbitrary, and set
$f(x,s)=f_{\iota(x)}(s)$.
Let $u$ solve \eqref{eq:tree-cauchy} with
$u_0:\Tq\to[0,1]$.
\begin{enumerate}
\item[\textup{(i)}] \textbf{Uniform extinction: $a\le\rho_q$.}
Assume $u_0\in\ell^2(\Tq)$. Then
\[
 \norm{u(t,\cdot)}_\infty\to0
 \quad\text{ as }t\to\infty.
\]
If $a<\rho_q$, then
\[
 \norm{u(t,\cdot)}_\infty
 \le
 e^{(a-\rho_q)t}\norm{u_0}_2
 \quad\text{ for }t\ge0.
\]

\item[\textup{(ii)}] \textbf{Spreading: $a>\rho_q$.}
If $u_0$ is nonzero and finitely supported, then for every
$\theta\in(0,1)$ there are $C_\theta,T_\theta>0$ such that
\begin{equation}
\label{eq:common-balls}
 \begin{aligned}
 B\left(o,c_*(a)t-\frac3{2\lambda_*(a)}\log(1+t)-C_\theta\right)
 &\subseteq \{x:u(t,x)\ge\theta\}\\
 &\subseteq
 B\left(o,c_*(a)t-\frac3{2\lambda_*(a)}\log(1+t)+C_\theta\right)
 \end{aligned}
\end{equation}
for $t\ge T_\theta$.
\end{enumerate}
\end{theorem}

\begin{corollary}\label{cor:common-step}
Assume \eqref{eq:hyp} and $f_i'(0)=a>0$ for $0\le i<P$. Fix
integers $A\le B$. Consider the following two initial-value problems.
\begin{enumerate}
\item[\textup{(i)}] Let $w$ solve \eqref{eq:lattice} with
$\ind_{\{j\le A\}}\le w_j(0)\le\ind_{\{j\le B\}}$ for $j\in\Z$,
and put $H_\theta(t)=\max\{j:w_j(t)\ge\theta\}$.
\item[\textup{(ii)}] Let $u$ solve \eqref{eq:tree-cauchy} with
$f(x,s)=f_{b(x)}(s)$ and
$\ind_{\{b\le A\}}\le u_0\le\ind_{\{b\le B\}}$,
and put $H_\theta(t)=\max\{b(x):u(t,x)\ge\theta\}$.
\end{enumerate}
In each case, for every $\theta\in(0,1)$, the maximum defining
$H_\theta(t)$ exists for all sufficiently large $t$, and
\begin{equation}
\label{eq:common-step}
 H_\theta(t)
 =
 c_*(a)t-\frac3{2\lambda_*(a)}\log t
 +O_{\theta,A,B}(1)\quad\text{ as } t\to\infty.
\end{equation}
In addition, equation \eqref{eq:lattice} admits an entire solution
$W(t,j)$, defined
for $(t,j)\in\R\times\Z$, with $0<W<1$, $W(t,j+P)\le W(t,j)$,
and, for every $\theta\in(0,1)$, a constant
$C_\theta$ such that
\begin{equation}
\label{eq:entire-levels}
 \{j:j\le c_*(a)t-C_\theta\}
 \subseteq\{j:W(t,j)\ge\theta\}
 \subseteq\{j:j\le c_*(a)t+C_\theta\}
 \quad\text{ for } t\in\R.
\end{equation}
{The function $W$ is the limit of a sequence of space and time
translates of the solution of \eqref{eq:lattice} with
$w_j(0)=\ind_{\{j\le0\}}$; convergence holds in}
$\ell^\infty(\Z)$ uniformly on bounded time intervals.
\end{corollary}

\paragraph{Moving and stationary fronts.}
For $c\in\R\setminus\{0\}$, a solution $w$ of
\eqref{eq:lattice} defined for $(t,j)\in\R\times\Z$ is a
pulsating front of speed $c$ if $w_j(t)=\Phi_j(j-ct)$ for functions
$\Phi_j:\R\to(0,1)$, $j\in\Z$, satisfying $\Phi_{j+P}=\Phi_j$ and
$\lim_{z\to-\infty}\max_{0\le i<P}|\Phi_i(z)-1|=0$ and
$\lim_{z\to+\infty}\max_{0\le i<P}\Phi_i(z)=0$. A stationary front is a stationary solution
$p=(p_j)_{j\in\Z}$ of \eqref{eq:lattice} satisfying
$0<p_j<1$, $p_{j+P}<p_j$, $p_j\to1$ as $j\to-\infty$, and
$p_j\to0$ as $j\to+\infty$.

\begin{theorem}\label{thm:all-speeds}

Assume \eqref{eq:hyp} and $f(x,u)=f_{b(x)}(u)$.
\begin{enumerate}
\item[\textup{(i)}] 
For every $c\in[c_+,\infty)\setminus\{0\}$, equation
\eqref{eq:lattice} has a pulsating front $w_j(t)=\Phi_j(j-ct)$
with each $\Phi_j$ strictly decreasing. Its lift
\[
 u(t,x)=w_{b(x)}(t)=\Phi_{b(x)}(b(x)-ct)\quad\text{for }(t,x)\in\R\times\Tq
\]
solves the differential equation in \eqref{eq:tree-cauchy}.

\item[\textup{(ii)}] 
Equation \eqref{eq:lattice} has a stationary front $p=(p_j)_{j\in\Z}$
exactly when $\mu\le0$; its lift $x\mapsto p_{b(x)}$ is a
stationary solution of the differential equation in
\eqref{eq:tree-cauchy}.

\item[\textup{(iii)}] 
For the orientation with profile limits one at $-\infty$ and
zero at $+\infty$, the speeds of pulsating fronts of
\eqref{eq:lattice}, together with zero when a stationary front
exists, form exactly $[c_+,\infty)$.
\end{enumerate}

\end{theorem}

Theorems~\ref{thm:wave-stability} and~\ref{thm:tree-wave-stability}
concern a fixed front from Theorem~\ref{thm:all-speeds} and assume
\begin{equation}
\label{eq:derivative-kpp}
 f_i'(s)\le a_i
 \quad\text{ for } 0\le s\le1,\ 0\le i<P.
\end{equation}

\begin{theorem}\label{thm:wave-stability}
Assume \eqref{eq:hyp}, \eqref{eq:derivative-kpp}, and $\mu>0$.
Let $u=(u_j)$ solve \eqref{eq:lattice}, and fix a front
$W_j(t)=\Phi_j(j-ct)$ from Theorem~\ref{thm:all-speeds}
with $c\ge c_+$. Suppose that, for some $A\ge0$,
\[
 \Phi_j(j+A)\le u_j(0)\le\Phi_j(j-A)
 \quad\text{ for }j\in\Z,
\]
and that, for some $1\le p<\infty$,
\[
 B_p
 =\|e^{\lambda_+j}(u_j(0)-W_j(0))\|_{\ell^p(\Z)}
 <\infty.
\]
Then the following conclusions hold.
\begin{enumerate}
\item[\textup{(i)}] \textbf{Uniform convergence.}
\begin{equation}
\label{eq:wave-stability}
 \sup_{j\in\Z}|u_j(t)-W_j(t)|\to0
 \quad\text{ as }t\to\infty.
\end{equation}

\item[\textup{(ii)}] \textbf{Weighted estimate.}
There is a constant $C_p<\infty$ such that
\begin{equation}
\label{eq:wave-stability-weighted}
 \sup_{j\in\Z}
 e^{\lambda_+(j-ct)}|u_j(t)-W_j(t)|
 \le
 C_p(1+t)^{-1/(2p)}
 e^{-\lambda_+(c-c_+)t}B_p
 \quad\text{ for }t\ge0.
\end{equation}
\end{enumerate}
\end{theorem}

\begin{theorem}\label{thm:tree-wave-stability}
Assume \eqref{eq:hyp}, \eqref{eq:derivative-kpp}, and
$f(x,u)=f_{b(x)}(u)$. Let $w$ be a front of \eqref{eq:lattice} from
Theorem~\ref{thm:all-speeds} with speed $c$, where $c=0$
is allowed for a stationary front, and set $W(t,x)=w_{b(x)}(t)$.
Let $u$ solve \eqref{eq:tree-cauchy} with initial datum $u_0$.
Suppose that, for some integer $m\ge0$,
\begin{equation}
\label{eq:tree-front-bracket}
 w_{b(x)+mP}(0)\le u_0(x)\le w_{b(x)-mP}(0)
 \quad\text{ for }x\in\Tq,
\end{equation}
and that, for some $1\le p<\infty$,
\[
 B_p
 =\|e^{\lambda_+b(x)}(u_0(x)-W(0,x))\|_{\ell^p(\Tq)}
 <\infty.
\]
Then the following conclusions hold.
\begin{enumerate}
\item[\textup{(i)}] \textbf{Uniform convergence.}
\begin{equation}
\label{eq:tree-front-stability}
 \|u(t,\cdot)-W(t,\cdot)\|_{\ell^\infty(\Tq)}\to0
 \quad\text{ as }t\to\infty.
\end{equation}

\item[\textup{(ii)}] \textbf{Weighted estimate.}
There is a constant $C_p<\infty$ such that
\[
 \sup_{x\in\Tq}
 e^{\lambda_+(b(x)-ct)}|u(t,x)-W(t,x)|
 \le
 C_p(1+t)^{-1/(2p)}
 e^{-\lambda_+(c-c_+)t}B_p
 \quad\text{ for }t\ge0.
\]
\end{enumerate}
\end{theorem}

\paragraph{Positive stationary solutions near a geodesic.}
For uniqueness and convergence when $\mu\le0<\sigma_{\mathcal G}$,
we impose the additional condition
\begin{equation}
\label{eq:strict-quotient}
 \frac{f_i(r)}r>\frac{f_i(s)}s
 \quad\text{ for } 0<r<s\le1,\ 0\le i<P.
\end{equation}

\Needspace{16\baselineskip}
\begin{theorem}\label{thm:geodesic-critical-convergence}
Assume \eqref{eq:hyp}, \eqref{eq:strict-quotient},
$f(x,u)=f_{s(x)}(u)$, and $\mu\le0<\sigma_{\mathcal G}$.
Set $h_r=q^{-r/2}\phi_r(\lambda_0)$ for $r\in\N_0$.
Let $p_*$ be the smallest positive stationary solution from
Theorem~\ref{thm:geodesic-threshold}, and write $p_{*,r}=p_*(x)$
when $s(x)=r$. Then the following conclusions hold.
\begin{enumerate}
\item[\textup{(i)}] 
The function $p_*$ is the only stationary solution $p$ of the
differential equation in \eqref{eq:tree-cauchy} satisfying
\[
 0<p(x)\le1\quad\text{ for }x\in\Tq,
 \qquad
 \sup_{x\in\Tq}\frac{p(x)}{h_{s(x)}}<\infty.
\]

\item[\textup{(ii)}] 
Let $u$ solve \eqref{eq:tree-cauchy} with a nonzero initial datum
$u_0\colon\Tq\to[0,1]$ satisfying
$u_0(x)\le C_0h_{s(x)}$ for every $x\in\Tq$ and some $C_0<\infty$.
For every $0\le M<\infty$,
\begin{equation}
\label{eq:critical-geodesic-convergence}
 \lim_{t\to\infty}
 \sup_{\{x:|\ell_{\mathcal G}(x)|\le M\}}
 |u(t,x)-p_*(x)|
 =0.
\end{equation}
In particular, \eqref{eq:critical-geodesic-convergence} holds for
every nonzero finitely supported $u_0\colon\Tq\to[0,1]$.

\item[\textup{(iii)}] 
If $\mu=0$, there is a constant $A_*>0$ such that
\[
 p_{*,r}
 =A_*q^{-r/2}\phi_r(\lambda_0)+O(q^{-r})
 \quad\text{ as }r\to\infty.
\]
\end{enumerate}
\end{theorem}

For Theorem~\ref{thm:geodesic-convergence}, define $\nu_0=1$ and
$\nu_r=(q-1)q^{r-1}$ for $r\ge1$, and set
\[
 \mathcal H_\nu
 =\left\{v=(v_r)_{r\ge0}:v_r\in\R,\
                 \sum_{r\ge0}\nu_rv_r^2<\infty\right\}.
\]

\begin{theorem}\label{thm:geodesic-convergence}
Assume \eqref{eq:hyp}, \eqref{eq:strict-quotient},
$f(x,u)=f_{s(x)}(u)$, and $\mu<0<\sigma_{\mathcal G}$.
Let $\lambda_{\mathrm R}$ be the larger positive root of $k(\lambda)=0$,
and fix $\lambda\in(\lambda_0,\lambda_{\mathrm R})$.
Let $p_*$ be the smallest positive stationary solution from
Theorem~\ref{thm:geodesic-threshold}, and write
$p_{*,r}=p_*(x)$ when $s(x)=r$. Then the following conclusions hold.
\begin{enumerate}
\item[\textup{(i)}] 
The function $p_*$ is the only stationary solution
$p:\Tq\to(0,1]$ of the differential equation in
\eqref{eq:tree-cauchy} satisfying
\[
 \sup_{x\in\Tq} e^{\lambda s(x)}p(x)<\infty.
\]

\item[\textup{(ii)}] 
The sequence $(p_{*,r})_{r\ge0}$ belongs to $\mathcal H_\nu$
and is the only sequence $(p_r)_{r\ge0}\in\mathcal H_\nu$
with $0<p_r\le1$ for $r\in\N_0$ such that
$x\mapsto p_{s(x)}$ is a stationary solution of the differential
equation in \eqref{eq:tree-cauchy}.

\item[\textup{(iii)}] 
Let $u$ solve \eqref{eq:tree-cauchy} with a nonzero initial datum
$u_0:\Tq\to[0,1]$. If there is $C_0<\infty$ such that
$u_0(x)\le C_0e^{-\lambda s(x)}$ for every $x\in\Tq$, then
\begin{equation}
 \label{eq:geodesic-local-convergence}
 \lim_{t\to\infty}
 \sup_{\{x:|\ell_{\mathcal G}(x)|\le M\}}
           |u(t,x)-p_*(x)|=0\quad\text{ for } 0\le M<\infty.
\end{equation}

\item[\textup{(iv)}] 
Let $u$ solve \eqref{eq:tree-cauchy} with a nonzero finitely
supported datum $u_0:\Tq\to[0,1]$. Then $u(t,x)\to p_*(x)$
as $t\to\infty$ for every fixed $x\in\Tq$, and
\begin{equation}
 \label{eq:no-global-geodesic-convergence}
 \|u(t,\cdot)-p_*(\cdot)\|_{\ell^\infty(\Tq)}
       \ge p_*(\gamma(0))>0\quad\text{ for }t\ge0.
\end{equation}
\end{enumerate}
\end{theorem}

\begin{theorem}\label{thm:stationary-family}
Assume \eqref{eq:hyp}, \eqref{eq:strict-quotient},
$f(x,u)=f_{s(x)}(u)$ and $\mu\le0<\sigma_{\mathcal G}$.
Let $p_*$ be the smallest positive stationary solution from
Theorem~\ref{thm:geodesic-threshold}. Write $p_{*,r}=p_*(x)$
when $s(x)=r$, and put $h_r=q^{-r/2}\phi_r(\lambda_0)$
for $r\in\N_0$. Then the following conclusions hold.
\begin{enumerate}
\item[\textup{(i)}] 
For each $\alpha\in[p_{*,0},1]$, there is exactly one sequence
$U^\alpha=(U_r^\alpha)_{r\ge0}$ such that
$0<U_r^\alpha\le1$ for $r\in\N_0$, $U_0^\alpha=\alpha$, and
$x\mapsto U_{s(x)}^\alpha$ is a stationary solution of the
differential equation in \eqref{eq:tree-cauchy}.
For $r\in\N_0$, the endpoint values are
$U_r^{p_{*,0}}=p_{*,r}$ and $U_r^1=1$.

\item[\textup{(ii)}] 
For each $\alpha\in(p_{*,0},1)$,
\begin{equation}\label{eq:stationary-family-intro}
 p_{*,r}<U_r^\alpha<1,\qquad
 U_{r+1}^\alpha<U_r^\alpha
 \quad\text{ for }r\in\N_0,
\end{equation}
and
\begin{equation}\label{eq:stationary-family-decay-intro}
 \lim_{r\to\infty}U_r^\alpha=0,\qquad
 \sup_{r\ge0}\frac{U_r^\alpha}{h_r}=\infty.
\end{equation}

\item[\textup{(iii)}] 
If $p:\Tq\to(0,1]$ is a stationary solution of the differential
equation in \eqref{eq:tree-cauchy} and depends only on $s(x)$,
then $p(x)=U_{s(x)}^\alpha$ for every $x\in\Tq$,
where $\alpha=p(\gamma(0))\in[p_{*,0},1]$.
\end{enumerate}
\end{theorem}

To determine propagation in the coordinate $\ell_{\mathcal G}$,
assume \eqref{eq:hyp} and $\mu\le0<\sigma_{\mathcal G}$.
For $\lambda\in\R$, define the bounded self-adjoint operator
$\mathcal J_{\parallel}(\lambda)$ on $\ell^2(\N_0)$ and its upper
spectral endpoint $\Lambda_{\parallel}(\lambda)$ by
\[
 \mathcal J_{\parallel}(\lambda)
 =
 \mathcal J_{\mathcal G}
 +
 2D(\cosh\lambda-1)
 \langle\delta_0,\,\cdot\,\rangle\delta_0,
 \qquad
 \Lambda_{\parallel}(\lambda)
 =
 \sup\spec(\mathcal J_{\parallel}(\lambda)).
\]
Set
\begin{equation}
\label{eq:tangential-speed-intro}
 c_{\parallel}
 =
 \min_{\lambda>0}
 \frac{\Lambda_{\parallel}(\lambda)}{\lambda}.
\end{equation}
The additional diagonal term in $\mathcal J_{\parallel}(\lambda)$
comes from the two neighbors on $\mathcal G$; edges outside
$\mathcal G$ preserve $\ell_{\mathcal G}$.
By Lemma~\ref{lem:tangential-dispersion}, the minimum in
\eqref{eq:tangential-speed-intro} is attained at a unique
$\lambda_{\parallel}>0$.

\begin{theorem}\label{thm:tangential-spreading}

Assume \eqref{eq:hyp}, $f(x,u)=f_{s(x)}(u)$, and
$\mu\le0<\sigma_{\mathcal G}$.
Let $u$ solve \eqref{eq:tree-cauchy} with a nonzero finitely
supported initial datum $u_0\colon\Tq\to[0,1]$, and let $p_*$
be the smallest positive stationary solution from
Theorem~\ref{thm:geodesic-threshold}. Then $c_{\parallel}>0$,
and the following conclusions hold.
\begin{enumerate}
\item[\textup{(i)}] 
For every $c>c_{\parallel}$,
\begin{equation}
\label{eq:tangential-outer}
 \lim_{t\to\infty}
 \sup_{\{x:|\ell_{\mathcal G}(x)|\ge ct\}}
 u(t,x)=0.
\end{equation}

\item[\textup{(ii)}] 
For every $0<\theta<p_*(\gamma(0))$, the set
$\{n\in\Z:u(t,\gamma(n))\ge\theta\}$ is nonempty and bounded
for all sufficiently large $t$. For these times, set
\[
 L_\theta(t)
 =\max\{|n|:n\in\Z,\ u(t,\gamma(n))\ge\theta\}.
\]
Then
\begin{equation}
\label{eq:tangential-level-speed}
 \frac{L_\theta(t)}t\to c_{\parallel}
 \quad\text{ as }t\to\infty.
\end{equation}

\item[\textup{(iii)}] 
If \eqref{eq:strict-quotient} also holds, then, for every
$0<c<c_{\parallel}$,
\begin{equation}
\label{eq:tangential-interior-profile}
 \lim_{t\to\infty}
 \sup_{\{x:|\ell_{\mathcal G}(x)|\le ct\}}
 |u(t,x)-p_*(x)|=0.
\end{equation}
The supremum in \eqref{eq:tangential-interior-profile} includes
every $x\in\Tq$ satisfying $|\ell_{\mathcal G}(x)|\le ct$,
with every transverse distance $s(x)\in\N_0$.
\end{enumerate}

\end{theorem}

The limit $L_\theta(t)/t\to c_{\parallel}$ as $t\to\infty$ in
\eqref{eq:tangential-level-speed} has the following logarithmic refinement.

\Needspace{24\baselineskip}
\begin{theorem}\label{thm:tangential-bramson}

Assume \eqref{eq:hyp}, $f(x,u)=f_{s(x)}(u)$, and
$\mu\le0<\sigma_{\mathcal G}$.
Let $u$ solve \eqref{eq:tree-cauchy} with a nonzero finitely
supported initial datum $u_0\colon\Tq\to[0,1]$, and let $p_*$
be the smallest positive stationary solution from
Theorem~\ref{thm:geodesic-threshold}. For $t\ge0$, set
\[
 m_{\parallel}(t)
 =c_{\parallel}t
 -\frac{3}{2\lambda_{\parallel}}\log(1+t).
\]
For every $0<\theta<p_*(\gamma(0))$, the following conclusions hold.
\begin{enumerate}
\item[\textup{(i)}] 
There are $C_\theta,T_\theta<\infty$ such that
\[
 \{n\in\Z:|n|\le m_{\parallel}(t)-C_\theta\}
 \subseteq
 \{n\in\Z:u(t,\gamma(n))\ge\theta\}
 \subseteq
 \{n\in\Z:|n|\le m_{\parallel}(t)+C_\theta\}
\]
for every $t\ge T_\theta$.

\item[\textup{(ii)}] 
The quantity
$L_\theta(t)=\max\{|n|:n\in\Z,\ u(t,\gamma(n))\ge\theta\}$ satisfies
\[
 L_\theta(t)
 =c_{\parallel}t
 -\frac{3}{2\lambda_{\parallel}}\log(1+t)
 +O_\theta(1)\quad\text{ as }t\to\infty.
\]
\end{enumerate}

\end{theorem}

For each of the three reactions in \eqref{eq:reaction-assignments},
set $\mathcal L=D\Delta_{\Tq}+\partial_u f(\cdot,0)$ on
$\ell^2(\Tq)$.
The threshold proofs combine quadratic-form estimates with positive
subsolutions supported on finite sets. For $f(x,u)=f_{r(x)}(u)$,
the equation at $r=0$ has diffusion coefficient $D(q+1)$; for
$f(x,u)=f_{s(x)}(u)$, the equation at $s=0$ has coefficient $D(q-1)$. Positive-function identities for graph
quadratic forms are developed in
\cite{KellerPinchoverPogorzelski}.
At $\mu=0$ in the Busemann and radial cases, and at
$\sigma_{\mathcal G}=0$ in the geodesic case,
Lemmas~\ref{lem:horo-spectrum}, \ref{lem:radial-spectrum}
and~\ref{lem:geodesic-spectrum} give
$\sup\spec(\mathcal L)=0$ and $\ker\mathcal L=\{0\}$.
Since $\mathcal L$ is bounded and self-adjoint, the spectral theorem
and dominated convergence give $\|e^{t\mathcal L}v\|_2\to0$ as
$t\to\infty$ for every $v\in\ell^2(\Tq)$.
By the KPP inequality in \eqref{eq:hyp}, comparison therefore yields
$\|u(t,\cdot)\|_\infty\le\|e^{t\mathcal L}u_0\|_2\to0$ as
$t\to\infty$ for $u_0\in\ell^2(\Tq)$ with $0\le u_0\le1$.

For the logarithmic positions, positive eigenfunctions of the linearized
equations define continuous-time Markov chains. Killing a path at its
first crossing of a moving boundary, including a crossing between jumps,
gives the linear comparison functions. Forward and reversed survival
estimates retain the distances of both endpoints from the boundary.
By \eqref{eq:taylor}, $0\le a_i u-f_i(u)\le C_2u^2$ for
$0\le i<P$ and $0\le u\le1$.
The conditional integral bounds in Lemmas~\ref{lem:bridge-occupation}
and~\ref{lem:tangential-backward-log} control the contribution of
$C_2u^2$ along paths conditioned on survival and their endpoints,
yielding lower bounds for solutions of the nonlinear equation.
For the Markov and additive-functional formulas, see
\cite{Norris,EthierKurtz,FitzsimmonsPitman}; for walks conditioned to stay
positive, see
\cite{Caravenna,DenisovSakhanenkoWachtel,GramaLauvergnatLePage2018,GramaLauvergnatLePage}.
For longitudinal propagation, assume $\mu\le0<\sigma_{\mathcal G}$.
Let $I\subset(0,\infty)$ be a sufficiently small compact interval
containing $\lambda_{\parallel}$.
For each $\lambda\in I$, Lemma~\ref{lem:tangential-regeneration}
constructs a Markov chain $(N_t,R_t)$ on $\Z\times\N_0$, whose
coordinates record longitudinal position and transverse distance.
Started with $R_0=0$, the pairs of increments of $N_t$ and waiting
times between successive jumps of $N_t$ are independent and
identically distributed. The waiting time $T$ until the first such
jump satisfies $\sup_{\lambda\in I}\mathbb E^\lambda e^{\eta T}<\infty$
for some $\eta>0$, where $\mathbb E^\lambda$ denotes
expectation for $(N_t,R_t)$ started at $(N_0,R_0)=(0,0)$.
Longitudinal translation preserves the law of $T$, so the same
bound holds when $(N_0,R_0)=(j,0)$ for any $j\in\Z$. Lemma~\ref{lem:tangential-regeneration} proves the exponential-moment
bound for $T$ by killing the transverse chain at the first longitudinal
jump and using
$\inf_{\lambda\in I}(\Lambda_{\parallel}(\lambda)-\sigma_{\mathcal G})>0$.

\paragraph{Notation.}
The following table records the spectral parameters and level positions
used in the main results.
\begin{center}
\small
\renewcommand{\arraystretch}{1.15}
\begin{tabular}{@{}p{0.28\textwidth}p{0.68\textwidth}@{}}
\hline
Notation & Definition \\
\hline
$r(x),b(x)$ & Distance from $o$ and the integer-valued Busemann coordinate. \\
$s(x),\ell_{\mathcal G}(x)$ & Distance from $\mathcal G$ and the index of the closest vertex on $\mathcal G$. \\
$M_\lambda,k(\lambda)$ & The periodic matrix in \eqref{eq:matrix} and its principal eigenvalue. \\
$\lambda_0,\mu$ & $\lambda_0=\tfrac12\log q$ and $\mu=k(\lambda_0)$. \\
$c_\pm,\lambda_\pm$ & The extremal quotients in \eqref{eq:parameters} and their unique extremizers. \\
$\Psi_a,\rho_q$; $c_*(a),\lambda_*(a)$ & The homogeneous dispersion function, extinction threshold, speed and minimizing exponent in \eqref{eq:homogeneous-parameters}. \\
$\mathcal J_{\mathcal G},\sigma_{\mathcal G}$ & The operator on $\ell^2(\N_0)$ in \eqref{eq:geodesic-jacobi} and $\sigma_{\mathcal G}=\sup\spec(\mathcal J_{\mathcal G})$. \\
$p_*$ & The pointwise smallest stationary solution $p:\Tq\to(0,1]$ of the differential equation in \eqref{eq:tree-cauchy} with $f(x,u)=f_{s(x)}(u)$ when $\sigma_{\mathcal G}>0$. \\

$H_\theta,H_\theta^b,R_\theta^o,R_\theta^{\mathcal G}$ & The largest lattice index, Busemann coordinate, and distances from $o$ and $\mathcal G$, respectively, at level $\theta$. \\
$m_*(t)$ & $c_+t-3\log(1+t)/(2\lambda_+)$. \\
$\Lambda_{\parallel},c_{\parallel},\lambda_{\parallel}$ & The upper spectral endpoint of $\mathcal J_{\parallel}(\lambda)$, the minimum in \eqref{eq:tangential-speed-intro}, and its unique minimizer. \\
$m_{\parallel}(t),L_\theta(t)$ & $c_{\parallel}t-3\log(1+t)/(2\lambda_{\parallel})$ and the largest $|n|$ with $u(t,\gamma(n))\ge\theta$. \\
\hline
\end{tabular}
\end{center}

\section{Comparison and periodic kernel estimates}\label{sec:comparison}

The comparison principle applies to each of the three coordinate
reductions and to the equation on the tree. We prove it on the common
space $\ell^\infty(S)$, with $S=\Tq$ for \eqref{eq:tree-cauchy},
$S=\Z$ for \eqref{eq:lattice}, and $S=\N_0$ for the half-line
systems. Throughout this section the reaction family satisfies
\eqref{eq:hyp}. In statements concerning the tree, we also allow
$f(x,s)=f_{\iota(x)}(s)$ for an arbitrary fixed map
$\iota:\Tq\to\{0,\ldots,P-1\}$.

A solution of one of these systems defined on $\R\times S$ is called
an entire solution.
Fix $S\in\{\Tq,\Z,\N_0\}$, let $I\subseteq\R$ be an
interval, and let $v:I\times S\to\R$. For functions $v_n:D_n\to\R$
with $D_n\subseteq\R\times S$, we say that $v_n$ converges locally
to $v$ if, for every compact interval $J\subseteq I$ and every
nonempty finite $F\subseteq S$, the inclusion $J\times F\subseteq D_n$
holds for all sufficiently large $n$ and
\[
 \sup_{t\in J}\sup_{x\in F}|v_n(t,x)-v(t,x)|\to0
 \quad\text{as }n\to\infty.
\]
In particular, for $I=\R$ one may take $J=[-T,T]$ for every $T>0$.

For each of \eqref{eq:tree-cauchy}, \eqref{eq:lattice}, and the
half-line systems, write $Q$ for its diffusion operator and $f(x,s)$
for its reaction at spatial coordinate $x\in S$ and value
$s\in[0,1]$. For $i\in\{0,\ldots,P-1\}$, fix the globally
Lipschitz extension
\[
 \widetilde f_i(s)=
 \begin{cases}
  a_i s,&s<0,\\
  f_i(s),&0\le s\le1,\\
  f_i'(1)(s-1),&s>1.
 \end{cases}
\]
Here $s\in\R$. For $x\in S$, define
$\widetilde f(x,s)=\widetilde f_i(s)$ whenever $f(x,\cdot)=f_i$. On $[0,1]$, the extensions
agree with the original reactions.
A subsolution on an interval $I$ is a locally Lipschitz map
$v:I\to\ell^\infty(S)$ satisfying
\[
 \partial_t v(t,x)\le (Qv(t,\cdot))(x)+\widetilde f(x,v(t,x))
 \quad\text{for every }x\in S\text{ and almost every }t\in I.
\]
A supersolution $z:I\to\ell^\infty(S)$ is locally Lipschitz and
satisfies
\[
 \partial_t z(t,x)\ge (Qz(t,\cdot))(x)+\widetilde f(x,z(t,x))
 \quad\text{for every }x\in S\text{ and almost every }t\in I.
\]
Since $S$ is countable, these coordinate inequalities hold outside a
common null set of times.
If $w$ solves \eqref{eq:lattice} and $m\in\Z$, then the shifted sequence
$v_j(t)=w_{j+m}(t)$ solves the same diffusion equation with reaction
$f_{j+m}(v_j)$, where the reaction indices are reduced modulo $P$.
For reactions $f^{(n)}(x,s)$ and $f^\infty(x,s)$ on
$S\times[0,1]$, local convergence means
$\sup_{s\in[0,1]}|f^{(n)}(x,s)-f^\infty(x,s)|\to0$ as
$n\to\infty$ for every fixed $x\in S$. Their bounds and Lipschitz
constants in $s\in[0,1]$ are required to be independent of $n$ and
$x$. Assignments from the fixed finite family $\{f_0,\ldots,f_{P-1}\}$
have these uniform bounds and Lipschitz constants.

For a nonempty finite set $\Omega\subset\Tq$ and a function
$v:\Omega\to\C$, extend $v$ by zero outside $\Omega$ and set
$a(x)=\partial_u f(x,0)$. The Dirichlet linearization is
\[
 (\mathcal L_\Omega v)(x)
 =D\sum_{\substack{y\sim x\\y\in\Omega}}v(y)
   +(a(x)-D(q+1))v(x),\quad\text{ for }x\in\Omega.
\]
The diagonal term $-D(q+1)v(x)$ counts all edges incident to $x$,
including those ending outside $\Omega$. For adjacent $x,y\in\Omega$,
the two matrix entries are both $D$, so $\mathcal L_\Omega$ is real
symmetric. We call its largest eigenvalue the principal eigenvalue.
{When $\Omega$ is connected, the matrix
$\mathcal L_\Omega+D(q+1)I$ is nonnegative and irreducible,
because its diagonal entries are $a(x)>0$ and its off-diagonal
entries along the edges of $\Omega$ are $D>0$.} The
Perron--Frobenius theorem then gives simplicity of the principal
eigenvalue and a strictly positive eigenvector
\cite[Chapter~8]{HornJohnson}.
For a bounded operator $\mathcal A$, $\spec(\mathcal A)$ denotes
the set of $z\in\C$ for which $zI-\mathcal A$ has no bounded inverse
on the stated space. For $z\notin\spec(\mathcal A)$ the inverse
$(zI-\mathcal A)^{-1}$ is its resolvent. We use the norm-convergent
series $e^{t\mathcal A}=\sum_{n\ge0}t^n\mathcal A^n/n!$.
For an isolated simple eigenvalue $E$, let $\Gamma$ be a positively
oriented closed contour in $\C\setminus\spec(\mathcal A)$ enclosing
$E$ and no other spectral value. The Riesz projection associated with $E$ is
\[
 \Pi_E=\frac1{2\pi\mathrm i}
     \int_\Gamma(zI-\mathcal A)^{-1}\,\dd z.
\]
See \cite[Chapters~III and VII]{Kato}. If $\mathcal A$
is self-adjoint on a Hilbert space $\mathcal H$, then
\[
 \sup\spec(\mathcal A)
 =\sup_{0\ne v\in\mathcal H}
   \frac{\langle v,\mathcal Av\rangle_{\mathcal H}}
        {\|v\|_{\mathcal H}^2};
\]
see \cite[Theorem~2.19]{Teschl}. Comparison methods for locally finite
graphs appear in \cite{HuWang2026}. We give a direct proof on
$\ell^\infty$ using the uniform bound on the jump rates.

\begin{lemma}\label{lem:comparison}

Let $S\in\{\Tq,\Z,\N_0\}$ be the spatial domain, with
$S=\Tq$ for \eqref{eq:tree-cauchy}, $S=\Z$ for \eqref{eq:lattice},
and $S=\N_0$ for \eqref{eq:radial}. Each of these equations has a
unique global solution for each initial datum with values in $[0,1]$.
The solution belongs to $C^1([0,\infty);\ell^\infty(S))$ and takes values in
$[0,1]$. If the initial datum is nonzero, the solution is strictly positive
at each vertex or index for every $t>0$. These conclusions also
hold for the half-line system with the interior equation in
\eqref{eq:radial} and boundary equation
\[
 \dot U_0=b_0(U_1-U_0)+f_0(U_0),\quad\text{ for }t\ge0,\ b_0>0.
\]
The choices $b_0=D(q+1)$ and $b_0=D(q-1)$ give the root-distance and
geodesic-distance systems, respectively. If the initial datum is not
identically one, the solution is strictly below one at every positive
time and every vertex or index.
For comparison, let $S$ be the spatial domain of one of these systems,
let $Q$ be its diffusion operator, and let $\widetilde f(x,s)$ be
the corresponding reaction extension fixed above. Fix $T>0$ and an
arbitrary nonempty set $F\subseteq S$, which may be infinite or
disconnected. Suppose that $v,z:[0,T]\to\ell^\infty(S)$ are locally
Lipschitz and satisfy, for every $x\in F$ and almost every $t\in(0,T)$,
\[
 \begin{aligned}
 \partial_t v(t,x)&\le (Qv(t,\cdot))(x)+\widetilde f(x,v(t,x)),\\
 \partial_t z(t,x)&\ge (Qz(t,\cdot))(x)+\widetilde f(x,z(t,x)).
 \end{aligned}
\]
Assume also that
\[
 \begin{aligned}
 v(0,x)&\le z(0,x) &&\text{for every }x\in F,\\
 v(t,y)&\le z(t,y) &&\text{for every }(t,y)\in[0,T]\times(S\setminus F).
 \end{aligned}
\]
Then $v(t,x)\le z(t,x)$ for every $(t,x)\in[0,T]\times F$.
Taking $F=S$ gives comparison on the full spatial domain.

\end{lemma}
\begin{proof}

\noindent\textbf{Step 1: Existence and uniqueness.}
Fix one of the systems in the statement and let $S$ be its spatial
domain. We work in the real Banach space $\ell^\infty(S)$. Its diffusion
operator can be written as
\begin{equation}\label{eq:kern-bounded-rates}
 \begin{aligned}
 (Qw)(x)&=\sum_{y\in S}d_{xy}(w(y)-w(x))
       &&\text{for }x\in S,\ w\in\ell^\infty(S),\\
 d_{xy}&\ge0 &&\text{for }x,y\in S,\\
 \sigma_x&:=\sum_{y\in S}d_{xy}\le\Lambda
       &&\text{for }x\in S.
 \end{aligned}
\end{equation}
For the tree, $d_{xy}=D$ when $y\sim x$ and is zero otherwise.
For the lattice, $d_{j,j-1}=D$ and $d_{j,j+1}=Dq$.
For the half-line, these rates hold at $j\ge1$, while $d_{0,1}=b_0$;
all other rates are zero. One may take $\Lambda=D(q+1)$ for the tree
and lattice and $\Lambda=\max\{D(q+1),b_0\}$ for the half-line.
It follows from \eqref{eq:kern-bounded-rates} that
\begin{equation}\label{eq:kern-diffusion-norm}
 \|Qw\|_{\ell^\infty(S)}\le2\Lambda\|w\|_{\ell^\infty(S)}\quad\text{for }w\in\ell^\infty(S).
\end{equation}
Write $\widetilde f(x,s)$ for the reaction extension fixed above at
coordinate $x$ and value $s$. Set
$L_0=\max_{0\le i<P}\max_{0\le s\le1}|f_i'(s)|<\infty$.
The definition of $\widetilde f_i$ then gives
\[
 |\widetilde f(x,s)-\widetilde f(x,r)|\le L_0|s-r|
 \quad\text{ for }x\in S,\ s,r\in\R.
\]
Define $\mathcal F:\ell^\infty(S)\to\ell^\infty(S)$ by
\[
 (\mathcal F(w))(x)=(Qw)(x)+\widetilde f(x,w(x))\quad\text{for }w\in\ell^\infty(S),\ x\in S.
\]
Since $\widetilde f(x,0)=0$, the map $\mathcal F$ takes bounded
functions to bounded functions. By \eqref{eq:kern-diffusion-norm} and
the uniform Lipschitz bound for $\widetilde f$, it satisfies
\begin{equation}\label{eq:kern-field-lipschitz}
 \|\mathcal F(w_1)-\mathcal F(w_2)\|_{\ell^\infty(S)}
 \le(2\Lambda+L_0)\|w_1-w_2\|_{\ell^\infty(S)}\quad\text{for }w_1,w_2\in\ell^\infty(S).
\end{equation}
Fix $w^0\in\ell^\infty(S)$ and choose
$T_{\rm loc}=[2(2\Lambda+L_0)]^{-1}$. On the Banach space
\[
 X=C([0,T_{\rm loc}];\ell^\infty(S)),\qquad
 \|w\|_X=\sup_{0\le t\le T_{\rm loc}}
                   \|w(t,\cdot)\|_{\ell^\infty(S)},
\]
define $\mathcal T:X\to X$ by
\[
 (\mathcal Tw)(t,x)=w^0(x)+\int_0^t
                  (\mathcal F(w(s,\cdot)))(x)\,\dd s\quad\text{for }w\in X,\ (t,x)\in[0,T_{\rm loc}]\times S.
\]
By \eqref{eq:kern-field-lipschitz}, the integrand is continuous in
$\ell^\infty(S)$, so the integral is a Banach-space integral. The same
estimate gives
\begin{equation}\label{eq:kern-contraction}
 \|\mathcal Tw_1-\mathcal Tw_2\|_X
 \le T_{\rm loc}(2\Lambda+L_0)\|w_1-w_2\|_X
 =\tfrac12\|w_1-w_2\|_X
 \text{ for }w_1,w_2\in X.
\end{equation}
By \eqref{eq:kern-contraction}, the Banach fixed-point theorem gives
a unique $w\in X$ with $\mathcal Tw=w$. Applying
\eqref{eq:kern-field-lipschitz} along $w$ shows that
$s\mapsto\mathcal F(w(s,\cdot))$ is continuous. We may therefore
differentiate the integral equation in $\ell^\infty(S)$.
Thus $w\in C^1([0,T_{\rm loc}];\ell^\infty(S))$ and
$\partial_t w(t,x)=(\mathcal F(w(t,\cdot)))(x)$ at every coordinate.
The value of $T_{\rm loc}$ depends only on $\Lambda$ and $L_0$.
Thus the same construction applies successively on every interval
$[nT_{\rm loc},(n+1)T_{\rm loc}]$, with the value already obtained at
its left endpoint. Uniqueness on each interval makes the constructions
agree at their common endpoints. Their derivatives also agree there,
by continuity of $\mathcal F$. Hence
$w\in C^1([0,\infty);\ell^\infty(S))$.

\medskip
\noindent\textbf{Step 2: Comparison on an arbitrary subset.}
Fix $T>0$ and a nonempty set $F\subseteq S$. Let
$v,z:[0,T]\to\ell^\infty(S)$ be the locally Lipschitz functions
satisfying the subsolution and supersolution inequalities on $F$
in the statement. Their initial and exterior
values satisfy
\[
 \begin{aligned}
 v(0,x)&\le z(0,x) &&\text{for }x\in F,\\
 v(t,y)&\le z(t,y) &&\text{for }(t,y)\in[0,T]\times(S\setminus F).
 \end{aligned}
\]
Choose $\beta\ge\Lambda+L_0$ and define
\[
 H(x,s)=(\beta-\sigma_x)s+\widetilde f(x,s)\quad(x\in S,\ s\in\R),\qquad
 M_F(t)=\sup_{x\in F}(v(t,x)-z(t,x))_+\quad\text{for }t\in[0,T],
\]
where $r_+=\max\{r,0\}$ for $r\in\R$.
For $s\ge r$, the Lipschitz bound for $\widetilde f$ gives
\[
 0\le(\beta-\sigma_x-L_0)(s-r)
 \le H(x,s)-H(x,r)
 \le(\beta-\sigma_x+L_0)(s-r)\quad\text{for }x\in S.
\]
Consequently $H(x,\cdot)$ is nondecreasing and has Lipschitz
constant at most $\beta-\sigma_x+L_0$.
Moreover,
\[
 |M_F(t)-M_F(s)|
 \le\|v(t,\cdot)-v(s,\cdot)\|_\infty
     +\|z(t,\cdot)-z(s,\cdot)\|_\infty\quad\text{for }s,t\in[0,T],
\]
so $M_F$ is continuous on $[0,T]$. The coordinate functions of
$v,z$ are absolutely continuous. Since $S$ is countable, their
differential inequalities hold simultaneously at every $x\in F$
outside one null set of times.

At each $x\in F$, linearity of $Q$ gives
\[
 \begin{aligned}
 &(Qv(t,\cdot))(x)-(Qz(t,\cdot))(x)\\
 &=\sum_{y\in F}d_{xy}(v(t,y)-z(t,y))
   +\sum_{y\in S\setminus F}d_{xy}(v(t,y)-z(t,y))\\
 &\quad-\sigma_x(v(t,x)-z(t,x))\quad\text{for }t\in[0,T].
 \end{aligned}
\]
The exterior ordering makes the second sum nonpositive for every
$t\in[0,T]$. Subtracting the subsolution and supersolution
inequalities therefore yields
\[
 \begin{aligned}
 &\partial_t(v(t,x)-z(t,x))+\beta(v(t,x)-z(t,x))\\
 &\quad\le\sum_{y\in F}d_{xy}(v(t,y)-z(t,y))
       +H(x,v(t,x))-H(x,z(t,x))\quad\text{for }x\in F\text{ and a.e. }t\in(0,T).
 \end{aligned}
\]
For every $x\in F$, the sum over $F$ is at most
$\sigma_xM_F(t)$, while monotonicity and the Lipschitz bound for $H$
give
\[
 H(x,v(t,x))-H(x,z(t,x))
 \le(\beta-\sigma_x+L_0)(v(t,x)-z(t,x))_+\quad\text{for }t\in[0,T].
\]
Using the displayed bound on
$H(x,v(t,x))-H(x,z(t,x))$ and summing the coefficients
$d_{xy}$ in \eqref{eq:kern-bounded-rates}, we obtain
\begin{equation}\label{eq:kern-coordinate-comparison}
 \partial_t(v(t,x)-z(t,x))+\beta(v(t,x)-z(t,x))
 \le(\beta+L_0)M_F(t)
\end{equation}
for each $x\in F$ and almost every $t\in(0,T)$.
Multiplying \eqref{eq:kern-coordinate-comparison} by $e^{\beta t}$
and integrating at a fixed $x\in F$, we find
\[
 \begin{aligned}
 v(t,x)-z(t,x)
 &\le e^{-\beta t}(v(0,x)-z(0,x))\\
 &\quad+(\beta+L_0)\int_0^t e^{-\beta(t-s)}M_F(s)\,\dd s\quad\text{for }t\in[0,T].
 \end{aligned}
\]
Taking positive parts and then the supremum over $x\in F$ gives
\begin{equation}\label{eq:kern-supremum-comparison}
 M_F(t)\le e^{-\beta t}M_F(0)
       +(\beta+L_0)\int_0^t e^{-\beta(t-s)}M_F(s)\,\dd s\quad\text{for }t\in[0,T].
\end{equation}
The initial ordering implies $M_F(0)=0$. Applying Gronwall's
inequality to \eqref{eq:kern-supremum-comparison}, after multiplication
by $e^{\beta t}$, yields $M_F(t)=0$ on $[0,T]$.
Thus $v(t,x)\le z(t,x)$ on $[0,T]\times F$. Since $T>0$ is
arbitrary, this conclusion applies on every finite interval where
the stated hypotheses hold. Translation of time gives comparison
with any prescribed initial time.

\medskip
\noindent\textbf{Step 3: Values in $[0,1]$.}
Let $u$ be the global solution constructed in Step~1 with initial
datum $u_0:S\to[0,1]$. Since $Q\mathbf1=0$ and
$\widetilde f(x,0)=\widetilde f(x,1)=0$, both constant functions
zero and one solve the extended equation. Applying Step~2 with
$F=S$ gives $0\le u(t,x)\le1$ for all $t\ge0$ and $x\in S$.
On this range $\widetilde f(x,u(t,x))=f(x,u(t,x))$, so $u$ solves
the original system.

\medskip
\noindent\textbf{Step 4: Strict positivity.}
By \eqref{eq:kern-bounded-rates}, the operator $Q+\Lambda I$ has
nonnegative matrix entries. The norm-convergent series
\[
 e^{tQ}=e^{-\Lambda t}
       \sum_{n=0}^\infty\frac{t^n}{n!}(Q+\Lambda I)^n\quad\text{for }t\ge0
\]
therefore defines an operator preserving nonnegativity.
For every $x,y\in S$, there are an integer $n\ge0$ and vertices
$x=x_0,x_1,\ldots,x_n=y$ such that
$d_{x_{k-1},x_k}>0$ for $1\le k\le n$. The corresponding product
of matrix entries is positive; when $n=0$, use the identity operator.
For $y\in S$, define $\delta_y:S\to\{0,1\}$ by
$\delta_y(x)=\mathbf1_{\{x=y\}}$ for $x\in S$. Consequently,
\[
 ((Q+\Lambda I)^n\delta_y)(x)>0.
\]
By the exponential series, $(e^{tQ}\delta_y)(x)>0$ for every
$t>0$. If $u_0\ge0$ is nonzero, choose $y$ with $u_0(y)>0$.
Positivity of the semigroup gives
\[
 (e^{tQ}u_0)(x)\ge u_0(y)(e^{tQ}\delta_y)(x)>0\quad\text{for }t>0,\ x\in S.
\]
By Step~3 and \eqref{eq:hyp}, $f(x,u(t,x))\ge0$ at every coordinate.
Variation of constants and positivity of $e^{tQ}$ therefore give
\begin{equation}\label{eq:kern-positive-solution}
 u(t,\cdot)=e^{tQ}u_0+
    \int_0^t e^{(t-s)Q}f(\cdot,u(s,\cdot))\,\dd s
 \ge e^{tQ}u_0\quad\text{for }t\ge0.
\end{equation}
The strict lower bound for $(e^{tQ}u_0)(x)$ proved above now gives
$u(t,x)>0$ for $t>0$ and $x\in S$.

\medskip
\noindent\textbf{Step 5: Values strictly below one.}
Define $h(t,x)=1-u(t,x)$ for $(t,x)\in[0,\infty)\times S$.
By Step~3, $h\ge0$, and the equation
for $u$ gives
\[
 \partial_t h(t,x)=(Qh(t,\cdot))(x)-f(x,1-h(t,x))\quad\text{for }t\ge0,\ x\in S.
\]
Since $f(x,1)=0$ and $f(x,\cdot)$ has Lipschitz constant at most
$L_0$ on $[0,1]$,
\[
 -f(x,1-h(t,x))\ge-L_0h(t,x)\quad\text{for }t\ge0,\ x\in S.
\]
Variation of constants for the operator $Q-L_0I$, whose semigroup
is $e^{-L_0t}e^{tQ}$, consequently gives
\[
 h(t,\cdot)\ge e^{-L_0t}e^{tQ}(1-u_0)\quad\text{for }t\ge0.
\]
If $u_0\not\equiv1$, Step~4 applied to the nonzero function
$1-u_0\ge0$ makes this lower bound strictly positive at every
coordinate for $t>0$. Therefore $u(t,x)<1$ for every $t>0$ and
$x\in S$.

\medskip
\noindent\textbf{Step 6: The half-line boundary equation.}
For $S=\N_0$ and $w\in\ell^\infty(\N_0)$, the rates in Step~1 give
\[
 (Qw)(0)=b_0(w(1)-w(0)),\qquad
 (Qw)(j)=D(w(j-1)+qw(j+1)-(q+1)w(j))
 \text{ for }j\ge1.
\]
Thus the equation at zero is exactly
$\dot U_0=b_0(U_1-U_0)+f_0(U_0)$.
For each fixed $b_0>0$, the rate bound
$\Lambda=\max\{D(q+1),b_0\}$ is finite, and both directions of every
half-line edge have positive rates. The norm bound \eqref{eq:kern-diffusion-norm}, the contraction
\eqref{eq:kern-contraction}, and the comparison estimate
\eqref{eq:kern-supremum-comparison} consequently apply with this
choice of $\Lambda$. Positive rates in both directions along every
edge give the path condition used in Step~4. Hence Steps~1--5 prove
existence, uniqueness, comparison, and the strict bounds for the stated
half-line boundary equation.\qedhere

\end{proof}

\begin{remark}\label{rem:finite-mass-levels}

Let $u$ solve \eqref{eq:tree-cauchy} with
$u_0\in\ell^1(\Tq)$ and $0\le u_0\le1$. The KPP bound and the
chosen extension give $\widetilde f_i(s)\le a_{\max}s$ for $s\ge0$.
Thus comparison with the linear solution gives
\[
 0\le u(t,\cdot)\le e^{a_{\max}t}P_tu_0,
 \qquad P_t=e^{tD\Delta_{\Tq}}\quad\text{for }t\ge0.
\]
Write $P_t(x,y)=(P_t\delta_y)(x)$ for the heat kernel.
Define the adjacency operator $A$ on $\ell^\infty(\Tq)$ by
$(Av)(x)=\sum_{y\sim x}v(y)$. The expansion
\[
 P_t=e^{-D(q+1)t}\sum_{n\ge0}\frac{(Dt)^n}{n!}A^n
 \quad\text{ for }t\ge0
\]
gives nonnegative symmetric entries. Moreover,
$\Delta_{\Tq}\mathbf1=0$ gives $\sum_{y\in\Tq}P_t(x,y)=1$, so symmetry
also gives $\sum_{x\in\Tq}P_t(x,y)=1$. Tonelli's theorem therefore yields,
for every finite $t\ge0$,
\begin{equation}
 \label{eq:finite-mass-level-count}
 \|u(t,\cdot)\|_1\le e^{a_{\max}t}\|u_0\|_1,
 \qquad
 \#\{x:u(t,x)\ge\theta\}
 \le\theta^{-1}e^{a_{\max}t}\|u_0\|_1
 \quad\text{ for }0<\theta<1.
\end{equation}
By \eqref{eq:finite-mass-level-count}, every positive level set is
finite at each finite time. Horoballs and nonempty geodesic
neighborhoods contain infinitely many vertices. Accordingly, the
initial-data assumptions in Theorems~\ref{thm:periodic-bramson-step}
and~\ref{thm:periodic-bramson-distance} specify positive lower bounds
on these infinite sets in their horoball and geodesic-neighborhood assertions.

\end{remark}

\begin{lemma}\label{lem:liouville}
{Let $S$ be the spatial domain of one of the systems in
Lemma~\ref{lem:comparison}, including any spatial phase of
\eqref{eq:lattice}, and let $W:\R\times S\to[0,1]$ be an
entire solution of that system. If, for some constant $\delta>0$,}
\[
 \delta\le W(t,x)\le1
 {\quad\text{for }(t,x)\in\R\times S,}
\]
then $W\equiv1$.
\end{lemma}

\begin{proof}
Set $g(s)=\min_{0\le i<P}f_i(s)$ for $s\in[0,1]$. By
\eqref{eq:hyp} and finiteness of the phase set, $g$ is Lipschitz,
$g(0)=g(1)=0$, and $g(s)>0$ for $0<s<1$. If
$0<\delta<1$, let $y:[0,\infty)\to[\delta,1]$ solve
$y'=g(y)$ with $y(0)=\delta$. Positivity of $g$ on $(0,1)$ and
$g(1)=0$ imply that $y(t)$ increases to one as $t\to\infty$.
{Fix $t_0\in\R$ and $\tau>0$. The spatially constant function
$y(t-t_0+\tau)$, defined for $t\ge t_0-\tau$, is a subsolution of the
system solved by $W$: its diffusion term is zero and $g\le f_i$
for every $i\in\{0,\ldots,P-1\}$. At time $t_0-\tau$, its value is
$y(0)=\delta\le W(t_0-\tau,x)$. Lemma~\ref{lem:comparison} therefore
gives $W(t_0,x)\ge y(\tau)$. Letting $\tau\to\infty$ and using $W\le1$}
yields $W(t_0,x)=1$. Since $t_0$ and $x$ were arbitrary, $W\equiv1$.
If $\delta=1$, the conclusion follows directly from the assumed bounds.\qedhere
\end{proof}

\begin{lemma}\label{lem:compactness}

Fix $S\in\{\Tq,\Z,\N_0\}$ and a diffusion operator $Q$ from
Lemma~\ref{lem:comparison} on $\ell^\infty(S)$. Let $I_n$ be time
intervals that eventually contain every bounded interval of $\R$,
and let $f^{(n)}(x,\cdot)$ be the reaction at coordinate $x$ in the
$n$-th equation. Suppose $u_n:I_n\times S\to[0,1]$ solves
\[
 \partial_tu_n(t,x)=(Qu_n(t,\cdot))(x)
                    +f^{(n)}(x,u_n(t,x))
 \quad\text{for }(t,x)\in I_n\times S.
\] Assume that, for every $x\in S$,
$f^{(n)}(x,\cdot)$ converges uniformly on $[0,1]$ to $f^\infty(x,\cdot)$, and that all
$f^{(n)}(x,\cdot)$ have common bounds and Lipschitz constants on $[0,1]$.
Then, along a subsequence as $n\to\infty$, the solutions converge
locally to an entire solution $u$ of
$\partial_tu(t,x)=(Qu(t,\cdot))(x)+f^\infty(x,u(t,x))$, with
$u\in C^1(\R;\ell^\infty(S))$.

The same conclusion, with limiting domain $\Z$, holds for solutions
$u_n:I_n\times J_n\to[0,1]$, where $J_n\subseteq\Z$ eventually contains every finite
subset of $\Z$, provided the lattice equation holds wherever both
nearest neighbors belong to $J_n$ and the reactions converge locally
with the same uniform bounds and Lipschitz constants. In particular,
let $w$ be a $[0,1]$-valued solution of \eqref{eq:lattice} for
$t\ge0$. For any integer sequence $(m_n)$ and any sequence
$t_n\to\infty$, the functions $v_n(t,j)=w_{j+m_n}(t+t_n)$ have a
local entire subsequential limit after restricting to a subsequence
on which $m_n\bmod P$ is constant. If this residue is $m$, the
limiting equation has reaction $f_{j+m}$.

For any fixed system in Lemma~\ref{lem:comparison}, a solution in
$[0,1]$ that is monotone in time has a pointwise limit as
$t\to\infty$, and that limit is stationary.
\end{lemma}
\begin{proof}
Let $F_0$ and $K$ be a common bound and a common Lipschitz constant
for the reactions on $[0,1]$, respectively, and let $\Lambda$ bound
the total jump rates of the fixed diffusion operator. The limiting
reactions inherit the bounds $F_0$ and $K$. On a fixed domain,
\[
 |\partial_tu_n(t,x)|\le\Lambda+F_0
\]
uniformly in $n$, $t$ and $x$ whenever the equation is defined.
The derivative bound gives equicontinuity at each coordinate,
uniformly in $n$ on any fixed bounded time interval contained in
$I_n$ for all sufficiently large $n$. Enumerate $S$ and apply the
Arzel\`a--Ascoli theorem successively on $[-m,m]$ at its first $m$
coordinates. Taking the diagonal subsequence gives the claimed local
convergence to a function $u:\R\times S\to[0,1]$.

At a fixed coordinate $x$, only finitely many neighboring values
enter $Q$. The reaction term converges uniformly on each bounded
time interval, because
$|f^{(n)}(x,u_n(t,x))-f^\infty(x,u(t,x))|
\le K|u_n(t,x)-u(t,x)|+
\|f^{(n)}(x,\cdot)-f^\infty(x,\cdot)\|_{C([0,1])}$.
The uniform convergence just established permits passage to the limit
in the coordinatewise integral equation at each $x\in S$. Passing to
the limit also in the time increment estimate gives
$|u(t,x)-u(s,x)|\le(\Lambda+F_0)|t-s|$ for every $x\in S$.
Taking the supremum over $x$ shows that
$t\mapsto u(t,\cdot)$ is locally Lipschitz in $\ell^\infty(S)$.
For $v\in\ell^\infty(S)$ with $0\le v\le1$, define
$\mathcal F_\infty(v)(x)=(Qv)(x)+f^\infty(x,v(x))$.
By \eqref{eq:kern-diffusion-norm} and the common Lipschitz bound $K$,
$\mathcal F_\infty$ has Lipschitz constant at most $2\Lambda+K$
on this set. Hence $t\mapsto\mathcal F_\infty(u(t,\cdot))$ is
continuous in $\ell^\infty(S)$, and the coordinatewise integral
equation is the Banach-space integral equation for $\mathcal F_\infty$. It follows that
$u\in C^1(\R;\ell^\infty(S))$.

For the domains $J_n$, fix a finite set $F\subset\Z$ and put
$F^{\pm}=F\cup\{j-1,j+1:j\in F\}$. For all sufficiently large $n$,
$F^{\pm}\subset J_n$, so the lattice equation holds at every
$j\in F$. Use $f_j^{(n)}$ and $f_j^\infty$ for the reactions and their
limits, as above. At each $j\in F$, the equation gives
$|\partial_tu_n(t,j)|\le D(q+1)+F_0$. Together with
$0\le u_n\le1$, this gives equicontinuity at each fixed index for all
sufficiently large $n$. Arzel\`a--Ascoli on $[-m,m]$ for $m\in\N$ and a diagonal extraction
over $(m,j)\in\N\times\Z$ therefore give, along a
subsequence as $n\to\infty$, local convergence to a function
$u:\R\times\Z\to[0,1]$. The coordinatewise integral equation reads
\[
 \begin{aligned}
 u_n(t,j)-u_n(s,j)
 =\int_s^t\Bigl[&D\bigl(u_n(r,j-1)+q u_n(r,j+1)
                    -(q+1)u_n(r,j)\bigr)\\
              &+f_j^{(n)}(u_n(r,j))\Bigr]\,\dd r.
 \end{aligned}
\]
For fixed $s,t\in\R$, the interval between them lies in $I_n$ for
all sufficiently large $n$. Local convergence on $F^{\pm}$ and
uniform convergence of $f_j^{(n)}$ on $[0,1]$ therefore give the
lattice integral equation at each $j\in\Z$. The derivative bound yields
$\|u(t,\cdot)-u(s,\cdot)\|_\infty\le(D(q+1)+F_0)|t-s|$.
On $\{v\in\ell^\infty(\Z):0\le v\le1\}$, the map
$v\mapsto (D(v_{j-1}+qv_{j+1}-(q+1)v_j)+f_j^\infty(v_j))_{j\in\Z}$
has Lipschitz constant at most $2D(q+1)+K$. Thus the integral equation gives
$u\in C^1(\R;\ell^\infty(\Z))$. For
$v_n(t,j)=w_{j+m_n}(t+t_n)$, the time domains $[-t_n,\infty)$
eventually contain every bounded interval of $\R$. The reaction at
$j$ is $f_{j+m_n}$, which equals $f_{j+m}$ along a subsequence with
$m_n\equiv m\pmod P$. The fixed-domain conclusion therefore applies.
More generally, for assignments from the fixed finite reaction
family, a diagonal subsequence makes the reaction eventually
constant at each fixed coordinate of the countable spatial domain.

Finally, fix a time-monotone solution $u$ of a fixed system and
write $Q$ for its diffusion operator and $f(x,s)$ for its reaction at
coordinate $x$ and value $s$. Since $0\le u\le1$, the limit
$p(x)=\lim_{t\to\infty}u(t,x)$ exists at every coordinate.
For $h>0$, write the integral equation between $t$ and $t+h$ and
substitute $s=t+r$, with $0\le r\le h$. At $x$ and each of its
finitely many neighbors, $u(t+r,\cdot)$ converges pointwise to $p$
as $t\to\infty$. The integrands are uniformly bounded, so bounded
convergence gives
\[
 0=\lim_{t\to\infty}\bigl(u(t+h,x)-u(t,x)\bigr)
   =h\bigl((Qp)(x)+f(x,p(x))\bigr).
\]
Thus $p$ is stationary.
\end{proof}

\subsection{The periodic matrix and the associated chain}\label{sec:periodic-chain}

The matrix $M_\lambda$ in \eqref{eq:matrix} determines the growth
rate $k(\lambda)$ and the weights $e^{-\lambda j}\phi_j(\lambda)$
for $j\in\Z$, with $\phi$ extended periodically. Recall that
$a_i=f_i'(0)$ and $\bar a=P^{-1}\sum_{i=0}^{P-1}a_i$. We first prove
the required properties of its principal eigenvalue $k(\lambda)$.
We then use the positive eigenvector $\phi(\lambda)$ to define a
continuous-time Markov chain. Its transition kernel gives the
linearized evolution, while the corresponding absorbed kernels give
comparison functions near the moving level sets. A path constraint must
hold at every time in the prescribed interval, including between
jumps. Constants may depend on the fixed equation and on the
parameters indicated in each statement, but are independent of the
initial and terminal sites. Lemma~\ref{lem:uniform-tilt} establishes
uniformity in $\lambda$ on compact intervals.
For a $P\times P$ complex matrix $A$, the matrix norm in the Fourier
estimates is $\|A\|_\infty=\max_{0\le i<P}\sum_{h=0}^{P-1}|A_{ih}|$.
An estimate $r(\theta)=O(|\theta|^m)$ as $\theta\to0$ means that
$|r(\theta)|\le C|\theta|^m$ in a specified neighborhood of zero;
for a matrix remainder, absolute value is replaced by this norm.

\begin{lemma}\label{lem:matrix}

The maps $\lambda\mapsto k(\lambda)$ and
$\lambda\mapsto\phi(\lambda)$ are real analytic on $\R$. Moreover,
\begin{equation}
\label{eq:convex-symmetry}
 k''(\lambda)>0,\qquad k(\lambda)=k(\log q-\lambda)\quad(\lambda\in\R),
 \qquad \min_{\lambda\in\R}k(\lambda)=\mu.
\end{equation}
Let $d_q(\lambda)=D(e^\lambda+qe^{-\lambda}-(q+1))$. Then
\begin{equation}
\label{eq:bounds}
 a_{\min}+d_q(\lambda)\le k(\lambda)\le a_{\max}+d_q(\lambda),
 \qquad k(\lambda)\ge\bar a+d_q(\lambda)\quad\text{for }\lambda\in\R.
\end{equation}
The bound $k(\lambda)\ge\bar a+d_q(\lambda)$ is strict for every
real $\lambda$ if the $a_i$ are not all equal. The minimum over $\lambda>0$ and the maximum over $\lambda<0$ in
\eqref{eq:parameters} are attained uniquely, and
\begin{equation}
\label{eq:critical-pairs}
 k(\lambda_\pm)=\lambda_\pm k'(\lambda_\pm),\qquad
 c_\pm=k'(\lambda_\pm),\qquad c_-<0,\qquad c_-<c_+.
\end{equation}
The sign of $c_+$ is the sign of $\mu$. For every $c\in(c_-,c_+)$ there
is a unique $\lambda_c\in(\lambda_-,\lambda_+)$ such that
\begin{equation}
\label{eq:interior-speed}
 k'(\lambda_c)=c,\qquad k(\lambda_c)-c\lambda_c>0.
\end{equation}

\end{lemma}
\begin{proof}

\emph{Analyticity and elementary bounds.}
Fix $\lambda\in\R$. For $P\ge2$, the off-diagonal entries of
$M_\lambda$ are nonnegative and their directed graph is strongly
connected. Since $a_i>0$, the matrix $M_\lambda+D(q+1)I$ has
positive diagonal and is nonnegative and irreducible. The
Perron--Frobenius theorem gives a simple real eigenvalue $k(\lambda)$,
strictly larger in real part than the other eigenvalues, with positive
right and left eigenvectors; see \cite[Chapter~8]{HornJohnson}. For $P=1$, \eqref{eq:matrix} gives
$k(\lambda)=a_0+D(e^\lambda+qe^{-\lambda}-(q+1))$
and $\phi_0(\lambda)=1$.
The identity
$e^{tM_\lambda}=e^{-D(q+1)t}e^{t(M_\lambda+D(q+1)I)}$
and the exponential series for the nonnegative irreducible matrix
$M_\lambda+D(q+1)I$ show that $e^{tM_\lambda}$ has strictly positive
entries for every $t>0$.
Since the entries of $M_\lambda$ are analytic in $\lambda$ and its
principal eigenvalue is simple, analytic perturbation theory gives
local analytic eigenvalue and eigenvector branches
\cite[Chapter~II]{Kato}. The normalization
$\sum_{i=0}^{P-1}\phi_i(\lambda)=1$ makes the positive right eigenvector unique,
so the local analytic branches agree on overlaps and define a real
analytic map on all of $\R$.

Apply the Collatz--Wielandt inequalities to
$M_\lambda+D(q+1)I$, and then subtract $D(q+1)$.
For $\mathbf1=(1,\ldots,1)^T\in\R^P$, this gives
\[
 \min_{0\le i<P} (M_\lambda\mathbf1)_i
 \le k(\lambda)\le
 \max_{0\le i<P} (M_\lambda\mathbf1)_i.
\]
Since $(M_\lambda\mathbf1)_i=a_i+d_q(\lambda)$, this proves the first
two bounds in \eqref{eq:bounds}.

\emph{Strict convexity and symmetry.}
Fix $\lambda\in\R$ and put $\phi_i=\phi_i(\lambda)$. In this
proof, indices are reduced modulo $P$, and unqualified sums and
products range over $i=0,\ldots,P-1$. Define the rates $\alpha_i,\beta_i$
and the operator $Q$ on $\R^P$ by the formulas below, where
$0\le i<P$ and $v\in\R^P$:
\[
 \alpha_i=De^\lambda\frac{\phi_{i-1}}{\phi_i},\qquad
 \beta_i=Dqe^{-\lambda}\frac{\phi_{i+1}}{\phi_i},\qquad
 (Qv)_i=\alpha_i(v_{i-1}-v_i)+\beta_i(v_{i+1}-v_i).
\]
Choose a positive left eigenvector $\psi$ of $M_\lambda$ with
$\sum_{i=0}^{P-1}\psi_i\phi_i=1$, and set $\pi_i=\psi_i\phi_i$. The left and
right eigenvalue equations give $\pi^TQ=0$; also $\pi_i>0$ and
$\sum_{i=0}^{P-1}\pi_i=1$. Analyticity allows us to define
$h_i=\partial_\lambda\log\phi_i$ and
$\ell_i=\partial_\lambda^2\log\phi_i$. Dividing the eigenvalue equation
$(M_\lambda\phi)_i=k(\lambda)\phi_i$ by $\phi_i$ and differentiating
twice in $\lambda$ gives
\[
 k''(\lambda)=(Q\ell)_i+
 \alpha_i(1+h_{i-1}-h_i)^2+
 \beta_i(-1+h_{i+1}-h_i)^2.
\]
Multiplying by $\pi_i$ and summing over one period removes the
$Q\ell$ term, since $\pi^TQ=0$. Thus
\begin{equation}\label{eq:kern-eigenvalue-curvature}
 k''(\lambda)=\sum_{i=0}^{P-1}\pi_i\Bigl[
 \alpha_i(1+h_{i-1}-h_i)^2+
 \beta_i(-1+h_{i+1}-h_i)^2\Bigr]\ge0.
\end{equation}
Equality in \eqref{eq:kern-eigenvalue-curvature} would require
$1+h_{i-1}-h_i=0$ for $0\le i<P$, because $\pi_i\alpha_i>0$.
The sum of these equalities is $P=0$, whereas $P\ge1$. Hence
\eqref{eq:kern-eigenvalue-curvature} gives $k''(\lambda)>0$.

From \eqref{eq:matrix} we have
$M_\lambda^T=M_{\log q-\lambda}$. Equality of their principal
eigenvalues gives the symmetry in \eqref{eq:convex-symmetry}.
Differentiating that symmetry at $\lambda_0=\frac12\log q$ yields
$k'(\lambda_0)=0$. By the strict convexity just proved, $\lambda_0$
is the unique minimizer, and therefore $\min_\R k=k(\lambda_0)=\mu$.

\emph{The arithmetic-mean bound.}
Divide the eigenvalue equation by
$\phi_i$ and average over one period:
\[
 k(\lambda)=\bar a-D(q+1)
    +De^\lambda\frac1P\sum_{i=0}^{P-1}\frac{\phi_{i-1}}{\phi_i}
    +Dqe^{-\lambda}\frac1P\sum_{i=0}^{P-1}\frac{\phi_{i+1}}{\phi_i}.
\]
The products $\prod_{i=0}^{P-1}\phi_{i-1}/\phi_i$ and
$\prod_{i=0}^{P-1}\phi_{i+1}/\phi_i$ both equal one. The
arithmetic--geometric mean inequality therefore makes each average at
least one and yields
$k(\lambda)\ge\bar a+d_q(\lambda)$. Equality forces all the ratios to
be one, hence $\phi$ to be constant. The eigenvalue equation then
forces all $a_i$ to be equal. Thus the inequality is strict for every
$\lambda$ whenever the coefficients $a_i$ are not all equal.

\emph{Critical parameters.}
By \eqref{eq:bounds}, $k(0)\ge a_{\min}>0$ and
$k(\lambda)$ grows at least exponentially as $|\lambda|\to\infty$.
Consequently
\[
 \frac{k(\lambda)}{\lambda}\to+\infty
 \quad\text{ as }\lambda\downarrow0\text{ or }\lambda\to+\infty,
\]
and
\[
 \frac{k(\lambda)}{\lambda}\to-\infty
 \quad\text{ as }\lambda\uparrow0\text{ or }\lambda\to-\infty.
\]
These limits imply that the minimum on $(0,\infty)$ and the maximum
on $(-\infty,0)$ are attained at interior points. Put $G(\lambda)=\lambda k'(\lambda)-k(\lambda)$. By
\eqref{eq:bounds}, $G(0)=-k(0)<0$, while
\eqref{eq:kern-eigenvalue-curvature} gives
$G'(\lambda)=\lambda k''(\lambda)$. Consequently $G$ decreases
strictly on $(-\infty,0)$ and increases strictly on $(0,\infty)$. Since
\[
 \left(\frac{k(\lambda)}{\lambda}\right)'
   =\frac{G(\lambda)}{\lambda^2}\quad\text{for }\lambda\in\R\setminus\{0\},
\]
the two interior extrema are the unique zeros $\lambda_-<0<\lambda_+$
of $G$. This proves their uniqueness and
$k(\lambda_\pm)=\lambda_\pm k'(\lambda_\pm)$, hence
$c_\pm=k'(\lambda_\pm)$.

Since $q\ge2$, we have $\lambda_0>0$. Strict convexity and
$k'(\lambda_0)=0$ imply that $k$ is strictly decreasing on
$(-\infty,\lambda_0)$. Hence $k(\lambda)>k(0)>0$ for $\lambda<0$,
and therefore $c_-<0$. Since $k'$ is strictly increasing and
$\lambda_-<\lambda_+$, one also has $c_-<c_+$. If $\mu<0$, evaluating $k(\lambda)/\lambda$ at $\lambda_0>0$
gives $c_+<0$. If $\mu=0$, then $k(\lambda)\ge0$ for all
$\lambda$ and $k(\lambda_0)=0$ with $\lambda_0>0$, so $c_+=0$.
If $\mu>0$, the attained minimum
$k(\lambda_+)/\lambda_+$ is positive. Hence $\operatorname{sgn}c_+=
\operatorname{sgn}\mu$.

Finally, $k'$ maps $(\lambda_-,\lambda_+)$ bijectively onto
$(c_-,c_+)$. Thus, for each $c$ in this interval, there is a unique
$\lambda_c$ with $k'(\lambda_c)=c$. Moreover,
$G(\lambda)<0$ for every
$\lambda\in(\lambda_-,\lambda_+)$, so
$k(\lambda_c)-c\lambda_c=-G(\lambda_c)>0$. This proves
\eqref{eq:interior-speed}.
\end{proof}

\subsection{The continuous-time chain}

Fix $\lambda\in\R$. Extend the rates $\alpha_i,\beta_i$ and the
weights $\pi_i$ defined in the proof of Lemma~\ref{lem:matrix}
periodically to $\Z$. Write $[j]$ for the residue of $j$ in
$\{0,\ldots,P-1\}$. We use the same symbol $Q$ for the bounded
operator on $\ell^\infty(\Z)$ given by
\[
 (Qv)_j=\alpha_j(v_{j-1}-v_j)+\beta_j(v_{j+1}-v_j)
 \quad\text{for }j\in\Z,\ v\in\ell^\infty(\Z).
\]
{After identifying a $P$-periodic vector with its entries at
$0,\ldots,P-1$, the restriction of $Q$ to these vectors is the
phase matrix from Lemma~\ref{lem:matrix}. On this $P$-dimensional
space,}
\[
 Q=\Phi^{-1}M_\lambda\Phi-k(\lambda)I,
 \qquad \Phi=\diag(\phi_0(\lambda),\ldots,\phi_{P-1}(\lambda)).
\]
In the whole-line generator, $j\to j-1$ and $j\to j+1$ remain
distinct jumps. This convention also applies when their residues
coincide for $P=1$ or $P=2$. Throughout the kernel arguments, the
first index is the starting site and the second is the terminal site;
$Q$ therefore acts on the first index.

Set $c=k'(\lambda)$ and
$h_i=\partial_\lambda\log\phi_i(\lambda)$ for
$i\in\{0,\ldots,P-1\}$.
Differentiating the eigenvalue equation in \eqref{eq:matrix} and
then multiplying by the invariant vector $\pi$ gives
\begin{equation}\label{eq:kern-drift-identity}
 k'(\lambda)=\alpha_i-\beta_i+(Qh)_i\quad\text{ for }0\le i<P,
 \qquad \sum_{i=0}^{P-1}\pi_i(\beta_i-\alpha_i)=-c.
\end{equation}
Thus, when the phase has law $\pi$, the chain with generator $Q$
has mean infinitesimal drift $-c$: at phase $i$, the left and right
jumps contribute $-\alpha_i$ and $\beta_i$, respectively.

To construct the chain, choose
$\Lambda\ge\max_{0\le i<P}(\alpha_i+\beta_i)$ with $\Lambda>0$, and set
$R=I+Q/\Lambda$ on $\Z$. The entries of $R$ are nonnegative and each
row sums to one. Let $(Y_n)_{n\ge0}$ be a discrete-time chain on $\Z$
with transition matrix $R$, and let $(N_t)_{t\ge0}$ be an independent
Poisson process of rate $\Lambda$. Define $X_t=Y_{N_t}$. Its paths
are right-continuous with left limits and have finitely many jumps
on each bounded time interval, because $N_t<\infty$ almost surely.
Write $\mathbb P_j$ for the law with $X_0=Y_0=j$ and $\mathbb E_j$
for its expectation. The transition and survival estimates use
$X_t-X_0+ct$, where $-ct$ is the mean displacement when the
initial phase has law $\pi$, by \eqref{eq:kern-drift-identity}.
Conditioning on $N_t$ gives
\begin{equation}
 \label{eq:poisson-kernel}
 p_t(j,l):=\mathbb P_j(X_t=l)
  =e^{-\Lambda t}\sum_{n=0}^\infty
       \frac{(\Lambda t)^n}{n!}R^n(j,l),\qquad
 (e^{tQ}v)_j=\sum_{l\in\Z} p_t(j,l)v_l.
\end{equation}
Here $t\ge0$ and $j,l\in\Z$.
The series converges in the operator norm on $\ell^\infty(\Z)$,
and the semigroup formula holds for $v\in\ell^\infty(\Z)$. It gives
$p_t\ge0$, $\sum_{l\in\Z}p_t(j,l)=1$, and
$\sum_{m\in\Z}p_s(j,m)p_t(m,l)=p_{s+t}(j,l)$ for $s,t\ge0$ and $j,l\in\Z$. A ring of the Poisson clock
may leave the site unchanged. Since $\Lambda(R-I)=Q$, the constructed
chain has generator $Q$; see
\cite[Sections~2.4 and 2.6--2.8]{Norris} for Poisson processes,
jump chains, nonexplosion and transition equations.

By construction, $[X_t]$ has the phase generator in
Lemma~\ref{lem:matrix}. Its stationary probability vector is therefore
$\pi$. {Periodic extension of the phase balance equations gives
$\sum_{i\in\Z}\pi_iQ(i,j)=0$ at each $j\in\Z$, where only
$i\in\{j-1,j,j+1\}$ contribute.} Hence $(\pi_j)_{j\in\Z}$ is an invariant measure for
the whole-line chain. Its mass on each period is one, so its total
mass on $\Z$ is infinite. Define the reversed rates and generator by
\begin{equation}
 \label{eq:dual-rates}
\begin{gathered}
 \widehat\alpha_i=\frac{\pi_{i-1}\beta_{i-1}}{\pi_i},\qquad
 \widehat\beta_i=\frac{\pi_{i+1}\alpha_{i+1}}{\pi_i},\\
 (\widehat Qv)_i=\widehat\alpha_i(v_{i-1}-v_i)
                  +\widehat\beta_i(v_{i+1}-v_i)\quad\text{for }i\in\Z,\ v\in\ell^\infty(\Z).
\end{gathered}
\end{equation}
The periodic rates in \eqref{eq:dual-rates} define a chain $\widehat X$
on $\Z$ with generator $\widehat Q$. Write
$\widehat{\mathbb P}_j$ and $\widehat{\mathbb E}_j$ for its law and
expectation when $\widehat X_0=j$, and set
$\widehat p_t(j,l)=\widehat{\mathbb P}_j(\widehat X_t=l)$.
Set $\lambda^\dagger=\log q-\lambda$. The identity
$M_\lambda^T=M_{\lambda^\dagger}$ makes
$\phi(\lambda^\dagger)$ a positive left eigenvector of $M_\lambda$.
The definition $\pi_i=\psi_i\phi_i(\lambda)$ and its normalization
therefore give
\[
 \pi_i=
 \frac{\phi_i(\lambda)\phi_i(\lambda^\dagger)}
 {\sum_{r=0}^{P-1}\phi_r(\lambda)\phi_r(\lambda^\dagger)}\quad\text{for }i\in\Z.
\]
Substitution in \eqref{eq:dual-rates} gives
\[
 \widehat\alpha_i=
 De^{\lambda^\dagger}
 \frac{\phi_{i-1}(\lambda^\dagger)}{\phi_i(\lambda^\dagger)},
 \qquad
 \widehat\beta_i=
 Dqe^{-\lambda^\dagger}
 \frac{\phi_{i+1}(\lambda^\dagger)}{\phi_i(\lambda^\dagger)}\quad\text{for }i\in\Z.
\]
Put $\Phi^\dagger=\diag(\phi_i(\lambda^\dagger))_{i=0}^{P-1}$.
The formulas for $\widehat\alpha_i$ and $\widehat\beta_i$ above
give, on $P$-periodic vectors,
\[
 \widehat Q=(\Phi^\dagger)^{-1}M_{\lambda^\dagger}\Phi^\dagger
                      -k(\lambda^\dagger)I.
\]
Thus $\widehat Q$ is obtained by the same construction at parameter
$\lambda^\dagger$. Differentiating
$k(\lambda)=k(\log q-\lambda)$ gives
$k'(\lambda^\dagger)=-k'(\lambda)$. The mean drift of $\widehat X$
under its stationary phase law is therefore
$-k'(\lambda^\dagger)=k'(\lambda)$.
If $\lambda$ ranges over a compact interval $I=[r_1,r_2]\subset\R$,
then $\lambda^\dagger$ ranges
over $[\log q-r_2,\log q-r_1]$. Any estimate for the chains constructed from $Q$ that is uniform on
$I\cup[\log q-r_2,\log q-r_1]$ therefore also applies to the
reversed chain for $\lambda\in I$.

The periodically extended phase balance equations read
$\pi_{i-1}\beta_{i-1}+\pi_{i+1}\alpha_{i+1}
 =\pi_i(\alpha_i+\beta_i)$ for $i\in\Z$.
Together with \eqref{eq:dual-rates}, they give
$\widehat\alpha_i+\widehat\beta_i=\alpha_i+\beta_i$.
Use the same Poisson rate $\Lambda$ for $Q$ and $\widehat Q$, and
put $\widehat R=I+\widehat Q/\Lambda$. The rate identities give
$\pi_jR(j,l)=\pi_l\widehat R(l,j)$. Multiplication and summation over
intermediate sites give the same identity for $R^n$ and
$\widehat R^n$. The power series in \eqref{eq:poisson-kernel} then gives
\begin{equation}
 \label{eq:dual-kernel}
 \pi_jp_t(j,l)=\pi_l\widehat p_t(l,j),\qquad
 \sum_{j\in\Z}\pi_jp_t(j,l)=\pi_l
 \quad\text{for }t\ge0,\ j,l\in\Z.
\end{equation}
Summing the first identity in \eqref{eq:dual-kernel} over $j\in\Z$
and using $\sum_{j\in\Z}\widehat p_t(l,j)=1$ proves the second. For a fixed
horizon $t>0$, reverse a path by reversing its successive states and
holding intervals, taking the right-continuous version at the reversed
jump times. The probability of a jump exactly at $t$ is zero. If $E\in\sigma(X_s:0\le s\le t)$ is a path event and
$\widehat E$ is the event obtained by reversing the states and holding
intervals of the paths in $E$, then
\begin{equation}\label{eq:kern-path-reversal}
 \pi_j\mathbb P_j(E,X_t=l)
   =\pi_l\widehat{\mathbb P}_l(\widehat E,\widehat X_t=j)
 \quad\text{for }j,l\in\Z.
\end{equation}
To prove \eqref{eq:kern-path-reversal}, multiply the jump rates
along a path with finitely many jumps. By \eqref{eq:dual-rates}, the
factors involving $\pi$ cancel at every intermediate state. The
holding-time factors agree because the forward and reversed total
rates agree. Summing over the number of jumps and the successive
sites, and integrating over the ordered jump times in $(0,t)$, gives
\eqref{eq:kern-path-reversal} by nonexplosion. This is path reversal
with respect to the periodically extended invariant measure; see
\cite[Sections~3.5 and 3.7]{Norris} for invariant distributions
and time reversal.
Let $\mathcal G_t=\sigma(X_s:0\le s\le t)$ be the natural path
sigma-algebra. For stopping arguments, $\mathcal F_t$ denotes its
completed, right-continuous augmentation. When laws with the same
starting site and different parameters are compared, fix a horizon
$T<\infty$ and use the sets in $\mathcal G_T$ that are null under
all these laws to perform the completion. Lemma~\ref{lem:uniform-tilt}
shows that the finite-horizon laws have positive densities with
respect to one another, so their null sets agree. A random time $T_{\mathrm{stop}}$
is a stopping time if $\{T_{\mathrm{stop}}\le t\}\in\mathcal F_t$ for every $t\ge0$.
For $T_{\mathrm{stop}}\le T$, the sigma-algebra $\mathcal F_{T_{\mathrm{stop}}}$ consists of the events
$E\in\mathcal F_T$ such that
$E\cap\{T_{\mathrm{stop}}\le t\}\in\mathcal F_t$ for every $0\le t\le T$.
The notation $\mathbb E[Z;E]$ means $\mathbb E[Z\mathbf1_E]$.
For an integrable random variable $Z$ and an event $E$ with positive
probability, $\mathbb E[Z\mid E]=\mathbb E[Z;E]/\mathbb P(E)$.
If $E$ requires the path to remain in a prescribed allowed region
throughout $[0,t]$, the probabilities $\mathbb P_j(E,X_t=l)$ are the
entries of the absorbed kernel. Their sum over $l$ is $\mathbb P_j(E)$.

We use the strong Markov property in the following form: if $T_{\mathrm{stop}}<\infty$
is a stopping time, then, conditional on $\mathcal F_{T_{\mathrm{stop}}}$, the path
$(X_{T_{\mathrm{stop}}+r})_{r\ge0}$ has the law of a new chain started at $X_{T_{\mathrm{stop}}}$.
For example, for bounded $g:\Z\to\R$ and $r\ge0$,
\[
 \mathbb E_j[g(X_{T_{\mathrm{stop}}+r})\mid\mathcal F_{T_{\mathrm{stop}}}]
       =\sum_{l\in\Z} p_r(X_{T_{\mathrm{stop}}},l)g(l)\quad\text{for }j\in\Z.
\]
See \cite[Theorem~2.8.1 and Section~6.5]{Norris}.
For a prescribed moving boundary $b:[0,\infty)\to\R$, the future
boundary after a restart at $T_{\mathrm{stop}}$ is $r\mapsto b(T_{\mathrm{stop}}+r)$ for $r\ge0$.
Estimates after the restart must retain the stopping time $T_{\mathrm{stop}}$,
the site $X_{T_{\mathrm{stop}}}$, and the difference $X_{T_{\mathrm{stop}}}-b(T_{\mathrm{stop}})$.

\begin{lemma}\label{lem:kernel}
Fix $\lambda\in\R$ and let $p_t(j,l)$ be the kernel of the
whole-line generator $Q$ defined above. There is $C_\lambda<\infty$
such that
\begin{equation}
\label{eq:kernel-bound}
 0\le p_t(j,l)\le C_\lambda(1+t)^{-1/2}
 \text{ for } t\ge0,\ j,l\in\Z.
\end{equation}
\end{lemma}
\begin{proof}

By \eqref{eq:poisson-kernel}, $p_t(j,l)\ge0$ and
$\sum_{l\in\Z}p_t(j,l)=1$. For the decay estimate, let $z:\Z\to\C$ be finitely
supported and define its Fourier transform over one period by
\[
 \widehat z_i(\theta)=\sum_{n\in\Z}z_{i+nP}e^{-\mathrm i\theta(i+nP)},
 \quad 0\le i<P,\qquad -\pi/P\le\theta\le\pi/P.
\]
The inverse formula is
\[
 z_j=\frac{P}{2\pi}\int_{-\pi/P}^{\pi/P}
          e^{\mathrm i\theta j}\widehat z_{[j]}(\theta)\,\dd\theta\quad\text{for }j\in\Z.
\]
For $-\pi/P\le\theta\le\pi/P$, define
$Q(\theta):\C^P\to\C^P$ by the following formula, with indices
reduced modulo $P$:
\[
 (Q(\theta)v)_i=
 \alpha_i(e^{-\mathrm i\theta}v_{i-1}-v_i)
 +\beta_i(e^{\mathrm i\theta}v_{i+1}-v_i)\quad\text{for }0\le i<P,\ v\in\C^P.
\]
Changing the summation index by one in the definition of $\widehat z$
gives $\widehat{Qz}(\theta)=Q(\theta)\widehat z(\theta)$.
Since
$\|Qz\|_1\le2(\max_{0\le i<P}\alpha_i+\max_{0\le i<P}\beta_i)\|z\|_1$,
this identity extends to $\ell^1(\Z)$ by density and then to the
norm-convergent exponential series. Thus
$\widehat{e^{tQ}z}(\theta)=e^{tQ(\theta)}\widehat z(\theta)$,
and Fourier inversion applies in particular to $z=\delta_l$.
Extend $M_\lambda$ to complex parameters by the same entry formula
\eqref{eq:matrix}. With $\Phi=\diag(\phi_i(\lambda))$,
\[
 Q(\theta)=\Phi^{-1}M_{\lambda-\mathrm i\theta}\Phi-k(\lambda)I.
\]
The phase matrix $Q(0)$ is similar to $M_\lambda-k(\lambda)I$.
Lemma~\ref{lem:matrix} therefore implies that zero is simple and
that all other eigenvalues of $Q(0)$ have negative real parts. Near
$\lambda$ in $\C$, let $k$ denote the analytic continuation of the
simple eigenvalue $k(\lambda)$ of $M_\lambda$. The corresponding
eigenvalue of $Q(\theta)$ is
\[
 \zeta(\theta)=k(\lambda-\mathrm i\theta)-k(\lambda)
 =-\mathrm i k'(\lambda)\theta
   -\tfrac12 k''(\lambda)\theta^2+O(|\theta|^3)
 \quad\text{ as } \theta\to0.
\]
By \eqref{eq:convex-symmetry}, $k''(\lambda)>0$. The Taylor
expansion therefore permits a choice of $0<\epsilon<\pi/P$ and $d>0$
such that $\re\zeta(\theta)\le-d\theta^2$ for
$|\theta|\le\epsilon$. Analyticity bounds the spectral projection
onto $\zeta(\theta)$ on this interval. Since the remaining eigenvalues
of $Q(0)$ have strictly negative real parts, continuity keeps the
complementary spectrum in a fixed strict left half-plane after
$\epsilon$ is decreased. Hence, for some $\gamma>0$ and $C<\infty$ depending on
the fixed $\lambda$,
\begin{equation}
\label{eq:fourier-small}
 \|e^{tQ(\theta)}\|_\infty\le Ce^{-dt\theta^2}+Ce^{-\gamma t}
 \text{ for } |\theta|\le\epsilon,\ t\ge0.
\end{equation}
The term $Ce^{-\gamma t}$ in \eqref{eq:fourier-small} follows from
the contour formula for the exponential on the complementary spectral subspace, after
decreasing $\gamma$ if necessary.

We next show that every eigenvalue of $Q(\theta)$ has strictly
negative real part when $0<|\theta|\le\pi/P$. If
$Q(\theta)v=\sigma v$ with $v\ne0$, choose $i$ with
$|v_i|=\max_{0\le h<P}|v_h|>0$. Dividing the $i$th equation by $v_i$
and taking real parts gives $\re\sigma\le0$.
If $\re\sigma=0$, the two ratios
$e^{-\mathrm i\theta}v_{i-1}/v_i$ and
$e^{\mathrm i\theta}v_{i+1}/v_i$ have modulus at most one and
real part one, so both equal one. In particular,
$v_{i+1}=e^{-\mathrm i\theta}v_i$ and $|v_{i+1}|=|v_i|$.
Applying the eigenvalue equation at successive components of maximal modulus gives
$v_{i+P}=e^{-\mathrm iP\theta}v_i=v_i$. Hence
$e^{\mathrm iP\theta}=1$, which in the stated interval requires
$\theta=0$.
By continuity of the spectrum and compactness, after decreasing
$\gamma>0$ the spectra for $\epsilon\le|\theta|\le\pi/P$ lie in
$\{z:\re z\le-2\gamma\}$. Choose a common rectangular contour
$\Gamma$ surrounding these spectra and contained in
$\{z:\re z\le-\gamma\}$. The resolvent
$(zI-Q(\theta))^{-1}$ is continuous and uniformly bounded on the
compact set $\Gamma\times\{\theta:\epsilon\le|\theta|\le\pi/P\}$.
The contour formula for the exponential therefore gives
$\|e^{tQ(\theta)}\|_\infty\le Ce^{-\gamma t}$ on this frequency range. This reasoning also covers $P=1$ and $P=2$, when repeated
indices in the matrix are added.

Fourier inversion bounds each kernel entry by a constant multiple of
$\int_{-\pi/P}^{\pi/P}\|e^{tQ(\theta)}\|_\infty\,\dd\theta$.
By \eqref{eq:fourier-small}, the contribution of
$|\theta|\le\epsilon$ is at most $Ct^{-1/2}+Ce^{-\gamma t}$ for
$t\ge1$. The complementary frequencies contribute at most
$Ce^{-\gamma t}$ by the contour estimate. Combining these bounds
with $p_t(j,l)\le1$ for $0\le t\le1$ proves
\eqref{eq:kernel-bound}.
\end{proof}

\begin{lemma}\label{lem:local-gaussian}
Fix $\lambda\in\R$, set $c=k'(\lambda)$ and
$\varsigma^2=k''(\lambda)>0$, and let $X$ have generator $Q$.
Write $p_t(j,l)=\mathbb P_j(X_t=l)$ for its transition
probabilities and $\pi$ for its stationary phase probability vector.
Then
\begin{equation}
 \label{eq:local-gaussian}
 \lim_{t\to\infty}\sup_{j,l\in\Z}
 \left|\sqrt t\,p_t(j,l)
 -\frac{P\pi_{[l]}}{\sqrt{2\pi\varsigma^2}}
       \exp\!\left(-\frac{(l-j+ct)^2}{2\varsigma^2t}\right)
 \right|=0.
\end{equation}
Consequently, for every $A>0$ there are $b_A>0$ and $t_A<\infty$,
depending on $A$ and the fixed $\lambda$, such that
$p_t(j,l)\ge b_A t^{-1/2}$ for $j,l\in\Z$ whenever $t\ge t_A$ and
$|l-j+ct|\le A\sqrt t$.
\end{lemma}

\begin{proof}

For the Fourier matrices $Q(\theta)$ defined in the proof of
Lemma~\ref{lem:kernel}, choose $\epsilon,d,\gamma>0$ as in
\eqref{eq:fourier-small}. Set
$\zeta(\theta)=k(\lambda-\mathrm i\theta)-k(\lambda)$ for
$|\theta|\le\epsilon$, and let $\Pi(\theta)$ be its spectral
projection for $Q(\theta)$. On the space of
$P$-periodic vectors, $Q\mathbf1=0$, $\pi^TQ=0$ and
$\pi^T\mathbf1=1$. Thus $\Pi(0)=\mathbf1\pi^T$, and analyticity gives
$\Pi(\theta)=\Pi(0)+O(|\theta|)$ as $\theta\to0$.
For $i=[j]$ and $h=[l]$, Fourier inversion gives
\[
 p_t(j,l)=\frac{P}{2\pi}\int_{-\pi/P}^{\pi/P}
       e^{\mathrm i\theta(j-l)}[e^{tQ(\theta)}]_{ih}\,\dd\theta\quad\text{for }t\ge0,\ j,l\in\Z.
\]
The contour formula on the complementary spectral subspace in
Lemma~\ref{lem:kernel} gives
\[
 \|e^{tQ(\theta)}(I-\Pi(\theta))\|_\infty\le Ce^{-\gamma t}
 \quad\text{for }|\theta|\le\epsilon,\ t\ge0.
\]
On $\epsilon\le|\theta|\le\pi/P$, Lemma~\ref{lem:kernel} also gives
$\|e^{tQ(\theta)}\|_\infty\le Ce^{-\gamma t}$ for $t\ge0$.
These two norm bounds make the complementary spectral integral and
the integral over $|\theta|\ge\epsilon$ exponentially small,
uniformly in $j,l$. In the remaining integral put $\eta=\sqrt t\,\theta$.
The function
\[
 \mathbf1_{\{|\eta|\le\epsilon\sqrt t\}}
 e^{t\zeta(\eta/\sqrt t)+\mathrm i c\eta\sqrt t}
       \Pi(\eta/\sqrt t)
\]
converges as $t\to\infty$ in $L^1(\R)$, entry by entry, to
$e^{-\varsigma^2\eta^2/2}\mathbf1\pi^T$.
Pointwise convergence follows from the Taylor expansion of $\zeta$.
Moreover, $\re(t\zeta(\eta/\sqrt t))\le-d\eta^2$ on the
integration range, and $\Pi(\eta/\sqrt t)$ is uniformly bounded.
Dominated convergence therefore gives the asserted $L^1$ convergence.
The factor $\exp(-\mathrm i\eta(l-j+ct)/\sqrt t)$ has modulus
one. Hence the $L^1$ error bounds the inversion error independently
of $j,l$. Evaluating the Fourier transform of
$e^{-\varsigma^2\eta^2/2}$ now gives \eqref{eq:local-gaussian},
with uniformity over both integer indices. The factor $P$ comes from
integration over the phase Fourier interval, and the terminal factor
$\pi_{[l]}$ comes from the left eigenvector. For fixed $A>0$, take
$b_A=\frac{P\min_{0\le i<P}\pi_i}{2\sqrt{2\pi\varsigma^2}}
 e^{-A^2/(2\varsigma^2)}>0$.
On $|l-j+ct|\le A\sqrt t$, the Gaussian term in
\eqref{eq:local-gaussian} is at least $2b_A$. By the uniform convergence
in \eqref{eq:local-gaussian}, choose $t_A\ge1$ so that its error is at
most $b_A$ for every $t\ge t_A$ and every $j,l\in\Z$.
Then $p_t(j,l)\ge b_A t^{-1/2}$ whenever
$t\ge t_A$ and $|l-j+ct|\le A\sqrt t$.
For the Fourier method in local limit estimates, compare
\cite[Chapter~2]{LawlerLimic}.\qedhere

\end{proof}

\subsection{Absorption at a straight line}

For an integer-valued path $X=(X_s)_{s\ge0}$ with $X_0=j$,
$y>0$ and $v\in\R$, survival up to $t\ge0$ means
$y+X_s-j-vs>0$ for every $s\in[0,t]$. This condition includes
every interval between successive jumps. The next lemma estimates
the probability of this event for a chain with positive periodic
nearest-neighbor jump rates.

\begin{lemma}\label{lem:periodic-survival}
Let $(A_i)_{i\in\Z}$ and $(B_i)_{i\in\Z}$ be positive
$P$-periodic rates. Define $\mathcal Q$ on
$\ell^\infty(\Z)$ by
\[
 (\mathcal Qg)_j=A_j(g_{j-1}-g_j)+B_j(g_{j+1}-g_j)
 \quad\text{ for }g\in\ell^\infty(\Z),\ j\in\Z.
\]
{Let $X$ be the continuous-time chain with generator $\mathcal Q$.
Let $\varpi=(\varpi_i)_{i=0}^{P-1}$ be the stationary probability
vector of its phase process $[X_t]$ on $\{0,\ldots,P-1\}$. Set}
$v=\sum_{i=0}^{P-1}\varpi_i(B_i-A_i)$.
{Write $\mathbb P_j$ and $\mathbb E_j$ for the law and expectation
when $X_0=j$. For every}
$j\in\Z$ and $y>0$,
\begin{equation}
 \label{eq:periodic-survival}
 \mathbb P_j\{y+X_s-j-vs>0\text{ for all }0\le s\le t\}
       \le C\frac{1+y}{\sqrt{1+t}},\quad\text{ for }t\ge0.
\end{equation}
The constant $C$ depends only on the rates
$\{A_i,B_i:0\le i<P\}$ and is independent of $j$, $y$ and $t$.
\end{lemma}

\begin{proof}

Put $m_i=B_i-A_i$ for $0\le i<P$. On $P$-periodic sequences,
identify $\mathcal Q$ with the $P\times P$ matrix obtained by
reducing the spatial indices modulo $P$ and adding terms whose
indices coincide. Put $m=(m_i)_{i=0}^{P-1}$ and
$\mathbf1=(1,\ldots,1)^T\in\R^P$. We seek
$d=(d_i)_{i=0}^{P-1}\in\R^P$ satisfying
\begin{equation}
 \label{eq:corrector}
 \mathcal Qd=v\mathbf1-m,\qquad
       \sum_{i=0}^{P-1}\varpi_i d_i=0.
\end{equation}
By the definition of $v$, the vector $v\mathbf1-m$ has zero
$\varpi$-mean. Irreducibility gives
$\ker\mathcal Q=\operatorname{span}\{\mathbf1\}$, while
$\varpi^T\mathcal Q=0$ and rank $P-1$ give
$\operatorname{Ran}\mathcal Q=
\{u\in\R^P:\sum_{i=0}^{P-1}\varpi_i u_i=0\}$. Thus
\eqref{eq:corrector} has a unique solution with zero $\varpi$-mean.
For $P=1$, the solution is $d=0$.
Extend $d$ periodically and define $F:\Z\to\R$ by
$F(j)=j+d_{[j]}$. Evaluating the finite-neighbor formula for
$\mathcal Q$ at $F$ gives $(\mathcal QF)(j)=v$ for every
$j\in\Z$. Let $(\mathcal F_t)_{t\ge0}$ be the completed, right-continuous natural
filtration of $X$. The process $M_t=F(X_t)-F(j)-vt$ satisfies
$\mathbb E_jM_t^2<\infty$ and
$\mathbb E_j[M_t\mid\mathcal F_s]=M_s$ for $0\le s\le t$;
hence it is a square-integrable martingale. Its predictable quadratic variation is
\begin{equation}
 \label{eq:quadratic-rate}
 \langle M\rangle_t=\int_0^t V_{[X_s]}\,\dd s\quad\text{ for }t\ge0,\qquad
 V_i=A_i(-1+d_{i-1}-d_i)^2+B_i(1+d_{i+1}-d_i)^2\quad\text{ for }0\le i<P.
\end{equation}
Here $\langle M\rangle$ denotes the continuous increasing process,
starting at zero, for which $M_t^2-\langle M\rangle_t$ is a martingale.
We justify both martingale identities by stopping the Poisson
construction on a finite state space. In the whole-line construction
\eqref{eq:poisson-kernel}, replace $Q$ by $\mathcal Q$ and take
$\Lambda=\max_{0\le i<P}(A_i+B_i)>0$. Let $(Y_n)_{n\ge0}$
have transition matrix $I+\mathcal Q/\Lambda$ on $\Z$ and
$Y_0=j$. Let $(N_t)_{t\ge0}$ be an independent Poisson process
of rate $\Lambda$, so $X_t=Y_{N_t}$. Let $(\mathcal F_t^{Y,N})_{t\ge0}$ be the
completed, right-continuous natural filtration of $(N_t,Y_{N_t})$.
For $N\in\N$, the time $\tau_N=\inf\{t\ge0:N_t\ge N\}$ is an
$(\mathcal F_t^{Y,N})$-stopping time. The stopped pair $(N_{t\wedge\tau_N},X_{t\wedge\tau_N})$ takes values
in the finite set $\{0,\ldots,N\}\times\{j-N,\ldots,j+N\}$.
Before $\tau_N$, its generator applied to $F(x)$ and $F(x)^2$ gives
$\mathcal QF=v$ and $\mathcal Q(F^2)=2vF+V_{[x]}$.
The finite-state generator formula therefore gives both stopped
martingale identities in $(\mathcal F_t^{Y,N})$.
The bound $|F(j\pm1)-F(j)|\le C$ and the finite moments of the Poisson
count justify passage to the limit as $N\to\infty$ in $L^2$ on each
bounded time interval. Since $M$ and $M^2-\langle M\rangle$ are adapted
to $(\mathcal F_t)_{t\ge0}$, the tower property gives the two martingale
identities with respect to $(\mathcal F_t)_{t\ge0}$. This is the jump-process form of the generator
martingale formula; see \cite[Chapter~4]{EthierKurtz}.
For the stopping and maximal inequalities below, first use finite
dyadic time grids and then let the mesh tend to zero. Right-continuity
gives the pathwise limit, while domination of the jump count by a
Poisson random variable supplies an integrable bound on each finite
time interval. Thus the martingale stopping and
maximal inequalities of \cite[Sections~4.4 and 4.8]{Durrett2019}
apply after this passage to the limit; the continuous-time formulation
is also given in \cite[Chapter~2]{EthierKurtz}.

Set $V_-=\min_{0\le i<P}V_i$ and
$V_+=\max_{0\le i<P}V_i$. We verify that
$0<V_-\le V_i\le V_+<\infty$ for $0\le i<P$.
The upper bound follows from the finite phase set. For $v\ne0$,
Cauchy--Schwarz applied to $\mathcal QF=v$ and
\eqref{eq:quadratic-rate} gives $V_i\ge v^2/(A_i+B_i)>0$.
If $v=0$, set $e_i=1+d_{i+1}-d_i$ for $0\le i<P$,
with indices reduced modulo $P$. The equation for $d$ becomes
$B_i e_i=A_i e_{i-1}$. If $V_i$ vanished at one phase, this recurrence
would force every $e_i$ to vanish, contrary to
$\sum_{i=0}^{P-1}e_i=P$. The finite phase set now gives $V_->0$.
Choose $B\ge1+2\|d\|_\infty$; then every jump of $M$ has
absolute size at most $B$.

For $0<a<R$, stop $Z_t=a+M_t$ at
$\sigma_R=\inf\{t\ge0:Z_t\notin(0,R)\}$, with
$\inf\varnothing=\infty$. The paths of $Z$ are affine between
successive jumps, so $\sigma_R$ is a stopping time. At a finite exit,
a jump overshoots either endpoint by at most $B$, whereas a continuous
exit reaches the endpoint. Hence $Z_{t\wedge\sigma_R}\in[-B,R+B]$ for $t\ge0$. Define
$g_R(z)=(z+B)(R+B-z)$, so $g_R\ge0$ on the stopped range.
By \eqref{eq:quadratic-rate}, the generator of $(Z_t,[X_t])$
applied to $g_R$ equals $-V_i$ before exit at phase $i$. Indeed, put
$e_i^-=-1+d_{i-1}-d_i$ and $e_i^+=1+d_{i+1}-d_i$.
Since $A_ie_i^-+B_ie_i^+=v$, the continuous drift and jump terms
in the generator of $(Z_t,[X_t])$ give
\[
 -v g_R'(z)
 +A_i\bigl[g_R(z+e_i^-)-g_R(z)\bigr]
 +B_i\bigl[g_R(z+e_i^+)-g_R(z)\bigr]=-V_i\quad\text{for }z\in(0,R),\ 0\le i<P.
\]
By the generator formula stopped at $t\wedge\sigma_R$,
\[
 V_-\mathbb E_j(t\wedge\sigma_R)
    \le g_R(a)\le(a+B)(R+B)\quad\text{for }t\ge0.
\]
Letting $t\to\infty$ and using monotone convergence gives
$\mathbb E_j\sigma_R\le(a+B)(R+B)/V_-<\infty$, so
$\sigma_R<\infty$ almost surely. Bounded optional stopping gives
$\mathbb E_jZ_{t\wedge\sigma_R}=a$, and the bound on the stopped
values permits dominated convergence as $t\to\infty$. Thus
$\mathbb E_jZ_{\sigma_R}=a$. Upper exits have value at least $R$,
and lower exits have value at least $-B$, so
$a\ge R\mathbb P_j(Z_{\sigma_R}\ge R)-B\mathbb P_j(Z_{\sigma_R}<R)$.
It follows that
\begin{equation}\label{eq:kern-upper-exit}
 \mathbb P_j\{Z_{\sigma_R}\ge R\}
      \le\frac{a+B}{R+B}\quad\text{for }0<a<R,\ j\in\Z.
\end{equation}
On $\{a+M_s>0\text{ for every }0\le s\le t\}$, either
$\sigma_R>t$ or $Z_{\sigma_R}\ge R$. Consequently, the union bound, \eqref{eq:kern-upper-exit}, and
Markov's inequality applied to the exit-time estimate give, for $t>0$,
\[
 \begin{aligned}
 \mathbb P_j\{a+M_s>0\text{ for }0\le s\le t\}
 &\le\mathbb P_j(Z_{\sigma_R}\ge R)+\mathbb P_j(\sigma_R>t)\\
 &\le\frac{a+B}{R+B}
       +\frac{(a+B)(R+B)}{V_-t}.
 \end{aligned}
\]
For $t\ge1$ and $a<\sqrt t$, choose $R=\sqrt t$; the right-hand
side is at most $C(1+a)/\sqrt{1+t}$. If $a\ge\sqrt t$ or
$0\le t\le1$, the bound follows from a probability being at most
one after increasing $C$.
Put $D_0=\|d\|_\infty$ and $a=y+2D_0$. Since
$a+M_s=y+X_s-j-vs+2D_0+d_{[X_s]}-d_{[j]}$ and
$|d_{[X_s]}-d_{[j]}|\le2D_0$, we obtain
\[
 \begin{aligned}
 &\mathbb P_j\{y+X_s-j-vs>0\text{ for all }0\le s\le t\}\\
 &\qquad\le\mathbb P_j\{a+M_s>0\text{ for all }0\le s\le t\}
 \le C\frac{1+y+2D_0}{\sqrt{1+t}}
 \le C(1+2D_0)\frac{1+y}{\sqrt{1+t}}\\
 &\hspace{1em}\text{for }y>0,\ t\ge0,\ j\in\Z.
 \end{aligned}
\]
This proves \eqref{eq:periodic-survival}.\qedhere

\end{proof}

Return to the generator $Q$ with rates $\alpha_i,\beta_i$ at the
fixed parameter $\lambda$. In the following formulas,
$\mathbb P_j$ and $\mathbb E_j$ denote its law and expectation
from $j$. Put $h_i=\partial_\lambda\log\phi_i(\lambda)$
for $0\le i<P$.
Differentiating the eigenvector equation gives
$ (Qh)_i=k'(\lambda)+\beta_i-\alpha_i$ for $0\le i<P$.
Multiplication by
$\pi^T$ therefore gives
$v=\sum_{i=0}^{P-1}\pi_i(\beta_i-\alpha_i)=-k'(\lambda)$.
Define $d=(d_i)_{i=0}^{P-1}$ by
$d_i=-h_i+\sum_{r=0}^{P-1}\pi_rh_r$ for $0\le i<P$. Then $d$ satisfies
$ (Qd)_i=-k'(\lambda)-\beta_i+\alpha_i$ for $0\le i<P$
and $\pi^Td=0$,
so $d$ satisfies \eqref{eq:corrector} with $\mathcal Q=Q$ and $\varpi=\pi$.
Extend $d$ periodically. In the notation of Lemma~\ref{lem:periodic-survival},
take $A_i=\alpha_i$ and $B_i=\beta_i$ for $0\le i<P$, and set
$V_+=\max_{0\le i<P}V_i$ with $V_i$ given by \eqref{eq:quadratic-rate}.
For $X_0=j$, the corresponding martingale is
$M_t=X_t-j+d_{[X_t]}-d_{[j]}+k'(\lambda)t$ for $t\ge0$.
Substitution in \eqref{eq:quadratic-rate} gives
\[
 \begin{aligned}
 \sum_{i=0}^{P-1}\pi_iV_i
 &=\sum_{i=0}^{P-1}\pi_i\bigl[
   \alpha_i(1+h_{i-1}-h_i)^2
   +\beta_i(-1+h_{i+1}-h_i)^2\bigr]\\
 &=k''(\lambda),
 \end{aligned}
\]
where the last equality is \eqref{eq:kern-eigenvalue-curvature}. With $c=k'(\lambda)$,
we also have
\[
 \begin{aligned}
 \mathbb E_j(X_t-j+ct)
    &=-\mathbb E_j(d_{[X_t]}-d_{[j]}),\\
 \mathbb E_j|X_t-j+ct|^2
    &\le 2V_+t+8\|d\|_\infty^2\quad\text{for }t\ge0,\ j\in\Z.
 \end{aligned}
\]
The expectation identity follows from $\mathbb E_jM_t=0$, and the
second-moment bound follows from \eqref{eq:quadratic-rate}. For every $j\in\Z$ and $t>0$, dividing the second-moment bound by $t^2$
gives
\[
 \mathbb E_j\left|\frac{X_t-j}{t}+c\right|^2
 \le \frac{2V_+}{t}+\frac{8\|d\|_\infty^2}{t^2}
 \to0\quad\text{ as }t\to\infty.
\]
Thus $(X_t-j)/t\to-c$ in $L^2(\mathbb P_j)$ as $t\to\infty$.

\begin{lemma}\label{lem:killed-upper}

Fix $\lambda\in\R$, set $c=k'(\lambda)$, and let $X$ have generator
$Q$ and transition probabilities $p_t(j,l)=\mathbb P_j(X_t=l)$.
For $t\ge0$, $j,l\in\Z$ and $y>0$, define
\begin{equation}
 \label{eq:killing-definition}
 \begin{aligned}
 T_{j,y}&=\inf\{s\ge0:y+X_s-j+cs\le0\},\qquad X_0=j,\\
 K_t^y(j,l)&=\mathbb P_j\{X_t=l,\ T_{j,y}>t\}.
 \end{aligned}
\end{equation}
There is a constant $C_\lambda<\infty$, independent of
$t,j,l,y$, such that, for $z=y+l-j+ct>0$,
\begin{equation}
 \label{eq:killed-upper}
 K_t^y(j,l)\le C_\lambda
        \frac{(1+y)(1+z)}{(1+t)^{3/2}}.
\end{equation}
Here $\inf\varnothing=\infty$. If $z\le0$, the kernel entry is zero.

\end{lemma}
\begin{proof}
The distance $y+X_s-j+cs$ is affine between successive jumps.
Thus $T_{j,y}$ is a stopping time, and $\{T_{j,y}>t\}$ is the event
that this distance stays positive throughout $[0,t]$. The chain $X$
has velocity $-c$, so Lemma~\ref{lem:periodic-survival} gives
$\mathbb P_j(T_{j,y}>t)\le C(1+y)(1+t)^{-1/2}$.
The reversed rates in \eqref{eq:dual-rates} are positive and periodic,
have stationary phase vector $\pi$, and have velocity $c$, because
\[
 \sum_{i=0}^{P-1}\pi_i(\widehat\beta_i-\widehat\alpha_i)
     =\sum_{i=0}^{P-1}\pi_i(\alpha_i-\beta_i)=c.
\]
Thus Lemma~\ref{lem:periodic-survival} also applies to
$z+\widehat X_s-l-cs$.

For a deterministic time
$s\in[0,t]$, the Markov property gives
\[
 K_t^y(j,l)=\sum_{\substack{m\in\Z\\h_m>0}}K_s^y(j,m)K_{t-s}^{h_m}(m,l),
 \qquad h_m=y+m-j+cs\quad(m\in\Z).
\]
At time $s$, the height at site $m$ is $h_m$. After restarting
from $m$, the height at elapsed time $r\in[0,t-s]$ is
$h_m+X_{s+r}-m+cr$. Thus the second factor imposes exactly the
original path constraint on $[s,t]$. We assign value zero to a kernel
entry whose initial height is nonpositive.

Suppose $t\ge3$ and set $a=t/3$. Retain absorption in the first and
last thirds, but discard it in the middle third. The Markov property
at the two deterministic dividing times and positivity give
\[
 K_t^y(j,l)\le\sum_{m,n\in\Z}
      K_a^y(j,m)\,p_a(m,n)\,
      K_a^{\,y+n-j+2ca}(n,l),
\]
where an entry with a nonpositive initial distance is defined as zero.
By \eqref{eq:kernel-bound}, $p_a(m,n)\le C(1+a)^{-1/2}$.
By Lemma~\ref{lem:periodic-survival}, the sum
$\sum_{m\in\Z} K_a^y(j,m)$ is at most $C(1+y)(1+a)^{-1/2}$.

For the reversed chain, define, for $r\ge0$, $z>0$ and $l,n\in\Z$,
\[
 \widehat K_r^{\,z}(l,n)
 =\widehat{\mathbb P}_l\{\widehat X_r=n,
      z+\widehat X_s-l-cs>0\text{ for all }0\le s\le r\}.
\]
A surviving path in the last third runs from $n$ to $l$.
Its reversal starts at $l$ with distance $z$, and its distance at
time $s\in[0,a]$ is $z+\widehat X_s-l-cs$. Since
$z+n-l-ca=y+n-j+2ca$, the path-reversal identity
\eqref{eq:kern-path-reversal} gives
\[
 \pi_n K_a^{\,y+n-j+2ca}(n,l)
   =\pi_l\widehat K_a^{\,z}(l,n).
\]
The path identity applies to the whole survival event. If $c\ne0$,
contact with the line at a specified site can occur at only one
deterministic time. There are countably many sites, and a jump time
has a continuous distribution. Almost surely, therefore, no jump
occurs exactly when its pre-jump or post-jump position lies on the
line. Reversing the holding intervals preserves strict positivity,
and the choice of right-continuous versions changes the event only
on this null set. If $c=0$, the forbidden region is a fixed set of
sites, and reversal preserves avoidance of those sites. The
finite-path jump-density identity and the Poisson construction
therefore justify the displayed identity in both cases.
Summing and using the upper and lower positive bounds for $\pi$ gives
\[
 \sum_{n\in\Z} K_a^{\,y+n-j+2ca}(n,l)
  \le \frac{\max_{0\le i<P}\pi_i}{\min_{0\le i<P}\pi_i}
      \widehat{\mathbb P}_l
       \{z+\widehat X_s-l-cs>0\text{ for }0\le s\le a\}
  \le C(1+z)(1+a)^{-1/2}.
\]
Substitution of these bounds into the sum over $m,n$ gives
$K_t^y(j,l)\le C(1+y)(1+z)(1+a)^{-3/2}$.
Since $a=t/3$, this proves \eqref{eq:killed-upper} for $t\ge3$. For $0\le t<3$, increase $C_\lambda$ to at least $8$. Since
$(1+t)^{3/2}\le8$, $K_t^y(j,l)\le1$ and
$(1+y)(1+z)\ge1$, \eqref{eq:killed-upper} also holds on this interval. Since survival requires the terminal height $z$ to be positive,
$K_t^y(j,l)=0$ when $z\le0$.\qedhere

\end{proof}

Conditioned local limit theorems for walks with independent increments
are proved in \cite{Caravenna}. Survival and conditioned limit theorems
for sums of a fixed real-valued function along a discrete-time Markov
chain are proved, under the specified spectral hypotheses, in
\cite{GramaLauvergnatLePage2018,GramaLauvergnatLePage}.

\subsection{A lower bound away from the boundary}

We derive lower bounds for $K_t^y(j,l)$, first when both endpoint
heights are comparable to $\sqrt t$ and then when they exceed a fixed
positive constant. At stopping times we use the strong Markov property
of the chain constructed in \eqref{eq:poisson-kernel}; see
\cite[Theorem~2.8.1 and Section~6.5]{Norris}. When restarting $X$ at time $s$, we use its current height
$y+X_s-j+cs$ above the absorbing line as the new initial height.

\Needspace{9\baselineskip}
\begin{lemma}\label{lem:killed-interior-lower}

Fix $\lambda\in\R$, and use $c$, $T_{j,y}$ and $K_t^y(j,l)$ from
Lemma~\ref{lem:killed-upper}. For every $0<a\le b<\infty$ there are
$\kappa_{a,b}>0$ and $t_{a,b}<\infty$ such that
\[
 K_t^y(j,l)\ge\kappa_{a,b}t^{-1/2}
\]
whenever $j,l\in\Z$, $t\ge t_{a,b}$ and
$a\sqrt t\le y,z\le b\sqrt t$, where $z=y+l-j+ct$.

\end{lemma}
\begin{proof}

We first treat $y,z\ge H\sqrt t$ and $|z-y|\le\sqrt t$, with
$H$ to be chosen. Constants denoted by $C$ in this part of the proof
may change, but depend only on the fixed chain and are independent
of $H$ and the endpoint heights and sites. Let $\widehat X$ start at
$l$ with generator
$\widehat Q$ from \eqref{eq:dual-rates}. Its height above the reversed absorbing line is
$z+\widehat X_s-l-cs$. Let $d$ and $\widehat d$ be the normalized
$P$-periodic solutions of \eqref{eq:corrector}, with phase weights
$\varpi=\pi$, for $Q$ with velocity $-c$ and $\widehat Q$ with
velocity $c$, respectively. Put
$D_0=\max\{\|d\|_\infty,\|\widehat d\|_\infty\}$ and define
\[
 \begin{aligned}
 M_s&=X_s-j+cs+d_{[X_s]}-d_{[j]},\\
 \widehat M_s&=\widehat X_s-l-cs
                  +\widehat d_{[\widehat X_s]}-\widehat d_{[l]}.
 \end{aligned}
\]
By \eqref{eq:quadratic-rate},
$\mathbb E_jM_s^2+\widehat{\mathbb E}_l\widehat M_s^2\le Cs$ for $s\ge0$,
uniformly for $j,l\in\Z$. Applying the maximal inequality to
$M^2$ and $\widehat M^2$ on dyadic grids and then using right-continuity
therefore gives, for $r>0$,
\[
 \mathbb P_j\{\sup_{0\le s\le t/2}|M_s|\ge r\}
 +\widehat{\mathbb P}_l
       \{\sup_{0\le s\le t/2}|\widehat M_s|\ge r\}
 \le Ct/r^2\quad\text{for }t\ge0;
\]
see \cite[Section~4.4]{Durrett2019}.

Suppose $t\ge1$ and $y,z\ge H\sqrt t$, with $H\ge4D_0+1$.
By the definitions of $M$, $\widehat M$ and $D_0$, a forward contact
at time $s\le t/2$ implies
$M_s\le-y+2D_0\le-H\sqrt t/2$, and a reversed contact implies
$\widehat M_s\le-z+2D_0\le-H\sqrt t/2$.
The maximal estimate with $r=H\sqrt t/2$ therefore bounds each
crossing probability by $C/H^2$.
For paths making contact with the line by time $t/2$, the Markov
property at $t/2$ and \eqref{eq:kernel-bound} give
\[
 \begin{aligned}
 \mathbb P_j(X_t=l,\ T_{j,y}\le t/2)
 &=\mathbb E_j\bigl[\mathbf1_{\{T_{j,y}\le t/2\}}
                        p_{t/2}(X_{t/2},l)\bigr]\\
 &\le C t^{-1/2}\mathbb P_j(T_{j,y}\le t/2)
 \le C H^{-2}t^{-1/2}.
 \end{aligned}
\]
The first inequality uses
$p_{t/2}(m,l)\le Ct^{-1/2}$ from \eqref{eq:kernel-bound}, uniformly
in $m,l$ for $t\ge1$; the second uses the crossing probability
$\mathbb P_j(T_{j,y}\le t/2)\le C/H^2$ proved above.
For the reversed chain, define
$\widehat T=\inf\{s\ge0:z+\widehat X_s-l-cs\le0\}$ and let
$\widehat p$ be its unrestricted kernel. On a path from $j$ to $l$, multiplication of the reversed jump
rates contributes the factor $\pi_l/\pi_j$, and the holding factors
agree. Reversal sends a contact in $[t/2,t]$ to a contact in
$[0,t/2]$. Thus
\[
 \begin{aligned}
 &\mathbb P_j\{X_t=l,\ y+X_s-j+cs\le0
                       \text{ for some }s\in[t/2,t]\}\\
 &=\frac{\pi_l}{\pi_j}
       \widehat{\mathbb P}_l(\widehat X_t=j,\ \widehat T\le t/2)\\
 &\le\frac{\pi_l}{\pi_j}C t^{-1/2}
       \widehat{\mathbb P}_l(\widehat T\le t/2)
 \le C H^{-2}t^{-1/2}.
 \end{aligned}
\]
The penultimate inequality follows by conditioning at $t/2$ and
using $\widehat p_{t/2}\le Ct^{-1/2}$; the last uses
$\pi_l/\pi_j\le\max_i\pi_i/\min_i\pi_i$ and the bound
$\widehat{\mathbb P}_l\{\sup_{0\le s\le t/2}|\widehat M_s|\ge H\sqrt t/2\}\le C/H^2$
proved above. In particular, applying \eqref{eq:kern-path-reversal} to the
survival event gives
$\pi_jK_t^y(j,l)=\pi_l\widehat K_t^{\,z}(l,j)$, where
$\widehat K_t^{\,z}(l,j)=\widehat{\mathbb P}_l\{\widehat X_t=j,
z+\widehat X_s-l-cs>0\text{ for }0\le s\le t\}$.
A path counted by $p_t(j,l)-K_t^y(j,l)$ starts at $j$, ends at
$l$ at time $t$, and satisfies $y+X_s-j+cs\le0$ for some
$s\in[0,t]$. Such a time belongs to $[0,t/2]$ or $[t/2,t]$. The union bound therefore yields
\[
 0\le p_t(j,l)-K_t^y(j,l)\le C H^{-2}t^{-1/2}
 \quad\text{for }t\ge1,\ j,l\in\Z,\ H\ge4D_0+1.
\]
Here $z=y+l-j+ct$, $y,z\ge H\sqrt t$, and $|z-y|\le\sqrt t$.

Since $l-j+ct=z-y$, \eqref{eq:local-gaussian} gives
$p_t(j,l)\ge b_1t^{-1/2}$ for all large $t$ when
$|z-y|\le\sqrt t$. {Choose
$H\ge\max\{4D_0+1,\sqrt{2C/b_1}\}$, where $C$ is the constant
in the bound for $p_t(j,l)-K_t^y(j,l)$ above.} We have proved, uniformly in the two sites,
\begin{equation}
 \label{eq:high-distance-lower}
 K_t^y(j,l)\ge(b_1/2)t^{-1/2}
 \quad\text{if}\quad y,z\ge H\sqrt t,
 \quad |z-y|\le\sqrt t,
 \quad\text{ for }t\ge t_1.
\end{equation}

Now fix $0<a\le b<\infty$ as in the statement. Set $\eta=1/8$
and choose an integer
\[
 N\ge\max\left\{2,\left(\frac{H+\eta}{a}\right)^2,
                 \left(\frac{b-a}{1-2\eta}\right)^2\right\}.
\]
Then $a\sqrt N-\eta\ge H$ and
$(b-a)/\sqrt N+2\eta\le1$.
Write $\tau=t/N$. Define $\xi_k$ by the following formula for
$0\le k\le N$, and define $I_k$ for $1\le k<N$:
\[
 \begin{aligned}
 \xi_k&=j-ck\tau+\frac{k}{N}(z-y)
       &&\text{ for }0\le k\le N,\\
 I_k&=\{m\in\Z:|m-\xi_k|\le\eta\sqrt\tau\}
       &&\text{ for }1\le k<N.
 \end{aligned}
\]
The interpolation satisfies $\xi_0=j$ and $\xi_N=l$. Choose $m_k\in I_k$ for $1\le k<N$,
and set $m_0=j$ and $m_N=l$. For $0\le k\le N$, let
$y_k=y+m_k-j+ck\tau$ be the height of $m_k$ above the line at time $k\tau$. These distances differ from the linear interpolation between
$y$ and $z$ by at most $\eta\sqrt\tau$. Since
$a\sqrt t\le y,z\le b\sqrt t$, we have
$y_k\ge(a\sqrt N-\eta)\sqrt\tau$ and
$|y_{k+1}-y_k|\le((b-a)/\sqrt N+2\eta)\sqrt\tau$.
The choice of $N$ therefore gives
\[
 y_k\ge H\sqrt\tau\quad(0\le k\le N),\qquad |y_{k+1}-y_k|\le\sqrt\tau\quad(0\le k<N).
\]
Choose $t_{a,b}$ so that $\tau=t/N\ge\max\{t_1,\eta^{-2}\}$
whenever $t\ge t_{a,b}$. By \eqref{eq:high-distance-lower}, each
$K_\tau^{y_k}(m_k,m_{k+1})$ is then at least
$(b_1/2)\tau^{-1/2}$. Moreover, each $I_k$ contains at least
$2\eta\sqrt\tau-1\ge\eta\sqrt\tau$ integers. Applying the Markov
property at $\tau,2\tau,\ldots,(N-1)\tau$ and retaining only
$m_k\in I_k$ at time $k\tau$ gives
\[
 K_t^y(j,l)
 \ge (\eta\sqrt\tau)^{N-1}
           \bigl((b_1/2)\tau^{-1/2}\bigr)^N
 =\eta^{N-1}(b_1/2)^N\sqrt N\,t^{-1/2}.
\]
The choice of $N$ depends only on $a,b$ and the fixed chain.
Thus $\eta^{N-1}(b_1/2)^N\sqrt N$ is a fixed positive lower
constant, which proves the assertion for every stated choice of
$j,l,y,z$ and $t\ge t_{a,b}$.\qedhere

\end{proof}

\Needspace{9\baselineskip}
\begin{lemma}\label{lem:entrance-height}

Let $X$ have generator $\mathcal Q$ with positive periodic rates
$A_i,B_i$ and velocity $v$ as in Lemma~\ref{lem:periodic-survival}.
Write $\mathbb P_j$ and $\mathbb E_j$ for its law and expectation
when $X_0=j$, and let $(\mathcal F_s)_{s\ge0}$ be its completed,
right-continuous natural filtration. There are $y_0\ge1$ and
$\delta>0$, depending only
on $A_i,B_i$, with the following property. For each $A>0$ there are
$B_A<\infty$, $\kappa_A>0$ and $t_A<\infty$ such that, for every
$j\in\Z$, if $t\ge t_A$ and $y_0\le y\le A\sqrt t$,
there are an $(\mathcal F_s)$-stopping time $T_{\mathrm{stop}}$
and an event $E\in\mathcal F_{T_{\mathrm{stop}}}$ satisfying
\[
 \begin{gathered}
 E\subseteq\{T_{\mathrm{stop}}\le t/3,\ y+X_s-j-vs>0\text{ for }0\le s\le T_{\mathrm{stop}}\},\\
 \delta\sqrt t\le y+X_{T_{\mathrm{stop}}}-j-vT_{\mathrm{stop}}\le B_A\sqrt t\quad\text{on }E,
 \qquad \mathbb P_j(E)\ge\kappa_A y/\sqrt t.
 \end{gathered}
\]
The same $y_0,\delta$ may be chosen for a chain and its reversed chain.

\end{lemma}
\begin{proof}

Let $d$ be the zero-$\varpi$-mean solution of \eqref{eq:corrector},
extended periodically to $\Z$, and let
$M_s=X_s-j-vs+d_{[X_s]}-d_{[j]}$ for $s\ge0$. Let $B\ge1$ bound its jumps,
and let $V_->0$ be a lower bound for the variance rates in
\eqref{eq:quadratic-rate}. We first choose constants that also work for the reversed chain.
Extend $\varpi$ periodically to $\Z$. Define the reversed rates for $i\in\Z$
and the reversed generator on $\ell^\infty(\Z)$ by
\[
 \begin{aligned}
 \widehat A_i&=\frac{\varpi_{i-1}B_{i-1}}{\varpi_i},\qquad
 \widehat B_i=\frac{\varpi_{i+1}A_{i+1}}{\varpi_i},\\
 (\widehat{\mathcal Q}g)_j
 &=\widehat A_j(g_{j-1}-g_j)+\widehat B_j(g_{j+1}-g_j)
 \quad\text{ for }g\in\ell^\infty(\Z),\ j\in\Z.
 \end{aligned}
\]
These rates are positive and periodic, their invariant phase vector
is $\varpi$, and their velocity is
$\sum_{i=0}^{P-1}\varpi_i(\widehat B_i-\widehat A_i)=-v$.
Let $\widehat X$ have this generator. Write
$\widehat{\mathbb P}_j$ and $\widehat{\mathbb E}_j$ for its law and
expectation when $\widehat X_0=j$, and let
$(\widehat{\mathcal F}_s)_{s\ge0}$ be its completed,
right-continuous natural filtration.
Let $\widehat d$ be the normalized corrector satisfying
$(\widehat{\mathcal Q}\widehat d)_i=-v-\widehat B_i+\widehat A_i$
and $\sum_{i=0}^{P-1}\varpi_i\widehat d_i=0$.
The right-hand side $(-v-\widehat B_i+\widehat A_i)_{i=0}^{P-1}$
has zero $\varpi$-mean because the reversed velocity is $-v$.
The same range calculation as for \eqref{eq:corrector}, applied to
$\widehat{\mathcal Q}$, therefore gives the unique normalized
corrector $\widehat d$. Put
$\widehat B=1+2\|\widehat d\|_\infty$ and define
\[
 \begin{aligned}
 \widehat V_i&=
 \widehat A_i(-1+\widehat d_{i-1}-\widehat d_i)^2
 +\widehat B_i(1+\widehat d_{i+1}-\widehat d_i)^2\quad\text{ for }0\le i<P,\\
 \widehat V_-&=\min_{0\le i<P}\widehat V_i.
 \end{aligned}
\]
The triangle inequality bounds the reversed corrected jumps by
$\widehat B$. If $v\ne0$, Cauchy--Schwarz gives
$\widehat V_i\ge v^2/(\widehat A_i+\widehat B_i)>0$.
If $v=0$, put $\widehat e_i=1+\widehat d_{i+1}-\widehat d_i$.
The corrector equation gives
$\widehat B_i\widehat e_i=\widehat A_i\widehat e_{i-1}$ and
$\sum_{i=0}^{P-1}\widehat e_i=P$. If $\widehat V_i=0$ for some $0\le i<P$, then
$\widehat e_{i-1}=\widehat e_i=0$, and the recurrence would force
$\widehat e_h=0$ for every $0\le h<P$. This contradicts
$\sum_{h=0}^{P-1}\widehat e_h=P$. Thus
$\widehat V_->0$ in both cases.
Set
\[
 D_*=\max\{\|d\|_\infty,\|\widehat d\|_\infty\},\qquad
 B_*=\max\{B,\widehat B\},\qquad
 V_*=\min\{V_-,\widehat V_-\}>0.
\]
Choose $y_0\ge\max\{1,4D_*,2D_*+B_*\}$,
$\delta>0$ with $\delta^2\le V_*/48$, and
$t_A\ge\max\{1,(B_*/\delta)^2\}$.

We give the entrance construction for $X$ using these common bounds.
Fix $A>0$, $j\in\Z$, $t\ge t_A$ and $y_0\le y\le A\sqrt t$.
Put $D_0=\|d\|_\infty$, $a=y-2D_0$ and $R=\delta\sqrt t$.
For $y\ge y_0$ and $t\ge t_A$, we have $a\ge y/2$,
$a\ge B_*$ and $R\ge B_*$. If $a\ge R$, take $T_{\mathrm{stop}}=0$ and
$E$ equal to the entire sample space. The initial height then lies
in $[R,A\sqrt t]$, and $1\ge y/(A\sqrt t)$.
Suppose that $a<R$. Put $Z_s=a+M_s$ and let
$\sigma_R=\inf\{s\ge0:Z_s\notin(0,R)\}$.
For $s\ge0$, $Z_{s\wedge\sigma_R}\in[-B_*,R+B_*]$.
For $g_R(z)=(z+B_*)(R+B_*-z)$, set
$e_i^-=-1+d_{i-1}-d_i$ and $e_i^+=1+d_{i+1}-d_i$.
Before exit, the generator of $(Z_s,[X_s])$ applied to $g_R$ is
\[
 -v g_R'(z)
 +A_i[g_R(z+e_i^-)-g_R(z)]
 +B_i[g_R(z+e_i^+)-g_R(z)]
 =-V_i\quad\text{for }z\in(0,R),\ 0\le i<P.
\]
Indeed, $A_ie_i^-+B_ie_i^+=v$ by \eqref{eq:corrector},
and $g_R''(z)=-2$.
Bounded stopping and nonnegativity of $g_R$ on the stopped range give
$V_*\mathbb E_j(u\wedge\sigma_R)\le(a+B_*)(R+B_*)$.
Monotone convergence as $u\to\infty$ gives a finite expected exit
time. Optional stopping for the bounded stopped martingale, followed
by dominated convergence, also gives
\[
 \mathbb E_j\sigma_R\le (a+B_*)(R+B_*)/V_*,\qquad
 \mathbb E_jZ_{\sigma_R}=a.
\]
On a lower exit $Z_{\sigma_R}\le0$, whereas on an upper exit
$Z_{\sigma_R}\le R+B_*$. Hence
$\mathbb P_j\{Z_{\sigma_R}\ge R\}\ge a/(R+B_*)$. Take
\[
 T_{\mathrm{stop}}=\sigma_R\wedge(t/3),\qquad
 E=\{\sigma_R\le t/3,\ Z_{\sigma_R}\ge R\}.
\]
Then $T_{\mathrm{stop}}$ is bounded, $E\in\mathcal F_{T_{\mathrm{stop}}}$, and $T_{\mathrm{stop}}=\sigma_R$ on $E$.
Indeed, for $r<t/3$, the event $E\cap\{T_{\mathrm{stop}}\le r\}$ is
$\{\sigma_R\le r,Z_{\sigma_R}\ge R\}\in\mathcal F_r$,
and $E\in\mathcal F_{t/3}$. Since $a+B_*\le2a$,
$R+B_*\le2R$ and $\delta^2\le V_*/48$, Markov's inequality yields
\[
 \begin{aligned}
 \mathbb P_j(E)
 &\ge\frac{a}{R+B_*}-\frac{3(a+B_*)(R+B_*)}{V_*t}\\
 &\ge\frac{a}{R+B_*}
          \left(1-\frac{6(R+B_*)^2}{V_*t}\right)
 \ge\frac{a}{2(R+B_*)}
 \ge\frac{y}{8\delta\sqrt t}.
 \end{aligned}
\]
For $0\le s<\sigma_R$,
\[
 y+X_s-j-vs=Z_s+2D_0-d_{[X_s]}+d_{[j]}\ge Z_s>0.
\]
At $s=\sigma_R$ on $E$, the same identity places the height in
$[R,R+B_*+4D_0]$. Since $R>0$, the required strict inequality also
holds at the stopping time itself. Hence $E$ is contained in the
survival event through $T_{\mathrm{stop}}$.

For the reversed chain, the corrected martingale is
$\widehat M_s=\widehat X_s-j+vs+
\widehat d_{[\widehat X_s]}-\widehat d_{[j]}$.
Its corrector has norm at most $D_*$, its jumps have absolute size
at most $B_*$, and its variance rates are at least $V_*$.
Put $\widehat a=y-2\|\widehat d\|_\infty$. If
$\widehat a\ge R$, choose $\widehat{T}_{\mathrm{stop}}=0$ and let
$\widehat E$ be the entire path space. If $\widehat a<R$, set
$\widehat Z_s=\widehat a+\widehat M_s$,
$\widehat\sigma_R=\inf\{s\ge0:\widehat Z_s\notin(0,R)\}$,
$\widehat{T}_{\mathrm{stop}}=\widehat\sigma_R\wedge(t/3)$, and
$\widehat E=\{\widehat\sigma_R\le t/3,\ 
\widehat Z_{\widehat\sigma_R}\ge R\}$. The stopped martingale and
quadratic-variation identities give
$\widehat{\mathbb E}_j\widehat\sigma_R
\le(\widehat a+B_*)(R+B_*)/V_*$ and
$\widehat{\mathbb P}_j(\widehat Z_{\widehat\sigma_R}\ge R)
\ge\widehat a/(R+B_*)$. Since
$\widehat a\ge y/2$, $\widehat a\ge B_*$, $R\ge B_*$ and
$\delta^2\le V_*/48$, subtracting the bound
$3\widehat{\mathbb E}_j\widehat\sigma_R/t$ yields
$\widehat{\mathbb P}_j(\widehat E)\ge y/(8\delta\sqrt t)$.
In the case $\widehat a<R$, the identity
$y+\widehat X_s-j+vs
=\widehat Z_s+2\|\widehat d\|_\infty
-\widehat d_{[\widehat X_s]}+\widehat d_{[j]}$ shows positivity
for $0\le s\le\widehat{T}_{\mathrm{stop}}$ on $\widehat E$. At a successful upper exit its
terminal height lies in $[R,R+B_*+4D_*]$ and its success probability
is at least $y/(8\delta\sqrt t)$; in the time-zero case its height
is $y\in[R,A\sqrt t]$ and its success probability is one.
Choose any $B_A>\max\{A,\delta+B_*+4D_*\}$.
Since $t\ge1$, this gives the required upper height bound in both
chains, with the common probability constant
$\kappa_A=\min\{A^{-1},(8\delta)^{-1}\}$ and the same $y_0,\delta$.

\end{proof}

\subsection{Changing the exponential parameter}

Write $Q^\lambda$, $p_t^\lambda$ and $\mathbb P_j^\lambda$ for the
generator, transition kernel and law at parameter $\lambda$, and
$\mathbb E_j^\lambda$ for expectation under $\mathbb P_j^\lambda$.
For each starting site $j$, take the laws $\mathbb P_j^\lambda$
on the same space of integer-valued right-continuous paths with
left limits, with coordinate process $X$. Write
$\widehat{\mathbb P}_j^\lambda$ for the law started at $j$ of the
reversed chain whose generator is given by \eqref{eq:dual-rates}
for $Q^\lambda$.
Retain $\mathcal G_t=\sigma(X_s:0\le s\le t)$ for the natural
path sigma-algebra and $\mathcal F_t$ for the stopping filtration
specified above. On each fixed horizon $T<\infty$, this filtration
uses the common null sets in $\mathcal G_T$. The change-of-parameter
identity is first proved on $\mathcal G_t$ and extends to
$\mathcal F_t$ for $t<T$ by the finite-horizon density and its
right-continuity. Choosing $T>t$ gives the identity at every finite $t$.

\begin{lemma}\label{lem:uniform-tilt}

On every compact interval $I\subset\R$, the constants in
Lemmas~\ref{lem:kernel},~\ref{lem:local-gaussian},
\ref{lem:periodic-survival},~\ref{lem:killed-upper},
\ref{lem:killed-interior-lower} and~\ref{lem:entrance-height}, for
the process with generator $Q^\lambda$ and its reversal, may be chosen uniformly for
$\lambda\in I$; the convergence in \eqref{eq:local-gaussian} is
also uniform on $I$. For $t\ge0$, an event $E\in\mathcal F_t$ and $j,l\in\Z$,
\begin{equation}
 \label{eq:path-tilt-identity}
 \mathbb P_j^\lambda(E,\ X_t=l)
 =e^{(k(\eta)-k(\lambda))t+(\eta-\lambda)(l-j)}
  \frac{\phi_j(\eta)\phi_l(\lambda)}
       {\phi_j(\lambda)\phi_l(\eta)}
  \mathbb P_j^\eta(E,\ X_t=l),\qquad \lambda,\eta\in\R.
\end{equation}

\end{lemma}
\begin{proof}

Put $J=I\cup\{\log q-r:r\in I\}$. This compact set contains the
parameters of both the original chains and their reversals. Let $\mathsf Q^\lambda$ be the $P\times P$ phase matrix of
the whole-line generator $Q^\lambda$, and let $\pi^\lambda$
be its stationary probability vector. The rates $\alpha_i^\lambda,\beta_i^\lambda$
and the weights $\pi_i^\lambda$ are positive and continuous, so
they have common positive lower and finite upper bounds on $J$.
For each $\lambda\in J$, put
$m=(m_i)_{i=0}^{P-1}$ with $m_i=\beta_i^\lambda-\alpha_i^\lambda$,
and $v=\sum_{i=0}^{P-1}\pi_i^\lambda m_i=-k'(\lambda)$.
Then the matrix
$\mathsf Q^\lambda-\mathbf1(\pi^\lambda)^T$ is invertible. Indeed,
$\mathsf Q^\lambda$ has kernel spanned by $\mathbf1$ and range
$\{z:(\pi^\lambda)^Tz=0\}$; multiplication by $(\pi^\lambda)^T$
shows that a vector in the kernel of
$\mathsf Q^\lambda-\mathbf1(\pi^\lambda)^T$ has zero phase mean and
then is zero. Consequently
\[
 d^\lambda=
  \bigl(\mathsf Q^\lambda-\mathbf1(\pi^\lambda)^T\bigr)^{-1}
       (v\mathbf1-m)\quad\text{for }\lambda\in J
\]
satisfies $\mathsf Q^\lambda d^\lambda=v\mathbf1-m$ and
$(\pi^\lambda)^Td^\lambda=0$. The inverse is continuous on $J$,
so $D_I:=\sup_{\lambda\in J}\|d^\lambda\|_\infty<\infty$.
Set $B_I=1+2D_I$. Since $\|d^\lambda\|_\infty\le D_I$, the quantities
$-1+d_{i-1}^\lambda-d_i^\lambda$ and $1+d_{i+1}^\lambda-d_i^\lambda$
have absolute size at most $B_I$.
The rates
\[
 V_i^\lambda=
 \alpha_i^\lambda(-1+d_{i-1}^\lambda-d_i^\lambda)^2
 +\beta_i^\lambda(1+d_{i+1}^\lambda-d_i^\lambda)^2\quad\text{for }0\le i<P,\ \lambda\in J
\]
are continuous. To verify their positivity, set
$F^\lambda(j)=j+d_{[j]}^\lambda$. Evaluating the finite-neighbor formula for $Q^\lambda$ at
$F^\lambda$ gives $(Q^\lambda F^\lambda)(j)=v$ for every $j\in\Z$. If one $V_i^\lambda$ vanished, positivity of $\alpha_i^\lambda$ and
$\beta_i^\lambda$ would give
$-1+d_{i-1}^\lambda-d_i^\lambda=1+d_{i+1}^\lambda-d_i^\lambda=0$,
so $Q^\lambda F^\lambda=v$ would give $v=0$.
Put $e_h=1+d_{h+1}^\lambda-d_h^\lambda$ for $0\le h<P$,
with indices taken modulo $P$. At the vanishing phase $i$,
$e_{i-1}=e_i=0$. The relation
$\beta_h^\lambda e_h=\alpha_h^\lambda e_{h-1}$ for $0\le h<P$ and positivity
of the rates then force $e_h=0$ for every $0\le h<P$, contrary to
$\sum_{h=0}^{P-1}e_h=P$. Hence there are common constants
$0<V_I\le V_i^\lambda\le V_I^+<\infty$ for $0\le i<P$ and $\lambda\in J$.

For $\lambda\in J$ and $\theta\in[-\pi/P,\pi/P]$,
let $\xi\in\C^P\setminus\{0\}$ be an eigenvector of
$Q^\lambda(\theta)$ with eigenvalue $\sigma$. Choose
$i\in\{0,\ldots,P-1\}$ such that
$|\xi_i|=\max_{0\le h<P}|\xi_h|$.
At this maximizing coordinate $i$, the eigenvalue equation gives
\[
 \re\sigma
 =\alpha_i^\lambda\left(
       \re\frac{e^{-\mathrm i\theta}\xi_{i-1}}{\xi_i}-1\right)
 +\beta_i^\lambda\left(
       \re\frac{e^{\mathrm i\theta}\xi_{i+1}}{\xi_i}-1\right)
 \le0.
\]
If $\re\sigma=0$, both ratios in this display have modulus at
most one and real part one. Hence
$e^{-\mathrm i\theta}\xi_{i-1}=\xi_i$ and
$e^{\mathrm i\theta}\xi_{i+1}=\xi_i$. The neighboring components
also have maximal modulus, so the same equality propagates through
all $P$ phases and gives $e^{\mathrm iP\theta}=1$. This is impossible
for $0<|\theta|\le\pi/P$. For each $0<\epsilon<\pi/P$,
compactness of $J\times\{\theta:\epsilon\le|\theta|\le\pi/P\}$
therefore gives a common strictly negative spectral upper bound. Choose $\gamma_I>0$ so that these spectra lie in
$\{w:\re w\le-2\gamma_I\}$. A common enclosing
rectangle with right side $\re w=-\gamma_I$ has a
uniformly bounded resolvent on its boundary. The contour formula
therefore gives $\|e^{tQ^\lambda(\theta)}\|_\infty
\le C_Ie^{-\gamma_It}$ for $\epsilon\le|\theta|\le\pi/P$.

Near $\theta=0$, let $\zeta_\lambda(\theta)$ and
$\Pi_\lambda(\theta)$ continue the simple eigenvalue zero and its
spectral projection $\mathbf1(\pi^\lambda)^T$.
After decreasing $\epsilon>0$ if necessary, analytic perturbation
and compactness of $J$ give the following uniform Taylor bounds for
$|\theta|\le\epsilon$ and $\lambda\in J$:
\[
 \begin{aligned}
 \zeta_\lambda(\theta)
 &=-\mathrm i k'(\lambda)\theta
        -\tfrac12 k''(\lambda)\theta^2+R_\lambda(\theta),\\
 |R_\lambda(\theta)|&\le C_I|\theta|^3,\\
 \|\Pi_\lambda(\theta)-\mathbf1(\pi^\lambda)^T\|_\infty
 &\le C_I|\theta|.
 \end{aligned}
\]
At $\theta=0$ the other eigenvalues have negative real parts;
compactness keeps them in a common strict left half-plane for small
$|\theta|$. The contour formula bounds the complementary semigroup
by $C_Ie^{-\gamma_It}$. Since
$\min_{\lambda\in J}k''(\lambda)>0$, decreasing the common
neighborhood gives constants $d_I,\gamma_I>0$ such that
\[
 \|e^{tQ^\lambda(\theta)}\|_\infty
 \le C_Ie^{-d_It\theta^2}+C_Ie^{-\gamma_It}
 \quad\text{ for }|\theta|\le\epsilon,\ t\ge0,\ \lambda\in J.
\]
Fourier inversion and integration in $\theta$ give
$p_t^\lambda(j,l)\le C_I(1+t)^{-1/2}$ uniformly for $\lambda\in J$.
After the substitution $r=\sqrt t\,\theta$, the matrix
\[
 \mathbf1_{\{|r|\le\epsilon\sqrt t\}}
 e^{t\zeta_\lambda(r/\sqrt t)+\mathrm i k'(\lambda)r\sqrt t}
       \Pi_\lambda(r/\sqrt t)
\]
converges as $t\to\infty$, uniformly on every bounded $r$-interval
and for $\lambda\in J$, to
$e^{-k''(\lambda)r^2/2}\mathbf1(\pi^\lambda)^T$.
The displayed matrix has norm at most $C_Ie^{-d_Ir^2}$ for $r\in\R$,
$\lambda\in J$ and all sufficiently large $t$. The limit
$e^{-k''(\lambda)r^2/2}\mathbf1(\pi^\lambda)^T$ also has an
integrable Gaussian bound uniform for $\lambda\in J$, because
$\min_{\lambda\in J}k''(\lambda)>0$ and
$\sup_{\lambda\in J}\|\mathbf1(\pi^\lambda)^T\|_\infty<\infty$.
The integrals of the norms of both matrices over $\{|r|>R\}$ therefore
tend to zero as $R\to\infty$, uniformly for $\lambda\in J$ and all sufficiently large $t$.
The convergence on each bounded $r$-interval is uniform in $\lambda$.
Fourier inversion gives the asserted uniformity in
\eqref{eq:local-gaussian}, because the remaining oscillatory factor
in the inversion integral has modulus one for every $j,l$.

For the chain with law $\mathbb P_j^\lambda$, put
$M_t^\lambda=X_t-j+k'(\lambda)t
+d_{[X_t]}^\lambda-d_{[j]}^\lambda$. By \eqref{eq:quadratic-rate},
its predictable quadratic variation is
$\langle M^\lambda\rangle_t=\int_0^t V_{[X_s]}^\lambda\,\dd s$.
For $0<a<R$, define
$\sigma_R=\inf\{s\ge0:a+M_s^\lambda\notin(0,R)\}$.
The jump bound $B_I$ and quadratic-variation lower bound $V_I$ give
\[
 \mathbb E_j^\lambda\sigma_R\le\frac{(a+B_I)(R+B_I)}{V_I},\qquad
 \mathbb P_j^\lambda(a+M_{\sigma_R}^\lambda\ge R)
       \le\frac{a+B_I}{R+B_I}\quad\text{for }\lambda\in J,\ j\in\Z,\ 0<a<R.
\]
For $y>0$, put $a=y+2D_I$. Then
$a+M_s^\lambda\ge y+X_s-j+k'(\lambda)s$ for all $s\ge0$.
Taking $R=\sqrt t$ in the exit estimate, and using the bound one
when $a\ge\sqrt t$ or $t\le1$, therefore gives
$\mathbb P_j^\lambda\{y+X_s-j+k'(\lambda)s>0
\text{ for }0\le s\le t\}\le C_I(1+y)/\sqrt{1+t}$.
The process $\widehat X$ has parameter $\log q-\lambda\in J$, so
Lemma~\ref{lem:periodic-survival} has the same constant $C_I$ for $\widehat X$. For
$z=y+l-j+k'(\lambda)t>0$, the Markov property at $t/3$ and $2t/3$
combines the forward and reversed survival bounds
$C_I(1+y)(1+t)^{-1/2}$ and $C_I(1+z)(1+t)^{-1/2}$ with the
bound $p_{t/3}^\lambda(m,n)\le C_I(1+t)^{-1/2}$ for
$m,n\in\Z$ on the middle time interval $[t/3,2t/3]$.
Since
$\max_{0\le i<P}\pi_i^\lambda/\min_{0\le i<P}\pi_i^\lambda$
is uniformly bounded for $\lambda\in J$, we obtain
$K_t^y(j,l)\le C_I(1+y)(1+z)/(1+t)^{3/2}$.

Fix $0<a_0\le b_0<\infty$ for the interior lower bound. The
uniform version of \eqref{eq:local-gaussian} gives
$p_t^\lambda(j,l)\ge b_1t^{-1/2}$ on $|z-y|\le\sqrt t$, with a
common $b_1>0$ and time threshold. The bound $V_i^\lambda\le V_I^+$ gives a common second-moment
bound for the forward and reversed corrected martingales. If
$t\ge1$, $y,z\ge H\sqrt t$ and $H\ge4D_I+1$, a contact of
$y+X_s-j+k'(\lambda)s$ or
$z+\widehat X_s-l-k'(\lambda)s$ with zero before $t/2$
requires the corresponding corrected martingale to have absolute
value at least $H\sqrt t-2D_I\ge H\sqrt t/2$.
The maximal inequalities therefore bound each contact probability
by $C_IH^{-2}$. For contacts in $[0,t/2]$, conditioning at $t/2$
multiplies this bound by
$p_{t/2}^\lambda(m,l)\le C_It^{-1/2}$, uniformly in $m,l$.
For contacts in $[t/2,t]$, reverse the path and condition at $t/2$
under the reversed law. The same kernel bound and the uniformly
bounded invariant-weight ratio then give
$p_t^\lambda(j,l)-K_t^y(j,l)\le C_IH^{-2}t^{-1/2}$.
Choose $H\ge\max\{4D_I+1,\sqrt{2C_I/b_1}\}$, where $C_I$
is the constant in the bound for
$p_t^\lambda(j,l)-K_t^y(j,l)$, and then an integer $N\ge2$
with $a_0\sqrt N-1/8\ge H$ and
$(b_0-a_0)/\sqrt N+1/4\le1$. For endpoint heights in $[a_0\sqrt t,b_0\sqrt t]$, set
$\tau=t/N$, $c=k'(\lambda)$ and
$\xi_k=j-ck\tau+(k/N)(z-y)$ for $0\le k\le N$.
For $1\le k<N$, put
$I_k=\{m\in\Z:|m-\xi_k|\le\sqrt\tau/8\}$.
Choose $m_k\in I_k$ for $1\le k<N$, set $m_0=j$ and $m_N=l$, and define
$y_k=y+m_k-j+ck\tau$ for $0\le k\le N$.
The choice of $N$ then gives $y_k\ge H\sqrt{t/N}$ and
$|y_{k+1}-y_k|\le\sqrt{t/N}$ for $0\le k<N$. Thus each factor $K_{t/N}^{y_k}(m_k,m_{k+1})$ is at
least $(b_1/2)(t/N)^{-1/2}$ once $t/N$ exceeds the common time
threshold in Lemma~\ref{lem:local-gaussian} and is at least one.
Each $I_k$ has radius $\sqrt{t/N}/8$ and, for $t/N\ge64$, contains
at least $\sqrt{t/N}/8$ integers. Iterating the killed-kernel
composition rule and summing over these intermediate sites gives
the common lower constant $(1/8)^{N-1}(b_1/2)^N\sqrt N$ in front
of $t^{-1/2}$.
For Lemma~\ref{lem:entrance-height}, choose
\[
 y_0\ge\max\{1,4D_I,2D_I+B_I\},\qquad
 \delta>0,\qquad \delta^2\le V_I/48,\qquad
 t_A\ge\max\{1,(B_I/\delta)^2\}.
\]
Fix $A>0$, $j\in\Z$, $t\ge t_A$ and $y_0\le y\le A\sqrt t$.
Set $a=y-2D_I$ and
$R=\delta\sqrt t$; then $a\ge y/2$ and $a\ge B_I$.
If $a\ge R$, take $T_{\mathrm{stop}}=0$ and let $E$ be the entire path space.
If $a<R$, define
\[
 \begin{aligned}
 \sigma_R&=\inf\{s\ge0:a+M_s^\lambda\notin(0,R)\},\\
 T_{\mathrm{stop}}&=\sigma_R\wedge(t/3),\\
 E&=\{\sigma_R\le t/3,\ a+M_{\sigma_R}^\lambda\ge R\}.
 \end{aligned}
\]
The exit-time and upper-exit probability estimates in the proof of
Lemma~\ref{lem:entrance-height}, with jump bound $B_I$ and variance
lower bound $V_I$, give $\mathbb P_j^\lambda(E)\ge
\min\{A^{-1},(8\delta)^{-1}\}y/\sqrt t$ in both cases.
On $E$, the distance $y+X_{T_{\mathrm{stop}}}-j+k'(\lambda)T_{\mathrm{stop}}$ lies between
$\delta\sqrt t$ and $B_A\sqrt t$ for any
$B_A>\max\{A,\delta+B_I+4D_I\}$.
All choices are uniform on $J$, which proves the six uniform
assertions for the chains with generator $Q^\lambda$ and their
reversals.

To prove \eqref{eq:path-tilt-identity}, fix $t\ge0$, $j\in\Z$,
$\lambda,\eta\in\R$ and an arbitrary event $E\in\mathcal G_t$. For a path from $j$ to $l$
with finitely many jumps, the ratio of its jump factors under parameters
$\lambda$ and $\eta$ telescopes to
\[
 e^{(\eta-\lambda)(l-j)}
 \frac{\phi_l(\lambda)\phi_j(\eta)}
      {\phi_j(\lambda)\phi_l(\eta)}.
\]
The eigenvector equation gives the total rate at phase $i$ as
$\alpha_i^\lambda+\beta_i^\lambda=k(\lambda)-a_i+D(q+1)$.
The ratio of the holding-time factors is consequently
$e^{(k(\eta)-k(\lambda))t}$. Integration over jump times and summation
over the number of jumps prove the identity, since the number is
finite almost surely. Since $E$ is an arbitrary event in $\mathcal G_t$, the identity also
applies when $E$ imposes absorption throughout $[0,t]$. The density on
$\mathcal G_t$ is
\[
 \left.\frac{\dd\mathbb P_j^\lambda}
                  {\dd\mathbb P_j^\eta}\right|_{\mathcal G_t}
 =e^{(k(\eta)-k(\lambda))t+(\eta-\lambda)(X_t-j)}
   \frac{\phi_j(\eta)\phi_{X_t}(\lambda)}
        {\phi_j(\lambda)\phi_{X_t}(\eta)}\quad\text{for }t\ge0,\ j\in\Z.
\]
Integrating over all paths gives expectation one for this positive
density under $\mathbb P_j^\eta$. Interchanging $\lambda$ and $\eta$
gives mutual absolute continuity on $\mathcal G_t$ for every finite $t$.
To pass to the stopping filtration, fix $T>t$ and let $Z_s$ denote
the displayed density at time $s\le T$. Restricting the density
identity from $\mathcal G_T$ to $\mathcal G_s$ gives
$Z_s=\mathbb E_j^\eta[Z_T\mid\mathcal G_s]$, so
$(Z_s)_{0\le s\le T}$ is uniformly integrable. For
$E\in\mathcal F_t$ and $t<r<T$, the event $E$ is measurable in
$\mathcal G_r$ completed by the common null sets in $\mathcal G_T$.
The density identity therefore gives
$\mathbb P_j^\lambda(E)=\mathbb E_j^\eta[Z_r;E]$.
The explicit density has right-continuous paths. Letting $r\downarrow t$
and using uniform integrability yields
$\mathbb P_j^\lambda(E)=\mathbb E_j^\eta[Z_t;E]$.
Intersecting $E$ with $\{X_t=l\}$ proves
\eqref{eq:path-tilt-identity} on $\mathcal F_t$.
The finite-path calculation uses the jump-time construction in
\cite[Section~2.6]{Norris}.\qedhere

\end{proof}

\begin{lemma}\label{lem:killed-straight-lower}
Fix $\lambda\in\R$ and use the killed kernel $K_t^y$ from
\eqref{eq:killing-definition}. There exists $y_0\ge1$ such that, for
every $A>0$, there are $c_A>0$ and $t_A<\infty$ satisfying
\begin{equation}
 \label{eq:killed-straight-lower}
 K_t^y(j,l)\ge c_A\frac{yz}{t^{3/2}},
 \qquad z=y+l-j+k'(\lambda)t,
\end{equation}
whenever $j,l\in\Z$, $t\ge t_A$ and
\[
 y_0\le y\le A\sqrt t,\qquad y_0\le z\le A\sqrt t.
\]
The constants are independent of $j,l\in\Z$. If
$\lambda$ ranges in a fixed compact interval, $y_0,c_A,t_A$ may be
chosen uniformly.
\end{lemma}

\begin{proof}

The forward chain has velocity $-k'(\lambda)$, and its reversal
has velocity $k'(\lambda)$. Apply Lemma~\ref{lem:entrance-height}
to $X$ with initial site and distance $(j,y)$ and to $\widehat X$
with initial site and distance $(l,z)$. Choose the common $y_0$
and $\delta>0$ from that lemma, and take $t$ above both entrance
time thresholds for the fixed $A$. Let $T_{\mathrm{stop}},\widehat{T}_{\mathrm{stop}}$ and
$E,\widehat E$ be the resulting stopping times and events.
Define finite measures
\[
\begin{gathered}
 F(ds,m)=\mathbb P_j^\lambda(E,T_{\mathrm{stop}}\in ds,X_{T_{\mathrm{stop}}}=m)\quad(m\in\Z),\\
 \widehat F(du,n)=
 \widehat{\mathbb P}_l^\lambda(\widehat E,\widehat{T}_{\mathrm{stop}}\in du,
                         \widehat X_{\widehat{T}_{\mathrm{stop}}}=n)\quad(n\in\Z).
\end{gathered}
\]
Their time supports lie in $[0,t/3]$, and their total masses satisfy
\[
 \sum_{m\in\Z}\int_{[0,t/3]} F(ds,m)\ge \kappa_A\frac{y}{\sqrt t},
 \qquad
 \sum_{n\in\Z}\int_{[0,t/3]} \widehat F(du,n)\ge
 \widehat\kappa_A\frac{z}{\sqrt t}.
\]
For $(s,m)$ in the support of $F$ and $(u,n)$ in the support of
$\widehat F$, put
\[
 y_s=y+m-j+k'(\lambda)s,
 \qquad
 z_u=z+n-l-k'(\lambda)u.
\]
Taking $B_A$ to be the larger of the two upper height constants in
Lemma~\ref{lem:entrance-height}, one has
\[
 \delta\sqrt t\le y_s,z_u\le B_A\sqrt t.
\]

By $E\in\mathcal F_{T_{\mathrm{stop}}}$, the strong Markov property permits
integration over the joint measure $F(ds,m)$ of the initial path
segment. For each resulting pair $(s,m)$, reverse the remaining
path from $m$ to $l$ by \eqref{eq:kern-path-reversal}; its weight
changes by $\pi_l/\pi_m$.
The event $\widehat E$ is measurable at $\widehat{T}_{\mathrm{stop}}$ for the
reversed chain, so its strong Markov property separates the entrance
segment from $l$ to $n$. The remaining reversed segment runs from $n$ to $m$; returning
it to the forward law contributes $\pi_m/\pi_n$.
The product of the two reversal factors is $\pi_l/\pi_n$.
On the support of $F(ds,m)\widehat F(du,n)$, we have
$s,u\le t/3$. Hence the initial and terminal segments occupy
disjoint time intervals, and the remaining duration $t-s-u$ is
at least $t/3$. Each entrance event and its stopping time are
determined by its own outer path segment, since
$E\in\mathcal F_{T_{\mathrm{stop}}}$ and $\widehat E\in\widehat{\mathcal F}_{\widehat{T}_{\mathrm{stop}}}$.
Inserting a middle path therefore does not change the forward stopping time
$T_{\mathrm{stop}}$ or the reversed stopping time $\widehat{T}_{\mathrm{stop}}$,
and each concatenated path is counted at most once.
Integrating the strong Markov identities in the
order just described, and applying Tonelli's theorem to their
nonnegative terms, gives
\begin{equation}
 \label{eq:technical-double-entrance}
 K_t^y(j,l)\ge
 \sum_{m,n\in\Z}\int_{[0,t/3]}\!\int_{[0,t/3]}
 F(ds,m)\widehat F(du,n)
 \frac{\pi_l}{\pi_n}
 K_{t-s-u}^{y_s}(m,n).
\end{equation}
The terminal distance in $K_{t-s-u}^{y_s}(m,n)$ is
\[
 y_s+n-m+k'(\lambda)(t-s-u)
   =z+n-l-k'(\lambda)u=z_u,
\]
so the constraints defining $F$ and $\widehat F$ agree with the constraint in
$K_{t-s-u}^{y_s}(m,n)$ at the two ends of the middle interval.
The remaining time $L=t-s-u$ lies in $[t/3,t]$, and hence
\[
 \delta\sqrt L\le y_s,z_u\le\sqrt3 B_A\sqrt L.
\]
By Lemma~\ref{lem:killed-interior-lower}, if
$t\ge3t_{\delta,\sqrt3 B_A}$, then
\[
 K_L^{y_s}(m,n)
 \ge\kappa_{\delta,\sqrt3 B_A}L^{-1/2}
 \ge\kappa_{\delta,\sqrt3 B_A}t^{-1/2}.
\]
Furthermore,
$\pi_l/\pi_n\ge\min_{0\le i<P}\pi_i/\max_{0\le i<P}\pi_i>0$. The measures $F$ and
$\widehat F$ have masses at least $\kappa_A y/\sqrt t$ and
$\widehat\kappa_A z/\sqrt t$, respectively. Hence
\eqref{eq:technical-double-entrance} gives
\[
 K_t^y(j,l)\ge
 \kappa_{\delta,\sqrt3 B_A}\,\frac{\min_{0\le i<P}\pi_i}{\max_{0\le i<P}\pi_i}
       \kappa_A\widehat\kappa_A\frac{yz}{t^{3/2}},
\]
which proves \eqref{eq:killed-straight-lower}.

For $\lambda$ in a fixed compact interval,
Lemma~\ref{lem:uniform-tilt} gives common choices of
$y_0,\delta,B_A$ and the time thresholds. It also gives positive
lower bounds for the three constants
$\kappa_{\delta,\sqrt3 B_A},\kappa_A,\widehat\kappa_A$ used in
\eqref{eq:technical-double-entrance}. Positivity and continuity of $\pi_i$ give a common
positive lower bound for $\min_{0\le i<P}\pi_i/\max_{0\le i<P}\pi_i$.
Consequently, the coefficient
\[
 \kappa_{\delta,\sqrt3 B_A}
 \frac{\min_{0\le i<P}\pi_i}{\max_{0\le i<P}\pi_i}
 \kappa_A\widehat\kappa_A
\]
in the lower bound is uniformly positive,
which proves the uniform assertion.\qedhere

\end{proof}

\subsection{A decreasing logarithmic boundary}

For the chain and velocity $v$ in Lemma~\ref{lem:periodic-survival},
we estimate the probability that $y+X_s-j-vs$ stays above
$-\kappa\log(1+s/(1+a))$, uniformly in the time offset $a\ge0$.
Here $y>0$ and $\kappa\ge0$. We use the corrected martingale
$M_s=X_s-j-vs+d_{[X_s]}-d_{[j]}$ for $s\ge0$, where $d$ is the
periodic extension of the normalized solution of \eqref{eq:corrector}.

\begin{lemma}\label{lem:descending-log-survival}
Use the positive periodic rates and velocity $v$ of
Lemma~\ref{lem:periodic-survival}. For every $\kappa\ge0$ there is
$C_\kappa<\infty$ such that, for all $a\ge0$, $y>0$, $j\in\Z$ and $t\ge0$,
\begin{equation}
 \label{eq:descending-log-survival}
 \mathbb P_j\left\{y+X_s-j-vs>
       -\kappa\log\left(1+\frac{s}{1+a}\right)
                    \text{ for all }0\le s\le t\right\}
 \le C_\kappa\frac{1+y}{\sqrt{1+t}}.
\end{equation}
More explicitly, for fixed $P$, $\kappa$ and
$0<r_-\le r_+<\infty$, one constant $C_\kappa$ works for every choice of
$P$-periodic rates with $r_-\le A_i,B_i\le r_+$ for $0\le i<P$.
\end{lemma}

\begin{proof}

Fix $j\in\Z$. Throughout this proof,
$\mathbb P=\mathbb P_j$ and $\mathbb E=\mathbb E_j$ refer to
the chain with generator $\mathcal Q$ and rates $A_i,B_i$
in Lemma~\ref{lem:periodic-survival}, started at $X_0=j$.
Let $d$ solve \eqref{eq:corrector} and set
$M_s=X_s-j-vs+d_{[X_s]}-d_{[j]}$. Choose
$B\ge1+2\max_{0\le i<P}|d_i|$; then every jump of $M$
has absolute value at most $B$. Put
$V_-=\min_{0\le i<P}V_i>0$, with $V_i$ defined in
\eqref{eq:quadratic-rate}.
For $0<x<R$, define
$\sigma_R=\inf\{s\ge0:x+M_s\notin(0,R)\}$.
For every $s\ge0$, $x+M_{s\wedge\sigma_R}\in[-B,R+B]$.
At a point $(z,i)\in(0,R)\times\{0,\ldots,P-1\}$,
the generator of $(x+M_s,[X_s])$, including the continuous drift
of $M_s$, sends the function $(z,i)\mapsto(z+B)(R+B-z)$
to $-V_i$. Bounded stopping,
followed by monotone convergence for the exit time and dominated
convergence for the stopped martingale, gives
\[
 \mathbb E\sigma_R\le\frac{(x+B)(R+B)}{V_-},\qquad
 \mathbb E(x+M_{\sigma_R})=x,\qquad
 \mathbb P(x+M_{\sigma_R}\ge R)\le\frac{x+B}{R+B}\quad\text{for }0<x<R.
\]
For $t>0$, a path staying positive up to $t$ either has not left
$(0,R)$ by $t$ or leaves through its upper endpoint. Thus
\[
 \mathbb P\{x+M_s>0\text{ for }0\le s\le t\}
 \le\frac{x+B}{R+B}+\frac{(x+B)(R+B)}{V_-t}\quad\text{for }0<x<R.
\]
Taking $R=\sqrt t$ when $t\ge1$ and $x<\sqrt t$, and using the
bound one in the other cases, gives
\[
 \mathbb P\{x+M_s>0\text{ for }0\le s\le t\}
       \le C(1+x)(1+t)^{-1/2},\quad\text{ for }x>0,\ t\ge0,
\]
uniformly for $j\in\Z$.
Fix $a\ge0$ and $y>0$. Put $D_0=\|d\|_\infty$, $x=y+2D_0$ and
$g(s)=\kappa\log(1+s/(1+a))$ for $s\ge0$.

Since $|d_{[X_s]}-d_{[j]}|\le2D_0$, the process
$Z_s=x+M_s+g(s)$ is at least $y+X_s-j-vs+g(s)$.
Consequently,
\[
 \{y+X_s-j-vs+g(s)>0\text{ for every }s\in[0,t]\}
 \subseteq\{\tau>t\},\qquad
 \tau=\inf\{s\ge0:Z_s\le0\}.
\]
Write $F(t)=\mathbb P(\tau>t)$ for $t\ge0$.
On $\{\tau<\infty\}$, $-B\le Z_\tau\le0$. Fix $t\ge0$ and put
$T_{\mathrm{stop}}=t\wedge\tau$. To justify optional stopping at $T_{\mathrm{stop}}$,
let $N_t$ be the number of jumps of $X$ on $[0,t]$. Then
\[
 \sup_{0\le s\le t}|M_s|\le BN_t+|v|t,\qquad
 \mathbb E\sup_{0\le s\le t}|M_s|^2<\infty\quad\text{for }t\ge0,
\]
because $N_t$ is dominated by a Poisson random variable with finite
mean. For $n\in\N$, the finite-grid stopping times
$T_{\mathrm{stop},n}=(2^{-n}\lceil2^nT_{\mathrm{stop}}\rceil)\wedge t$ take values in $[0,t]$, and optional stopping gives
$\mathbb EM_{T_{\mathrm{stop},n}}=0$. Right-continuity and the displayed integrable
bound imply $\mathbb EM_{T_{\mathrm{stop}}}=0$ as $n\to\infty$; see
\cite[Chapter~2]{EthierKurtz} and
\cite[Section~4.8]{Durrett2019}.
Moreover, $g(T_{\mathrm{stop}})=\int_0^t g'(r)\mathbf1_{\{r<T_{\mathrm{stop}}\}}\,\dd r$,
so Tonelli's theorem and $g'\ge0$ give
$\mathbb Eg(T_{\mathrm{stop}})=\int_0^t g'(r)F(r)\,\dd r$.
Since $Z_{T_{\mathrm{stop}}}=x+M_{T_{\mathrm{stop}}}+g(T_{\mathrm{stop}})$, these identities yield
\[
 \mathbb E[Z_t;\tau>t]+\mathbb E[Z_\tau;\tau\le t]
  =x+\int_0^t g'(s)F(s)\,\dd s\quad\text{for }t\ge0.
\]
Using $Z_\tau\ge-B$ on $\{\tau\le t\}$, we obtain
\begin{equation}
 \label{eq:stopped-first-moment}
 \mathbb E[Z_t;\tau>t]
   \le x+B+\int_0^t g'(s)F(s)\,\dd s\quad\text{for }t\ge0.
\end{equation}

On $\{\tau>t/2\}$, condition on $Z_{t/2}=z>0$ and the
site $X_{t/2}$. For $0\le r\le t/2$,
$0\le g(t/2+r)-g(t/2)\le\kappa\log2$, and
\[
 M_{t/2+r}-M_{t/2}
 =X_{t/2+r}-X_{t/2}-vr
      +d_{[X_{t/2+r}]}-d_{[X_{t/2}]}.
\]
The Markov property identifies these increments with the corrected
martingale for the chain started at $X_{t/2}$. Survival up to $t$
requires
$z+\kappa\log2+M_{t/2+r}-M_{t/2}>0$ for $0\le r\le t/2$.
Applying the bound for
$\mathbb P\{x+M_s>0\text{ for }0\le s\le t\}$ proved at the start of this proof,
with initial height $z+\kappa\log2$, gives the conditional upper bound
$C(1+z+\kappa\log2)(1+t/2)^{-1/2}$. Integration over
$\{\tau>t/2\}$ therefore yields
\[
 \begin{aligned}
 F(t)
 &\le\frac{C_\kappa}{\sqrt{1+t}}
       \bigl(F(t/2)+\mathbb E[Z_{t/2};\tau>t/2]\bigr)\\
 &\le\frac{C_\kappa}{\sqrt{1+t}}
      \left(1+x+B+\int_0^t g'(s)F(s)\,\dd s\right).
 \end{aligned}
\]
The last step is \eqref{eq:stopped-first-moment} evaluated at
$t/2$, followed by $F(t/2)\le1$ and $g'F\ge0$.
Increasing $C_\kappa$ includes bounded $t$ in this estimate.
Put $G(t)=1+x+B+\int_0^t g'(s)F(s)\,\dd s$. By the bound just
proved, $F(t)\le C_\kappa G(t)/\sqrt{1+t}$. Since
$g'(t)=\kappa/(1+a+t)$, we obtain, for almost every $t\ge0$, 
\[
 G'(t)\le\frac{C_\kappa g'(t)}{\sqrt{1+t}}G(t),\qquad
 \int_0^\infty\frac{g'(t)}{\sqrt{1+t}}\,\dd t
   \le\kappa\int_0^\infty(1+t)^{-3/2}\,\dd t=2\kappa.
\]
Applying Gronwall's inequality with the integrable coefficient
$C_\kappa g'(t)/\sqrt{1+t}$ gives $G(t)\le C(1+x)$, with $C$
independent of $a,t$. Since $x=y+2D_0$ and
$\mathbb P_j\{y+X_s-j-vs+g(s)>0\text{ for }0\le s\le t\}
\le F(t)$, this bound on $G$ gives
\eqref{eq:descending-log-survival}.
On the compact rate set $r_-\le A_i,B_i\le r_+$, the phase
stationary probability vector $\varpi$ depends continuously on the
rates. The normalized system \eqref{eq:corrector} has a unique
solution, so its solution $d$ also depends continuously on the
rates. Consequently the quantities $V_i$ in
\eqref{eq:quadratic-rate} are continuous, and the increments
$-1+d_{i-1}-d_i$ and $1+d_{i+1}-d_i$ have a common bound.
Lemma~\ref{lem:periodic-survival} gives $V_i>0$ at every point of this
set. Compactness and the finite phase set therefore give a common
positive lower bound for $V_i$, so the exit and survival estimates
have common constants on the stated rate set.\qedhere

\end{proof}

\paragraph{The backward logarithmic boundary.}

For $\kappa\ge0$, $a\ge0$, $t\ge0$, $\lambda\in\R$,
$j,l\in\Z$ and $y>0$, define $b_{t,a}(s)$ for $0\le s\le t$
and the kernel $\mathcal K_{t,a}^{y,\lambda}(j,l)$ by
\begin{equation}
 \label{eq:backward-boundary-definition}
 \begin{aligned}
 b_{t,a}(s)&=\kappa\log\frac{1+a+t}{1+a+t-s},\\
 \mathcal K_{t,a}^{y,\lambda}(j,l)
  &=\mathbb P_j^\lambda\{X_t=l,\ y+X_s-j+k'(\lambda)s>b_{t,a}(s)
                       \text{ for }0\le s\le t\}.
 \end{aligned}
\end{equation}
Here $a$ denotes the time offset. The terminal height for a path from $j$ to $l$ is
$z=y+l-j+k'(\lambda)t-b_{t,a}(t)$. We set
$\mathcal K_{t,a}^{y,\lambda}(j,l)=0$ when $y\le0$ or $z\le0$.
At $t=0$, we set $\mathcal K_{0,a}^{y,\lambda}(j,l)=\mathbf1_{\{j=l\}}$
for $y>0$.

For $0\le s\le t$, the Markov property gives the composition rule
\begin{equation}
 \label{eq:offset-kernel-composition}
 \mathcal K_{t,a}^{y,\lambda}(j,l)
 ={\sum_{\substack{n\in\Z\\h_n>0}}
   }\mathcal K_{s,a+t-s}^{y,\lambda}(j,n)
   \mathcal K_{t-s,a}^{h_n,\lambda}(n,l),\qquad
 h_n=y+n-j+k'(\lambda)s-b_{t,a}(s).
\end{equation}
Indeed, $b_{s,a+t-s}(r)=b_{t,a}(r)$ for $0\le r\le s$, and
$b_{t,a}(s+r)-b_{t,a}(s)=b_{t-s,a}(r)$ for $0\le r\le t-s$.
These identities give the terminal height
$h_n+l-n+k'(\lambda)(t-s)-b_{t-s,a}(t-s)=z$
for $\mathcal K_{t-s,a}^{h_n,\lambda}(n,l)$.
The survival events defining the two kernel factors in
\eqref{eq:offset-kernel-composition} therefore both require $h_n>0$
at the joining time $s$.

\begin{lemma}\label{lem:backward-log-kernel}

Fix $\lambda\in\R$ and $\kappa\ge0$. There is $C<\infty$ such that
\begin{equation}
 \label{eq:backward-log-upper}
 \mathcal K_{t,a}^{y,\lambda}(j,l)
       \le C\frac{(1+y)(1+z)}{(1+t)^{3/2}}
\end{equation}
for all $a,t\ge0$, $j,l\in\Z$ and $y>0$ for which
$z=y+l-j+k'(\lambda)t-b_{t,a}(t)>0$.
For every $A>0$ there are $Y\ge1$, $T<\infty$ and $c_A>0$ such that
\begin{equation}
 \label{eq:backward-log-lower}
 \mathcal K_{t,a}^{y,\lambda}(j,l)\ge c_A\frac{yz}{t^{3/2}}
\end{equation}
when $a\ge0$, $j,l\in\Z$, $t\ge T$ and $Y\le y,z\le A\sqrt t$.
The constants are independent of $a,j,l$ and can be chosen uniformly
for $\lambda$ in a compact interval.

\end{lemma}
\begin{proof}

Write $c=k'(\lambda)$, $Q=Q^\lambda$ and
$\pi_i=\pi_i^\lambda$ for the periodically extended stationary
phase weights. Let $\widehat X$ have the reversed generator
\eqref{eq:dual-rates}, and write $\widehat{\mathbb P}_l^\lambda$ for
its law starting at $l$. The finite-path reversal calculation gives
\begin{equation}
 \label{eq:backward-kernel-reversal}
 \begin{aligned}
 \pi_j\mathcal K_{t,a}^{y,\lambda}(j,l)
 &=\pi_l\widehat{\mathbb P}_l^\lambda
 \biggl\{\widehat X_t=j,\\
 &\hspace{12mm}z+\widehat X_r-l-cr>
       -\kappa\log\left(1+\frac r{1+a}\right)
       \text{ for }0\le r\le t\biggr\}.
 \end{aligned}
\end{equation}
To verify \eqref{eq:backward-kernel-reversal}, take a path with sites
$x_0=j,\ldots,x_n=l$ and jump times
$0<s_1<\cdots<s_n<t$. Its jump factors satisfy
\[
 \pi_j\prod_{h=1}^n Q(x_{h-1},x_h)
   =\pi_l\prod_{h=1}^n\widehat Q(x_h,x_{h-1}).
\]
The holding factors also agree, since
$\widehat\alpha_i+\widehat\beta_i=\alpha_i+\beta_i$.
On each open holding interval, with reversed time $r=t-s$,
\[
 \begin{aligned}
 y+X_{t-r}-j+c(t-r)-b_{t,a}(t-r)
 &=z+X_{t-r}-l-cr\\
 &\quad+\kappa\log\left(1+\frac r{1+a}\right).
 \end{aligned}
\]
Thus the two positivity constraints agree between every pair of
jumps. For each fixed site, contacts solve
$X_s=j-y-cs+b_{t,a}(s)$ with $X_s$ fixed at that site. If
$\kappa>0$, the right-hand side is strictly convex and meets that
site at most twice; if $\kappa=0$ and $c\ne0$, it meets it at most
once. These contact times are deterministic, so a jump occurs at
one of them with probability zero. When $\kappa=c=0$, reversal
preserves avoidance of the fixed forbidden set of sites. The constraints therefore agree up to a null
set on the full closed time interval, proving
\eqref{eq:backward-kernel-reversal} for the required all-time event.

The reversal formula starts the reversed chain at height $z$
above the decreasing boundary $-\kappa\log(1+r/(1+a))$.
To prove \eqref{eq:backward-log-upper} for $t\ge3$, split $[0,t]$ into
$[0,t/3]$, $[t/3,2t/3]$ and $[2t/3,t]$, and discard absorption only
on $[t/3,2t/3]$. On $[0,t/3]$ the boundary $b_{t,a}$ is nonnegative, so the sum
over its surviving endpoints is at most
$C(1+y)(1+t)^{-1/2}$ by Lemma~\ref{lem:periodic-survival}.
By \eqref{eq:kernel-bound}, the unrestricted transition probability
on the middle interval is at most $C(1+t)^{-1/2}$. To sum the last part over its initial point, reverse
its paths with the periodic invariant weights. Relative to the final
height $z$, the reversed boundary at elapsed time $r$ is exactly
\[
 b_{t,a}(t-r)-b_{t,a}(t)
       =-\kappa\log\left(1+\frac r{1+a}\right)\quad\text{for }0\le r\le t.
\]
Put $h_n=y+n-j+2ct/3-b_{t,a}(2t/3)$ for $n\in\Z$.
By \eqref{eq:backward-kernel-reversal}, the velocity $c$ of
$\widehat X$ and Lemma~\ref{lem:descending-log-survival} give
\[
 \sum_{\substack{n\in\Z\\h_n>0}}
 \mathcal K_{t/3,a}^{h_n,\lambda}(n,l)
 \le\frac{\max_{0\le i<P}\pi_i}{\min_{0\le i<P}\pi_i}
 C(1+z)(1+t)^{-1/2}.
\] By \eqref{eq:offset-kernel-composition}, multiplying the first, middle
and reversed terminal bounds proves \eqref{eq:backward-log-upper} for $t\ge3$; bounded times follow by
enlarging the constant.

To prove \eqref{eq:backward-log-lower}, set
$d_{t,a}=\kappa t^{-1}\log(1+t/(1+a))$.
Since $0\le d_{t,a}\le\kappa\log(1+t)/t$, for all sufficiently
large $t$ the value $c-d_{t,a}$ lies in the image under $k'$ of a
fixed neighborhood of $\lambda$. Choose $\eta=\eta_{t,a}$ there
with $k'(\eta)=c-d_{t,a}$. The positive lower bound for $k''$ on
that neighborhood, supplied by Lemma~\ref{lem:matrix}, gives
$|\eta-\lambda|\le C\log(1+t)/t$, uniformly in $a$.
Since $b_{t,a}$ is convex, $b_{t,a}(0)=0$ and
$b_{t,a}(t)=d_{t,a}t$, we have
\[
 b_{t,a}(s)-d_{t,a}s\le0\quad\text{ for }0\le s\le t.
\]
On the common path space used in Lemma~\ref{lem:uniform-tilt}, put
$H_s^\eta=y+X_s-j+k'(\eta)s$ for $0\le s\le t$. Since
\[
 y+X_s-j+cs-b_{t,a}(s)
  =H_s^\eta+d_{t,a}s-b_{t,a}(s)\ge H_s^\eta\quad\text{for }0\le s\le t,
\]
the event $H_s^\eta>0$ for all $s\in[0,t]$ is contained in the
survival event in \eqref{eq:backward-boundary-definition}. Moreover,
$H_0^\eta=y$ and, on $X_t=l$, $H_t^\eta=z$, because
$d_{t,a}t=b_{t,a}(t)$. By Lemmas~\ref{lem:killed-straight-lower} and~\ref{lem:uniform-tilt},
the event $\{X_t=l,\ H_s^\eta>0\text{ for }0\le s\le t\}$ has
probability at least $c_Ayz/t^{3/2}$ under $\mathbb P_j^\eta$.
The constant is uniform for $\eta$ in a fixed compact neighborhood
of $\lambda$.

Apply \eqref{eq:path-tilt-identity} with parameters $\lambda,\eta$
to $\{X_t=l,\ H_s^\eta>0\text{ for all }0\le s\le t\}$,
and put $\delta=\eta-\lambda$. Since
$l-j=z-y-k'(\eta)t$, the logarithm of the exponential factor is
\[
 [k(\eta)-k(\lambda)-k'(\eta)\delta]t+\delta(z-y).
\]
Taylor's theorem bounds its absolute value by
$Ct\delta^2+2A\sqrt t\,|\delta|$. The bound
$|\delta|\le C\log(1+t)/t$ makes both terms tend to zero as
$t\to\infty$, uniformly in $a$.
The ratio $\phi_j(\eta)\phi_l(\lambda)/(\phi_j(\lambda)\phi_l(\eta))$
in \eqref{eq:path-tilt-identity} has positive lower and finite upper bounds,
uniformly for $j,l\in\Z$ on the fixed compact parameter interval containing
$\lambda$ and $\eta$. Combining these bounds with
$c_Ayz/t^{3/2}$ under $\mathbb P_j^\eta$ proves
\eqref{eq:backward-log-lower}. If $\lambda$ ranges over a fixed compact
interval, choose a common compact neighborhood containing $\eta$.
Lemmas~\ref{lem:matrix} and~\ref{lem:uniform-tilt} give a common
positive lower bound for $k''$, common positive lower and finite upper bounds for
$\phi_i(\eta)$ and $\phi_i(\lambda)$, and common constants in \eqref{eq:killed-straight-lower} there. The choices of the time threshold
and the lower-bound constant are therefore uniform on that interval.\qedhere

\end{proof}

\begin{lemma}\label{lem:bridge-occupation}
Fix $\lambda>0$ and $\kappa\ge0$. For $t>0$, $j,l\in\Z$ and $y>0$, put
\[
 \begin{aligned}
 Y_s&=y+X_s-j+k'(\lambda)s-b_{t,0}(s)
       &&\text{for }0\le s\le t,\\
 z&=y+l-j+k'(\lambda)t-\kappa\log(1+t).
 \end{aligned}
\]
Suppose that both endpoint heights $y,z$ belong to a fixed compact
interval $[Y,B]\subset(0,\infty)$, where $Y$ is at least the threshold in
\eqref{eq:backward-log-lower} with $A=1$. There are a constant $C_B$ and
a time threshold, both depending only on $\lambda,\kappa,Y,B$ and the fixed
model parameters, such that, for all $j,l\in\Z$, $y,z\in[Y,B]$ and $t$
above that threshold,
\begin{equation}
 \label{eq:bridge-occupation}
 \mathbb E_j^\lambda\left[
   \int_0^t(1+Y_s)e^{-\lambda Y_s}\,\dd s
   \,\middle|\,X_t=l,\ Y_r>0\text{ for }0\le r\le t\right]\le C_B.
\end{equation}
By \eqref{eq:backward-log-lower}, the conditioning event has positive
probability for all sufficiently large $t$.
\end{lemma}

\begin{proof}
For $s\in[0,t]$ and $n\in\Z$, put
$h=h(s,n)=y+n-j+k'(\lambda)s-b_{t,0}(s)$. By \eqref{eq:offset-kernel-composition},
the boundary on $[0,s]$ is $b_{s,t-s}$, and
$b_{t,0}(s+r)-b_{t,0}(s)=b_{t-s,0}(r)$ for $0\le r\le t-s$.
Conditioning at time $s$ and using Tonelli's theorem therefore
express the numerator of the conditional expectation in
\eqref{eq:bridge-occupation} as
\begin{equation}
 \label{eq:bridge-convolution}
 {\int_0^t\sum_{\substack{n\in\Z\\h>0}}(1+h)e^{-\lambda h}}
       \mathcal K_{s,t-s}^{y,\lambda}(j,n)
       \mathcal K_{t-s,0}^{h,\lambda}(n,l)\,\dd s.
\end{equation}
In $\mathcal K_{t-s,0}^{h,\lambda}(n,l)$ the terminal height is $z$, because
$b_{t,0}(s)+b_{t-s,0}(t-s)=b_{t,0}(t)$.
The upper bound \eqref{eq:backward-log-upper} shows that
\eqref{eq:bridge-convolution} is at most
\[
 C(1+y)(1+z)\int_0^t
 \frac{{\sum_{\substack{n\in\Z\\h>0}}(1+h)^3e^{-\lambda h}}}
 {(1+s)^{3/2}(1+t-s)^{3/2}}\,\dd s
 \le C_B(1+t)^{-3/2}.
\]
For each fixed $s\in[0,t]$,
$\{h(s,n):n\in\Z,\ h(s,n)>0\}=\rho+\N_0$
for some $\rho=\rho(s)\in(0,1]$. Hence
\[
 {\sum_{\substack{n\in\Z\\h>0}}(1+h)^3e^{-\lambda h}
 \le}\sum_{m=0}^\infty(m+2)^3e^{-\lambda m}<\infty.
\]
The series converges because $\lambda>0$, and the bound is uniform
in $\rho$. For the time integral, split at $t/2$. On $[0,t/2]$,
use $(1+t-s)^{-3/2}\le 2^{3/2}(1+t)^{-3/2}$ and integrate
$(1+s)^{-3/2}$; on $[t/2,t]$, interchange $s$ and $t-s$.
For $t$ above both $B^2$ and the time threshold in
\eqref{eq:backward-log-lower}, the endpoint heights satisfy
$Y\le y,z\le\sqrt t$. The conditioning probability is therefore
at least $c_1yz/t^{3/2}$, by \eqref{eq:backward-log-lower} with $A=1$.
Dividing the bound for \eqref{eq:bridge-convolution} by the lower
bound $c_1yz/t^{3/2}$ gives at most
$C_Bt^{3/2}(1+t)^{-3/2}/(yz)$. Since $yz\ge Y^2$ and
$t^{3/2}(1+t)^{-3/2}\le1$, the resulting constant is uniform over
$y,z\in[Y,B]$. This proves \eqref{eq:bridge-occupation}. General formulas for additive
functionals are discussed in \cite{FitzsimmonsPitman}.\qedhere
\end{proof}

\subsection{The linearized equation}

Combining the kernel estimate \eqref{eq:kernel-bound} with the KPP
bound in \eqref{eq:hyp} gives a pointwise upper bound for solutions
of \eqref{eq:lattice}.
Boutillon~\cite{Boutillon2024} uses estimates for the linearized
equation to obtain upper bounds for a nonlocal Fisher--KPP equation.

\begin{lemma}\label{lem:linear-bound}
Let $w$ be a solution of \eqref{eq:lattice} with values in $[0,1]$.
If $\lambda>0$ and $\sum_{j\in\Z} e^{\lambda j}w_j(0)<\infty$, then
\begin{equation}
\label{eq:linear-bound}
 w_j(t)\le C_\lambda(1+t)^{-1/2}
 e^{-\lambda j+k(\lambda)t}
 \sum_{l\in\Z} e^{\lambda l}w_l(0).
\end{equation}
The estimate holds for $t\ge0$ and $j\in\Z$. The constant
$C_\lambda$ is independent of $j$, $t$ and the initial datum;
it may depend on $\lambda$, $D$, $q$, $P$, and $(a_i)_{i=0}^{P-1}$.
\end{lemma}
\begin{proof}

Fix $\lambda>0$ as in the statement and set
$H_j=e^{-\lambda j}\phi_j(\lambda)>0$ for $j\in\Z$.
For $v\in\ell^\infty(\Z)$, define $\mathcal L_{\Z}v$ by
\[
 (\mathcal L_{\Z}v)_j
 =D(v_{j-1}+qv_{j+1}-(q+1)v_j)+a_jv_j
 \quad\text{ for }j\in\Z.
\]
Let $Q$ be the generator with rates $DH_{j-1}/H_j$ from $j$ to
$j-1$ and $DqH_{j+1}/H_j$ from $j$ to $j+1$; its kernel is
$p_t(j,l)$ in \eqref{eq:poisson-kernel}.
Periodicity of $a_j$ and $\phi_j(\lambda)$ makes
$\mathcal L_{\Z}$ and $Q$ bounded on $\ell^\infty(\Z)$. We first establish the kernel identity using finitely
supported sequences. For every such $z:\Z\to\C$, \eqref{eq:matrix}
and the definition of $Q$ give
$(\mathcal L_{\Z}(Hz))_j=H_j(Qz+k(\lambda)z)_j$.
{Set $C=D(q+1)$, and let $I$ denote the identity on $\ell^\infty(\Z)$.
The diagonal entries of both
$\mathcal L_{\Z}+CI$ and $Q+(k(\lambda)+C)I$ are $a_j>0$,
and their off-diagonal entries are nonnegative.}
For the operators below, write $\mathcal A(j,l)=(\mathcal A\delta_l)_j$,
where $\delta_l(j)=\mathbf1_{\{j=l\}}$.
For every integer $n\ge0$, multiplication along the finitely many
length-$n$ nearest-neighbor or holding paths gives
\[
 (\mathcal L_{\Z}+CI)^n(j,l)
 =\frac{H_j}{H_l}
       \bigl(Q+(k(\lambda)+C)I\bigr)^n(j,l),
 \quad\text{for }n\in\N_0,\ j,l\in\Z,
\]
since the intermediate factors $H$ cancel. Summing the
operator-norm-convergent exponential series therefore yields
\[
 (e^{t\mathcal L_{\Z}})(j,l)
   =e^{k(\lambda)t}\frac{H_j}{H_l}p_t(j,l)
 \quad\text{for }t\ge0,\ j,l\in\Z.
\]

The exponential series for $\mathcal L_{\Z}+CI$ has
nonnegative entries. Define
$U_j(t)=(e^{t\mathcal L_{\Z}}w(0,\cdot))_j$ for
$t\ge0$ and $j\in\Z$, where $w(0,\cdot)=(w_l(0))_{l\in\Z}$.
Then $U_j(t)\ge0$ and
$U\in C^1([0,\infty);\ell^\infty(\Z))$.
The KPP bound in \eqref{eq:hyp}, together with the fixed extensions
of the reactions outside $[0,1]$, gives
$\widetilde f_j(s)\le a_js$ for all $s\ge0$. Since $\partial_tU_j=(\mathcal L_{\Z}U)_j$,
we obtain
\[
 \partial_tU_j(t)\ge
 D(U_{j-1}(t)+qU_{j+1}(t)-(q+1)U_j(t))
 +\widetilde f_j(U_j(t)),\qquad U_j(0)=w_j(0)\quad\text{for }t\ge0,\ j\in\Z.
\]
Applied to $w(0,\cdot)\ge0$, the nonnegative exponential series
for $\mathcal L_{\Z}+CI$ permits summation over its matrix entries
by Tonelli's theorem. Applying Lemma~\ref{lem:comparison} to $w$
and $U$ therefore gives
\[
 w_j(t)\le\sum_{l\in\Z}(e^{t\mathcal L_{\Z}})(j,l)w_l(0)
 =e^{-\lambda j+k(\lambda)t}\phi_j(\lambda)
   \sum_{l\in\Z} p_t(j,l)\frac{e^{\lambda l}}{\phi_l(\lambda)}w_l(0)\quad\text{for }t\ge0,\ j\in\Z.
\]
The series over $l$ is absolutely convergent because
$\sum_{j\in\Z} e^{\lambda j}w_j(0)<\infty$ and
$\min_{0\le i<P}\phi_i(\lambda)>0$. By \eqref{eq:kernel-bound} and the bound
$\phi_j(\lambda)/\phi_l(\lambda)\le
\max_{0\le i<P}\phi_i(\lambda)/\min_{0\le i<P}\phi_i(\lambda)$, we obtain
\eqref{eq:linear-bound}.\qedhere

\end{proof}

\section{Subsolutions and moving boundaries}\label{sec:moving-boundaries}

Throughout this section, assume \eqref{eq:hyp}.
Lemma~\ref{lem:bump} constructs subsolutions $v(t,j)$ of
\eqref{eq:lattice} supported in $\{j\in\Z:|j-ct|<L_0\}$ for
any prescribed $c\in(c_-,c_+)$. Lemma~\ref{lem:local-general}
gives a criterion, in terms of principal eigenvectors on finite
subsets of $\Tq$, for $u(t,x)\to1$ as $t\to\infty$ at every
fixed $x\in\Tq$.

For each $b\in\Z$, let $w$ solve \eqref{eq:lattice} with
$w_j(0)=\mathbf1_{\{j\le b\}}$ for $j\in\Z$.
Lemmas~\ref{lem:sharp-step-upper} and~\ref{lem:sharp-step-lower}
locate $H_\theta(t)=\max\{j\in\Z:w_j(t)\ge\theta\}$ at
$c_+t-(3/(2\lambda_+))\log(1+t)+O(1)$ as $t\to\infty$,
for each $\theta\in(0,1)$.
Their proofs use the kernel bounds \eqref{eq:backward-log-upper}
and \eqref{eq:backward-log-lower}, the conditional expectation bound
\eqref{eq:bridge-occupation}, and the comparison formulas in
Lemma~\ref{lem:moving-boundary-representation}.
When $\mu>0$, Lemma~\ref{lem:periodic-halfline-log} gives analogous
bounds for solutions $U=(U_r(t))_{r\in\N_0}$ with finite nonzero
initial support and $U_0(t)\to1$ as $t\to\infty$, for the interior
equation in \eqref{eq:radial} and boundary equation
$\dot U_0=b_0(U_1-U_0)+f_0(U_0)$ with fixed $b_0>0$.
Throughout this section, define $\mathcal A$ and $\mathcal L_{\Z}$
on real or complex sequences $v=(v_j)_{j\in\Z}$ by the finite
coordinate formulas
\[
 (\mathcal Av)_j=D(v_{j-1}+qv_{j+1}-(q+1)v_j),\qquad
 (\mathcal L_{\Z}v)_j=(\mathcal Av)_j+a_jv_j
 \quad\text{for }j\in\Z.
\]

\begin{lemma}\label{lem:bump}

For every $c\in(c_-,c_+)$ and $R>0$, there are $L_0>R$,
$\delta>0$ and a subsolution $v:\R\times\Z\to[0,1]$ of
\eqref{eq:lattice} such that, for all $t\in\R$ and $j\in\Z$,
\[
 v(t,j)=0\ \text{ for }|j-ct|\ge L_0,\qquad
 v(t,j)\ge\delta\ \text{ for }|j-ct|\le R.
\]
The map $t\mapsto v(t,\cdot)$ is locally Lipschitz into
$\ell^\infty(\Z)$, and
$\partial_tv(t,j)\le(\mathcal Av(t,\cdot))_j+f_j(v(t,j))$
for every $j\in\Z$ and almost every $t\in\R$.
Every fixed multiple $\vartheta v$, $0<\vartheta\le1$, is also
a subsolution. When $c=0$, each $\vartheta v$ is independent of $t$
and satisfies
$0<(\mathcal A(\vartheta v))_j+f_j(\vartheta v(t,j))$
wherever $v(t,j)>0$. The same $L_0$ and $\delta$ may be used when
$f_j$ is replaced by $f_{j+m}$ for any $m\in\Z$.

\end{lemma}
\begin{proof}

Fix $c\in(c_-,c_+)$ and $R>0$. By \eqref{eq:interior-speed},
there is a unique $\lambda_c\in(\lambda_-,\lambda_+)$ with
$k'(\lambda_c)=c$ and
$\sigma_0:=k(\lambda_c)-c\lambda_c>0$.
For complex $z$ near $\lambda_c$, define $M_z$ by
\eqref{eq:matrix}, with $\lambda$ replaced by $z$ and
$v\in\C^P$. Each entry of $M_z$ is holomorphic in $z$.
Since $k(\lambda_c)$ is a simple eigenvalue of $M_{\lambda_c}$,
analytic perturbation theory gives holomorphic continuations of the
eigenvalue $k$ and eigenvector $\phi$, satisfying
$M_z\phi(z)=k(z)\phi(z)$ and
$\sum_{i=0}^{P-1}\phi_i(z)=1$; see \cite{Kato}.
Extend $\phi(z)$ periodically by
$\phi_{j+P}(z)=\phi_j(z)$ for $j\in\Z$.
For real $\lambda$ near $\lambda_c$ and small real $\eta$, define
\[
 F(\lambda,\eta)=
 \begin{cases}
  \displaystyle\frac{\im k(\lambda+\mathrm i\eta)}{\eta},&\eta\ne0,\\[2mm]
  k'(\lambda),&\eta=0.
 \end{cases}
\]
By holomorphy and the reality of $k$ on the real axis,
\[
 F(\lambda,\eta)=\int_0^1
       \re k'(\lambda+\mathrm i\eta r)\,\dd r.
\]
Thus $F$ is smooth near $(\lambda_c,0)$, and
\eqref{eq:interior-speed} and \eqref{eq:convex-symmetry} give
\[
 F(\lambda_c,0)=c,
 \qquad
 \partial_\lambda F(\lambda_c,0)=k''(\lambda_c)>0.
\]
The implicit function theorem therefore gives a real function
$\lambda(\eta)$ with $\lambda(\eta)\to\lambda_c$ as $\eta\downarrow0$
and $F(\lambda(\eta),\eta)=c$ for all sufficiently small $\eta>0$.
Put
\[
 z=\lambda(\eta)+\mathrm i\eta,
 \qquad
 \sigma=k(z)-cz.
\]
Since $F(\lambda(\eta),\eta)=c$, $\sigma$ is real. Moreover,
$\sigma\to\sigma_0>0$ as $\eta\downarrow0$, so $\sigma>0$ for all
sufficiently small $\eta>0$.

For $(t,j)\in\R\times\Z$, set
$Z(t,j)=e^{-z(j-ct)}\phi_j(z)$.
The equation $M_z\phi(z)=k(z)\phi(z)$ and the definition of
$\mathcal L_{\Z}$ give
$(\mathcal L_{\Z}Z(t,\cdot))_j=k(z)Z(t,j)$, while
$\partial_tZ(t,j)=czZ(t,j)$.
Since $\sigma=k(z)-cz$, we obtain
\[
 Z_t=\mathcal L_{\Z}Z-\sigma Z
 \qquad \text{for }t\in\R,\ j\in\Z.
\]
Define $V(t,j)=\re Z(t,j)$ for $(t,j)\in\R\times\Z$.
Taking real parts gives
$V_t=\mathcal L_{\Z}V-\sigma V$, since $\sigma$ is real.
For each $i\in\{0,\ldots,P-1\}$, holomorphy of $\phi_i$ gives
$\phi_i(z)\to\phi_i(\lambda_c)>0$ as $\eta\downarrow0$.
Use the principal argument $\arg\zeta\in(-\pi,\pi]$ for
$\zeta\in\C\setminus\{0\}$. Since $\phi_i(\lambda_c)>0$
for $0\le i<P$, we may decrease $\eta$ so that
\[
 |\arg\phi_i(z)|<\frac{\pi}{8}
 \quad\text{ for } 0\le i<P.
\]
Writing $s=j-ct$, we have
\[
 V(t,j)=e^{-\lambda(\eta)s}|\phi_j(z)|
        \cos\bigl(\arg\phi_j(z)-\eta s\bigr)
 \qquad \text{for }t\in\R,\ j\in\Z,\ s=j-ct.
\]
If $3\pi/(4\eta)\le|s|\le\pi/\eta$, then
$5\pi/8\le|\arg\phi_j(z)-\eta s|\le9\pi/8$.
Hence $\cos(\arg\phi_j(z)-\eta s)<0$, so
\[
 V(t,j)<0
 \quad\text{ whenever }\quad
 \frac{3\pi}{4\eta}\le |j-ct|\le\frac{\pi}{\eta}
 \qquad t\in\R,\ j\in\Z.
\]
Choose $\eta$ still smaller so that $\pi/(4\eta)>2$ and
$R<\pi/(4\eta)$, and fix this value for the rest of the proof.
Put $L_0=\pi/\eta>R$. For every $s\in[-R,R]$ and
$i\in\{0,\ldots,P-1\}$,
$|\arg\phi_i(z)-\eta s|<3\pi/8$, so
\[
 \re(e^{-zs}\phi_i(z))
 \ge e^{-|\lambda(\eta)|R}\min_{0\le h<P}|\phi_h(z)|
          \cos(3\pi/8)=:m_R>0.
\]
For $\varepsilon>0$, to be chosen below, define, for
$(t,j)\in\R\times\Z$,
\[
 v(t,j)=\varepsilon (V(t,j))_+
                   \ind_{\{|j-ct|<L_0\}}.
\]
For every $i\in\{0,\ldots,P-1\}$, the function
\[
 s\mapsto\varepsilon\bigl(\re(e^{-zs}\phi_i(z))\bigr)_+
                  \ind_{\{|s|<L_0\}}
 \qquad s\in\R.
\]
is zero near $s=\pm L_0$ and is globally Lipschitz on $\R$.
Let $K$ be the maximum of the Lipschitz constants of these $P$
functions. Then
$\|v(t_1,\cdot)-v(t_2,\cdot)\|_\infty\le |c|K|t_1-t_2|$
for all $t_1,t_2\in\R$. On every bounded interval $I\subset\R$,
$v(t,j)=0$ for all $t\in I$ outside the finite set
$\{j\in\Z:|j|\le L_0+|c|\sup_{t\in I}|t|\}$.
The finitely many remaining coordinate functions are Lipschitz in
$t$, so $t\mapsto v(t,\cdot)$ is differentiable in
$\ell^\infty(\Z)$ for almost every $t\in I$.
Fix such a time $t\in\R$ and a coordinate $j\in\Z$.
If $v(t,j)>0$, then $|j-ct|<L_0$. Since $V(t,j)<0$ when
$3\pi/(4\eta)\le |j-ct|\le L_0$, we have
$|j-ct|<3\pi/(4\eta)$. Since $\pi/(4\eta)>2$, both neighbors
$j-1$ and $j+1$ lie in $\{l\in\Z:|l-ct|<L_0\}$.
At $j$, $v(t,j)=\varepsilon V(t,j)$ and
$\partial_tv(t,j)=\varepsilon\partial_tV(t,j)$; for
$l\in\{j-1,j+1\}$, $v(t,l)\ge\varepsilon V(t,l)$.
The coefficients of $v(t,j-1)$ and $v(t,j+1)$ in
$(\mathcal L_{\Z}v(t,\cdot))_j$ are $D>0$ and $Dq>0$,
respectively. Therefore
\[
 (\mathcal L_{\Z}v)_j-\partial_t v(t,j)
 \ge\varepsilon\bigl((\mathcal L_{\Z}V)_j-\partial_tV(t,j)\bigr)
 =\sigma v(t,j).
\]
If $v(t,j)=0$, nonnegativity and differentiability of
$v(\,\cdot\,,j)$ give $\partial_t v(t,j)=0$, while
$(\mathcal L_{\Z}v(t,\cdot))_j\ge0$. Thus, for both
$v(t,j)>0$ and $v(t,j)=0$,
\begin{equation}\label{eq:kern-compact-subsolution-linear}
 v_t\le\mathcal L_{\Z}v-\sigma v
\end{equation}
at every coordinate for almost every time.

For $|j-ct|<L_0$, the formula for $V$ gives
$(V(t,j))_+\le e^{|\lambda(\eta)|L_0}\max_{0\le i<P}|\phi_i(z)|$.
Choose
\[
 0<\varepsilon\le
 \frac{e^{-|\lambda(\eta)|L_0}}{\max_{0\le i<P}|\phi_i(z)|}
 \min\left\{1,\frac{\sigma}{2(1+C_2)}\right\}.
\]
Then $0\le v(t,j)\le1$ and $C_2v(t,j)\le\sigma/2$
for all $(t,j)\in\R\times\Z$.
For every $j\in\Z$ and $r\in[0,1]$, periodicity and
\eqref{eq:taylor} give $f_j(r)\ge a_jr-C_2r^2$.
Combining this bound with \eqref{eq:kern-compact-subsolution-linear}
yields, for $0<\vartheta\le1$, every $j\in\Z$, and almost every
$t\in\R$, with $v_j=v(t,j)$,
\[
 \begin{aligned}
 (\mathcal A(\vartheta v))_j+f_j(\vartheta v_j)
       -\vartheta\partial_t v_j
 &\ge \vartheta v_j(\sigma-C_2\vartheta v_j)\\
 &\ge (\sigma/2)\vartheta v_j.
 \end{aligned}
\]
Thus $\vartheta v$ is a subsolution for every $0<\vartheta\le1$.
When $c=0$, $v(t,j)$ is independent of $t$, and
$(\mathcal A(\vartheta v(t,\cdot)))_j+f_j(\vartheta v(t,j))>0$
for every $j\in\Z$ with $v(t,j)>0$ and $0<\vartheta\le1$.
Set $\delta=\varepsilon m_R>0$. For $|j-ct|\le R$, the bound
$\re(e^{-z(j-ct)}\phi_j(z))\ge m_R$ gives $v(t,j)\ge\delta$.
For $m\in\Z$, replacing $f_i$ by $f_{i+m}$ replaces
$\phi_i(z)$ by $\phi_{i+m}(z)$, with indices modulo $P$,
and leaves $k(z)$ unchanged. Hence the same $z$, $\sigma$ and
$L_0$ may be used for every $m\in\Z$.
The minimum and maximum of $|\phi_i(z)|$ over $0\le i<P$
and the bounds on $|\arg\phi_i(z)|$ are unchanged by this
permutation. Therefore the same $\varepsilon$, $m_R$ and
$\delta=\varepsilon m_R$ work for every $m\in\Z$.\qedhere

\end{proof}

\begin{lemma}\label{lem:local-general}

Consider \eqref{eq:tree-cauchy} with
$f(x,s)=f_{\iota(x)}(s)$ for $x\in\Tq$ and $s\in[0,1]$,
where $\iota:\Tq\to\{0,\ldots,P-1\}$ is fixed and
$\{f_i\}_{i=0}^{P-1}$ satisfies \eqref{eq:hyp}.
Suppose there are $R_0<\infty$,
$\eta_0>0$ and $m_0>0$ such that every vertex $x$ belongs to a finite
connected set $\Omega_x\subseteq B(x,R_0)$, where
$B(x,R_0)=\{y\in\Tq:d(x,y)\le R_0\}$.
For every $x\in\Tq$, assume that the zero-exterior linearized
operator $\mathcal L_{\Omega_x}$ has an eigenpair
$(\eta_x,\varphi_x)$ with $\eta_x\ge\eta_0$,
$\varphi_x(y)>0$ for $y\in\Omega_x$,
$\mathcal L_{\Omega_x}\varphi_x=\eta_x\varphi_x$,
$\max_{y\in\Omega_x}\varphi_x(y)=1$ and
$\varphi_x(x)\ge m_0$, where $\eta_x$ is its principal eigenvalue.
Then every solution $u:[0,\infty)\times\Tq\to[0,1]$ of
\eqref{eq:tree-cauchy} with $u(0,x)=u_0(x)$ for $x\in\Tq$
and $u_0\not\equiv0$ satisfies $u(t,x)\to1$ as $t\to\infty$
for every fixed $x\in\Tq$.

\end{lemma}
\begin{proof}

Put $a(y)=\partial_s f(y,0)=a_{\iota(y)}$ for $y\in\Tq$,
and choose
$0<\varepsilon\le\min\{1/2,\eta_0/[2(1+C_2)]\}$.
Then $\varepsilon<1$ and $C_2\varepsilon<\eta_0$.
For each $x\in\Tq$, define $z_x:\Tq\to[0,\varepsilon]$ by
\[
 z_x(y)=
 \begin{cases}
  \varepsilon\varphi_x(y),&y\in\Omega_x,\\
  0,&y\notin\Omega_x.
 \end{cases}
\]
On $\Omega_x$, $0<z_x(y)\le\varepsilon$ and
$\mathcal L_{\Omega_x}z_x=\eta_xz_x$.
For $y\in\Tq$ and $r\in[0,1]$, \eqref{eq:taylor} gives
$f(y,r)\ge a(y)r-C_2r^2$.
Consequently, for every $0<\tau\le1$ and $y\in\Omega_x$,
\[
 \begin{aligned}
 D\Delta_{\Tq}(\tau z_x)(y)+f(y,\tau z_x(y))
 &\ge\tau\bigl(D\Delta_{\Tq}z_x(y)+a(y)z_x(y)\bigr)
       -C_2\tau^2z_x(y)^2\\
 &=\tau z_x(y)(\eta_x-C_2\tau z_x(y))>0.
 \end{aligned}
\]
For $y\in\Tq\setminus\Omega_x$, one has $z_x(y)=0$ and
$f(y,0)=0$, so
$D\Delta_{\Tq}(\tau z_x)(y)+f(y,\tau z_x(y))
 =D\tau\sum_{w\sim y}z_x(w)\ge0$. Thus
$\tau z_x$ is a stationary subsolution of \eqref{eq:tree-cauchy}
for every $0<\tau\le1$, and
$D\Delta_{\Tq}(\tau z_x)(y)+f(y,\tau z_x(y))>0$
for every $y\in\Omega_x=\{y\in\Tq:z_x(y)>0\}$.

Let $u:[0,\infty)\times\Tq\to[0,1]$ solve
\eqref{eq:tree-cauchy} with $u(0,\cdot)=u_0\not\equiv0$.
Fix $x_0\in\Tq$ and $t_0>0$. Lemma~\ref{lem:comparison}
gives $u(t_0,y)>0$ for every $y\in\Tq$.
Since $\Omega_{x_0}$ is finite and $z_{x_0}(y)>0$ there,
\[
 \tau=\min\left\{1,
       \min_{y\in\Omega_{x_0}}\frac{u(t_0,y)}{z_{x_0}(y)}\right\}>0,
 \qquad u(t_0,\cdot)\ge\tau z_{x_0}.
\]
Let $U:[0,\infty)\times\Tq\to[0,1]$ solve
\eqref{eq:tree-cauchy} with $U(0,\cdot)=\tau z_{x_0}$.
Since $\tau z_{x_0}$ is a stationary subsolution,
Lemma~\ref{lem:comparison} gives
$U(s,\cdot)\ge\tau z_{x_0}=U(0,\cdot)$ for $s\ge0$.
For each $s\ge0$, autonomy and Lemma~\ref{lem:comparison}, applied
to the initial values $U(s,\cdot)\ge U(0,\cdot)$, give
$U(t+s,\cdot)\ge U(t,\cdot)$ for $t\ge0$.
For every $x\in\Tq$, the function $t\mapsto U(t,x)$ is
nondecreasing and lies in $[0,1]$. Hence
$p(x):=\lim_{t\to\infty}U(t,x)$ exists.
The stationary-limit assertion of Lemma~\ref{lem:compactness},
applied to this bounded solution which is monotone in time, gives
$D\Delta_{\Tq}p(x)+f(x,p(x))=0$ for every $x\in\Tq$.
Since $\tau z_{x_0}\not\equiv0$, Lemma~\ref{lem:comparison}
gives $p(x)\ge U(t,x)>0$ for every $x\in\Tq$ and $t>0$.

Fix $x\in\Tq$. Since $\Omega_x$ is finite and
$p(y),z_x(y)>0$ for all $y\in\Omega_x$, the number
\[
 s_*:=\min\left\{1,\min_{y\in\Omega_x}\frac{p(y)}{z_x(y)}\right\}
 \in(0,1]
\]
satisfies $p(w)\ge s_*z_x(w)$ for every $w\in\Tq$.
If $s_*<1$, there is $y\in\Omega_x$ at which the minimum of
$p(w)/z_x(w)$ over $w\in\Omega_x$ is attained, so
$p(y)=s_*z_x(y)>0$. Since
$p(w)-s_*z_x(w)\ge0$ for every $w\in\Tq$ and
$p(y)-s_*z_x(y)=0$,
\[
 \Delta_{\Tq}p(y)\ge\Delta_{\Tq}(s_*z_x)(y).
\]
Since $p(y)=s_*z_x(y)$, we have
$f(y,p(y))=f(y,s_*z_x(y))$. Thus
\[
 0=D\Delta_{\Tq}p(y)+f(y,p(y))
 \ge D\Delta_{\Tq}(s_*z_x)(y)+f(y,s_*z_x(y))>0,
\]
a contradiction. Thus $s_*=1$, and in particular
$p(x)\ge z_x(x)\ge\varepsilon m_0$ for every vertex $x$.
Thus $\varepsilon m_0\le p(x)\le1$ for every $x\in\Tq$.
The function $(t,x)\mapsto p(x)$ is an entire solution of
\eqref{eq:tree-cauchy} with the positive lower bound
$\varepsilon m_0$. Lemma~\ref{lem:liouville} therefore yields
$p(x)=1$ for every $x\in\Tq$.
Since $u(t_0,\cdot)\ge U(0,\cdot)$,
Lemma~\ref{lem:comparison} gives
$1\ge u(t_0+t,x)\ge U(t,x)\to1$ as $t\to\infty$
for every fixed $x\in\Tq$. Hence $u(t,x)\to1$ as $t\to\infty$
for every fixed $x\in\Tq$.\qedhere

\end{proof}

\subsection{Comparison at a moving boundary}

Let $\mathbb P_j^{\mathrm d}$ denote the law of the continuous-time
chain $X$ started at $j\in\Z$ with generator $\mathcal A$: from each
$i\in\Z$, it jumps to $i-1$ at rate $D$ and to $i+1$ at rate $Dq$. Write
$\mathbb E_j^{\mathrm d}$ for expectation under
$\mathbb P_j^{\mathrm d}$. For each $t\ge0$, the jump count of $X$ on
$[0,t]$ is Poisson with mean $D(q+1)t$. For $\mathbb P_j^{\mathrm d}$ and the laws $\mathbb P_j^\lambda$
defined in Subsection~\ref{sec:periodic-chain}, take a common path space
of integer-valued right-continuous paths with finitely many jumps on
every bounded interval. Let $X$ be its coordinate process and set
\[
 \mathcal F_s=\bigcap_{r>s}\sigma(X_u:0\le u\le r),\quad\text{for }s\ge0.
\]
In this subsection, $\mathcal F_s$ denotes this common right-continuous
filtration for every law. For the construction, transition equations
and strong Markov property, see
\cite[Sections~2.4, 2.6--2.8 and~6.5]{Norris}.

For a family of sets $\Omega_t\subseteq\Z$, indexed by $t\ge0$,
fix a final time $t\ge0$ and a site $j\in\Omega_t$. Under $\mathbb P_j^{\mathrm d}$, the path starts at $X_0=j$.
At elapsed time $s\in[0,t]$, the spatial index is $X_s$ and the
calendar time is $t-s$. We retain only paths satisfying
$X_s\in\Omega_{t-s}$ for every $s\in[0,t]$.
The following lemma constructs the corresponding linear evolution,
including its zero assignments when $\Omega_t$ changes, and compares
it with solutions of \eqref{eq:lattice}.

\begin{lemma}\label{lem:moving-boundary-representation}

Let $m:[0,\infty)\to\R$ be continuously differentiable and suppose that
$\{s\in[0,T]:m(s)\in\Z\}$ is finite for every $T<\infty$.
Put $\Omega_t=\{j\in\Z:j>m(t)\}$ for $t\ge0$.
Let $V:[0,\infty)\times\Z\to[0,\infty)$ be uniformly bounded
on $[0,T]\times\Z$ for every $T<\infty$. For each $T<\infty$,
there is a finite set $E_T\subseteq[0,T]$ such that
$t\mapsto V(t,j)$ is continuous on $[0,T]\setminus E_T$
for every $j\in\Z$.
Let $v_0\in\ell^\infty(\Z)$ be nonnegative, with
$v_0(j)=0$ for $j\in\Z\setminus\Omega_0$.
Let $v(t,j)$ denote the linear evolution from $v(0,j)=v_0(j)$,
with $v(t,j)=0$ for $j\in\Z\setminus\Omega_t$. For $t>0$,
a site entering $\Omega_t$ is initialized at zero, and a site leaving
it is reset to zero, including at an isolated contact that removes
the site only at that instant. Between consecutive times at which $m(t)\in\Z$, the function $v$
satisfies $\partial_t v(t,j)=(\mathcal L_{\Z}v(t,\cdot))_j-V(t,j)v(t,j)$
for $j\in\Omega_t$ and almost every $t>0$. For every $t\ge0$
and $j\in\Z$, this evolution satisfies
\begin{equation}
 \label{eq:moving-feynman-kac}
 v(t,j)=\mathbb E_j^{\mathrm d}\left[
  \exp\left(\int_0^t[a_{X_s}-V(t-s,X_s)]\,\dd s\right)v_0(X_t);
  X_s>m(t-s)\text{ for }0\le s\le t\right].
\end{equation}
If $0\le v(t,j)\le1$ and $V(t,j)\ge C_2v(t,j)$ for every
$t\ge0$ and $j\in\Z$, then $v(t,j)\le w_j(t)$ for every
$t\ge0$ and $j\in\Z$, whenever $w$ solves \eqref{eq:lattice}
and $w_j(0)\ge v_0(j)$ for every $j\in\Z$.

For a solution $w$ of \eqref{eq:lattice} with
$0\le w_j(r)\le1$ for every $r\ge0$ and $j\in\Z$, fix
$t\ge0$ and define
$\tau=\inf\{s\in[0,t]:X_s\le m(t-s)\}$, with $\tau=\infty$ if this
set is empty. If $w_j(0)=0$ for $j>m(0)$, then, for $j>m(t)$,
\begin{equation}
 \label{eq:stopped-boundary-comparison}
 w_j(t)\le
 \mathbb E_j^{\mathrm d}\left[\exp\left(\int_0^\tau a_{X_s}\,\dd s\right);
                       \tau\le t\right].
\end{equation}

\end{lemma}
\begin{proof}

\emph{The linear evolution.}
Fix $T<\infty$ and put $\Lambda=D(q+1)$ and
\[
 C_T=\sup_{0\le r\le T}\sup_{i\in\Z}|a_i-V(r,i)|<\infty.
\]
First fix $\Omega\subseteq\Z$, use the initial datum
$v_0\mathbf1_\Omega$, and impose zero values on $\Z\setminus\Omega$. Sum over $n\in\N_0$, the two possible
directions of each jump, and the ordered jump times in $[0,t]$. Weight each path by
$\exp(\int_0^t[a_{X_s}-V(t-s,X_s)]\,\dd s)v_0(X_t)$ and retain
only paths satisfying $X_s\in\Omega$ for every $s\in[0,t]$.
For $0\le t\le T$, the absolute value
of the sum of the terms with $n$ jumps is at most
\[
 \|v_0\|_\infty e^{C_TT}e^{-\Lambda t}
             \frac{(\Lambda t)^n}{n!}
 \qquad n\in\N_0.
\]
This bound gives absolute convergence of the sums and integrals.
Conditioning on the first short path interval gives
$\partial_t v(t,j)=(\mathcal L_{\Z}v(t,\cdot))_j-V(t,j)v(t,j)$
for $j\in\Omega$ and almost every $t\in(0,T)$.
The path sum also gives
$\|v(t,\cdot)\|_\infty\le\|v_0\|_\infty e^{C_TT}$.
By \eqref{eq:kern-diffusion-norm},
$\|\mathcal Av(t,\cdot)\|_\infty\le2\Lambda\|v(t,\cdot)\|_\infty$.
The bound from the path sum therefore bounds every coordinate
derivative between switching times by
$(2\Lambda+C_T)\|v_0\|_\infty e^{C_TT}$. Integrating at each
coordinate and taking the supremum shows that $v$ is locally
Lipschitz in $\ell^\infty$ on each interval between switches. For two evolutions $v$ and $\widetilde v$ with the
same initial value, the coordinatewise integral equations give
\[
 \|v(t,\cdot)-\widetilde v(t,\cdot)\|_\infty
 \le(2\Lambda+C_T)\int_0^t
       \|v(r,\cdot)-\widetilde v(r,\cdot)\|_\infty\,\dd r
 \qquad 0\le t\le T.
\]
Gronwall's inequality therefore gives uniqueness on each interval with
a fixed active set, after translating its initial time to zero.

The set $\{r\in[0,T]:m(r)\in\Z\}$ is finite by hypothesis.
At each $r\in\{r\in[0,T]:m(r)\in\Z\}$, retain the current
values on $\Omega_r$ and assign zero at every
$j\in\Z\setminus\Omega_r$. At an isolated
contact, make this restriction before continuing on the same active
set as before the contact. Composing the fixed-set evolutions with
these restrictions gives \eqref{eq:moving-feynman-kac}: the path must
satisfy $X_s\in\Omega_{t-s}$ at every time, including each switching
time. At a site
entering at time $r>0$, a surviving path evaluated at time $r+h$ must
jump during its first $h$ units of path time to reach a site active
at calendar time $r$. Its contribution is therefore bounded by
\[
 \|v_0\|_\infty e^{C_TT}(1-e^{-\Lambda h}),
 \qquad 0<h\le T-r.
\]
Hence $v(r+h,j)\to0$ as $h\downarrow0$ at every site $j$
entering at time $r$. At a removed site $j$, the path condition
$X_s\in\Omega_{t-s}$ sets $v(r,j)=0$ at the switching time $r$.
At an isolated contact, the same condition gives $v(r,j)=0$;
the bound $\|v_0\|_\infty e^{C_TT}(1-e^{-\Lambda h})$ gives
$v(r+h,j)\to0$ as $h\downarrow0$ after that contact. The unique bounded evolutions on the intervening intervals
therefore give exactly the stated resets and representation.
Formula \eqref{eq:moving-feynman-kac} is the Feynman--Kac formula
for the killed process with the calendar-time coordinate $t-s$
adjoined; see
\cite{FitzsimmonsPitman} for its additive-functional formulation.

\emph{Nonlinear comparison.}
For $j\in\Omega_t$, the hypotheses $V(t,j)\ge C_2v(t,j)$
and $0\le v(t,j)\le1$,
together with \eqref{eq:taylor}, give almost everywhere in time
\[
 \begin{aligned}
 \partial_t v(t,j)
 &= (\mathcal Av(t,\cdot))_j+a_jv(t,j)-V(t,j)v(t,j)\\
 &\le (\mathcal Av(t,\cdot))_j+a_jv(t,j)-C_2v(t,j)^2\\
 &\le (\mathcal Av(t,\cdot))_j+f_j(v(t,j)).
 \end{aligned}
\]
Between consecutive times with $m(t)\in\Z$, every $j\notin\Omega_t$
satisfies $v(t,j)=0$, $\partial_t v(t,j)=0$ and
$(\mathcal Av(t,\cdot))_j=D(v(t,j-1)+qv(t,j+1))\ge0$,
while $f_j(0)=0$.
Thus the zero extension of $v$ is a subsolution on the full lattice
on each such interval. On each interval between switches, $v$ is locally Lipschitz
with values in $\ell^\infty(\Z)$, as proved above, and $w$ belongs
to $C^1([0,\infty);\ell^\infty(\Z))$ by Lemma~\ref{lem:comparison}.
Thus that lemma applies on the full lattice. Starting from
$v_0\le w(0,\cdot)$, we obtain $v\le w$ up to the first switch.
At every switch, removed and entering values are zero, including
the reset and reentry values at an isolated contact. Since $w\ge0$,
these assignments preserve $v\le w$. Induction over the finitely
many switches on $[0,T]$ therefore gives $v\le w$ there. Since $T$ is arbitrary, $v(t,j)\le w_j(t)$ for every
$t\ge0$ and $j\in\Z$.

\emph{The stopped identity.}
For the stopped estimate \eqref{eq:stopped-boundary-comparison},
define $\mathcal R_j(r)=a_jr-f_j(r)\ge0$ for $j\in\Z$ and $r\in[0,1]$, and put
$T_{\mathrm{stop}}=t\wedge\tau$. The first exit occurs at a jump time or at one of
the finitely many deterministic times $s$ with $m(t-s)\in\Z$.
Consequently $\{\tau\le r\}$ is determined by $(X_s)_{0\le s\le r}$,
so $T_{\mathrm{stop}}$ is a bounded stopping time for $\mathcal F_s$.
The coordinate equation
$\dot w_j=(\mathcal L_{\Z}w)_j-\mathcal R_j(w_j)$ then gives
\begin{equation}
 \label{eq:stopped-residual-identity}
 \begin{aligned}
 w_j(t)
 &=\mathbb E_j^{\mathrm d}\left[
     e^{\int_0^{T_{\mathrm{stop}}} a_{X_r}\,\dd r}w_{X_{T_{\mathrm{stop}}}}(t-T_{\mathrm{stop}})\right]\\
 &\quad-\mathbb E_j^{\mathrm d}\left[
     \int_0^{T_{\mathrm{stop}}} e^{\int_0^s a_{X_r}\,\dd r}
       \mathcal R_{X_s}\bigl(w_{X_s}(t-s)\bigr)\,\dd s\right].
 \end{aligned}
 \qquad t\ge0,\ j\in\Z,\ j>m(t).
\end{equation}
More explicitly, for $0\le s\le t$ the process
\[
 e^{\int_0^s a_{X_r}\,\dd r}w_{X_s}(t-s)-w_j(t)
 -\int_0^s e^{\int_0^v a_{X_r}\,\dd r}
       \mathcal R_{X_v}(w_{X_v}(t-v))\,\dd v
\]
is a martingale. Indeed, Lemma~\ref{lem:comparison} gives the
$C^1$ time regularity of $w$, and its equation yields
\[
 G(s,i)=w_i(t-s),
 \qquad
 (\partial_s+\mathcal A)G(s,i)
       =-a_iw_i(t-s)+\mathcal R_i(w_i(t-s))
 \qquad 0\le s\le t,\ i\in\Z.
\]
Multiplication by $e^{\int_0^s a_{X_r}\,\dd r}$ cancels the
$-a_iw_i(t-s)$ term. Compensating the left and right jumps, of rates
$D$ and $Dq$, gives the stated martingale in the natural filtration.
By \eqref{eq:hyp}, $0\le\mathcal R_j(r)\le a_{\max}$ for
$0\le r\le1$. Together with $0\le w\le1$, this bounds the
martingale's absolute value on $[0,t]$ by
$e^{a_{\max}t}(1+a_{\max}t)+1$. Its right-continuity and this bound
extend the martingale identity to $\mathcal F_s$ by letting the
conditioning time $s+\varepsilon$ decrease to $s$. To justify stopping in
continuous time, approximate $T_{\mathrm{stop}}$ from above by the finite-grid
stopping times $T_{\mathrm{stop},n}=(2^{-n}\lceil2^nT_{\mathrm{stop}}\rceil)\wedge t$.
Discrete-time optional stopping, as in
\cite[Section~4.8]{Durrett2019}, applies on each such grid.
Right-continuity and bounded convergence as $n\to\infty$ then give the identity at
$T_{\mathrm{stop}}$, proving \eqref{eq:stopped-residual-identity}. For continuous-time stopping
conventions, see \cite[Sections~4.1 and~6.5]{Norris}. On $\{\tau>t\}$, the path ends
at $X_t>m(0)$, where $w_{X_t}(0)=0$ by the initial-data hypothesis.
On $\{\tau\le t\}$, one has $T_{\mathrm{stop}}=\tau$ and
$w_{X_\tau}(t-\tau)\le1$. Substituting these bounds in
\eqref{eq:stopped-residual-identity} and using
$\mathcal R_j(r)\ge0$ for $j\in\Z$ and $r\in[0,1]$ gives
\eqref{eq:stopped-boundary-comparison}.\qedhere

\end{proof}

Fix $b\in\Z$. For the rest of this subsection, write
\[
 \lambda=\lambda_+,
 \qquad c=c_+,
 \qquad \kappa=\frac{3}{2\lambda},
 \qquad m_b(t)=b+ct-\kappa\log(1+t),\quad\text{for }t\ge0.
\]
By \eqref{eq:parameters} and \eqref{eq:critical-pairs},
$k(\lambda)=c\lambda$, $k'(\lambda)=c$ and $\kappa\lambda=3/2$. The derivative
$m_b'(t)=c-\kappa/(1+t)$ is strictly increasing, so $m_b$ meets
each integer level at most twice. Its bounded range on $[0,T]$
contains finitely many integers. Hence
$\{t\in[0,T]:m_b(t)\in\Z\}$ is finite for every $T<\infty$, as
required by Lemma~\ref{lem:moving-boundary-representation}.

\begin{lemma}\label{lem:sharp-step-upper}

Assume \eqref{eq:hyp}, and let $w$ be the $[0,1]$-valued solution of
\eqref{eq:lattice} with $0\le w_j(0)\le\mathbf1_{\{j\le b\}}$
for an integer $b$. Then there is $C<\infty$ such that
\begin{equation}
 \label{eq:sharp-step-upper}
 w_j(t)\le C(1+y)e^{-\lambda y},\qquad
 y=j-m_b(t)>0,\quad\text{ for }t\ge0,\ j\in\Z.
\end{equation}
The constant $C$ is independent of $b\in\Z$.

\end{lemma}
\begin{proof}

At $t=0$, \eqref{eq:sharp-step-upper} follows from $w_j(0)=0$ for
$j>b$. Fix $t>0$ and $j>m_b(t)$, and put
\[
 \begin{aligned}
 y&=j-m_b(t),\\
 Y_s&=X_s-m_b(t-s)=y+X_s-j+cs-b_{t,0}(s),\quad 0\le s\le t,\\
 \tau&=\inf\{s\in[0,t]:Y_s\le0\},
 \end{aligned}
\]
with $\tau=\infty$ when the set is empty. Here
$b_{t,0}(s)=\kappa\log((1+t)/(1+t-s))$ for $0\le s\le t$,
by \eqref{eq:backward-boundary-definition}. Since $m_b(0)=b$ and $w_i(0)=0$ for every integer $i>b$,
the zero initial-data condition in
Lemma~\ref{lem:moving-boundary-representation} holds, and
\eqref{eq:stopped-boundary-comparison} applies to this $\tau$.
Put $H(i)=e^{-\lambda i}\phi_i(\lambda)$ for $i\in\Z$.
For each stopping time
$0\le T_{\mathrm{stop}}\le t$ and each bounded $\mathcal F_{T_{\mathrm{stop}}}$-measurable random
variable $Z$, the change-of-measure identity is
\begin{equation}
 \label{eq:stopped-change-of-law}
 \mathbb E_j^{\mathrm d}\left[e^{\int_0^{T_{\mathrm{stop}}} a_{X_r}\,\dd r}Z\right]
 =\mathbb E_j^\lambda\left[
       e^{k(\lambda)T_{\mathrm{stop}}}\frac{H(j)}{H(X_{T_{\mathrm{stop}}})}Z\right].
\end{equation}
At each deterministic time $s\ge0$, the density of $\mathbb P_j^{\mathrm d}$
relative to $\mathbb P_j^\lambda$ on $\mathcal F_s$ is
\[
 D_s=\exp\left(k(\lambda)s-\int_0^s a_{X_r}\,\dd r\right)
                         \frac{H(j)}{H(X_s)}.
\]
Under $\mathbb P_j^\lambda$, the left and right jump rates from
$i\in\Z$ are, respectively,
\[
 \alpha_i^\lambda=D\frac{H(i-1)}{H(i)},\qquad
 \beta_i^\lambda=Dq\frac{H(i+1)}{H(i)},\qquad
 \alpha_i^\lambda+\beta_i^\lambda
       =D(q+1)+k(\lambda)-a_i.
\]
For a path with finitely many jumps, the ratios of the jump factors
telescope to $H(j)/H(X_s)$, and the ratio of the exponential holding
factors is $\exp(\int_0^s[k(\lambda)-a_{X_r}]\,\dd r)$.
Integrating over the ordered jump times proves that $D_s$ is the
Radon--Nikodym density of $\mathbb P_j^{\mathrm d}$ with respect to
$\mathbb P_j^\lambda$ on $\sigma(X_r:0\le r\le s)$.
Let $N_t$ be the actual number of jumps up to $t$. Periodicity and
positivity of $\phi$ give
\[
 \sup_{0\le s\le t}D_s\le C\exp(Ct+C N_t),\qquad
 \mathbb E_j^\lambda
       \bigl[\sup_{0\le s\le t}D_s^p\bigr]<\infty\quad\text{ for } p>1,
\]
because $N_t$ is bounded by the number of rings of a Poisson clock
of rate $\max_i(\alpha_i^\lambda+\beta_i^\lambda)$ in the construction
of Subsection~\ref{sec:periodic-chain}. The bound $\mathbb E_j^\lambda[\sup_{0\le s\le t}D_s^p]<\infty$
for $p>1$ also extends the density identity to $\mathcal F_s$ and
to the common completion of the two laws on every fixed finite horizon. The finite-horizon
densities are strictly positive, so the two laws have the same null
sets and the identity extends to their completion. Indeed, if
$A\in\mathcal F_s$, then, for every $\varepsilon>0$, $A$ agrees,
up to a common null set, with an event in
$\sigma(X_r:0\le r\le s+\varepsilon)$.
Right-continuity of $(D_r)_{r\ge0}$ and
$\mathbb E_j^\lambda[\sup_{0\le r\le s+1}D_r^p]<\infty$ for $p>1$ give
\[
 \mathbb P_j^{\mathrm d}(A)
 =\lim_{\varepsilon\downarrow0}
       \mathbb E_j^\lambda[D_{s+\varepsilon}\mathbf1_A]
 =\mathbb E_j^\lambda[D_s\mathbf1_A].
\]
Thus $\mathbb E_j^\lambda D_s=1$ and
$\mathbb E_j^\lambda[D_t\mid\mathcal F_s]=D_s$ for $s\le t$.
The bound
$\mathbb E_j^\lambda[\sup_{0\le r\le t}D_r^p]<\infty$ for $p>1$
gives uniform integrability of the stopped densities. Apply the
conditional-expectation form of
discrete-time optional stopping from
\cite[Section~4.8]{Durrett2019} at
$T_{\mathrm{stop},n}=(2^{-n}\lceil2^nT_{\mathrm{stop}}\rceil)\wedge t$. Since the process $s\mapsto D_s$ is
right-continuous, $D_{T_{\mathrm{stop},n}}\to D_{T_{\mathrm{stop}}}$ as $n\to\infty$, almost surely and in $L^1$, by
the displayed $L^p$ bound. Testing against bounded
$\mathcal F_{T_{\mathrm{stop}}}$-measurable variables and passing to the limit gives
$D_{T_{\mathrm{stop}}}=\mathbb E_j^\lambda[D_t\mid\mathcal F_{T_{\mathrm{stop}}}]$. Testing this identity
against $e^{\int_0^{T_{\mathrm{stop}}} a_{X_r}\,\dd r}Z$, which is bounded for fixed $t$,
proves \eqref{eq:stopped-change-of-law}. On $\{\tau\le t\}$, $-1\le Y_\tau\le0$: a downward jump of $X$
has size one, and a first exit between jumps has $Y_\tau=0$. Apply \eqref{eq:stopped-change-of-law} with
$T_{\mathrm{stop}}=t\wedge\tau$ and $Z=\mathbf1_{\{\tau\le t\}}$.
Using $0<\min_{0\le i<P}\phi_i(\lambda)
\le\max_{0\le i<P}\phi_i(\lambda)<\infty$ in
\eqref{eq:stopped-boundary-comparison} gives
\begin{equation}
 \label{eq:exit-weight-bound}
 w_j(t)\le C e^{-\lambda y}
 \mathbb E_j^\lambda\left[
       \left(\frac{1+t}{1+t-\tau}\right)^{3/2};\ \tau\le t\right].
\end{equation}
The power $3/2$ in \eqref{eq:exit-weight-bound} follows from
$c\tau+X_\tau-j=-y+Y_\tau+\kappa\log((1+t)/(1+t-\tau))$
and $\kappa\lambda=3/2$.

An exit at a jump at time $s\in(0,t]$ requires a downward jump from
\[
 i(s)=\lfloor m_b(t-s)\rfloor+1,
 \qquad 0<i(s)-m_b(t-s)\le1.
\]
For $0\le r\le s$, the identity
\[
 b_{s,t-s}(r)=\kappa\log\frac{1+t}{1+t-r}=b_{t,0}(r)
\]
shows that
$\mathbb P_j^\lambda\{\tau>s,\ X_s=i(s)\}
=\mathcal K_{s,t-s}^{y,\lambda}(j,i(s))$.
By \eqref{eq:backward-log-upper} and
$0<i(s)-m_b(t-s)\le1$,
\[
 \mathcal K_{s,t-s}^{y,\lambda}(j,i(s))
 \le C(1+y)(1+s)^{-3/2}
 \qquad 0<s\le t.
\]
Let $F:[0,t]\to[0,\infty)$ be Borel measurable. Summing the ordered paths
whose first exit is a downward jump gives
\[
 \begin{aligned}
 &\mathbb E_j^\lambda[F(\tau);\ \tau\le t,
                     \ \tau\text{ is a jump time of }X]\\
 &\quad=\int_0^t F(s)\alpha_{i(s)}^\lambda
           \mathcal K_{s,t-s}^{y,\lambda}(j,i(s))\,\dd s\\
 &\quad\le C(1+y)\int_0^t F(s)(1+s)^{-3/2}\,\dd s.
 \end{aligned}
\]
Indeed, the path before the exiting jump survives and ends at
$i(s)$, and the exiting jump contributes
$\alpha_{i(s)}^\lambda\,\dd s$. Jumps at the finitely many prescribed
boundary-contact times have probability zero. The upper bound $\max_{0\le i<P}\alpha_i^\lambda$ for the
downward-jump rate is finite. Bounded rates justify the
ordered-jump expansion and the summation over all jump counts;
see \cite[Section~2.6]{Norris}.
An exit between jumps can occur only at a deterministic time in
\[
 \mathcal C_t=\{s\in(0,t]:m_b(t-s)\in\Z\}.
\]
The number of points of $\mathcal C_t$ in any interval of length
one is at most $2(\lceil |c|+\kappa\rceil+2)$. Indeed,
$s\mapsto m_b(t-s)$ has at most two monotone pieces and derivative
of magnitude at most $|c|+\kappa$. On either piece within a unit
interval, its range has length at most $|c|+\kappa$ and therefore
contains at most $\lceil |c|+\kappa\rceil+2$ integer levels.

For $s\in\mathcal C_t$, write $k=m_b(t-s)\in\Z$ and define
$E_s=\{\tau=s,\ X_{s-}=X_s\}$. Thus $E_s$ is the event of a
first exit at $s$ with no jump of $X$ at that time. Such an exit requires
$m_b(t-r)<k$ for all $r<s$ sufficiently close to $s$. For a time $s\in\mathcal C_t$ satisfying
$m_b(t-r)<k$ for all $r<s$ sufficiently close to $s$,
choose deterministic $0\le s_n<s$ with $s_n\uparrow s$ as $n\to\infty$.
Then $z_n=k-m_b(t-s_n)>0$ for all large $n$ and $z_n\to0$.
The chain has finitely many jumps and is almost surely continuous
at the deterministic time $s$. Hence, on $E_s$, it is constant at $k$
immediately before $s$, which gives
\[
 E_s\subseteq\liminf_{n\to\infty}
        \{\tau>s_n,\ X_{s_n}=k\}
 \quad\text{up to a null set}.
\]
Since $b_{s_n,t-s_n}(r)=b_{t,0}(r)$ for $0\le r\le s_n$,
the probability of $\{\tau>s_n,\ X_{s_n}=k\}$ is
$\mathcal K_{s_n,t-s_n}^{y,\lambda}(j,k)$.
Fatou's inequality and \eqref{eq:backward-log-upper}, with
duration $s_n$, offset $t-s_n$ and terminal height $z_n$, give
\[
 \mathbb P_j^\lambda(E_s)
 \le\liminf_{n\to\infty}\mathcal K_{s_n,t-s_n}^{y,\lambda}(j,k)
 \le C(1+y)(1+s)^{-3/2}.
\]
The bound $\mathbb P_j^\lambda(E_s)\le C(1+y)(1+s)^{-3/2}$
holds at every possible first-exit time between jumps, including
$s=t$, because $s_n<t$ and $t-s_n\ge0$. If
$m_b(t-r)<k$ does not hold for all $r<s$ sufficiently close to $s$,
then $\mathbb P_j^\lambda(E_s)=0$.
Consequently, for every Borel measurable $F:[0,t]\to[0,\infty)$,
\begin{equation}
 \label{eq:exit-measure-bound}
 \mathbb E_j^\lambda[F(\tau);\tau\le t]
 \le C(1+y)\left[
     \int_0^t\frac{F(s)}{(1+s)^{3/2}}\,\dd s
     +\sum_{s\in\mathcal C_t}\frac{F(s)}{(1+s)^{3/2}}
                         \right].
\end{equation}
Up to a null set, $\{\tau\le t\}$ is the union of
$\{\tau\le t,\ \tau\text{ is a jump time of }X\}$ and
$\bigcup_{s\in\mathcal C_t}E_s$: the path $X$ has finitely many
jumps on $[0,t]$ and $Y_s$ is continuous between jumps.

Take $F(s)=((1+t)/(1+t-s))^{3/2}$ for $0\le s\le t$. For $0\le s\le t$,
\[
 \frac{(1+t)^{3/2}}{(1+s)^{3/2}(1+t-s)^{3/2}}
 \le 2^{3/2}\bigl((1+s)^{-3/2}+(1+t-s)^{-3/2}\bigr).
\]
For $s\le t/2$, use $1+t-s\ge(1+t)/2$; for $s\ge t/2$, use
$1+s\ge(1+t)/2$. Each term on the right has an integral over
$[0,t]$ bounded independently of $t$.
Put $N_c=2(\lceil |c|+\kappa\rceil+2)$. By the bound on the number
of times in $\mathcal C_t$, each interval of length one contains at
most $N_c$ points of $\mathcal C_t$. Partitioning $[0,t]$ from each endpoint into such
intervals gives
\[
 \sum_{s\in\mathcal C_t}
 \bigl((1+s)^{-3/2}+(1+t-s)^{-3/2}\bigr)
 \le 2N_c\sum_{n=0}^\infty(1+n)^{-3/2}<\infty.
\]
The preceding integral and series bounds, substituted into
\eqref{eq:exit-measure-bound} with
$F(s)=((1+t)/(1+t-s))^{3/2}$, give
$\mathbb E_j^\lambda[F(\tau);\tau\le t]\le C(1+y)$. By \eqref{eq:exit-weight-bound}, this proves
\eqref{eq:sharp-step-upper}. The bound on the number of times in
$\mathcal C_t$ depends only on $c,\kappa$. The estimate
\eqref{eq:backward-log-upper} and the positive upper and lower bounds
for the entries of $\phi(\lambda)$ are uniform in the starting phase.
Thus $C$ is independent of $b$.\qedhere

\end{proof}

\begin{lemma}\label{lem:sharp-step-lower}

Assume \eqref{eq:hyp}, and let $w$ solve \eqref{eq:lattice} from
$w_j(0)=\mathbf1_{\{j\le b\}}$, where $b\in\Z$. For every
$\theta\in(0,1)$ there are $C_\theta,T_\theta<\infty$, independent of
$b$, such that
\begin{equation}
 \label{eq:sharp-step-lower}
 w_j(t)\ge\theta\qquad
 \text{if }t\ge T_\theta\text{ and }
 j\in\Z,\ j\le m_b(t)-C_\theta.
\end{equation}

\end{lemma}
\begin{proof}

\emph{A lower bound near $m_b(t)$.}
Let $Y_*$ be the lower height threshold in
\eqref{eq:backward-log-lower} with $A=1$. Choose a fixed integer
$z_0\ge\max\{Y_*,1\}$, put
\[
 m(t)=m_{b-z_0}(t),
 \qquad j_0=b
 \qquad t\ge0.
\]
Recall $H(i)=e^{-\lambda i}\phi_i(\lambda)$ for $i\in\Z$.
Since $j_0-m(0)=z_0>0$, define
$\delta_{j_0}(j)=\mathbf1_{\{j=j_0\}}$ for $j\in\Z$, and let $W$ be the evolution
$W_t=\mathcal L_{\Z}W$ on $j>m(t)$, with initial datum
$\delta_{j_0}$ and the zero exterior values and resets in
Lemma~\ref{lem:moving-boundary-representation}.
For $t\ge0$ and $j\in\Z$ with $y=j-m(t)>0$, formulas \eqref{eq:moving-feynman-kac} with
$V=0$ and \eqref{eq:stopped-change-of-law} with $T_{\mathrm{stop}}=t$ apply on
$\{X_t=j_0,\ X_s>m(t-s)\text{ for }0\le s\le t\}$.
On $\{X_t=j_0,\ X_s>m(t-s)\text{ for }0\le s\le t\}$,
the factor $e^{k(\lambda)t}H(j)/H(X_t)$ in
\eqref{eq:stopped-change-of-law} is
\[
 e^{k(\lambda)t-\lambda(j-j_0)}
       \frac{\phi_j(\lambda)}{\phi_{j_0}(\lambda)}
 =\frac{\phi_j(\lambda)}{\phi_{j_0}(\lambda)}
       e^{\lambda z_0}(1+t)^{3/2}e^{-\lambda y},
\]
since $j-j_0=ct-\kappa\log(1+t)+y-z_0$ and
$k(\lambda)=c\lambda$. Therefore
\begin{equation}
 \label{eq:linear-boundary-seed}
 W(t,j)=
 \frac{\phi_j(\lambda)}{\phi_{j_0}(\lambda)}
 e^{\lambda z_0}(1+t)^{3/2}e^{-\lambda y}
 \mathcal K_{t,0}^{y,\lambda}(j,j_0)
 \qquad t\ge0,\ j\in\Z,\ y=j-m(t)>0.
\end{equation}
Indeed, $c=k'(\lambda)$ by \eqref{eq:critical-pairs}, and
$m(t)-m(t-s)=cs-b_{t,0}(s)$. Hence the height above the moving
boundary is
\[
 X_s-m(t-s)=y+X_s-j+k'(\lambda)s-b_{t,0}(s)
 \qquad 0\le s\le t.
\]
On $\{X_t=j_0\}$, the terminal height above the boundary is
$j_0-m(0)=z_0$.
Thus
$\mathbb P_j^\lambda\{X_t=j_0,\ X_s>m(t-s)\text{ for }0\le s\le t\}
=\mathcal K_{t,0}^{y,\lambda}(j,j_0)$ as defined in
\eqref{eq:backward-boundary-definition}.
Define $\omega(y)=(1+y)e^{-\lambda y}$ for $y>0$.
By \eqref{eq:backward-log-upper}, the positive upper and lower bounds
for the entries of $\phi(\lambda)$ and the fixed value of $z_0$ give
a constant $C_W$, independent of $b$, such that
\[
 0\le W(t,j)\le C_W\omega(y)
\]
for every $t\ge0$ and $j\in\Z$ with $y=j-m(t)>0$.
Set $W(t,j)=0$ for $t\ge0$ and $j\le m(t)$.

Choose $\varepsilon>0$ so that
\[
 \varepsilon\le1,
 \qquad
 \varepsilon C_W\sup_{y>0}\omega(y)\le1.
\]
On $j>m(t)$, let $v$ satisfy
$v_t=\mathcal L_{\Z}v-Vv$, with initial datum
$\varepsilon\delta_{j_0}$, the zero exterior values and resets of
Lemma~\ref{lem:moving-boundary-representation}, and the potential
\[
 V(t,j)=C_2\varepsilon C_W\omega(j-m(t)),
 \quad\text{ for }t\ge0,\ j\in\Z,\ j>m(t).
\]
Set $V(t,j)=0$ for $j\le m(t)$.
The potential $V$ is bounded and continuous between the finitely
many switching times on each bounded time interval, so
\eqref{eq:moving-feynman-kac} applies. Since $V\ge0$, \eqref{eq:moving-feynman-kac} gives
$0\le v(t,j)\le\varepsilon W(t,j)\le1$ for every $t\ge0$ and
$j\in\Z$. Consequently $V(t,j)\ge C_2v(t,j)$ for every such
$(t,j)$, and
\eqref{eq:taylor} gives, for $j\in\Z$ with $j>m(t)$ and
almost every $t\ge0$,
\[
 v_t-\mathcal Av-f_j(v)
 =a_jv-Vv-f_j(v)\le C_2v^2-Vv\le0.
\]
Since $j_0=b$ and $\varepsilon\le1$, the initial values satisfy
$\varepsilon\delta_{j_0}\le w(0,\cdot)$. The potential inequality
$V\ge C_2v$ and the bound $0\le v\le1$ verify all hypotheses of
the nonlinear comparison in
Lemma~\ref{lem:moving-boundary-representation}. That lemma includes
the zero assignments at every switch, so it gives $v(t,j)\le w_j(t)$
for all $t\ge0$ and $j\in\Z$.

Set $Y=z_0$, $B=z_0+P+1$ and
$Y_s=X_s-m(t-s)$ for $0\le s\le t$. If $Y\le y\le B$,
then both endpoint heights $y,z_0$ lie in $[Y,B]$, with $Y\ge Y_*$.
Take $t\ge B^2$ above the time thresholds in
Lemmas~\ref{lem:backward-log-kernel} and~\ref{lem:bridge-occupation}.
By \eqref{eq:backward-log-lower} with $A=1$,
\[
 \mathcal K_{t,0}^{y,\lambda}(j,j_0)
       \ge c_1yz_0t^{-3/2}>0
 \qquad j\in\Z,\ y=j-m(t)\in[Y,B],
\]
where $c_1$ is independent of $t,j,b$. On the event
$X_t=j_0$ and $Y_s>0$ for every $s\in[0,t]$, one has
$V(t-s,X_s)=C_2\varepsilon C_W\omega(Y_s)$.
Applying \eqref{eq:moving-feynman-kac} and
\eqref{eq:stopped-change-of-law} at $T_{\mathrm{stop}}=t$ and dividing by
$\varepsilon W(t,j)>0$ cancels the factor
$\varepsilon e^{k(\lambda)t}H(j)/H(j_0)$ and gives
\[
\begin{gathered}
 \frac{v(t,j)}{\varepsilon W(t,j)}
 =\mathbb E_j^\lambda\left[
  \exp\left(-C_2\varepsilon C_W\int_0^t\omega(Y_s)\,\dd s\right)
       \,\middle|\, X_t=j_0,\ Y_s>0\text{ for }0\le s\le t\right]
 \\\text{ for }j\in\Z,\ y=j-m(t)\in[Y,B].
\end{gathered}
\]
Since $\lambda=\lambda_+>0$ and $y,z_0\in[Y,B]$ with $Y\ge Y_*$,
\eqref{eq:bridge-occupation}, applied with terminal site $l=j_0$,
bounds the conditional mean of $\int_0^t\omega(Y_s)\,\dd s$ by
$C_B$, independently of $t,j,b$. Jensen's inequality for the convex
function $x\mapsto e^{-C_2\varepsilon C_Wx}$ now gives
$v(t,j)\ge\varepsilon e^{-C\varepsilon}W(t,j)$, where
$C=C_2C_WC_B$. In view of \eqref{eq:linear-boundary-seed}, the kernel lower bound
\eqref{eq:backward-log-lower} and $(1+t)^{3/2}t^{-3/2}\ge1$ imply
\[
 v(t,j)\ge
 \varepsilon e^{-C\varepsilon}
 \frac{\min_{0\le i<P}\phi_i(\lambda)}{\max_{0\le i<P}\phi_i(\lambda)}
 e^{\lambda z_0-\lambda B}c_1Yz_0>0
 \qquad j\in\Z,\ y=j-m(t)\in[Y,B].
\]
Choose $\delta\in(0,1)$ such that
\[
 \delta<\varepsilon e^{-C\varepsilon}
 \frac{\min_{0\le i<P}\phi_i(\lambda)}{\max_{0\le i<P}\phi_i(\lambda)}
 e^{\lambda z_0-\lambda B}c_1Yz_0.
\]
Choose $T\ge B^2$ so that the lower estimate in
Lemma~\ref{lem:backward-log-kernel} and the conditional occupation
estimate in Lemma~\ref{lem:bridge-occupation} both hold for
$t\ge T$, terminal height $z_0$, and every initial height
$y\in[Y,B]$. We obtain
\begin{equation}
 \label{eq:positive-at-logarithmic-position}
 w_j(t)\ge v(t,j)\ge\delta
 \quad\text{if }Y\le j-m(t)\le Y+P+1,
 \quad\text{ for }t\ge T,\ j\in\Z.
\end{equation}
The choices of $z_0,Y,B$ and $\varepsilon$ are independent of $b$.
The kernel and occupation estimates are uniform in $j,j_0$, and the
ratios of the periodic vector $\phi(\lambda)$ have fixed positive
bounds. Thus $\delta$ and $T$ are independent of $b$.

\emph{Extension by periodicity.}
Periodicity and Lemma~\ref{lem:comparison} imply
$w_{j+P}(t)\le w_j(t)$ for every $t\ge0$ and $j\in\Z$. Indeed, the shifted function
$j\mapsto w_{j+P}(t)$ solves the same equation and starts from the
smaller initial datum $\mathbf1_{\{j\le b-P\}}$. Given $j\in\Z$ with $j\le m(t)+Y=m_b(t)$, let $n$ be the
smallest nonnegative integer for which $j+nP\ge m(t)+Y$.
This choice gives
\[
 Y\le j+nP-m(t)\le Y+P.
\]
By \eqref{eq:positive-at-logarithmic-position},
$w_{j+nP}(t)\ge\delta$ for $t\ge T$. Repeated comparison with the
period shift gives $w_j(t)\ge w_{j+nP}(t)$, and hence
\[
 w_j(t)\ge\delta
 \quad\text{ for }j\in\Z,\ j\le m_b(t),\ t\ge T.
\]

\emph{Raising the level.}
To raise $\delta$ to the prescribed level $\theta$, put
$g(s)=\min_i f_i(s)$ for $s\in[0,1]$, where the minimum is over
$i\in\{0,\ldots,P-1\}$. The function $g$ is Lipschitz on $[0,1]$,
is positive on $(0,1)$, and vanishes at $0$ and $1$. Let $z:[0,\infty)\to(0,1)$ solve
$z'(s)=g(z(s))$, $z(0)=\delta$. Then $z(t)\uparrow1$ as $t\to\infty$, so choose a fixed
$\tau>0$ with $z(\tau)>\theta$. Extend $g$ by setting $g(s)=0$ for $s\notin[0,1]$. Since $g$
vanishes at both endpoints, this extension is globally Lipschitz.
For $J\in\Z$, let $V^J$ solve
$\dot V_j^J(s)=(\mathcal AV^J(s))_j+g(V_j^J(s))$ for $j\in\Z$ and $s\ge0$,
with initial datum $V_j^J(0)=\delta\mathbf1_{\{j\le J\}}$.
Let $L_g$ be a Lipschitz constant of $g$ on $[0,1]$.
The operator $\mathcal A$ is bounded on $\ell^\infty(\Z)$, and the
extension of $g$ is globally Lipschitz. Picard iteration for the
coordinate integral equation therefore gives a unique global solution
$V^J\in C^1([0,\infty);\ell^\infty(\Z))$.
Since the off-diagonal coefficients of $\mathcal A$ are
nonnegative and $g(s)+L_gs$ is nondecreasing on $[0,1]$,
the integrating-factor comparison in Lemma~\ref{lem:comparison}
applies. The function $(s,j)\mapsto z(s)$ solves the lattice
equation with reaction $g$, because $\mathcal A\mathbf1=0$ and
$z'=g(z)$. Its initial value is $\delta$ at every $j\in\Z$.
Comparison with the zero solution and with the spatially constant
solution $(s,j)\mapsto z(s)$ therefore gives
$0\le V_j^J(s)\le z(s)\le1$ for $s\ge0$ and $j\in\Z$. Consequently,
$h_j(s)=z(s)-V_j^J(s)$, for $s\ge0$ and $j\in\Z$, satisfies
\[
 \partial_s h_j\le(\mathcal Ah)_j+L_gh_j,
 \qquad h_j(0)=\delta\mathbf1_{\{j>J\}}
 \qquad j\in\Z,\ s\ge0.
\]
Linear comparison with the positive semigroup $e^{s\mathcal A}$
therefore gives
\[
 0\le z(\tau)-V_j^J(\tau)
 \le\delta e^{L_g\tau}
       (e^{\tau\mathcal A}\mathbf1_{\{\cdot>J\}})_j
 =\delta e^{L_g\tau}\mathbb P_j^{\mathrm d}(X_\tau>J)
 \qquad j\in\Z.
\]
For each integer $R\ge1$ and $j\le J-R$,
$\mathbb P_j^{\mathrm d}(X_\tau>J)$ is bounded by the probability
that $X$ makes at least $R$ jumps in $[0,\tau]$.
That number has the Poisson distribution with mean $D(q+1)\tau$;
see \cite[Section~2.4]{Norris}. Choose an integer $R\ge1$ such that
\[
 e^{-D(q+1)\tau}\sum_{n=R}^{\infty}
       \frac{(D(q+1)\tau)^n}{n!}
 <\frac{z(\tau)-\theta}{\delta e^{L_g\tau}}.
\]
The left side tends to zero as $R\to\infty$, and the right side
is positive by $z(\tau)>\theta$. For $t\ge T$, take $J=\lfloor m_b(t)\rfloor$.
Then $w_j(t)\ge V_j^J(0)$ for every $j\in\Z$: for $j\le J$
this follows from $w_j(t)\ge\delta$, and for $j>J$ from
$w_j(t)\ge0$. Since $f_i(r)\ge g(r)$ for $0\le i<P$ and
$r\in[0,1]$, Lemma~\ref{lem:comparison} gives
$w_j(t+\tau)\ge V_j^J(\tau)$ for every $j\in\Z$. Hence
\[
 w_j(t+\tau)\ge\theta
 \quad\text{ for }j\in\Z,\ j\le J-R,\ t\ge T.
\]
Finally,
\[
 |m_b(t+\tau)-m_b(t)|
 \le |c|\tau+\kappa\log(1+\tau)
 \qquad t\ge0.
\]
Thus one may take
\[
 T_\theta=T+\tau,\qquad
 C_\theta=R+1+|c|\tau+\kappa\log(1+\tau).
\]
Given a final time $s\ge T_\theta$, put $t=s-\tau\ge T$ and
$J=\lfloor m_b(t)\rfloor$. If $j\le m_b(s)-C_\theta$, then
\[
 j\le m_b(t)-R-1\le\lfloor m_b(t)\rfloor-R=J-R.
\]
Since $j\le J-R$ gives $w_j(s)=w_j(t+\tau)\ge\theta$,
\eqref{eq:sharp-step-lower} follows, with constants independent of $b$.\qedhere

\end{proof}

\subsection{A half-line equation}

\begin{lemma}\label{lem:periodic-halfline-log}

Assume \eqref{eq:hyp} and $\mu>0$. Let
$U\in C^1([0,\infty);\ell^\infty(\N_0))$ satisfy the interior equation
in \eqref{eq:radial}, namely
\[
 \dot U_r(t)=D\bigl(U_{r-1}(t)+qU_{r+1}(t)-(q+1)U_r(t)\bigr)
              +f_r(U_r(t)),\quad\text{for }r\in\N,\ t\ge0,
\]
and
$\dot U_0(t)=b_0(U_1(t)-U_0(t))+f_0(U_0(t))$ for $t\ge0$,
where $b_0>0$.
Suppose $0\le U_r(t)\le1$ for $r\in\N_0$ and $t\ge0$, the initial
datum $(U_r(0))_{r\in\N_0}$ is nonzero and supported at finitely
many indices, and $U_0(t)\to1$ as $t\to\infty$.
Set
\[
 m(t)=c_+t-\frac{3}{2\lambda_+}\log(1+t),\quad\text{for }t\ge0.
\]
For every $\theta\in(0,1)$ there are $C_\theta,T_\theta\ge0$,
both finite, such that
\begin{equation}
 \label{eq:periodic-halfline-log}
 \{r\in\N_0:r\le m(t)-C_\theta\}
 {\subseteq\{r\in\N_0:U_r(t)\ge\theta\}}
 \subseteq\{r\in\N_0:r\le m(t)+C_\theta\}
                  \quad\text{ for } t\ge T_\theta.
\end{equation}

\end{lemma}
\begin{proof}

\emph{An upper reaction family.}
Fix $\theta\in(0,1)$. We construct a family $(F_i)_{i\in\Z}$ of functions $F_i:[0,1]\to\R$,
$P$-periodic in $i$, with $F_i'(0)=a_i$, $f_i\le F_i$ on $[0,1]$, and
$F_i(s)/s$ nonincreasing for $0<s\le1$. Choose
\[
 0<\varepsilon\le\min\left\{\frac12,\frac{a_{\min}}{2L}\right\}
\]
and, for $i\in\{0,\ldots,P-1\}$, define
$H_0:[0,\infty)\to[0,1]$ and $F_i:[0,1]\to\R$ by
\[
 H_0(z)=
 \begin{cases}
  1-(1-z)^3,&0\le z\le1,\\
  1,&z\ge1,
 \end{cases}
 \qquad
 F_i(s)=a_i s H_0((1-s)/\varepsilon).
\]
For $0\le z\le1$, one has $H_0'(z)=3(1-z)^2\ge0$ and
$H_0(z)-z=z(1-z)(2-z)\ge0$. Since the first two derivatives of
$H_0$ match at $z=1$, each $F_i$ belongs to $C^2([0,1])$.
Moreover, $F_i(0)=F_i(1)=0$ and $0<F_i(s)\le a_i s$ for $0<s<1$.
For $0\le s\le1-\varepsilon$, the definition of $F_i$ gives
$F_i(s)=a_i s\ge f_i(s)$ by \eqref{eq:hyp}.
For $1-\varepsilon\le s\le1$, we have $s\ge1/2$ and
$H_0((1-s)/\varepsilon)\ge(1-s)/\varepsilon$. Also,
$f_i(1)=0$ and $|f_i'|\le L$ in \eqref{eq:taylor} give
$f_i(s)\le L(1-s)$. Hence
\[
 F_i(s)\ge \frac{a_{\min}}{2\varepsilon}(1-s)
          \ge L(1-s)\ge f_i(s).
\]
Moreover,
\[
 F_i'(0)=a_i,
 \qquad F_i'(1)=-\frac{3a_i}{\varepsilon}<0
 \qquad i\in\{0,\ldots,P-1\},
\]
and $F_i(s)/s=a_iH_0((1-s)/\varepsilon)$ is nonincreasing on
$(0,1]$. Extend $F_i$ by $F_{i+P}=F_i$ for $i\in\Z$.
The family $(F_i)_{i\in\Z}$ satisfies \eqref{eq:hyp}, and
$f_i(s)\le F_i(s)$ for $i\in\Z$ and $s\in[0,1]$.

\emph{An upper comparison satisfying the half-line boundary inequality.}
Let $W\in C^1([0,\infty);\ell^\infty(\Z))$ solve
\eqref{eq:lattice} with reactions $F_i$ and initial datum
$W_j(0)=\mathbf1_{\{j\le0\}}$ for $j\in\Z$.
Lemma~\ref{lem:comparison} gives $0\le W_j(t)\le1$ for
$j\in\Z$ and $t\ge0$. The matrix $M_\lambda$ in
\eqref{eq:matrix} depends on the reactions through their derivatives
at zero. Since $F_i'(0)=a_i$, the families $(F_i)$ and $(f_i)$
have the same $M_\lambda$, $c_+$ and $\lambda_+$.
By Lemma~\ref{lem:matrix}, $\mu>0$ gives $c_+>0$, so
$m(t)\to\infty$ as $t\to\infty$. For each $\alpha\in(0,1)$,
Lemma~\ref{lem:sharp-step-lower} therefore gives $W_0(t)\ge\alpha$
for all sufficiently large $t$. In view of $W_0\le1$, we obtain
$W_0(t)\to1$ as $t\to\infty$.

We choose a positive multiple $R$ of $P$ such that
$W_{-R}(t)\ge1/2$ for every $t\ge0$ and
$U_r(0)=0$ for every $r\in\N_0$ with $r>R$.
Indeed, $F_{j+P}=F_j$ and $W_{j+P}(0)\le W_j(0)$ give
$W_{j+P}(t)\le W_j(t)$ for $j\in\Z$ and $t\ge0$ by
Lemma~\ref{lem:comparison}. Thus $W_{-R}(t)\ge W_0(t)$
whenever $R$ is a positive multiple of $P$.
Since $W_0(t)\to1$, choose $T>0$ such that
$W_0(t)\ge1/2$ for $t\ge T$.
For $0\le t\le T$, the inequalities $F_j(W_j(t))\ge0$
for $j\in\Z$ and the identity $\mathcal A\mathbf1=0$ give
\[
 \partial_t(1-W)\le\mathcal A(1-W),
 \qquad (1-W_j(0))=\mathbf1_{\{j>0\}}.
\]
Positivity of the diffusion semigroup consequently yields
\[
 1-W_{-R}(t)
 \le\bigl(e^{t\mathcal A}\mathbf1_{\{j>0\}}\bigr)_{-R}
 =\mathbb P_{-R}^{\mathrm d}(X_t>0)
 \qquad 0\le t\le T.
\]
Under $\mathbb P_{-R}^{\mathrm d}$, a path of $X$ from $-R$
to a positive index makes at least $R+1$ jumps. Its jump count by
time $T$ has the Poisson distribution with mean $D(q+1)T$.
Consequently, uniformly for $0\le t\le T$,
$\mathbb P_{-R}^{\mathrm d}(X_t>0)\le\mathbb P(N_T\ge R+1)$,
where $N_T$ has this Poisson distribution.
Choose a positive multiple $R$ of $P$ with
$R\ge\max\{r\in\N_0:U_r(0)>0\}$ and
\[
 e^{-D(q+1)T}\sum_{n=R+1}^{\infty}
       \frac{(D(q+1)T)^n}{n!}\le\frac12.
\]
The set $\{r\in\N_0:U_r(0)>0\}$ is nonempty and finite by
hypothesis, and the Poisson tail in the preceding display tends to
zero as $R\to\infty$. The chosen $R$ therefore gives
$W_{-R}(t)\ge1/2$ for $0\le t\le T$; for $t\ge T$, use
$W_{-R}(t)\ge W_0(t)\ge1/2$.

Set
\[
 \overline U_r(t)=\min\{1,2W_{r-R}(t)\},
 \quad\text{for }r\in\N_0,\ t\ge0.
\]
Since $R$ is a multiple of $P$, $F_{r-R}=F_r$ for every
$r\in\N_0$. Fix $r\in\N$ and a time $t\ge0$ at which
$\overline U_r(t)<1$. Then $\overline U_r(t)=2W_{r-R}(t)$,
and the monotonicity of $F_r(s)/s$ gives
\[
 f_r(\overline U_r)
 \le F_r(2W_{r-R})
 \le2F_r(W_{r-R}).
\]
For the sequence $\overline U$ on $\N_0$, write
$(\mathcal A\overline U)_r
=D(\overline U_{r-1}+q\overline U_{r+1}-(q+1)\overline U_r)$
only at interior indices $r\ge1$. Since
$\overline U_{r\pm1}\le2W_{r-R\pm1}$ and
$\overline U_r=2W_{r-R}$, the diffusion terms satisfy
$(\mathcal A\overline U)_r\le2(\mathcal AW)_{r-R}$.
Using the diffusion inequality together with
$f_r(2W_{r-R})\le2F_r(W_{r-R})$ and the equation solved by $W$, we
obtain
\[
 \partial_t\overline U_r
 =2(\mathcal AW)_{r-R}+2F_r(W_{r-R})
 \ge(\mathcal A\overline U)_r+f_r(\overline U_r).
\]
For each $r\in\N$, the derivative
$\partial_t\overline U_r(t)$ is zero for almost every time at which
$\overline U_r(t)=1$. At such times,
$(\mathcal A\overline U)_r\le0$ because
$\overline U_{r-1},\overline U_{r+1}\le1$, and $f_r(1)=0$.
At $r=0$, $W_{-R}(t)\ge1/2$ gives $\overline U_0(t)=1$
for all $t\ge0$, and therefore
\[
 b_0(\overline U_1-\overline U_0)+f_0(\overline U_0)\le0
 \qquad t\ge0.
\]
The solution $W$ belongs to $C^1([0,\infty);\ell^\infty(\Z))$,
and the map $s\mapsto\min\{1,2s\}$ is Lipschitz. Hence
$t\mapsto(\overline U_r(t))_{r\in\N_0}$ is locally Lipschitz in
$\ell^\infty(\N_0)$. For every $r\in\N$, the preceding
interior inequalities hold for almost every $t\ge0$; the boundary
inequality holds for every $t\ge0$.
For $0\le r\le R$, $W_{r-R}(0)=1$, so
$U_r(0)\le1=\overline U_r(0)$. For $r>R$,
$U_r(0)=0\le\overline U_r(0)$. Lemma~\ref{lem:comparison},
applied on each finite time interval with spatial domain $\N_0$,
therefore gives $U_r(t)\le\overline U_r(t)$ for
$r\in\N_0$ and $t\ge0$.
By \eqref{eq:sharp-step-upper} with reactions $(F_i)$ and step
location $b=0$,
$W_{r-R}(t)\le C_F(1+y)e^{-\lambda_+y}$ whenever
$t\ge0$ and $y=r-R-m(t)>0$, where $C_F$ depends on $(F_i)$.
Choose $C>0$ such that $2C_F(1+y)e^{-\lambda_+y}<\theta$ for $y>C$.
Then $U_r(t)<\theta$ for $t\ge0$ and $r\in\N_0$ with
$r>m(t)+R+C$. Thus $U_r(t)\ge\theta$ implies
$r\le m(t)+R+C$, proving the second inclusion in
\eqref{eq:periodic-halfline-log} with $C_\theta\ge R+C$.

\emph{A lower comparison with prescribed boundary values.}
Choose $\beta\in(\theta,1)$. Since $U_0(t)\to1$ as
$t\to\infty$, choose $T_0\ge0$ such that
\[
 U_0(T_0+t)\ge\beta\quad\text{ for } t\ge0.
\]
For $i\in\Z$ and $0\le s\le\beta$, define
$f_{i,\beta}(s)=f_i(s)(1-s/\beta)$. Then
\[
 f_{i,\beta}'(0)=a_i,
 \qquad
 f_{i,\beta}'(\beta)=-\frac{f_i(\beta)}{\beta}<0,
 \qquad
 0<f_{i,\beta}(s)\le f_i(s)\le a_i s
 \quad\text{ for } 0<s<\beta.
\]
For $0\le z\le1$, define
\[
 \widetilde f_{i,\beta}(z)
 =\beta^{-1}f_{i,\beta}(\beta z)
 =\beta^{-1}f_i(\beta z)(1-z)
 \qquad i\in\Z.
\]
For every $i\in\Z$, the identities
$f_{i,\beta}'(0)=a_i$ and
$f_{i,\beta}'(\beta)=-f_i(\beta)/\beta$ give
\[
 \widetilde f_{i,\beta}'(0)=a_i,\qquad
 \widetilde f_{i,\beta}'(1)=-f_i(\beta)/\beta<0,\qquad
 0<\widetilde f_{i,\beta}(z)\le a_i z\quad\text{ for } 0<z<1
 \qquad i\in\Z.
\]
The functions $\widetilde f_{i,\beta}$ are $P$-periodic in $i$,
belong to $C^2([0,1])$, and vanish at zero and one.
Together with the preceding derivative and positivity bounds, these
properties verify \eqref{eq:hyp} for $(\widetilde f_{i,\beta})$.
Since $\widetilde f_{i,\beta}'(0)=a_i$, this family has matrix
$M_\lambda$ in \eqref{eq:matrix} and the same principal-eigenvalue
function $k(\lambda)$ and parameters $c_+,\lambda_+$ as $(f_i)$.

Let $V\in C^1([0,\infty);\ell^\infty(\Z))$ solve
\[
 \dot V_j(t)=(\mathcal AV(t))_j+f_{j,\beta}(V_j(t)),
 \quad\text{for }j\in\Z,\ t\ge0,
 \qquad V_j(0)=\beta\mathbf1_{\{j\le0\}}.
\]
Indeed, $V/\beta$ solves
\eqref{eq:lattice} with reactions $\widetilde f_{i,\beta}$ and
initial datum $\mathbf1_{\{j\le0\}}$.
Lemma~\ref{lem:comparison} gives existence of this normalized
solution and hence of $V$, with $0\le V_j(t)\le\beta$ for
$j\in\Z$ and $t\ge0$.
For every $r\in\N$, $V_r(0)=0\le U_r(T_0)$, and
$V_0(t)\le\beta\le U_0(T_0+t)$ for $t\ge0$.
For $r\in\N$ and $t\ge0$, $0\le V_r(t)\le\beta$ gives
\[
 \dot V_r-(\mathcal AV)_r-f_r(V_r)
 =f_{r,\beta}(V_r)-f_r(V_r)\le0.
\]
Apply Lemma~\ref{lem:comparison} on each finite time interval,
with spatial domain $S=\N_0$, subset $F=\N$, and the half-line
diffusion operator whose boundary coefficient is $b_0$.
The restrictions $(V_r(t))_{r\in\N_0}$ and
$(U_r(T_0+t))_{r\in\N_0}$ are locally Lipschitz trajectories in
$\ell^\infty(\N_0)$. On $F$, $V$ satisfies the preceding
subsolution inequality and $U(T_0+\cdot)$ satisfies the interior
equation. On $S\setminus F=\{0\}$,
$V_0(t)\le\beta\le U_0(T_0+t)$ for $t\ge0$.
Together with $V_r(0)\le U_r(T_0)$ for $r\in F$, these
inequalities yield
\[
 V_r(t)\le U_r(T_0+t),\quad\text{for }r\in\N,\ t\ge0.
\]
The same inequality holds at $r=0$ by the boundary order.
Since $0<\theta/\beta<1$, \eqref{eq:sharp-step-lower}, applied to
$V/\beta$ at level $\theta/\beta$ with step location $b=0$,
gives finite constants $C',T'\ge0$ such that
\[
 U_r(T_0+t)\ge V_r(t)\ge\theta
 \quad\text{ for }t\ge T'\text{ and }0\le r\le m(t)-C'
 \qquad r\in\N_0.
\]
For $t\ge0$,
\[
 \left|m(t+T_0)-m(t)\right|
 \le c_+T_0+\frac{3}{2\lambda_+}\log(1+T_0).
\]
Choose
$C_\theta\ge\max\{R+C,C'+c_+T_0+(3/(2\lambda_+))\log(1+T_0)\}$
and $T_\theta\ge T_0+T'$.
For $s\ge T_\theta$, put $t=s-T_0\ge T'$.
If $r\in\N_0$ and $r\le m(s)-C_\theta$, the displayed bound
on $m(t+T_0)-m(t)$ gives $r\le m(t)-C'$.
Thus $U_r(s)\ge V_r(t)\ge\theta$ by the preceding lower bound.
This proves the first inclusion in \eqref{eq:periodic-halfline-log};
the choice $C_\theta\ge R+C$ also preserves the second inclusion.\qedhere

\end{proof}
\section{Tree spectra and fronts}\label{sec:tree-spectra}

Recall $\lambda_0=\tfrac12\log q$ and $\mu=k(\lambda_0)$ from
\eqref{eq:parameters}. When $a_i=a$ for $0\le i<P$,
$\mu=a-\rho_q$, where $\rho_q=D(\sqrt q-1)^2$.
We compute the spectral upper endpoint for the potentials
$a(x)=a_{b(x)}$ and, under $\mu\le0$, $a(x)=a_{r(x)}$.
When $\mu>0$, we construct finite connected sets and principal
eigenvectors satisfying Lemma~\ref{lem:local-general}, with
$R_0,\eta_0,m_0$ independent of the vertex. For $f(x,s)=f_{b(x)}(s)$ or $f(x,s)=f_{r(x)}(s)$,
if $\mu>0$, Lemma~\ref{lem:local-general} gives $u(t,x)\to1$
as $t\to\infty$ at every fixed vertex for the solution of
\eqref{eq:tree-cauchy} with $0\le u_0\le1$ and $u_0\not\equiv0$.
Every finite restriction below acts on functions extended by zero to
the whole tree, so its diffusion diagonal is $-D(q+1)$ at every
vertex of the finite set.

\begin{lemma}\label{lem:horo-spectrum}
Assume \eqref{eq:hyp}. For the potential $a(x)=a_{b(x)}$, let
$\mathcal L=D\Delta_{\Tq}+a(x)$ on $\ell^2(\Tq)$. Then
\begin{equation}
\label{eq:spectral-top}
 \sup\spec(\mathcal L)=\mu=k(\lambda_0).
\end{equation}
The eigenspace $\ker(\mu-\mathcal L)$ is $\{0\}$. If $\mu>0$, the
hypotheses of Lemma~\ref{lem:local-general} hold with constants
independent of the vertex.
\end{lemma}
\begin{proof}

\emph{The upper spectral bound.}
The degree bound and the bounded real potential make $\mathcal L$
a bounded self-adjoint operator on $\ell^2(\Tq)$; indeed,
$\|\mathcal Lz\|_2\le(2D(q+1)+a_{\max})\|z\|_2$.
Put $p_j=\phi_j(\lambda_0)$ for $j\in\Z$ and
$H(x)=q^{-b(x)/2}p_{b(x)}>0$. The formula defining $\mathcal L$ also acts pointwise on $H$,
since each vertex has finitely many neighbors. At a vertex with $b(x)=j$, the
parent has coordinate $j-1$ and the $q$ children have coordinate
$j+1$. By \eqref{eq:matrix} and $e^{\lambda_0}=\sqrt q$,
\[
 (\mathcal LH)(x)
 =q^{-j/2}\bigl[D\sqrt q\,(p_{j-1}+p_{j+1})
                +(a_j-D(q+1))p_j\bigr]
 =\mu H(x)
 \qquad x\in\Tq,\ j=b(x).
\]

Since $p$ is positive and $P$-periodic, the ratios $H(y)/H(x)$
for $x\sim y$ have positive upper and lower bounds.
Fix a finitely supported function $z:\Tq\to\C$.
Expanding one square for each unordered edge gives
\eqref{eq:ground-state}. By $\mathcal LH=\mu H$, the diagonal
contribution at $x$ is
\[
 D\frac{|z(x)|^2}{H(x)}\sum_{y\sim x}H(y)
   =(\mu+D(q+1)-a(x))|z(x)|^2
 \qquad x\in\Tq.
\]
The mixed contribution is
$-2D\re\sum_{\{x,y\}\in E}\overline{z(x)}z(y)$, with each edge
counted once. Combining these terms gives
\begin{equation}
\label{eq:ground-state}
 \langle z,(\mu-\mathcal L)z\rangle
 =D\sum_{\{x,y\}\in E}H(x)H(y)
       \left|\frac{z(x)}{H(x)}-\frac{z(y)}{H(y)}\right|^2.
\end{equation}
To extend \eqref{eq:ground-state} to every $z\in\ell^2(\Tq)$,
choose an orientation for each edge and define
\[
 (\mathcal B_Hz)_{\{x,y\}}
 =\sqrt{DH(x)H(y)}
       \left(\frac{z(x)}{H(x)}-\frac{z(y)}{H(y)}\right)
 \qquad \{x,y\}\in E.
\]
The adjacent-ratio bounds and the degree bound give
\[
 \|\mathcal B_Hz\|_{\ell^2(E)}^2
 \le 2D\sum_{\{x,y\}\in E}
 \left(\frac{H(y)}{H(x)}|z(x)|^2+
       \frac{H(x)}{H(y)}|z(y)|^2\right)
 \le C\|z\|_2^2.
\]
Thus $\mathcal B_H:\ell^2(\Tq)\to\ell^2(E)$ is bounded. For
$z\in\ell^2(\Tq)$, choose finitely supported $z_n$ with
$\|z_n-z\|_2\to0$ as $n\to\infty$. Then
$\mathcal B_Hz_n\to\mathcal B_Hz$ in $\ell^2(E)$, and boundedness of
$\mathcal L$ permits passage to the limit on both sides of
\eqref{eq:ground-state}. Thus the edge sum converges and the identity
holds for every $z\in\ell^2(\Tq)$. Its nonnegative right side gives
$\langle z,\mathcal Lz\rangle\le\mu\norm{z}_2^2$.
The variational characterization of the spectral upper endpoint
therefore yields $\sup\spec(\mathcal L)\le\mu$; see
\cite[Theorem~2.19]{Teschl}.
Identity \eqref{eq:ground-state} is the ground-state transform;
its graph formulation is also treated in \cite{KellerPinchoverPogorzelski}.

\emph{Cone test functions.}
For $x\in\Tq$ and $N\in\N_0$, define
\[
 \mathcal C_N(x)=\{y\in\Tq:d(x,y)\le N,
                         \ b(y)=b(x)+d(x,y)\}.
\]
Each path from $x$ to a vertex of $\mathcal C_N(x)$ follows children,
so this set consists of $x$ and its descendants through depth $N$.
Put $s=b(x)$ and $z_N=H\ind_{\mathcal C_N(x)}$. At depth $n$ there
are $q^n$ vertices, each with squared value $q^{-s-n}p_{s+n}^2$.
Hence
\begin{equation}\label{eq:horo-cone-norm}
 \norm{z_N}_2^2=q^{-s}\sum_{n=0}^Np_{s+n}^2,
 \qquad
 c_0q^{-s}(N+1)\le\norm{z_N}_2^2\le C_0q^{-s}(N+1)
 \qquad x\in\Tq,\ N\in\N_0,\ s=b(x),
\end{equation}
where $c_0=\min_{0\le i<P}p_i^2>0$ and
$C_0=\max_{0\le i<P}p_i^2$. In \eqref{eq:ground-state}, the ratio
$z_N/H$ equals one on $\mathcal C_N(x)$ and zero elsewhere. The
contributing edges are therefore the single parent edge at $x$ and
the $q^{N+1}$ edges from depth $N$ to depth $N+1$. Their sum is
\begin{equation}\label{eq:horo-cone-boundary}
 \langle z_N,(\mu-\mathcal L)z_N\rangle
 =D\sqrt q\,q^{-s}
       \bigl(p_{s-1}p_s+p_{s+N}p_{s+N+1}\bigr)
 \qquad x\in\Tq,\ N\in\N_0,\ s=b(x).
\end{equation}
By \eqref{eq:horo-cone-norm} and \eqref{eq:horo-cone-boundary},
each boundary product is at most $C_0$, whereas
$\norm{z_N}_2^2\ge c_0q^{-s}(N+1)$. Hence
\begin{equation}
\label{eq:cone-rayleigh}
 \frac{\langle z_N,\mathcal Lz_N\rangle}{\norm{z_N}_2^2}
 \ge\mu-\frac{C}{N+1}
 \qquad x\in\Tq,\ N\in\N_0,
\end{equation}
with $C=2D\sqrt q\,C_0/c_0$, independently of $x$ and $N$.
The Rayleigh characterization gives
$\sup\spec(\mathcal L)\ge\mu-C/(N+1)$ for every $N$.
Letting $N\to\infty$ and combining this with
$\sup\spec(\mathcal L)\le\mu$ proves \eqref{eq:spectral-top}.

\emph{The eigenspace at $\mu$.}
Suppose that $\mathcal Lz=\mu z$ for some $z\in\ell^2(\Tq)$.
By \eqref{eq:ground-state}, the sum of nonnegative edge terms is
zero. Thus $z(x)/H(x)=z(y)/H(y)$ whenever $x\sim y$, and
connectedness gives $z=cH$ for a constant $c\in\C$.
For any fixed vertex $x$, the generation count above yields
\[
 \sum_{y\in\mathcal C_N(x)}H(y)^2
   =q^{-b(x)}\sum_{n=0}^Np_{b(x)+n}^2
   \ge c_0q^{-b(x)}(N+1)
 \qquad N\in\N_0.
\]
As $N\to\infty$, this forces $c=0$ whenever $cH\in\ell^2(\Tq)$.
Hence $z=0$, and therefore $\ker(\mu-\mathcal L)=\{0\}$.

\emph{Uniform finite-domain eigenvectors.}
Assume $\mu>0$. Let $\eta_N(x)$ be the largest eigenvalue of the
zero-exterior restriction of $\mathcal L$ to $\mathcal C_N(x)$.
The finite-dimensional Rayleigh characterization and
\eqref{eq:cone-rayleigh} give
$\eta_N(x)\ge\mu-C/(N+1)$ for every $x$. Fix $N$ with
$C/(N+1)\le\mu/2$. Adding $D(q+1)$ times the identity makes
each restricted matrix nonnegative with diagonal $a_{b(y)}>0$
at a vertex $y$, and irreducible because $\mathcal C_N(x)$ is connected.
The Perron--Frobenius theorem, followed by subtraction of the scalar
shift, shows that $\eta_N(x)\ge\mu/2$ is simple and has a strictly
positive eigenvector $\varphi_x$. Normalize it by
$\max_{y\in\mathcal C_N(x)}\varphi_x(y)=1$.
Fix the integer $N$ chosen above. If $b(x)-b(y)\in P\Z$,
a rooted isomorphism from $\mathcal C_N(x)$ to $\mathcal C_N(y)$
preserves every potential value and every diffusion coefficient.
It consequently identifies the two restricted matrices. Uniqueness of the normalized positive eigenvector then
shows that the root values $\varphi_x(x)$ range over at most $P$
positive numbers. Set
$m_0=\min_{x\in\Tq}\varphi_x(x)>0$.
Since $x\in\mathcal C_N(x)\subseteq B(x,N)$, the choices
$R_0=N$, $\eta_0=\mu/2$ and $\Omega_x=\mathcal C_N(x)$ verify the
hypotheses of Lemma~\ref{lem:local-general}.

\end{proof}

\begin{lemma}\label{lem:radial-local}
Assume \eqref{eq:hyp} and $f(x,s)=f_{r(x)}(s)$. If $\mu>0$, the hypotheses of
Lemma~\ref{lem:local-general} hold with constants independent of the
vertex. Consequently, every solution $u$ of \eqref{eq:tree-cauchy}
with $0\le u_0\le1$ and $u_0\not\equiv0$ converges to one as
$t\to\infty$ at each fixed vertex.
\end{lemma}
\begin{proof}

Fix $N\in\N$. For $x\ne o$, define
$\mathcal C_N(x)=\{y\in\Tq:d(x,y)\le N,
 r(y)=r(x)+d(x,y)\}$. At $x=o$, choose $q$ distinct neighbors
$y_1,\ldots,y_q$ of $o$ and define
\[
 \mathcal C_N(o)=\{o\}\cup\bigcup_{i=1}^{q}
 \{y\in\Tq:1\le r(y)\le N,\ d(y_i,y)=r(y)-1\}.
\]
Extend functions on
these cones by zero outside them. Every restricted diffusion matrix
then has diagonal $-D(q+1)$, including at the cone root. Put $s=r(x)$;
the potential at depth $n$ is $a_{s+n}$.

Write $p_j=\phi_j(\lambda_0)$ for $j\in\Z$. Define
$z_N(y)=q^{-n/2}p_{s+n}$ when $y\in\mathcal C_N(x)$ has depth $n$
from $x$, and set $z_N=0$ elsewhere. The $q^n$ vertices at depth $n$
give
\[
 S_N(s):=\|z_N\|_2^2=\sum_{n=0}^Np_{s+n}^2
 \qquad s=r(x),\ x\in\Tq,\ N\in\N.
\]
Between depths $n$ and $n+1$ there are $q^{n+1}$ edges.
Each such edge contributes twice to the quadratic form, so their
total contribution is $2D\sqrt q\,p_{s+n}p_{s+n+1}$.
The Rayleigh numerator is therefore
\[
 Q_N(s)=
 \sum_{n=0}^N(a_{s+n}-D(q+1))p_{s+n}^2
 +2D\sqrt q\sum_{n=0}^{N-1}p_{s+n}p_{s+n+1}
 \qquad s=r(x),\ x\in\Tq,\ N\in\N.
\]
At $\lambda_0=\frac12\log q$, the periodic eigenvector equation reads
\[
 D\sqrt q\,(p_{j-1}+p_{j+1})
 +(a_j-D(q+1))p_j=\mu p_j
 \qquad j\in\Z.
\]
Multiplying by $p_j$ and summing from $j=s$ to $s+N$ counts each
interior product twice and leaves the two endpoint products:
\begin{equation}
\label{eq:radial-cone-rayleigh}
 Q_N(s)=\mu S_N(s)
 -D\sqrt q\bigl(p_{s-1}p_s+p_{s+N}p_{s+N+1}\bigr)
 \qquad s=r(x),\ x\in\Tq,\ N\in\N.
\end{equation}
By positivity and periodicity of $p$,
$S_N(s)\ge(N+1)(\min_{0\le i<P}p_i)^2$, and each boundary product in
\eqref{eq:radial-cone-rayleigh} is at most $(\max_{0\le i<P}p_i)^2$.
It follows that
\[
 \frac{Q_N(s)}{S_N(s)}\ge\mu-\frac{C}{N+1},
 \qquad C=2D\sqrt q
 \left(\frac{\max_{0\le i<P}p_i}{\min_{0\le i<P}p_i}\right)^2
 \qquad s=r(x),\ x\in\Tq,\ N\in\N,
\]
independently of $s$ and $N$.

Choose $N$ with $C/(N+1)\le\mu/2$. By the Rayleigh characterization,
the largest eigenvalue of every zero-exterior cone operator is at
least $\mu/2$. Adding $D(q+1)$ times the identity makes each
such matrix nonnegative with diagonal $a_{r(y)}>0$ at a vertex $y$,
and irreducible because the cone is connected. The Perron--Frobenius theorem, followed by
subtraction of the scalar shift, gives a strictly positive principal
eigenvector $\varphi_x$, unique under the normalization
$\max_{y\in\mathcal C_N(x)}\varphi_x(y)=1$.
For this fixed $N$, the cone operators are determined, up to a
root-preserving isomorphism, by $r(x)$ modulo $P$. At $o$, the chosen
$q$ branches and the retained diagonal $-D(q+1)$ give the same matrix
as any other cone whose root has phase zero. The normalized root
values therefore range over at most $P$ positive numbers. Set
$m_0=\min_{x\in\Tq}\varphi_x(x)>0$.
Taking $\Omega_x=\mathcal C_N(x)\subseteq B(x,N)$, $R_0=N$ and
$\eta_0=\mu/2$ verifies the hypotheses of
Lemma~\ref{lem:local-general}, which gives the stated convergence.

\end{proof}

\begin{lemma}\label{lem:decreasing-ground}
Assume \eqref{eq:hyp} and $\mu\le0$. Define
$h_j=q^{-j/2}\phi_j(\lambda_0)$ for $j\in\Z$. Then
$h_{j+1}<h_j$ for every $j\in\Z$.
\end{lemma}
\begin{proof}

Put $d_j=h_j-h_{j+1}$. Since $\phi(\lambda_0)$ is $P$-periodic,
\[
 h_{j-P}=q^{P/2}h_j,
 \qquad d_{j-P}=q^{P/2}d_j
 \qquad j\in\Z.
\]
The eigenvector equation for $M_{\lambda_0}$, written in terms of $h$,
is
\[
 D(h_{j-1}+qh_{j+1}-(q+1)h_j)+a_jh_j=\mu h_j
 \qquad j\in\Z.
\]
Since $h_{j-1}-h_j=d_{j-1}$ and $h_j-h_{j+1}=d_j$, this equation
becomes $D(d_{j-1}-qd_j)+(a_j-\mu)h_j=0$. Hence
\begin{equation}\label{eq:decreasing-ground-recursion}
 d_j=q^{-1}d_{j-1}+\frac{a_j-\mu}{Dq}h_j
 \qquad j\in\Z.
\end{equation}
Iterating \eqref{eq:decreasing-ground-recursion} through one period gives
\[
 d_j=q^{-P}d_{j-P}
 +\frac1D\sum_{m=0}^{P-1}q^{-m-1}(a_{j-m}-\mu)h_{j-m}
 \qquad j\in\Z.
\]
Using $d_{j-P}=q^{P/2}d_j$, we obtain
\begin{equation}\label{eq:decreasing-ground-period}
 (1-q^{-P/2})d_j
 =\frac1D\sum_{m=0}^{P-1}q^{-m-1}(a_{j-m}-\mu)h_{j-m}
 \qquad j\in\Z.
\end{equation}
Every summand on the right of \eqref{eq:decreasing-ground-period}
is strictly positive because $a_i>0$,
$\mu\le0$, and $h_j>0$. Since $q>1$, it follows that $d_j>0$, which is
exactly $h_{j+1}<h_j$.\qedhere

\end{proof}

\begin{lemma}\label{lem:radial-spectrum}
Assume \eqref{eq:hyp} and $\mu\le0$. For the radial potential, let
$\mathcal L_r=D\Delta_{\Tq}+a_{r(x)}$ on $\ell^2(\Tq)$. Then
\begin{equation}
\label{eq:radial-spectral-edge}
 \sup\spec(\mathcal L_r)=\mu,
 \qquad \ker(\mu-\mathcal L_r)=\{0\}.
\end{equation}
\end{lemma}
\begin{proof}

The bounded degree and the bounded real potential make
$\mathcal L_r$ bounded and self-adjoint on $\ell^2(\Tq)$.
Set $h_j=q^{-j/2}\phi_j(\lambda_0)$ for $j\in\Z$ and
$H(x)=h_{r(x)}$. Evaluating the formula for $\mathcal L_r$ pointwise
at $H$ and using the recurrence
\[
 D(h_{j-1}+qh_{j+1}-(q+1)h_j)+a_jh_j=\mu h_j
 \qquad j\in\Z,
\]
we obtain $\mathcal L_rH(x)=\mu H(x)$ for $x\ne o$.
At $o$, the $q+1$ neighbors all have radial index one, so
$\mathcal L_rH(o)=D(q+1)(h_1-h_0)+a_0h_0$.
Subtracting this expression from the recurrence at $j=0$ gives
\begin{equation}\label{eq:radial-root-remainder}
 (\mu-\mathcal L_r)H(o)=D(h_{-1}-h_1)>0,
\end{equation}
where Lemma~\ref{lem:decreasing-ground} gives $h_{-1}>h_0>h_1$.
Put $\kappa_o=D(h_{-1}-h_1)/h_0>0$. For finitely supported
complex-valued $z$, expanding the squares gives
\[
 \begin{aligned}
 &D\sum_{\{x,y\}\in E}H(x)H(y)
       \left|\frac{z(x)}{H(x)}-\frac{z(y)}{H(y)}\right|^2\\
 &=D\sum_{x\in\Tq}\frac{\sum_{y\sim x}H(y)}{H(x)}|z(x)|^2
       -2D\re\sum_{\{x,y\}\in E}\overline{z(x)}z(y).
 \end{aligned}
\]
The recurrence and the root calculation show that the coefficient
$D\sum_{y\sim x}H(y)/H(x)$ equals
$\mu-a_{r(x)}+D(q+1)-\kappa_o\mathbf1_{\{o\}}(x)$.
Substitution into this expansion yields
\begin{equation}
\label{eq:radial-ground-identity}
 \langle z,(\mu-\mathcal L_r)z\rangle
 =D\sum_{\{x,y\}\in E}H(x)H(y)
       \left|\frac{z(x)}{H(x)}-\frac{z(y)}{H(y)}\right|^2
       +\kappa_o|z(o)|^2.
\end{equation}
Positivity and periodicity of $\phi(\lambda_0)$ bound the ratios
$H(y)/H(x)$ above and below by positive constants for $x\sim y$.
Hence
\[
 \begin{aligned}
 D\sum_{\{x,y\}\in E}H(x)H(y)
 \left|\frac{z(x)}{H(x)}-\frac{z(y)}{H(y)}\right|^2
 &\le2D\sum_{\{x,y\}\in E}
 \left[\frac{H(y)}{H(x)}|z(x)|^2+
       \frac{H(x)}{H(y)}|z(y)|^2\right]\\
 &\le C\|z\|_2^2.
 \end{aligned}
 \qquad z\in\ell^2(\Tq).
\]
For $z\in\ell^2(\Tq)$ and $n\in\N$, set
$z_n=z\mathbf1_{B(o,n)}$, so $\|z_n-z\|_2\to0$ as $n\to\infty$.
Choose one orientation for each edge and define
\[
 \xi_n(\{x,y\})=\sqrt{DH(x)H(y)}
       \left(\frac{z_n(x)}{H(x)}-\frac{z_n(y)}{H(y)}\right)
 \qquad \{x,y\}\in E.
\]
The adjacent-ratio bounds and the degree $q+1$ give
$\|\xi_n-\xi_m\|_{\ell^2(E)}^2\le C\|z_n-z_m\|_2^2$.
By completeness, there is $\xi\in\ell^2(E)$ such that
$\xi_n\to\xi$ in $\ell^2(E)$. Coordinatewise convergence of $z_n$
identifies
\[
 \xi(\{x,y\})=\sqrt{DH(x)H(y)}
       \left(\frac{z(x)}{H(x)}-\frac{z(y)}{H(y)}\right)
 \qquad \{x,y\}\in E.
\]
Thus $\|\xi_n\|_{\ell^2(E)}^2\to\|\xi\|_{\ell^2(E)}^2$.
Together with $z_n(o)\to z(o)$ and boundedness of $\mathcal L_r$,
the convergence $\|\xi_n\|_{\ell^2(E)}^2\to\|\xi\|_{\ell^2(E)}^2$
extends \eqref{eq:radial-ground-identity} to all $z\in\ell^2(\Tq)$.
Since the right side of \eqref{eq:radial-ground-identity} is
nonnegative, the Rayleigh characterization gives
$\sup\spec(\mathcal L_r)\le\mu$.

For the reverse inequality, fix $x\ne o$, write $s=r(x)$ and
$p_j=\phi_j(\lambda_0)$. For $N\in\N$, define
$\mathcal C_N(x)=\{y\in\Tq:d(x,y)\le N,
 r(y)=r(x)+d(x,y)\}$. Define $z_N(y)=q^{-n/2}p_{s+n}$ at depth $0\le n\le N$
in this cone, and set $z_N=0$ elsewhere. The $q^n$ vertices at depth
$n$ give $\|z_N\|_2^2=\sum_{n=0}^Np_{s+n}^2$.
There are $q^{n+1}$ internal edges between depths $n$ and $n+1$,
so their quadratic-form contribution is
$2D\sqrt q\,p_{s+n}p_{s+n+1}$. Summing the periodic equation
$D\sqrt q(p_{j-1}+p_{j+1})+(a_j-D(q+1))p_j=\mu p_j$ over
$j=s,\ldots,s+N$ after multiplication by $p_j$ therefore gives
\[
 \frac{\langle z_N,\mathcal L_rz_N\rangle}{\|z_N\|_2^2}
 =\mu-
 \frac{D\sqrt q\,(p_{s-1}p_s+p_{s+N}p_{s+N+1})}
      {\sum_{n=0}^Np_{s+n}^2}
 \ge\mu-\frac{C}{N+1}
 \qquad N\in\N.
\]
Since $p$ has a positive minimum and a finite maximum over one period,
the last inequality holds with $C$ independent of $N$ and of $x$. Letting $N\to\infty$ gives
$\sup\spec(\mathcal L_r)\ge\mu$ and proves the spectral equality in
\eqref{eq:radial-spectral-edge}.

Finally suppose $\mathcal L_rz=\mu z$ for some $z\in\ell^2(\Tq)$.
Then every nonnegative term in \eqref{eq:radial-ground-identity} must
vanish. Since $\kappa_o>0$, we first obtain $z(o)=0$. The edge terms
then imply that $z/H$ is constant along every edge, and connectedness
forces $z\equiv0$. This proves $\ker(\mu-\mathcal L_r)=\{0\}$ and completes
\eqref{eq:radial-spectral-edge}.
\end{proof}

\subsection{Comparison reactions with a common derivative at zero}\label{sec:common-derivative}

For a scalar reaction with derivative $a>0$ at zero, the constants
$c_*(a)$ and $\lambda_*(a)$ are defined in
\eqref{eq:homogeneous-parameters}. We first prove level-position
estimates on $\Z$ and, when $a>\rho_q$, on $\Tq$, with center
$c_*(a)t-(3/(2\lambda_*(a)))\log t$, and then construct scalar
reactions bounding a finite family with $f_i'(0)=a$ for every $i$.
The identity $\log(1+t)-\log t=\log(1+1/t)\to0$ as $t\to\infty$
permits either logarithm in these estimates.

\begin{lemma}\label{lem:graham}

Let $h\in C^2([0,K])$, where $K>0$, satisfy
\[
 h(0)=h(K)=0,\qquad h'(0)=a>0,\qquad h'(K)<0,
 \qquad 0<h(s)\le as\quad\text{ for } 0<s<K.
\]
Let $W:[0,\infty)\times\Z\to[0,K]$ be the solution of
\[
 \dot W_j=D(W_{j-1}+qW_{j+1}-(q+1)W_j)+h(W_j),
 \qquad W_j(0)=K\ind_{\{j\le0\}}.
\]
The differential equation holds for $j\in\Z$ and $t>0$, and the
initial condition holds for every $j\in\Z$.
For each $t\ge0$, the sequence $j\mapsto W_j(t)$ is nonincreasing,
with $W_j(t)\to K$ as $j\to-\infty$ and $W_j(t)\to0$ as
$j\to+\infty$. In particular, for every $\theta\in(0,K)$ the
following maximum is well defined, and
\begin{equation}
\label{eq:hom-lattice-bramson}
 \max\{j:W_j(t)\ge\theta\}
 =c_*(a)t-\frac{3}{2\lambda_*(a)}\log t+O_\theta(1)
 \quad\text{ as } t\to\infty.
\end{equation}
The conclusion holds for every $a>0$, with the sign of $c_*(a)$
determined by $a-\rho_q$.

\end{lemma}
\begin{proof}

Define $\widetilde W_j(t)=W_j(t)/K$ and
$\widetilde h(s)=K^{-1}h(Ks)$ for $0\le s\le1$.
Then $\widetilde W$ solves \eqref{eq:lattice} with
$f_j=\widetilde h$ for every $j\in\Z$ and
$\widetilde W_j(0)=\ind_{\{j\le0\}}$.
By the hypotheses on $h$, the function $\widetilde h$ satisfies
\eqref{eq:hyp}, with $\widetilde h'(0)=a$ and
$\widetilde h'(1)=h'(K)<0$. Thus it suffices to prove the result for
$K=1$, with the level $\theta$ replaced by $\theta/K$.

\emph{Monotonicity and the spatial tails.}
Lemma~\ref{lem:comparison}, with $P=1$, gives $0\le W\le1$.
Translation invariance and the same comparison principle, applied to
$W_{j+1}(0)\le W_j(0)$, give $W_{j+1}(t)\le W_j(t)$.
Let $S_t=e^{t\mathcal A}$ on $\ell^\infty(\Z)$, where
\[
 (\mathcal Av)_j=D(v_{j-1}+qv_{j+1}-(q+1)v_j)
 \qquad v\in\ell^\infty(\Z),\ j\in\Z.
\]
The semigroup $S_t$ preserves nonnegativity. Since
$0\le h(W)\le aW$, variation of constants gives
\begin{equation}\label{eq:homogeneous-step-linear-bounds}
 \begin{aligned}
 0\le W_j(t)&\le e^{at}(S_t\ind_{\{j\le0\}})_j,\\
 0\le1-W_j(t)&\le(S_t\ind_{\{j>0\}})_j.
 \end{aligned}
 \qquad j\in\Z,\ t\ge0.
\end{equation}
The bound on $1-W_j(t)$ in \eqref{eq:homogeneous-step-linear-bounds}
follows by linear comparison from
$\partial_t(1-W)=\mathcal A(1-W)-h(W)\le\mathcal A(1-W)$.
{Let $X=(X_t)_{t\ge0}$ be the continuous-time chain on $\Z$
with generator $\mathcal A$, and let $\mathbb P_j^{\mathrm d}$ and
$\mathbb E_j^{\mathrm d}$ denote its law and expectation when $X_0=j$.
For fixed $t\ge0$, let $N_t$ be its number of jumps in $[0,t]$.
The holding times of $X$ are independent exponential variables of rate}
$D(q+1)$, so $N_t$ is Poisson with mean $D(q+1)t$;
see \cite[Sections~2.4 and 2.6]{Norris}. The identity
$(S_tv)_j=\mathbb E_j^{\mathrm d}[v_{X_t}]$ for
$v\in\ell^\infty(\Z)$, together with the fact that each jump has size
one, turns \eqref{eq:homogeneous-step-linear-bounds} into
$W_j(t)\le e^{at}\mathbb P_j^{\mathrm d}(N_t\ge j)$ for $j\ge1$ and
$1-W_j(t)\le\mathbb P_j^{\mathrm d}(N_t\ge1-j)$ for $j\le0$.
The Poisson tails therefore imply $W_j(t)\to0$ as $j\to+\infty$
and $W_j(t)\to1$ as $j\to-\infty$. For $0<\theta<1$,
 $\{j\in\Z:W_j(t)\ge\theta\}$ is therefore a nonempty subset
of $\Z$ that is bounded above and hence has a largest element.

For $P=1$ and $f_0=h$, \eqref{eq:matrix} and
\eqref{eq:homogeneous-parameters} give
$k(\lambda)=\Psi_a(\lambda)$, $\lambda_+=\lambda_*(a)$, and
$c_+=c_*(a)$. Define
\[
 m(t)=c_*(a)t-\frac{3}{2\lambda_*(a)}\log(1+t)
 \qquad t\ge0.
\]
Fix $\theta\in(0,1)$, and let $C$ be the constant in
Lemma~\ref{lem:sharp-step-upper} for this reaction. Since
$(1+y)e^{-\lambda_*(a)y}\to0$ as $y\to\infty$, choose
$Y_\theta>0$ such that
$C(1+y)e^{-\lambda_*(a)y}<\theta$ for every $y\ge Y_\theta$.
Applying Lemmas~\ref{lem:sharp-step-upper} and
\ref{lem:sharp-step-lower} with $b=0$, we obtain constants
$C_\theta,T_\theta\in(0,\infty)$ such that
\[
 \lfloor m(t)-C_\theta\rfloor
 \le\max\{j:W_j(t)\ge\theta\}<m(t)+Y_\theta
 \text{ for }t\ge T_\theta.
\]
All these constants are independent of $t$. The quantity
$\lfloor m(t)-C_\theta\rfloor$ differs from $m(t)-C_\theta$ by less than one. In view of
\[
 0\le\log(1+t)-\log t=\log(1+t^{-1})\le t^{-1}
 \text{ for }t\ge1,
\]
the two bounds prove \eqref{eq:hom-lattice-bramson} for $K=1$.
The initial normalization then gives the same conclusion for general
$K$, by applying the bounds to $\widetilde W$ at level $\theta/K$.

\emph{An alternative derivation from the nonlocal equation.}
We again work with $K=1$.
For $z\in\R$, let $\delta_z$ denote the unit point mass at $z$. Set
\begin{equation}\label{eq:homogeneous-jump-measure}
\begin{gathered}
 \nu=D(q+1),\qquad
 J=\frac1{q+1}\delta_1+\frac q{q+1}\delta_{-1},\\
 \qquad
 U(t,x)=W_j(t)\quad\text{ for }x\in(j-1,j]
 \qquad j\in\Z,\ t\ge0.
\end{gathered}
\end{equation}
Write $(J*U)(t,x)=\int_\R U(t,x-z)J(\dd z)$.
Integer translations preserve the intervals $(j-1,j]$, so $U$
satisfies
\[
 U_t=\nu(J*U-U)+h(U),\qquad U(0,x)=\ind_{(-\infty,0]}(x)
 \qquad t\ge0,\ x\in\R.
\]
Indeed,
\[
 \nu\bigl((J*U)(t,x)-U(t,x)\bigr)
 =D\bigl(U(t,x-1)+qU(t,x+1)-(q+1)U(t,x)\bigr)
 \qquad t\ge0,\ x\in\R,
\]
and the integral equations agree on each interval. If $U$ and $V$
are bounded solutions with the same initial datum and values in
$[0,1]$, set $d(t)=\sup_{x\in\R}|U(t,x)-V(t,x)|$. Their integral
equations give
\[
 d(t)\le\bigl(2\nu+\|h'\|_{L^\infty(0,1)}\bigr)
              \int_0^t d(s)\,\dd s
 \qquad t\ge0.
\]
Gronwall's inequality gives $d\equiv0$, so $U$ is the unique bounded
$[0,1]$-valued solution of the real-line initial-value problem.

By its definition, $J$ is a Borel probability measure supported
in $\{-1,1\}$, with $J(\{0\})=0$ and
$J((0,\infty))=1/(q+1)>0$. Moreover, $h\in C^2([0,1])$ gives
$C^{1,\gamma}$ regularity near zero for every $0<\gamma<1$, while
$h(0)=h(1)=0$ and $0<h(s)\le as$ for $0<s<1$.
The datum $U(0,\cdot)=\ind_{(-\infty,0]}$ takes values in $[0,1]$,
equals one for $x<-1$, and equals zero for $x>1$.
These measure, reaction, and initial-data assumptions satisfy the
hypotheses of Graham's logarithmic-correction theorem \cite{Graham}.
For $\lambda>0$, its variational function is
\begin{equation}\label{eq:homogeneous-graham-function}
 \Gamma(\lambda)
 =\frac{\nu\int_\R e^{\lambda x}J(\dd x)-\nu+a}{\lambda}
 =\frac{a+D(e^\lambda+qe^{-\lambda}-(q+1))}{\lambda}
 =\frac{\Psi_a(\lambda)}{\lambda}.
\end{equation}
The unique minimizer and minimum of $\Gamma$ are $\lambda_*(a)$ and
$c_*(a)$;
see \eqref{eq:homogeneous-parameters} and
Lemma~\ref{lem:matrix} with $P=1$. Graham's theorem therefore gives
\[
 \sup\{x:U(t,x)\ge\theta\}
 =c_*(a)t-\frac{3}{2\lambda_*(a)}\log t+O_\theta(1)
 \quad\text{ as } t\to\infty.
\]
With the chosen interval convention, this real-line supremum equals
$\max\{j:W_j(t)\ge\theta\}$. The conclusion of Graham's theorem is stated for $c_*(a)\in\R$.
It therefore applies when $a-\rho_q$ is positive, zero, or negative. Returning to
$W=K\widetilde W$ and the original level $\theta$ proves
\eqref{eq:hom-lattice-bramson}.
\end{proof}

\begin{lemma}\label{lem:hom-tree}

Let $h$ satisfy \eqref{eq:hyp} with $h'(0)=a>\rho_q$.
Let $u$ solve \eqref{eq:tree-cauchy} with $f(x,s)=h(s)$, so that
\[
 u_t=D\Delta_{\Tq}u+h(u)
 \qquad t\ge0,\ x\in\Tq,
\]
with $0\le u_0\le1$, $u_0\not\equiv0$, and $u_0$ finitely
supported. For $t>0$, put
\[
 m_a(t)=c_*(a)t-\frac{3}{2\lambda_*(a)}\log t.
\]
For every $\theta\in(0,1)$ there are $C_\theta,T_\theta>0$ such that
\begin{equation}
\label{eq:hom-tree-balls}
 B(o,m_a(t)-C_\theta)
 \subseteq\{x\in\Tq:u(t,x)\ge\theta\}
 \subseteq B(o,m_a(t)+C_\theta),
 \quad\text{ for }t\ge T_\theta.
\end{equation}
The constants depend only on $\theta$, $h$, $D$, $q$, $u_0$, and $o$.

\end{lemma}
\begin{proof}

Fix $\theta\in(0,1)$. Since $\mu=a-\rho_q>0$,
Lemmas~\ref{lem:horo-spectrum} and~\ref{lem:local-general}, applied
to the potential $a(x)\equiv a$, give $u(t,x)\to1$ as $t\to\infty$ at
every fixed vertex. Lemma~\ref{lem:matrix} also gives $c_*(a)>0$.

\emph{The outer bound for radial data.}
First suppose that $u_0$ is radial. The equation and initial datum
are invariant under every tree automorphism fixing $o$. These
automorphisms act transitively on each set $\{x:r(x)=r\}$, so
uniqueness in Lemma~\ref{lem:comparison} gives $u(t,x)=U_{r(x)}(t)$.
The sequence $U$ satisfies \eqref{eq:radial} with $f_r=h$.
Choose an integer $R\ge0$ such that $U_r(0)=0$ for $r>R$, and let
$W^R$ solve \eqref{eq:lattice} with $f_j=h$ and
$W_j^R(0)=\ind_{\{j\le R\}}$. Translation invariance and
Lemma~\ref{lem:graham} give
\[
 W_{j+1}^R(t)\le W_j^R(t)\quad(j\in\Z,\ t\ge0),\qquad
 \max\{j:W_j^R(t)\ge\theta\}=R+m_a(t)+O_\theta(1)
 \quad\text{ as }t\to\infty.
\]
The restriction of $W^R$ to $r\ge0$ satisfies the equation for
$U_r$ in \eqref{eq:radial} for every $r\ge1$. At the root,
\begin{equation}\label{eq:homogeneous-radial-upper-root}
 \dot W_0^R-D(q+1)(W_1^R-W_0^R)-h(W_0^R)
 =D(W_{-1}^R-W_1^R)\ge0
 \qquad t\ge0.
\end{equation}
By \eqref{eq:homogeneous-radial-upper-root} and
Lemma~\ref{lem:comparison}, the initial
ordering $U_r(0)\le W_r^R(0)$ gives $U_r(t)\le W_r^R(t)$ for every
$r\ge0$. If $U_r(t)\ge\theta$, then $W_r^R(t)\ge\theta$, and
\eqref{eq:hom-lattice-bramson}, applied to $W^R$ with the shift $R$,
gives $r\le m_a(t)+C_\theta$ for
all sufficiently large $t$, after including the fixed shift $R$ in
$C_\theta$. This proves the outer inclusion in
\eqref{eq:hom-tree-balls} for radial data.

\emph{The inner bound for radial data.}
Fix $\theta<\beta<1$. The convergence at the root gives $T_0>0$
such that $U_0(T_0+t)\ge\beta$ for every $t\ge0$. Define
\[
 h_\beta(s)=h(s)(1-s/\beta),\qquad 0\le s\le\beta.
\]
$h_\beta$ belongs to $C^2([0,\beta])$ and satisfies
\[
 \begin{gathered}
 h_\beta(0)=h_\beta(\beta)=0,\qquad
 h_\beta'(0)=a,\qquad h_\beta'(\beta)=-h(\beta)/\beta<0,\\
 0<h_\beta(s)\le h(s)\le as\quad\text{ for } 0<s<\beta.
 \end{gathered}
\]
Define $\widetilde h_\beta(s)=\beta^{-1}h_\beta(\beta s)$ for
$0\le s\le1$. The preceding properties of $h_\beta$ give
$\widetilde h_\beta'(0)=a$,
$\widetilde h_\beta'(1)=-h(\beta)/\beta<0$, and
$0<\widetilde h_\beta(s)\le as$ for $0<s<1$.
Together with $\widetilde h_\beta\in C^2([0,1])$ and
$\widetilde h_\beta(0)=\widetilde h_\beta(1)=0$, the derivative
and positivity bounds verify \eqref{eq:hyp}. Apply Lemma~\ref{lem:comparison}
to the lattice equation with reaction $\widetilde h_\beta$ and
initial datum $\ind_{\{j\le0\}}$, and multiply its solution by
$\beta$. We obtain the solution $V$ of \eqref{eq:lattice} with
reaction $h_\beta$ and initial datum
$V_j(0)=\beta\ind_{\{j\le0\}}$, with $0\le V\le\beta$.
Since $V_r(0)=0$ for $r\ge1$ and $0\le V_0(t)\le\beta$,
the initial and boundary values satisfy
\begin{equation}\label{eq:homogeneous-radial-lower-data}
 V_r(0)=0\le U_r(T_0)\quad\text{ for }r\in\N,\qquad
 V_0(t)\le\beta\le U_0(T_0+t)\quad\text{ for }t\ge0.
\end{equation}
For $r\ge1$, the equation for $V$ gives
\[
 \dot V_r-D(V_{r-1}+qV_{r+1}-(q+1)V_r)-h(V_r)
 =h_\beta(V_r)-h(V_r)\le0
 \qquad r\in\N,\ t\ge0.
\]
Thus $(V_r)_{r\in\N_0}$ is a subsolution of the interior equation
in \eqref{eq:radial} with reaction $h$ on $\{1,2,\ldots\}$.
The only exterior index is $r=0$, where
$V_0(t)\le\beta\le U_0(T_0+t)$ for every $t\ge0$ by
\eqref{eq:homogeneous-radial-lower-data}. The solutions $V$ and $U(T_0+\cdot)$ belong to
$C^1([0,\infty);\ell^\infty)$, so their restrictions satisfy the
time-regularity assumptions of Lemma~\ref{lem:comparison}. Applying that lemma on
$\{r\ge1\}$ gives
\begin{equation}\label{eq:homogeneous-radial-lower-comparison}
 V_r(t)\le U_r(T_0+t)\quad\text{ for }r\in\N,\ t\ge0.
\end{equation}
The endpoint values and derivatives of $h_\beta$, together with
$h_\beta\in C^2([0,\beta])$ and $0<h_\beta(s)\le as$ for
$0<s<\beta$, verify the hypotheses of Lemma~\ref{lem:graham} with
$K=\beta$. Since $h_\beta'(0)=a$, its parameters are
$c_*(a)$ and $\lambda_*(a)$.
By \eqref{eq:hom-lattice-bramson}, applied to $V$ with $K=\beta$,
$\max\{j:V_j(t)\ge\theta\}\ge m_a(t)-C_\theta$ for all sufficiently
large $t$. Since $V$ is nonincreasing in $j$, we have
$V_r(t)\ge\theta$ for every integer $r\le m_a(t)-C_\theta$.
Combining this bound with \eqref{eq:homogeneous-radial-lower-comparison},
we obtain
\[
 U_r(T_0+t)\ge\theta
 \quad\text{ for }1\le r\le m_a(t)-C_\theta
 \qquad r\in\N,
\]
for all sufficiently large $t$. At $r=0$, the lower bound follows
from $U_0(T_0+t)\ge\beta>\theta$. For $t\ge2T_0$,
\[
 \begin{aligned}
 m_a(t-T_0)-m_a(t)
 &=-c_*(a)T_0-\frac{3}{2\lambda_*(a)}\log(1-T_0/t),\\
 |m_a(t-T_0)-m_a(t)|
 &\le c_*(a)T_0+\frac{3}{2\lambda_*(a)}\log2.
 \end{aligned}
\]
Replacing $t$ by $t-T_0$ in the lower estimate and increasing
$C_\theta,T_\theta$ therefore proves the inner inclusion in
\eqref{eq:hom-tree-balls} for radial data.

\emph{General finitely supported data.}
Choose an integer $R\ge0$ such that $u_0\le\ind_{B(o,R)}$, and let
$u^R$ solve $u_t^R=D\Delta_{\Tq}u^R+h(u^R)$ with
$u^R(0,x)=\ind_{B(o,R)}(x)$. The equation and this initial datum are
invariant under tree automorphisms fixing $o$, so uniqueness in
Lemma~\ref{lem:comparison} makes $u^R$ radial. The same lemma gives
$u(t,x)\le u^R(t,x)$.
The radial upper estimate consequently gives
$\{x:u^R(t,x)\ge\theta\}\subseteq B(o,m_a(t)+C_\theta)$ for large $t$.
For the lower bound, fix $\tau>0$. Strict positivity in
Lemma~\ref{lem:comparison} gives $u(\tau,o)>0$, so choose
$0<\varepsilon\le u(\tau,o)$. Let $u^\varepsilon$ solve
$u_t^\varepsilon=D\Delta_{\Tq}u^\varepsilon+h(u^\varepsilon)$ from
$u^\varepsilon(0,x)=\varepsilon\ind_{\{o\}}(x)$. The initial ordering
and Lemma~\ref{lem:comparison} give
$u^\varepsilon(t,x)\le u(\tau+t,x)$. By uniqueness,
$u^\varepsilon$ is radial, so its radial lower estimate gives
\[
 u(\tau+t,x)\ge u^\varepsilon(t,x)\ge\theta
 \quad\text{ for }r(x)\le m_a(t)-C_\theta
 \qquad x\in\Tq,
\]
for all sufficiently large $t$. For $t\ge2\tau$,
$|m_a(t-\tau)-m_a(t)|\le c_*(a)\tau+
(3/(2\lambda_*(a)))\log2$. Hence the lower bound for
$u(\tau+t,x)$, with $t$ replaced by $t-\tau$, gives
$u(t,x)\ge\theta$ for $r(x)\le m_a(t)-C_\theta$ after increasing
$C_\theta$.
Choose $C_\theta$ and $T_\theta$ large enough for both comparisons.
This proves \eqref{eq:hom-tree-balls} for the prescribed datum $u_0$.\qedhere

\end{proof}

\begin{lemma}\label{lem:envelopes}
Suppose a finite family $f_i$ satisfies \eqref{eq:hyp} and
$f_i'(0)=a$ for every $i$. Then there are
$h_-,h_+\in C^2([0,1])$, each satisfying \eqref{eq:hyp} with
derivative $a$ at zero, such that
\begin{equation}
\label{eq:scalar-envelopes}
 h_-(s)\le f_i(s)\le h_+(s)
 \quad\text{ for }0\le s\le1
 \text{ and }i\in\{0,\ldots,P-1\}.
\end{equation}
\end{lemma}
\begin{proof}

\emph{Construction of $h_-$.}
For each $i\in\{0,\ldots,P-1\}$, define
\[
 R_i(s)=\frac{f_i(s)}{as(1-s)}
 \quad\text{ for }0<s<1.
\]
By \eqref{eq:hyp}, Taylor's formula gives
$f_i(s)=as+O(s^2)$ as $s\downarrow0$ and
$f_i(s)=-f_i'(1)(1-s)+O((1-s)^2)$ as $s\uparrow1$.
Consequently, each $R_i$ extends continuously to $[0,1]$ and is
strictly positive there, with
\[
 R_i(0)=1,\qquad R_i(1)=-f_i'(1)/a>0
 \qquad i\in\{0,\ldots,P-1\}.
\]
The expansion at zero gives $R_i(s)=1+O(s)$ uniformly in $i$.
Hence $R_i(s)\ge1/2$ for all sufficiently small $s>0$, uniformly in
$i$, and
\[
 \left|\frac{R_i(s)^{-1}-1}{s}\right|
 =\frac{|R_i(s)-1|}{sR_i(s)}\le C
 \qquad i\in\{0,\ldots,P-1\},\ \text{for all sufficiently small }s>0.
\]
On each $[\delta,1]$ with $0<\delta<1$, continuity, positivity and
finiteness of the family bound this quotient uniformly in $i$.
Consequently
\[
 K_0=\max\left\{0,
       \sup_{\substack{0\le i<P\\0<s\le1}}
             \frac{R_i(s)^{-1}-1}{s}\right\}<\infty.
\]
By the definition of $K_0$, $R_i(s)\ge(1+K_0s)^{-1}$ for
$0<s\le1$. Multiplying by $as(1-s)$ and using the endpoint values
therefore gives
\begin{equation}\label{eq:scalar-lower-reaction}
 h_-(s):=\frac{as(1-s)}{1+K_0s}\le f_i(s)
 \quad\text{ for } 0\le s\le1
 \qquad i\in\{0,\ldots,P-1\}.
\end{equation}
The denominator is positive on $[0,1]$, so $h_-$ is $C^2$ and
has zero endpoint values. Moreover,
\[
 h_-'(0)=a,\qquad h_-'(1)=-\frac{a}{1+K_0}<0,
 \qquad 0<h_-(s)\le as\quad\text{ for } 0<s<1.
\]

\emph{Construction of $h_+$.}
Let $L=\max_{0\le i<P}\|f_i'\|_{L^\infty(0,1)}>0$.
Since $f_i(1)=0$, the fundamental theorem of calculus gives
$f_i(s)=-\int_s^1 f_i'(r)\,\dd r\le L(1-s)$ for $0\le s\le1$.
Choose
\[
 0<\varepsilon\le\min\left\{\frac12,\frac{a}{2L}\right\},
\]
and define $H_0:[0,\infty)\to[0,1]$ and
$h_+:[0,1]\to\R$ by
\begin{equation}\label{eq:scalar-upper-reaction}
 H_0(z)=
 \begin{cases}
 1-(1-z)^3,&0\le z\le1,\\
 1,&z\ge1,
 \end{cases}
 \qquad
 h_+(s)=asH_0\left(\frac{1-s}{\varepsilon}\right).
\end{equation}
For $0\le z\le1$, one has
$H_0'(z)=3(1-z)^2$ and $H_0''(z)=-6(1-z)$.
Thus both derivatives vanish at $z=1$, where $H_0(1)=1$, so
$h_+\in C^2([0,1])$. Moreover, $0\le H_0\le1$ on $[0,\infty)$ and
$H_0(z)-z=z(1-z)(2-z)\ge0$ on $[0,1]$.
For $0\le s\le1-\varepsilon$, \eqref{eq:hyp} gives
$h_+(s)=as\ge f_i(s)$. For $1-\varepsilon\le s\le1$, the bounds
$s\ge1/2$ and $H_0((1-s)/\varepsilon)\ge(1-s)/\varepsilon$
give
\[
 h_+(s)\ge\frac{a}{2\varepsilon}(1-s)
          \ge L(1-s)\ge f_i(s)
 \qquad i\in\{0,\ldots,P-1\}.
\]
Since $\varepsilon\le1/2$, we have $h_+(s)=as$ near zero.
The definition of $h_+$ gives zero endpoint values.
Differentiating $h_+(s)=asH_0((1-s)/\varepsilon)$ at $s=1$ and
using $H_0(0)=0$ and $H_0'(0)=3$ yields
\[
 h_+'(0)=a,\qquad h_+'(1)=-\frac{3a}{\varepsilon}<0,
 \qquad 0<h_+(s)\le as\quad\text{ for } 0<s<1.
\]
Thus $h_-$ and $h_+$ satisfy \eqref{eq:hyp}, and
$h_-(s)\le f_i(s)\le h_+(s)$ for every $i$ and $0\le s\le1$,
which proves \eqref{eq:scalar-envelopes}. Since
$h_-'(0)=h_+'(0)=a$, the $P=1$ matrix in \eqref{eq:matrix} has the
single entry $\Psi_a(\lambda)=a+D(e^\lambda+qe^{-\lambda}-(q+1))$
for either reaction. Its principal eigenvalue is therefore
$\Psi_a(\lambda)$. The definitions in
\eqref{eq:homogeneous-parameters} give the same speed $c_*(a)$ and
minimizing exponent $\lambda_*(a)$ for both reactions.\qedhere

\end{proof}

\subsection{Front equations}\label{sec:waves}

Throughout this section we assume \eqref{eq:hyp}. We first construct
moving fronts by monotone iteration and then treat speed zero,
including the critical case $\mu=0$. Lemmas~\ref{lem:difference}
and~\ref{lem:tree-difference} assume \eqref{eq:derivative-kpp}, which
bounds every secant slope of $f_i$ by $a_i$.

For $c\ne0$, a pulsating front from one to zero is a solution of
\eqref{eq:lattice} of the form $w_j(t)=\Phi_j(j-ct)$ for
$t\in\R$ and $j\in\Z$, where
$\Phi_i\in C^1(\R;(0,1))$, $\Phi_{j+P}=\Phi_j$, and
\[
 \lim_{z\to-\infty}\Phi_i(z)=1,\qquad
 \lim_{z\to+\infty}\Phi_i(z)=0
\]
uniformly for $i\in\{0,\ldots,P-1\}$. The profiles satisfy
\begin{equation}
\label{eq:profile}
 -c\Phi_i'(z)
 =
 D\bigl(\Phi_{i-1}(z-1)+q\Phi_{i+1}(z+1)
 -(q+1)\Phi_i(z)\bigr)+f_i(\Phi_i(z)).
\end{equation}
Equation~\eqref{eq:profile} holds for every $i\in\Z$ and $z\in\R$.
Periodicity in the phase index gives, for $c\ne0$,
\[
 w_{j+P}(t+P/c)=w_j(t)
 \qquad j\in\Z,\ t\in\R.
\]
At speed zero, we use $w_j(t)=p_j$, with $p:\Z\to(0,1)$
satisfying
\[
 \begin{gathered}
 D(p_{j-1}+qp_{j+1}-(q+1)p_j)+f_j(p_j)=0
 \quad\text{for }j\in\Z,\\
 \lim_{j\to-\infty}p_j=1,
 \qquad \lim_{j\to+\infty}p_j=0.
 \end{gathered}
\]

For $c\in\R$ and a $P$-periodic family $V=(V_i)_{i\in\Z}$ of
locally Lipschitz functions $V_i:\R\to[0,1]$, write
\begin{equation}
\label{eq:profile-residual}
 (\mathcal R_cV)_i(z)
 :=
 cV_i'(z)
 +D\bigl(V_{i-1}(z-1)+qV_{i+1}(z+1)-(q+1)V_i(z)\bigr)
 +f_i(V_i(z))
 \qquad i\in\Z,\ \text{a.e. }z\in\R.
\end{equation}
The profile equation is $\mathcal R_cV=0$. We call $V$ a profile
subsolution if $(\mathcal R_cV)_i(z)\ge0$ for every $i\in\Z$ and
almost every $z\in\R$, and a profile supersolution if
$(\mathcal R_cV)_i(z)\le0$ for every $i\in\Z$ and almost every
$z\in\R$.

\begin{lemma}\label{lem:periodic-equilibria}
Suppose $\xi=(\xi_i)_{i\in\Z}$ is $P$-periodic, $0\le\xi_i\le1$, and
\[
 D(\xi_{i-1}+q\xi_{i+1}-(q+1)\xi_i)+f_i(\xi_i)=0
 \quad\text{ for } i\in\Z.
\]
Then either $\xi_i=0$ for every $i$ or $\xi_i=1$ for every $i$.
\end{lemma}

\begin{proof}

Sum the stationary equations for $i=0,\ldots,P-1$. Periodicity gives
$\sum_{i=0}^{P-1}\xi_{i-1}=\sum_{i=0}^{P-1}\xi_i
=\sum_{i=0}^{P-1}\xi_{i+1}$, so the diffusion sum is zero and
\[
 \sum_{i=0}^{P-1}f_i(\xi_i)=0.
\]
By \eqref{eq:hyp}, $f_i(s)\ge0$ for $s\in[0,1]$, with equality
if and only if $s\in\{0,1\}$. Thus every $\xi_i$ belongs to $\{0,1\}$. If
$\xi_i=0$ for one index, its equation gives
$\xi_{i-1}=\xi_{i+1}=0$. Repeating the equation at the indices
$i+1,\ldots,i+P-1$ gives $\xi_j=0$ for $0\le j<P$, and hence for $j\in\Z$ by periodicity. If all
components are positive, their membership in $\{0,1\}$ gives
$\xi_j=1$ for every $j$.\qedhere

\end{proof}

\subsection{Moving fronts}

\begin{lemma}\label{lem:wave-construction}

Assume \eqref{eq:hyp} and $\mu>0$. For every $c\ge c_+$ there are
$C^1$ profiles $(\Phi_i)_{i\in\Z}$ satisfying
$\Phi_{i+P}=\Phi_i$ and \eqref{eq:profile}, with
$\Phi_i'(z)<0$ for all $i$ and $z$, and with $\Phi_i(z)\to1$ as $z\to-\infty$ and
$\Phi_i(z)\to0$ as $z\to+\infty$.
If $c>c_+$, let $\lambda$ be the smaller positive root of
\[
 k(\lambda)=c\lambda.
\]
Then $(\Phi_i)$ may be chosen so that
\begin{equation}
\label{eq:front-tail}
 e^{\lambda z}\Phi_i(z)\to A_c\phi_i(\lambda)
 \quad\text{ as } z\to+\infty
\end{equation}
for some $A_c>0$, uniformly over the finite phase set. No front with
these end states has nonzero speed $c<c_+$, and no stationary sequence
with these end states exists. The speed exclusions require only the prescribed end limits.

\end{lemma}

\begin{proof}

\emph{Exponential barriers.}
Fix $c>c_+$. By Lemma~\ref{lem:matrix}, $k(0)>0$,
$k(\lambda_+)-c\lambda_+=\lambda_+(c_+-c)<0$, and
$k(r)-cr\to+\infty$ as $r\to+\infty$. Strict convexity therefore
gives exactly two positive zeros $\lambda<\lambda_{\rm r}$ of
$k(r)-cr$. Choose
\[
 \lambda<\lambda'
 <
 \min\{2\lambda,\lambda_{\rm r}\},
\]
and put $\delta=c\lambda'-k(\lambda')>0$.
Define
\[
 E_i(z)=e^{-\lambda z}\phi_i(\lambda),
 \qquad
 \widetilde E_i(z)=e^{-\lambda' z}\phi_i(\lambda')
 \qquad i\in\Z,\ z\in\R.
\]
With the linearized profile operator on $P$-periodic families
of $C^1$ functions $V_i:\R\to\R$ given by
\[
 (\mathcal P_cV)_i
 :=
 cV_i'
 +D\bigl(V_{i-1}(\,\cdot-1)+qV_{i+1}(\,\cdot+1)
 -(q+1)V_i\bigr)+a_iV_i
 \qquad i\in\Z,
\]
the eigenvalue equation in \eqref{eq:matrix} gives
\[
 \mathcal P_cE=0,
 \qquad
 \mathcal P_c\widetilde E=-\delta\widetilde E.
\]

For $A\ge1$, to be fixed below, set
\[
 \overline\Phi_i(z)=\min\{1,E_i(z)\},
 \qquad
 \underline\Phi_i(z)
 =
 \bigl(E_i(z)-A\widetilde E_i(z)\bigr)_+
 \qquad i\in\Z,\ z\in\R.
\]
We choose $A$ so that $0\le\underline\Phi_i\le1$ and
$\mathcal R_c\underline\Phi\ge0$ almost everywhere. Put
$d=\lambda'-\lambda>0$ and
$r_i=\phi_i(\lambda')/\phi_i(\lambda)$. Define $z_i(A)=(\log A+\log r_i)/d$. Then
$\underline\Phi_i(z)>0$ exactly when $z>z_i(A)$.
For $z>z_i(A)$, positivity of $\phi(\lambda)$ and
$\phi(\lambda')$ and the inequality $\lambda'<2\lambda$ give
constants $C,c_0>0$, independent of $A\ge1$ and
$i\in\{0,\ldots,P-1\}$, such that
\[
 E_i(z)\le C A^{-\lambda/d},\qquad
 \frac{A\delta\widetilde E_i(z)}{E_i(z)^2}
       \ge c_0 A^{\lambda/d}.
\]
Indeed, the second ratio equals
$A\delta\phi_i(\lambda')\phi_i(\lambda)^{-2}
 e^{(2\lambda-\lambda')z}$. Since $2\lambda-\lambda'>0$, its
infimum for $z>z_i(A)$ is its limit at $z_i(A)$, which is a positive
constant depending on $i$ times $A^{\lambda/d}$. Taking the minimum
over $i\in\{0,\ldots,P-1\}$ gives $c_0>0$ independent of $A$.
Fix $A$ large enough that
\[
 C A^{-\lambda/d}\le1,\qquad c_0 A^{\lambda/d}\ge C_2.
\]
Then $0\le\underline\Phi\le\overline\Phi\le1$ and
$\underline\Phi\not\equiv0$.

At a point where $\underline\Phi_i(z)>0$,
\[
 \mathcal P_c(E-A\widetilde E)_i=A\delta\widetilde E_i.
\]
Replacing $E_{i-1}(z-1)-A\widetilde E_{i-1}(z-1)$ and
$E_{i+1}(z+1)-A\widetilde E_{i+1}(z+1)$ by their positive parts
increases their diffusion contributions, whose coefficients are
$D>0$ and $Dq>0$. Since $0<\underline\Phi_i(z)\le E_i(z)\le1$
at the point under consideration, \eqref{eq:taylor} yields
\[
 (\mathcal R_c\underline\Phi)_i
 \ge A\delta\widetilde E_i-C_2\underline\Phi_i^2
 \ge A\delta\widetilde E_i-C_2E_i^2\ge0.
\]
Where $\underline\Phi_i=0$, its profile derivative vanishes almost
everywhere, the reaction vanishes, and the neighboring contributions
are nonnegative. Thus $\mathcal R_c\underline\Phi\ge0$ almost
everywhere. The functions $\underline\Phi$ and $\overline\Phi$ are bounded
and locally absolutely continuous.

If $\overline\Phi_i=E_i<1$, truncation decreases the neighboring
values and $f_i(E_i)\le a_iE_i$, so
$\mathcal R_c\overline\Phi\le0$. If $\overline\Phi_i=1$, its
derivative vanishes almost everywhere, $f_i(1)=0$, and the neighboring
values are at most one. Hence $(\mathcal R_c\overline\Phi)_i\le0$
also where $\overline\Phi_i=1$, at every point of differentiability.

\emph{Monotone iteration.}
With $L$ as in \eqref{eq:taylor}, choose
$\beta\ge D(q+1)+L$. For a bounded measurable $P$-periodic family
$V=(V_i)_{i\in\Z}$ with $V_i:\R\to[0,1]$, define
\[
 \mathcal H_i(V)(z)
 :=
 DV_{i-1}(z-1)+DqV_{i+1}(z+1)
 +(\beta-D(q+1))V_i(z)+f_i(V_i(z))
 \qquad i\in\Z,\ z\in\R.
\]
The order on profile vectors is the pointwise order:
$V\le W$ means $V_i(z)\le W_i(z)$ for every $i,z$.
Since $\beta-D(q+1)+f_i'(s)\ge0$ for $0\le s\le1$,
$\mathcal H$ preserves this order. Whenever $0\le V_i(z)\le1$
for every $i,z$, one has $0\le\mathcal H_i(V)(z)\le\beta$
for every $i,z$. If every component $V_i$ is nonincreasing in $z$, then
$z\mapsto\mathcal H_i(V)(z)$ is nonincreasing for every $i$. Since $\mu>0$, Lemma~\ref{lem:matrix} gives
$c_+>0$, and hence $c>0$. For a bounded measurable $P$-periodic
family $V$ with values in $[0,1]$, define $\mathcal TV$ by
\begin{equation}\label{eq:positive-front-iteration}
 (\mathcal TV)_i(z)
 =
 \frac1c\int_z^\infty
 e^{-\beta(s-z)/c}\mathcal H_i(V)(s)\,\dd s
 \qquad i\in\Z,\ z\in\R.
\end{equation}
For either barrier $V$, one has
$\mathcal R_cV=cV'-\beta V+\mathcal H(V)$. Multiplication by
$e^{-\beta(s-z)/c}$ and integration over $s\in[z,\infty)$ give
\[
 V_i(z)-(\mathcal TV)_i(z)
 =-\frac1c\int_z^\infty e^{-\beta(s-z)/c}
                    (\mathcal R_cV)_i(s)\,\dd s
 \qquad i\in\Z,\ z\in\R.
\]
The boundary term at infinity is zero because $V$ is bounded.
The inequalities $\mathcal R_c\underline\Phi\ge0$ and
$\mathcal R_c\overline\Phi\le0$ therefore imply
\begin{equation}\label{eq:positive-front-barrier-order}
 \mathcal T\underline\Phi\ge\underline\Phi,
 \qquad
 \mathcal T\overline\Phi\le\overline\Phi.
\end{equation}
By \eqref{eq:positive-front-iteration}, positivity of the integral
kernel and the order property of $\mathcal H$ imply that
$\mathcal T$ preserves order. The substitution $s=z+r$ in that integral shows
that $\mathcal T$ also preserves nonincreasing dependence on $z$:
for each fixed $r\ge0$, $\mathcal H_i(V)(z+r)$ is nonincreasing
whenever all the $V_i$ are nonincreasing. Starting from
$V^0=\overline\Phi$ and setting $V^{n+1}=\mathcal TV^n$, we therefore
obtain nonincreasing profiles with
$\underline\Phi\le V^{n+1}\le V^n\le\overline\Phi$ for every $n\ge0$.
Define $\Phi_i(z)=\lim_{n\to\infty}V_i^n(z)$. By continuity of $f_i$ and convergence of each of the three
components entering $\mathcal H_i$,
$\mathcal H_i(V^n)(s)\to\mathcal H_i(\Phi)(s)$ for every $i,s$.
The integrand defining $(\mathcal TV^n)_i(z)$, including the
prefactor $1/c$, is bounded by $(\beta/c)e^{-\beta(s-z)/c}$,
which is integrable for $s\ge z$.
Dominated convergence therefore gives $\Phi=\mathcal T\Phi$. Since
$0\le\mathcal H_i(\Phi)\le\beta$, the integral defining $\mathcal T$
is locally absolutely continuous and gives
$-c\Phi_i'+\beta\Phi_i=\mathcal H_i(\Phi)$ almost everywhere and
$|\Phi_i'|\le2\beta/c$. Thus the profiles are locally Lipschitz,
$\mathcal H(\Phi)$ is continuous, and \eqref{eq:profile} gives
$\Phi\in C^1$. For $z_1<z_2$, the inequality $V_i^n(z_1)\ge V_i^n(z_2)$
passes to the limit. Thus every $\Phi_i$ is nonincreasing.

\emph{End states and the noncritical tail.}
By boundedness and monotonicity, the limits
$\ell_i^\pm=\lim_{z\to\pm\infty}\Phi_i(z)$ exist.
Put $\ell^\pm=(\ell_i^\pm)_{i=0}^{P-1}$.
In view of \eqref{eq:profile}, $\Phi_i'$ has a limit at each end,
because each shifted profile and each reaction term has a limit. The identity
$\Phi_i(z+1)-\Phi_i(z)=\int_0^1\Phi_i'(z+r)\,\dd r$, together with
the finite end limits of $\Phi_i$, shows that both derivative limits
are zero. Passing to the limit in \eqref{eq:profile} therefore gives
\[
 D(\ell_{i-1}^\pm+q\ell_{i+1}^\pm-(q+1)\ell_i^\pm)
       +f_i(\ell_i^\pm)=0
 \quad\text{ for }i\in\{0,\ldots,P-1\}.
\]
The vectors $\ell^-$ and $\ell^+$ thus satisfy the system in
Lemma~\ref{lem:periodic-equilibria}.  The upper exponential barrier
gives $\ell^+=0$.  For each $i$ and any $z>z_i(A)$,
$\ell_i^-\ge\Phi_i(z)\ge\underline\Phi_i(z)>0$, so $\ell^-$ is nonzero.
Lemma~\ref{lem:periodic-equilibria} therefore gives $\ell^-=1$.

By the definitions of $\underline\Phi$ and $\overline\Phi$,
the inequalities $\underline\Phi\le\Phi\le\overline\Phi$ give
\[
 E_i(z)-A\widetilde E_i(z)
 \le \Phi_i(z)\le E_i(z)
 \qquad i\in\Z,
\]
for all sufficiently large $z$. Dividing these inequalities by $e^{-\lambda z}$ and using
$\lambda'>\lambda$, we obtain
\[
 e^{\lambda z}\Phi_i(z)\to\phi_i(\lambda)
 \quad\text{ as } z\to+\infty
 \qquad i\in\Z.
\]
For $h\in\R$, replacing $\Phi_i(z)$ by $\Phi_i(z+h)$ multiplies the limiting
coefficient by $e^{-\lambda h}$. Since $e^{-\lambda h}$ ranges over $(0,\infty)$, this
gives \eqref{eq:front-tail} with any prescribed $A_c>0$.

\emph{The critical speed.}
Since $\mu>0$,
Lemma~\ref{lem:matrix} gives $c_+>0$.  Choose $c_n\downarrow c_+$ as $n\to\infty$, with $c_n>c_+$, and
translate the corresponding profiles so that $\Phi_0^n(0)=1/2$.
Since $c_n\ge c_+>0$ and $0\le\Phi_i^n\le1$, \eqref{eq:profile} gives
\[
 \sup_{i\in\Z,\,z\in\R}|(\Phi_i^n)'(z)|
 \le\frac{D(q+1)+\max_{0\le i<P}\|f_i\|_{L^\infty([0,1])}}{c_n}\le C.
\]
For each $m\in\N$, apply Arzel\`a--Ascoli to the finitely many
components on $[-m,m]$. A diagonal choice gives indices $n_k\to\infty$
such that $\Phi_i^{n_k}\to\Phi_i$ uniformly on every bounded interval
for $0\le i<P$. Integrate \eqref{eq:profile} between fixed
$z_1<z_2$. Uniform convergence on $[z_1-1,z_2+1]$, continuity of
$f_i$ and $c_{n_k}\to c_+$ permit passage to the limit in that
identity. The resulting integral equation has a continuous
integrand, so $\Phi\in C^1$ and satisfies \eqref{eq:profile} with
$c=c_+$. The inequalities $\Phi_i^{n_k}(z_1)\ge\Phi_i^{n_k}(z_2)$
for $z_1<z_2$ pass to the limit, and uniform convergence near zero
gives $\Phi_0(0)=1/2$.
Define $\ell_i^\pm=\lim_{z\to\pm\infty}\Phi_i(z)$, which exist
by monotonicity. Integrating \eqref{eq:profile} over $[z,z+1]$ gives
\[
\begin{gathered}
 -c_+\bigl(\Phi_i(z+1)-\Phi_i(z)\bigr)
 =\int_z^{z+1}
 \bigl[D(\Phi_{i-1}(r-1)+q\Phi_{i+1}(r+1)-(q+1)\Phi_i(r))
       +f_i(\Phi_i(r))\bigr]\,\dd r
 \\\text{ for }i\in\Z,\ z\in\R.
\end{gathered}
\]
As $z\to\pm\infty$, the left side tends to zero and the right side
tends to
$D(\ell_{i-1}^\pm+q\ell_{i+1}^\pm-(q+1)\ell_i^\pm)+f_i(\ell_i^\pm)$.
Thus the vectors $\ell^\pm=(\ell_i^\pm)_{i=0}^{P-1}$ satisfy
Lemma~\ref{lem:periodic-equilibria}. Since
$\ell_0^-\ge\Phi_0(0)=1/2$ and $\ell_0^+\le\Phi_0(0)=1/2$,
that lemma gives $\ell^-=1$ and $\ell^+=0$.

\emph{Strict monotonicity.}
Since $c\ne0$ and the right-hand side of \eqref{eq:profile} is
$C^1$ for $C^1$ profiles, \eqref{eq:profile} gives $\Phi_i\in C^2$.  For the entire solution
$w_j(t)=\Phi_j(j-ct)$ set
\[
 Y_j(t)=\partial_tw_j(t)=-c\Phi_j'(j-ct)\ge0
 \qquad j\in\Z,\ t\in\R.
\]
Then
\[
 \dot Y_j
 =
 D(Y_{j-1}+qY_{j+1}-(q+1)Y_j)
 +f_j'(w_j(t))Y_j
 \qquad j\in\Z,\ t\in\R.
\]
Since $0\le w_j(t)\le1$, \eqref{eq:lattice} gives
\[
 \|Y(t,\cdot)\|_\infty
 \le D(q+1)+\max_{0\le i<P}\|f_i\|_{L^\infty([0,1])}
 \qquad t\in\R.
\]
Set $K=D(q+1)+L$. For fixed $l,j\in\Z$ and $s<t$, the equation for
$e^{K(t-s)}Y_j(t)$ has nonnegative coefficients, since its diagonal
coefficient is $L+f_j'(w_j(t))\ge0$ and its two off-diagonal
coefficients are $D$ and $Dq$. Iterating the integral equation along
a nearest-neighbor path from $l$ to $j$, of length $n=|j-l|$, gives
\begin{equation}\label{eq:positive-front-path}
 Y_j(t)\ge e^{-K(t-s)}\frac{D^n(t-s)^n}{n!}Y_l(s)
 \quad\text{ for }s<t.
\end{equation}
The factor $(t-s)^n/n!$ is the volume of the ordered jump-time
simplex, and each of the $n$ off-diagonal factors is at least $D$.
Suppose $j_0\in\Z$, $t_0\in\R$ and $Y_{j_0}(t_0)=0$.
By \eqref{eq:positive-front-path}, $Y_l(s)=0$ for every $l$ and
every $s<t_0$. Time continuity then gives $Y(t_0,\cdot)=0$.
The linear equation for $Y$ has uniformly bounded coefficients;
Gronwall's inequality in $\ell^\infty$ gives uniqueness for its
integral equation on every bounded forward time interval. Hence
$Y(t,\cdot)=0$ also for $t\ge t_0$, and $Y\equiv0$ on $\R$.
Since $c\ge c_+>0$, $j-ct$ ranges over $\R$ at each fixed index $j$;
hence $Y\equiv0$ would make every $\Phi_i$ constant, contradicting
its two end limits. Therefore $Y_j(t)>0$ and $\Phi_i'(z)<0$.
The limits one and zero now give $0<\Phi_i(z)<1$.

\emph{Exclusion of lower speeds.}
Let $w$ be a putative moving front of nonzero speed $c<c_+$, or a
stationary sequence with the prescribed end states, in which case
set $c=0$. Since $\mu>0$, both cases have $c<c_+$. Choose
\[
 \max\{c,c_-\}<s<c_+.
\]
Let $v$, $L_0$ and $\delta$ be supplied by Lemma~\ref{lem:bump}
for speed $s$ and $R=1$. Fix $0<\vartheta\le1/2$. The left end
limit gives an integer $J$ with $w_j(0)>1/2$ for $j\le J$.
Choose $\ell\in P\Z$ with $\ell+L_0\le J$ and define
\[
 v^\ell(t,j)=\vartheta v(t,j-\ell)
 \qquad t\in\R,\ j\in\Z.
\]
Lemma~\ref{lem:bump} states that $\vartheta v$ is a subsolution
of \eqref{eq:lattice} for $0<\vartheta\le1$. Since
$\ell\in P\Z$ gives $f_{j-\ell}=f_j$, the shifted function
$v^\ell$ also satisfies that subsolution inequality for every
$j\in\Z$ and almost every $t\in\R$. At time zero it vanishes outside
$\{|j-\ell|<L_0\}$ and is at most $1/2$ on that set, so
$v^\ell(0,j)\le w_j(0)$ everywhere. Lemma~\ref{lem:comparison}
gives $w_j(t)\ge v^\ell(t,j)$ for all $t\ge0$.

Choose $j_t\in\Z$ with $|j_t-(\ell+st)|\le1/2$. The lower bound
in Lemma~\ref{lem:bump} for $R=1$ gives
$w_{j_t}(t)\ge v^\ell(t,j_t)\ge\vartheta\delta$.
For a moving candidate,
\[
 j_t-ct=(s-c)t+\ell+O(1)\to+\infty
 \quad\text{ as }t\to\infty,
\]
so its prescribed right end limit gives $w_{j_t}(t)\to0$.
For a stationary candidate, $c=0$ and $s>0$, so $j_t\to+\infty$
and its sequence limit also gives $w_{j_t}\to0$.
Both conclusions contradict $w_{j_t}(t)\ge\vartheta\delta$.
\end{proof}

\begin{lemma}\label{lem:signed-fronts}
Assume \eqref{eq:hyp}. For every $c\ge c_+$ with $c\ne0$,
equation \eqref{eq:profile} has profiles $\Phi_i\in C^1(\R;(0,1))$
with $\Phi_{i+P}=\Phi_i$ such that $\Phi_i'(z)<0$ for every $i,z$,
$\Phi_i(z)\to1$ as $z\to-\infty$, and $\Phi_i(z)\to0$ as
$z\to+\infty$.
If $c>c_+$, the profiles may be chosen to satisfy
\eqref{eq:front-tail}, where $\lambda$ is the smaller positive root
of $k(\lambda)=c\lambda$.  No front with these limiting states has a
nonzero speed below $c_+$.
\end{lemma}

\begin{proof}

First let $c>c_+$ with $c\ne0$. By Lemma~\ref{lem:matrix}, the function
\[
 \lambda\mapsto k(\lambda)-c\lambda
\]
is strictly convex, is positive at $\lambda=0$, and tends to
$+\infty$ as $\lambda\to+\infty$. Its value at $\lambda_+$ is
$\lambda_+(c_+-c)<0$. It therefore has exactly two positive zeros,
$\lambda<\lambda_{\rm r}$. Choose
$\lambda<\lambda'<\min\{2\lambda,\lambda_{\rm r}\}$ and put
$\delta=c\lambda'-k(\lambda')>0$. For $A\ge1$, to be fixed below, define
\[
\begin{gathered}
 E_i(z)=e^{-\lambda z}\phi_i(\lambda),\qquad
 \widetilde E_i(z)=e^{-\lambda'z}\phi_i(\lambda'),\\
 \overline\Phi_i=\min\{1,E_i\},\qquad
 \underline\Phi_i=(E_i-A\widetilde E_i)_+
 \\\text{ for }i\in\Z,\ z\in\R.
\end{gathered}
\]
Writing $d=\lambda'-\lambda$, the set where
$\underline\Phi_i>0$ is
$z>(\log A+\log(\phi_i(\lambda')/\phi_i(\lambda)))/d$.
The exponential formulas give constants $C,c_0>0$, independent of
the phase $i$ and $A\ge1$, such that on this set
\[
 E_i(z)\le C A^{-\lambda/d},\qquad
 \frac{A\delta\widetilde E_i(z)}{E_i(z)^2}
       \ge c_0A^{\lambda/d}.
\] Choose $A$ so large that
$CA^{-\lambda/d}\le1$ and $c_0A^{\lambda/d}\ge C_2$.
For a $P$-periodic family $V=(V_i)_{i\in\Z}$ of locally Lipschitz
functions $V_i:\R\to\R$, define the linear operator
\[
 (\mathcal P_cV)_i(z)=cV_i'(z)
 +D\bigl(V_{i-1}(z-1)+qV_{i+1}(z+1)-(q+1)V_i(z)\bigr)
 +a_iV_i(z)
 \qquad i\in\Z,\ \text{a.e. }z\in\R.
\]
By \eqref{eq:matrix}, $\mathcal P_cE=0$ and
$\mathcal P_c\widetilde E=-\delta\widetilde E$, so
$\mathcal P_c(E-A\widetilde E)=A\delta\widetilde E$.
At any $(i,z)$ with $\underline\Phi_i(z)>0$, replacing the negative
neighboring values of $E-A\widetilde E$ by zero increases the diffusion
term. Since $0\le\underline\Phi_i\le E_i\le1$ there, the bound
$f_i(r)\ge a_ir-C_2r^2$ in \eqref{eq:taylor} gives
\[
 (\mathcal R_c\underline\Phi)_i
 \ge A\delta\widetilde E_i-C_2E_i^2\ge0
 \quad\text{where }\underline\Phi_i>0.
\]
Where $\underline\Phi_i=0$, its derivative is zero almost everywhere
and the neighboring terms are nonnegative. For $\overline\Phi_i=E_i<1$,
truncation lowers the neighboring terms and $f_i(E_i)\le a_iE_i$,
so $(\mathcal R_c\overline\Phi)_i\le0$. Where
$\overline\Phi_i=1$, its derivative and reaction vanish almost
everywhere and its diffusion term is nonpositive. Thus
$0\le\underline\Phi\le\overline\Phi\le1$,
$\mathcal R_c\underline\Phi\ge0$ and
$\mathcal R_c\overline\Phi\le0$ almost everywhere, for either sign of $c$.

Choose $\beta\ge D(q+1)+L$. For a bounded measurable $P$-periodic
family $V=(V_i)_{i\in\Z}$ with $V_i:\R\to[0,1]$, define
\[
 \mathcal H_i(V)(z)=DV_{i-1}(z-1)+DqV_{i+1}(z+1)
       +(\beta-D(q+1))V_i(z)+f_i(V_i(z))
 \qquad i\in\Z,\ z\in\R.
\]
We use the pointwise order: $V\le W$ means $V_i(z)\le W_i(z)$ for
every $i,z$. By \eqref{eq:taylor},
$\beta-D(q+1)+f_i'(r)\ge0$ for $0\le r\le1$, so $\mathcal H$ preserves
this order. Since $\mathcal H(0)=0$ and $\mathcal H(1)=\beta$,
we obtain $0\le\mathcal H_i(V)\le\beta$ whenever $0\le V\le1$.
If every $V_i$ is nonincreasing, then for $z_1<z_2$ the inequalities
$V_i(z_1)\ge V_i(z_2)$ also hold at the shifted arguments
$z_1\pm1,z_2\pm1$ for every $i$.
Together with the nonnegative neighboring coefficients and
$\beta-D(q+1)+f_i'\ge0$ on $[0,1]$, they show that
$\mathcal H(V)_i(z_1)\ge\mathcal H(V)_i(z_2)$.
By \eqref{eq:profile}, a profile satisfies
\[
 (-c\partial_z+\beta)V=\mathcal H(V).
\]
For a bounded measurable $P$-periodic family $V=(V_i)_{i\in\Z}$ with
$V_i:\R\to[0,1]$, define $\mathcal T_cV$ by
\begin{equation}\label{eq:signed-front-iteration}
 (\mathcal T_cV)_i(z)
 =
 \begin{cases}
 \displaystyle
 \frac1c\int_z^\infty
 e^{-\beta(s-z)/c}\mathcal H_i(V)(s)\,\dd s,
 &c>0,\\[7pt]
 \displaystyle
 \frac1{|c|}\int_{-\infty}^z
 e^{-\beta(z-s)/|c|}\mathcal H_i(V)(s)\,\dd s,
 &c<0.
 \end{cases}
 \qquad i\in\Z,\ z\in\R.
\end{equation}
For either sign of $c$, the integration kernel in
\eqref{eq:signed-front-iteration} is nonnegative and has total mass
$1/\beta$. More generally, replacing $\mathcal H_i(V)$ by a bounded
measurable function $F:\R\to\R$ in
\eqref{eq:signed-front-iteration} defines a positive linear map on
bounded measurable functions with the supremum norm. Its operator
norm is $1/\beta$: integration gives the upper bound
$\|F\|_\infty/\beta$, and $F\equiv1$ gives equality.
If $h$ denotes the output of this linear map applied to $F$, then $h$
is bounded and locally absolutely continuous and satisfies
$-ch'+\beta h=F$ almost everywhere. Taking $F=\mathcal H_i(V)$
therefore gives
$-c\partial_z(\mathcal T_cV)_i+\beta(\mathcal T_cV)_i
=\mathcal H_i(V)$ almost everywhere.
For a bounded locally Lipschitz profile vector $V$ with values in
$[0,1]$, the inequality $\mathcal R_cV\le0$ gives
$-cV_i'+\beta V_i\ge\mathcal H_i(V)$ almost everywhere. For $c>0$, integrate this
inequality on $[z,R]$ against $e^{-\beta(s-z)/c}/c$.
The left side becomes $V_i(z)-e^{-\beta(R-z)/c}V_i(R)$.
For $c<0$, integrate on $[-R,z]$ against
$e^{-\beta(z-s)/|c|}/|c|$; the left side becomes
$V_i(z)-e^{-\beta(z+R)/|c|}V_i(-R)$.
Since $V$ is bounded, letting $R\to\infty$ gives
$V\ge\mathcal T_cV$ in both cases. Reversing the residual inequality
gives $V\le\mathcal T_cV$ when $\mathcal R_cV\ge0$. Therefore
\begin{equation}\label{eq:signed-front-barrier-order}
 \underline\Phi
 \le\mathcal T_c\underline\Phi
 \le\mathcal T_c\overline\Phi
 \le\overline\Phi.
\end{equation}
Set $V^0=\overline\Phi$ and $V^{n+1}=\mathcal T_cV^n$.
Positivity of the integral kernel and the order property of
$\mathcal H$ give
$\underline\Phi\le V^{n+1}\le V^n\le\overline\Phi$.
In \eqref{eq:signed-front-iteration}, substitute $s=z+r$ for $c>0$
and $s=z-r$ for $c<0$. In both cases the resulting kernel on
$r\ge0$ is independent of $z$, so the monotonicity property of
$\mathcal H$ shows that $\mathcal T_c$ preserves nonincreasing
dependence on $z$. Thus each $V^n$ is componentwise nonincreasing in $z$. Define
$\Phi_i(z)=\lim_{n\to\infty}V_i^n(z)$, which exists by the displayed
ordering. Continuity of $f_i$ and the finitely many shifted components
give $\mathcal H_i(V^n)(s)\to\mathcal H_i(\Phi)(s)$ for every $s$.
The bound $0\le\mathcal H_i(V^n)\le\beta$ and the integrable exponential
kernel therefore permit dominated convergence in $\mathcal T_c$,
yielding $\Phi=\mathcal T_c\Phi$. Differentiating this integral of a
bounded measurable function gives
$-c\Phi_i'+\beta\Phi_i=\mathcal H_i(\Phi)$ almost everywhere and
$|\Phi_i'|\le2\beta/|c|$. Thus $\Phi$ is locally Lipschitz,
$\mathcal H_i(\Phi)$ is continuous, and the differential identity gives
$\Phi\in C^1$ satisfying \eqref{eq:profile}.

By monotonicity, $\ell_i^\pm=\lim_{z\to\pm\infty}\Phi_i(z)$ exist.
Put $\ell^\pm=(\ell_i^\pm)_{i=0}^{P-1}$.
Integrating \eqref{eq:profile} over $[z,z+1]$ and letting
$z\to\pm\infty$ gives
\[
 D(\ell_{i-1}^\pm+q\ell_{i+1}^\pm-(q+1)\ell_i^\pm)
       +f_i(\ell_i^\pm)=0
 \qquad i\in\Z,
\]
since $\Phi_i(z+1)-\Phi_i(z)\to0$. The upper bound
$\Phi_i\le E_i$ gives $\ell^+=0$, and
$\Phi\ge\underline\Phi\not\equiv0$ makes $\ell^-\ne0$.
Lemma~\ref{lem:periodic-equilibria} therefore gives $\ell^-=1$.
For all sufficiently large $z$, the barrier bounds give
\[
 1-Ae^{-(\lambda'-\lambda)z}
       \frac{\phi_i(\lambda')}{\phi_i(\lambda)}
 \le\frac{\Phi_i(z)}{e^{-\lambda z}\phi_i(\lambda)}\le1
 \qquad i\in\Z.
\]
Since $\lambda'>\lambda$ and the phase set is finite, this quotient
tends uniformly to one as $z\to+\infty$. Hence \eqref{eq:front-tail}
holds with $A_c=1$. For any $\zeta\in\R$, the profile
$z\mapsto\Phi_i(z+\zeta)$ still solves \eqref{eq:profile}, and
$e^{\lambda z}\Phi_i(z+\zeta)\to e^{-\lambda\zeta}\phi_i(\lambda)$
as $z\to+\infty$. For a prescribed $A_c>0$, choose
$\zeta=-\lambda^{-1}\log A_c$; then
$e^{-\lambda\zeta}=A_c$, which gives \eqref{eq:front-tail}.

If $c=c_+\ne0$, choose $c_n>c_+$ with $c_n\to c_+$ and
$|c_n|\ge|c_+|/2$. Translate the profiles of speed $c_n$ so that
$\Phi_0^n(0)=1/2$. Since $0\le\Phi_i^n\le1$, \eqref{eq:profile} gives
\[
 \sup_{i\in\Z,\,z\in\R}|(\Phi_i^n)'(z)|
 \le\frac{D(q+1)+\max_{0\le i<P}\|f_i\|_{L^\infty([0,1])}}{|c_n|}\le C.
\]
The uniform derivative bound and $0\le\Phi_i^n\le1$ permit
Arzel\`a--Ascoli on every $[-m,m]$, $m\in\N$. A diagonal subsequence,
common to the finitely many phases, converges to profiles $\Phi_i$
uniformly on every bounded interval. For fixed $a<b$, integrate
\eqref{eq:profile} for this subsequence. Uniform convergence on
$[a-1,b+1]$, continuity of $f_i$ on $[0,1]$, and $c_n\to c_+$
permit passage to the limit, giving
\[
 \begin{aligned}
 -c_+\bigl(\Phi_i(b)-\Phi_i(a)\bigr)
 &=\int_a^b\bigl[D(\Phi_{i-1}(z-1)+q\Phi_{i+1}(z+1)
                      -(q+1)\Phi_i(z))\\
 &\hspace{35mm}+f_i(\Phi_i(z))\bigr]\,\dd z.
 \end{aligned}
 \qquad i\in\Z.
\]
The integrand is continuous and $c_+\ne0$, so this identity makes
$\Phi_i\in C^1$ and proves \eqref{eq:profile} with speed $c_+$.
Monotonicity and $\Phi_0(0)=1/2$ pass to the limit.
Define $\ell_i^\pm=\lim_{z\to\pm\infty}\Phi_i(z)$, which exist
by monotonicity. Integrate \eqref{eq:profile} with $c=c_+$ over
$[z,z+1]$ and let $z\to\pm\infty$. Since
$\Phi_i(z+1)-\Phi_i(z)\to0$ at either end, we obtain
\[
 D(\ell_{i-1}^\pm+q\ell_{i+1}^\pm-(q+1)\ell_i^\pm)
       +f_i(\ell_i^\pm)=0
 \qquad i\in\Z.
\]
The inequalities $\ell_0^-\ge1/2$ and $\ell_0^+\le1/2$ now allow
Lemma~\ref{lem:periodic-equilibria} to give $\ell^-=1$ and $\ell^+=0$.

Define $w_j(t)=\Phi_j(j-ct)$ for $t\in\R$ and $j\in\Z$.
Since each $\Phi_i$ is nonincreasing, put
\[
 Y_j(t)=\operatorname{sgn}(c)\,\partial_tw_j(t)
       =-|c|\Phi_j'(j-ct)\ge0.
\]
Since the constructed profiles are $C^1$, the right side of
\eqref{eq:profile} is $C^1$. Division by $c\ne0$ then gives
$\Phi_i\in C^2$. By \eqref{eq:lattice} and $0\le w_j(t)\le1$,
\[
 \|Y(t,\cdot)\|_\infty
 \le D(q+1)+\max_{0\le i<P}\|f_i\|_{L^\infty([0,1])}
 \qquad t\in\R.
\]
By differentiating \eqref{eq:lattice}, $Y$ satisfies
\[
 \dot Y_j=D(Y_{j-1}+qY_{j+1}-(q+1)Y_j)
                  +f_j'(w_j(t))Y_j
 \qquad j\in\Z,\ t\in\R.
\]
Set $K=D(q+1)+L$. The equation for $e^{K(t-s)}Y_j(t)$ has
nonnegative diagonal $L+f_j'(w_j(t))$ and off-diagonal
coefficients $D,Dq$. For fixed $l,j\in\Z$ and $s<t$, choose the
nearest-neighbor path from $l$ to $j$, of length $n=|j-l|$.
Each of its off-diagonal factors is at least $D$, and integration over
$s<r_1<\cdots<r_n<t$ contributes $(t-s)^n/n!$. Iterating the positive
integral equation along this path therefore gives
\[
 Y_j(t)\ge e^{-K(t-s)}\frac{D^n(t-s)^n}{n!}Y_l(s).
\]
If $Y_{j_0}(t_0)=0$, the path lower bound gives $Y_l(s)=0$ for every
$l\in\Z$ and $s<t_0$. Time continuity gives $Y(t_0,\cdot)=0$.
The two off-diagonal coefficients are $D,Dq$, and the diagonal
coefficient has absolute value at most $D(q+1)+L$. Gronwall's
inequality in $\ell^\infty(\Z)$ therefore gives $Y(t,\cdot)=0$ for
$t\ge t_0$.
Since $c\ne0$, the argument $j-ct$ ranges over $\R$ as $t$ varies
for each fixed index $j$. Thus every $\Phi_i$ would be constant,
contradicting its end limits. Thus $Y_j(t)>0$ and $\Phi_i'(z)<0$.

Let $w_j(t)=\Psi_j(j-ct)$ be a front with nonzero speed $c<c_+$,
where $\Psi_i\in C^1(\R;(0,1))$ are $P$-periodic in $i$ and satisfy
$\Psi_i(-\infty)=1$ and $\Psi_i(+\infty)=0$. Choose
$\max\{c,c_-\}<s<c_+$ and take $v,L_0,\delta$ from
Lemma~\ref{lem:bump} with speed $s$ and $R=1$.
Fix $0<\vartheta\le1/2$. Since $w_j(0)\to1$ as $j\to-\infty$,
choose $J$ such that $w_j(0)>1/2$ for $j\le J$, and then choose
$\ell\in P\Z$ with $\ell+L_0\le J$. The function
$v^\ell(t,j)=\vartheta v(t,j-\ell)$ is a subsolution of
\eqref{eq:lattice}: Lemma~\ref{lem:bump} gives the subsolution
inequality for $\vartheta v$, and $f_{j-\ell}=f_j$ preserves it
under the shift by $\ell\in P\Z$. Moreover,
$v^\ell(0,j)\le w_j(0)$ for every $j\in\Z$.
Indeed, $v^\ell(0,j)=0$ outside $|j-\ell|<L_0$, whereas on that
set $j\le J$ and $v^\ell(0,j)\le\vartheta\le1/2<w_j(0)$.
Lemma~\ref{lem:comparison} therefore gives
$w_j(t)\ge v^\ell(t,j)$ for every $t\ge0$ and $j\in\Z$.
For $j_t\in\Z$ with $|j_t-(\ell+st)|\le1/2$, the bound in
Lemma~\ref{lem:bump} for $R=1$ yields
$w_{j_t}(t)\ge v^\ell(t,j_t)\ge\vartheta\delta$. On the other hand,
$j_t-ct=(s-c)t+\ell+O(1)\to+\infty$ as $t\to\infty$.
Since $\Psi_i(z)\to0$ as $z\to+\infty$ uniformly over the finite
phase set, we obtain $w_{j_t}(t)\to0$ as $t\to\infty$, contradicting
$w_{j_t}(t)\ge\vartheta\delta$. This excludes every nonzero speed
$c<c_+$.\qedhere

\end{proof}

\Needspace{12\baselineskip}
\subsection{Stationary fronts}

\begin{lemma}\label{lem:stationary-fronts}
Under \eqref{eq:hyp}, there exists a stationary sequence
$p=(p_j)_{j\in\Z}$ satisfying
\begin{equation}
\label{eq:stationary-front}
 \begin{gathered}
 D(p_{j-1}+qp_{j+1}-(q+1)p_j)+f_j(p_j)=0,\\
 0<p_j<1,\qquad p_{j+P}<p_j,\\
 \lim_{j\to-\infty}p_j=1,\qquad
 \lim_{j\to+\infty}p_j=0
 \end{gathered}
\end{equation}
with the equation and inequalities holding for every $j\in\Z$,
if and only if $\mu\le0$.
\end{lemma}

\begin{proof}

\emph{The case $\mu<0$.}
By Lemma~\ref{lem:matrix}, $k$ is strictly convex,
$k(0)>0$, $k(\lambda_0)=\mu<0$, and $k(\lambda)\to+\infty$ as
$\lambda\to+\infty$. Thus $k(\lambda)=0$ has exactly two positive
roots, which we denote by $\lambda<\lambda_{\rm r}$. Choose
\[
 \lambda<\lambda'
 <
 \min\{2\lambda,\lambda_{\rm r}\}.
\]
By strict convexity, $k(\lambda')<0$.  For $A\ge1$, to be fixed below, the sequences
\[
 \overline p_j
 =
 \min\{1,e^{-\lambda j}\phi_j(\lambda)\},
 \qquad
 \underline p_j
 =
 \bigl(
 e^{-\lambda j}\phi_j(\lambda)
 -Ae^{-\lambda'j}\phi_j(\lambda')
 \bigr)_+
 \qquad j\in\Z.
\]
satisfy the stationary supersolution and subsolution inequalities,
respectively, once $A$ is chosen as follows. Put
$E_j=e^{-\lambda j}\phi_j(\lambda)$,
$\widetilde E_j=e^{-\lambda'j}\phi_j(\lambda')$ and
$\delta=-k(\lambda')>0$. Here
$(\mathcal Av)_j=D(v_{j-1}+qv_{j+1}-(q+1)v_j)$ and
$(\mathcal L_{\Z}v)_j=(\mathcal Av)_j+a_jv_j$.
The matrix equations give $\mathcal L_{\Z}E=0$ and
$\mathcal L_{\Z}\widetilde E=-\delta\widetilde E$.
Put $d=\lambda'-\lambda>0$. The condition
$E_j>A\widetilde E_j$ is equivalent to
\[
 j>\frac{\log A+\log(\phi_j(\lambda')/\phi_j(\lambda))}{d}.
\]
For every $j\in\Z$ satisfying $E_j>A\widetilde E_j$,
positivity and periodicity of $\phi(\lambda)$ and
$\phi(\lambda')$, together with $2\lambda-\lambda'>0$,
give constants $C,c_0>0$, independent of $A\ge1$ and $j$, such that
\[
 E_j\le CA^{-\lambda/d},\qquad
 \frac{A\delta\widetilde E_j}{E_j^2}
 =A\delta\frac{\phi_j(\lambda')}{\phi_j(\lambda)^2}
           e^{(2\lambda-\lambda')j}
 \ge c_0A^{\lambda/d}.
\]
Choose
$A\ge\max\{1,C^{d/\lambda},(C_2/c_0)^{d/\lambda}\}$.
Then $CA^{-\lambda/d}\le1$ and $c_0A^{\lambda/d}\ge C_2$.
At a site $j$ with $\underline p_j>0$, one has
$0<\underline p_j\le E_j\le1$. Replacing negative neighboring
values of $E-A\widetilde E$ by zero increases the diffusion term.
Using $f_j(\underline p_j)\ge a_j\underline p_j-C_2\underline p_j^2$
from \eqref{eq:taylor} and $\underline p_j\le E_j$ therefore gives
\[
 (\mathcal A\underline p)_j+f_j(\underline p_j)
 \ge A\delta\widetilde E_j-C_2E_j^2\ge0.
\]
At a zero site the reaction is zero and the neighboring contributions
are nonnegative. Where $\overline p_j=E_j<1$, truncating the neighbors
and using $f_j(E_j)\le a_jE_j$ gives
$(\mathcal A\overline p)_j+f_j(\overline p_j)\le0$.
Where $\overline p_j=1$, the reaction vanishes and the diffusion
term is nonpositive. Thus
$0\le\underline p\le\overline p\le1$ and the required stationary
inequalities hold at every site.

Choose $\beta\ge D(q+1)+L$ and define
$\mathcal S:[0,1]^\Z\to[0,1]^\Z$ by
\[
 (\mathcal Sv)_j
 =
 \frac1\beta
 \bigl[
 Dv_{j-1}+Dqv_{j+1}
 +(\beta-D(q+1))v_j+f_j(v_j)
 \bigr]
 \qquad j\in\Z.
\]
By $\beta-D(q+1)+f_j'(r)\ge0$ for $0\le r\le1$, the map
$\mathcal S$ is order preserving. Since $\mathcal S0=0$ and
$\mathcal S\mathbf1=\mathbf1$, it maps $[0,1]^\Z$ into itself.
Define $p^0=\overline p$ and $p^{n+1}=\mathcal Sp^n$ for $n\in\N_0$.
The inequalities
$\mathcal S\underline p\ge\underline p$ and
$\mathcal S\overline p\le\overline p$, together with
order preservation of $\mathcal S$, give, by induction,
\[
 \underline p\le p^{n+1}\le p^n\le\overline p
 \qquad n\in\N_0.
\]
Thus $p_j=\lim_{n\to\infty}p_j^n$ exists at every $j$.
Continuity of $f_j$ permits passage to the limit in the three
coordinate terms defining $\mathcal S$, so $p=\mathcal Sp$ and $p$
is stationary. Periodicity of $\phi$ gives
$\overline p_{j+P}\le\overline p_j$, and $\mathcal S$ commutes with a
shift by $P$. Induction and passage to the limit therefore yield
\[
 p_{j+P}\le p_j
 \qquad j\in\Z.
\]
For each $i\in\{0,\ldots,P-1\}$, define
$\ell_i^\pm=\lim_{m\to\pm\infty}p_{i+mP}$.
Since $0\le p_j\le1$ and $p_{j+P}\le p_j$, both limits exist for $0\le i<P$.
Put $\ell^\pm=(\ell_i^\pm)_{i=0}^{P-1}$, and extend $\ell_i^\pm$
periodically in $i$. Passing to the stationary equation along each
residue class shows that these periodic sequences satisfy the system
in Lemma~\ref{lem:periodic-equilibria}. That lemma gives either
$\ell_i^\pm=0$ for every $i$ or $\ell_i^\pm=1$ for every $i$,
for each choice of the sign. The bound $p_j\le e^{-\lambda j}\phi_j(\lambda)$ gives
$\ell^+=0$. Since $\lambda'>\lambda$, $\underline p_j>0$
for all sufficiently large $j$; hence $p_j>0$ at some site.
Monotonicity under a shift by $P$ then gives a positive component of
$\ell^-$, so Lemma~\ref{lem:periodic-equilibria} gives $\ell^-=1$.
Since there are only $P$ residue classes, the limits
$\ell_i^-=1$ and $\ell_i^+=0$ for $0\le i<P$ imply
$p_j\to1$ as $j\to-\infty$ and $p_j\to0$ as $j\to+\infty$.

\emph{The critical case $\mu=0$.}
Lemma~\ref{lem:matrix} gives $c_+=0$.
Choose $c_n\downarrow0$ as $n\to\infty$ and take the positive-speed fronts from
Lemma~\ref{lem:signed-fronts}, normalized by
\[
 \Phi_0^n(0)=\frac12.
\]
Let $\mathbb Q$ denote the rational numbers. Since
$\Phi_i^n(r)\in[0,1]$ for every $i\in\{0,\ldots,P-1\}$ and
$r\in\mathbb Q$, diagonal extraction gives a subsequence on which
all these values converge. Relabel this subsequence by $n$ and set
$\alpha_i(r)=\lim_{n\to\infty}\Phi_i^n(r)$.
Define
\[
 \Phi_i(z)=\sup\{\alpha_i(r):r\in\mathbb Q,\ r>z\},\quad\text{ for }z\in\R
 \qquad i\in\{0,\ldots,P-1\}.
\]
Extend $\Phi_i$ periodically by $\Phi_{i+P}(z)=\Phi_i(z)$ for
$i\in\Z$ and $z\in\R$. The functions $\alpha_i$ are nonincreasing
on $\mathbb Q$, so each $\Phi_i$ is nonincreasing and right-continuous.
For each $i\in\{0,\ldots,P-1\}$ and each continuity point $z$ of $\Phi_i$,
rational points on either side of $z$ bound $\Phi_i^n(z)$ and give
$\Phi_i^n(z)\to\Phi_i(z)$ as $n\to\infty$.
A nonincreasing function has at most countably many discontinuities.
Since $0\le\Phi_i^n\le1$, bounded convergence on each compact
interval therefore gives
\[
 \Phi_i^n\to\Phi_i
 \quad\text{a.e.\ and in }L^1_{\rm loc}(\R)\quad\text{ as } n\to\infty
 \qquad i\in\{0,\ldots,P-1\}.
\]
Translation by $-1$ or $1$ preserves Lebesgue null sets and
maps compact intervals to compact intervals. Hence
$\Phi_{i-1}^n(\,\cdot-1)\to\Phi_{i-1}(\,\cdot-1)$ and
$\Phi_{i+1}^n(\,\cdot+1)\to\Phi_{i+1}(\,\cdot+1)$
almost everywhere and in $L^1_{\rm loc}(\R)$ as $n\to\infty$.  Since every $f_i$ is Lipschitz on
$[0,1]$,
\[
 f_i(\Phi_i^n)\to f_i(\Phi_i)
 \quad\text{in }L^1_{\rm loc}(\R)\quad\text{ as } n\to\infty
 \qquad i\in\{0,\ldots,P-1\}.
\]
For $\zeta\in C_c^\infty(\R)$, integration by parts gives
\[
\begin{gathered}
 \left|
 \int_\R c_n(\Phi_i^n)'(z)\zeta(z)\,\dd z
 \right|
 =
 \left|
 c_n\int_\R\Phi_i^n(z)\zeta'(z)\,\dd z
 \right|
 \le
 c_n\|\zeta'\|_{L^1}\to0\quad\text{ as } n\to\infty
 \\\text{ for }i\in\{0,\ldots,P-1\}.
\end{gathered}
\]
The local $L^1$ convergence of the shifted profiles and reaction
terms, together with $c_n\|\zeta'\|_{L^1}\to0$, permits passage
to the limit in the weak form of \eqref{eq:profile}. Since the
limiting diffusion and reaction terms belong to $L^1_{\rm loc}(\R)$,
the resulting stationary identity holds almost everywhere:
\begin{equation}
\label{eq:zero-speed-profile-ae}
 D\bigl(
 \Phi_{i-1}(z-1)+q\Phi_{i+1}(z+1)-(q+1)\Phi_i(z)
 \bigr)+f_i(\Phi_i(z))=0
\end{equation}
for every $i\in\Z$ and almost every $z\in\R$.

For each $i\in\{0,\ldots,P-1\}$, choose a Lebesgue null set
$E_i\subseteq\R$ such that \eqref{eq:zero-speed-profile-ae} holds
at every $z\in\R\setminus E_i$. Here $\Phi_i$ denotes the
right-continuous representative fixed in its definition. Monotonicity gives limits
$\ell_i^\pm=\lim_{z\to\pm\infty}\Phi_i(z)$ for $0\le i<P$;
extend these limits periodically by $\ell_{i+P}^\pm=\ell_i^\pm$.
Choose sequences tending to either end in
$\R\setminus\bigcup_{i=0}^{P-1}E_i$.
The shifted arguments have the same end limits, and continuity of
$f_i$ gives
\[
 D(\ell_{i-1}^\pm+q\ell_{i+1}^\pm-(q+1)\ell_i^\pm)
       +f_i(\ell_i^\pm)=0
 \qquad i\in\Z.
\]
For every continuity point $z$ of $\Phi_0$, normalization and
monotonicity of $\Phi_0^n$ give
\[
 \Phi_0(z)\ge\frac12\quad\text{ for } z<0,
 \qquad
 \Phi_0(z)\le\frac12\quad\text{ for } z>0.
\]
Consequently $\ell_0^-\ge1/2$ and $\ell_0^+\le1/2$.
Applying Lemma~\ref{lem:periodic-equilibria} to
$\ell^\pm=(\ell_i^\pm)_{i=0}^{P-1}$ yields $\ell^-=1$ and
$\ell^+=0$.

For $0\le i<P$, let $D_i\subseteq\R$ be the countable set of
discontinuities of $\Phi_i$.  Define
\begin{equation}\label{eq:stationary-exceptional-set}
 \mathcal N
 =
 \bigcup_{\substack{0\le i<P\\m\in\Z}}
 \bigl((E_i-m)\cup(D_i-m)\bigr).
\end{equation}
By countable subadditivity, $\mathcal N$ has Lebesgue measure zero. Choose
$\zeta_0\in\R\setminus\mathcal N$ and define
\[
 p_j=\Phi_{j\bmod P}(j+\zeta_0)
 \qquad j\in\Z.
\]
For $0\le i<P$ and $m\in\Z$, the point $m+\zeta_0$
lies outside both $E_i$ and $D_i$. In particular, substituting
$z=j+\zeta_0$ in \eqref{eq:zero-speed-profile-ae} for phase
$j\bmod P$ gives the stationary equation for $p_j$, with
neighboring terms $p_{j-1}$ and $p_{j+1}$. The single choice of $\zeta_0$ therefore makes the stationary
equation valid at every $j\in\Z$. By monotonicity of $\Phi_i$,
$p_{j+P}\le p_j$, and the end limits give
$p_j\to1$ as $j\to-\infty$ and $p_j\to0$ as $j\to+\infty$.

\emph{Strict inequalities.}
For the sequence $p$ constructed when $\mu<0$ and for the
sequence $p_j=\Phi_{j\bmod P}(j+\zeta_0)$ constructed when
$\mu=0$, one has $0\le p_j\le1$ and $p_{j+P}\le p_j$
for every $j\in\Z$. We prove $0<p_j<1$ and $p_{j+P}<p_j$.  If $p_j=0$ at one site, its stationary equation gives
$p_{j-1}=p_{j+1}=0$, and this propagates to all sites, contradicting
the left limit.  If $p_j=1$, then $p_{j-1},p_{j+1}\le1$ and the
stationary equation forces both neighbors to equal one; propagation
contradicts the right limit.  Hence
\[
 0<p_j<1
 \quad\text{ for } j\in\Z.
\]
Set $d_j=p_j-p_{j+P}\ge0$ and define
\[
 b_j=\int_0^1 f_j'\bigl(p_{j+P}+\tau(p_j-p_{j+P})\bigr)\,\dd\tau,
 \qquad |b_j|\le L
 \qquad j\in\Z.
\]
The fundamental theorem of calculus gives
$f_j(p_j)-f_j(p_{j+P})=b_jd_j$. Subtracting the stationary equations
at $j$ and $j+P$, and using $f_{j+P}=f_j$, therefore yields
\[
 D(d_{j-1}+qd_{j+1}-(q+1)d_j)+b_jd_j=0
 \qquad j\in\Z.
\]
If $d_j=0$ at one site, the displayed
equation gives $d_{j-1}=d_{j+1}=0$, and propagation gives $d\equiv0$.
That would make $p$ $P$-periodic, which is incompatible with its two
different end limits.  Therefore
\[
 p_{j+P}<p_j
 \quad\text{ for } j\in\Z.
\]

\emph{Nonexistence when $\mu>0$.}
Suppose $p$ satisfies the stationary equation in
\eqref{eq:stationary-front}, with values in $[0,1]$ and limits one
as $j\to-\infty$ and zero as $j\to+\infty$.
By Lemma~\ref{lem:matrix}, $c_-<0<c_+$. Choose $s\in(0,c_+)$ and
take $v,L_0,\delta$ from Lemma~\ref{lem:bump} with speed $s$ and
$R=1$. Fix $0<\vartheta\le1/2$. The left limit gives $J\in\Z$
such that $p_j>1/2$ for $j\le J$; choose $\ell\in P\Z$ with
$\ell+L_0\le J$. Define $v^\ell(t,j)=\vartheta v(t,j-\ell)$.
Periodicity gives $f_{j-\ell}=f_j$, so Lemma~\ref{lem:bump} makes
$v^\ell$ a subsolution. At time zero its support lies in
$\{j:j\le J\}$ and its values are at most $\vartheta\le1/2$.
Thus $v^\ell(0,j)\le p_j$, and Lemma~\ref{lem:comparison} gives
$p_j\ge v^\ell(t,j)$ for every $t\ge0$. Choose $j_t\in\Z$ with
$|j_t-(\ell+st)|\le1/2$. By the lower bound from
Lemma~\ref{lem:bump} on $|j-st|\le1$, one has
$p_{j_t}\ge\vartheta\delta$. Since $s>0$, $j_t\to+\infty$ as
$t\to\infty$, while the right limit of $p$ gives $p_{j_t}\to0$.
This contradiction proves nonexistence when $\mu>0$.
\end{proof}

\subsection{Weighted difference estimates}

\begin{lemma}\label{lem:difference}

Assume \eqref{eq:hyp} and \eqref{eq:derivative-kpp}. Let $u,v$ be
two solutions of \eqref{eq:lattice} with values in $[0,1]$.
Let $1\le p<\infty$, $\lambda>0$, and suppose
\[
 B_p
 :=
 \|e^{\lambda j}(u_j(0)-v_j(0))\|_{\ell^p(\Z)}
 <\infty.
\]
Then, for every real $c$ and $t\ge0$,
\begin{equation}
\label{eq:difference}
 \sup_{j\in\Z}
 e^{\lambda(j-ct)}|u_j(t)-v_j(t)|
 \le
 C_{\lambda,p}(1+t)^{-1/(2p)}
 e^{(k(\lambda)-c\lambda)t}B_p.
\end{equation}

\end{lemma}

\begin{proof}

Set $z_j(t)=u_j(t)-v_j(t)$ for $j\in\Z$ and $t\ge0$.
By Lemma~\ref{lem:comparison}, $z\in C^1([0,\infty);\ell^\infty(\Z))$.
The fundamental theorem of calculus and \eqref{eq:derivative-kpp} give
\[
 \dot z_j
 =
 D(z_{j-1}+qz_{j+1}-(q+1)z_j)+b_j(t)z_j
 \qquad j\in\Z,\ t\ge0,
\]
where
\[
 b_j(t)
 =
 \int_0^1
 f_j'\bigl(v_j(t)+\tau(u_j(t)-v_j(t))\bigr)\,\dd\tau
 \le a_j
 \qquad j\in\Z,\ t\ge0.
\]
By \eqref{eq:taylor}, $|b_j(t)|\le L$. Put $\gamma=D(q+1)+L$ and
$y_j(t)=e^{\gamma t}z_j(t)$. For $w\in\ell^\infty(\Z)$, define
\[
 (K(t)w)_j=D(w_{j-1}+qw_{j+1})+(b_j(t)+L)w_j
 \qquad j\in\Z,\ t\ge0.
\]
The operators $K(t)$ are nonnegative and satisfy
$\|K(t)\|_{\ell^\infty\to\ell^\infty}\le D(q+1)+2L$.
Since $u$ and $v$ are continuous in the $\ell^\infty$ norm and the
finitely many functions $f_i'$ are uniformly continuous on $[0,1]$,
$t\mapsto K(t)$ is continuous in operator norm. Consequently, the
integral equation
$y(t,\cdot)=z(0,\cdot)+\int_0^tK(s)y(s,\cdot)\,\dd s$
has an operator-norm-convergent Volterra series. The operator in its
$n$-th term has norm at most $[t(D(q+1)+2L)]^n/n!$ and preserves
nonnegativity. Since $b_j(t)\le a_j$, each term applied to nonnegative initial data
increases when $b_j(t)$ is replaced by $a_j$. Multiplication by
$e^{-\gamma t}$ then gives domination by the semigroup of the bounded
operator
\[
 (\mathcal L_{\Z}w)_j
 =D(w_{j-1}+qw_{j+1}-(q+1)w_j)+a_jw_j
 \qquad w\in\ell^\infty(\Z),\ j\in\Z.
\]
Thus, with $z_l^0=|z_l(0)|$, positivity and this termwise comparison give
\begin{equation}\label{eq:lattice-difference-domination}
 |u_j(t)-v_j(t)|\le(e^{t\mathcal L_{\Z}}z^0)_j
 \quad\text{ for }j\in\Z,\ t\ge0.
\end{equation}

Put $H_j=e^{-\lambda j}\phi_j(\lambda)$. The entrywise identity
proved in Lemma~\ref{lem:linear-bound} gives
\[
 (e^{t\mathcal L_{\Z}})(j,l)
 =e^{k(\lambda)t}\frac{H_j}{H_l}p_t^\lambda(j,l)
 \qquad j,l\in\Z,\ t\ge0.
\]
Define
\[
 G_l=\frac{e^{\lambda l}|z_l(0)|}{\phi_l(\lambda)},
 \qquad \|G\|_p\le\frac{B_p}{\min_{0\le i<P}\phi_i(\lambda)}
 \qquad l\in\Z.
\]
For $N\in\N$, put $z^{0,N}=z^0\mathbf1_{\{-N,\ldots,N\}}$.
Applying the entrywise kernel identity to this finitely supported
sequence gives
\[
 e^{\lambda(j-ct)}(e^{t\mathcal L_{\Z}}z^{0,N})_j
 =e^{(k(\lambda)-c\lambda)t}\phi_j(\lambda)
       \sum_{|l|\le N}p_t^\lambda(j,l)G_l
 \qquad j\in\Z,\ t\ge0.
\]
Since $G\in\ell^p\subseteq\ell^\infty$ and
$\sum_{l\in\Z} p_t^\lambda(j,l)=1$, the sums over $|l|\le N$ are bounded by
$\|G\|_\infty$. Set $C=D(q+1)$. The identity
$e^{t\mathcal L_{\Z}}=e^{-Ct}e^{t(\mathcal L_{\Z}+CI)}$ represents
its action on bounded sequences by a nonnegative kernel: the diagonal
entries of $\mathcal L_{\Z}+CI$ are $a_j>0$, every power has finitely
many entries in each row, and the exponential series converges in
operator norm. The kernel has finite row sums. Applying monotone
convergence to this kernel as $N\to\infty$, and using
\eqref{eq:lattice-difference-domination}, gives
\[
 e^{\lambda(j-ct)}|z_j(t)|
 \le e^{(k(\lambda)-c\lambda)t}\phi_j(\lambda)
            \sum_{l\in\Z} p_t^\lambda(j,l)G_l
 \qquad j\in\Z,\ t\ge0.
\]

The rows of $p_t^\lambda$ have sum one, while
\eqref{eq:kernel-bound} gives
\[
 \sup_{j,l\in\Z}p_t^\lambda(j,l)
 \le C_\lambda(1+t)^{-1/2},\qquad t\ge0.
\]
If $1<p<\infty$ and $p'=p/(p-1)$, then
\[
 \sum_{l\in\Z} p_t^\lambda(j,l)^{p'}
 \le\bigl(\sup_l p_t^\lambda(j,l)\bigr)^{p'-1}
              \sum_{l\in\Z} p_t^\lambda(j,l)
 \qquad j\in\Z,\ t\ge0.
\]
Thus, also including $p=1$ with $p'=\infty$,
\[
 \|p_t^\lambda(j,\cdot)\|_{\ell^{p'}}
 \le
 C_{\lambda,p}(1+t)^{-1/(2p)}
 \qquad j\in\Z,\ t\ge0.
\]
For $p=1$ this is the row-maximum estimate. For $p>1$,
$(p'-1)/p'=1/p$, so the row bound follows with the factor
$C_\lambda^{1/p}$. Applying H\"older's inequality to
$\sum_{l\in\Z} p_t^\lambda(j,l)G_l$, then taking the supremum over $j$,
proves \eqref{eq:difference}. The constant can be chosen as
$C_{\lambda,p}=C_\lambda^{1/p}\max_i\phi_i(\lambda)/\min_i\phi_i(\lambda)$,
independently of $c$, $t$, $u$ and $v$.\qedhere

\end{proof}

\begin{lemma}\label{lem:tree-difference}

Assume \eqref{eq:hyp}, \eqref{eq:derivative-kpp}, and
$f(x,u)=f_{b(x)}(u)$. Let $u,v$ be two solutions of the differential equation in
\eqref{eq:tree-cauchy} with values in $[0,1]$. If
$1\le p<\infty$, $\lambda>0$, and
\[
 B_p
 =
 \left(
 \sum_{x\in\Tq}
 e^{p\lambda b(x)}
 |u(0,x)-v(0,x)|^p
 \right)^{1/p}
 <\infty,
\]
then, for every real $c$ and $t\ge0$,
\begin{equation}
\label{eq:tree-difference}
 \sup_{x\in\Tq}
 e^{\lambda(b(x)-ct)}|u(t,x)-v(t,x)|
 \le
 C_{\lambda,p}(1+t)^{-1/(2p)}
 e^{(k(\lambda)-c\lambda)t}B_p.
\end{equation}
The assumptions apply to arbitrary $[0,1]$-valued solutions on
$\Tq$, including solutions that vary along the sets $\{x:b(x)=j\}$.

\end{lemma}

\begin{proof}

Set $z(t,x)=u(t,x)-v(t,x)$ for $t\ge0$ and $x\in\Tq$.
By Lemma~\ref{lem:comparison},
$z\in C^1([0,\infty);\ell^\infty(\Tq))$. Define
\[
 b_{u,v}(t,x)=\int_0^1
 f_{b(x)}'\bigl(v(t,x)+\tau(u(t,x)-v(t,x))\bigr)\,\dd\tau
 \qquad t\ge0,\ x\in\Tq.
\]
Then $z_t=D\Delta_{\Tq}z+b_{u,v}(t,x)z$. The bounds
$|f_i'(s)|\le L$ in \eqref{eq:taylor} and $f_i'(s)\le a_i$ in
\eqref{eq:derivative-kpp} give, respectively,
$|b_{u,v}(t,x)|\le L$ and $b_{u,v}(t,x)\le a_{b(x)}$. Multiplication by
$e^{(D(q+1)+L)t}$ gives nonnegative diagonal coefficients
$L+b_{u,v}(t,x)$ and positive edge coefficients $D$.
For $w\in\ell^\infty(\Tq)$, define
\[
 (K_{\Tq}(t)w)(x)
 =D\sum_{y\sim x}w(y)+(L+b_{u,v}(t,x))w(x)
 \qquad t\ge0,\ x\in\Tq.
\]
Then $Y(t,x)=e^{(D(q+1)+L)t}z(t,x)$ satisfies
$\partial_tY=K_{\Tq}(t)Y$, and
$\|K_{\Tq}(t)\|_{\ell^\infty(\Tq)\to\ell^\infty(\Tq)}
\le D(q+1)+2L$. The map $t\mapsto K_{\Tq}(t)$
is continuous in operator norm
because $u,v$ are continuous in $\ell^\infty$ and the finite family
$(f_i')$ is uniformly continuous on $[0,1]$. Thus the $n$-th iterated
operator in the integrating-factor integral equation has norm at most
$[t(D(q+1)+2L)]^n/n!$, and the series converges in operator norm on
every bounded time interval. Each iterated term has nonnegative entries and sums over finitely
many nearest-neighbor or holding paths from each initial vertex.
Its row sums are bounded by its operator norm, so the series represents
the evolution on bounded sequences by a nonnegative kernel with
finite row sums. Since $b_{u,v}\le a_{b(x)}$, its
entries are at most those obtained by replacing $b_{u,v}$ with
$a_{b(x)}$. Applying the positive kernels to $z(0,\cdot)$ and
$|z(0,\cdot)|$ gives
\[
 |u(t,\cdot)-v(t,\cdot)|
 \le
 e^{t\mathcal L}|u(0,\cdot)-v(0,\cdot)|,
 \qquad
 \mathcal L=D\Delta_{\Tq}+a_{b(\cdot)}
 \qquad t\ge0.
\]
Set
\[
 H_\lambda(x)
 =
 e^{-\lambda b(x)}\phi_{b(x)}(\lambda)
 \qquad x\in\Tq.
\]
By \eqref{eq:matrix} and the one-parent, $q$-children neighbor counts,
$\mathcal LH_\lambda=k(\lambda)H_\lambda$ pointwise. For finitely
supported $\zeta:\Tq\to\C$, define
$Q_{\Tq}\zeta=H_\lambda^{-1}\mathcal L(H_\lambda\zeta)-k(\lambda)\zeta$.
The coordinate formula for this operator is
\begin{equation}
\label{eq:tree-tilted-generator}
 (Q_{\Tq}\zeta)(x)
 =
 \alpha_{b(x)}(\zeta(x^-)-\zeta(x))
 +
 \frac{\beta_{b(x)}}q
 \sum_{\substack{y\sim x\\b(y)=b(x)+1}}
 (\zeta(y)-\zeta(x))
 \qquad x\in\Tq,
\end{equation}
where
\[
 \alpha_i=De^\lambda\frac{\phi_{i-1}(\lambda)}{\phi_i(\lambda)},
 \qquad
 \beta_i=Dqe^{-\lambda}\frac{\phi_{i+1}(\lambda)}{\phi_i(\lambda)}
 \qquad i\in\Z.
\]
These rates are positive and periodic. Set
$r_* =\max_i(\alpha_i+\beta_i)$. Formula
\eqref{eq:tree-tilted-generator} defines $Q_{\Tq}$ on
$\ell^\infty(\Tq)$, with $\|Q_{\Tq}\|\le2r_*$. Take a Poisson clock
of rate $r_*$. At a ring at vertex $x$, move to $x^-$ with probability
$\alpha_{b(x)}/r_*$, to each of its $q$ children with probability
$\beta_{b(x)}/(qr_*)$, and remain at $x$ with probability
$1-(\alpha_{b(x)}+\beta_{b(x)})/r_*\ge0$.
This defines the chain with generator \eqref{eq:tree-tilted-generator};
the finite number of clock rings on every bounded interval gives
nonexplosion. Denote its transition probabilities by
$P_t^\lambda(x,y)$; see \cite[Sections~2.6--2.8]{Norris} for jump-chain
constructions and nonexplosion at bounded jump rates.

Let $Q$ be the generator on $\Z$ from
\eqref{eq:poisson-kernel}, with the same parameter $\lambda$.
For every bounded sequence $g$ on $\Z$, the neighbor counts give
\[
 Q_{\Tq}(g\circ b)=(Qg)\circ b.
\]
Induction gives $Q_{\Tq}^n(g\circ b)=(Q^ng)\circ b$ for every
$n\ge0$. Since both generators are bounded on their
$\ell^\infty$ spaces, summing their norm-convergent exponential
series gives
$e^{tQ_{\Tq}}(g\circ b)=(e^{tQ}g)\circ b$.
Taking $g=\ind_{\{l\}}$ gives the exact projection identity
\begin{equation}
\label{eq:tree-projection-kernel}
 \sum_{\substack{y\in\Tq\\b(y)=l}}P_t^\lambda(x,y)
 =
 p_t^\lambda(b(x),l).
\end{equation}
Equation~\eqref{eq:tree-projection-kernel} holds for
$x\in\Tq$, $l\in\Z$ and $t\ge0$.
By \eqref{eq:tree-projection-kernel},
$P_t^\lambda(x,y)\le p_t^\lambda(b(x),b(y))$ for every
$x,y\in\Tq$ and $t\ge0$. Together with
\eqref{eq:kernel-bound} and the row sums of $P_t^\lambda$, this gives
\[
 0\le P_t^\lambda(x,y)
 \le p_t^\lambda(b(x),b(y))
 \le C_\lambda(1+t)^{-1/2},
 \qquad
 \sum_{y\in\Tq}P_t^\lambda(x,y)=1
 \qquad x,y\in\Tq,\ t\ge0.
\]
For $1<p<\infty$, put $p'=p/(p-1)$. Using the row sum and the pointwise kernel bound, we obtain
\[
 \sum_{y\in\Tq}P_t^\lambda(x,y)^{p'}
 \le\left(\sup_yP_t^\lambda(x,y)\right)^{p'-1}
        \sum_{y\in\Tq}P_t^\lambda(x,y)
 \le\bigl[C_\lambda(1+t)^{-1/2}\bigr]^{p'-1}
 \qquad x\in\Tq,\ t\ge0.
\]
Taking the $p'$-th root and using $(p'-1)/p'=1/p$ gives
the estimate below for $1<p<\infty$. For $p=1$, the same estimate
follows from the row maximum, with $p'=\infty$:
\[
 \|P_t^\lambda(x,\cdot)\|_{\ell^{p'}(\Tq)}
 \le
 C_{\lambda,p}(1+t)^{-1/(2p)}
 \qquad x\in\Tq,\ t\ge0.
\]

The semigroup conjugation is
\begin{equation}\label{eq:tree-entry-conjugation}
 (e^{t\mathcal L})(x,y)
 =
 e^{k(\lambda)t}
 \frac{H_\lambda(x)}{H_\lambda(y)}
 P_t^\lambda(x,y)
 \qquad x,y\in\Tq,\ t\ge0.
\end{equation}
{To verify \eqref{eq:tree-entry-conjugation} entry by entry, take
$C=D(q+1)$. The diagonal entries of both $\mathcal L+CI$ and
$Q_{\Tq}+(k(\lambda)+C)I$ are then $a_{b(x)}>0$, by
$k(\lambda)-\alpha_{b(x)}-\beta_{b(x)}=a_{b(x)}-D(q+1)$.
Their off-diagonal entries are nonnegative.} For every finitely supported $\zeta:\Tq\to\C$, the defining identity
$\mathcal L(H_\lambda\zeta)=H_\lambda(Q_{\Tq}\zeta+k(\lambda)\zeta)$
and multiplication along finite paths give, for every $n\ge0$,
\[
 (\mathcal L+CI)^n(x,y)
 =\frac{H_\lambda(x)}{H_\lambda(y)}
   \bigl(Q_{\Tq}+(k(\lambda)+C)I\bigr)^n(x,y)
 \qquad x,y\in\Tq.
\]
Each entry involves only finitely many paths of length $n$.
Summing the nonnegative exponential series and cancelling $e^{Ct}$
{proves \eqref{eq:tree-entry-conjugation}. Define}
\[
 G(y)=\frac{e^{\lambda b(y)}|u(0,y)-v(0,y)|}
                  {\phi_{b(y)}(\lambda)}
 \qquad y\in\Tq.
\]
Then $\|G\|_p\le B_p/\min_i\phi_i(\lambda)$, and hence
$G\in\ell^\infty$. For $N\in\N$, apply the entrywise semigroup identity to the
finitely supported sequence
$|u(0,\cdot)-v(0,\cdot)|\mathbf1_{B(o,N)}$.
The finite-path exponential series represents $e^{t\mathcal L}$
on bounded sequences by its nonnegative kernel. Linear comparison
with the constant potential $a_{\max}$ bounds its row sums by
$e^{a_{\max}t}$. Hence monotone convergence as $N\to\infty$
identifies the limit of these semigroup values with
$e^{t\mathcal L}|u(0,\cdot)-v(0,\cdot)|$.
For the sums against $P_t^\lambda$, the bound
$\sum_{y\in\Tq}P_t^\lambda(x,y)G(y)\le\|G\|_\infty$ makes the limit finite.
Combining \eqref{eq:tree-entry-conjugation} with
$|u(t,\cdot)-v(t,\cdot)|\le
e^{t\mathcal L}|u(0,\cdot)-v(0,\cdot)|$ in the limiting kernel
representation gives
\[
 \begin{aligned}
 e^{\lambda(b(x)-ct)}|u(t,x)-v(t,x)|
 &\le e^{(k(\lambda)-c\lambda)t}\phi_{b(x)}(\lambda)
                  \sum_{y\in\Tq}P_t^\lambda(x,y)G(y).
 \end{aligned}
 \qquad x\in\Tq,\ t\ge0.
\]
Applying H\"older's inequality to the sum over $y\in\Tq$ gives
\[
 \begin{aligned}
 \sum_{y\in\Tq}P_t^\lambda(x,y)G(y)
 &\le\|P_t^\lambda(x,\cdot)\|_{\ell^{p'}(\Tq)}
       \|G\|_{\ell^p(\Tq)}\\
 &\le C_{\lambda,p}(1+t)^{-1/(2p)}
       \frac{B_p}{\min_i\phi_i(\lambda)}.
 \end{aligned}
 \qquad x\in\Tq,\ t\ge0.
\]
Using $\phi_{b(x)}(\lambda)\le\max_i\phi_i(\lambda)$ in the preceding
bound for $e^{\lambda(b(x)-ct)}|u(t,x)-v(t,x)|$, and taking the supremum
over $x\in\Tq$, proves \eqref{eq:tree-difference}.\qedhere

\end{proof}

\Needspace{30\baselineskip}
\section{Distance from a geodesic}\label{sec:geodesic}

Recall $\mathcal G=\gamma(\Z)$ and
$s(x)=\min_{n\in\Z}d(x,\gamma(n))$ from the Introduction.
Throughout this section assume \eqref{eq:hyp} and
$f(x,u)=f_{s(x)}(u)$. We identify the spectral upper endpoint of
$D\Delta_{\Tq}+a_{s(x)}$ on $\ell^2(\Tq)$ and study stationary
solutions of the differential equation in \eqref{eq:tree-cauchy}
of the form $p(x)=p_{s(x)}$, where $(p_r)_{r\ge0}$ is a sequence
on $\N_0$.
The quotient condition \eqref{eq:strict-quotient} is imposed in
Lemmas~\ref{lem:weighted-steady-uniqueness}
and~\ref{lem:cutoff-steady-uniqueness}.
Every vertex with $s(x)\ge1$ has one neighbor at level $s(x)-1$ and
$q$ neighbors at level $s(x)+1$.  A vertex on $\mathcal G$ has two
neighbors on $\mathcal G$ and $q-1$ neighbors at distance one.
Consequently, for $f(x,u)=f_{s(x)}(u)$, a solution
$u(t,x)=U_{s(x)}(t)$ satisfies
\[
 \dot U_r(t)=D(U_{r-1}(t)+qU_{r+1}(t)-(q+1)U_r(t))
             +f_r(U_r(t))
 \quad\text{for }t>0,\ r\in\N,
\]
and the boundary equation
\begin{equation}
\label{eq:geodesic-boundary}
 \dot U_0(t)=D(q-1)(U_1(t)-U_0(t))+f_0(U_0(t)),
 \quad\text{ for }t>0.
\end{equation}
In \eqref{eq:geodesic-boundary}, the two neighbors on $\mathcal G$
have the same solution value as the vertex itself. The remaining
$q-1$ neighbors therefore give the coefficient $D(q-1)$.

For each $x\in\Tq$, denote its unique closest vertex on
$\mathcal G$ by
$\pi_{\mathcal G}(x)=\gamma(\ell_{\mathcal G}(x))$,
where $\ell_{\mathcal G}(x)\in\Z$. Indeed, the unique path from $x$
to any vertex of $\mathcal G$ first meets $\mathcal G$ at
$\pi_{\mathcal G}(x)$; a second first meeting point would give a cycle.
For $R\ge0$, the set
$\mathcal N_R(\mathcal G)=\{x\in\Tq:s(x)\le R\}$ contains $\mathcal G$
and is therefore infinite.

\subsection{The spectral upper endpoint}

Recall $\lambda_0=\tfrac12\log q$. The number
$\mu=k(\lambda_0)$ is the principal eigenvalue of the matrix
$M_{\lambda_0}$ in \eqref{eq:matrix}.  To incorporate the
boundary equation \eqref{eq:geodesic-boundary}, put
\[
 \tau=D\sqrt q,\qquad
 \tau_0=D\sqrt{q-1},\qquad
 d_r=a_r-D(q+1)\text{ for }r\in\N_0.
\]
For $v=(v_r)_{r\in\N_0}\in\ell^2(\N_0)$, define the Jacobi operator
$\mathcal J_{\mathcal G}$ by
\begin{equation}
\label{eq:geodesic-jacobi}
 \begin{aligned}
 (\mathcal J_{\mathcal G}v)_0
   &=(a_0-D(q-1))v_0+\tau_0v_1,\\
 (\mathcal J_{\mathcal G}v)_1
   &=\tau_0v_0+d_1v_1+\tau v_2,\\
 (\mathcal J_{\mathcal G}v)_r
   &=\tau v_{r-1}+d_rv_r+\tau v_{r+1},
   \quad\text{ for }r\ge2.
 \end{aligned}
\end{equation}
A Jacobi operator is a symmetric tridiagonal difference operator; see
\cite[Chapters~1--2]{TeschlJacobi}. The coefficients in
\eqref{eq:geodesic-jacobi} are real and bounded, since $(a_r)$ is
periodic. Each row has at most three nonzero entries with a common bound.
Since the matrix is symmetric, its row and column sums of absolute
values are uniformly bounded. It therefore defines a bounded
symmetric operator on all of $\ell^2(\N_0)$, and hence a
self-adjoint operator.
All scalar products and sequence norms in this subsection use counting
measure. Put
\begin{equation}
\label{eq:geodesic-sigma}
 \sigma_{\mathcal G}
 =
 \sup\spec(\mathcal J_{\mathcal G}),
 \qquad
 \mathcal L_{\mathcal G}
 =
 D\Delta_{\Tq}+a_{s(\cdot)}.
\end{equation}
Since $a_{s(x)}$ is real and bounded and every vertex has degree
$q+1$, $\mathcal L_{\mathcal G}$ is also a bounded self-adjoint
operator on $\ell^2(\Tq)$.
Each layer $\{x:s(x)=r\}$ is infinite. To compare the spectral
upper endpoints on $\ell^2(\Tq)$ and $\ell^2(\N_0)$, we use the
$\ell^2$ norms on individual layers and finitely supported functions
constructed below.

\begin{lemma}\label{lem:geodesic-spectrum}

Under \eqref{eq:hyp},
\begin{equation}
\label{eq:geodesic-spectrum}
 \sup\spec_{\ell^2(\Tq)}(\mathcal L_{\mathcal G})
 =
 \sigma_{\mathcal G}\ge\mu,
 \qquad
 \ker_{\ell^2(\Tq)}(\sigma_{\mathcal G}-\mathcal L_{\mathcal G})
 =
 \{0\}.
\end{equation}

\end{lemma}

\Needspace{7\baselineskip}
\begin{proof}

\emph{Upper bound.}
For $z\in\ell^2(\Tq)$, define $F(z)=(F_r(z))_{r\ge0}$ by
\[
 F_r(z)
 =
 \left(
 \sum_{\substack{x\in\Tq\\s(x)=r}}|z(x)|^2
 \right)^{1/2},
 \text{ for }r\in\N_0,
\]
and write $F=F(z)$. The layers partition $\Tq$, so
$\|F\|_{\ell^2(\N_0)}=\|z\|_{\ell^2(\Tq)}$.
First suppose that $z$ has finite support. For $r\ge1$, every vertex
in layer $r+1$ has one predecessor in layer $r$, while a vertex in
layer $r$ has $q$ successors. Thus
\[
 \sum_{\substack{x\sim y\\s(x)=r,\ s(y)=r+1}}|z(x)|^2=qF_r^2,
 \qquad
 \sum_{\substack{x\sim y\\s(x)=r,\ s(y)=r+1}}|z(y)|^2=F_{r+1}^2,
 \quad\text{ for }r\in\N.
\]
For every $r\in\N$, Cauchy--Schwarz gives
\[
 \left|
 \sum_{\substack{x\sim y\\s(x)=r,\ s(y)=r+1}}
 \overline{z(x)}z(y)
 \right|
 \le
 \sqrt q\,F_rF_{r+1}.
\]
Between layers zero and one, Cauchy--Schwarz gives
\[
 \begin{aligned}
 \left|\sum_{\substack{x\sim y\\s(x)=0,\ s(y)=1}}
                \overline{z(x)}z(y)\right|
 &\le
 \left(\sum_{\substack{x\sim y\\s(x)=0,\ s(y)=1}}|z(x)|^2\right)^{1/2}
 \left(\sum_{\substack{x\sim y\\s(x)=0,\ s(y)=1}}|z(y)|^2\right)^{1/2}\\
 &=\sqrt{q-1}\,F_0F_1.
 \end{aligned}
\]
Here each $x\in\mathcal G$ occurs $q-1$ times in the first sum,
and each vertex with $s(y)=1$ occurs once in the second.
Along the geodesic,
\[
\begin{aligned}
 2\re\sum_{n\in\Z}
 \overline{z(\gamma(n))}z(\gamma(n+1))
 &=
 2F_0^2
 -
 \sum_{n\in\Z}
 |z(\gamma(n+1))-z(\gamma(n))|^2.
\end{aligned}
\]
The diagonal contribution to
$\langle z,\mathcal L_{\mathcal G}z\rangle$ is
$\sum_{r\ge0}d_rF_r^2$. The geodesic identity adds $2DF_0^2$ to this
term, giving $(a_0-D(q-1))F_0^2$ at $r=0$. Combining this identity
with the Cauchy--Schwarz estimates between layers $r$ and $r+1$
for $r\in\N_0$ and \eqref{eq:geodesic-jacobi} yields, for every
finitely supported $z\in\ell^2(\Tq)$ with $F=F(z)$,
\begin{equation}
\label{eq:geodesic-layer-bound}
 \langle z,\mathcal L_{\mathcal G}z\rangle
 \le
 \langle F,\mathcal J_{\mathcal G}F\rangle
 -
 D\sum_{n\in\Z}
 |z(\gamma(n+1))-z(\gamma(n))|^2.
\end{equation}
For general $z\in\ell^2(\Tq)$, approximate $z$ in $\ell^2(\Tq)$ by
finitely supported sequences. The map
$z\mapsto(z(\gamma(n+1))-z(\gamma(n)))_{n\in\Z}$ is bounded
from $\ell^2(\Tq)$ to $\ell^2(\Z)$, because
\[
 \sum_{n\in\Z}|z(\gamma(n+1))-z(\gamma(n))|^2
 \le 4\|z\|_{\ell^2(\Tq)}^2.
\]
For $z,\widetilde z\in\ell^2(\Tq)$, the reverse triangle
inequality on each layer also gives
\begin{equation}\label{eq:geodesic-layer-continuity}
 \|F(z)-F(\widetilde z)\|_{\ell^2(\N_0)}
 \le
 \|z-\widetilde z\|_{\ell^2(\Tq)}.
\end{equation}
Since $\mathcal L_{\mathcal G}$ on $\ell^2(\Tq)$ and
$\mathcal J_{\mathcal G}$ on $\ell^2(\N_0)$ are bounded,
\eqref{eq:geodesic-layer-continuity}
and continuity of the geodesic edge sum pass
\eqref{eq:geodesic-layer-bound} to the $\ell^2$ limit. In view of $\|F\|_2=\|z\|_2$ and
$\langle F,\mathcal J_{\mathcal G}F\rangle
\le\sigma_{\mathcal G}\|F\|_2^2$,
\eqref{eq:geodesic-layer-bound} gives
$\langle z,\mathcal L_{\mathcal G}z\rangle
\le\sigma_{\mathcal G}\|z\|_2^2$.
The Rayleigh characterization of the spectral upper endpoint therefore yields
\[
 \sup\spec_{\ell^2(\Tq)}(\mathcal L_{\mathcal G})
 \le
 \sigma_{\mathcal G}.
\]

\emph{Reverse bound and comparison with $\mu$.}
Set $\nu_0=1$ and $\nu_r=(q-1)q^{r-1}$ for $r\ge1$.
For every $n\in\Z$ and $r\in\N_0$,
$\nu_r=\#\{x\in\Tq:s(x)=r,\ \ell_{\mathcal G}(x)=n\}$. Let $v\ne0$ be a finitely supported nonnegative sequence
on $\N_0$ and let $w$ be a finitely supported real sequence on $\Z$ with
$\|w\|_2=1$.  Define
\[
 z(x)
 =
 w_{\ell_{\mathcal G}(x)}
 \frac{v_{s(x)}}{\sqrt{\nu_{s(x)}}},
 \quad\text{ for }x\in\Tq.
\]
The supports of $v$ and $w$ are finite, and every set counted by
$\nu_r$ is finite, so $z$ has finite support. Summing first over the
vertices counted by $\nu_r$, and then over $n\in\Z$ and
$r\in\N_0$, gives $\|z\|_2^2=
\sum_{n\in\Z}w_n^2\sum_{r\in\N_0}v_r^2=\|v\|_2^2$.
Every edge joining layers $r$ and $r+1$ has endpoints with the same
nearest point on $\mathcal G$. There are $\nu_{r+1}$ such edges over
each $\gamma(n)$, for $n\in\Z$; hence, for every $r\in\N_0$,
\[
 \sum_{\substack{x\sim y\\s(x)=r,\ s(y)=r+1}}z(x)z(y)
 =\sqrt{\frac{\nu_{r+1}}{\nu_r}}\,v_rv_{r+1}
       \sum_{n\in\Z}w_n^2.
\]
The square-root factor is $\sqrt{q-1}$ for $r=0$ and $\sqrt q$
for $r\ge1$. Along $\mathcal G$, the edge contribution is
$2Dv_0^2\sum_{n\in\Z}w_nw_{n+1}$; the identity
$2\sum_{n\in\Z}w_nw_{n+1}=2-\sum_{n\in\Z}(w_{n+1}-w_n)^2$
changes the layer-zero diagonal from $a_0-D(q+1)$ to
$a_0-D(q-1)$. Thus, for every finitely supported nonnegative
$v\ne0$ on $\N_0$ and finitely supported real $w$ on $\Z$ with
$\|w\|_2=1$, and for $z$ defined by
$z(x)=w_{\ell_{\mathcal G}(x)}v_{s(x)}/\sqrt{\nu_{s(x)}}$ for $x\in\Tq$,
\begin{equation}
\label{eq:geodesic-product-test}
 \langle z,\mathcal L_{\mathcal G}z\rangle
 =
 \langle v,\mathcal J_{\mathcal G}v\rangle
 -
 Dv_0^2
 \sum_{n\in\Z}|w_{n+1}-w_n|^2.
\end{equation}
For $M\in\N_0$, take
\[
 w_n=(2M+1)^{-1/2}\ind_{\{|n|\le M\}},
 \quad\text{ for }n\in\Z.
\]
Only the two boundary differences are nonzero, so
$\sum_{n\in\Z}|w_{n+1}-w_n|^2=2/(2M+1)$, which tends to zero as
$M\to\infty$. Since the off-diagonal entries of
$\mathcal J_{\mathcal G}$ are positive, replacing a finite test vector
by its absolute value preserves its norm and gives a quadratic form
at least as large. Moreover, finitely supported vectors are dense in
$\ell^2(\N_0)$ and $\mathcal J_{\mathcal G}$ is bounded. Its Rayleigh
supremum can therefore be taken over finitely supported nonnegative
vectors. For fixed $v$, divide \eqref{eq:geodesic-product-test} by
$\|z\|_2^2=\|v\|_2^2$ and let $M\to\infty$. Taking the supremum
over these $v$ then gives
\[
 \sup\spec_{\ell^2(\Tq)}(\mathcal L_{\mathcal G})
 \ge
 \sigma_{\mathcal G}.
\]

To prove $\sigma_{\mathcal G}\ge\mu$, recall that
$\phi(\lambda_0)=(\phi_r(\lambda_0))_{r\in\Z}$ is the positive
$P$-periodic extension of the eigenvector of $M_{\lambda_0}$,
normalized by $\sum_{i=0}^{P-1}\phi_i(\lambda_0)=1$.
Fix an integer $R\ge2$ and, for $N\in\N_0$, define
\[
 v^{R,N}_r=\phi_r(\lambda_0)\ind_{\{R\le r\le R+N\}}
 \text{ for }r\in\N_0.
\]
By \eqref{eq:matrix} at $\lambda_0$, the positive $P$-periodic vector
$\phi(\lambda_0)$ satisfies, for every $r\in\N_0$,
\[
 \tau\bigl(\phi_{r-1}(\lambda_0)+\phi_{r+1}(\lambda_0)\bigr)
 +d_r\phi_r(\lambda_0)
 =
 \mu\phi_r(\lambda_0).
\]
Multiplying by $\phi_r(\lambda_0)$ and summing over $R\le r\le R+N$
therefore gives, for every $N\in\N_0$,
\[
 \begin{aligned}
 &\langle v^{R,N},\mathcal J_{\mathcal G}v^{R,N}\rangle
       -\mu\|v^{R,N}\|_2^2\\
 &\quad=-\tau\bigl(\phi_{R-1}(\lambda_0)\phi_R(\lambda_0)
       +\phi_{R+N}(\lambda_0)\phi_{R+N+1}(\lambda_0)\bigr).
 \end{aligned}
\]
Since $\phi(\lambda_0)$ is positive and periodic,
$\|v^{R,N}\|_2^2\ge(N+1)\min_{0\le i<P}\phi_i(\lambda_0)^2$,
whereas the absolute value of the displayed right-hand side is at most
$2\tau\max_{0\le i<P}\phi_i(\lambda_0)^2$. The Rayleigh quotients of
$v^{R,N}$ consequently tend to $\mu$ as $N\to\infty$, and the
Rayleigh characterization yields $\sigma_{\mathcal G}\ge\mu$.

\emph{The kernel at the spectral upper endpoint.}
Suppose
\[
 \mathcal L_{\mathcal G}z
 =
 \sigma_{\mathcal G}z
\]
for some $z\in\ell^2(\Tq)$, and set $F=F(z)$. Since
$\mathcal J_{\mathcal G}\le\sigma_{\mathcal G}$ in the quadratic-form
sense, \eqref{eq:geodesic-layer-bound} gives
\[
\begin{aligned}
 0
 &=\langle z,(\sigma_{\mathcal G}-\mathcal L_{\mathcal G})z\rangle\\
 &\ge
 \langle F,
 (\sigma_{\mathcal G}-\mathcal J_{\mathcal G})F\rangle\\
 &\qquad+
 D\sum_{n\in\Z}
 |z(\gamma(n+1))-z(\gamma(n))|^2
 \ge0.
\end{aligned}
\]
The quadratic form of $\sigma_{\mathcal G}-\mathcal J_{\mathcal G}$
at $F$ and the geodesic edge sum are nonnegative, so both vanish.
The vanishing edge sum gives
$z(\gamma(n+1))=z(\gamma(n))$ for every $n\in\Z$.
Since $(z(\gamma(n)))_{n\in\Z}\in\ell^2(\Z)$, the common value
is zero. Hence $F_0=0$.  Also
\[
 \langle F,
 (\sigma_{\mathcal G}-\mathcal J_{\mathcal G})F\rangle=0.
\]
Since $\sigma_{\mathcal G}-\mathcal J_{\mathcal G}$ is bounded,
self-adjoint, and nonnegative, the vanishing quadratic form is
\[
 \| (\sigma_{\mathcal G}-\mathcal J_{\mathcal G})^{1/2}F\|_2^2=0.
\]
Applying $(\sigma_{\mathcal G}-\mathcal J_{\mathcal G})^{1/2}$
to $(\sigma_{\mathcal G}-\mathcal J_{\mathcal G})^{1/2}F=0$ gives
\[
 \mathcal J_{\mathcal G}F
 =
 \sigma_{\mathcal G}F.
\]
The equation at $r=0$ gives $\tau_0F_1=0$. Since $F_0=0$ and
$\tau_0,\tau>0$, the recurrence in \eqref{eq:geodesic-jacobi}
gives $F_r=0$ successively for every $r\ge0$. Hence $z=0$. This proves the kernel identity in
\eqref{eq:geodesic-spectrum}.\qedhere

\end{proof}

For each sequence index set $I$ used below and each $j\in I$, let
$\delta_j\in\ell^2(I)$ be the coordinate vector defined by
$\delta_j(r)=\mathbf1_{\{r=j\}}$ for $r\in I$.

\begin{lemma}\label{lem:geodesic-resolvent}

For $j\in\N$, define $\mathcal B_j$ on
$\ell^2(\{j,j+1,\ldots\})$ by
$(\mathcal B_jv)_r=\tau v_{r-1}+d_rv_r+\tau v_{r+1}$ for
$v\in\ell^2(\{j,j+1,\ldots\})$ and integer $r\ge j$, with $v_{j-1}=0$.
Then $\sup\spec(\mathcal B_j)=\mu$. For $j\in\N$ and $E>\mu$ put
\[
 m_j(E)
 =
 \langle\delta_j,(E-\mathcal B_j)^{-1}\delta_j\rangle
\]
and, for $E>\mu$, put
\[
 \mathfrak{s}(E)
 =
 E-a_0+D(q-1)-D^2(q-1)m_1(E).
\]
For every $j\ge1$ and $E>\mu$, one has $m_j(E)>0$ and
\begin{equation}
\label{eq:periodic-resolvent}
 m_{j+P}(E)=m_j(E),
 \qquad
 m_j(E)
 =
 \frac{1}
 {E-d_j-D^2q\,m_{j+1}(E)}.
\end{equation}
For each fixed $E>\mu$ and $j\ge1$, $m_j(E)$ is the increasing
limit of the diagonal resolvent values for the restrictions of
$\mathcal B_j$ to $\{j,\ldots,R\}$ with zero exterior values,
as $R\to\infty$. For every $j\in\N$, the function $m_j$ satisfies
\[
 Em_j(E)\to1
 \quad\text{ as } E\to\infty.
\]
Moreover, $\mathfrak{s}$ is strictly increasing and
$\mathfrak{s}(E)\to+\infty$ as $E\to\infty$.
Write $\mathfrak{s}(\mu+)=\lim_{E\downarrow\mu}\mathfrak{s}(E)\in[-\infty,\infty)$.
Then
\begin{equation}
\label{eq:geodesic-scalar-test}
 \sigma_{\mathcal G}
 =
 \begin{cases}
 \mu,&\mathfrak{s}(\mu+)\ge0,\\
 E_{\mathcal G},&\mathfrak{s}(\mu+)<0,
 \end{cases}.
\end{equation}
In the second case, $E_{\mathcal G}$ is uniquely determined by
$E_{\mathcal G}>\mu$ and $\mathfrak{s}(E_{\mathcal G})=0$. If $\mu<0$, then
\begin{equation}
\label{eq:geodesic-sign-test}
 \operatorname{sgn}\sigma_{\mathcal G}
 =
 \operatorname{sgn}\left(
 a_0-D(q-1)+D^2(q-1)m_1(0)
 \right).
\end{equation}
When $\mu=0$, the equality $\sigma_{\mathcal G}=0$ is equivalent
to $\mathfrak{s}(0+)\ge0$.

\end{lemma}

\begin{proof}

\emph{The tail resolvent.}
Put $p_r=\phi_r(\lambda_0)$ for $r\in\Z$, extend the sequence
$(d_r)_{r\ge0}$ periodically to $\Z$, and let
$\mathcal B_{\Z}$ be the full-line Jacobi operator on $\ell^2(\Z)$
with diagonal $d_r$ for $r\in\Z$ and off-diagonal entries $\tau$. The bounded real coefficients and
symmetric entries make $\mathcal B_{\Z}$ and every $\mathcal B_j$
bounded and self-adjoint on their respective $\ell^2$ spaces.
By \eqref{eq:matrix}, the positive periodic sequence $p$ satisfies,
for every $r\in\Z$,
\[
 \tau(p_{r-1}+p_{r+1})+d_rp_r=\mu p_r.
\]
For every finitely supported $v:\Z\to\C$, expansion of the squares
and the identity $\tau(p_{r-1}+p_{r+1})+d_rp_r=\mu p_r$ give
\[
 \langle v,(\mu-\mathcal B_{\Z})v\rangle
 =\tau\sum_{r\in\Z}p_rp_{r+1}
       \left|\frac{v_r}{p_r}-\frac{v_{r+1}}{p_{r+1}}\right|^2
 \ge0.
\]
For each $j\in\N$, extending a finitely supported vector on
$\{j,j+1,\ldots\}$ by zero to $\Z$ and using density gives
\[
 \sup\spec(\mathcal B_j)\le\mu.
\]
For the reverse inequality, fix $j\in\N$ and, for integers $R\ge j$, set
$v_r^{(R)}=p_r\mathbf1_{\{j,\ldots,R\}}(r)$ on
$\{j,j+1,\ldots\}$. The displayed recurrence
for $p$ gives the exact quotient, for every integer $R\ge j$,
\[
 \frac{\langle v^{(R)},\mathcal B_jv^{(R)}\rangle}
      {\|v^{(R)}\|_2^2}
 =\mu-\tau\frac{p_{j-1}p_j+p_Rp_{R+1}}
                  {\sum_{r=j}^Rp_r^2}.
\]
The numerator in the last fraction is at most $2(\max_{0\le i<P}p_i)^2$,
whereas its denominator is at least $(R-j+1)(\min_{0\le i<P}p_i)^2$.
Letting $R\to\infty$ in the Rayleigh characterization therefore gives
$\sup\spec(\mathcal B_j)\ge\mu$. Together with the upper bound, this yields, for every $j\in\N$,
\[
 \sup\spec(\mathcal B_j)=\mu.
\]

For $j\in\N$ and $E>\mu$, one has $E-\mathcal B_j\ge(E-\mu)I>0$, so the
resolvent exists. The spectral bound just proved gives
$\|e^{-Et}e^{t\mathcal B_j}\|\le e^{-(E-\mu)t}$ for $t\ge0$,
where the norm is the operator norm on $\ell^2(\{j,j+1,\ldots\})$.
Hence
\[
 (E-\mathcal B_j)^{-1}
 =
 \int_0^\infty e^{-Et}e^{t\mathcal B_j}\,\dd t
\]
{in operator norm. Take $C=D(q+1)$, so that $\mathcal B_j+CI$
has diagonal entries $a_r>0$ for $r\ge j$ and off-diagonal
entries $\tau>0$.} For integers $r,s\ge j$, the path of
length $|r-s|$ from $r$ to $s$ contributes to the $(r,s)$ entry of
$e^{t(\mathcal B_j+CI)}$ the positive factor
$\tau^{|r-s|}t^{|r-s|}/|r-s|!$ for $t>0$. Since
$e^{t\mathcal B_j}=e^{-Ct}e^{t(\mathcal B_j+CI)}$, every semigroup
entry is positive for $t>0$, and the integral gives $m_j(E)>0$.
The operators
$\mathcal B_{j+P}$ and $\mathcal B_j$ are unitarily equivalent by a
shift through one period, which gives
$m_{j+P}(E)=m_j(E)$.  For $j\in\N$ and $E>\mu$, separating the coordinate at $j$
from $\ell^2(\{j+1,j+2,\ldots\})$ gives the block decomposition
on $\C\delta_j\oplus\ell^2(\{j+1,j+2,\ldots\})$
\[
 E-\mathcal B_j=
 \begin{pmatrix}
 E-d_j&-\tau\langle\delta_{j+1},\,\cdot\,\rangle\\
 -\tau\delta_{j+1}&E-\mathcal B_{j+1}
 \end{pmatrix}.
\]
By $E-\mathcal B_j\ge(E-\mu)I$, its first-site Schur complement
$E-d_j-\tau^2m_{j+1}(E)$ is positive. Indeed, apply the quadratic
form to the vector whose first coordinate is one and whose remaining
coordinates are $\tau(E-\mathcal B_{j+1})^{-1}\delta_{j+1}$;
the resulting value is precisely this Schur complement. Taking the first diagonal entry in the block inverse formula gives
$m_j(E)^{-1}=E-d_j-D^2q\,m_{j+1}(E)$, as asserted in
\eqref{eq:periodic-resolvent}.

Fix $j\in\N$ and $E>\mu$. For integers $R\ge j$ and
$j\le r\le R$, let
$\mathcal B_r^{(R)}$ be the restriction of $\mathcal B_r$ to
$\{r,\ldots,R\}$ with zero exterior values, and define
\[
 m_r^{(R)}(E)
 =\langle\delta_r,(E-\mathcal B_r^{(R)})^{-1}\delta_r\rangle.
\]
Every matrix $E-\mathcal B_r^{(R)}$ is bounded below by
$(E-\mu)I$. Applying the first-site Schur complement successively
therefore computes these values from
\[
 m_{R+1}^{(R)}(E)=0
\]
by the recurrence in \eqref{eq:periodic-resolvent}, applied for
$r=R,R-1,\ldots,j$; every denominator
is positive. Completing the square for $E-\mathcal B_j^{(R)}$ gives
\[
 m_j^{(R)}(E)
 =
 \sup_{\substack{v\in\ell^2(\{j,j+1,\ldots\})\\
                    \supp v\subseteq\{j,\ldots,R\}}}
 \left\{
 2\re v_j-\langle v,(E-\mathcal B_j)v\rangle
 \right\}.
\]
Put $A=E-\mathcal B_j$. On the full half-line space, completing the
square gives, for every $v\in\ell^2(\{j,j+1,\ldots\})$,
\[
 2\re v_j-\langle v,Av\rangle
 =m_j(E)-\|A^{1/2}(v-A^{-1}\delta_j)\|_2^2.
\]
The subspaces of vectors supported in $\{j,\ldots,R\}$ increase
with $R$, and their union is dense in $\ell^2(\{j,j+1,\ldots\})$.
Since $A^{1/2}$ is bounded, they approximate $A^{-1}\delta_j$ in
this identity. Thus the finite suprema increase to $m_j(E)$.
The bound $A\ge(E-\mu)I$ also gives
\[
 0<m_j^{(R)}(E)\uparrow m_j(E)\le\frac1{E-\mu}
 \quad\text{ as } R\to\infty,\ R\in\N,\ R\ge j.
\]
For $E>\|\mathcal B_j\|$, the operator-norm-convergent Neumann
series is
$(E-\mathcal B_j)^{-1}=E^{-1}\sum_{n=0}^{\infty}(\mathcal B_j/E)^n$.
Taking its $(j,j)$ entry gives
\[
 m_j(E)=E^{-1}+O(E^{-2})\quad\text{ as } E\to+\infty.
\]

\emph{The boundary Schur complement.}
Fix $E>\mu$, put $R_E=(E-\mathcal B_1)^{-1}$, and write
$v=(v_0,v')\in\ell^2(\N_0)$ with respect
to the decomposition $\C\delta_0\oplus\ell^2(\{1,2,\ldots\})$.
Completing the square gives
\begin{equation}
\label{eq:geodesic-schur}
 \begin{aligned}
 \langle v,(E-\mathcal J_{\mathcal G})v\rangle
 &=
 \left\langle
 v'-\tau_0v_0R_E\delta_1,
 (E-\mathcal B_1)
 (v'-\tau_0v_0R_E\delta_1)
 \right\rangle\\
 &\quad+
 \mathfrak{s}(E)|v_0|^2.
 \end{aligned}
\end{equation}
If $\mathfrak{s}(E)\ge0$, both terms in \eqref{eq:geodesic-schur} are
nonnegative. Conversely, choose $v_0=1$ and
$v'=\tau_0R_E\delta_1$; the first term then vanishes, so nonnegativity
of $E-\mathcal J_{\mathcal G}$ forces $\mathfrak{s}(E)\ge0$.
By the Rayleigh characterization, this proves
\begin{equation}\label{eq:geodesic-schur-order}
 \sigma_{\mathcal G}\le E
 \quad\Longleftrightarrow\quad
 \mathfrak{s}(E)\ge0
 \quad\text{ for } E>\mu.
\end{equation}
The resolvent identity gives norm differentiability for $E>\mu$ and
\[
 \frac{\dd}{\dd E}m_1(E)
 =
 -\|(E-\mathcal B_1)^{-1}\delta_1\|_2^2,
 \quad\text{ for }E>\mu,
\]
and therefore
\begin{equation}\label{eq:geodesic-schur-derivative}
 \mathfrak{s}'(E)
 =
 1+D^2(q-1)
 \|(E-\mathcal B_1)^{-1}\delta_1\|_2^2
 >0,
 \quad\text{ for }E>\mu.
\end{equation}
Since $m_1(E)=O(E^{-1})$ as $E\to+\infty$, we have
$\mathfrak{s}(E)\to+\infty$ as $E\to+\infty$. By \eqref{eq:geodesic-schur-derivative}, $\mathfrak{s}$ is strictly increasing
and has a one-sided limit $\mathfrak{s}(\mu+)$ in
$[-\infty,\infty)$. If $\mathfrak{s}(\mu+)\ge0$, then $\mathfrak{s}(E)\ge0$ for every
$E>\mu$. By \eqref{eq:geodesic-schur}, this gives
$\sigma_{\mathcal G}\le E$ for every $E>\mu$. Letting
$E\downarrow\mu$ yields $\sigma_{\mathcal G}\le\mu$. Lemma~\ref{lem:geodesic-spectrum}
gives the reverse inequality, so $\sigma_{\mathcal G}=\mu$.
If $\mathfrak{s}(\mu+)<0$, continuity, strict increase and $\mathfrak{s}(E)\to+\infty$
give a unique $E_{\mathcal G}>\mu$ with $\mathfrak{s}(E_{\mathcal G})=0$.
At this value, define
$v=(1,\tau_0R_{E_{\mathcal G}}\delta_1)\in\ell^2(\N_0)$.
In \eqref{eq:geodesic-jacobi}, the tail equations follow from
$(E_{\mathcal G}-\mathcal B_1)R_{E_{\mathcal G}}\delta_1=\delta_1$,
and the zeroth equation reduces to $\mathfrak{s}(E_{\mathcal G})=0$. Hence
$\mathcal J_{\mathcal G}v=E_{\mathcal G}v$.
It follows that $\sigma_{\mathcal G}\ge E_{\mathcal G}$, while
\eqref{eq:geodesic-schur} gives $\sigma_{\mathcal G}\le E_{\mathcal G}$.
These two cases prove \eqref{eq:geodesic-scalar-test}.

If $\mu<0$, the point $E=0$ lies above the spectrum of every
$\mathcal B_j$. If $\mathfrak{s}(0)>0$, continuity gives $E\in(\mu,0)$ with
$\mathfrak{s}(E)>0$, and hence $\sigma_{\mathcal G}\le E<0$.
If $\mathfrak{s}(0)=0$, strict increase gives $\mathfrak{s}(E)<0$ for $\mu<E<0$,
so $\mathfrak{s}(\mu+)<0$. The unique root in
\eqref{eq:geodesic-scalar-test} is then zero, and
$\sigma_{\mathcal G}=0$.
If $\mathfrak{s}(0)<0$, \eqref{eq:geodesic-schur-order} at $E=0$ gives
$\sigma_{\mathcal G}>0$. Consequently, when $\mu<0$,
\[
 \sigma_{\mathcal G}<0,\ =0,\ >0
 \quad\Longleftrightarrow\quad
 \mathfrak{s}(0)>0,\ =0,\ <0,
\]
respectively.  Since
\[
 -\mathfrak{s}(0)
 =
 a_0-D(q-1)+D^2(q-1)m_1(0),
\]
this is \eqref{eq:geodesic-sign-test}.  When $\mu=0$, the equivalence
$\sigma_{\mathcal G}\le E\iff \mathfrak{s}(E)\ge0$ for $E>0$, together
with $\sigma_{\mathcal G}\ge0$, gives the assertion on letting
$E\downarrow0$.\qedhere

\end{proof}

For a constant growth rate $a_r=a$ for every $r\in\N_0$,
put $d=a-D(q+1)$.
For $E>\mu=d+2D\sqrt q$, all $m_j(E)$ are equal, and
\eqref{eq:periodic-resolvent} becomes
$D^2q\,m_1(E)^2-(E-d)m_1(E)+1=0$.
Here $E-d>2D\sqrt q$, so the two roots are distinct.
The condition $Em_1(E)\to1$ as $E\to\infty$ selects the smaller
root, and continuity of the resolvent keeps that choice throughout
$E>\mu$. Thus
\[
 m_1(E)
 =
 \frac{
 E-d-\sqrt{(E-d)^2-4D^2q}
 }{2D^2q},
 \quad\text{ for }E>\mu,
\]
where the square root is the positive real root for $E>\mu$.
Since
\[
 \mu=d+2D\sqrt q,
 \qquad
 m_1(\mu+)=\frac1{D\sqrt q},
\]
substitution in the definition of $\mathfrak{s}$ gives
\[
 \mathfrak{s}(\mu+)
 =
 D\left(\sqrt q+\frac1{\sqrt q}-2\right)>0.
\]
By \eqref{eq:geodesic-scalar-test},
\[
 \sigma_{\mathcal G}
 =
 \mu
 =
 a-\rho_q
\]
for constant $a_r=a$. In this case the potentials $a_{b(x)}$,
$a_{r(x)}$, and $a_{s(x)}$ all equal $a$, and their common spectral
upper endpoint is $a-\rho_q$.

\subsection{A uniqueness identity}

In addition to \eqref{eq:hyp},
Lemmas~\ref{lem:weighted-steady-uniqueness}
and~\ref{lem:cutoff-steady-uniqueness} assume
\eqref{eq:strict-quotient}: for every $i\in\{0,\ldots,P-1\}$, the function
$s\mapsto f_i(s)/s$ is strictly decreasing on $(0,1]$. For example, the reactions $f_i(s)=a_is(1-s)$, with $a_i>0$ for
$i\in\{0,\ldots,P-1\}$, satisfy both \eqref{eq:hyp} and
\eqref{eq:strict-quotient}. Quotient conditions also appear in the uniqueness
theory for sublinear elliptic equations \cite{BrezisOswald}.
We prove the identities for stationary solutions of the interior
equation in \eqref{eq:radial} and the boundary equation
\eqref{eq:geodesic-boundary} directly.

With the weights $\nu_0=1$ and $\nu_r=(q-1)q^{r-1}$ for
$r\ge1$, define
\[
 \mathcal H_\nu
 =
 \left\{
 v:\N_0\to\R:
 \sum_{r\ge0}\nu_rv_r^2<\infty
 \right\}.
\]
For $v,w\in\mathcal H_\nu$, put
\[
 \langle v,w\rangle_\nu=\sum_{r\ge0}\nu_rv_rw_r,
 \qquad
 \|v\|_\nu^2=\langle v,v\rangle_\nu.
\]
For every $r\in\N_0$, the weight $\nu_r$ counts the vertices
at distance $r$ whose nearest vertex on $\mathcal G$ is $\gamma(0)$:
\[
 \nu_r=\#\{x\in\Tq:s(x)=r,\
                    \pi_{\mathcal G}(x)=\gamma(0)\}.
\]
Consequently, $\sum_{r\ge0}\nu_rv_r^2$ is the sum of
$|v_{s(x)}|^2$ over vertices whose nearest geodesic vertex is
$\gamma(0)$. For every $n\in\Z$ and $r\in\N_0$, the function
$x\mapsto v_{s(x)}$ equals $v_r$ on
$\{x\in\Tq:s(x)=r,\ \pi_{\mathcal G}(x)=\gamma(n)\}$. For a real sequence $v$ on $\N_0$, define
\[
 \begin{aligned}
 (\mathcal Dv)_0
   &=D(q-1)(v_1-v_0),\\
 (\mathcal Dv)_r
   &=D(v_{r-1}+qv_{r+1}-(q+1)v_r),
   \quad\text{ for }r\ge1,
 \end{aligned}
\]
and, for $v\in\mathcal H_\nu$, set
\[
 \mathcal E(v)
 =
 \sum_{r\ge0}
 D\nu_{r+1}(v_{r+1}-v_r)^2.
\]
The map $v\mapsto(\sqrt{\nu_r}v_r)_{r\ge0}$ is an isometric
bijection from $\mathcal H_\nu$ onto $\ell^2(\N_0;\R)$. Substitution in
\eqref{eq:geodesic-jacobi} shows that it conjugates
$v\mapsto\mathcal Dv+(a_rv_r)_{r\ge0}$ to
the restriction of $\mathcal J_{\mathcal G}$ to
$\ell^2(\N_0;\R)$. Since $\mathcal J_{\mathcal G}$ is bounded
and $\sup_{r\in\N_0}|a_r|<\infty$, the multiplication map
$v\mapsto(a_rv_r)_{r\in\N_0}$ is bounded on $\mathcal H_\nu$,
and therefore so is $\mathcal D$.
For $r=0$, $\nu_1=q-1\le q\nu_0$, and for $r\ge1$,
$\nu_{r+1}=q\nu_r$. Consequently, for every $v\in\mathcal H_\nu$,
\[
 \mathcal E(v)
 \le 2D\sum_{r\ge0}\nu_{r+1}(v_{r+1}^2+v_r^2)
 \le 2D(q+1)\|v\|_\nu^2.
\]
The triangle inequality applied to
$(\sqrt{D\nu_{r+1}}(v_{r+1}-v_r))_{r\ge0}$ and the bound
$\mathcal E(v)\le2D(q+1)\|v\|_\nu^2$ show that
$v\mapsto\sqrt{\mathcal E(v)}$ is continuous on $\mathcal H_\nu$.
For finitely supported $v$, collecting the two contributions of
each edge gives $\mathcal E(v)=\langle v,-\mathcal Dv\rangle_\nu$.
For $v\in\mathcal H_\nu$ and integers $N\ge0$, define
$v^{(N)}_r=v_r\mathbf1_{\{0,\ldots,N\}}(r)$ for $r\in\N_0$.
Then $v^{(N)}\to v$ in $\mathcal H_\nu$ as $N\to\infty$. By the continuity of $\mathcal E$ and boundedness
of $\mathcal D$, passage to the limit therefore yields
\[
 \mathcal E(v)=\langle v,-\mathcal Dv\rangle_\nu
 \quad\text{ for }v\in\mathcal H_\nu.
\]

\begin{lemma}\label{lem:weighted-steady-uniqueness}
Assume \eqref{eq:hyp} and \eqref{eq:strict-quotient}.  There is at most
one sequence $U\in\mathcal H_\nu$ with $0<U_r\le1$ satisfying
\[
 (\mathcal DU)_r+f_r(U_r)=0,
 \quad\text{ for }r\ge0.
\]
\end{lemma}

\begin{proof}

Let $U=(U_r)_{r\ge0}$ satisfy $(\mathcal DU)_r+f_r(U_r)=0$
and $0<U_r\le1$ for every $r\ge0$, and let $w$ be a finitely
supported real sequence on $\N_0$. Put $b_r=D\nu_{r+1}$ for
$r\in\N_0$.
The weights satisfy $\nu_0D(q-1)=b_0$ and
$\nu_rDq=\nu_{r+1}D=b_r$ for $r\ge1$. Multiplying the stationary
equations by $\nu_r$ therefore gives
\[
 \begin{aligned}
 \nu_0f_0(U_0)&=b_0(U_0-U_1),\\
 \nu_rf_r(U_r)&=b_{r-1}(U_r-U_{r-1})
                 +b_r(U_r-U_{r+1})
                 \quad\text{ for }r\ge1.
 \end{aligned}
\]
For each edge joining $r$ to $r+1$, with $r\in\N_0$, expansion
of the squares gives
\[
 \begin{aligned}
 &b_r(w_{r+1}-w_r)^2
 -b_rU_rU_{r+1}
       \left(\frac{w_r}{U_r}-\frac{w_{r+1}}{U_{r+1}}\right)^2\\
 &\qquad=b_r(U_r-U_{r+1})
       \left(\frac{w_r^2}{U_r}-\frac{w_{r+1}^2}{U_{r+1}}\right).
 \end{aligned}
\]
For the positive stationary sequence $U$ and every finitely supported
real sequence $w$ on $\N_0$, summing the edge equality over
$r\in\N_0$ and substituting the preceding stationary equations gives
the discrete Picone identity
\begin{equation}
\label{eq:weighted-picone}
 \begin{aligned}
 &\mathcal E(w)
 -
 \sum_{r\ge0}
 \nu_r\frac{f_r(U_r)}{U_r}w_r^2\\
 &\qquad=
 \sum_{r\ge0}
 b_rU_rU_{r+1}
 \left(
 \frac{w_r}{U_r}
 -
 \frac{w_{r+1}}{U_{r+1}}
 \right)^2
 \ge0.
 \end{aligned}
\end{equation}
Every denominator in \eqref{eq:weighted-picone} is positive,
and finite support of $w$ makes both sums finite. By \eqref{eq:hyp},
we have, for every $r\in\N_0$,
\[
 0\le\frac{f_r(U_r)}{U_r}\le a_{\max},
\]
so the left side of \eqref{eq:weighted-picone} is continuous in
$w\in\mathcal H_\nu$. For $w\in\mathcal H_\nu$ and an integer $N\ge0$, put
$w_r^{(N)}=w_r\mathbf1_{\{0,\ldots,N\}}(r)$. Since
$\|w-w^{(N)}\|_\nu^2=\sum_{r>N}\nu_rw_r^2\to0$ as
$N\to\infty$, the continuity of $\mathcal E$ and boundedness of
$f_r(U_r)/U_r$ let us pass to the limit in the nonnegativity of
the left side of
\eqref{eq:weighted-picone}. We obtain
\[
 \mathcal E(w)\ge
 \sum_{r\ge0}\nu_r\frac{f_r(U_r)}{U_r}w_r^2
 \quad\text{ for }w\in\mathcal H_\nu.
\]

{Let $U,V\in\mathcal H_\nu$ satisfy
$(\mathcal DU)_r+f_r(U_r)=0$ and
$(\mathcal DV)_r+f_r(V_r)=0$, with $0<U_r,V_r\le1$,
for every $r\ge0$.} By $\mathcal E(W)=\langle W,-\mathcal DW\rangle_\nu$ for
$W\in\mathcal H_\nu$ and the stationary equations for $U,V$,
\[
 \mathcal E(V)
 =
 \sum_{r=0}^\infty\nu_rf_r(V_r)V_r,
 \qquad
 \mathcal E(U)
 =
 \sum_{r=0}^\infty\nu_rf_r(U_r)U_r.
\]
Apply $\mathcal E(w)\ge\sum_{r\ge0}\nu_r f_r(U_r)w_r^2/U_r$
with $w=V$, and then interchange $U$ and $V$ and take $w=U$.
Adding the two inequalities gives
\[
 0
 \le
 \sum_{r\ge0}
 \nu_r
 \left(
 \frac{f_r(U_r)}{U_r}
 -
 \frac{f_r(V_r)}{V_r}
 \right)
 (U_r^2-V_r^2)
 \le0.
\]
The absolute value of the $r$th summand is at most
$a_{\max}\nu_r(U_r^2+V_r^2)$, so the series converges absolutely.
By \eqref{eq:strict-quotient}, each summand is nonpositive and
vanishes exactly when $U_r=V_r$. Since the series has sum zero,
we obtain $U_r=V_r$ for every $r\ge0$.\qedhere

\end{proof}

\subsection{Uniqueness under an exponential bound}

We next prove uniqueness for positive sequences satisfying
$(\mathcal DU)_r+f_r(U_r)=0$ for $r\in\N_0$ and bounded by a
constant multiple of $q^{-r/2}\phi_r(\lambda_0)$. Put
\[
 h_r=q^{-r/2}\phi_r(\lambda_0),
 \qquad
 b_r=D\nu_{r+1},\quad\text{for }r\in\N_0.
\]
Then
\begin{equation}
\label{eq:bounded-transformed-conductance}
 b_rh_rh_{r+1}
 =
 D(q-1)q^{-1/2}
 \phi_r(\lambda_0)\phi_{r+1}(\lambda_0),
 \quad\text{for }r\in\N_0.
\end{equation}
The positive sequence $\phi_r(\lambda_0)$ is $P$-periodic.
Its positive minimum and finite maximum therefore give positive
constants $B_-,B_+$ such that
$B_-\le b_rh_rh_{r+1}\le B_+$ for every $r\in\N_0$.

\begin{lemma}\label{lem:cutoff-steady-uniqueness}

Assume \eqref{eq:hyp} and \eqref{eq:strict-quotient}.  There is at most
one sequence $U=(U_r)_{r\in\N_0}$ with $0<U_r\le1$ for every
$r\in\N_0$ satisfying
$(\mathcal DU)_r+f_r(U_r)=0$ for $r\ge0$ and
\[
 U_r\le Ch_r
 \quad\text{ for } r\ge0
\]
for some $C<\infty$.

\end{lemma}

\begin{proof}

{Let $U,V$ satisfy
$(\mathcal DU)_r+f_r(U_r)=0$ and
$(\mathcal DV)_r+f_r(V_r)=0$ for every $r\ge0$, with
$0<U_r,V_r\le1$. Choose $C<\infty$ so that
$U_r,V_r\le Ch_r$ for every $r\ge0$.} For $r\in\N_0$ and $s\in(0,1]$, put $g_r(s)=f_r(s)/s$.
For a positive sequence $W$ on $\N_0$ and a finitely supported real
sequence $\chi$ on $\N_0$, define
\[
 \mathcal R_\chi(W)
 =
 \sum_{r\ge0}
 b_rW_rW_{r+1}(\chi_{r+1}-\chi_r)^2.
\]
Apply \eqref{eq:weighted-picone} with stationary sequence $V$ and
finitely supported test sequence $w=\chi V$. Since
$w_r/V_r=\chi_r$, we obtain
\begin{equation}
\label{eq:cutoff-product-identity}
 \mathcal E(\chi V)
 =
 \sum_{r=0}^\infty\nu_rg_r(V_r)\chi_r^2V_r^2
 +
 \mathcal R_\chi(V),
\end{equation}
while \eqref{eq:weighted-picone}, with stationary sequence $U$ and
test sequence $\chi V$, gives
\[
 \mathcal E(\chi V)\ge\sum_{r=0}^\infty\nu_rg_r(U_r)\chi_r^2V_r^2 .
\]
Combining
$\mathcal E(\chi V)\ge\sum_{r=0}^\infty\nu_rg_r(U_r)\chi_r^2V_r^2$
with \eqref{eq:cutoff-product-identity} gives
$\sum_{r=0}^\infty\nu_r\chi_r^2(g_r(U_r)-g_r(V_r))V_r^2
\le\mathcal R_\chi(V)$. Applying
\eqref{eq:weighted-picone} with stationary sequence $U$ and test
sequence $\chi U$ gives
\[
 \mathcal E(\chi U)=\sum_{r=0}^\infty\nu_rg_r(U_r)\chi_r^2U_r^2
+\mathcal R_\chi(U) .
\] With stationary sequence $V$ and test sequence
$\chi U$, the nonnegative right-hand side of
\eqref{eq:weighted-picone} gives
\[
 \mathcal E(\chi U)\ge\sum_{r=0}^\infty\nu_rg_r(V_r)\chi_r^2U_r^2 .
\]
Consequently,
\[
 \begin{aligned}
 \sum_{r=0}^\infty\nu_r\chi_r^2(g_r(U_r)-g_r(V_r))V_r^2
     &\le\mathcal R_\chi(V),\\
 \sum_{r=0}^\infty\nu_r\chi_r^2(g_r(V_r)-g_r(U_r))U_r^2
     &\le\mathcal R_\chi(U).
 \end{aligned}
\]
Adding these bounds and using \eqref{eq:strict-quotient} gives
\begin{equation}
\label{eq:cutoff-uniqueness-inequality}
 \begin{aligned}
 0
 &\le
 \sum_{r=0}^\infty\nu_r\chi_r^2
 (g_r(U_r)-g_r(V_r))(V_r^2-U_r^2)\\
 &\le
 \mathcal R_\chi(U)+\mathcal R_\chi(V).
 \end{aligned}
\end{equation}
For every $r\in\N_0$, the summand
$\nu_r\chi_r^2(g_r(U_r)-g_r(V_r))(V_r^2-U_r^2)$ in
\eqref{eq:cutoff-uniqueness-inequality} is nonnegative because
$g_r$ is strictly decreasing on $(0,1]$.

For each integer $N\ge1$, define the finitely supported real sequence
$\chi=\chi^{(N)}=(\chi_r)_{r\in\N_0}$ by
\[
 \chi_r=
 \begin{cases}
 1,&0\le r\le N,\\
 2-r/N,&N<r<2N,\\
 0,&r\ge2N.
 \end{cases}
\]
The differences $\chi_{r+1}-\chi_r$ equal $-1/N$ for
$N\le r\le2N-1$ and zero for all other $r\in\N_0$. Put
\[
 C_1=2C^2D(q-1)q^{-1/2}
       \max_{0\le i<P}
       \bigl(\phi_i(\lambda_0)\phi_{i+1}(\lambda_0)\bigr).
\]
The constant $C_1$ is finite and independent of $N$. For every
integer $N\ge1$ and $\chi=\chi^{(N)}$, the bounds
$U_r,V_r\le Ch_r$ for $r\in\N_0$ and
\eqref{eq:bounded-transformed-conductance} give
\[
 \mathcal R_\chi(U)+\mathcal R_\chi(V)
 \le
 \frac{2C^2}{N^2}\sum_{r=N}^{2N-1}b_rh_rh_{r+1}
 \le\frac{C_1}{N}.
\]
For a fixed $r\in\N_0$ and $N\ge\max\{1,r\}$, one has $\chi_r=1$.
For this fixed $r$ and every integer $N\ge\max\{1,r\}$,
nonnegativity of the summands in
\eqref{eq:cutoff-uniqueness-inequality} therefore gives
\[
 0\le\nu_r(g_r(U_r)-g_r(V_r))(V_r^2-U_r^2)
 \le\frac{C_1}{N}.
\]
Letting $N\to\infty$ in the preceding pointwise bound and using
strict decrease of $g_r$ on $(0,1]$, we obtain
$U_r=V_r$. Since $r\in\N_0$ was arbitrary, $U=V$.\qedhere

\end{proof}

\begin{lemma}\label{lem:critical-stationary-tail}
Assume \eqref{eq:hyp} and $\mu=0$. Let $U=(U_r)_{r\ge0}$ satisfy
$(\mathcal DU)_r+f_r(U_r)=0$ for every $r\ge0$. Suppose that one constant $C<\infty$ satisfies, for every $r\ge0$,
\[
 0<U_r\le1,
 \qquad
 U_r\le Ch_r.
\]
Then there is a constant $A_U>0$ such that
\begin{equation}
\label{eq:critical-stationary-tail}
 U_r
 =
 A_Uq^{-r/2}\phi_r(\lambda_0)
 +
 O(q^{-r})
 \quad\text{ as } r\to\infty.
\end{equation}
Moreover, $\sum_{r=0}^{\infty}\nu_rU_r^2=\infty$.
\end{lemma}

\begin{proof}

Set $V_r=U_r/h_r$ and $B_r=b_rh_rh_{r+1}$ for $r\ge0$.
Define $B_-=\min_{0\le i<P}b_ih_ih_{i+1}$ and
$B_+=\max_{0\le i<P}b_ih_ih_{i+1}$.
By \eqref{eq:bounded-transformed-conductance}, $B_r$ is positive
and $P$-periodic. Together with $U_r\le Ch_r$, this gives $B_->0$ and
\[
 0<V_r\le C,
 \qquad B_-\le B_r\le B_+\text{ for }r\ge0.
\]
Define
\[
 J_r=B_r(V_{r+1}-V_r)\text{ for }r\ge0,
 \qquad
 F_r=\nu_rh_r\bigl(a_rU_r-f_r(U_r)\bigr)\text{ for }r\ge1.
\]
For a real sequence $W$ on $\N_0$, write $(aW)_r=a_rW_r$ for
$r\in\N_0$, so that $\mathcal D+a$ acts on real sequences on $\N_0$.
For $r\ge1$, the periodic eigenvector equation and $\mu=0$ give
$((\mathcal D+a)h)_r=0$, while stationarity gives
$((\mathcal D+a)U)_r=a_rU_r-f_r(U_r)$. Since
$\nu_rD=b_{r-1}$ and $\nu_rDq=b_r$ for $r\ge1$, expansion gives,
for every integer $r\ge1$,
\[
 \begin{aligned}
 F_r
 &=\nu_r\bigl[h_r((\mathcal D+a)U)_r
              -U_r((\mathcal D+a)h)_r\bigr]\\
 &=b_rh_rh_{r+1}(V_{r+1}-V_r)
   -b_{r-1}h_{r-1}h_r(V_r-V_{r-1}).
 \end{aligned}
\]
Consequently,
\begin{equation}
\label{eq:critical-flux}
 J_r-J_{r-1}=F_r\quad\text{ for } r\ge1,
 \qquad J_r=B_r(V_{r+1}-V_r)\quad\text{ for } r\ge0,
\end{equation}
with $J_0$ defined on the edge from $0$ to $1$. Moreover,
\eqref{eq:taylor} gives $0\le a_rU_r-f_r(U_r)\le C_2U_r^2$.
Using $U_r=h_rV_r$ and
$\nu_rh_r^3=(q-1)q^{-1}q^{-r/2}\phi_r(\lambda_0)^3$ for $r\ge1$,
we obtain
\[
 0\le F_r\le C_2\nu_rh_rU_r^2
             \le C_Fq^{-r/2}V_r^2\text{ for }r\ge1,
\]
where
$C_F=C_2(q-1)q^{-1}\max_{0\le i<P}\phi_i(\lambda_0)^3$.

By \eqref{eq:critical-flux} and $F_r\ge0$ for $r\ge1$, the
sequence $(J_r)_{r\in\N_0}$ is nondecreasing. If $J_r>0$ for some
$r\in\N_0$, then for every integer $s\ge r$,
\[
 V_{s+1}-V_s=\frac{J_s}{B_s}\ge\frac{J_r}{B_+}>0.
\]
For each integer $m>r$, summing over $s=r,\ldots,m-1$ gives
$V_m\ge V_r+(m-r)J_r/B_+\to+\infty$ as $m\to\infty$,
contrary to $V_m\le C$.
Thus $J_r\le0$ for every $r\ge0$. Put $J_\infty=\lim_{r\to\infty}J_r$;
this finite limit exists because $J_0\le J_r\le0$. If $J_\infty<0$,
the inequalities $J_r\le J_\infty<0$ and $B_r\le B_+$ give
\[
 V_{r+1}-V_r=\frac{J_r}{B_r}
 \le\frac{J_\infty}{B_+}<0\text{ for }r\ge0,
\]
which contradicts $V_r>0$ after summation. Hence
\[
 J_r\uparrow0\quad\text{ as }r\to\infty.
\]
Since $V_{r+1}-V_r=J_r/B_r\le0$, define
$A_U=\lim_{r\to\infty}V_r\ge0$; then
\[
 V_r\downarrow A_U\quad\text{ as }r\to\infty.
\]
For every $r\in\N_0$, sum \eqref{eq:critical-flux} from $r+1$
to an integer $N>r$ and let $N\to\infty$. Since $J_N\to0$, this gives
\[
 -J_r
 =
 \sum_{n=r+1}^\infty F_n.
\]
For every $r\in\N_0$, the nonnegativity of $F_n$ for $n\ge1$
permits Tonelli's theorem to interchange the sums over integers
$s\ge r$ and $n\ge s+1$. Put
$C_V=(C_F/B_-)q^{-1/2}(1-q^{-1/2})^{-2}$.
For every $r\in\N_0$, using $B_s\ge B_-$ for $s\ge r$ and
$F_n\le C_Fq^{-n/2}V_n^2$ for $n\ge s+1$ gives
\begin{equation}\label{eq:critical-relative-tail-estimate}
\begin{aligned}
 0\le V_r-A_U
 &=\sum_{s=r}^\infty\frac{-J_s}{B_s}
   =\sum_{s=r}^\infty\frac1{B_s}\sum_{n=s+1}^\infty F_n\\
 &\le B_-^{-1}\sum_{n=r+1}^\infty(n-r)F_n\\
 &\le (C_F/B_-)\sum_{n=r+1}^\infty(n-r)q^{-n/2}V_n^2
 \le C_Vq^{-r/2}V_r^2.
\end{aligned}
\end{equation}
For $r\in\N_0$, one has $V_n\le V_r$ for every integer $n\ge r$,
and the geometric series evaluates to
\[
 \sum_{n=r+1}^\infty(n-r)q^{-n/2}
 =q^{-r/2}\frac{q^{-1/2}}{(1-q^{-1/2})^2}.
\]
If $A_U=0$, dividing $V_r\le C_Vq^{-r/2}V_r^2$ by $V_r>0$ would give
\[
 1\le C_Vq^{-r/2}V_r\le C_VCq^{-r/2}\text{ for }r\ge0.
\]
The upper bound $C_VCq^{-r/2}$ tends to zero as $r\to\infty$,
so $A_U>0$.
By \eqref{eq:critical-relative-tail-estimate} and boundedness of $V$,
we conclude that
\[
 V_r=A_U+O(q^{-r/2})\quad\text{ as } r\to\infty,
\]
and multiplication by
$h_r=q^{-r/2}\phi_r(\lambda_0)$ proves
\eqref{eq:critical-stationary-tail}.

Finally,
\[
 \nu_rh_r^2
 =
 (q-1)q^{-1}\phi_r(\lambda_0)^2,
 \quad\text{ for }r\ge1.
\]
Since the periodic sequence $\phi_r(\lambda_0)$ has a positive minimum,
there is $c>0$ such that $\nu_rh_r^2\ge c$ for $r\ge1$.
The monotonicity $V_r\downarrow A_U>0$ therefore gives
\[
 \nu_rU_r^2=\nu_rh_r^2V_r^2\ge cA_U^2\text{ for }r\ge1.
\]
Summing over $r\ge1$ yields $\sum_{r=0}^\infty\nu_rU_r^2=\infty$ and hence
$U\notin\mathcal H_\nu$.
\end{proof}

\subsection{Propagation and an eight-periodic example}

\begin{lemma}\label{lem:geodesic-local}
Assume \eqref{eq:hyp}, $f(x,u)=f_{s(x)}(u)$, and $\mu>0$.
Let $u$ solve \eqref{eq:tree-cauchy} with
$u_0:\Tq\to[0,1]$ and $u_0\not\equiv0$. Then
$u(t,x)\to1$ as $t\to\infty$ for every fixed $x\in\Tq$.
In particular,
the solution $U$ of the interior equation in \eqref{eq:radial} and
the boundary equation \eqref{eq:geodesic-boundary} with initial values
\[
 U_r(0)=\eta\ind_{\{r=0\}},
 \quad\text{ for }r\in\N_0,\qquad 0<\eta\le1,
\]
satisfies
\[
 U_r(t)\to1\quad\text{ as } t\to\infty
 \qquad\text{for each fixed }r\in\N_0.
\]
\end{lemma}

\begin{proof}

Put $p_j=\phi_j(\lambda_0)$ for $j\in\Z$, and let
$p_{\min}=\min_{0\le j<P}p_j>0$ and
$p_{\max}=\max_{0\le j<P}p_j<\infty$. By periodicity,
$p_{\min}\le p_j\le p_{\max}$ for every $j\in\Z$.
For $x\in\Tq$ with $s(x)\ge1$ and $N\in\N_0$, define
\[
 \mathcal C_N(x)=
 \{y\in\Tq:d(x,y)\le N,\ s(y)=s(x)+d(x,y)\}.
\]
Write $r_0=s(x)$ and define $z:\Tq\to\R$ by
\[
 z(y)=
 \begin{cases}
 q^{-d(x,y)/2}p_{r_0+d(x,y)},&y\in\mathcal C_N(x),\\
 0,&y\notin\mathcal C_N(x),
 \end{cases}
 \quad\text{for }y\in\Tq.
\]
The cone $\mathcal C_N(x)$ has $q^m$ vertices at distance $m$ from $x$
for $0\le m\le N$, each with potential $a_{r_0+m}$,
and $q^{m+1}$ edges between distances $m$ and $m+1$ for $0\le m<N$.
Consequently, with norms and scalar products on $\ell^2(\Tq)$,
\[
 \|z\|_2^2=\sum_{m=0}^Np_{r_0+m}^2,
 \qquad
 \langle z,\mathcal L_{\mathcal G}z\rangle
 =\sum_{m=0}^N(a_{r_0+m}-D(q+1))p_{r_0+m}^2
   +2D\sqrt q\sum_{m=0}^{N-1}p_{r_0+m}p_{r_0+m+1}.
\]
Multiplying the periodic eigenvector equation
\[
 D\sqrt q(p_{j-1}+p_{j+1})+(a_j-D(q+1))p_j=\mu p_j,
 \quad\text{for }j\in\Z,
\]
by $p_j$ and summing from $j=r_0$ to $j=r_0+N$ gives
\[
 \langle z,\mathcal L_{\mathcal G}z\rangle
 =\mu\sum_{m=0}^Np_{r_0+m}^2
 -D\sqrt q\bigl(p_{r_0-1}p_{r_0}
                 +p_{r_0+N}p_{r_0+N+1}\bigr).
\]
Since $p_{\min}\le p_j\le p_{\max}$ for $j\in\Z$, for every
$x\in\Tq$ with $r_0=s(x)\ge1$ and every $N\in\N_0$ the quotient
$\langle z,\mathcal L_{\mathcal G}z\rangle/\|z\|_2^2$ is at least
\[
 \mu-\frac{2D\sqrt q(p_{\max}/p_{\min})^2}{N+1}.
\]
{Choose an integer $N\ge0$ satisfying
$N+1\ge4D\sqrt q(p_{\max}/p_{\min})^2/\mu$.
Then $\langle z,\mathcal L_{\mathcal G}z\rangle/\|z\|_2^2\ge\mu/2>0$
for every $r_0\in\N$.}

For $s(x)\ge1$, set $\Omega_x=\mathcal C_N(x)$. If $x\in\mathcal G$,
choose a neighbor $y$ with $s(y)=1$ and set
$\Omega_x=\{x\}\cup\mathcal C_N(y)$; such a neighbor exists because
$q-1\ge1$. The restriction of $\mathcal L_{\mathcal G}$ to
$\Omega_x$ with zero exterior values retains the diffusion diagonal
$-D(q+1)$. Extending the function $z$ constructed on $\mathcal C_N(y)$
by zero at $x$ preserves its $\ell^2$ norm and quadratic form.
By the Rayleigh principle, the principal eigenvalue of the zero-exterior
restriction to $\Omega_x$ is at least $\mu/2$.
For $s(x)\ge1$, the same eigenvalue bound follows from
$\langle z,\mathcal L_{\mathcal G}z\rangle/\|z\|_2^2\ge\mu/2$.

For every $x\in\Tq$, the set $\Omega_x$ is finite and connected,
contains $x$, and is contained in $B(x,N+1)$.
The matrix of the zero-exterior restriction of $\mathcal L_{\mathcal G}$
to $\ell^2(\Omega_x)$ is real symmetric, with entries $D$ along the edges.
{Adding $D(q+1)I$ to this restriction matrix gives a matrix with
nonnegative entries, irreducible because $\Omega_x$ is connected.
The Perron--Frobenius theorem therefore gives a unique positive
principal eigenvector $\varphi_x:\Omega_x\to(0,1]$} normalized by
$\max_{y\in\Omega_x}\varphi_x(y)=1$.
For fixed $N$, the matrices for $s(x)\ge1$ depend only on
$s(x)\bmod P$ under vertex relabellings that preserve the
distinguished vertex $x$. For $x\in\mathcal G$, all choices of
the neighbor $y$ give the same additional matrix under relabellings
that also preserve $x$. Thus the values $\varphi_x(x)$, with the
normalization $\max_{y\in\Omega_x}\varphi_x(y)=1$, range over at most
$P+1$ positive numbers. Their minimum equals
$m_0:=\inf_{x\in\Tq}\varphi_x(x)>0$.
The sets $\Omega_x$ and eigenvectors $\varphi_x$ satisfy
Lemma~\ref{lem:local-general} with $R_0=N+1$, $\eta_0=\mu/2$ and
$m_0$, proving $u(t,x)\to1$ as $t\to\infty$ for every fixed $x\in\Tq$.

For the initial sequence $U_r(0)=\eta\ind_{\{r=0\}}$ for $r\in\N_0$,
put $u(t,x)=U_{s(x)}(t)$ for $t\ge0$ and $x\in\Tq$.
A vertex with $s(x)\ge1$ has one neighbor at distance $s(x)-1$ and
$q$ at distance $s(x)+1$ from $\mathcal G$. A vertex on $\mathcal G$
has two neighbors on $\mathcal G$ and $q-1$ at distance one.
The neighbor counts show that $u$ solves \eqref{eq:tree-cauchy} with initial
datum $\eta\mathbf1_{\mathcal G}$. Applying the convergence just proved
for solutions of \eqref{eq:tree-cauchy} with nonzero initial data gives
$U_r(t)=u(t,x)\to1$ as $t\to\infty$ for every $r\in\N_0$ and any
fixed vertex $x\in\Tq$ with $s(x)=r$.\qedhere

\end{proof}

\begin{lemma}\label{lem:boundary-log}

Let $h$ satisfy \eqref{eq:hyp} with $P=1$ and put
\[
 a=h'(0)>\rho_q.
\]
Consider
\[
 \dot U_r(t)
 =
 D(U_{r-1}(t)+qU_{r+1}(t)-(q+1)U_r(t))+h(U_r(t)),
 \quad\text{for }r\in\N,\ t>0,
\]
with a fixed coefficient $b_0>0$ and boundary equation
\[
 \dot U_0(t)=b_0(U_1(t)-U_0(t))+h(U_0(t)),
 \quad\text{for }t>0.
\]
Assume that the initial sequence $(U_r(0))_{r\in\N_0}$ takes values in
$[0,1]$ and has finite nonempty support,
that $0\le U_r(t)\le1$ for $r\in\N_0$ and $t\ge0$, and that
\[
 U_0(t)\to1\quad\text{ as } t\to\infty.
\]
For $t>0$, put
\[
 m_a(t)
 =
 c_*(a)t-\frac{3}{2\lambda_*(a)}\log t.
\]
Then, for every $\theta\in(0,1)$, there are
$C_\theta,T_\theta>0$ such that
\[
 \{0,\ldots,\lfloor m_a(t)-C_\theta\rfloor\}
 \subseteq
 \{r\in\N_0:U_r(t)\ge\theta\}
 \subseteq
 \{0,\ldots,\lfloor m_a(t)+C_\theta\rfloor\}
\]
for $t\ge T_\theta$, with $\{0,\ldots,\lfloor R\rfloor\}=\varnothing$ when $R<0$.

\end{lemma}

\begin{proof}

For $P=1$ and $f_0=h$, formula \eqref{eq:matrix} gives
$k(\lambda)=\Psi_a(\lambda)$, so
$\mu=a-\rho_q>0$, $c_+=c_*(a)$ and
$\lambda_+=\lambda_*(a)$. Since $0\le U_r(t)\le1$, the initial
sequence has finite support and is nonzero, and $U_0(t)\to1$ as
$t\to\infty$, Lemma~\ref{lem:periodic-halfline-log} applies to this
half-line system with $f_0=h$ and the given $b_0>0$.
Fix $\theta\in(0,1)$, let $C_\theta^0,T_\theta^0>0$ be constants in
Lemma~\ref{lem:periodic-halfline-log}, and, for $t\ge0$, put
\[
 \widetilde m_a(t)=c_*(a)t-\frac{3}{2\lambda_*(a)}\log(1+t).
\]
For $t\ge1$,
\[
 0\le m_a(t)-\widetilde m_a(t)
 =\frac{3}{2\lambda_*(a)}\log(1+t^{-1})
 \le\frac{3\log2}{2\lambda_*(a)}.
\]
Choose
\[
 C_\theta=C_\theta^0+\frac{3\log2}{2\lambda_*(a)},
 \qquad T_\theta=\max\{1,T_\theta^0\}.
\]
Then $m_a(t)-C_\theta\le\widetilde m_a(t)-C_\theta^0$ and
$\widetilde m_a(t)+C_\theta^0\le m_a(t)+C_\theta$ for
$t\ge T_\theta$. Substitution in \eqref{eq:periodic-halfline-log}
gives the asserted inclusions, since
$\{r\in\N_0:r\le R\}=\{0,\ldots,\lfloor R\rfloor\}$ for every real $R$,
with the stated convention when $R<0$.\qedhere

\end{proof}

\Needspace{12\baselineskip}
\begin{lemma}\label{lem:geodesic-certificate}

Take
\[
 q=2,\qquad D=1,\qquad P=8,
\]
and
\begin{equation}
\label{eq:geodesic-example-rates}
 a_0=\frac34,
 \qquad
 a_i=\frac1{100},
 \quad\text{ for }i\in\{1,\ldots,7\}.
\end{equation}
Then $\mu\le-97/9800<0$. For the potential $a_{s(x)}$, there is a
finite connected set $\Omega\subset\Tq$ such that the restriction
of $\Delta_{\Tq}+a_{s(x)}$ to $\Omega$, with zero values on
$\Tq\setminus\Omega$, has principal eigenvalue strictly greater
than $1/25$. In particular,
$\sigma_{\mathcal G}>1/25$.

\end{lemma}

\begin{proof}

At $\lambda_0=\frac12\log2$, choose
\[
 v=(50,39,32,28,27,28,32,39)^T.
\]
By \eqref{eq:matrix},
\[
 (M_{\lambda_0}v)_i=(a_i-3)v_i+\sqrt2(v_{i-1}+v_{i+1}),
 \quad\text{ for }i\in\{0,\ldots,7\},
\]
with indices taken modulo eight. Since $(99/70)^2-2=1/4900>0$, we have $\sqrt2<99/70$.
Dividing the component formula for $(M_{\lambda_0}v)_i$ by $v_i>0$
and using \eqref{eq:geodesic-example-rates} gives, for $i\in\{0,\ldots,7\}$,
\begin{equation}\label{eq:eight-periodic-matrix-bound}
 \frac{(M_{\lambda_0}v)_i}{v_i}
 \le
 \left(
 -\frac{153}{3500},
 -\frac{149}{9100},
 -\frac{323}{11200},
 -\frac{97}{9800},
 -\frac{17}{300},
 -\frac{97}{9800},
 -\frac{323}{11200},
 -\frac{149}{9100}
 \right)_i.
\end{equation}
Every entry on the right of \eqref{eq:eight-periodic-matrix-bound}
is at most $-97/9800$. Hence
\[
 M_{\lambda_0}v
 \le
 -\frac{97}{9800}v\quad\text{componentwise in }\R^8.
\]
Let $\psi\in(0,\infty)^8$ satisfy
$\psi^TM_{\lambda_0}=\mu\psi^T$.
Multiplication of the componentwise inequality by $\psi^T$ gives
\[
 \mu\,\psi^Tv
 =
 \psi^TM_{\lambda_0}v
 \le
 -\frac{97}{9800}\psi^Tv,
\]
and hence
\[
 \mu\le-\frac{97}{9800}<0.
\]

For $1\le m\le15$, let $y_m$ be the unique neighbor of
$\gamma(m)$ outside the geodesic.  Let $\Omega$ consist of the fifteen
vertices $\gamma(1),\ldots,\gamma(15)$ and the fifteen vertices
$y_1,\ldots,y_{15}$.  It is finite and connected.  Put
\[
 \theta=\frac{\pi}{16},
 \qquad
 z(\gamma(m))=\sin(m\theta),
 \qquad
 z(y_m)=\frac13\sin(m\theta),
 \quad\text{for }m\in\{1,\ldots,15\},
\]
and extend $z:\Tq\to\R$ by zero outside $\Omega$.
The zero-exterior restriction of $\Delta_{\Tq}$ retains diagonal $-3$
at every vertex of $\Omega$, including $\gamma(1),\ldots,\gamma(15)$.
Thus $\Delta_{\Tq}+a_{s(x)}$ has diagonal $-9/4$ at $\gamma(m)$
and $-299/100$ at $y_m$ for $1\le m\le15$. The squared norm is
\[
 \|z\|_{\ell^2(\Tq)}^2=\frac{10}{9}\sum_{m=1}^{15}\sin^2(m\theta).
\]
All scalar products below are also taken on $\ell^2(\Tq)$.
To compute $\langle z,(\Delta_{\Tq}+a_{s(x)})z\rangle$, put
$s_m=\sin(m\theta)$ for $m\in\{0,\ldots,16\}$.
Since $s_0=s_{16}=0$ and
$s_{m-1}+s_{m+1}=2\cos\theta\,s_m$ for $1\le m\le15$,
summation over $1\le m\le15$ gives
\[
 2\sum_{m=1}^{14}s_ms_{m+1}=2\cos\theta\,\sum_{m=1}^{15}\sin^2(m\theta).
\]
For $1\le m\le15$, the edge $\{\gamma(m),y_m\}$ contributes
$2s_m(s_m/3)$ to $\langle z,(\Delta_{\Tq}+a_{s(x)})z\rangle$,
so the sum over these fifteen edges is $(2/3)\sum_{m=1}^{15}\sin^2(m\theta)$.
Combining the edge contributions with the diagonal values $-9/4$
and $-299/100$ gives
\begin{equation}\label{eq:eight-periodic-rayleigh}
 \frac{
 \langle z,(\Delta_{\Tq}+a_{s(x)})z\rangle
 }{\|z\|_2^2}
 =
 \frac{
 -\frac94+2\cos\theta+\frac23-\frac{299}{900}
 }{10/9}.
\end{equation}
Since
\[
 \frac{\pi}{16}<\frac15,
 \qquad
 \cos\theta
 >
 1-\frac{\theta^2}{2}
 >
 \frac{49}{50},
\]
the numerator in \eqref{eq:eight-periodic-rayleigh} satisfies
\[
 -\frac94+2\cos\theta+\frac23-\frac{299}{900}
 >-\frac94+\frac{49}{25}+\frac23-\frac{299}{900}
 =\frac2{45}.
\]
Consequently, the quotient in \eqref{eq:eight-periodic-rayleigh} is strictly greater than
\[
 \frac{2}{45}\frac9{10}
 =
 \frac1{25}.
\]
Since $z$ vanishes outside $\Omega$, its $\ell^2(\Tq)$ norm equals
the $\ell^2(\Omega)$ norm of its restriction, and
$\langle z,(\Delta_{\Tq}+a_{s(x)})z\rangle$ equals the quadratic form
of the zero-exterior restriction to $\Omega$ evaluated at $z|_\Omega$.
The Rayleigh principle therefore makes the principal eigenvalue on
$\Omega$ strictly greater than $1/25$. By
Lemma~\ref{lem:geodesic-spectrum}, which identifies
$\sigma_{\mathcal G}$ with the full-tree spectral supremum, $z$ also
gives $\sigma_{\mathcal G}>1/25$.\qedhere

\end{proof}

\subsection{Heat estimates}

\begin{lemma}\label{lem:tree-lp}

For each $r\in[1,\infty]$, define $P_t$ on $\ell^r(\Tq)$ by
\[
 P_t=e^{tD\Delta_{\Tq}},\qquad
 P_t(x,y)=(P_t\delta_y)(x),\qquad
 \delta_y(x)=\mathbf1_{\{x=y\}}.
\]
These definitions apply for $t\ge0$ and $x,y\in\Tq$.
For $1\le p<\infty$, put
\[
 \gamma_p=
 \begin{cases}
 \rho_q,&1\le p\le2,\\
 2\rho_q/p,&2\le p<\infty.
 \end{cases}
\]
Then, for every $v\in\ell^p(\Tq)$ and $t\ge0$,
\begin{equation}
\label{eq:tree-lp}
 \|P_tv\|_\infty
 \le
 C_p(1+t)^{-3/(2p)}
 e^{-\gamma_pt}
 \|v\|_p.
\end{equation}
Here the norms are taken on $\Tq$ with counting measure,
and $C_p$ depends only on $p$, $q$, and $D$.

\end{lemma}

\begin{proof}

The spectral measure of the adjacency operator on $\ell^2(\Tq)$
at a vertex determines $P_t(x,x)$. We compute this measure and
interpolate bounds for $P_t(x,\cdot)$ to prove \eqref{eq:tree-lp}.
For the corresponding tree spectral and heat-kernel theory, see
\cite{Kesten,Woess,CowlingMedaSetti}.

\emph{The spectral measure.}
Let $T_{\mathrm f}$ be the rooted tree with root $o$ in which every
vertex has $q$ children. Let $A$ and $A_{\mathrm f}$ be the adjacency
operators on $\ell^2(\Tq)$ and $\ell^2(T_{\mathrm f})$,
respectively, both with counting measure. Each operator acts by
summing the function over the neighbors of the argument in its tree.
Both trees have degree at most $q+1$, so $A$ and $A_{\mathrm f}$
are bounded and self-adjoint on their respective spaces. In particular,
$\|A\|,\|A_{\mathrm f}\|\le q+1$.
For real $z>q+1$, define
\[
 F(z)=\langle\delta_o,(z-A_{\mathrm f})^{-1}\delta_o\rangle.
\]
Deleting $o$ leaves $q$ copies of $T_{\mathrm f}$. For
$h=(z-A_{\mathrm f})^{-1}\delta_o$, the equation on each copy has
forcing $h(o)$ at its root and zero elsewhere. After identifying
the copy with $T_{\mathrm f}$, applying $(z-A_{\mathrm f})^{-1}$
gives the value $F(z)h(o)$ at the root of the copy. Since $h(o)=F(z)$, the
equation at $o$ is $zF(z)-qF(z)^2=1$. Thus, for real $z>q+1$,
\[
 F(z)=\frac1{z-qF(z)}.
\]
The resolvent expansion gives $zF(z)\to1$ as $z\to+\infty$
on the real axis, which selects, for real $z>q+1$,
\[
 F(z)
 =
 \frac{z-\sqrt{z^2-4q}}{2q}.
\]
For $x\in\Tq$ and real $z>q+1$, deleting $x$ leaves $q+1$
copies of $T_{\mathrm f}$. Solving the equation
$(z-A)h=\delta_x$ on these copies and then at $x$ gives
\begin{equation}\label{eq:tree-diagonal-resolvent}
 \begin{aligned}
 G(z)&=\langle\delta_x,(z-A)^{-1}\delta_x\rangle
       =\frac1{z-(q+1)F(z)}\\
     &=\frac{2q}{(q-1)z+(q+1)\sqrt{z^2-4q}}.
 \end{aligned}
\end{equation}
Let $\varrho_q$ denote the spectral measure of $A$ at $\delta_x$,
characterized, for real $z>q+1$, by
\[
 G(z)=\int_\R\frac{\dd\varrho_q(\xi)}{z-\xi}.
\]
Vertex transitivity of $\Tq$ makes $\varrho_q$ independent of $x$.
Choose the branch of $\sqrt{z^2-4q}$ on
$\C\setminus[-2\sqrt q,2\sqrt q]$ for which
$\sqrt{z^2-4q}/z\to1$ as $|z|\to\infty$.
The functions
$\langle\delta_o,(z-A_{\mathrm f})^{-1}\delta_o\rangle$ and
$\langle\delta_x,(z-A)^{-1}\delta_x\rangle$, and their
algebraic expressions above, are analytic on the connected open set
$\{z:\im z>0\}\cup\{z:\re z>q+1\}$.
Their agreement for real $z>q+1$ and the identity theorem therefore
give the formulas for $F(z)$ and $G(z)$ throughout the upper
half-plane. For $-2\sqrt q<\xi<2\sqrt q$, the chosen square root
has boundary value
$\lim_{\varepsilon\downarrow0}\sqrt{(\xi+\mathrm i\varepsilon)^2-4q}
=\mathrm i\sqrt{4q-\xi^2}$. Writing
$G(\xi+\mathrm i0)=\lim_{\varepsilon\downarrow0}
G(\xi+\mathrm i\varepsilon)$, we obtain, for $-2\sqrt q<\xi<2\sqrt q$,
\[
 -\im G(\xi+i0)
 =\frac{(q+1)\sqrt{4q-\xi^2}}
 {2((q+1)^2-\xi^2)}.
\]
Outside $[-2\sqrt q,2\sqrt q]$, the two terms in the denominator
of $G$ have the same sign, so $G$ has finite real boundary values.
At the endpoints, $G(\pm2\sqrt q)=\pm\sqrt q/(q-1)$ is finite
because $q\ge2$. Thus $\lim_{\varepsilon\downarrow0}\im[-G(\xi+\mathrm i\varepsilon)]$
exists and is finite for every $\xi\in\R$. The function $-G$ has
positive imaginary part in the upper half-plane, so
\cite[Lemma~B.7]{TeschlJacobi} shows that $\varrho_q$ is absolutely
continuous. Stieltjes inversion in the
sign convention of \cite[Theorem~B.2]{TeschlJacobi} therefore gives
on $[-2\sqrt q,2\sqrt q]$, with zero measure outside this interval,
\begin{equation}\label{eq:tree-kesten-density}
 \dd\varrho_q(\xi)
 =
 \frac{q+1}{2\pi}
 \frac{\sqrt{4q-\xi^2}}
 {(q+1)^2-\xi^2}
 \ind_{[-2\sqrt q,\,2\sqrt q]}(\xi)\,\dd\xi.
\end{equation}
This is the Kesten measure for the adjacency operator.

\emph{Heat-kernel row bounds.}
In the estimates below, $C$ depends only on $q$ and $D$.
By $D\Delta_{\Tq}=D(A-(q+1)I)$ and the spectral theorem,
$P_t(x,x)$ is an integral against $\varrho_q$ in
\eqref{eq:tree-kesten-density}. On $\supp\varrho_q
=[-2\sqrt q,2\sqrt q]$,
$ (q+1)^2-\xi^2\ge(q-1)^2$. Since
\[
 \rho_q=D(q+1-2\sqrt q),
\]
we obtain, for every $x\in\Tq$ and $t\ge0$,
\begin{equation}\label{eq:tree-diagonal-heat}
\begin{aligned}
 P_t(x,x)
 &=
 \int_{-2\sqrt q}^{2\sqrt q}
 e^{-D(q+1-\xi)t}\,\dd\varrho_q(\xi)\\
 &\le C e^{-\rho_qt}
       \int_0^{4\sqrt q}e^{-Dtu}\sqrt u\,\dd u\\
 &\le C(1+t)^{-3/2}e^{-\rho_qt}.
\end{aligned}
\end{equation}
Here $u=2\sqrt q-\xi\in[0,4\sqrt q]$ for
$\xi\in[-2\sqrt q,2\sqrt q]$, so
$\sqrt{4q-\xi^2}=\sqrt{u(4\sqrt q-u)}\le2q^{1/4}\sqrt u$.
For $t\ge1$, the substitution $v=Dtu$ bounds the integral by
$(Dt)^{-3/2}\int_0^\infty e^{-v}\sqrt v\,\dd v$.
For $0\le t\le1$, the integral is bounded by
$\int_0^{4\sqrt q}\sqrt u\,\dd u$. These bounds give the factor
$(1+t)^{-3/2}$ after increasing $C$.

For $t\ge0$, the adjacency matrix $A$ is real and symmetric, so
$P_{t/2}(x,y)=P_{t/2}(y,x)\in\R$ for $x,y\in\Tq$.
Since $P_{t/2}$ is bounded on $\ell^2(\Tq)$,
$P_{t/2}(x,\cdot)=P_{t/2}\delta_x$ belongs to
$\ell^2(\Tq)$ for every $x\in\Tq$. The semigroup identity at time
$t/2$, symmetry, and Cauchy--Schwarz therefore give, for
$x,y\in\Tq$ and $t\ge0$,
\[
 \begin{aligned}
 |P_t(x,y)|
 &=\left|\sum_{w\in\Tq}P_{t/2}(x,w)P_{t/2}(y,w)\right|\\
 &\le\sqrt{P_t(x,x)P_t(y,y)}
 \le C(1+t)^{-3/2}e^{-\rho_qt},
 \end{aligned}
\]
and, for $x\in\Tq$ and $t\ge0$,
\[
 \|P_t(x,\cdot)\|_2^2
 =\sum_{y\in\Tq}P_t(x,y)^2
 =P_{2t}(x,x)
 \le
 C(1+t)^{-3/2}e^{-2\rho_qt}.
\]
Thus, for $x\in\Tq$ and $t\ge0$,
\begin{equation}\label{eq:tree-heat-row-endpoints}
 \|P_t(x,\cdot)\|_\infty
 \le
 C(1+t)^{-3/2}e^{-\rho_qt},
 \qquad
 \|P_t(x,\cdot)\|_2
 \le
 C(1+t)^{-3/4}e^{-\rho_qt}.
\end{equation}
For $1\le r<\infty$, $v\in\ell^r(\Tq)$ and $x\in\Tq$,
convexity gives
$|Av(x)|^r\le(q+1)^{r-1}\sum_{y\sim x}|v(y)|^r$.
Summing over $x$, each $y$ occurs $q+1$ times, so
$\|Av\|_r\le(q+1)\|v\|_r$. The same bound for $r=\infty$
follows by taking the supremum in the neighbor sum. Thus the
exponential series converges in operator norm on each $\ell^r(\Tq)$,
$1\le r\le\infty$, and gives, for $t\ge0$,
\[
 P_t=e^{-D(q+1)t}\sum_{n=0}^{\infty}\frac{(Dt)^n}{n!}A^n.
\]
For every integer $n\ge0$, $A^n$ has nonnegative matrix entries and
\[
 \sum_{y\in\Tq}A^n(x,y)=(q+1)^n
 \quad\text{for }x\in\Tq,
\]
because $A\mathbf1=(q+1)\mathbf1$ for the function
$\mathbf1(x)=1$ on $\Tq$. The coordinate entries of the exponential
series agree on every $\ell^r(\Tq)$, $1\le r\le\infty$.
Tonelli's theorem therefore gives
$\sum_{y\in\Tq}P_t(x,y)
=e^{-D(q+1)t}\sum_{n=0}^\infty(Dt)^n(q+1)^n/n!=1$
for $x\in\Tq$ and $t\ge0$.
Consequently, $\|P_t(x,\cdot)\|_1=1$ for $x\in\Tq$ and $t\ge0$.

\emph{Interpolation.}
Fix $1\le p<\infty$ and let $1/p+1/p'=1$, with $p'=\infty$ when $p=1$.
If $1\le p\le2$, then $p'\in[2,\infty]$. Set
\[
 \frac1{p'}=\frac{1-\vartheta}{2},
 \qquad
 \vartheta=\frac2p-1.
\]
Interpolating the two estimates in \eqref{eq:tree-heat-row-endpoints}
gives, for $x\in\Tq$, $t\ge0$ and $1\le p\le2$,
\[
 \begin{aligned}
 \|P_t(x,\cdot)\|_{p'}
 &\le\|P_t(x,\cdot)\|_2^{1-\vartheta}
       \|P_t(x,\cdot)\|_\infty^{\vartheta}\\
 &\le C_p(1+t)^{-3/(2p)}e^{-\rho_qt}.
 \end{aligned}
\]
If $2\le p<\infty$, then $1\le p'\le2$, and interpolation between
$\|P_t(x,\cdot)\|_1=1$ and the $\ell^2$ bound in
\eqref{eq:tree-heat-row-endpoints} yields, for $x\in\Tq$ and $t\ge0$,
\[
 \begin{aligned}
 \|P_t(x,\cdot)\|_{p'}
 &\le\|P_t(x,\cdot)\|_1^{1-2/p}
       \|P_t(x,\cdot)\|_2^{2/p}\\
 &\le C_p(1+t)^{-3/(2p)}e^{-2\rho_qt/p}.
 \end{aligned}
\]
For every finitely supported $v\in\ell^p(\Tq)$, $1\le p<\infty$,
$x\in\Tq$ and $t\ge0$, H\"older's inequality therefore gives
\[
 |P_tv(x)|
 \le\sum_{y\in\Tq}P_t(x,y)|v(y)|
 \le\|P_t(x,\cdot)\|_{p'}\|v\|_p.
\]
Taking the supremum over $x\in\Tq$ gives \eqref{eq:tree-lp} for
finitely supported $v$. Fix $1\le p<\infty$, $t\ge0$ and
$v\in\ell^p(\Tq)$, and choose finitely supported $v^n$ with
$\|v^n-v\|_p\to0$ as $n\to\infty$.
Since $P_t$ is bounded on $\ell^p(\Tq)$ and
$\|w\|_\infty\le\|w\|_p$ for $w\in\ell^p(\Tq)$,
$P_tv^n\to P_tv$ uniformly on $\Tq$ as $n\to\infty$.
Passing to the limit in \eqref{eq:tree-lp} for $v^n$ proves
\eqref{eq:tree-lp} for $v$.\qedhere

\end{proof}
\section{Proofs of the main results}\label{sec:proofs-main}

Under \eqref{eq:hyp}, for nonzero finitely supported initial data
$u_0:\Tq\to[0,1]$, extinction and persistence under either assignment
$f(x,u)=f_{b(x)}(u)$ or $f(x,u)=f_{r(x)}(u)$ are determined by the sign of
\[
 \mu=k\left(\frac12\log q\right),
\]
where $k(\lambda)$ is the principal eigenvalue of $M_\lambda$ in
\eqref{eq:matrix}, with $a_j=f_j'(0)$.
For $f(x,u)=f_{s(x)}(u)$, the threshold
is the upper spectral endpoint $\sigma_{\mathcal G}$ of the
self-adjoint half-line operator $\mathcal J_{\mathcal G}$ in \eqref{eq:geodesic-jacobi} on
$\ell^2(\N_0)$. The quantity $c_+$ in \eqref{eq:parameters} is the
signed propagation speed in \eqref{eq:lattice}; by
Lemma~\ref{lem:matrix}, it has the same sign as $\mu$.

\Needspace{12\baselineskip}
\subsection{Extinction and spreading}

\Needspace{7\baselineskip}
\begin{proof}[Proof of Theorem~\ref{thm:threshold}]

Set $\mathcal L=D\Delta_{\Tq}+a_{b(\cdot)}$ on
$\ell^2(\Tq)$, with counting measure. This operator is bounded and
self-adjoint, and Lemma~\ref{lem:horo-spectrum} gives
$\sup\spec(\mathcal L)=\mu$. The same coordinate formula defines a bounded operator on
$\ell^\infty(\Tq)$. For data in $\ell^2(\Tq)$, its exponential
series in the two spaces agree at every vertex. By the inequality
$f_{b(x)}(s)\le a_{b(x)}s$ in \eqref{eq:hyp} and
Lemma~\ref{lem:comparison},
\begin{equation}\label{eq:threshold-semigroup-comparison}
 0\le u(t,x)\le(e^{t\mathcal L}u_0)(x)
 \quad\text{ for }t\ge0\text{ and }x\in\Tq.
\end{equation}
If $\mu<0$, the spectral theorem gives
$\|e^{t\mathcal L}u_0\|_2\le e^{\mu t}\|u_0\|_2$.
Since $\|z\|_\infty\le\|z\|_2$ for $z\in\ell^2(\Tq)$ with counting measure,
\eqref{eq:threshold-semigroup-comparison} proves
\eqref{eq:exp-extinction}.
For $\mu=0$, let $\nu_{u_0}$ be the scalar spectral measure of
$u_0$ for $\mathcal L$.
Then
\[
 \norm{e^{t\mathcal L}u_0}_2^2
 =\int_{(-\infty,0]}e^{2ts}\,\dd\nu_{u_0}(s)\to0
 \quad\text{ as } t\to\infty.
\]
By Lemma~\ref{lem:horo-spectrum}, $\ker\mathcal L=\{0\}$,
and hence $\nu_{u_0}(\{0\})=0$. For $s<0$, the integrand tends
to zero and is bounded by one. Dominated convergence applies because
$\nu_{u_0}(\R)=\|u_0\|_2^2<\infty$.
Combining this limit with \eqref{eq:threshold-semigroup-comparison}
proves uniform extinction for every $u_0\in\ell^2(\Tq)$ with $0\le u_0\le1$. For
$\mu>0$, Lemma~\ref{lem:horo-spectrum} supplies finite connected
sets $\Omega_x\subseteq B(x,R_0)$ and principal eigenvectors
$\varphi_x$ of their zero-exterior linearized operators, normalized
by $\max_{y\in\Omega_x}\varphi_x(y)=1$, with eigenvalues at least $\mu/2$ and
$\varphi_x(x)\ge m_0>0$. The constants $R_0,m_0$ are independent of
$x$. The bounds on the eigenvalues, the support radii and
$\varphi_x(x)$ verify the three uniform hypotheses of
Lemma~\ref{lem:local-general}, with $\eta_0=\mu/2$.
That lemma gives $u(t,x)\to1$ as $t\to\infty$ at each fixed
vertex for every nonzero datum in $[0,1]$.\qedhere

\end{proof}

\begin{proof}[Proof of Theorem~\ref{thm:radial-threshold}]

Let $\mathcal L_r=D\Delta_{\Tq}+a_{r(\cdot)}$ on
$\ell^2(\Tq)$ with counting measure. The bounded degree and real
bounded potential make $\mathcal L_r$ bounded and self-adjoint.
The same coordinate formula is bounded on $\ell^\infty(\Tq)$, and
its exponential series agrees pointwise on $\ell^2$ data with the
Hilbert-space series. For $\mu\le0$, Lemma~\ref{lem:radial-spectrum} gives
$\sup\spec(\mathcal L_r)=\mu$ and $\ker(\mu-\mathcal L_r)=\{0\}$.
Since $f_{r(x)}(u)\le a_{r(x)}u$ on $[0,1]$, comparison yields
$0\le u(t,\cdot)\le e^{t\mathcal L_r}u_0$. If $\mu<0$, then
\[
 \|u(t,\cdot)\|_\infty
 \le\|e^{t\mathcal L_r}u_0\|_2
 \le e^{\mu t}\|u_0\|_2
 \qquad t\ge0.
\]
If $\mu=0$, let $\nu_{u_0}^{\,r}$ be the scalar spectral measure of
$u_0\in\ell^2(\Tq)$ for $\mathcal L_r$. It is supported in
$(-\infty,0]$, has total mass $\|u_0\|_2^2$, and satisfies
$\nu_{u_0}^{\,r}(\{0\})=0$ because $\ker\mathcal L_r=\{0\}$. Hence
\[
 \|e^{t\mathcal L_r}u_0\|_2^2
 =\int_{(-\infty,0]}e^{2tE}\,\dd\nu_{u_0}^{\,r}(E)\to0
 \quad\text{ as }t\to\infty
\]
by dominated convergence: $e^{2tE}\le1$ for $E\le0$ and tends to zero
for $E<0$. The inequality
$\|u(t,\cdot)\|_\infty\le\|e^{t\mathcal L_r}u_0\|_2$
proves uniform extinction. For $\mu>0$,
Lemma~\ref{lem:radial-local}, applied to the nonzero datum
$0\le u_0\le1$, gives convergence to one at each fixed vertex.\qedhere

\end{proof}

\begin{theorem}\label{thm:step}

Assume \eqref{eq:hyp}. Suppose the initial datum of \eqref{eq:lattice}
satisfies $\ind_{\{j\le A\}}\le w_j(0)\le\ind_{\{j\le B\}}$ for
every $j\in\Z$, where $A\le B$ are integers. Then
\begin{equation}
\label{eq:step-limits}
 \lim_{t\to\infty}\sup_{\substack{j\in\Z\\j\ge ct}}w_j(t)=0\quad\text{ for } c>c_+,\qquad
 \lim_{t\to\infty}\inf_{\substack{j\in\Z\\j\le ct}}w_j(t)=1\quad\text{ for } c<c_+.
\end{equation}
For every $\theta\in(0,1)$, the level position
$H_\theta(t)=\max\{j\in\Z:w_j(t)\ge\theta\}$ exists for all sufficiently large $t$,
and there is a constant $C_\theta<\infty$ such that
\begin{equation}
\label{eq:step-position}
 \frac{H_\theta(t)}{t}\to c_+\quad\text{ as } t\to\infty,
 \qquad
 H_\theta(t)\le c_+t-\frac1{2\lambda_+}\log(1+t)+C_\theta\quad\text{for all sufficiently large }t.
\end{equation}
The same conclusions hold for solutions of \eqref{eq:tree-cauchy} with
$f(x,s)=f_{b(x)}(s)$ and
$\ind_{\{b\le A\}}\le u_0\le\ind_{\{b\le B\}}$, with $j$ replaced
by $b(x)$ and with level position
\[
 H_\theta(t)=\max\{b(x):x\in\Tq,\ u(t,x)\ge\theta\}.
\]
The initial datum may vary within each horocycle.
The speed can have either sign or be zero.

\end{theorem}
\begin{proof}

For every $\lambda>0$, the bound $w_j(0)\le\ind_{\{j\le B\}}$ gives
\[
 \sum_{j\in\Z}e^{\lambda j}w_j(0)
 \le\sum_{j\le B}e^{\lambda j}
 =\frac{e^{\lambda B}}{1-e^{-\lambda}}<\infty.
\]
Thus Lemma~\ref{lem:linear-bound} applies with $\lambda=\lambda_+$.
Since $k(\lambda_+)=c_+\lambda_+$, it gives
\begin{equation}
\label{eq:step-upper}
 w_j(t)\le C(1+t)^{-1/2}e^{-\lambda_+(j-c_+t)}
 \qquad j\in\Z,\ t\ge0.
\end{equation}
For $j\ge ct$ with $c>c_+$, \eqref{eq:step-upper} tends to zero
uniformly, proving the first limit in \eqref{eq:step-limits}.
For $\theta\in(0,1)$, $t\ge0$ and $j\in\Z$ with
$w_j(t)\ge\theta$, \eqref{eq:step-upper} gives
\[
 j\le c_+t-\frac{\log(1+t)}{2\lambda_+}
          +\frac1{\lambda_+}\log\frac C\theta.
\]
This is the upper bound in \eqref{eq:step-position} whenever
$\{j:w_j(t)\ge\theta\}$ is nonempty.

For the lower limit, fix $c<c_+$. Since $c_-<c_+$ by
Lemma~\ref{lem:matrix}, choose
$\max\{c,c_-\}<c_2<c_+$ and then $c<c_1<c_2$.
Let $v(t,j)$ be the subsolution given by Lemma~\ref{lem:bump}
at speed $c_2$ with $R=P$.
Translation by a multiple of $P$ preserves \eqref{eq:lattice}.
Let $L_0$ be the support radius supplied by Lemma~\ref{lem:bump},
so that $v(0,j)=0$ for $|j|\ge L_0$. Set
$K_0=\lceil(L_0-A)/P\rceil$. Then the support of
$j\mapsto v(0,j+K_0P)$ lies in $\{j\in\Z:j\le A\}$. This remains true
for every integer $k\ge K_0$, and comparison gives
$w_j(t)\ge v(t,j+kP)$ for $j\in\Z$ and $t\ge0$. For $j\le c_1t$, take
\[
 k=\left\lceil\frac{c_2t-j}{P}\right\rceil.
\]
Then $0\le j+kP-c_2t<P$ and
$k\ge(c_2-c_1)t/P\ge K_0$ for all sufficiently large $t$,
uniformly in $j\le c_1t$. The bound
$v(t,j+kP)\ge\delta$ from Lemma~\ref{lem:bump} on
$|j+kP-c_2t|\le P$ therefore gives
\begin{equation}
\label{eq:step-positive}
 w_j(t)\ge\delta\quad\text{ for }j\in\Z,\ j\le c_1t,\ t\ge t_0.
\end{equation}
The construction uses $c_2-c_1>0$ and is valid for every sign
of $c_2$.

If the lower limit in \eqref{eq:step-limits} failed, there would be
$t_n\to\infty$ as $n\to\infty$, $j_n\le ct_n$ and $\epsilon_0>0$ with
$w_{j_n}(t_n)\le1-\epsilon_0$. Pass to a subsequence with
$j_n\equiv\rho\pmod P$ for one fixed $\rho\in\{0,\ldots,P-1\}$.
The translated reactions are then exactly $f_{l+\rho}$.
After a further subsequence, Lemma~\ref{lem:compactness} gives
an entire solution $V:\R\times\Z\to[0,1]$ that is the limit of
$w_{j_n+l}(t_n+s)$, uniformly on
bounded $s$-intervals at each
$l\in\Z$, satisfying
\[
 \partial_sV(s,l)
 =D\bigl(V(s,l-1)+qV(s,l+1)-(q+1)V(s,l)\bigr)
       +f_{l+\rho}(V(s,l))
 \qquad s\in\R,\ l\in\Z.
\]
For fixed $s,l$,
$c_1(t_n+s)-(j_n+l)\ge(c_1-c)t_n+c_1s-l\to\infty$
as $n\to\infty$. Thus \eqref{eq:step-positive} gives
$\delta\le V(s,l)\le1$ at every space--time point.
The functions $f_{l+\rho}$ satisfy \eqref{eq:hyp}, so
Lemma~\ref{lem:liouville} gives $V\equiv1$, contradicting
$V(0,0)\le1-\epsilon_0$. This proves the lower limit.

For each $\theta\in(0,1)$, the lower limit in \eqref{eq:step-limits}
makes $\{j:w_j(t)\ge\theta\}$ nonempty for large $t$, while
\eqref{eq:step-upper} bounds this set above. A nonempty subset of $\Z$
bounded above has a maximum. For each $c<c_+$, the lower limit gives
$H_\theta(t)\ge\lfloor ct\rfloor$ for all sufficiently large $t$;
for each $c>c_+$, the upper limit gives $H_\theta(t)<ct$.
Letting these fixed speeds approach $c_+$ proves
$H_\theta(t)/t\to c_+$ as $t\to\infty$.
For the tree equation, let $w^A,w^B$ solve \eqref{eq:lattice} with
initial values $\ind_{\{j\le A\}}$ and $\ind_{\{j\le B\}}$,
respectively, and use their lifts
$x\mapsto w^A_{b(x)}(t)$ and $x\mapsto w^B_{b(x)}(t)$ as lower
and upper solutions. Each vertex has one neighbor at Busemann coordinate $b(x)-1$
and $q$ at coordinate $b(x)+1$. The two lifted solutions therefore
satisfy the differential equation in \eqref{eq:tree-cauchy}, and
Lemma~\ref{lem:comparison} gives the two limits and the
level-position bound for $u(t,x)$ under the
stated bounds on $u_0$.\qedhere

\end{proof}

\begin{theorem}\label{thm:radial}

Assume \eqref{eq:hyp}, $f(x,s)=f_{r(x)}(s)$, and
$\mu>0$. For every nonzero finitely supported initial datum
$0\le u_0\le1$,
\begin{equation}
\label{eq:radial-limits}
 \lim_{t\to\infty}\sup_{\substack{x\in\Tq\\r(x)\ge ct}}u(t,x)=0\quad\text{ for } c>c_+,\qquad
 \lim_{t\to\infty}\inf_{\substack{x\in\Tq\\r(x)\le ct}}u(t,x)=1\quad\text{ for } 0<c<c_+.
\end{equation}
Thus, for every $0<\theta<1$ and $0<\epsilon<c_+$, for all sufficiently
large $t$,
\[
 B(o,(c_+-\epsilon)t)\subseteq\{x\in\Tq:u(t,x)\ge\theta\}
 \subseteq B(o,(c_++\epsilon)t).
\]
Consequently, $R_\theta^o(t)/t\to c_+$ as $t\to\infty$ for every
$\theta\in(0,1)$, where
$R_\theta^o(t)=\max\{r(x):u(t,x)\ge\theta\}$.

\end{theorem}
\begin{proof}

Since $\mu>0$, Lemma~\ref{lem:matrix} gives $c_+>0$ and
$k(\lambda_+)=c_+\lambda_+>0$. Extend the positive eigenvector
$\phi(\lambda_+)$ periodically to $r\in\Z$. Since $u_0$ has finite
{support, choose
\[
 C\ge\max\left\{\frac1{\phi_0(\lambda_+)},
 \max_{x\in\supp u_0}
 \frac{e^{\lambda_+r(x)}u_0(x)}{\phi_{r(x)}(\lambda_+)}
 \right\}.
\]
The maximum is finite, and this choice gives
$Ce^{-\lambda_+r(x)}\phi_{r(x)}(\lambda_+)\ge u_0(x)$ for every
$x\in\Tq$ and $C\phi_0(\lambda_+)\ge1$.} Define
\[
 \overline U_r(t)=\min\{1,Ce^{-\lambda_+r+k(\lambda_+)t}
                                      \phi_r(\lambda_+)\}
 \qquad r\in\N_0,\ t\ge0.
\]
Its lift satisfies $\overline U_{r(x)}(0)\ge u_0(x)$ for every
$x\in\Tq$, and $\overline U_0(t)=1$ for all $t\ge0$.
Where the exponential term is below one and $r\ge1$, the eigenvector
equation \eqref{eq:matrix} makes it a solution of the equation
linearized at zero. The inequality $f_r(s)\le a_rs$ and truncation
of the neighboring values from above can only decrease the nonlinear
right-hand side. Where the function equals one, its derivative is
zero and its diffusion term is nonpositive. At $r=0$, $\overline U_0=1$, so the boundary right-hand side is
$D(q+1)(\overline U_1-1)+f_0(1)\le0$. Where $\overline U_r(t)<1$, its derivative is
$k(\lambda_+)\overline U_r(t)\le k(\lambda_+)$; where
$\overline U_r(t)=1$, its derivative vanishes almost everywhere.
Thus
\[
 \sup_{r\in\N_0}|\overline U_r(t)-\overline U_r(s)|
 \le k(\lambda_+)|t-s|
 \quad\text{ for }s,t\ge0.
\]
In particular, $t\mapsto(\overline U_r(t))_{r\in\N_0}$ is locally Lipschitz
in $\ell^\infty(\N_0)$. The supersolution inequality holds
almost everywhere, including the intervals adjacent to a switching time. The version
of comparison in Lemma~\ref{lem:comparison} therefore gives
$u(t,x)\le C(\max_{0\le i<P}\phi_i(\lambda_+))
e^{-\lambda_+(c-c_+)t}$ for $t\ge0$, $c>c_+$ and
$x\in\Tq$ with $r(x)\ge ct$. The right-hand side tends to zero
as $t\to\infty$, uniformly over these vertices. This proves the upper limit in
\eqref{eq:radial-limits}.

For the lower limit, fix $\tau>0$ and choose
$0<\varepsilon\le u(\tau,o)$ by Lemma~\ref{lem:comparison}.
Let $U$ solve \eqref{eq:radial} from
$U_r(0)=\varepsilon\mathbf1_{\{r=0\}}$.
Comparison gives $U_{r(x)}(t-\tau)\le u(t,x)$ for $t\ge\tau$.
It suffices to prove the lower limit for $U$ at every speed in
$(0,c_+)$: given $0<c<c'<c_+$, one has $ct\le c'(t-\tau)$
for all sufficiently large $t$. Lemma~\ref{lem:radial-local} gives $U_r(t)\to1$ as $t\to\infty$ for every fixed $r$.
Fix $0<c<c_1<c_2<c_+$. Since $c_-<0<c_2<c_+$,
let $v(t,j)$ be the subsolution given by Lemma~\ref{lem:bump}
at speed $c_2$, with $R=1$, and denote its support radius by $L_0$. Choose a fixed multiple $j_0$ of $P$
with $j_0>L_0+2$. Replace $v$ by $v/2$ and $\delta$ by $\delta/2$, which is permitted
by Lemma~\ref{lem:bump}; then $0\le v\le1/2$. The set of indices
supporting $v(0,\cdot-j_0)$ is finite. Local convergence of
$U_r(s)$ as $s\to\infty$ therefore
gives a time $s_0$ such that
\[
 U_r(s)\ge v(0,r-j_0)\quad\text{ for }r\in\N_0,\ s\ge s_0.
\]
For every launch time $s\ge s_0$, translation by $j_0\in P\Z$
preserves the interior reaction phases. For $t\ge s$, the support of
$r\mapsto v(t-s,r-j_0)$ lies in
$|r-j_0-c_2(t-s)|<L_0$, hence in $r>2$ because $c_2>0$.
The function $r\mapsto v(t-s,r-j_0)$ vanishes at $r=0,1$, so its
boundary derivative and the boundary right-hand side are both zero.
Thus $r\mapsto v(t-s,r-j_0)$ is a subsolution of the radial half-line
equation, and comparison holds for every such launch time.

For integers $r$ with $j_0\le r\le c_1t$, take
$s=t-(r-j_0)/c_2$. For all sufficiently large $t$ this satisfies
$s\ge s_0$, uniformly in that range of $r$. Since
$r-j_0-c_2(t-s)=0$, the lower bound for $v$ gives $U_r(t)\ge\delta$. At the finitely many
$r<j_0$, local convergence gives the same conclusion after making
$\delta$ smaller. Thus $U_r(t)\ge\delta$ for $0\le r\le c_1t$ and large $t$.

If convergence to one on $r\le ct$ failed, there would be
$\epsilon_0>0$, $t_n\to\infty$, and $0\le r_n\le ct_n$ with
$U_{r_n}(t_n)\le1-\epsilon_0$. A bounded subsequence of $r_n$
would contradict convergence at each fixed index. Thus we may take
$r_n\to\infty$ and $r_n\equiv\rho\pmod P$ along a subsequence,
for one $\rho\in\{0,\ldots,P-1\}$.
After passing to a further subsequence, Lemma~\ref{lem:compactness}
gives convergence of $U_{r_n+j}(t_n+s)$, uniformly on bounded
$s$-intervals at each $j\in\Z$, to an entire solution
$W:\R\times\Z\to[0,1]$ of the lattice equation with reaction
$f_{j+\rho}$ at index $j$;
the shifted domains contain every fixed integer
interval because $r_n\to\infty$. For each fixed $(s,j)$,
\[
 0\le r_n+j\le c_1(t_n+s)
\]
for all sufficiently large $n$, since $r_n\to\infty$,
$r_n\le ct_n$ and $c<c_1$.
Consequently $\delta\le W(s,j)\le1$, whereas the chosen sequence gives
$W(0,0)\le1-\epsilon_0$. Lemma~\ref{lem:liouville} implies
$W\equiv1$, a contradiction. This proves
\eqref{eq:radial-limits}. The level-set assertions follow by taking
speeds on either side of $c_+$; the upper exponential estimate also
ensures finiteness of the level sets at each positive time.\qedhere

\end{proof}

\subsection{Periodic linear growth at zero}

We prove Theorems~\ref{thm:periodic-bramson-step},
\ref{thm:horo-compact-bramson}, and~\ref{thm:periodic-bramson-distance}
by applying the moving-boundary
estimates of Section~\ref{sec:moving-boundaries}.
 Constants may depend
on the initial datum and the level. Recall
\[
 m_*(t)=c_+t-\frac{3}{2\lambda_+}\log(1+t),\quad\text{for }t\ge0.
\]

\begin{proof}[Proof of Theorem~\ref{thm:periodic-bramson-step}]

Let $w^A$ and $w^B$ solve \eqref{eq:lattice} with initial values
$\mathbf1_{\{j\le A\}}$ and $\mathbf1_{\{j\le B\}}$,
respectively. Comparison gives $w^A_j(t)\le w_j(t)\le w^B_j(t)$.
The centers in Lemmas~\ref{lem:sharp-step-lower} and
\ref{lem:sharp-step-upper} are $m_A(t)=A+m_*(t)$ and
$m_B(t)=B+m_*(t)$, with exponent $\lambda=\lambda_+$.
Both lemmas apply at every integer step location and for
positive, zero or negative $c_+$.
By Lemma~\ref{lem:sharp-step-lower}, for every $\theta\in(0,1)$
there is $C_\theta^-$ such that
$w^A_j(t)\ge\theta$ whenever $j\le A+m_*(t)-C_\theta^-$
and $t$ is sufficiently large.
Lemma~\ref{lem:sharp-step-upper} gives
$w^B_j(t)\le C(1+y)e^{-\lambda_+y}$ for $y=j-B-m_*(t)>0$.
Choose $Y_\theta>0$ so that
$\sup_{y\ge Y_\theta}C(1+y)e^{-\lambda_+y}<\theta$.
Comparison therefore yields
\[
 \{j\in\Z:j\le A+m_*(t)-C_\theta^-\}
 \subseteq\{j\in\Z:w_j(t)\ge\theta\}
 \subseteq\{j\in\Z:j<B+m_*(t)+Y_\theta\}
\]
for all sufficiently large $t$. The two inclusions make
$\{j\in\Z:w_j(t)\ge\theta\}$ nonempty and bounded above, and hence give
$H_\theta(t)=m_*(t)+O_{\theta,A,B}(1)$ after absorbing the fixed shifts
$A,B$ into the constants.

A vertex $x$ with $b(x)=j$ has one neighbor at coordinate $j-1$
and $q$ neighbors at coordinate $j+1$. Consequently,
$x\mapsto w^A_{b(x)}(t)$ and $x\mapsto w^B_{b(x)}(t)$ satisfy the
differential equation in \eqref{eq:tree-cauchy}, by
\eqref{eq:lattice}.
Their initial values bracket $u_0$, and tree comparison gives
$w^A_{b(x)}(t)\le u(t,x)\le w^B_{b(x)}(t)$.
Thus $u(t,x)\ge\theta$ for $b(x)\le A+m_*(t)-C_\theta^-$ and
$u(t,x)<\theta$ for $b(x)\ge B+m_*(t)+Y_\theta$, for all
sufficiently large $t$. Absorbing $A,B$ into the constants proves
the claimed tree bounds.
\end{proof}

\begin{proof}[Proof of Theorem~\ref{thm:horo-compact-bramson}]
Fix $\theta\in(0,1)$. By \eqref{eq:finite-mass-level-count},
\[
 \#\{x\in\Tq:u(t,x)\ge\theta\}
 \le\theta^{-1}e^{a_{\max}t}\|u_0\|_1<\infty
 \quad\text{ for every finite }t\ge0.
\]
Since $\mu>0$, Theorem~\ref{thm:threshold} gives
$u(t,x)\to1$ as $t\to\infty$ for each fixed $x\in\Tq$.
Thus the level set is eventually nonempty and $H_\theta^b(t)$ is defined.

For the upper bound, let $B=\max\{b(x):x\in\supp u_0\}$,
and let $w^B$ solve \eqref{eq:lattice} from
$w^B_j(0)=\mathbf1_{\{j\le B\}}$ for $j\in\Z$.
Every vertex with coordinate $j$ has one neighbor with coordinate
$j-1$ and $q$ neighbors with coordinate $j+1$. Hence
$x\mapsto w^B_{b(x)}(t)$ solves the tree equation, and
Lemma~\ref{lem:comparison} gives $u(t,x)\le w^B_{b(x)}(t)$ for
$t\ge0$ and $x\in\Tq$.
By Lemma~\ref{lem:sharp-step-upper}, there is $C>0$ such that
\[
 w^B_j(t)\le C(1+y)e^{-\lambda_+y},\qquad
 y=j-B-m_*(t)>0,\quad\text{ for }t\ge0,\ j\in\Z.
\]
Choose $Y_\theta\ge1$ with
$\sup_{y\ge Y_\theta}C(1+y)e^{-\lambda_+y}<\theta$.
The comparison then yields
$H_\theta^b(t)\le m_*(t)+|B|+Y_\theta$ whenever the maximum exists.

For the lower bound, fix $x_*\in\Tq$ and write $b_*=b(x_*)$.
Choose $\beta\in(\theta,1)$. By Theorem~\ref{thm:threshold},
there is $T_0\ge0$ such that
$u(T_0+t,x_*)\ge\beta$ for every $t\ge0$.
Define
\[
 \mathcal C_r=\{x\in\Tq:d(x_*,x)=r,\ b(x)=b_*+r\}
 \quad\text{ for }r\in\N_0,
 \qquad F=\bigcup_{r\ge1}\mathcal C_r.
\]
Then $\mathcal C_0=\{x_*\}$ and $\#\mathcal C_r=q^r$.
Every vertex in $\mathcal C_r$, $r\ge1$, has one neighbor in
$\mathcal C_{r-1}$ and $q$ neighbors in $\mathcal C_{r+1}$.

For $j\in\Z$ and $z\in[0,1]$, set
\begin{equation}\label{eq:horo-compact-reactions}
 \widehat f_j(z)=\beta^{-1}f_{b_*+j}(\beta z)(1-z).
\end{equation}
The functions $\widehat f_j$ are $P$-periodic in $j$ and $C^2$ in $z$.
From \eqref{eq:hyp} and \eqref{eq:horo-compact-reactions},
\[
 \widehat f_j(0)=\widehat f_j(1)=0,\quad
 \widehat f_j'(0)=a_{b_*+j},\quad
 \widehat f_j'(1)=-\beta^{-1}f_{b_*+j}(\beta)<0
 \quad\text{ for }j\in\Z,
\]
and $0<\widehat f_j(z)\le a_{b_*+j}z$ for
$j\in\Z$ and $0<z<1$. Thus \eqref{eq:hyp} holds for
$(\widehat f_j)$. Cyclic permutation of the periodic indices
conjugates its linearization matrices to $M_\lambda$; consequently
its functions $k(\lambda)$ and its constants $c_+,\lambda_+$
coincide with those of $(f_j)$.

Let $W$ solve \eqref{eq:lattice} with reactions $\widehat f_j$
and initial values $W_j(0)=\mathbf1_{\{j\le0\}}$ for $j\in\Z$.
Define $V_j(t)=\beta W_j(t)$ for $t\ge0$ and $j\in\Z$.
By Lemma~\ref{lem:comparison}, $0\le V_j(t)\le\beta$ and
\[
 \dot V_j(t)=D\big(V_{j-1}(t)+qV_{j+1}(t)-(q+1)V_j(t)\big)
 +f_{b_*+j}(V_j(t))\big(1-V_j(t)/\beta\big)
 \quad\text{ for }t\ge0,\ j\in\Z.
\]
Applying Lemma~\ref{lem:sharp-step-lower} to $W$ at level
$\theta/\beta\in(0,1)$ gives $C^-,T^-\ge0$ such that
\[
 V_r(t)\ge\theta\quad\text{ for }t\ge T^-,\quad
 r\in\N_0\text{ with }r\le m_*(t)-C^-.
\]

Define $v:[0,\infty)\times\Tq\to[0,\beta]$ by
\[
 v(t,x)=\begin{cases}
 V_r(t),&x\in\mathcal C_r,\ r\in\N_0,\\
 0,&x\notin\bigcup_{r\ge0}\mathcal C_r.
 \end{cases}
\]
Since $t\mapsto(V_j(t))_{j\in\Z}$ belongs to $C^1([0,\infty);\ell^\infty(\Z))$ and
$\|v(t,\cdot)-v(s,\cdot)\|_\infty
\le\sup_{j\in\Z}|V_j(t)-V_j(s)|$ for $s,t\ge0$,
$v$ is locally Lipschitz with values in $\ell^\infty(\Tq)$.
The neighbor counts in $\mathcal C_r$ and the equation for $V$ give
\[
 \partial_t v(t,x)-D\Delta_{\Tq}v(t,x)-f_{b(x)}(v(t,x))
 =-\beta^{-1}v(t,x)f_{b(x)}(v(t,x))\le0
 \quad\text{ for }t\ge0,\ x\in F.
\]
At $x_*$ we have $v(t,x_*)\le\beta\le u(T_0+t,x_*)$ for
$t\ge0$, and $v(t,x)=0\le u(T_0+t,x)$ for all other
$x\in\Tq\setminus F$. Moreover, $v(0,x)=0\le u(T_0,x)$
for $x\in F$. Applying Lemma~\ref{lem:comparison} on $F$, with
these inequalities on $\Tq\setminus F$ throughout each finite time
interval, yields
\begin{equation}\label{eq:horo-compact-cone-comparison}
 u(T_0+t,x)\ge V_r(t)
 \quad\text{ for }t\ge0,\ r\ge1,\ x\in\mathcal C_r.
\end{equation}
Here and below $r$ ranges over integers.

Since $\mu>0$ implies $c_+>0$, the integer
$r_t=\lfloor m_*(t)-C^-\rfloor$ is at least one for all sufficiently
large $t$. The set $\mathcal C_{r_t}$ is nonempty, so
\eqref{eq:horo-compact-cone-comparison} and the lower bound for $V$
give
\[
 H_\theta^b(T_0+t)\ge b_*+r_t
 \ge m_*(t)+b_*-C^--1
 \quad\text{ for all sufficiently large }t.
\]
Finally, the definition of $m_*$ gives
\[
 |m_*(t+T_0)-m_*(t)|
 \le c_+T_0+\frac{3}{2\lambda_+}\log(1+T_0)
 \quad\text{ for }t\ge0.
\]
Combining the last two estimates gives the lower bound in
\eqref{eq:horo-compact-bramson} with
$C_\theta^-=|b_*|+C^-+1+c_+T_0+
3\log(1+T_0)/(2\lambda_+)$.
Together with the upper bound, this proves the result with
$C_\theta=\max\{|B|+Y_\theta,C_\theta^-\}$ and a sufficiently
large $T_\theta$ such that $t-T_0\ge T^-$ and
$r_{t-T_0}\ge1$ for every $t\ge T_\theta$.
\end{proof}

\Needspace{7\baselineskip}
\begin{proof}[Proof of Theorem~\ref{thm:periodic-bramson-distance}]

First suppose $f(x,u)=f_{r(x)}(u)$. Choose an integer $R\ge0$
with $\supp u_0\subseteq B(o,R)$.
Let $U^+$ solve \eqref{eq:radial} from
$U_r^+(0)=\mathbf1_{\{r\le R\}}$. Fix $\tau>0$ and choose
$0<\varepsilon\le u(\tau,o)$, using strict positivity from
Lemma~\ref{lem:comparison}. Let $U^-$ solve \eqref{eq:radial} from
$U_r^-(0)=\varepsilon\mathbf1_{\{r=0\}}$. Their tree lifts and
comparison give
\[
 U^-_{r(x)}(t-\tau)\le u(t,x)\le U^+_{r(x)}(t)
 \quad\text{ for }t\ge\tau,\ x\in\Tq.
\]
Both coordinate data are nonzero, take values in $[0,1]$, and have
finite support in $\N_0$. The boundary equation for each of $U^-$ and $U^+$ is
$\dot U_0^\pm=D(q+1)(U_1^\pm-U_0^\pm)+f_0(U_0^\pm)$,
so the coefficient $b_0$ in Lemma~\ref{lem:periodic-halfline-log}
equals $D(q+1)>0$. Since $\mu>0$,
Lemma~\ref{lem:radial-local} gives $U_0^\pm(t)\to1$ as
$t\to\infty$. Thus both solutions satisfy the finite-support condition and
the boundary convergence hypothesis of
Lemma~\ref{lem:periodic-halfline-log}. Apply the outer inclusion in \eqref{eq:periodic-halfline-log} to
$U^+$ and the inner inclusion to $U^-$. For the fixed comparison time $\tau>0$ and $t>\tau$,
\[
 m_*(t-\tau)-m_*(t)
 =-c_+\tau+\frac{3}{2\lambda_+}
            \log\frac{1+t}{1+t-\tau}=O(1)
 \quad\text{ as }t\to\infty.
\]
Thus the lower bound centered at $m_*(t-\tau)$ and the upper bound
centered at $m_*(t)$ give \eqref{eq:periodic-bramson-balls}.

Now suppose $f(x,u)=f_{s(x)}(u)$. Let $V^-,V^+$ solve the interior
equation in \eqref{eq:radial} for $r\ge1$ and the boundary equation
\eqref{eq:geodesic-boundary}, with initial values $V_r^-(0)=\eta\mathbf1_{\{r=0\}}$ and
$V_r^+(0)=\mathbf1_{\{r\le R_2\}}$, respectively. On
$\mathcal G$, the two neighboring values along the geodesic equal
$V_0^\pm$ and cancel from the diffusion sum; the other $q-1$
values equal $V_1^\pm$. Thus the boundary coefficient of each lift
is $b_0=D(q-1)>0$. Since $R_1\ge0$,
$\mathcal G\subseteq\mathcal N_{R_1}(\mathcal G)$, so the bounds
$\eta\mathbf1_{\mathcal N_{R_1}(\mathcal G)}\le u_0
\le\mathbf1_{\mathcal N_{R_2}(\mathcal G)}$ give
\[
 V^-_{s(x)}(0)\le u_0(x)\le V^+_{s(x)}(0)
 \qquad x\in\Tq,
\]
so comparison preserves these inequalities at all $t\ge0$.
The two coordinate data are nonzero and finitely supported in
$\N_0$, although their tree lifts have infinite support.
Lemma~\ref{lem:geodesic-local}, with $\mu>0$, gives
$V_0^\pm(t)\to1$ as $t\to\infty$. Apply
Lemma~\ref{lem:periodic-halfline-log} to both, use the inner bound
for $V^-$ and the outer bound for $V^+$, and substitute $r=s(x)$.
This proves \eqref{eq:periodic-bramson-tubes}.\qedhere

\end{proof}

\subsection{Common linear growth at zero}

\begin{proof}[Proof of Theorem~\ref{thm:common}]

Suppose first that $a>\rho_q$ and that $u_0$ is nonzero and
finitely supported. Lemma~\ref{lem:envelopes} gives
$h_-,h_+\in C^2([0,1])$ satisfying
\[
 h_-(s)\le f_i(s)\le h_+(s)\quad\text{ for }0\le s\le1,\ 0\le i<P,
 \qquad h_-'(0)=h_+'(0)=a.
\]
Let $u^-$ and $u^+$ solve \eqref{eq:tree-cauchy} with reactions
$h_-$ and $h_+$, respectively, and initial datum $u_0$.
Comparison gives $u^-(t,x)\le u(t,x)\le u^+(t,x)$.
By Lemma~\ref{lem:envelopes}, both comparison reactions satisfy
\eqref{eq:hyp}. Lemma~\ref{lem:hom-tree} therefore applies to
$u^-$ and $u^+$ with the same nonzero finitely supported datum.
Their derivatives at zero both equal $a$, so
\eqref{eq:homogeneous-parameters} gives the same dispersion function
$\Psi_a$, speed $c_*(a)$ and exponent $\lambda_*(a)$.
Consequently, for each $\theta\in(0,1)$, the sets
$\{x\in\Tq:u^\pm(t,x)\ge\theta\}$ satisfy the two ball
inclusions in \eqref{eq:hom-tree-balls}, with radii
$c_*(a)t-3\log t/(2\lambda_*(a))$ up to bounded errors.
The inner bound for $u^-$ and the outer bound for $u^+$ give
\eqref{eq:common-balls}, since $\log(1+t)-\log t\to0$
as $t\to\infty$.
For $a\le\rho_q$ and $u_0\in\ell^2(\Tq)$, use
$f_{\iota(x)}(s)\le as$ for $0\le s\le1$ and compare with
the constant-coefficient linear semigroup on $\ell^2(\Tq)$,
\[
 e^{t\mathcal L_a}u_0,
 \qquad
 \mathcal L_a:=D\Delta_{\Tq}+a
 \qquad t\ge0.
\]
Lemma~\ref{lem:horo-spectrum}, applied to the constant potential,
gives
\[
 \sup\spec(\mathcal L_a)=a-\rho_q
\]
and gives $\ker(\mathcal L_a-(a-\rho_q)I)=\{0\}$.
If $a<\rho_q$, the spectral theorem yields
\[
 \|u(t,\cdot)\|_\infty
 \le
 \|e^{t\mathcal L_a}u_0\|_2
 \le
 e^{(a-\rho_q)t}\|u_0\|_2
 \qquad t\ge0.
\]
If $a=\rho_q$, let $\nu_{u_0}$ be the scalar spectral measure of
$u_0$ for $\mathcal L_a$. It is supported in $(-\infty,0]$ and has
finite mass $\nu_{u_0}(\R)=\|u_0\|_2^2$. Since $\ker\mathcal L_a=\{0\}$, one has
$\nu_{u_0}(\{0\})=0$, and
\[
 \|e^{t\mathcal L_a}u_0\|_2^2
 =
 \int_{(-\infty,0]}e^{2ts}\,\dd\nu_{u_0}(s)
 \to0\quad\text{ as } t\to\infty
\]
by dominated convergence.  Again
$\|u(t,\cdot)\|_\infty\le\|e^{t\mathcal L_a}u_0\|_2$, which proves the
critical extinction assertion for every assignment $\iota$.\qedhere

\end{proof}

\begin{proof}[Proof of Corollary~\ref{cor:common-step}]

Choose $h_-,h_+$ from Lemma~\ref{lem:envelopes}, so that
$h_-\le f_i\le h_+$ on $[0,1]$ and $h_-'(0)=h_+'(0)=a$.
Let $W^-,W^+$ solve
$\dot W_j^\pm=D(W_{j-1}^\pm+qW_{j+1}^\pm-(q+1)W_j^\pm)
+h_\pm(W_j^\pm)$ for $j\in\Z$, with initial values
$W^-_j(0)=W^+_j(0)=\mathbf1_{\{j\le0\}}$.
The initial bounds at $A,B$ and comparison give
\[
 W^-_{j-A}(t)\le w_j(t)\le W^+_{j-B}(t)
 \qquad j\in\Z,\ t\ge0.
\]
By Lemma~\ref{lem:graham}, $W^\pm$ are nonincreasing in $j$ and
their level positions equal
$m_a(t)+O_\theta(1)$, where
$m_a(t)=c_*(a)t-3\log t/(2\lambda_*(a))$ for $t>0$.
Thus the lower comparison bounds $H_\theta(t)$ below by
$A+m_a(t)-C_\theta$, and the upper comparison bounds it above by
$B+m_a(t)+C_\theta$. This proves \eqref{eq:common-step}.
Every vertex $x\in\Tq$ has one neighbor with Busemann coordinate
$b(x)-1$ and $q$ neighbors with coordinate $b(x)+1$.
Consequently, the functions $W^-_{b(x)-A}(t)$ and $W^+_{b(x)-B}(t)$
solve the differential equation in \eqref{eq:tree-cauchy} with reactions
$h_-$ and $h_+$, respectively. For
$\ind_{\{b\le A\}}\le u_0\le\ind_{\{b\le B\}}$, comparison gives
\[
 W^-_{b(x)-A}(t)\le u(t,x)\le W^+_{b(x)-B}(t)
 \qquad x\in\Tq,\ t\ge0.
\]
The scalar level positions give the stated bounds with $j=b(x)$.
Lemma~\ref{lem:graham} permits positive, zero and negative $c_*(a)$.

For the entire solution, take the initial datum
$w_j(0)=\ind_{\{j\le0\}}$. Comparison with the scalar solutions
constructed above gives $W_j^-(t)\le w_j(t)\le W_j^+(t)$ for
$j\in\Z$ and $t\ge0$. Comparison with the period translate also
gives $w_{j+P}(t)\le w_j(t)$.
Let $t_n=n$ for $n\in\N$ and
$l_n=P\lfloor m_a(t_n)/P\rfloor$. Since
$l_n\in P\Z$, the functions
$w^{(n)}_j(t)=w_{j+l_n}(t+t_n)$ solve exactly \eqref{eq:lattice}
for $t>-t_n$, with the original reaction $f_j$ at index $j$.
Lemma~\ref{lem:compactness} gives a subsequence converging to an
entire solution $W(t,j)$, uniformly on bounded time intervals at
each fixed index. Since $0\le m_a(t_n)-l_n<P$, pass to a further
subsequence on which this difference converges, and relabel the
resulting subsequence by $n$. For bounded $t$,
\[
 m_a(t_n+t)-l_n
 =c_*(a)t+m_a(t_n)-l_n+o(1)\quad\text{ as } n\to\infty
\]
uniformly. Fix $\theta\in(0,1)$. Use the lower scalar comparison at
level $(1+\theta)/2$ and the upper scalar comparison at level
$\theta/2$. Their level-position constants are independent of $n$
and $t$, once $t_n+t$ is sufficiently large. The bound $0\le m_a(t_n)-l_n<P$ and the preceding asymptotic
identity give a constant $C_\theta$ such that, for each fixed
$t\in\R$ and $j\in\Z$, passage to the local limit yields
\[
 \begin{aligned}
 W(t,j)&\ge(1+\theta)/2>\theta
       &&\text{if }j\le c_*(a)t-C_\theta,\\
 W(t,j)&\le\theta/2<\theta
       &&\text{if }j>c_*(a)t+C_\theta.
 \end{aligned}
\]
This proves \eqref{eq:entire-levels} with one constant for all
$t\in\R$. Applying these bounds at levels approaching one and zero
gives $W(t,j)\to1$ as $j-c_*(a)t\to-\infty$ and
$W(t,j)\to0$ as $j-c_*(a)t\to+\infty$, uniformly for $t\in\R$.
At every time $t-1$, these limits make both $W(t-1,\cdot)$ and
$1-W(t-1,\cdot)$ nonzero. Strict positivity in
Lemma~\ref{lem:comparison} applied over $[t-1,t]$ gives $W(t,j)>0$.
To prove $W(t,j)<1$, set $Z(t,j)=1-W(t,j)$ for
$(t,j)\in\R\times\Z$ and define the operator $\mathcal A$ on
$\ell^\infty(\Z)$ by
$(\mathcal Az)_j=D(z_{j-1}+qz_{j+1}-(q+1)z_j)$ for $j\in\Z$.
With the common Lipschitz constant $L$ in \eqref{eq:taylor},
the equation gives $\partial_tZ\ge\mathcal AZ-LZ$.
Comparison over $[t-1,t]$ gives
$Z(t,\cdot)\ge e^{-L}e^{\mathcal A}Z(t-1,\cdot)$.
The sequence $Z(t-1,\cdot)$ is nonnegative and nonzero. Every pair
$j,l\in\Z$ is joined by a finite path with positive jump
coefficients, so every entry of the kernel of $e^{\mathcal A}$ is
positive. Hence $Z(t,j)>0$ for every $j\in\Z$.
Thus $0<W<1$, and $W(t,j+P)\le W(t,j)$ follows by passage to the
limit in the corresponding inequality for $w^{(n)}$.

Fix $T>0$ and $0<\epsilon<1/2$.
For $|t|\le T$ and $t_n>2T$,
\[
 m_a(t_n+t)-l_n
 =c_*(a)t+m_a(t_n)-l_n
  -\frac{3}{2\lambda_*(a)}\log\left(1+\frac t{t_n}\right).
\]
Here $0\le m_a(t_n)-l_n<P$, and
$|\log(1+t/t_n)|\le2T/t_n$. The centers are therefore bounded
uniformly for $|t|\le T$. The scalar bounds at levels $\epsilon$
and $1-\epsilon$ give an integer $J\ge1$ such that, for all
sufficiently large $n$ and $|t|\le T$,
\[
 \begin{aligned}
 0\le w^{(n)}_j(t),W(t,j)&\le\epsilon&&\text{if }j\ge J,\\
 1-\epsilon\le w^{(n)}_j(t),W(t,j)&\le1&&\text{if }j\le-J.
 \end{aligned}
\]
The bounds for $W$ follow by passage to the local limit.
Thus $|w^{(n)}_j(t)-W(t,j)|\le\epsilon$ outside $[-J,J]$.
On that finite interval the local convergence is uniform for
$|t|\le T$, proving convergence in $\ell^\infty(\Z)$ uniformly
on $[-T,T]$.

\end{proof}

\subsection{Moving and stationary fronts}

\begin{theorem}\label{thm:pulsating}

Assume \eqref{eq:hyp}, $f(x,s)=f_{b(x)}(s)$ and $\mu>0$.
Equation \eqref{eq:lattice} admits a strictly decreasing pulsating front
for each $c\ge c_+$, including $c=c_+$. Every front of \eqref{eq:lattice} with state one as $j-ct\to-\infty$
and state zero as $j-ct\to+\infty$ has speed at least $c_+$. The lift of
each front with $c\ge c_+$ to the differential equation in
\eqref{eq:tree-cauchy} is
\[
 u(t,x)=\Phi_{b(x)}(b(x)-ct)
 \qquad t\in\R,\ x\in\Tq,
\]
with profiles satisfying \eqref{eq:profile} and
$w_{j+P}(t+P/c)=w_j(t)$. The profiles can be normalized by
$\Phi_0(0)=1/2$.

\end{theorem}
\begin{proof}

Since $\mu>0$, Lemma~\ref{lem:matrix} gives $c_+>0$, so every
speed considered here is nonzero. Lemma~\ref{lem:wave-construction}
gives the $C^1$ profiles, their strict decrease, both end limits,
and the exclusion of speeds below $c_+$. There is a unique
$z_0\in\R$ with $\Phi_0(z_0)=1/2$; replacing every profile by
$z\mapsto\Phi_i(z+z_0)$ preserves \eqref{eq:profile} and gives the
stated normalization. Put $w_j(t)=\Phi_j(j-ct)$.
For $b(x)=j$, one neighbor has coordinate $j-1$ and $q$ have
coordinate $j+1$, so
\[
 D\Delta_{\Tq}\bigl[w_{b(\cdot)}(t)\bigr](x)
 =D\bigl(w_{j-1}(t)+qw_{j+1}(t)-(q+1)w_j(t)\bigr)
 \qquad t\in\R,\ x\in\Tq,\ j=b(x).
\]
Since $w$ satisfies \eqref{eq:lattice}, this identity proves that
$u(t,x)=w_{b(x)}(t)$ solves the tree equation.
By $\Phi_{j+P}=\Phi_j$ and
$(j+P)-c(t+P/c)=j-ct$, we obtain
$w_{j+P}(t+P/c)=w_j(t)$.\qedhere

\end{proof}

\Needspace{7\baselineskip}
\begin{proof}[Proof of Theorem~\ref{thm:all-speeds}]

Lemma~\ref{lem:signed-fronts} proves the nonzero-speed assertions, and
Lemma~\ref{lem:stationary-fronts} proves the zero-speed assertion.
The sign identity in Lemma~\ref{lem:matrix} says that zero belongs to
$[c_+,\infty)$ exactly when $\mu\le0$. For each resulting coordinate solution $w=(w_j(t))_{j\in\Z}$,
regard a stationary sequence as independent of $t$. The one neighbor
at $b(x)-1$ and the $q$ neighbors at $b(x)+1$ give, for $t\in\R$,
\[
 \begin{aligned}
 D\Delta_{\Tq}[w_{b(\cdot)}(t)](x)+f_{b(x)}(w_{b(x)}(t))
 &=D(w_{b(x)-1}(t)+qw_{b(x)+1}(t)-(q+1)w_{b(x)}(t))\\
 &\quad+f_{b(x)}(w_{b(x)}(t))=\dot w_{b(x)}(t).
 \end{aligned}
 \qquad x\in\Tq.
\]
Hence $u(t,x)=w_{b(x)}(t)$ solves the differential equation in
\eqref{eq:tree-cauchy} for every stated speed.\qedhere

\end{proof}

\begin{remark}

For negative speed, the transition moves towards decreasing Busemann
coordinate with state one on the left and state zero on the right.
Since stationary and moving fronts have coordinate values in $(0,1)$
at every $j\in\Z$, their lifts are positive at every vertex of $\Tq$. For $\mu\le0$, these infinite-support solutions coexist with
the uniform extinction of the finitely supported solutions in
Theorem~\ref{thm:threshold}.
At $\mu=0$, the stationary sequence in
\eqref{eq:stationary-front} is constructed in
Lemma~\ref{lem:stationary-fronts}.

\end{remark}

\subsection{Stability}

\begin{proof}[Proof of Theorem~\ref{thm:wave-stability}]

Translation of the profile argument preserves \eqref{eq:profile},
so $\Phi_j(j-ct+A)$ and $\Phi_j(j-ct-A)$ solve
\eqref{eq:lattice}. The initial bracket and
Lemma~\ref{lem:comparison} give
\[
 \Phi_j(j-ct+A)\le u_j(t)\le\Phi_j(j-ct-A)
 \quad\text{ for }j\in\Z,\ t\ge0.
\]
In particular, $u$ and $W$ take values in $[0,1]$.
Apply Lemma~\ref{lem:difference} to the pair $(u,W)$ with the
exponent $\lambda=\lambda_+$, the given $1\le p<\infty$, and the
front speed $c$. The weighted initial error
$e^{\lambda_+j}(u_j(0)-W_j(0))$ has $\ell^p(\Z)$ norm $B_p$, and
\[
 k(\lambda_+)-c\lambda_+
 =-\lambda_+(c-c_+).
\]
This proves \eqref{eq:wave-stability-weighted}. In particular,
\[
\begin{gathered}
 \sup_{\substack{j\in\Z\\|j-ct|\le R}}|u_j(t)-W_j(t)|
 \le C_p e^{\lambda_+R}(1+t)^{-1/(2p)}
       e^{-\lambda_+(c-c_+)t}B_p\to0
 \quad\text{ as } t\to\infty
 \\\text{ for each fixed }R>0,
\end{gathered}
\]
since $p<\infty$ and $c\ge c_+$: the exponential factor is at most
one, and $(1+t)^{-1/(2p)}\to0$ also at $c=c_+$. Given
$\varepsilon>0$, the end limits of the finitely many profiles allow
$R$ to be chosen so that both profile translates and $W$ lie in
$[1-\varepsilon,1]$ for $j-ct<-R$
and in $[0,\varepsilon]$ for $j-ct>R$. The preserved profile bounds
$\Phi_j(j-ct+A)\le u_j(t)\le\Phi_j(j-ct-A)$ then bound
$|u_j(t)-W_j(t)|$ by $\varepsilon$ on these two regions,
uniformly in $t$. Together with the convergence on
$|j-ct|\le R$, this proves \eqref{eq:wave-stability}. If $u_j(0)=W_j(0)$ except at finitely many indices, then $B_p<\infty$.
If every modified value belongs to $(0,1)$, strict monotonicity and the
two limits of each profile give one common $A$ for the initial bracket.\qedhere

\end{proof}

\begin{proof}[Proof of Theorem~\ref{thm:tree-wave-stability}]

The sequences $j\mapsto w_{j\pm mP}(t)$ solve
\eqref{eq:lattice}, because $f_{j\pm mP}=f_j$. Each vertex has one neighbor at Busemann coordinate $b(x)-1$
and $q$ neighbors at coordinate $b(x)+1$. By \eqref{eq:lattice},
the corresponding functions on the tree therefore solve
\eqref{eq:tree-cauchy}. If $c\ne0$, periodicity gives
\[
 w_{j\pm mP}(t)=\Phi_j(j-ct\pm mP)
 \qquad j\in\Z,\ t\in\R.
\]
Their order with $w_j(t)$ follows from $\Phi_j'<0$ for either sign
of $c$. If $c=0$, the order follows from $p_{j+P}<p_j$.
Thus \eqref{eq:tree-front-bracket} and comparison give
\[
 w_{b(x)+mP}(t)\le u(t,x)\le w_{b(x)-mP}(t)
 \quad\text{ for }x\in\Tq,\ t\ge0.
\]
The end limits of the finitely many profiles, or the sequence limits
in \eqref{eq:stationary-front} when $c=0$, show that both functions
$w_{b(x)\pm mP}(t)$ tend to one as $b(x)-ct\to-\infty$ and to
zero as $b(x)-ct\to+\infty$, uniformly in time.

Apply Lemma~\ref{lem:tree-difference} to the pair $(u,W)$ with
$\lambda=\lambda_+$, the given $1\le p<\infty$, and speed $c$.
The weighted initial error
$e^{\lambda_+b(x)}(u_0(x)-W(0,x))$ has $\ell^p(\Tq)$ norm $B_p$,
with counting measure on all vertices. Since
$k(\lambda_+)-c\lambda_+=-\lambda_+(c-c_+)$, estimate
\eqref{eq:tree-difference} gives the weighted bound. On
$|b(x)-ct|\le R$, \eqref{eq:tree-difference} bounds the unweighted error by
$C_p e^{\lambda_+R}(1+t)^{-1/(2p)}
e^{-\lambda_+(c-c_+)t}B_p$, which tends to zero as $t\to\infty$
because Theorem~\ref{thm:all-speeds} gives $c\ge c_+$ and $p<\infty$.
This includes $c=c_+$, where the polynomial factor tends to zero,
and $c=0$, which is admissible only when $c_+\le0$.
Given $\varepsilon>0$, the two end limits allow $R$ to be chosen
so that both bracket functions and $W$ lie in $[1-\varepsilon,1]$
on $b(x)-ct<-R$ and in $[0,\varepsilon]$ on $b(x)-ct>R$,
uniformly in $t$. The bracket bounds the error by $\varepsilon$
on both regions. Combining this with convergence on the strip proves
\eqref{eq:tree-front-stability}.

For a finite vertex modification with all new values in $(0,1)$,
$B_p$ is finite. At each modified vertex $x$, $w_{b(x)+mP}(0)\to0$ and
$w_{b(x)-mP}(0)\to1$ as $m\to\infty$. Strict inequalities at the finite collection of
new values therefore allow one common $m$. At every unmodified vertex
the bracket holds by monotonicity. Thus the finite modification
satisfies both the bracket \eqref{eq:tree-front-bracket} and the
condition $B_p<\infty$.\qedhere

\end{proof}

\begin{remark}
The weighted norm in Theorem~\ref{thm:tree-wave-stability} uses
counting measure on all vertices. If $u_0-W(0,\cdot)$ takes a fixed nonzero value on one horocycle, then
$B_p=\infty$ for every $1\le p<\infty$: the weight is constant on that
horocycle, which contains infinitely many vertices. Finite vertex
modifications can vary within a horocycle.
Theorem~\ref{thm:wave-stability} uses the coordinate norm on $\Z$,
while Theorem~\ref{thm:tree-wave-stability} uses the full-tree norm.
\end{remark}

\begin{remark}\label{rem:stability-scope}
Theorems~\ref{thm:wave-stability} and~\ref{thm:tree-wave-stability}
give convergence to the front appearing in their initial-data brackets
and weighted norm assumptions. Selection of a front from finitely
supported data, or from data satisfying
$\ind_{\{j\le A\}}\le w_j(0)\le\ind_{\{j\le B\}}$ for integers $A\le B$,
is formulated as a further question in Section~\ref{sec:examples}. For traveling-wave convergence in hyperbolic space, see \cite{GGQ}.
\end{remark}

\subsection{Distance from a geodesic}

\begin{proof}[Proof of Theorem~\ref{thm:geodesic-threshold}]

Write $\mathcal L_{\mathcal G}=D\Delta_{\Tq}+a_{s(\cdot)}$ on
$\ell^2(\Tq)$ with counting measure. The bounded real potential and
symmetric tree adjacency make this operator bounded and self-adjoint.
The same coordinate formula defines a bounded operator on
$\ell^\infty(\Tq)$, with norm at most $2D(q+1)+a_{\max}$.
The exponential $e^{t\mathcal L_{\mathcal G}}$ therefore acts on
$\ell^\infty(\Tq)$ for $t\ge0$. The bound
$f_{s(x)}(v)\le a_{s(x)}v$ for $0\le v\le1$ in \eqref{eq:hyp} and
Lemma~\ref{lem:comparison} give
$0\le u(t,x)\le(e^{t\mathcal L_{\mathcal G}}u_0)(x)$ for every
$t\ge0$ and $x\in\Tq$.

For the extinction assertions, assume $u_0\in\ell^2(\Tq)$.
The exponential series for $\mathcal L_{\mathcal G}$ on
$\ell^2(\Tq)$ and $\ell^\infty(\Tq)$ agree at each vertex.
By Lemma~\ref{lem:geodesic-spectrum},
$\sup\spec(\mathcal L_{\mathcal G})=\sigma_{\mathcal G}$.
For $\sigma_{\mathcal G}<0$, the spectral theorem gives
\[
 \|e^{t\mathcal L_{\mathcal G}}u_0\|_2
 \le e^{\sigma_{\mathcal G}t}\|u_0\|_2
 \qquad t\ge0.
\]
Together with linear comparison and
$\|v\|_{\ell^\infty(\Tq)}\le\|v\|_{\ell^2(\Tq)}$ for
$v\in\ell^2(\Tq)$, this gives
$\|u(t,\cdot)\|_\infty\le e^{\sigma_{\mathcal G}t}\|u_0\|_2$
for $t\ge0$.

When $\sigma_{\mathcal G}=0$, let $\nu_{u_0}$ be the scalar spectral
measure of $u_0\in\ell^2(\Tq)$ for $\mathcal L_{\mathcal G}$.
By \eqref{eq:geodesic-spectrum}, it is supported in $(-\infty,0]$
and has no atom at zero. Since its total mass is $\|u_0\|_2^2$,
the bound $e^{2tE}\le1$ for $E\le0$ and convergence to zero for $E<0$
give
\[
 \|e^{t\mathcal L_{\mathcal G}}u_0\|_2^2
 =\int_{(-\infty,0]}e^{2tE}\,\dd\nu_{u_0}(E)\to0
 \quad\text{ as }t\to\infty
\]
by dominated convergence. Linear comparison then gives
$\|u(t,\cdot)\|_\infty\le\|e^{t\mathcal L_{\mathcal G}}u_0\|_2\to0$.

\emph{A positive stationary solution and its minimality.}
Suppose $\sigma_{\mathcal G}>0$. Since $\mathcal L_{\mathcal G}$
is bounded and finitely supported vectors are dense in $\ell^2(\Tq)$,
its Rayleigh supremum can be approximated by a nonzero finitely
supported vector $z$ with
$\langle z,\mathcal L_{\mathcal G}z\rangle>0$.
Enlarge its support to a finite connected set $\Omega$. The
zero-exterior operator $\mathcal L_\Omega$ on $\ell^2(\Omega)$ is
\[
 (\mathcal L_\Omega v)(x)
 =D\sum_{\substack{y\in\Omega\\y\sim x}}v(y)
   +(a_{s(x)}-D(q+1))v(x),\quad\text{ for }x\in\Omega.
\]
Thus the quadratic forms of $\mathcal L_\Omega$ and
$\mathcal L_{\mathcal G}$ agree on vectors supported in $\Omega$. The largest
eigenvalue $\eta$ of $\mathcal L_\Omega$ is positive by the Rayleigh
principle. The matrix $\mathcal L_\Omega+D(q+1)I$ has diagonal entries
$a_{s(x)}>0$ and positive entries $D$ along the edges of the connected
set $\Omega$. It is therefore nonnegative and irreducible, so an
eigenvector $\varphi$ for $\eta$ can be chosen strictly positive
on $\Omega$. Extend $\varphi$ to $\Tq$ by zero and normalize it by
$\max_{x\in\Omega}\varphi(x)=1$. Choose
$0<\varepsilon_0\le\min\{1/2,\eta/[2(1+C_2)]\}$, so that
$\varepsilon_0<1$ and $C_2\varepsilon_0<\eta$, and put
$z_0=\varepsilon_0\varphi$. For every $0<\vartheta\le1$, the lower Taylor bound in
\eqref{eq:taylor} and $\mathcal L_\Omega\varphi=\eta\varphi$
give
\begin{equation}\label{eq:geodesic-finite-subsolution}
 D\Delta_{\Tq}(\vartheta z_0)(x)
   +f(x,\vartheta z_0(x))
 \ge\vartheta z_0(x)(\eta-C_2\vartheta z_0(x))>0
 \quad\text{ for }x\in\Omega.
\end{equation}
Outside $\Omega$, $\vartheta z_0=0$ and every neighboring
contribution to its diffusion is nonnegative. Thus $\vartheta z_0$
is a stationary subsolution on $\Tq$. Denote the solution of
\eqref{eq:tree-cauchy} with initial value $z_0$ by $u^{z_0}$.
Comparison with the stationary subsolution gives
$u^{z_0}(s,\cdot)\ge z_0$ for $s\ge0$. Autonomy and comparison
then give $u^{z_0}(t+s,x)\ge u^{z_0}(t,x)$ for $s,t\ge0$;
comparison with the constant solution one gives $u^{z_0}(t,x)\le1$
for $t\ge0$ and $x\in\Tq$. Set
$p_*(x)=\lim_{t\to\infty}u^{z_0}(t,x)$ for $x\in\Tq$.
For fixed $x\in\Tq$ and $T>0$, write the integral equation as
$u^{z_0}(t+T,x)-u^{z_0}(t,x)
=\int_0^T[D\Delta_{\Tq}u^{z_0}(t+\tau,\cdot)(x)
+f(x,u^{z_0}(t+\tau,x))]\,\dd\tau$.
The difference on the left tends to zero as $t\to\infty$.
Bounded convergence applies to the reaction and finite neighbor sum
on the right, giving the limit
$T[D\Delta_{\Tq}p_*(x)+f(x,p_*(x))]$. Hence
$D\Delta_{\Tq}p_*(x)+f(x,p_*(x))=0$ for every $x\in\Tq$. Strict positivity of the evolution at any positive time
and monotonicity give $p_*(x)>0$ at every vertex.

To prove minimality, let $p$ be any stationary solution of the
differential equation in \eqref{eq:tree-cauchy} with $0<p\le1$.
Since $\Omega$ is finite and $z_0>0$ on $\Omega$, the number
\[
 t_*=
 \min_{x\in\Omega}\frac{p(x)}{z_0(x)}
\]
is positive.  If $t_*<1$, choose $x_*\in\Omega$ with
$p(x_*)=t_*z_0(x_*)$.
Then
\[
 p(x)-t_*z_0(x)\ge0\quad(x\in\Tq),\qquad
 (p-t_*z_0)(x_*)=0
\]
and therefore
\[
 \Delta_{\Tq}(p-t_*z_0)(x_*)\ge0.
\]
At the contact point the reaction values agree:
$f(x_*,p(x_*))=f(x_*,t_*z_0(x_*))$.  Subtracting \eqref{eq:geodesic-finite-subsolution}, with
$\vartheta=t_*$, from the stationary equation for $p$ therefore gives
\[
 D\Delta_{\Tq}(p-t_*z_0)(x_*)<0,
\]
{a contradiction. Hence $t_*\ge1$, so $p\ge z_0$. Comparison gives
$p(x)\ge u^{z_0}(t,x)$ for every $t\ge0$ and $x\in\Tq$.
Passing to the limit $t\to\infty$ yields $p\ge p_*$.
Thus $p_*$ is the smallest stationary solution of the differential
equation in \eqref{eq:tree-cauchy} with values in $(0,1]$.}

For $0<\vartheta\le1$, denote the solution of
\eqref{eq:tree-cauchy} from $\vartheta z_0$ by $u^{\vartheta z_0}$.
Inequality \eqref{eq:geodesic-finite-subsolution} and comparison make
$u^{\vartheta z_0}(t,x)$ nondecreasing for $t\ge0$ and bounded by one
for every $x\in\Tq$.
Define $P_\vartheta(x)=\lim_{t\to\infty}u^{\vartheta z_0}(t,x)$
for $x\in\Tq$. The monotone-limit conclusion of
Lemma~\ref{lem:compactness} gives
$D\Delta_{\Tq}P_\vartheta(x)+f(x,P_\vartheta(x))=0$ for $x\in\Tq$.
Strict positivity at any positive time gives $P_\vartheta(x)>0$
for every $x$. Since $\vartheta z_0\le z_0\le p_*$, comparison
gives $P_\vartheta\le p_*$, while minimality gives
$p_*\le P_\vartheta$. Thus $P_\vartheta=p_*$ and
$u^{\vartheta z_0}(t,x)\to p_*(x)$ as $t\to\infty$ at every vertex.
The data $\vartheta z_0$ are finitely supported and tend to zero in
$\ell^\infty(\Tq)$ as $\vartheta\downarrow0$, proving the assertion about
arbitrarily small data.

\emph{Persistence and coordinate symmetry.}
For a solution with nonzero initial datum, strict positivity at time
one and finiteness of $\Omega$ give
$\vartheta:=\min\{1,\min_{x\in\Omega}u(1,x)/z_0(x)\}>0$.
Thus $u(1,\cdot)\ge\vartheta z_0$ on the whole tree.
Comparison gives $u(1+t,x)\ge u^{\vartheta z_0}(t,x)$ for
$t\ge0$ and $x\in\Tq$. Since $u^{\vartheta z_0}(t,x)\to p_*(x)$
as $t\to\infty$, this proves \eqref{eq:geodesic-minimal-persistence}.
The tree automorphisms mapping $\mathcal G$ onto itself act
transitively on each set $\{x\in\Tq:s(x)=r\}$, $r\in\N_0$:
one can translate
the line and then identify the attached regular rooted branches.
If $g$ is a tree automorphism preserving $\mathcal G$, then
$s(gx)=s(x)$, so $p_*\circ g$ is a positive stationary solution.
Minimality gives $p_*(x)\le p_*(gx)$. Applying minimality to
$g^{-1}$ at $gx$ gives $p_*(gx)\le p_*(x)$.
Thus $p_*(gx)=p_*(x)$ for every such $g$, and transitivity on
$\{x:s(x)=r\}$ makes $p_*$ depend only on $s(x)$.
When $\mu>0$, Lemma~\ref{lem:geodesic-local} applied to the increasing
solution from $z_0$ implies $p_*\equiv1$.

\emph{Transverse confinement.}
Suppose $\mu\le0$. Put $h_r=q^{-r/2}\phi_r(\lambda_0)$ for
$r\in\N_0$, with $\phi_r(\lambda_0)$ extended periodically in $r$.
For $r\ge1$, its eigenvector equation is
\[
 D(h_{r-1}+qh_{r+1}-(q+1)h_r)+a_rh_r=\mu h_r.
\]
Lemma~\ref{lem:decreasing-ground} gives $h_{r+1}<h_r$, and
boundedness of the phase vector gives $h_r\to0$ as $r\to\infty$.
For $Ch_0\ge1$, set
$\overline p_C(x)=\min\{1,Ch_{s(x)}\}$. At an off-geodesic vertex
where $Ch_{s(x)}<1$, the bound $f_{s(x)}(v)\le a_{s(x)}v$
from \eqref{eq:hyp} and truncation of neighboring values give
\[
 D\Delta_{\Tq}\overline p_C(x)+f(x,\overline p_C(x))
 \le C\mu h_{s(x)}\le0
 \qquad x\in\Tq,\ s(x)\ge1,\ Ch_{s(x)}<1.
\]
At a vertex $x$ where $\overline p_C(x)=1$, the reaction vanishes and diffusion
is nonpositive. This includes every vertex of $\mathcal G$, with its
boundary coefficient $D(q-1)$. Thus $\overline p_C$ is a bounded
stationary supersolution. For any given finitely supported datum, choose $C\ge1/h_0$ with
$Ch_{s(x)}\ge u_0(x)$ on its finite support. Then
$u_0\le\overline p_C$, and comparison gives
$u(t,x)\le\min\{1,Ch_{s(x)}\}$ for $t\ge0$, proving
\eqref{eq:geodesic-transverse-confinement}.
If $\sigma_{\mathcal G}>0$, choose $0<\vartheta\le1$ so small that
$\vartheta z_0\le\overline p_{1/h_0}$ on $\Omega$.
Comparison with this stationary supersolution and convergence of the
solution from $\vartheta z_0$ give
\[
 p_*(x)\le\overline p_{1/h_0}(x)
       =\min\{1,h_{s(x)}/h_0\}\le h_{s(x)}/h_0
 \qquad x\in\Tq.
\] The bound $p_*(x)\le h_{s(x)}/h_0$ tends to zero as
$s(x)\to\infty$. If $p_*(x_0)=1$ at some vertex $x_0$,
stationarity and $f(x_0,1)=0$ give
$\sum_{y\sim x_0}(p_*(y)-1)=0$. Every summand is nonpositive,
so every neighbor also has value one. Repeating along finite paths
would give $p_*\equiv1$, contrary to
$p_*(x)\le h_{s(x)}/h_0\to0$ as $s(x)\to\infty$.
Hence $0<p_*(x)<1$ at every vertex.\qedhere

\end{proof}

\begin{proof}[Proof of Theorem~\ref{thm:geodesic-convergence}]

For real sequences $U=(U_r)_{r\in\N_0}$, use the coordinate
diffusion operator
\[
 (\mathcal DU)_0=D(q-1)(U_1-U_0),\qquad
 (\mathcal DU)_r=D(U_{r-1}+qU_{r+1}-(q+1)U_r),\quad\text{ for }r\in\N,
\]
and the weights $\nu_0=1$, $\nu_r=(q-1)q^{r-1}$ for $r\ge1$.
Thus $\mathcal H_\nu$ consists of real sequences with
$\sum_{r\ge0}\nu_rU_r^2<\infty$.
The row sums of $M_0$ are $a_i$, so its principal eigenvalue satisfies
$k(0)\ge a_{\min}>0$ by comparison with the positive constant vector.
Strict convexity, $k(\lambda_0)=\mu<0$, and
$k(\lambda)\to\infty$ as $\lambda\to\infty$ therefore give exactly
two positive zeros of $k$. Since $\lambda\in(\lambda_0,\lambda_{\mathrm R})$ as in the
statement, we have $k(\lambda)<0$. Put $h_r=e^{-\lambda r}\phi_r(\lambda)$ for
$r\ge0$. For a constant
$C$ with $Ch_0\ge1$, define
\[
 \overline U_r^C=\min\{1,Ch_r\}
 \qquad r\in\N_0.
\]
The sequence $\overline U^C$ is a stationary supersolution of
$(\mathcal DU)_r+f_r(U_r)=0$ for $r\in\N_0$. At an
index $r\ge1$ with value less than one, the KPP inequality and the
truncation of neighboring values give a right-hand side at most
$Ck(\lambda)h_r<0$. At an index with value one the right-hand side is
nonpositive. This includes $r=0$, where the correct expression is
$D(q-1)(\overline U_1^C-1)$. Also
\[
 \sum_{r\ge0}\nu_r(\overline U_r^C)^2<\infty,
\]
because, for $r\ge1$,
$\nu_r(\overline U_r^C)^2\le C^2(q-1)q^{-1}
\max_i\phi_i(\lambda)^2 e^{-(2\lambda-\log q)r}$ and
$2\lambda>\log q$. The resulting geometric series converges.

Let $(U_r^C(t))_{r\ge0}$ solve
$\dot U_r^C=(\mathcal DU^C)_r+f_r(U_r^C)$ for $r\ge0$, with
$U_r^C(0)=\overline U_r^C$. Comparison gives $U_r^C(t)\le\overline U_r^C$ for $r\in\N_0$ and $t\ge0$,
and the semigroup property gives $U_r^C(t+s)\le U_r^C(t)$
for $s,t\ge0$. Define $P_r^C=\lim_{t\to\infty}U_r^C(t)$.
For every fixed $r\in\N_0$ and $T>0$,
\[
 U_r^C(t+T)-U_r^C(t)
 =\int_0^T\bigl[(\mathcal DU^C(t+s))_r
                         +f_r(U_r^C(t+s))\bigr]\,\dd s
 \qquad t\ge0.
\]
As $t\to\infty$, the difference $U_r^C(t+T)-U_r^C(t)$ tends
to zero. Bounded convergence in the displayed integral, using the
bounded reaction and finite neighbor sum, gives
$0=T[(\mathcal DP^C)_r+f_r(P_r^C)]$.
Thus $P^C$ is stationary and belongs to $\mathcal H_\nu$ by
$0\le P_r^C\le\overline U_r^C$.
The initial tree function $x\mapsto\overline U_{s(x)}^C$ is
nonzero. Applying Theorem~\ref{thm:geodesic-threshold} to its
solution gives $P^C_{s(x)}\ge p_*(x)>0$. Since $p_*$ depends only on $s(x)$ and its coordinate values are
bounded above by $P^C\in\mathcal H_\nu$, the coordinate sequence
of $p_*$ also belongs to $\mathcal H_\nu$. Lemma~\ref{lem:weighted-steady-uniqueness}
gives $P_r^C=p_*(x)$ whenever $s(x)=r$. In particular,
$p_*(x)\le\overline U_{s(x)}^C$ gives the stated exponential bound. If $V\in\mathcal H_\nu$ satisfies
$(\mathcal DV)_r+f_r(V_r)=0$ and $0<V_r\le1$ for every $r$,
Lemma~\ref{lem:weighted-steady-uniqueness}, applied to $V$ and
the coordinate sequence of $p_*$, gives $V_r=p_*(x)$ whenever $s(x)=r$.

Given the nonzero datum $u_0:\Tq\to[0,1]$ with
$u_0(x)\le C_0e^{-\lambda s(x)}$ for $x\in\Tq$ as in the statement,
choose
\[
 C\ge\max\left\{\frac1{h_0},
                  \frac{C_0}{\min_i\phi_i(\lambda)}\right\}.
\]
Then $Ch_{s(x)}\ge C_0e^{-\lambda s(x)}\ge u_0(x)$.
Together with $u_0\le1$, this gives
$\overline U_{s(x)}^C=\min\{1,Ch_{s(x)}\}\ge u_0(x)$. The tree comparison principle gives
\[
 u(t,x)\le U_{s(x)}^C(t)\to p_*(x)\quad\text{ as } t\to\infty
 \quad\text{for each fixed }x\in\Tq.
\]
Together with
$\liminf_{t\to\infty}u(t,x)\ge p_*(x)$ from
\eqref{eq:geodesic-minimal-persistence}, the upper limit above gives
$u(t,x)\to p_*(x)$ as $t\to\infty$ at every fixed vertex. Both $u(t,x)$ and $p_*(x)$
are bounded by $Ce^{-\lambda s(x)}\max_i\phi_i(\lambda)$, uniformly
in $t$. For fixed $M\in[0,\infty)$ and $R\in[0,\infty)$, the set
$\{x\in\Tq:|\ell_{\mathcal G}(x)|\le M,\ s(x)\le R\}$ is finite. Pointwise convergence is uniform on this finite set. On the
remaining part of the strip $\{|\ell_{\mathcal G}(x)|\le M\}$,
we have $s(x)>R$, and the exponential estimate bounds the error by $2C\max_i\phi_i(\lambda)e^{-\lambda R}$.
Letting first $t\to\infty$ and then $R\to\infty$ proves
\eqref{eq:geodesic-local-convergence}.

Let $p$ be a stationary solution of the differential equation in
\eqref{eq:tree-cauchy} with $0<p\le1$ and
$\sup_x e^{\lambda s(x)}p(x)<\infty$. Choose
\[
 C\ge\max\left\{\frac1{h_0},
       \frac{\sup_xe^{\lambda s(x)}p(x)}
            {\min_i\phi_i(\lambda)}\right\}.
\]
Then $p(x)\le\min\{1,Ch_{s(x)}\}=\overline U_{s(x)}^C$.
The solution from $p$ is stationary, so tree comparison gives
$p(x)\le U_{s(x)}^C(t)$ for every $t\ge0$.
Letting $t\to\infty$ yields $p\le p_*$.
Minimality gives $p_*\le p$, and therefore $p=p_*$.

Finally let $u_0:\Tq\to[0,1]$ be nonzero and have finite support. The mass estimate
\eqref{eq:finite-mass-level-count} gives $u(t,\cdot)\in\ell^1(\Tq)$ at
every finite time. Hence $u(t,\gamma(n))\to0$ as $|n|\to\infty$.
On the other hand, $p_*(\gamma(n))=p_*(\gamma(0))>0$ by its coordinate
symmetry. Taking the supremum along $\mathcal G$ proves
\eqref{eq:no-global-geodesic-convergence}.\qedhere

\end{proof}

Assume \eqref{eq:hyp}, \eqref{eq:strict-quotient}, and
$f(x,u)=f_{s(x)}(u)$. Under $\mu<0<\sigma_{\mathcal G}$,
Theorem~\ref{thm:geodesic-convergence} proves uniqueness
of the stationary coordinate sequence in $\mathcal H_\nu$.
Theorem~\ref{thm:geodesic-critical-convergence} proves uniqueness among
positive stationary tree solutions $p\le1$ satisfying
$\sup_x p(x)/(q^{-s(x)/2}\phi_{s(x)}(\lambda_0))<\infty$ when
$\mu\le0<\sigma_{\mathcal G}$. At $\mu=0$, the coordinate sequence of
$p_*$ satisfies this bound and lies outside $\mathcal H_\nu$ by
Lemma~\ref{lem:critical-stationary-tail}.

\begin{proof}[Proof of Theorem~\ref{thm:geodesic-critical-convergence}]

Here $h_r=q^{-r/2}\phi_r(\lambda_0)$ for $r\in\N_0$, as in
the statement. We use the coordinate diffusion operator with
$(\mathcal DU)_0=D(q-1)(U_1-U_0)$ and
$(\mathcal DU)_r=D(U_{r-1}+qU_{r+1}-(q+1)U_r)$ for $r\ge1$.
For $C\ge1/h_0$ and $r\in\N_0$, let
$\overline U_r^C=\min\{1,Ch_r\}$.
For $r\ge1$, the eigenvector equation gives
$(\mathcal Dh)_r+a_rh_r=\mu h_r$. At an index with $Ch_r<1$,
the KPP bound and truncation of neighboring values give
\[
 (\mathcal D\overline U^C)_r+f_r(\overline U_r^C)
 \le C\bigl((\mathcal Dh)_r+a_rh_r\bigr)
 =C\mu h_r\le0
 \qquad r\in\N,\ Ch_r<1.
\]
At an index with $\overline U_r^C=1$, the reaction vanishes and
$(\mathcal D\overline U^C)_r\le0$.
In particular $\overline U_0^C=1$ and the boundary right-hand side is
$D(q-1)(\overline U_1^C-1)\le0$.
Thus $\overline U^C$ is a bounded stationary supersolution.

Let $U^C$ solve
$\dot U_r^C=(\mathcal DU^C)_r+f_r(U_r^C)$ for $r\ge0$, with
$U_r^C(0)=\overline U_r^C$. Comparison with this stationary
supersolution and time translation give
$0\le U_r^C(t+s)\le U_r^C(t)\le\overline U_r^C$ for $s,t\ge0$.
Define $P_r^C=\lim_{t\to\infty}U_r^C(t)$.
For fixed $r\in\N_0$ and $T>0$, the reaction $f_r$ is bounded on $[0,1]$
and each row of $\mathcal D$ has finitely many terms. Bounded convergence in
\[
 U_r^C(t+T)-U_r^C(t)
 =\int_0^T\bigl[(\mathcal DU^C(t+s))_r
                         +f_r(U_r^C(t+s))\bigr]\,\dd s
 \qquad t\ge0.
\]
gives $0=T[(\mathcal DP^C)_r+f_r(P_r^C)]$.
Hence $P^C$ is stationary and $P_r^C\le Ch_r$.
Applying Theorem~\ref{thm:geodesic-threshold} to the solution
$x\mapsto U_{s(x)}^C(t)$ gives $P^C_{s(x)}\ge p_*(x)>0$.
The bound \eqref{eq:geodesic-minimal-tail} gives
$p_*(x)\le h_{s(x)}/h_0$, and Theorem~\ref{thm:geodesic-threshold}
shows that $p_*$ depends only on $s(x)$.
The sequence $P^C$ and the coordinate values of $p_*$ are therefore
positive stationary sequences bounded by a constant multiple of $h$.
Lemma~\ref{lem:cutoff-steady-uniqueness} gives
$P_r^C=p_*(x)$ whenever $s(x)=r$.

For the nonzero initial datum satisfying $u_0(x)\le C_0h_{s(x)}$
for $x\in\Tq$, choose $C\ge\max\{C_0,1/h_0\}$. Comparison gives
$u(t,x)\le U^C_{s(x)}(t)$, whose limit as $t\to\infty$ is $p_*(x)$.
The pointwise lower bound in
\eqref{eq:geodesic-minimal-persistence} proves $u(t,x)\to p_*(x)$ as $t\to\infty$.
Moreover, comparison and $p_*(x)\le h_{s(x)}/h_0$ give
$u(t,x),p_*(x)\le Ch_{s(x)}$ for every $t\ge0$ and $x\in\Tq$,
because $C\ge C_0$ and $C\ge1/h_0$. Since $h_r\to0$ as $r\to\infty$, this bound
controls vertices with $s(x)>R$. For fixed $M,R\in[0,\infty)$,
the set $\{x\in\Tq:|\ell_{\mathcal G}(x)|\le M,\ s(x)\le R\}$
is finite. On $\{x\in\Tq:|\ell_{\mathcal G}(x)|\le M,\ s(x)>R\}$,
the error $|u(t,x)-p_*(x)|$ is bounded uniformly for $t\ge0$ by
$2C\sup_{r>R}h_r$, which tends to zero as $R\to\infty$.
Pointwise convergence on the finite set therefore proves
\eqref{eq:critical-geodesic-convergence}.

Let $p$ be a stationary solution of the differential equation in
\eqref{eq:tree-cauchy} with $0<p\le1$ and
$\sup_xp(x)/h_{s(x)}<\infty$. Choose
$C\ge\max\{1/h_0,\sup_xp(x)/h_{s(x)}\}$.
Then $p(x)\le\min\{1,Ch_{s(x)}\}=\overline U_{s(x)}^C$.
Since the tree solution from $p$ equals $p$ at every time,
comparison gives $p(x)\le U^C_{s(x)}(t)$ for $t\ge0$.
The limit $U^C_{s(x)}(t)\to p_*(x)$ yields $p\le p_*$;
minimality gives $p_*\le p$. Hence $p=p_*$.
For finitely supported data, \eqref{eq:finite-mass-level-count}
gives $u(t,\cdot)\in\ell^1(\Tq)$ at each finite time, and hence
$u(t,\gamma(n))\to0$ as $|n|\to\infty$.
Since $p_*(\gamma(n))=p_*(\gamma(0))>0$, taking the supremum along
$\mathcal G$ gives \eqref{eq:no-global-geodesic-convergence}.
Finally, when $\mu=0$, write $p_{*,r}$ for the common value of
$p_*(x)$ on $\{x:s(x)=r\}$. This positive stationary sequence
satisfies $p_{*,r}\le h_r/h_0$ by
\eqref{eq:geodesic-minimal-tail}.
Lemma~\ref{lem:critical-stationary-tail} therefore gives
$p_{*,r}=A_*q^{-r/2}\phi_r(\lambda_0)+O(q^{-r})$ as
$r\to\infty$, with $A_*>0$, and shows that the coordinate sequence
$(p_{*,r})_{r\ge0}$ lies outside $\mathcal H_\nu$.
\end{proof}

The rates in \eqref{eq:geodesic-example-rates} also give an example satisfying
$\mu=0<\sigma_{\mathcal G}$. Start with the rates in \eqref{eq:geodesic-example-rates},
let $\mu_{\rm old}<0$ denote their value of $k(\lambda_0)$. Replace every
logistic growth rate $a_i$ by $a_i-\mu_{\rm old}$.
The matrix $M_{\lambda_0}$ and operator $\mathcal J_{\mathcal G}$
both shift by $-\mu_{\rm old}I$. Thus the new values satisfy
$\mu=0$ and $\sigma_{\mathcal G}>0$; all growth rates remain positive.
\subsection{Tangential spreading along the geodesic}

Throughout this subsection, assume \eqref{eq:hyp},
$f(x,u)=f_{s(x)}(u)$, and $\mu\le0<\sigma_{\mathcal G}$.
The positive stationary solution $p_*$ from
Theorem~\ref{thm:geodesic-threshold} tends to zero as
$s(x)\to\infty$, by \eqref{eq:geodesic-minimal-tail}.
We determine the speed at which solutions approach $p_*$ along
$\mathcal G$.

Recall that $\mu=\min_{\lambda\in\R}k(\lambda)$ is the upper
spectral endpoint of the full-line periodic Jacobi operator with
diagonal $a_r-D(q+1)$ and off-diagonal entries $D\sqrt q$.
The half-line operator $\mathcal J_{\mathcal G}$ includes the
boundary equation at $r=0$, and
$\sigma_{\mathcal G}=\sup\spec(\mathcal J_{\mathcal G})$.
Thus $\sigma_{\mathcal G}>\mu$ under the assumptions of this subsection.

Use counting measure on $\ell^2(\N_0)$, and let $\delta_r$
denote its unit vector at $r\in\N_0$.
For $\lambda\in\R$, put $\chi(\lambda)=2D(\cosh\lambda-1)$ and define
\begin{equation}
\label{eq:tangential-fiber}
 \mathcal J_{\parallel}(\lambda)
 =
 \mathcal J_{\mathcal G}
 +
 \chi(\lambda)
 \langle\delta_0,\,\cdot\,\rangle\delta_0,
 \qquad
 \Lambda_{\parallel}(\lambda)
 =
 \sup\spec(\mathcal J_{\parallel}(\lambda)).
\end{equation}
Put $n=\ell_{\mathcal G}(x)$, $r=s(x)$, $\nu_0=1$ and
$\nu_r=(q-1)q^{r-1}$ for $r\ge1$. An edge changes $n$ precisely
when it lies on $\mathcal G$, so the additional term acts at $r=0$.
For a complex sequence $\psi=(\psi_r)_{r\in\N_0}$, put
$h_r=\psi_r/\sqrt{\nu_r}$ for $r\in\N_0$. At a vertex of $\mathcal G$, the
longitudinal neighbors have coordinates $n-1,n+1$, and the other
$q-1$ neighbors have coordinates $(n,1)$. For $x=\gamma(n)$ with $n\in\Z$, applying $\mathcal L_{\mathcal G}$
to $x\mapsto e^{-\lambda\ell_{\mathcal G}(x)}h_{s(x)}$ gives
\[
 e^{-\lambda n}
 \bigl[D(q-1)(h_1-h_0)+a_0h_0+\chi(\lambda)h_0\bigr].
\]
For $r\ge1$, there is one neighbor at $r-1$ and $q$ at $r+1$,
all with the same $n$. After factoring out $e^{-\lambda n}$,
multiplication of the remaining transverse expression by $\sqrt{\nu_r}$
gives $(\mathcal J_{\mathcal G}\psi)_r$.
It follows from these neighbor counts that
\begin{equation}\label{eq:tan-coordinate-conjugacy}
 \mathcal L_{\mathcal G}
       \bigl[e^{-\lambda\ell_{\mathcal G}(\cdot)}h_{s(\cdot)}\bigr](x)
 =\frac{e^{-\lambda n}}{\sqrt{\nu_r}}
       (\mathcal J_{\parallel}(\lambda)\psi)_r
 \quad\text{ for }n=\ell_{\mathcal G}(x),\ r=s(x)
 \qquad x\in\Tq,\ \lambda\in\R.
\end{equation}
Both sides of \eqref{eq:tan-coordinate-conjugacy} involve only
finitely many entries of $\psi$ at each vertex. The identity therefore holds pointwise for every complex sequence
$\psi$ on $\N_0$. For $\lambda\in\C$, define
$\mathcal J_{\parallel}(\lambda)=\mathcal J_{\mathcal G}
+2D(\cosh\lambda-1)\langle\delta_0,\,\cdot\,\rangle\delta_0$.
The finite neighbor sums then give \eqref{eq:tan-coordinate-conjugacy}
for every $\lambda\in\C$ and every complex sequence $\psi$ on $\N_0$.

\begin{lemma}\label{lem:tangential-dispersion}
Assume $\mu\le0<\sigma_{\mathcal G}$.  Then
$\Lambda_{\parallel}$ is even, real analytic, and strictly convex on
$\R$.  Moreover,
\[
 \Lambda_{\parallel}(0)=\sigma_{\mathcal G}>0,
 \qquad
 \frac{\Lambda_{\parallel}(\lambda)}{|\lambda|}
 \to\infty
 \quad\text{ as }|\lambda|\to\infty.
\]
For every $\lambda\in\R$,
\begin{equation}
\label{eq:tangential-schur}
 \mathfrak{s}(\Lambda_{\parallel}(\lambda))
 =
 2D(\cosh\lambda-1),
\end{equation}
where $\mathfrak{s}$ is the scalar function in
Lemma~\ref{lem:geodesic-resolvent}.  Consequently there is a unique
$\lambda_{\parallel}>0$ minimizing
$\Lambda_{\parallel}(\lambda)/\lambda$ on $(0,\infty)$, and
\begin{equation}
\label{eq:tangential-speed}
 c_{\parallel}
 :=
 \frac{\Lambda_{\parallel}(\lambda_{\parallel})}
      {\lambda_{\parallel}}
 =
 \Lambda_{\parallel}'(\lambda_{\parallel})>0.
\end{equation}
For every $c\in(0,c_{\parallel})$ there is a unique
$\lambda_c\in(0,\lambda_{\parallel})$ such that
\[
 \Lambda_{\parallel}'(\lambda_c)=c,
 \qquad
 \Lambda_{\parallel}(\lambda_c)-c\lambda_c>0.
\]
\end{lemma}

\begin{proof}

Let $\mathcal B_1$ be the zero-exterior Jacobi operator on
$\ell^2(\{1,2,\ldots\})$ with diagonal $a_r-D(q+1)$ and
off-diagonal entries $D\sqrt q$, and put $\tau_0=D\sqrt{q-1}$.
Lemma~\ref{lem:geodesic-resolvent} gives
$\sup\spec(\mathcal B_1)=\mu$.
For $E>\mu$, $\lambda\in\R$, and $v\in\ell^2(\N_0;\C)$,
write $v=(v_0,v')$,
where $v'=(v_r)_{r\ge1}$. Thus $E-\mathcal B_1$ is bounded below
by $(E-\mu)I>0$. Applying \eqref{eq:geodesic-schur} to
the operator in \eqref{eq:tangential-fiber} gives
\begin{equation}\label{eq:tan-schur-form}
 \begin{aligned}
 \langle v,(E-\mathcal J_{\parallel}(\lambda))v\rangle
 &=
 \left\langle
 v'-\tau_0v_0(E-\mathcal B_1)^{-1}\delta_1,\,
 (E-\mathcal B_1)
 \bigl(v'-\tau_0v_0(E-\mathcal B_1)^{-1}\delta_1\bigr)
 \right\rangle\\
 &\quad+
 \bigl(\mathfrak{s}(E)-\chi(\lambda)\bigr)|v_0|^2.
 \end{aligned}
 \qquad \lambda\in\R.
\end{equation}
Since $\sigma_{\mathcal G}>\mu$,
Lemma~\ref{lem:geodesic-resolvent} gives
$\mathfrak{s}(\sigma_{\mathcal G})=0$.
The same lemma gives $\mathfrak{s}'(E)>0$ and
$\mathfrak{s}(E)\to+\infty$ as $E\to+\infty$.
Since $\chi(\lambda)\ge0$, there is a unique
$E\ge\sigma_{\mathcal G}$ such that $\mathfrak{s}(E)=\chi(\lambda)$.
By \eqref{eq:tan-schur-form}, we have
$E-\mathcal J_{\parallel}(\lambda)\ge0$.
For $v_0=1$ and $v'=\tau_0(E-\mathcal B_1)^{-1}\delta_1$,
both terms in \eqref{eq:tan-schur-form} vanish. Thus $E$ is an eigenvalue and
an upper bound for the spectrum, which proves
$E=\Lambda_{\parallel}(\lambda)$ and \eqref{eq:tangential-schur}.
For $E=\Lambda_{\parallel}(\lambda)$, the eigenvalue equations
at $r=0$ and $r\ge1$ give the eigenvector
\begin{equation}\label{eq:tan-positive-eigenvector}
 \psi^\lambda
 =\bigl(1,\tau_0(E-\mathcal B_1)^{-1}\delta_1\bigr)
 \in\ell^2(\N_0)
 \qquad \lambda\in\R,\ E=\Lambda_{\parallel}(\lambda).
\end{equation}
For $r\ge1$, the off-diagonal entries $D\sqrt q$ of $\mathcal B_1$
are positive. The matrix $\mathcal B_1+D(q+1)I$ has diagonal
entries $a_r>0$ and off-diagonal entries $D\sqrt q>0$; its
power of order $r-1$ has a positive $(r,1)$ entry. Thus
$(e^{t\mathcal B_1})(r,1)>0$ for $t>0$, and
\[
 ((E-\mathcal B_1)^{-1}\delta_1)_r
 =\int_0^\infty e^{-Et}(e^{t\mathcal B_1})(r,1)\,\dd t>0
 \qquad r\in\N.
\]
The integral converges because $E>\mu=\sup\spec(\mathcal B_1)$.
Together with $\psi_0^\lambda=1$ and $\tau_0>0$, this proves
$\psi_r^\lambda>0$ for every $r\ge0$. The Schur complement also shows that $E$ is an isolated
simple eigenvalue. Indeed, the tail resolvent exists throughout a
neighborhood of $E>\mu$, and the scalar factor
$\mathfrak{s}(\widetilde E)-\chi(\lambda)$ has a simple zero at
$\widetilde E=E$ because $\mathfrak{s}'(E)>0$. The tail block equation gives
$v'=\tau_0v_0(E-\mathcal B_1)^{-1}\delta_1$, so the eigenspace is
one-dimensional. Self-adjointness excludes generalized eigenvectors.

For $E>\mu$, the resolvent $(E-\mathcal B_1)^{-1}$ depends
analytically on $E$. By $\mathfrak{s}(E)=E-a_0+D(q-1)
-\tau_0^2\langle\delta_1,(E-\mathcal B_1)^{-1}\delta_1\rangle$
in Lemma~\ref{lem:geodesic-resolvent},
\begin{equation}\label{eq:tan-schur-derivatives}
 \begin{aligned}
 \mathfrak{s}'(E)
 &=1+\tau_0^2\langle\delta_1,
             (E-\mathcal B_1)^{-2}\delta_1\rangle>0,\\
 \mathfrak{s}''(E)
 &=-2\tau_0^2\langle\delta_1,
             (E-\mathcal B_1)^{-3}\delta_1\rangle<0.
 \end{aligned}
\end{equation}
The strict inequalities follow from $E-\mathcal B_1>0$ and
$\delta_1\ne0$. In view of \eqref{eq:tan-schur-derivatives}, the
analytic implicit function theorem applies to
\eqref{eq:tangential-schur}. Differentiating $\mathfrak{s}(\Lambda_{\parallel}(\lambda))=\chi(\lambda)$
twice yields
\begin{equation}\label{eq:tan-fiber-curvature}
 \Lambda_{\parallel}''(\lambda)
 =\frac{\chi''(\lambda)
       -\mathfrak{s}''(\Lambda_{\parallel}(\lambda))
        (\Lambda_{\parallel}'(\lambda))^2}
       {\mathfrak{s}'(\Lambda_{\parallel}(\lambda))}>0,
 \qquad \chi''(\lambda)=2D\cosh\lambda>0
 \qquad \lambda\in\R.
\end{equation}
By \eqref{eq:tan-fiber-curvature}, $\Lambda_{\parallel}$ is analytic and strictly convex. Its
evenness follows from
the uniqueness of $E>\mu$ satisfying $\mathfrak{s}(E)=\chi(\lambda)$ and
the evenness of $\chi$.

Testing $\mathcal J_{\parallel}(\lambda)$ on $\delta_0$ gives
\[
 \Lambda_{\parallel}(\lambda)
 \ge
 a_0-D(q-1)+2D(\cosh\lambda-1)
 \qquad \lambda\in\R,
\]
so
$\Lambda_{\parallel}(\lambda)/|\lambda|\to\infty$ as $|\lambda|\to\infty$.
Since $\Lambda_{\parallel}(0)=\sigma_{\mathcal G}>0$, the quotient
$\Lambda_{\parallel}(\lambda)/\lambda$ tends to infinity also as
$\lambda\downarrow0$. Together with its divergence as
$\lambda\to\infty$, this gives an interior minimizer.  Put
\[
 G_{\parallel}(\lambda)
 =
 \lambda\Lambda_{\parallel}'(\lambda)
 -
 \Lambda_{\parallel}(\lambda)
 \qquad \lambda\in\R.
\]
Then
\[
 G_{\parallel}(0)
 =
 -\sigma_{\mathcal G}<0,
 \qquad
 G_{\parallel}'(\lambda)
 =
 \lambda\Lambda_{\parallel}''(\lambda)>0
 \quad\text{ for } \lambda>0.
\]
At every interior minimizer of
$\Lambda_{\parallel}(\lambda)/\lambda$, differentiation gives
$G_{\parallel}(\lambda)=0$. Such a minimizer exists by the two
limits already proved. Since $G_{\parallel}'(\lambda)>0$ for
$\lambda>0$, the zero is unique; denote it by
$\lambda_{\parallel}$. This proves \eqref{eq:tangential-speed}.
Evenness gives $\Lambda_{\parallel}'(0)=0$, while strict convexity
makes $\Lambda_{\parallel}'$ strictly increasing. It follows from
\eqref{eq:tangential-speed} that $\Lambda_{\parallel}'$ maps
$(0,\lambda_{\parallel})$ bijectively onto $(0,c_{\parallel})$.  At the point
$\lambda_c$ with derivative $c$,
\[
 \Lambda_{\parallel}(\lambda_c)-c\lambda_c
 =
 -G_{\parallel}(\lambda_c)>0
 \qquad c\in(0,c_{\parallel}).
\qedhere\]

\end{proof}

We construct a compactly supported subsolution moving along $\mathcal G$.
For $n\in\Z$ and $r\in\N_0$, put
$\mathcal C_{n,r}=\{x\in\Tq:\ell_{\mathcal G}(x)=n,\ s(x)=r\}$.
By the neighbor counts for the tree, $|\mathcal C_{n,r}|=\nu_r$,
independently of $n$. Translations and
reflections of $\mathcal G$ extend to tree automorphisms preserving
$s(x)$; these are the translations and reflections of data used below.

\begin{lemma}\label{lem:tangential-bump}

Assume $\mu\le0<\sigma_{\mathcal G}$.  For every
$c\in(0,c_{\parallel})$ there are $R\in\N$, $L>1$, and
$\delta>0$, together with a nonnegative entire subsolution
$v^+(t,x)$ of the differential equation in \eqref{eq:tree-cauchy}
with $f(x,u)=f_{s(x)}(u)$, such that
\[
 0\le v^+(t,x)\le1\quad(t\in\R,\ x\in\Tq),
\]
\[
 v^+(t,x)=0
 \quad\text{ if }s(x)>R\text{ or }|\ell_{\mathcal G}(x)-ct|\ge L
 \qquad t\in\R,\ x\in\Tq,
\]
and
\[
 v^+(t,x)\ge\delta
 \quad\text{whenever}\quad
 s(x)\le R,\ 
 |\ell_{\mathcal G}(x)-ct|\le1
 \qquad t\in\R,\ x\in\Tq.
\]
Reflection in the geodesic parameter gives a left-moving subsolution
$v^-$ satisfying the same bounds with
$|\ell_{\mathcal G}(x)-ct|$ replaced by $|\ell_{\mathcal G}(x)+ct|$.

\end{lemma}

\begin{proof}

Fix $c\in(0,c_{\parallel})$ and let
$\lambda_c\in(0,\lambda_{\parallel})$ satisfy
$\Lambda_{\parallel}'(\lambda_c)=c$ and
$\Lambda_{\parallel}(\lambda_c)-c\lambda_c>0$, as given by
Lemma~\ref{lem:tangential-dispersion}.  For $R\in\N$, let
$\mathcal J_{\parallel,R}(\lambda)$ be the zero-exterior restriction of
$\mathcal J_{\parallel}(\lambda)$ to
$\{0,\ldots,R\}$, retaining the full diffusion diagonal, and let
$\Lambda_R(\lambda)$ be its largest eigenvalue on $\C^{R+1}$
with the Euclidean scalar product. The adjacent off-diagonal entries
of $\mathcal J_{\parallel,R}(\lambda)$ are positive. After adding
$D(q+1)I$, its diagonal entries are $a_0+2D+\chi(\lambda)>0$
at $r=0$ and $a_r>0$ at $1\le r\le R$. The shifted matrix is
therefore nonnegative and irreducible, so $\Lambda_R(\lambda)$ is simple and has a strictly
positive eigenvector. The characteristic polynomial of
$\mathcal J_{\parallel,R}(\lambda)$ is analytic in $\lambda$, with a
simple root at $\Lambda_R(\lambda)$. The analytic
implicit function theorem shows that $\Lambda_R$ is real analytic.
For each $R\in\N$, choose the real analytic, positive unit
eigenvector $e_R=e_R(\lambda)\in\R^{R+1}$ associated with
$\Lambda_R(\lambda)$, and let $Q_R=Q_R(\lambda)$ be the orthogonal
projection onto $e_R(\lambda)^\perp$ in $\C^{R+1}$.
For $v\in\C^{R+1}$, put $P_0v=v_0\delta_0$.
The analytic choice of $e_R$ follows from simplicity and its positive
unit normalization. Differentiating
$\|e_R\|_2^2=1$ gives $\langle e_R,e_R'\rangle=0$.
The restriction of
$\Lambda_R(\lambda)-\mathcal J_{\parallel,R}(\lambda)$ to
$e_R(\lambda)^\perp$ is invertible, since $\Lambda_R(\lambda)$ is
the simple largest eigenvalue. Denote its inverse by
$(\Lambda_R(\lambda)-\mathcal J_{\parallel,R}(\lambda))^{-1}_{e_R^\perp}$.
The orthogonal component of the differentiated eigenvalue equation is
\[
 e_R'=\chi'(\lambda)
 (\Lambda_R-\mathcal J_{\parallel,R})^{-1}_{e_R^\perp}
 Q_RP_0e_R
 \qquad \lambda\in\R,\ R\in\N.
\]
Differentiating
$\Lambda_R'=\chi'(\lambda)\langle e_R,P_0e_R\rangle$
then gives the following formula (the general perturbation method is
in \cite[Chapter~II]{Kato}):
\begin{equation}\label{eq:tan-finite-curvature}
 \begin{aligned}
 \Lambda_R''(\lambda)
 &=\chi''(\lambda)|(e_R)_0|^2\\
 &\quad+2(\chi'(\lambda))^2
 \left\langle
 Q_R P_0 e_R,
 (\Lambda_R(\lambda)-\mathcal J_{\parallel,R}(\lambda))^{-1}_{e_R^\perp}
 Q_R P_0 e_R
 \right\rangle>0,
 \end{aligned}
 \qquad \lambda\in\R,\ R\in\N.
\end{equation}
The inverse on $e_R(\lambda)^\perp$ has positive quadratic form on
nonzero vectors in that space because $\Lambda_R(\lambda)$ is the
simple largest eigenvalue. The term
$\chi''(\lambda)|(e_R)_0|^2$ is strictly positive, since $\chi''(\lambda)>0$ and $(e_R)_0>0$.

Extending a vector in $\C^{R+1}$ by zero preserves the quadratic
form of $\mathcal J_{\parallel}(\lambda)$. These finite-support
spaces increase with $R$ and have dense union in $\ell^2(\N_0)$.
Boundedness of the operator and the Rayleigh principle therefore give
\begin{equation}\label{eq:tan-finite-spectral-limit}
 \Lambda_R(\lambda)\uparrow
 \Lambda_{\parallel}(\lambda)\quad\text{ as } R\to\infty
\end{equation}
for every real $\lambda$. For each compact interval $J\subset\R$, the functions
$\Lambda_R$ and $\Lambda_{\parallel}$ are continuous on $J$,
and $\Lambda_R$ increases pointwise to $\Lambda_{\parallel}$
by \eqref{eq:tan-finite-spectral-limit}. Dini's theorem therefore
gives $\sup_{\lambda\in J}|\Lambda_R(\lambda)-\Lambda_{\parallel}(\lambda)|\to0$
as $R\to\infty$. For $h>0$, convexity gives
\[
 \frac{\Lambda_R(\lambda)-\Lambda_R(\lambda-h)}h
 \le\Lambda_R'(\lambda)
 \le\frac{\Lambda_R(\lambda+h)-\Lambda_R(\lambda)}h
 \qquad \lambda\in\R,\ R\in\N.
\]
Letting first $R\to\infty$ and then $h\downarrow0$, and using the
differentiability of $\Lambda_{\parallel}$, yields
$\Lambda_R'(\lambda)\to\Lambda_{\parallel}'(\lambda)$ at each fixed
$\lambda$. Choose $0<\lambda_{\mathrm{lo}}<\lambda_c<\lambda_{\mathrm{hi}}$ so close to
$\lambda_c$ that
\[
 \Lambda_{\parallel}'(\lambda_{\mathrm{lo}})<c<
 \Lambda_{\parallel}'(\lambda_{\mathrm{hi}})
\]
and
\[
 \Lambda_{\parallel}(\lambda)-c\lambda>0
 \quad\text{ for } \lambda_{\mathrm{lo}}\le\lambda\le\lambda_{\mathrm{hi}}.
\]
Set
\[
 d_c=\min_{\lambda\in[\lambda_{\mathrm{lo}},\lambda_{\mathrm{hi}}]}
       (\Lambda_{\parallel}(\lambda)-c\lambda)>0.
\]
Choose $R\in\N$ so large that
$\sup_{\lambda\in[\lambda_{\mathrm{lo}},\lambda_{\mathrm{hi}}]}
|\Lambda_R(\lambda)-\Lambda_{\parallel}(\lambda)|<d_c/2$ and
$\Lambda_R'(\lambda_{\mathrm{lo}})<c<
\Lambda_R'(\lambda_{\mathrm{hi}})$. Then
$\Lambda_R(\lambda)-c\lambda\ge d_c/2$ throughout the interval.
Since $\Lambda_R''>0$, the equation $\Lambda_R'(\lambda)=c$ has a unique
$\lambda_R\in(\lambda_{\mathrm{lo}},\lambda_{\mathrm{hi}})$ with
\[
 \Lambda_R'(\lambda_R)=c,
 \qquad
 \Lambda_R(\lambda_R)-c\lambda_R>0.
\]

Let $\mathcal J_{\mathcal G,R}$ be the zero-exterior restriction
of $\mathcal J_{\mathcal G}$ to $\{0,\ldots,R\}$. For complex
$z$ near $\lambda_R$, write
\[
 \mathcal J_{\parallel,R}(z)
 =
 \mathcal J_{\mathcal G,R}
 +
 2D(\cosh z-1)
 \langle\delta_0,\,\cdot\,\rangle\delta_0.
\]
Its entries are holomorphic in $z$, and its principal eigenvalue
at the real point $\lambda_R$ is simple. The characteristic
polynomial therefore has nonzero derivative in the eigenvalue
variable there. The holomorphic implicit function theorem gives a
local eigenvalue branch $\Lambda_R(z)$. An invertible minor of the
eigenvector equations gives a holomorphic eigenvector branch
$\psi^R(z)$; this is the simple-eigenvalue case of analytic
perturbation theory \cite[Chapter~II]{Kato}.
Choose $\psi^R(\lambda_R)$ positive. The zeroth entry of $\psi^R(z)$
stays nonzero near $\lambda_R$, so normalize $\psi^R$ by $\psi^R_0(z)=1$.
The identity
$\mathcal J_{\parallel,R}(\overline z)=
\overline{\mathcal J_{\parallel,R}(z)}$ and uniqueness of the simple
eigenvalue branch give
$\Lambda_R(\overline z)=\overline{\Lambda_R(z)}$.
Consequently the imaginary part of
$\Lambda_R(\alpha+\mathrm i\eta)$ has Taylor expansion
$\eta\Lambda_R'(\alpha)+O(\eta^3)$ near $(\lambda_R,0)$.
The quotient $\im\Lambda_R(\alpha+\mathrm i\eta)/\eta$
therefore extends smoothly to $\eta=0$.
For real $(\alpha,\eta)$ near $(\lambda_R,0)$, define
\[
 F(\alpha,\eta)
 =
 \begin{cases}
 \displaystyle
 \frac{\im\Lambda_R(\alpha+\mathrm i\eta)}{\eta},
 &\eta\ne0,\\[2mm]
 \Lambda_R'(\alpha),&\eta=0.
 \end{cases}
\]
Since $F(\lambda_R,0)=c$ and
$\partial_\alpha F(\lambda_R,0)=\Lambda_R''(\lambda_R)>0$,
the implicit function theorem gives
$\alpha(\eta)$ with $\alpha(\eta)\to\lambda_R$ as $\eta\downarrow0$ and
\[
 F(\alpha(\eta),\eta)=c
 \qquad \text{for all sufficiently small }\eta>0.
\]
Choose $\eta>0$ sufficiently small and set
\[
 z=\alpha(\eta)+\mathrm i\eta,
 \qquad
 \rho=\Lambda_R(z)-cz.
\]
By $F(\alpha(\eta),\eta)=c$, we have $\im\rho=0$. Since
$\Lambda_R(\lambda_R)-c\lambda_R>0$, continuity makes $\rho>0$
for all sufficiently small $\eta>0$.

Set
\[
 h_r^R(z)
 =
 \frac{\psi_r^R(z)}{\sqrt{\nu_r}},
 \quad\text{ for }0\le r\le R.
\]
Put $\xi=\ell_{\mathcal G}(x)-ct$. For $s(x)\le R$, define
\[
 Z(t,x)
 =
 e^{-z\xi}h_{s(x)}^R(z)
 \qquad t\in\R,\ x\in\Tq,\ s(x)\le R.
\]
Set $Z(t,x)=0$ for $t\in\R$ and $s(x)>R$.
By $\mathcal J_{\parallel,R}(z)\psi^R(z)=\Lambda_R(z)\psi^R(z)$
and \eqref{eq:tan-coordinate-conjugacy}, we have
\[
 Z_t
 =
 \mathcal L_{\mathcal G}Z-\rho Z
 \qquad t\in\R,\ x\in\Tq,\ s(x)\le R.
\]
The exterior value at $s(x)=R+1$ is zero. Define
$V(t,x)=\re Z(t,x)$ for $t\in\R$ and $x\in\Tq$.
After shrinking $\eta$, all finitely many arguments of
$h_r^R(z)$ have absolute value less than $\pi/8$, and we may also
assume $0<\eta<\pi/4$. Hence
\[
 V(t,x)<0
 \quad\text{whenever}\quad
 s(x)\le R,\qquad
 \frac{3\pi}{4\eta}
 \le|\ell_{\mathcal G}(x)-ct|
 \le\frac{\pi}{\eta}
 \qquad t\in\R,\ x\in\Tq.
\]
Take $L=\pi/\eta$. For $\varepsilon>0$, to be chosen below, define
\[
 v^+(t,x)
 =
 \varepsilon(V(t,x))_+
 \ind_{\{s(x)\le R,\,
 |\ell_{\mathcal G}(x)-ct|<L\}}
 \qquad t\in\R,\ x\in\Tq.
\]
The intervals $3\pi/(4\eta)\le|\xi|\le L$, where
$V(t,x)<0$ for $s(x)\le R$, have width $\pi/(4\eta)>1$.
For each $0\le r\le R$, the function
$\xi\mapsto(\re(e^{-z\xi}h_r^R(z)))_+
\ind_{\{|\xi|<L\}}$ vanishes in these boundary intervals and is
Lipschitz on $\R$. The maximum of these Lipschitz constants over $r\in\{0,\ldots,R\}$
is finite. Since $\xi=\ell_{\mathcal G}(x)-ct$, it follows that
$t\mapsto v^+(t,\cdot)$ is locally Lipschitz in
$\ell^\infty(\Tq)$, as required by Lemma~\ref{lem:comparison}.
Where $v^+>0$, the longitudinal coordinate lies more than one step
from the cutoff boundaries. Every neighbor $y\sim x$ with $s(y)\le R$ therefore satisfies
$|\ell_{\mathcal G}(y)-ct|<L$, since an edge changes
$\ell_{\mathcal G}$ by at most one. Replacing $V(t,y)$ by
$(V(t,y))_+$ at these neighbors increases the off-diagonal diffusion
terms. At a neighbor $y\sim x$ with $s(y)>R$, both $Z(t,y)$ and
$v^+(t,y)$ are zero.  At a site where $v^+(t,x)=0$, the derivative $\partial_t v^+(t,x)$
vanishes for almost every such $t$, since $t\mapsto v^+(t,x)$ is
nonnegative and locally Lipschitz. All off-diagonal diffusion terms
there are nonnegative.
Therefore
\begin{equation}\label{eq:tan-moving-bump-residual}
 (v^+)_t
 \le
 D\Delta_{\Tq}v^+
 +a_{s(x)}v^+
 -\rho v^+
 \qquad x\in\Tq,\ \text{a.e. }t\in\R.
\end{equation}
The countability of $\Tq$ gives a
common full-measure set of such times. For the fixed $R,z,L$ and any $\varepsilon>0$,
$v^+/\varepsilon$ is bounded uniformly over $t\in\R$ and
$x\in\Tq$, since $|\xi|\le L$ on its support. Choose
$\varepsilon>0$ so small that $0\le v^+\le1$ and
$C_2v^+\le\rho$ everywhere. By the bound
$f_r(v)\ge a_rv-C_2v^2$ in \eqref{eq:taylor}, \eqref{eq:tan-moving-bump-residual}
then gives $(v^+)_t\le D\Delta_{\Tq}v^++f_{s(x)}(v^+)$.

Since $h_r^R(\lambda_R)>0$ for $0\le r\le R$, choose $\eta$
at the preceding stage small enough that
$m_0:=\min_{\substack{0\le r\le R\\|\xi|\le1}}
\re(e^{-z\xi}h_r^R(z))>0$. After fixing this $\eta$ and then
$\varepsilon$ as above, take $\delta=\varepsilon m_0$. Choose a tree automorphism
$\mathcal R$ with $\ell_{\mathcal G}(\mathcal R x)=-\ell_{\mathcal G}(x)$
and $s(\mathcal R x)=s(x)$, and put $v^-(t,x)=v^+(t,\mathcal R x)$.
Since $\mathcal R$ preserves adjacency and $s$, the residual of $v^-$ satisfies
\[
\begin{gathered}
 (v^-)_t-D\Delta_{\Tq}v^--f_{s(x)}(v^-)
 =\bigl((v^+)_t-D\Delta_{\Tq}v^+\bigr)(t,\mathcal R x)
   -f_{s(\mathcal R x)}(v^+(t,\mathcal R x))\le0
 \\\text{ for }x\in\Tq,\ \text{a.e. }t\in\R.
\end{gathered}
\]
The support and lower bound for $v^-$ follow by
$|\ell_{\mathcal G}(\mathcal R x)-ct|=|\ell_{\mathcal G}(x)+ct|$.\qedhere

\end{proof}

\begin{proof}[Proof of Theorem~\ref{thm:tangential-spreading}]

We first prove the outer estimate.  Fix $c>c_{\parallel}$.
By the definition of $c_{\parallel}$ there is $\lambda>0$ such that
\[
 \frac{\Lambda_{\parallel}(\lambda)}{\lambda}<c.
\]
Let $\psi^\lambda\in\ell^2(\N_0)$ be the positive eigenvector
constructed in the proof of Lemma~\ref{lem:tangential-dispersion},
normalized by $\psi^\lambda_0=1$, and set
\[
 h_r^\lambda
 =
 \frac{\psi_r^\lambda}{\sqrt{\nu_r}}
 \qquad r\in\N_0.
\]
Since $\nu_r\ge1$ for every $r\ge0$ and
$\psi^\lambda\in\ell^2(\N_0)$,
$\sup_r h_r^\lambda\le\|\psi^\lambda\|_2<\infty$. For $C>0$,
to be chosen below, $t\ge0$ and $x\in\Tq$, define
\[
 E^+(t,x)
 =
 C
 e^{-\lambda\ell_{\mathcal G}(x)
   +\Lambda_{\parallel}(\lambda)t}
 h_{s(x)}^\lambda
 \qquad t\ge0,\ x\in\Tq.
\]
The function $E^+$ solves $(E^+)_t=\mathcal L_{\mathcal G}E^+$ pointwise.
Because $u_0$ has finite support, choose $C>0$ with
\[
 C\ge\max_{x\in\supp u_0}
 \frac{e^{\lambda\ell_{\mathcal G}(x)}u_0(x)}{h_{s(x)}^\lambda}.
\]
This choice gives
\[
 u_0\le\min\{1,E^+(0,\cdot)\}.
\]
By \eqref{eq:hyp}, $\min\{1,E^+\}$ is a bounded
supersolution. Where $E^+(t,x)<1$, truncation decreases the
neighboring values, while $f_{s(x)}(E^+)\le a_{s(x)}E^+$. Where
$E^+(t,x)\ge1$, the truncated function equals one, its diffusion
term is nonpositive, and $f_{s(x)}(1)=0$. These
inequalities hold almost everywhere in time. Since
$\Lambda_{\parallel}(\lambda)>0$, $\min\{1,E^+\}$ is
Lipschitz in time in $\ell^\infty$, with Lipschitz constant at most
$\Lambda_{\parallel}(\lambda)$. It is therefore an admissible
supersolution for Lemma~\ref{lem:comparison}. Thus, on
$\ell_{\mathcal G}(x)\ge ct$,
\[
 u(t,x)
 \le
 C\sup_{r\in\N_0}h_r^\lambda\,
 e^{-(\lambda c-\Lambda_{\parallel}(\lambda))t}
 \qquad t\ge0,\ x\in\Tq,\ \ell_{\mathcal G}(x)\ge ct.
\]
Since $\lambda c-\Lambda_{\parallel}(\lambda)>0$, this bound tends
to zero as $t\to\infty$. For the other side, put
\[
 E^-(t,x)=C e^{\lambda\ell_{\mathcal G}(x)
                    +\Lambda_{\parallel}(\lambda)t}h_{s(x)}^\lambda
 \qquad t\ge0,\ x\in\Tq.
\]
The two longitudinal neighbors on $\mathcal G$ contribute
$D(e^\lambda+e^{-\lambda}-2)$, so the eigenvalue equation also gives
$(E^-)_t=\mathcal L_{\mathcal G}E^-$. Increase $C$ to obtain
$u_0\le\min\{1,E^-(0,\cdot)\}$. Where $E^-<1$, truncation lowers the neighboring values and
$f_{s(x)}(E^-)\le a_{s(x)}E^-$. Where $E^-\ge1$, the truncated
function has nonpositive diffusion and zero reaction. Moreover,
$\min\{1,E^-\}$ is Lipschitz in time in $\ell^\infty$, with
constant at most $\Lambda_{\parallel}(\lambda)$. It is therefore
a bounded supersolution, and Lemma~\ref{lem:comparison} yields
\[
 \sup_{\substack{x\in\Tq\\\ell_{\mathcal G}(x)\le-ct}}u(t,x)
 \le C\sup_{r\in\N_0}h_r^\lambda
       e^{-(\lambda c-\Lambda_{\parallel}(\lambda))t}\to0
 \quad\text{ as }t\to\infty.
\]
Combining the estimates on $\ell_{\mathcal G}(x)\le-ct$ and
$\ell_{\mathcal G}(x)\ge ct$ proves \eqref{eq:tangential-outer}.

For the lower bound along the geodesic, fix $0<c<c_2<c_{\parallel}$
and take the two subsolutions $v^+,v^-$ from
Lemma~\ref{lem:tangential-bump} with speed $c_2$.
The union of the supports of $v^+(0,\cdot)$ and $v^-(0,\cdot)$ is
finite. Since $p_*>0$, decrease the common multiplicative factor in the
construction of Lemma~\ref{lem:tangential-bump} until
$v^\pm(0,x)<p_*(x)$ throughout this union. The construction supplies $\rho>0$ such that
$(v^\pm)_t\le\mathcal L_{\mathcal G}v^\pm-\rho v^\pm$. For every $\vartheta\in(0,1]$,
\eqref{eq:taylor} therefore gives
\[
 D\Delta_{\Tq}(\vartheta v^\pm)
 +f_{s(x)}(\vartheta v^\pm)-\vartheta\partial_t v^\pm
 \ge\vartheta v^\pm(\rho-C_2\vartheta v^\pm)\ge0
 \qquad x\in\Tq,\ \text{a.e. }t\in\R,
\]
because the construction gives $C_2v^\pm\le\rho$.
For each vertex in the finite union of the two initial supports,
$\liminf_{t\to\infty}u(t,x)\ge p_*(x)>v^\pm(0,x)$ by
\eqref{eq:geodesic-minimal-persistence}. Taking the largest of the
finitely many resulting times gives $t_0$ such that
\[
 u(t,\cdot)\ge v^\pm(0,\cdot)
 \quad\text{ for } t\ge t_0.
\]
For the reduced amplitudes with $v^\pm(0,x)<p_*(x)$ on their
initial supports, let $\delta>0$ be the common lower bound from
Lemma~\ref{lem:tangential-bump} on
$\{x:s(x)\le R,\ |\ell_{\mathcal G}(x)\mp c_2t|\le1\}$ for
$v^\pm$, respectively. Choose $0<\varepsilon_0\le\delta$ and put
$\zeta=\varepsilon_0\ind_{\{\gamma(0)\}}$.
For $t\ge0$ and $m\in\Z$ satisfying $m=c_2t$ for $v^+$ or
$m=-c_2t$ for $v^-$, the bound on
$|\ell_{\mathcal G}(x)\mp c_2t|\le1$ gives
$v^\pm(t,\cdot)\ge\zeta_m$, where
$\zeta_m=\varepsilon_0\ind_{\{\gamma(m)\}}$. In particular,
$v^\pm(0,\cdot)\ge\zeta$.

Let $w$ solve \eqref{eq:tree-cauchy} with
$f(x,u)=f_{s(x)}(u)$ and $w(0,\cdot)=\zeta$.  By
Theorem~\ref{thm:geodesic-threshold},
\[
 \liminf_{T\to\infty}w(T,\gamma(0))
 \ge p_*(\gamma(0)).
\]
Fix $0<\theta<p_*(\gamma(0))$ and choose $T>0$ such that
\[
 w(T,\gamma(0))>\theta.
\]
Let $m\in\Z$ satisfy $|m|\le ct$.  If $m\ge0$, start comparison with
$v^+$ at time
\[
 s_m=t-T-\frac{|m|}{c_2};
\]
if $m<0$, use $v^-$ with the same $s_m$.
Uniformly for $|m|\le ct$,
\[
 s_m
 \ge
 t\left(1-\frac c{c_2}\right)-T\to\infty
 \quad\text{ as } t\to\infty
 \qquad \text{uniformly for }m\in\Z,\ |m|\le ct.
\]
Thus $s_m\ge t_0$ for all sufficiently large $t$. At $s_m$ we have
$u(s_m,\cdot)\ge v^\pm(0,\cdot)$, choosing $v^+$ for $m\ge0$ and
$v^-$ for $m<0$. Comparison for the elapsed time $|m|/c_2$ gives
\[
 u(t-T,\cdot)
 \ge v^\pm(|m|/c_2,\cdot)\ge\zeta_m,
\]
since the chosen subsolution is centered at $m$ at that elapsed time.
Translation invariance along $\mathcal G$ and another comparison over
the last time interval of length $T$ give
\[
 u(t,\gamma(m))
 \ge
 w(T,\gamma(0))
 >
 \theta.
\]
Hence
\[
 \{m\in\Z:|m|\le ct\}
 \subseteq
 \{m\in\Z:u(t,\gamma(m))\ge\theta\}
\]
for all sufficiently large $t$. In particular,
$\{m\in\Z:u(t,\gamma(m))\ge\theta\}$ is nonempty. For any fixed $c_3>c_{\parallel}$, the outer estimate
makes $u(t,\gamma(m))<\theta$ when $|m|\ge c_3t$ and $t$ is
sufficiently large, so the level set is bounded and $L_\theta(t)$
is defined. The two inclusions give
$\lfloor ct\rfloor\le L_\theta(t)\le c_3t$ for large $t$.
Letting $c\uparrow c_{\parallel}$ and
$c_3\downarrow c_{\parallel}$ proves
\eqref{eq:tangential-level-speed}.

Assume now \eqref{eq:strict-quotient} and retain the speeds
$0<c<c_2<c_{\parallel}$ and the datum
$\zeta=\varepsilon_0\ind_{\{\gamma(0)\}}$ chosen above.
Theorem~\ref{thm:geodesic-critical-convergence}, applied to the
solution $w$ from $\zeta$, gives
\[
 \sup_{\{x\in\Tq:\ell_{\mathcal G}(x)=0\}}
 |w(T,x)-p_*(x)|\to0
 \quad\text{ as } T\to\infty.
\]
Given $\varepsilon>0$, choose $T$ so that
$|w(T,x)-p_*(x)|<\varepsilon$ on $\{x:\ell_{\mathcal G}(x)=0\}$.
Choose a tree automorphism $\mathcal T$ such that
$\ell_{\mathcal G}(\mathcal Tx)=\ell_{\mathcal G}(x)+1$ and
$s(\mathcal Tx)=s(x)$, and put $\mathcal T_m=\mathcal T^m$ for
$m\in\Z$. Thus $\mathcal T_{-m}=\mathcal T_m^{-1}$.
For $m\in\Z$, $\tau\ge0$ and $x\in\Tq$, put
$w_m(\tau,x)=w(\tau,\mathcal T_{-m}x)$. Then $w_m$ solves
\eqref{eq:tree-cauchy} with initial datum $\zeta_m$, and
$p_*(\mathcal T_{-m}x)=p_*(x)$.
For $|m|\le ct$, let $s_m=t-T-|m|/c_2$. The bound
$s_m\ge t(1-c/c_2)-T$ makes $s_m\ge t_0$ uniformly for large $t$.
Starting $v^+$ or $v^-$ at $s_m$ gives
$u(t-T,\cdot)\ge\zeta_m$; comparison over the last interval of length
$T$ therefore yields, on $\ell_{\mathcal G}(x)=m$,
\[
 u(t,x)\ge w_m(T,x)
 =w(T,\mathcal T_{-m}x)\ge p_*(x)-\varepsilon.
\]
The choice of $T$ is independent of $m$, so this lower bound is
uniform on $|\ell_{\mathcal G}(x)|\le ct$ for all sufficiently large $t$.

For the upper bound on $u-p_*$, put
$h_r=q^{-r/2}\phi_r(\lambda_0)$ for $r\in\N_0$.
Choose
$C\ge\max\{1/h_0,\max_{x\in\supp u_0}u_0(x)/h_{s(x)}\}$.
Then
\[
 u_0(x)\le\overline U^C_{s(x)}
 :=
 \min\{1,Ch_{s(x)}\}
 \qquad x\in\Tq.
\]
Let $U^C(t,r)$ solve the interior equation in \eqref{eq:radial} and
the boundary equation \eqref{eq:geodesic-boundary}, with
$U^C(0,r)=\overline U_r^C$. Its lifted initial datum is nonzero
and bounded by $Ch_{s(x)}$.
Thus Theorem~\ref{thm:geodesic-critical-convergence} applies to
$x\mapsto U^C(t,s(x))$ with $M=0$. {Every transverse index $r\ge0$
occurs in $\{x:\ell_{\mathcal G}(x)=0\}$. Both
$x\mapsto U^C(t,s(x))$ and $p_*$ depend only on $s(x)$;
recall that $p_{*,r}=p_*(x)$ whenever $s(x)=r$.
Consequently \eqref{eq:critical-geodesic-convergence} with $M=0$ gives}
\[
 \sup_{r\in\N_0}|U^C(t,r)-p_{*,r}|\to0
 \quad\text{ as } t\to\infty.
\]
Tree comparison therefore yields
\[
 u(t,x)\le p_*(x)+\varepsilon
\]
for all $x\in\Tq$ and all sufficiently large $t$, with the time
threshold independent of $x$. On
$|\ell_{\mathcal G}(x)|\le ct$, the lower comparison gives
$u(t,x)-p_*(x)\ge-\varepsilon$, while this upper comparison gives
$u(t,x)-p_*(x)\le\varepsilon$. Hence
\[
 \limsup_{t\to\infty}
 \sup_{\substack{x\in\Tq\\|\ell_{\mathcal G}(x)|\le ct}}|u(t,x)-p_*(x)|
 \le\varepsilon.
\]
Letting $\varepsilon\downarrow0$ proves
\eqref{eq:tangential-interior-profile}.\qedhere

\end{proof}

\Needspace{8\baselineskip}
\subsection{Tangential Bramson correction}\label{sec:tangential-bramson}

Throughout this subsection, assume \eqref{eq:hyp},
$f(x,u)=f_{s(x)}(u)$, and $\mu\le0<\sigma_{\mathcal G}$.
By Theorem~\ref{thm:tangential-spreading}, solutions from nonzero
finitely supported data $u_0:\Tq\to[0,1]$ spread along $\mathcal G$ at speed
$c_{\parallel}$. We prove the displacement
$-3\log t/(2\lambda_{\parallel})$ by estimating a Markov chain
killed at a moving longitudinal boundary.

We combine bounds of order $t^{-1/2}$ for survival near the initial
and terminal times with a transition bound of order $t^{-1/2}$ for
the middle part of the path. Their product gives the factor $t^{-3/2}$.
Put $\kappa_{\parallel}=3/(2\lambda_{\parallel})$. By
Lemma~\ref{lem:tangential-dispersion},
\begin{equation}\label{eq:tan-critical-identities}
 \Lambda_{\parallel}(\lambda_{\parallel})
 =c_{\parallel}\lambda_{\parallel},
 \qquad
 \Lambda_{\parallel}'(\lambda_{\parallel})=c_{\parallel},
 \qquad
 \Lambda_{\parallel}''(\lambda_{\parallel})>0.
\end{equation}
Fix
$a_*=1+\kappa_{\parallel}/c_{\parallel}
+2\lambda_{\parallel}\kappa_{\parallel}/
(\Lambda_{\parallel}(\lambda_{\parallel})-\mu)$.
Then $a_*>0$ and
\[
 c_{\parallel}-\frac{\kappa_{\parallel}}{1+a_*}>0
\]
and
\[
 \frac{\lambda_{\parallel}\kappa_{\parallel}}{1+a_*}
 <
 \frac12\bigl(\Lambda_{\parallel}(\lambda_{\parallel})-\mu\bigr).
\]
The inequality
$\lambda_{\parallel}\kappa_{\parallel}/(1+a_*)
<(\Lambda_{\parallel}(\lambda_{\parallel})-\mu)/2$
controls the first hitting time of transverse state zero in
Lemma~\ref{lem:tangential-linear-seed}.
For $b\in\R$ and $t\ge0$, put
\[
 m_b(t)
 =
 b+c_{\parallel}t
 -\kappa_{\parallel}\log\frac{1+a_*+t}{1+a_*}.
\]
Thus $m_b'(t)>0$ for every $t\ge0$, and
\[
 m_b(t)
 =
 b+c_{\parallel}t-\kappa_{\parallel}\log(1+t)+O(1)
 \quad\text{ as } t\to\infty.
\]

For $(n,r)\in\Z\times\N_0$, recall the sets
\[
 \mathcal C_{n,r}
 =
 {\{x\in\Tq:\ell_{\mathcal G}(x)=n,\ s(x)=r\}}.
\]
If a function on $\Tq$ is constant on each $\mathcal C_{n,r}$,
denote its value on that set by $U_{n,r}$.
The linearization of \eqref{eq:tree-cauchy} at zero induces the following
operator on arrays $U=(U_{n,r})_{(n,r)\in\Z\times\N_0}$, with each
coordinate evaluated by a finite sum: for $n\in\Z$,
\[
 \begin{aligned}
 (\mathcal A_{\parallel} U)_{n,0}
 &=
 D\bigl(U_{n-1,0}+U_{n+1,0}+(q-1)U_{n,1}
       -(q+1)U_{n,0}\bigr)+a_0U_{n,0},\\
 (\mathcal A_{\parallel} U)_{n,r}
 &=
 D\bigl(U_{n,r-1}+qU_{n,r+1}-(q+1)U_{n,r}\bigr)
 +a_rU_{n,r},
 \quad\text{ for }r\ge1.
 \end{aligned}
\]
Recall that $\nu_0=1$ and $\nu_r=(q-1)q^{r-1}$ for $r\ge1$.
Choose a compact interval $I\Subset(0,\infty)$ containing
$\lambda_{\parallel}$ in its interior, sufficiently small for the
constructions below. Unless further dependence is stated, constants in
the chain estimates depend only on $D,q$, the periodic coefficient
sequence $(a_r)$, and $I$. For $\lambda\in I$, let
$\psi^\lambda>0$ be the positive $\ell^2(\N_0)$ eigenvector of
$\mathcal J_{\parallel}(\lambda)$ associated with
$\Lambda_{\parallel}(\lambda)$, normalized by
$\psi^\lambda_0=1$, and put
\[
 h_r^\lambda
 =
 \frac{\psi_r^\lambda}{\sqrt{\nu_r}},\qquad r\in\N_0.
\]
For $\lambda\in I$, define
$H_\lambda(n,r)=e^{-\lambda n}h_r^\lambda$ for
$(n,r)\in\Z\times\N_0$.
The eigenvalue equation
$\mathcal J_{\parallel}(\lambda)\psi^\lambda
=\Lambda_{\parallel}(\lambda)\psi^\lambda$ gives
$\mathcal A_{\parallel}H_\lambda
=\Lambda_{\parallel}(\lambda)H_\lambda$ pointwise.

For the chain $(N_t,R_t)$ constructed in the next lemma,
$\mathbb P^\lambda_{(j,r)}$ and $\mathbb E^\lambda_{(j,r)}$ denote
its law and expectation when $(N_0,R_0)=(j,r)$. Its filtration is the completed,
right-continuous natural filtration of $(N_t,R_t)$, as in the chain
conventions of Subsection~\ref{sec:periodic-chain}. We write
$(\widehat N_t,\widehat R_t)$ for the chain with the same transverse
rates and with horizontal rates $De^{-\lambda}$ to $(n-1,0)$ and
$De^\lambda$ to $(n+1,0)$ from $(n,0)$, for $n\in\Z$.
Lemma~\ref{lem:tangential-survival-entrance} proves the
time-reversal identity for these rates.

\begin{lemma}\label{lem:tangential-regeneration}

For $\lambda$ in a sufficiently small compact interval
$I\Subset(0,\infty)$ containing $\lambda_{\parallel}$, define,
for a bounded real function $F$ on $\Z\times\N_0$,
\[
 \mathcal Q_{\parallel}^\lambda F
 =
 (H_\lambda)^{-1}
 \mathcal A_{\parallel}(H_\lambda F)
 -
 \Lambda_{\parallel}(\lambda)F.
\]
Here $\mathcal A_{\parallel}(H_\lambda F)$ is evaluated pointwise
by the coordinate formula defining $\mathcal A_{\parallel}$, whose
sums contain finitely many terms. Then $\mathcal Q_{\parallel}^\lambda$ is the generator of a conservative
continuous-time Markov chain $(N_t,R_t)$ on $\Z\times\N_0$.  For $n\in\Z$, its nonzero transition rates are
\[
 \begin{aligned}
 &(n,0)\to(n-1,0)
 &&\text{at rate }De^\lambda,\\
 &(n,0)\to(n+1,0)
 &&\text{at rate }De^{-\lambda},\\
 &(n,0)\to(n,1)
 &&\text{at rate }D(q-1)\frac{h_1^\lambda}{h_0^\lambda},\\
 &(n,r)\to(n,r-1)
 &&\text{at rate }D\frac{h_{r-1}^\lambda}{h_r^\lambda},
 \quad\text{ for }r\ge1,\\
 &(n,r)\to(n,r+1)
 &&\text{at rate }Dq\frac{h_{r+1}^\lambda}{h_r^\lambda},
 \quad\text{ for }r\ge1.
 \end{aligned}
\]
The transverse chain $R_t$ is positive recurrent and reversible with
probability vector
\[
 \pi_r^\lambda
 =
 \frac{(\psi_r^\lambda)^2}
 {\sum_{j\ge0}(\psi_j^\lambda)^2},\qquad r\in\N_0.
\]

For a process started at $(j,0)$, let $\tau_1<\tau_2<\cdots$ be the
times at which $N_t$ changes, and put
\[
 \bigl(T_k,\xi_k\bigr)
 =
 \bigl(\tau_k-\tau_{k-1},
       N_{\tau_k}-N_{\tau_{k-1}}\bigr),
 \qquad \tau_0=0,\quad k\in\N,
\]
where $T_k$ is the duration and $\xi_k$ the longitudinal displacement
between these jumps. Here and below,
$\mathbb P^\lambda=\mathbb P^\lambda_{(0,0)}$ and
$\mathbb E^\lambda=\mathbb E^\lambda_{(0,0)}$ in assertions about
$(T_k,\xi_k)$. The interval $[\tau_{k-1},\tau_k]$ is the $k$th regeneration cycle.
The pairs $(T_k,\xi_k)$ are independent and identically distributed,
and $T_k$ is independent of $\xi_k$.  Moreover,
\[
 \mathbb P^\lambda(\xi_k=1)
 =
 \frac{e^{-\lambda}}{2\cosh\lambda},
 \qquad
 \mathbb P^\lambda(\xi_k=-1)
 =
 \frac{e^\lambda}{2\cosh\lambda},\qquad k\in\N.
\]
There is $\eta_0>0$, independent of $\lambda\in I$, such that
\[
 \sup_{\lambda\in I}\mathbb E^\lambda e^{\eta_0T_1}<\infty.
\]
Recall that, for $E>\mu$,
$\mathfrak{s}(E)=E-a_0+D(q-1)-D^2(q-1)m_1(E)$, where
$m_1(E)=\langle\delta_1,(E-\mathcal B_1)^{-1}\delta_1\rangle$
and $\mathcal B_1$ is the zero-exterior Jacobi operator on
$\ell^2(\{1,2,\ldots\})$ from Lemma~\ref{lem:geodesic-resolvent}.
The Laplace transform of $T_1$ is
\[
 \mathbb E^\lambda e^{-sT_1}
 =
 \frac{2D\cosh\lambda}
 {\mathfrak{s}(\Lambda_{\parallel}(\lambda)+s)+2D}
\]
for real $s$ in a neighborhood of zero that is independent of
$\lambda\in I$. Consequently
\[
 \mathbb E^\lambda T_1
 =
 \frac{\mathfrak{s}'(\Lambda_{\parallel}(\lambda))}
 {2D\cosh\lambda},
\]
and, with
\[
 c_\lambda
 =
 \Lambda_{\parallel}'(\lambda),
\]
one has
\[
 \mathbb E^\lambda
 \bigl(\xi_1+c_\lambda T_1\bigr)=0.
\]
At $\lambda=\lambda_{\parallel}$, $c_\lambda=c_{\parallel}$.
Moreover,
$\inf_{\lambda\in I}\operatorname{Var}^\lambda(\xi_1+c_\lambda T_1)>0$,
and there is $\eta_1>0$ such that
$\sup_{\lambda\in I}\mathbb E^\lambda
 e^{\eta_1|\xi_1+c_\lambda T_1|}<\infty$.

\end{lemma}

\begin{proof}

For a finitely supported function $F$ on $\Z\times\N_0$,
substitution of $H_\lambda(n,r)=e^{-\lambda n}h_r^\lambda$ gives
\[
 (\mathcal Q_{\parallel}^\lambda F)(x)
 =\sum_{x'\ne x}q_\lambda(x,x')\bigl(F(x')-F(x)\bigr),
 \qquad x=(n,r)\in\Z\times\N_0.
\]
Here $q_\lambda(x,x')$ is the transition rate from $x$ to $x'$
listed in the statement and is zero for all other distinct pairs.
The coefficient of $F(x)$ is $-\sum_{x'\ne x}q_\lambda(x,x')$ because
$\mathcal A_{\parallel}H_\lambda
=\Lambda_{\parallel}(\lambda)H_\lambda$; hence the row sums vanish. The equation
$\mathcal J_{\parallel}(\lambda)\psi^\lambda
=\Lambda_{\parallel}(\lambda)\psi^\lambda$ and positivity of
$\psi^\lambda$ give uniform upper and positive lower bounds for
$\psi_{r+1}^\lambda/\psi_r^\lambda$ for $r\in\N_0$ and $\lambda\in I$.
Indeed, for $r\ge2$,
\[
 D\sqrt q\,\psi_{r+1}^\lambda
 \le (\Lambda_{\parallel}(\lambda)-a_r+D(q+1))\psi_r^\lambda.
\]
For the same $r\ge2$, the equation at coordinate $r+1$ gives
\[
 D\sqrt q\,\psi_r^\lambda
 \le(\Lambda_{\parallel}(\lambda)-a_{r+1}+D(q+1))
       \psi_{r+1}^\lambda.
\]
The coefficients in these inequalities are uniformly bounded for
$\lambda\in I$, since the $a_r$ are periodic. The finitely many ratios
involving $r=0,1,2$ are positive and continuous on $I$.
Since $h_r^\lambda=\psi_r^\lambda/\sqrt{\nu_r}$ and
$\nu_{r+1}/\nu_r=q$ for $r\ge1$, the ratios
$h_{r+1}^\lambda/h_r^\lambda$ also have uniform upper and positive
lower bounds. Consequently there is $C_I<\infty$ such that
$\sup_{\lambda\in I,\,x}\sum_{x'\ne x}q_\lambda(x,x')\le C_I$.
The coordinate expression therefore defines a bounded operator on $\ell^\infty(\Z\times\N_0)$,
with norm at most $2C_I$. Choose a Poisson clock rate
$\Gamma>C_I$. At a clock ring in state $x$, jump
to $x'$ with probability
$q_\lambda(x,x')/\Gamma$, and otherwise retain $x$, where
$q_\lambda$ denotes the displayed transition rates. These probabilities
sum to one. Only finitely many rings occur on each bounded time
interval, and the infinitesimal jump rates are $q_\lambda$; this gives
the asserted conservative chain
\cite[Sections~2.6--2.8]{Norris}.
For fixed $\lambda\in I$ and distinct $r,s\in\N_0$, write
$q_R(r,s)$ for the rate at which $R$ jumps from $r$ to $s$,
suppressing the dependence on $\lambda$, as read from the rates of
$\mathcal Q_{\parallel}^\lambda$. For $r\ge1$, both sides of
$\pi_r^\lambda q_R(r,r+1)=\pi_{r+1}^\lambda q_R(r+1,r)$ equal
\[
 \frac{Dq\nu_r h_r^\lambda h_{r+1}^\lambda}
      {\sum_{j\ge0}(\psi_j^\lambda)^2},
\]
because $\nu_{r+1}=q\nu_r$. At $r=0$ they both equal
$D(q-1)h_0^\lambda h_1^\lambda/
\sum_{j\ge0}(\psi_j^\lambda)^2$, since $\nu_1=q-1$.
Thus $\pi_r^\lambda q_R(r,s)=\pi_s^\lambda q_R(s,r)$ for all
$r,s\in\N_0$ with $r\ne s$, proving detailed balance for
$\pi^\lambda$.
Furthermore, for $\lambda\in I\Subset(0,\infty)$,
$\chi(\lambda)>0=\mathfrak{s}(\sigma_{\mathcal G})$. Since $\mathfrak{s}$ is strictly
increasing, \eqref{eq:tangential-schur} gives
$\Lambda_{\parallel}(\lambda)>\sigma_{\mathcal G}$. Continuity on $I$
therefore yields
\[
 g:=\min_{\lambda\in I}
       (\Lambda_{\parallel}(\lambda)-\sigma_{\mathcal G})>0.
\]
To obtain a bound on $\psi_r^\lambda$ uniform in $r\in\N_0$ and
$\lambda\in I$, put $E=\Lambda_{\parallel}(\lambda)$.
For $r\ge2$, the eigenvalue equation gives the periodic recurrence
\[
 \binom{\psi_{r+1}^\lambda}{\psi_r^\lambda}
 =A_r(E)\binom{\psi_r^\lambda}{\psi_{r-1}^\lambda},\qquad
 A_r(E)=
 \begin{pmatrix}
 (E-a_r+D(q+1))/(D\sqrt q)&-1\\1&0
 \end{pmatrix}.
\]
The eigenvalues of the one-period matrix
$A_{P+1}(E)\cdots A_2(E)$ are the Floquet multipliers. For $E>\mu$, these multipliers are $\vartheta(E)$ and
$\vartheta(E)^{-1}$, with $0<\vartheta(E)<1$. The
square-summable solution belongs to the eigenspace of the one-period
matrix corresponding to $\vartheta(E)$; see
\cite[Section~7.1]{TeschlJacobi}. For $\lambda\in I$, the values
$E=\Lambda_{\parallel}(\lambda)$ range over a compact subset of
$(\mu,\infty)$. The stable multiplier $\vartheta(E)$ is continuous
on this compact set, so $\sup_{\lambda\in I}\vartheta(E)<1$.
The transfer relation therefore gives
\[
 \binom{\psi_{kP+2}^\lambda}{\psi_{kP+1}^\lambda}
 =\vartheta(E)^k\binom{\psi_2^\lambda}{\psi_1^\lambda}
 \quad\text{ for }k\ge0.
\]
The initial two components and the finitely many transfer matrices
within one period are uniformly bounded on $I$. Since
$\sup_{\lambda\in I}\vartheta(E)<1$, this identity gives constants
$C,\rho>0$ such that
\[
 \psi_r^\lambda\le Ce^{-\rho r},
 \qquad h_r^\lambda=\psi_r^\lambda/\sqrt{\nu_r}
\]
uniformly for $r\in\N_0$ and $\lambda\in I$. Hence
$1\le\sum_{r\ge0}(\psi_r^\lambda)^2\le C$, because
$\psi_0^\lambda=1$. Define the diagonal entries of the transverse
generator by $q_R(r,r)=-\sum_{s\ne r}q_R(r,s)$.
Summing the detailed-balance identities over the neighbors of each
$r$ gives $(\pi^\lambda)^Tq_R=0$. The uniformly bounded rates
justify the exponential series for $e^{tq_R}$ and give
$(\pi^\lambda)^Te^{tq_R}=(\pi^\lambda)^T$ for every $t\ge0$.
Thus $\pi^\lambda$ is invariant. The transverse chain is irreducible
and nonexplosive, so this invariant probability makes it positive
recurrent. In
particular, its mean return time to zero after leaving zero is finite;
see \cite[Sections~3.4--3.5]{Norris}.

Let $L_t^0=\int_0^t\mathbf1_{\{R_s=0\}}\,\dd s$ be the time
spent at transverse state zero. Independently of $R$, take Poisson
processes $P^+$ and $P^-$ with rates $De^{-\lambda}$ and
$De^\lambda$, respectively, and set
\[
 N_t=j+P^+(L_t^0)-P^-(L_t^0),\qquad t\ge0,\quad j\in\Z.
\]
A ring of $P^+$ or $P^-$ changes only $N_t$ and occurs in physical
time precisely while $R_t=0$. Thus the time-change construction
has the two horizontal rates in the statement and preserves the
transverse rates. Recurrence gives infinitely many visits of
$R$ to zero. The holding times there are independent exponentials with
a fixed positive rate for each $\lambda$, so $L_t^0\to\infty$
almost surely as $t\to\infty$. Consequently every $\tau_k$ is
finite almost surely. Each $\tau_k$ is a stopping time for the natural
filtration of $(N_t,R_t)$. At each $\tau_k$, the strong Markov property
implies that, conditionally on $\tau_k$ and $N_{\tau_k}$, the translated future
path $(N_{\tau_k+t}-N_{\tau_k},R_{\tau_k+t})_{t\ge0}$ is independent
of the past and has the law of the chain started from $(0,0)$.
Indeed, the state at $\tau_k$ has transverse coordinate zero, and
longitudinal translation preserves all jump rates; see
\cite[Section~6.5]{Norris} for the strong Markov property.
Iterating this conditional law at $\tau_1,\tau_2,\ldots$ proves
that the pairs $(T_k,\xi_k)$ are independent and identically
distributed. The superposition $P^++P^-$ has rate
$2D\cosh\lambda$ and independent direction marks with the displayed
probabilities. These marks are independent of the unmarked horizontal
ring times and of $R$, whereas $T_k$ is determined by those times
and the transverse path. Hence $T_k$ and $\xi_k$ are independent.

Kill the chain at its first horizontal jump, retaining the horizontal
rates as killing rates at $r=0$. The killed transverse generator is
\[
 (K^\lambda F)_r=(q_RF)_r
 -2D\cosh\lambda\,\mathbf1_{\{r=0\}}F_0,\qquad r\in\N_0,\quad F\in\ell^2(\pi^\lambda).
\]
It is a bounded self-adjoint operator on $\ell^2(\pi^\lambda)$.
Under the unitary map
$F_r\mapsto\sqrt{\pi_r^\lambda}F_r$, it becomes
\[
 \mathcal J_{\mathcal G}
 -
 2D\langle\delta_0,\,\cdot\,\rangle\delta_0
 -
 \Lambda_{\parallel}(\lambda).
\]
The quadratic-form inequality
$\mathcal J_{\mathcal G}-2D\langle\delta_0,\,\cdot\,\rangle\delta_0
\le\mathcal J_{\mathcal G}$ shows that the spectral upper endpoint
of $K^\lambda$ is at most
$\sigma_{\mathcal G}-\Lambda_{\parallel}(\lambda)\le-g<0$,
where $g>0$ is independent of $\lambda\in I$. For $t>0$, the density of $T_1$ is
\[
 2D\cosh\lambda\,
 \bigl(e^{tK^\lambda}\bigr)(0,0).
\]
Here $(e^{tK^\lambda})(0,0)$ is the killed transition
probability from zero to zero; equivalently, it is the scalar semigroup
entry at the unit vector $\delta_0/\sqrt{\pi_0^\lambda}$ in
$\ell^2(\pi^\lambda)$. It is therefore at most $e^{-gt}$.
Integration against $e^{\eta_0t}$, with $0<\eta_0<g$, gives the
asserted uniform exponential moment. By the block decomposition in \eqref{eq:geodesic-schur}, the diagonal
resolvent of
$\mathcal J_{\mathcal G}
 -2D\langle\delta_0,\,\cdot\,\rangle\delta_0$
at the origin, for $E>\sigma_{\mathcal G}$, is
\[
 \frac1{\mathfrak{s}(E)+2D}.
\]
For $|s|$ sufficiently small, uniformly on $I$, we have
$\Lambda_{\parallel}(\lambda)+s>\sigma_{\mathcal G}$.
The resolvent identity at
$E=\Lambda_{\parallel}(\lambda)+s$ and Laplace transformation give
\[
 \mathbb E^\lambda e^{-sT_1}
 =
 \frac{2D\cosh\lambda}
 {\mathfrak{s}(\Lambda_{\parallel}(\lambda)+s)+2D}.
\]
At $s=0$ the denominator equals $2D\cosh\lambda$ because
\[
 \mathfrak{s}(\Lambda_{\parallel}(\lambda))
 =
 2D(\cosh\lambda-1).
\]
Differentiating the displayed Laplace transform at $s=0$ gives
$\mathbb E^\lambda T_1=\mathfrak{s}'(\Lambda_{\parallel}(\lambda))/(2D\cosh\lambda)$.
Differentiating
$\mathfrak{s}(\Lambda_{\parallel}(\lambda))=2D(\cosh\lambda-1)$ in
$\lambda$ gives
\[
 \mathfrak{s}'(\Lambda_{\parallel}(\lambda))
 \Lambda_{\parallel}'(\lambda)
 =
 2D\sinh\lambda.
\]
Hence
\[
 c_\lambda\mathbb E^\lambda T_1
 =
 \tanh\lambda.
\]
Since
\[
 \mathbb E^\lambda\xi_1
 =
 -\tanh\lambda,
\]
we have $\mathbb E^\lambda(\xi_1+c_\lambda T_1)=0$.
Independence of $T_1$ and $\xi_1$
gives
\[
 \operatorname{Var}(\xi_1+c_\lambda T_1)
 \ge \operatorname{Var}(\xi_1)
 =\frac1{\cosh^2\lambda}>0,
\]
uniformly on $I$. Set $M_I=\sup_{\lambda\in I}|c_\lambda|<\infty$
and choose $\eta_1=\eta_0/[2(1+M_I)]>0$. Since
$|\xi_1+c_\lambda T_1|\le1+M_IT_1$ and
$\eta_1M_I<\eta_0$, the exponential-moment bound for $T_1$ gives
\[
 \sup_{\lambda\in I}\mathbb E^\lambda
 e^{\eta_1|\xi_1+c_\lambda T_1|}
 \le e^{\eta_1}\sup_{\lambda\in I}\mathbb E^\lambda
 e^{\eta_0T_1}<\infty.
\]\qedhere

\end{proof}

\begin{lemma}\label{lem:tangential-free-kernel}

For the chain $(N_t,R_t)$ of
Lemma~\ref{lem:tangential-regeneration}, for $t\ge0$, $n\in\Z$ and
$\lambda\in I$, define
\[
 p_t^\lambda(n)
 =\mathbb P^\lambda_{(0,0)}\{(N_t,R_t)=(n,0)\}.
\]
For every $t\ge0$ and $n\in\Z$, uniformly for $\lambda\in I$,
\[
 p_t^\lambda(n)
 \le
 C(1+t)^{-1/2}.
\]
Moreover, for every finite $A\ge0$ there are $c_A>0$ and a finite
$T_A\ge1$ such that
\[
 p_t^\lambda(n)
 \ge
 c_A t^{-1/2}
\]
whenever
\[
 t\ge T_A,
 \qquad
 |n+c_\lambda t|\le A\sqrt t.
\]
In fact,
\[
 p_t^\lambda(n)
 =
 \frac{\pi_0^\lambda}
 {\sqrt{2\pi\Lambda_{\parallel}''(\lambda)t}}
 \exp\left(
 -\frac{(n+c_\lambda t)^2}
 {2\Lambda_{\parallel}''(\lambda)t}
 \right)
 +o(t^{-1/2})\quad\text{ as } t\to\infty
\]
uniformly for $\lambda\in I$ and $|n+c_\lambda t|\le A\sqrt t$,
for each fixed finite $A\ge0$.

\end{lemma}

\begin{proof}

For $t\ge0$, $\lambda\in I$ and $\theta\in[-\pi,\pi]$, use the characteristic-function
convention
\[
 \Phi_t^\lambda(\theta)
 =\mathbb E^\lambda_{(0,0)}
   [e^{\mathrm i\theta N_t};R_t=0]
 =\sum_{n\in\Z}e^{\mathrm i n\theta}p_t^\lambda(n).
\]
The series converges absolutely because
$\sum_{n\in\Z}p_t^\lambda(n)\le1$.
For $z\in\C$, extend the operator family in
\eqref{eq:tangential-fiber} by
$\mathcal J_{\parallel}(z)
=\mathcal J_{\mathcal G}+2D(\cosh z-1)P_0$ on $\ell^2(\N_0)$,
where $P_0v=v_0\delta_0$ for $v\in\ell^2(\N_0)$.
Applying the coordinate formula for $\mathcal Q_{\parallel}^\lambda$
to $e^{\mathrm i\theta n}F_r$ gives a bounded operator on
$\ell^2(\pi^\lambda)$ acting on the transverse vector $F$.
The unitary map $F_r\mapsto\sqrt{\pi_r^\lambda}F_r$ from
$\ell^2(\pi^\lambda)$ onto $\ell^2(\N_0)$ conjugates this Fourier
operator to
\[
 \mathcal J_{\parallel}(\lambda-\mathrm i\theta)
 -
 \Lambda_{\parallel}(\lambda).
\]
The diagonal factors of the unitary map cancel in its $(0,0)$
semigroup entry. Thus
$\Phi_t^\lambda(\theta)=(e^{tA_\theta})_{00}$, where
$A_\theta=\mathcal J_{\parallel}(\lambda-\mathrm i\theta)
-\Lambda_{\parallel}(\lambda)$ on $\ell^2(\N_0)$.
The notation $A_\theta$ suppresses its dependence on
$\lambda\in I$. All operator and vector norms in this proof are taken
on $\ell^2(\N_0)$ after this conjugation. In particular,
\[
 p_t^\lambda(n)=\frac1{2\pi}\int_{-\pi}^{\pi}
 e^{-\mathrm i n\theta}(e^{tA_\theta})_{00}\,\dd\theta,\qquad t\ge0,\quad n\in\Z.
\]
The simple isolated eigenvalue $\Lambda_{\parallel}(\lambda)$ of
$\mathcal J_{\parallel}(\lambda)$ continues analytically for small
$\theta$, uniformly on $I$. Write
$\Lambda_{\parallel}(\lambda-\mathrm i\theta)$ for this
continuation. The corresponding eigenvalue of $A_\theta$ is
\[
 \zeta_\lambda(\theta)
 =
 \Lambda_{\parallel}(\lambda-\mathrm i\theta)
 -
 \Lambda_{\parallel}(\lambda).
\]
Taylor expansion gives
\[
 \zeta_\lambda(\theta)
 =
 -\mathrm i c_\lambda\theta
 -
 \frac12\Lambda_{\parallel}''(\lambda)\theta^2
 +
 O(|\theta|^3)\quad\text{ as } \theta\to0,
\]
uniformly for $\lambda\in I$.  The coefficient of $\theta^2$ is
strictly negative in real part.

Since $A_0=\mathcal J_{\mathcal G}
+2D(\cosh\lambda-1)P_0-\Lambda_{\parallel}(\lambda)$,
the min--max principle for a positive rank-one perturbation puts
all spectral values except its simple eigenvalue zero at most
$\sigma_{\mathcal G}-\Lambda_{\parallel}(\lambda)<0$.
Set $g=\min_{\lambda\in I}
(\Lambda_{\parallel}(\lambda)-\sigma_{\mathcal G})>0$.
The formula for $\mathcal J_{\parallel}(\lambda-\mathrm i\theta)$
gives $\|A_\theta-A_0\|\le C|\theta|$ for $|\theta|\le1$ and
$\lambda\in I$.
A circle of radius less than $g/3$ about zero therefore separates the
simple eigenvalue for every sufficiently small $|\theta|$ and every
$\lambda\in I$, by the resolvent Neumann series. Let $P_\lambda(\theta)$ be the Riesz projection of $A_\theta$
onto its simple eigenvalue $\zeta_\lambda(\theta)$. The remaining spectrum
is enclosed by a fixed bounded contour in $\{\re z\le-g/2\}$ after
reducing the allowed $|\theta|$. The resolvents on that contour are
uniformly bounded, again by the Neumann series from the self-adjoint
operators $A_0$. The Dunford integral over the contour in
$\{\re z\le-g/2\}$ gives
\[
 \|e^{tA_\theta}(I-P_\lambda(\theta))\|\le Ce^{-gt/2},\qquad t\ge0.
\]

For real $\theta\in[-\pi,\pi]$ and $\lambda\in I$, the
self-adjoint part of $A_\theta$ is
\[
 \frac{A_\theta+A_\theta^*}{2}
 = A_0-2D\cosh\lambda(1-\cos\theta)P_0.
\]
Put $e_\lambda=\psi^\lambda/\|\psi^\lambda\|_2$ and
$m_I=\inf_{\lambda\in I}(e_\lambda)_0>0$, using the bound on
$\|\psi^\lambda\|_2$ and $\psi_0^\lambda=1$. Fix $\theta_0>0$ and
write $d_\theta=2D\cosh\lambda(1-\cos\theta)$. For
$\theta_0\le|\theta|\le\pi$, $d_\theta\ge d_0>0$ uniformly on $I$.
For $\theta_0\le|\theta|\le\pi$ and $\lambda\in I$, decompose
$v\in\ell^2(\N_0)$ with $\|v\|_2=1$ as $v=a e_\lambda+w$,
where $a\in\C$ and $w\perp e_\lambda$.
The bound $\langle w,A_0w\rangle\le-g\|w\|_2^2$ gives
\[
 \re\langle v,A_\theta v\rangle
 \le-g\|w\|_2^2-d_\theta|v_0|^2.
\]
If $\|w\|_2\ge m_I/4$, the first term is at most
$-gm_I^2/16$. Otherwise $|a|=\sqrt{1-\|w\|_2^2}$ and
$|v_0|\ge |a|(e_\lambda)_0-\|w\|_2\ge m_I/2$, so the second term
is at most $-d_0m_I^2/4$. Thus the real part is bounded above by a
common negative constant $-\gamma_0$. Differentiating the squared
norm gives, for every $v\in\ell^2(\N_0)$,
\[
 \frac{\dd}{\dd t}\|e^{tA_\theta}v\|_2^2
 \le-2\gamma_0\|e^{tA_\theta}v\|_2^2,
 \qquad \|e^{tA_\theta}\|\le e^{-\gamma_0t},
 \quad t\ge0,\quad\theta_0\le|\theta|\le\pi,\quad\lambda\in I.
\]

Let $a_\lambda(\theta)=(P_\lambda(\theta))_{00}$. The contour
formula for $P_\lambda(\theta)$ and the uniform resolvent bounds give
$a_\lambda(\theta)=\pi_0^\lambda+O(|\theta|)$ uniformly on $I$,
because $P_\lambda(0)$ is the orthogonal projection onto $e_\lambda$.
The contour formula and uniform resolvent bounds also give a uniform
third-order Taylor remainder for $\zeta_\lambda(\theta)$. Since
$\inf_{\lambda\in I}\Lambda_{\parallel}''(\lambda)>0$, reducing
$\theta_0$ yields
$\re\zeta_\lambda(\theta)\le-b\theta^2$ for $|\theta|\le\theta_0$,
with $b>0$ independent of $\lambda$. Fourier inversion is therefore
\[
 p_t^\lambda(n)=\frac1{2\pi}\int_{-\theta_0}^{\theta_0}
 e^{-\mathrm i n\theta+t\zeta_\lambda(\theta)}
 a_\lambda(\theta)\,\dd\theta+O(e^{-bt}),
\]
after replacing $b$ by $\min\{b,g/2,\gamma_0\}$ to include
the bounds for $e^{tA_\theta}(I-P_\lambda(\theta))$ and for
$e^{tA_\theta}$ on $\theta_0\le|\theta|\le\pi$. Taking absolute values bounds the integral by
$C\int_\R e^{-bt\theta^2}\,\dd\theta$, proving the upper bound for
$t\ge1$; the probability bound $p_t^\lambda(n)\le1$ covers $t\le1$.
After substituting $\theta=\eta/\sqrt t$ in the integral over
$[-\theta_0,\theta_0]$, the leading expression for
$e^{-\mathrm i n\eta/\sqrt t+t\zeta_\lambda(\eta/\sqrt t)}
a_\lambda(\eta/\sqrt t)$ is
\[
 \pi_0^\lambda
 \exp\left(-\mathrm i\frac{n+c_\lambda t}{\sqrt t}\eta
           -\frac12\Lambda_{\parallel}''(\lambda)\eta^2\right).
\]
On each fixed bounded interval of $\eta$, the Taylor remainders
tend to zero as $t\to\infty$, uniformly in $\lambda\in I$ and
$|(n+c_\lambda t)/\sqrt t|\le A$. Outside that interval the
majorant $Ce^{-b\eta^2}$ has an arbitrarily small integral,
uniformly for $\lambda\in I$ and $|n+c_\lambda t|\le A\sqrt t$.
Evaluating the Gaussian integral therefore gives the claimed
expansion of $p_t^\lambda(n)$ with an error $o(t^{-1/2})$ uniform
for $\lambda\in I$ and $|n+c_\lambda t|\le A\sqrt t$. Finally,
$\inf_I\pi_0^\lambda>0$ and positive upper and lower bounds for
$\Lambda_{\parallel}''$ give the stated lower bound whenever
$|n+c_\lambda t|\le A\sqrt t$.\qedhere

\end{proof}

The next lemma bounds sums of independent centered increments.
The exit argument in Lemma~\ref{lem:centered-cycle-exit} uses the
two inequalities in \eqref{eq:s9tan-residual-bound}, which we verify
for $\xi_k+c_\lambda T_k$ and its negative in
Lemma~\ref{lem:tangential-survival-entrance}. The stopping
identities used in the proof are the discrete-time optional-stopping
identities of \cite[Section~4.8]{Durrett2019}; their integrability
conditions are checked below.

\begin{lemma}\label{lem:centered-cycle-exit}

Let $(Z_k)_{k\in\N}$ be independent identically distributed real
random variables with $\mathbb EZ_1=0$,
$\operatorname{Var}(Z_1)>0$, and
\[
 \mathbb E e^{\eta|Z_1|}<\infty
\]
for some $\eta>0$. Assume also that, for a constant $C_{\mathrm o}<\infty$,
\begin{equation}
 \label{eq:s9tan-residual-bound}
 \begin{aligned}
 \mathbb E[(Z_1-a)^2;Z_1>a]
   &\le C_{\mathrm o}\mathbb P(Z_1>a),\\
 \mathbb E[(-Z_1-a)^2;-Z_1>a]
   &\le C_{\mathrm o}\mathbb P(-Z_1>a)
 \end{aligned}
 \quad\text{ for } a\ge0.
\end{equation}
For a family of increment laws, assume common constants $\eta$,
$C_{\mathrm o}$, a common upper bound for
$\mathbb E e^{\eta|Z_1|}$, and a common positive lower bound for
$\operatorname{Var}(Z_1)$. There is $C<\infty$, depending
only on these bounds, such that, for
$y\ge0$ and $n\in\N_0$,
\[
 \mathbb P\left\{
 y+\sum_{j=1}^kZ_j>0,\ 1\le k\le n
 \right\}
 \le
 C\frac{1+y}{\sqrt{1+n}}.
\]
There are $y_0\ge1$ and $\delta>0$ such that, for every $A>0$, there
are $B_A<\infty$, $\kappa_A>0$, and $n_A<\infty$ with the following
property.  If
\[
 n\ge n_A,
 \qquad
 y_0\le y\le A\sqrt n,
\]
then there are a stopping time $J\le n/3$ for
$\mathcal F_k=\sigma(Z_1,\ldots,Z_k)$, $k\in\N_0$, and an event
$E\in\mathcal F_J$ with $\mathbb P(E)\ge\kappa_Ay/\sqrt n$ such that,
on $E$, $y+\sum_{j=1}^kZ_j>0$ for $0\le k\le J$ and
\[
 \delta\sqrt n
 \le
 y+\sum_{j=1}^JZ_j
 \le
 B_A\sqrt n.
\]

\end{lemma}

\begin{proof}

For $R\ge1$ and $y\ge0$, write
\[
 S_k=y+\sum_{j=1}^kZ_j,\quad k\in\N_0,
 \qquad
 \tau_R=\inf\{k\ge0:S_k\notin(0,R)\},
\]
and, for $R\ge1$ and $0<y<R$, define the upper and lower overshoots by
\[
 O_R^+=(S_{\tau_R}-R)\mathbf1_{\{S_{\tau_R}\ge R\}},
 \qquad
 O_R^-=-S_{\tau_R}\mathbf1_{\{S_{\tau_R}\le0\}}.
\]
Both variables are defined to be zero on $\{\tau_R=\infty\}$
until almost-sure exit has been proved. We use
$\inf\varnothing=\infty$ in the definition of $\tau_R$.
Let $\mathcal F_k=\sigma(Z_1,\ldots,Z_k)$. The processes $S_k$
and $S_k^2-k\sigma^2$, where $\sigma^2=\operatorname{Var}(Z_1)$,
are integrable martingales for this filtration; $\tau_R\wedge m$
is a bounded stopping time. On $\{\tau_R>k-1\}$ the distances
$R-S_{k-1}$ and $S_{k-1}$ are positive. At equality with either
endpoint the squared overshoot is zero, so the strict-tail assumptions
in \eqref{eq:s9tan-residual-bound} also bound exits with weak
endpoint inequalities. Independence and \eqref{eq:s9tan-residual-bound} give
\[
 \begin{aligned}
 &\mathbb E[(S_k-R)^2;
       \tau_R=k,\ S_k\ge R\mid\mathcal F_{k-1}]\\
 &\qquad\le C_{\mathrm o}
       \mathbb P(\tau_R=k,\ S_k\ge R\mid\mathcal F_{k-1}),
 \end{aligned}
 \qquad k\in\N.
\]
At the lower endpoint, the second inequality in
\eqref{eq:s9tan-residual-bound}, applied with $a=S_{k-1}$, gives
\[
 \mathbb E[S_k^2;\tau_R=k,\ S_k\le0\mid\mathcal F_{k-1}]
 \le C_{\mathrm o}\mathbb P(\tau_R=k,\ S_k\le0\mid\mathcal F_{k-1}),
 \quad k\in\N.
\]
Sum the two conditional estimates over $1\le k\le m$, take
expectations, and use disjointness of the events $\{\tau_R=k\}$.
Monotone convergence as $m\to\infty$ gives
\[
 \sup_{R\ge1,\ 0<y<R}
 \mathbb E\bigl[(O_R^+)^2+(O_R^-)^2;\tau_R<\infty\bigr]
 \le 2C_{\mathrm o}.
\]
To justify passage to the unbounded exit time, first stop the
quadratic martingale at $\tau_R\wedge m$. Before exit,
$S_k^2\le R^2$; at exit, the overshoot bounds control
$S_{\tau_R}^2$. Thus bounded optional stopping gives
\[
 \sigma^2\mathbb E(\tau_R\wedge m)
 =\mathbb ES_{\tau_R\wedge m}^2-y^2
 \le R^2+2R\sqrt{2C_{\mathrm o}}+2C_{\mathrm o},\qquad m\in\N.
\]
Thus $\mathbb E\tau_R<\infty$, and the exit occurs almost surely.
On $\{\tau_R>m\}$ we have $0<S_m<R$; on
$\{\tau_R\le m\}$ the stopped value equals $S_{\tau_R}$. Hence
$S_{\tau_R\wedge m}^2\le R^2+S_{\tau_R}^2$. The overshoot bounds
make the right side integrable. Dominated convergence therefore justifies
$\mathbb ES_{\tau_R}=y$ and
$\mathbb ES_{\tau_R}^2-y^2=\sigma^2\mathbb E\tau_R$.

For the remaining exit estimates take $R\ge1$, $0<y<R$, and put
$p_R=\mathbb P(S_{\tau_R}\ge R)$. Optional stopping for $S_k$,
first at bounded stopping times and then by the integrable domination
just proved, gives
\[
 y=Rp_R+\mathbb EO_R^+-\mathbb EO_R^-.
\]
The bound
$\mathbb E[(O_R^+)^2+(O_R^-)^2]\le2C_{\mathrm o}$
and Cauchy--Schwarz imply
$\mathbb EO_R^++\mathbb EO_R^-\le2\sqrt{2C_{\mathrm o}}$.
The identity $y=Rp_R+\mathbb EO_R^+-\mathbb EO_R^-$ therefore gives
\[
 p_R\le C\frac{1+y}{R},
 \qquad
 p_R\ge c\frac{y}{R}\quad\text{ for } y_0\le y<R,
\]
where one may take $y_0=\max\{1,4\sqrt{2C_{\mathrm o}}\}$
and $c=1/2$: the bound
$|\mathbb EO_R^+-\mathbb EO_R^-|\le2\sqrt{2C_{\mathrm o}}$
gives $Rp_R\ge y/2$ for $y\ge y_0$. Moreover,
\[
 \mathbb ES_{\tau_R}^2
 =R^2p_R+2R\mathbb EO_R^+
   +\mathbb E\bigl[(O_R^+)^2+(O_R^-)^2\bigr]
 \le C(1+y)R.
\]
Optional stopping for
$S_k^2-k\operatorname{Var}(Z_1)$ therefore yields
\[
 \mathbb E\tau_R\le C(1+y)R.
\]
Here the stopped quadratic identity is
$\sigma^2\mathbb E\tau_R=\mathbb ES_{\tau_R}^2-y^2$.
The assumed common positive lower bound for $\sigma^2$ makes the
constant in the exit-time bound uniform over the family.

If $S_k>0$ for $0\le k\le n$, then either
$S_{\tau_R}\ge R$ and $\tau_R\le n$, or $\tau_R>n$. Thus, when $R\ge1$, $0<y<R$, and $n\in\N$,
\[
 \mathbb P\{\min_{k\le n}S_k>0\}
 \le C\frac{1+y}{R}+\frac{C(1+y)R}{n}.
\]
Taking $R=\sqrt n$ proves the bound when $0<y<\sqrt n$. For
$y\ge\sqrt n$, enlarge $C$ so that the right-hand side is at least
one. For $y=0$, the event
$\{\sum_{j=1}^kZ_j>0\text{ for }1\le k\le n\}$ is contained in
$\{1+\sum_{j=1}^kZ_j>0\text{ for }1\le k\le n\}$, so the estimate for
initial height one gives the bound after enlarging $C$.
For $n=0$, the bound holds after taking $C\ge1$; this choice
preserves dependence only on the stated family bounds.

For the entrance bound choose $\delta>0$ below and take
$R=\delta\sqrt n\ge1$. For $y_0\le y<R$,
\[
 \mathbb P\{S_{\tau_R}\ge R,\ \tau_R\le n/3\}
 \ge c\frac{y}{R}-\frac{3C(1+y)R}{n}.
\]
Since $y\ge y_0\ge1$, choose $\delta>0$ with
$6C\delta^2\le c/2$. The last display is then at least
$cy/(2\delta\sqrt n)$. Choose an integer
\[
 n_A\ge\max\left\{1,\delta^{-2},
          \left(\frac{8C_{\mathrm o}\delta}{cy_0}\right)^2\right\}.
\]
For $n\ge n_A$, the second-moment bound gives
$\mathbb P(O_R^+>\sqrt n)\le2C_{\mathrm o}/n
\le cy/(4\delta\sqrt n)$. Define
\[
 J=\tau_R\wedge\lfloor n/3\rfloor,\qquad
 E=\{\tau_R\le\lfloor n/3\rfloor,\ S_{\tau_R}\ge R,
                          \ O_R^+\le\sqrt n\}.
\]
Then $E\in\mathcal F_J$,
$\mathbb P(E)\ge cy/(4\delta\sqrt n)$, and $J=\tau_R$ on $E$.
On $E$, the walk stays positive through $J$ and
$R\le S_J\le(\delta+1)\sqrt n$. If $y\ge R$, take $J=0$
and let $E$ be the entire sample space. In both cases one may take
$B_A=\max\{A,\delta+1\}$ and
$\kappa_A=\min\{A^{-1},c/(4\delta)\}$.
The common bounds on $\mathbb E e^{\eta|Z_1|}$,
$C_{\mathrm o}$ and $\operatorname{Var}(Z_1)$ make
$y_0,\delta,B_A,\kappa_A,n_A$ uniform over the family.\qedhere

\end{proof}

\Needspace{9\baselineskip}
\begin{lemma}\label{lem:tangential-survival-entrance}

Fix $\lambda\in I$ and write $c=c_\lambda$.
Let $(\mathcal F_s)_{s\ge0}$ be the completed, right-continuous
natural filtration of the chain in
Lemma~\ref{lem:tangential-regeneration}.
For $j\in\Z$ and $y>0$, start the chain from $(j,0)$ and set
\[
 Y_s=y+N_s-j+cs,\qquad\text{ for }s\ge0.
\]
There is $C<\infty$, independent of $j,y,t$, such that, for every
$t\ge0$, we have
\[
 \mathbb P^\lambda_{(j,0)}
 \{Y_s>0\text{ for }0\le s\le t\}
 \le
 C\frac{1+y}{\sqrt{1+t}}.
\]
There are $y_0\ge1$ and $\delta>0$ such that, for each $A>0$, one can
choose finite constants $B_A>0$, $\kappa_A>0$ and $T_A\ge0$
with the following property.  If
\[
 t\ge T_A,
 \qquad
 y_0\le y\le A\sqrt t,
\]
then there are an $(\mathcal F_s)$-stopping time $T_{\mathrm{stop}}\le t/3$
and an event $E\in\mathcal F_{T_{\mathrm{stop}}}$ of probability at least
$\kappa_Ay/\sqrt t$ such that, on $E$,
$Y_s>0$ for $0\le s\le T_{\mathrm{stop}}$, $R_{T_{\mathrm{stop}}}=0$, and
\[
 \delta\sqrt t
 \le
 Y_{T_{\mathrm{stop}}}
 \le
 B_A\sqrt t.
\]
The invariant measure of the full chain is counting measure in $n$
times $\pi^\lambda$ in $r$. Its reversed chain has the same transverse
rates and interchanges the two horizontal rates. For that chain use
$\widehat Y_s=y+\widehat N_s-j-cs$.
{Start the reversed chain from $(j,0)$, and denote its law by
$\widehat{\mathbb P}^{\lambda}_{(j,0)}$. For $y>0$ and $t\ge0$,
\[
 \widehat{\mathbb P}^{\lambda}_{(j,0)}
 \{\widehat Y_s>0\text{ for }0\le s\le t\}
 \le C\frac{1+y}{\sqrt{1+t}}.
\]
With the same $y_0$ and $\delta$, for each $A>0$ and
$t\ge T_A$, $y_0\le y\le A\sqrt t$, there are a stopping time
$\widehat{T}_{\mathrm{stop}}\le t/3$ for the completed, right-continuous natural
filtration $(\widehat{\mathcal F}_s)_{s\ge0}$ of
$(\widehat N,\widehat R)$ and an event
$\widehat E\in\widehat{\mathcal F}_{\widehat{T}_{\mathrm{stop}}}$ such that
\[
 \widehat{\mathbb P}^{\lambda}_{(j,0)}(\widehat E)
 \ge\kappa_Ay/\sqrt t.
\]
On $\widehat E$, one has $\widehat Y_s>0$ for
$0\le s\le\widehat{T}_{\mathrm{stop}}$, $\widehat R_{\widehat{T}_{\mathrm{stop}}}=0$, and
$\delta\sqrt t\le\widehat Y_{\widehat{T}_{\mathrm{stop}}}\le B_A\sqrt t$,
where $\widehat R$ is the transverse coordinate of the reversed chain.}
All constants may be chosen uniformly for $\lambda\in I$.

\end{lemma}

\begin{proof}

By Lemma~\ref{lem:tangential-dispersion}, evenness and strict convexity
of $\Lambda_{\parallel}$ give $c=\Lambda_{\parallel}'(\lambda)>0$
for $\lambda\in I\Subset(0,\infty)$. At the times $\tau_k$ of horizontal
jumps, the increments of $Y_{\tau_k}$ are
\[
 Z_k=\xi_k+cT_k,\qquad k\in\N.
\]
By Lemma~\ref{lem:tangential-regeneration}, the increments are
independent and identically distributed under each law
$\mathbb P^\lambda$. There is $\eta>0$ such that
\[
 \mathbb E^\lambda Z_1=0\text{ for }\lambda\in I,\qquad
 \inf_{\lambda\in I}\operatorname{Var}^\lambda(Z_1)>0,\qquad
 \sup_{\lambda\in I}\mathbb E^\lambda e^{\eta|Z_1|}<\infty.
\]
The cycle filtration is
$\mathcal H_k=\sigma(T_i,\xi_i:1\le i\le k)$ for $k\in\N_0$;
all discrete stopping indices below are taken for this filtration.
For expectations and probabilities involving $(T_i,\xi_i)$ alone,
write $\mathbb E=\mathbb E^\lambda$ and $\mathbb P=\mathbb P^\lambda$. Constants in this proof
are independent of $j$ and $\lambda\in I$.
We verify the additional hypothesis
\eqref{eq:s9tan-residual-bound} before applying
Lemma~\ref{lem:centered-cycle-exit}. Integrating the density of $T_1$ in
Lemma~\ref{lem:tangential-regeneration} by the spectral theorem gives
\[
 \mathbb P(T_1>s)=\int_{[g,\infty)}e^{-\beta s}\,\dd\omega_\lambda(\beta),
 \qquad s\ge0,
\]
where $g>0$ is uniform on $I$. The measure $\omega_\lambda$ is obtained
by multiplying the spectral measure of $-K^\lambda$ at the unit vector
$\delta_0/\sqrt{\pi_0^\lambda}$ in $\ell^2(\pi^\lambda)$ by
$2D\cosh\lambda/\beta$. Its total mass is one because $T_1<\infty$
almost surely, as proved in Lemma~\ref{lem:tangential-regeneration}.
Thus
\[
 \mathbb E[(T_1-s)^2;T_1>s]
 \le 2g^{-2}\mathbb P(T_1>s)\quad\text{ for } s\ge0.
\]
We next verify the corresponding upper residual estimate for $Z_1$:
\[
 \mathbb E[(Z_1-a)^2;Z_1>a]
 \le C\mathbb P(Z_1>a)\quad\text{ for } a\ge0.
\]
For $\xi\in\{-1,1\}$ and $a\ge0$, put
$s_{\xi,a}=(a-\xi)/c$. Independence of $T_1$ and $\xi_1$ means that
the law of $T_1$ is unchanged by conditioning on $\xi_1=\xi$.
If $s_{\xi,a}\ge0$, the bound
$\mathbb E[(T_1-s)^2;T_1>s]\le2g^{-2}\mathbb P(T_1>s)$
for $s\ge0$ gives
\[
 \begin{aligned}
 \mathbb E[(cT_1+\xi-a)^2;cT_1+\xi>a]
 &=c^2\mathbb E[(T_1-s_{\xi,a})^2;T_1>s_{\xi,a}]\\
 &\le2c^2g^{-2}\mathbb P(T_1>s_{\xi,a}).
 \end{aligned}
\]
If $s_{\xi,a}<0$, then $\xi=1$ and $0\le a<1$. The event
$cT_1+\xi>a$ has probability one and
\[
 \mathbb E(cT_1+\xi-a)^2\le2c^2\mathbb ET_1^2+2.
\]
Since $c$ and $\mathbb ET_1^2$ are bounded uniformly on $I$,
one constant $C$ works in both cases. Averaging over the two marks
gives $\mathbb E[(Z_1-a)^2;Z_1>a]\le C\mathbb P(Z_1>a)$
for $a\ge0$. Since $Z_1\ge-1$, one has $0<-Z_1-a\le1$
on $\{-Z_1>a\}$ for $a\ge0$, and hence
\[
 \mathbb E[(-Z_1-a)^2;-Z_1>a]
 \le\mathbb P(-Z_1>a)\quad\text{ for } a\ge0.
\]
We have therefore verified \eqref{eq:s9tan-residual-bound} for
$Z_1$. Replacing $Z_1$ by $-Z_1$ interchanges the two estimates:

\[
 \begin{aligned}
 \mathbb E[(-Z_1-a)^2;-Z_1>a]
 &\le C_{\mathrm o}\mathbb P(-Z_1>a),\\
 \mathbb E[(Z_1-a)^2;Z_1>a]
 &\le C_{\mathrm o}\mathbb P(Z_1>a)
 \end{aligned}
 \qquad\text{ for }a\ge0.
\]
Both hold with $C_{\mathrm o}=\max\{1,C\}$ by the displayed bounds
for $Z_1$ and $-Z_1$.

Between horizontal jumps, $Y_s$ increases at rate $c>0$. Thus
$Y_s>0$ throughout $[0,t]$ if and only if $Y_{\tau_k}>0$ at every
horizontal jump time $\tau_k\le t$.

Choose $\alpha>0$ so small that
$\alpha\sup_{\lambda\in I}\mathbb ET_1<1/2$. Put
$M(\eta)=\sup_{\lambda\in I}\mathbb E^\lambda e^{\eta T_1}$.
The bound $\sup_{\lambda\in I}\mathbb E^\lambda e^{\eta_0T_1}<\infty$
from Lemma~\ref{lem:tangential-regeneration} gives
$\log M(\eta)\le\eta\sup_{\lambda\in I}\mathbb E^\lambda T_1+C\eta^2$
for small
$\eta>0$. Choose $\eta$ so that $\alpha\log M(\eta)<3\eta/4$.
Exponential Markov inequality, using independence of the cycle lengths,
then gives
\begin{equation}
 \label{eq:s9tan-cycle-concentration}
 \mathbb P^\lambda\{\tau_{\lfloor\alpha t\rfloor}>t\}
 \le e^{-\eta t}M(\eta)^{\lfloor\alpha t\rfloor}
 \le e^{-\eta t/4},\qquad t\ge0.
\end{equation}
The estimate \eqref{eq:s9tan-cycle-concentration} is the
exponential Markov inequality applied to the independent variables
$T_1,\ldots,T_{\lfloor\alpha t\rfloor}$; see
\cite[Section~2.7]{Durrett2019}.
Put $n=\lfloor\alpha t\rfloor$. Since $Y_s$ increases between
horizontal jumps, survival throughout $[0,t]$, together with
$\tau_n\le t$, implies positivity of $y+\sum_{i=1}^kZ_i$ for
$1\le k\le n$. Lemma~\ref{lem:centered-cycle-exit} and
\eqref{eq:s9tan-cycle-concentration} therefore give
\[
 \mathbb P^\lambda_{(j,0)}\{Y_s>0\text{ for }0\le s\le t\}
 \le C\frac{1+y}{\sqrt{1+n}}+e^{-\eta t/4}
 \le C'\frac{1+y}{\sqrt{1+t}}.
\]
Increasing $C'$ covers bounded $t$.

For the entrance estimate, decrease the positive constant $\alpha$
used in \eqref{eq:s9tan-cycle-concentration} so that
$\alpha\sup_{\lambda\in I}\mathbb E^\lambda T_1<1/6$ and put
$n=\lfloor\alpha t\rfloor$. For sufficiently small $\eta>0$,
$\alpha\log M(\eta)<\eta/4$, and hence
\[
 \mathbb P^\lambda\{\tau_n>t/3\}
 \le e^{-\eta t/3}M(\eta)^n\le e^{-\eta t/12}.
\]
For $t\ge2/\alpha$, $\alpha t/2\le n\le\alpha t$. Thus the range
$y_0\le y\le A\sqrt t$ is contained in
$y_0\le y\le A\sqrt{2/\alpha}\sqrt n$.
Apply Lemma~\ref{lem:centered-cycle-exit} with this latter upper
constant, obtaining a successful entrance event $E_0$ at an index
$J\le n/3$. Set
$T_{\mathrm{stop}}=\tau_J\wedge(t/3)$ and $E=E_0\cap\{\tau_J\le t/3\}$.
The first $k$ cycle pairs are observable by time $\tau_k$.
Since $J$ is a stopping time for $(\mathcal H_k)$, $\tau_J$ and
$T_{\mathrm{stop}}=\tau_J\wedge(t/3)$ are stopping times for the natural filtration
of $(N_t,R_t)$. Since the successful entrance event satisfies $E_0\in\mathcal H_J$,
it is measurable at $\tau_J$. Restricting to $\{\tau_J\le t/3\}$ gives
$E\in\mathcal F_{T_{\mathrm{stop}}}$, as required in the statement. Since $\tau_J\le\tau_n$,
\[
 \mathbb P^\lambda(E)
 \ge c_A\frac{y}{\sqrt n}-e^{-\eta t/12}
 \ge c'_A\frac{y}{\sqrt t}
\]
for all sufficiently large $t$, uniformly for $y\ge y_0$.
On $E$, $T_{\mathrm{stop}}=\tau_J$, $R_{T_{\mathrm{stop}}}=0$, and Lemma~\ref{lem:centered-cycle-exit}
gives
\[
 \delta\sqrt{\alpha/2}\sqrt t
 \le y+\sum_{i=1}^JZ_i=Y_{T_{\mathrm{stop}}}
 \le B_{A\sqrt{2/\alpha}}\sqrt\alpha\sqrt t.
\] On $E$, the discrete entrance event gives $Y_{\tau_k}>0$ for
$0\le k\le J$. Since $Y_s$ increases at rate $c>0$ between these
jumps, $Y_s>0$ for every $0\le s\le T_{\mathrm{stop}}$.

Write $(\widehat N,\widehat R)$ for the reversed chain and
$\widehat{\mathbb P}^\lambda_{(j,0)}$ for its law from $(j,0)$.
For assertions about its cycle pairs alone, put
$\widehat{\mathbb P}^\lambda=\widehat{\mathbb P}^\lambda_{(0,0)}$
and let $\widehat{\mathbb E}^\lambda$ denote expectation under this law.
Put $\mathfrak m(n,r)=\pi_r^\lambda$ for
$(n,r)\in\Z\times\N_0$. The transverse balance
identities in Lemma~\ref{lem:tangential-regeneration}, together with
cancellation of the two incoming and outgoing horizontal rates at
each $(n,0)$, give $\mathfrak m\mathcal Q_{\parallel}^\lambda=0$.
Thus counting measure in $n$ times $\pi^\lambda$ is invariant.
Write $q_\lambda(x,x')$ for the transition rate from $x$ to $x'$
of $\mathcal Q_{\parallel}^\lambda$. The rates in
Lemma~\ref{lem:tangential-regeneration} are uniformly bounded on $I$,
so we may take
\[
 \Gamma=1+\sup_{\lambda\in I}\sup_{x\in\Z\times\N_0}
                    \sum_{x'\ne x}q_\lambda(x,x')<\infty.
\]
The nonnegative matrix $I+\mathcal Q_{\parallel}^\lambda/\Gamma$
satisfies $\mathfrak m(I+\mathcal Q_{\parallel}^\lambda/\Gamma)
=\mathfrak m$. Tonelli's theorem applies to its nonnegative entries,
so the matrix powers and their Poisson sum also preserve
$\mathfrak m$, that is,
$\mathfrak m e^{t\mathcal Q_{\parallel}^\lambda}=\mathfrak m$
for $t\ge0$.
The reversed rates are
\[
 \widehat q_\lambda(x,x')
 =\frac{\mathfrak m(x')}{\mathfrak m(x)}q_\lambda(x',x),
 \quad\text{ for }x,x'\in\Z\times\N_0,\ x\ne x'.
\]
At a horizontal jump the two endpoints have the same weight, so
reversal interchanges $De^\lambda$ and $De^{-\lambda}$.
Transverse reversal leaves the rates unchanged by detailed balance.
For later use, consider a path with states $x_0,\ldots,x_k$ and
ordered jump times $0<s_1<\cdots<s_k<t$. By the definition of the
reversed rates,
\begin{equation}\label{eq:tan-finite-path-reversal}
 \mathfrak m(x_0)\prod_{i=1}^kq_\lambda(x_{i-1},x_i)
 =\mathfrak m(x_k)\prod_{i=1}^k
                         \widehat q_\lambda(x_i,x_{i-1}).
\end{equation}
The total rates agree under reversal: the transverse rates are
unchanged and the sum of the two horizontal rates is
$2D\cosh\lambda$ for both chains. Hence the holding factors agree
as well. After reversing the ordered jump times,
\eqref{eq:tan-finite-path-reversal} gives the path-reversal identity
for every event described by finitely many jumps. Only finitely many jumps occur on $[0,t]$ almost surely because
the total rates are bounded. Integrating the finite-path densities
over jump times and summing over the jump count and state sequences
therefore gives the same reversal identity for measurable path
events on $[0,t]$.

When both endpoints have transverse state zero, the endpoint factor
in \eqref{eq:tan-finite-path-reversal} is
$\mathfrak m(x_0)/\mathfrak m(x_k)=1$.
For the strict boundary events used below, reversal preserves
positivity on every open holding interval. A discrepancy at a jump
could occur only if that jump coincided with a contact of an
integer site and the deterministic boundary. For each fixed
integer, the straight boundary has at most one such contact; the
logarithmic boundary has at most two. Jump times have continuous
laws, so this exceptional event has probability zero after taking
the countable union over the sites. Thus the same reversal identity
applies to the all-time strict survival events.

Write $(\widehat T_i,\widehat\xi_i)$ for the waiting times
and displacements between successive longitudinal jumps of
$(\widehat N,\widehat R)$,
$\widehat\tau_k=\sum_{i=1}^k\widehat T_i$, and
$\widehat Z_i=\widehat\xi_i-c\widehat T_i$. Their laws satisfy
$\widehat T_i\overset{\mathrm d}=T_i$ and
$\widehat Z_i\overset{\mathrm d}=-Z_i$; the waiting time and mark remain
independent. Consequently
\[
 \widehat Y_{\widehat\tau_k}
 =y+\sum_{i=1}^k\widehat Z_i,\qquad k\in\N_0.
\]
For $n=\lfloor\alpha t\rfloor$, continuous survival and
$\widehat\tau_n\le t$ imply positivity of all these first $n$ values.
Since $\widehat Z_1\overset{\mathrm d}=-Z_1$, the bounds just verified
for the law of $Z_1$ give, for every $a\ge0$,
\[
 \begin{aligned}
 \widehat{\mathbb E}^\lambda[(\widehat Z_1-a)^2;\widehat Z_1>a]
 &\le C_{\mathrm o}\widehat{\mathbb P}^\lambda(\widehat Z_1>a),\\
 \widehat{\mathbb E}^\lambda[(-\widehat Z_1-a)^2;-\widehat Z_1>a]
 &\le C_{\mathrm o}\widehat{\mathbb P}^\lambda(-\widehat Z_1>a).
 \end{aligned}
\]
Thus Lemma~\ref{lem:centered-cycle-exit} and
\eqref{eq:s9tan-cycle-concentration}, using the identical laws of
$\widehat T_i$ and $T_i$, give
\[
 \widehat{\mathbb P}^\lambda_{(j,0)}
 \{\widehat Y_s>0\text{ for }0\le s\le t\}
 \le C\frac{1+y}{\sqrt{1+n}}+e^{-\eta t/4}
 \le C'\frac{1+y}{\sqrt{1+t}}.
\]
For the entrance lower bound, increase $y_0$ so that, whenever
$y\ge y_0$, one has $y-1\ge y/2$ and $y-1$ exceeds the threshold
in Lemma~\ref{lem:centered-cycle-exit}. Apply that lemma to
$y-1+\sum_{i=1}^k\widehat Z_i$.
Let $\widehat J\le n/3$ be the stopping index provided by that lemma.
On its successful entrance event $\widehat E_0$,
$y-1+\sum_{i=1}^k\widehat Z_i>0$ for $0\le k\le\widehat J$.
Hence $\widehat Y_{\widehat\tau_k}>1$ for these indices. Their pre-jump heights satisfy
\[
 \widehat Y_{\widehat\tau_k-}
 =\widehat Y_{\widehat\tau_k}-\widehat\xi_k>0,
 \qquad 1\le k\le\widehat J,\quad\text{on }\widehat E_0,
\]
since $\widehat\xi_k\le1$. The height decreases between jumps, so
positivity of these pre-jump heights ensures survival on every entire
waiting interval, including the first.
Set $\widehat{T}_{\mathrm{stop}}=\widehat\tau_{\widehat J}\wedge(t/3)$ and
$\widehat E=\widehat E_0\cap\{\widehat\tau_{\widehat J}\le t/3\}$.
With
$\widehat{\mathcal H}_k=\sigma(\widehat T_i,\widehat\xi_i:1\le i\le k)$,
$\widehat J$ is a stopping index and $\widehat E_0$ is measurable
at $\widehat J$. Hence $\widehat\tau_{\widehat J}$ and
$\widehat{T}_{\mathrm{stop}}$ are stopping times for
$(\widehat{\mathcal F}_s)_{s\ge0}$, and
$\widehat E\in\widehat{\mathcal F}_{\widehat{T}_{\mathrm{stop}}}$.
The identical waiting-time law gives
$\widehat{\mathbb P}^\lambda(\widehat\tau_n>t/3)
\le e^{-\eta t/12}$. Hence
\[
 \widehat{\mathbb P}^\lambda(\widehat E)
 \ge c_A\frac{y-1}{\sqrt n}-e^{-\eta t/12}
 \ge c'_A\frac{y}{\sqrt t}
\]
for large $t$. On $\widehat E$, $\widehat R_{\widehat{T}_{\mathrm{stop}}}=0$ and
\[
 \widehat Y_{\widehat{T}_{\mathrm{stop}}}
 =1+\left(y-1+\sum_{i=1}^{\widehat J}\widehat Z_i\right).
\]
Since $\alpha t/2\le n\le\alpha t$, the discrete entrance bounds
give a lower bound $\delta\sqrt{\alpha/2}\sqrt t$ and, for $t\ge1$,
an upper bound
$(1+B_{A\sqrt{2/\alpha}}\sqrt\alpha)\sqrt t$ for this height. Choose $\delta$ and $\kappa_A$ as the minima of the forward and
reversed lower constants, and choose $B_A$ and $T_A$ as the maxima
of their upper-height constants and time thresholds, respectively. After increasing $y_0$ as above,
these choices give the same $y_0,\delta,B_A,\kappa_A,T_A$ for the
forward and reversed chains, uniformly for $\lambda\in I$.\qedhere

\end{proof}

For $t\ge0$, $y>0$, $j,l\in\Z$, and $\lambda\in I$, define
\[
 K_t^{y,\lambda}(j,l)
 =
 \mathbb P^\lambda_{(j,0)}
 \left\{
 (N_t,R_t)=(l,0),\
 y+N_s-j+c_\lambda s>0
 \text{ for }0\le s\le t
 \right\}.
\]

\begin{lemma}\label{lem:tangential-killed}

Use the chain $(N_t,R_t)$ from Lemma~\ref{lem:tangential-regeneration},
its reversal from Lemma~\ref{lem:tangential-survival-entrance}, and the
kernel $K_t^{y,\lambda}(j,l)$ defined immediately above. Recall that
$c_\lambda=\Lambda_{\parallel}'(\lambda)$. There is $C<\infty$,
uniform for $\lambda\in I$, such that, with
\[
 z=y+l-j+c_\lambda t,
\]
\[
 K_t^{y,\lambda}(j,l)
 \le
 C\frac{(1+y)(1+z)}{(1+t)^{3/2}}
\]
for every $t\ge0$, $j,l\in\Z$, $\lambda\in I$ and $y>0$
for which $z=y+l-j+c_\lambda t>0$. There is $y_0\ge1$ such that,
for every $A>0$, there
are $c_A>0$ and $T_A<\infty$ satisfying
\[
 K_t^{y,\lambda}(j,l)
 \ge
 c_A\frac{yz}{t^{3/2}}
\]
for $j,l\in\Z$ and $\lambda\in I$, whenever
\[
 t\ge T_A,
 \qquad
 y_0\le y,z\le A\sqrt t.
\]

\end{lemma}

\begin{proof}

It suffices first to consider $t\ge1$; the bound $K_t^{y,\lambda}(j,l)\le1$
covers $0\le t\le1$ after increasing $C$.
To prove the upper bound, start the forward chain at $(j,0)$.
Let $T_{\mathrm{stop}}$ be its first horizontal jump strictly after $t/4$, and let
$\widehat{T}_{\mathrm{stop}}$ be the first horizontal jump strictly after $t/4$
for the reversed chain started at $(l,0)$. Stopped paths include
their cut jump and are evaluated in its post-jump state.
We first estimate the two exceptional events
$\{T_{\mathrm{stop}}>t/3\}$ and $\{\widehat{T}_{\mathrm{stop}}>t/3\}$. Each requires an interval
of length $t/12$ without a horizontal jump.

To bound $\mathbb P^\lambda_{(j,0)}(T_{\mathrm{stop}}-t/4>t/12)$
and $\widehat{\mathbb P}^\lambda_{(l,0)}(\widehat{T}_{\mathrm{stop}}-t/4>t/12)$,
let $L_R=q_R$ be the
generator of the transverse process $R_t$ and use the transverse generator
$K^\lambda$ killed at the first horizontal jump, from
Lemma~\ref{lem:tangential-regeneration}. Thus, for
$r\in\N_0$ and $F\in\ell^2(\pi^\lambda)$,
$(K^\lambda F)_r=(L_RF)_r
-2D\cosh\lambda\,\mathbf1_{\{r=0\}}F_0$.
Reversibility and boundedness of the rates make $L_R$ and $K^\lambda$
bounded self-adjoint operators on $\ell^2(\pi^\lambda)$. The
transverse generator satisfies $L_R\le0$, and the estimate
$\sup\spec(K^\lambda)\le-g<0$ in Lemma~\ref{lem:tangential-regeneration}
gives $\langle F,K^\lambda F\rangle_{\pi^\lambda}
\le-\rho\|F\|_{\ell^2(\pi^\lambda)}^2$ for every
$F\in\ell^2(\pi^\lambda)$, with a common $\rho>0$ on $I$. For every $v,u\ge0$,
\begin{equation}\label{eq:tan-no-horizontal-jump}
 \begin{aligned}
 \mathbb P_{(j,0)}^\lambda
 \{\text{no horizontal jump on }[v,v+u]\}
 &=\left\langle e^{vL_R}\frac{\delta_0}{\pi_0^\lambda},
                   e^{uK^\lambda}\mathbf1\right\rangle_{\pi^\lambda}\\
 &\le (\pi_0^\lambda)^{-1/2}e^{-\rho u}.
 \end{aligned}
\end{equation}
Horizontal reversal changes only the sign of the jump, so its
transverse generator is $L_R$ and its total horizontal killing rate
remains $2D\cosh\lambda$. Its killed generator is therefore also
$K^\lambda$, and the displayed inner-product formula gives
\[
 \widehat{\mathbb P}_{(l,0)}^\lambda
 \{\text{no horizontal jump on }[v,v+u]\}
 \le(\pi_0^\lambda)^{-1/2}e^{-\rho u},\qquad v,u\ge0.
\]
The exponential bound for $\psi^\lambda$ in the proof of
Lemma~\ref{lem:tangential-regeneration} gives
$\pi_0^\lambda=1/\sum_{r=0}^\infty(\psi_r^\lambda)^2\ge C^{-1}$ uniformly on $I$.

On paths with $T_{\mathrm{stop}}-t/4\le t/12$ and
$\widehat{T}_{\mathrm{stop}}-t/4\le t/12$,
\[
 T_{\mathrm{stop}}\le t/3,
 \qquad
 \widehat{T}_{\mathrm{stop}}\le t/3,
\]
and both cut points have transverse state zero. Define the
subprobability measures $\alpha_y(ds,m)$ and
$\widehat\alpha_z(dr,n)$ for the surviving forward and reversed
segments, each stopped just after its cut jump and restricted to a
stopping time at most $t/3$. With
$Y_u=y+N_u-j+c_\lambda u$ and
$\widehat Y_u=z+\widehat N_u-l-c_\lambda u$, these measures are
\[
 \begin{aligned}
 \alpha_y(ds,m)
 &=\mathbb P^\lambda_{(j,0)}
   \{T_{\mathrm{stop}}\in ds,\ N_{T_{\mathrm{stop}}}=m,\ T_{\mathrm{stop}}\le t/3,\ Y_u>0\text{ for }0\le u\le T_{\mathrm{stop}}\},\\
 \widehat\alpha_z(dr,n)
 &=\widehat{\mathbb P}^\lambda_{(l,0)}
   \{\widehat{T}_{\mathrm{stop}}\in dr,\ \widehat N_{\widehat{T}_{\mathrm{stop}}}=n,
     \ \widehat{T}_{\mathrm{stop}}\le t/3,\ \widehat Y_u>0\text{ for }0\le u\le\widehat{T}_{\mathrm{stop}}\}.
 \end{aligned}
\]
By \eqref{eq:tan-finite-path-reversal}, the terminal segment of a
forward path ending at $(l,0)$ can be read as a path of the reversed
chain started at $(l,0)$.
The first forward cut and the first reversed cut are determined by
their respective outer segments. For cut durations $s,r\in[t/4,t/3]$, one has $s+r\le2t/3$, so
these initial and terminal segments are disjoint. Expanding the path densities and integrating over the
two cut times therefore bounds the surviving contribution by
\[
 \sum_{m,n\in\Z}\int_{[t/4,t/3]}\!\int_{[t/4,t/3]}
 \alpha_y(ds,m)\,p_{t-s-r}^\lambda(n-m)\,
 \widehat\alpha_z(dr,n).
\]
In forward time, the first cut is a post-jump state and the second is
a pre-jump state. Each outer measure includes its cut jump once.
Both cut states have transverse coordinate zero, so the reversal weight
ratio equals one. Discarding middle absorption gives the displayed bound.
Since $T_{\mathrm{stop}},\widehat{T}_{\mathrm{stop}}\ge t/4$, Lemma~\ref{lem:tangential-survival-entrance}
bounds the total masses of $\alpha_y$ and $\widehat\alpha_z$
by $C(1+y)t^{-1/2}$ and $C(1+z)t^{-1/2}$, respectively,
because each surviving outer segment survives at least until $t/4$.
Since $t-s-r\ge t/3$, Lemma~\ref{lem:tangential-free-kernel}
bounds the middle kernel by $Ct^{-1/2}$. Substituting this kernel bound
and the mass bounds for $\alpha_y$ and $\widehat\alpha_z$ into the
displayed integral gives $Ct^{-3/2}(1+y)(1+z)$. By \eqref{eq:tan-no-horizontal-jump} and its reversed version,
the contribution from $T_{\mathrm{stop}}>t/3$ or $\widehat{T}_{\mathrm{stop}}>t/3$ is at most
$Ce^{-\rho t/12}\le Ct^{-3/2}$. Since $(1+y)(1+z)\ge1$,
this contribution is also at most
$Ct^{-3/2}(1+y)(1+z)$.

We next prove an interior lower estimate. For a parameter $H\ge1$
to be fixed below, first require
\[
 y,z\ge H\sqrt t,\qquad |z-y|\le\sqrt t.
\]
By Lemma~\ref{lem:tangential-free-kernel},
$p_t^\lambda(l-j)\ge b_1t^{-1/2}$ for all sufficiently large $t$.
The constant $b_1>0$ and the threshold for $t$ are independent of
$H$: the centered displacement is $z-y$, so the lower bound in
Lemma~\ref{lem:tangential-free-kernel} applies with $A=1$.
The number of horizontal jumps is dominated by a Poisson process
whose rate is uniformly bounded on $I$. Recall from
Lemma~\ref{lem:tangential-regeneration} that the cycle pairs
$(T_i,\xi_i)$ are independent and identically distributed. Set
$Z_i=\xi_i+c_\lambda T_i$ for $i\in\N$; these increments are
centered and have uniformly bounded variances on $I$. For $L>0$ and
$m=\lceil Lt\rceil$, Kolmogorov's maximal inequality
\cite[Theorem~2.5.5]{Durrett2019} gives
\[
 \mathbb P^\lambda\left\{
 \max_{1\le k\le m}\left|\sum_{i=1}^kZ_i\right|\ge H\sqrt t
 \right\}
 \le \frac{m\operatorname{Var}(Z_1)}{H^2t}
 \le \frac{C_L}{H^2}.
\]
Choose $L$ larger than twice the common upper bound for the horizontal
jump rate. The exponential Markov bound for the dominating Poisson
variable bounds the event of more than $m$ jumps by $Ce^{-\rho t}$.
Let $A_+$ be the event that $Y_s=y+N_s-j+c_\lambda s$ reaches
zero or a negative value during $[0,t/2]$. Since $Y_s$ increases
between horizontal jumps, $A_+$ requires
$\sum_{i=1}^kZ_i\le-y\le-H\sqrt t$ at a horizontal jump in that
interval. The maximal inequality and the Poisson tail thus give
$\mathbb P^\lambda_{(j,0)}(A_+)\le C/H^2+Ce^{-\rho t}$.
For the reversed chain started at $(l,0)$, set
$\widehat Z_i=\widehat\xi_i-c_\lambda\widehat T_i$ for $i\in\N$,
where $(\widehat T_i,\widehat\xi_i)$ are its cycle lengths and marks.
For its height $\widehat Y_s=z+\widehat N_s-l-c_\lambda s$, let
$A_-$ be the event that this height reaches a nonpositive value in
$[0,t/2]$. A crossing occurs either at a jump or
before the next jump. In the latter case the height keeps decreasing
until that jump, so its post-jump value is at most one. Hence in
both cases some cycle ending no later than the first horizontal jump
after $t/2$ satisfies
\[
 \sum_{i=1}^k\widehat Z_i\le1-z\le1-H\sqrt t.
\]
On the event that at most $m$ horizontal jumps occur before $t/2$,
at most $m+1$ cycles enter this estimate. Apply Kolmogorov's
inequality to the first $m+1$ centered increments
$\widehat Z_i\overset{\mathrm d}=-Z_i$ and use the Poisson tail bound
for the complementary event. Since $(m+1)/t$ is bounded for $t\ge1$,
we obtain
$\widehat{\mathbb P}^\lambda_{(l,0)}(A_-)
\le C/H^2+Ce^{-\rho t}$ for $H\sqrt t\ge2$.

Let $T_{\mathrm{stop},*}$ be the first horizontal jump after $t/2$. By \eqref{eq:tan-no-horizontal-jump}, applied with $v=t/2$ and $u=t/6$, we have
$\mathbb P^\lambda_{(j,0)}(T_{\mathrm{stop},*}>2t/3)\le Ce^{-\rho t}$.
On $T_{\mathrm{stop},*}\le2t/3$, the strong Markov property, $R_{T_{\mathrm{stop},*}}=0$, and the
free-kernel upper bound yield
\[
 \mathbb P^\lambda_{(j,0)}
 \{A_+,(N_t,R_t)=(l,0)\}
 \le Ct^{-1/2}\mathbb P^\lambda_{(j,0)}(A_+)+Ce^{-\rho t}.
\]
Under horizontal reflection the reversed free kernel is
$\widehat p_u^\lambda(k)=p_u^\lambda(-k)$, so it has the identical
$Cu^{-1/2}$ bound. Let $\widehat{T}_{\mathrm{stop},*}$ be the first reversed horizontal jump after $t/2$.
The no-jump estimate on $[t/2,2t/3]$ gives, after decreasing
$\rho>0$ if necessary,
\[
 \widehat{\mathbb P}^\lambda_{(l,0)}
       \{\widehat{T}_{\mathrm{stop},*}>2t/3\}\le Ce^{-\rho t}.
\]
On $\{\widehat{T}_{\mathrm{stop},*}\le2t/3\}$, the cut has
$\widehat R_{\widehat{T}_{\mathrm{stop},*}}=0$ and the remaining duration satisfies
$t-\widehat{T}_{\mathrm{stop},*}\ge t/3$. The event $A_-$ is determined by time $t/2$.
Thus the strong Markov property at $\widehat{T}_{\mathrm{stop},*}$ and the bound
$\widehat p_u^\lambda(k)\le Cu^{-1/2}$ give
\[
 \begin{aligned}
 &\widehat{\mathbb P}^\lambda_{(l,0)}
   \{A_-,(\widehat N_t,\widehat R_t)=(j,0)\}\\
 &\quad\le
 Ct^{-1/2}\widehat{\mathbb P}^\lambda_{(l,0)}(A_-)+Ce^{-\rho t}\\
 &\quad\le
 Ct^{-1/2}\bigl(C/H^2+Ce^{-\rho t}\bigr)+Ce^{-\rho t}\\
 &\quad\le CH^{-2}t^{-1/2}+Ce^{-\rho t}.
 \end{aligned}
\]
The path-reversal identity \eqref{eq:tan-finite-path-reversal}
from $(j,0)$ to $(l,0)$ has invariant-measure
weight ratio $\pi_0^\lambda/\pi_0^\lambda=1$ and sends crossings in
the last half to $A_-$. Hence the displayed estimate bounds their
contribution to the forward endpoint probability.
A path counted by $p_t^\lambda(l-j)-K_t^{y,\lambda}(j,l)$ has a
crossing in one of the two half-intervals. Consequently
\[
 0\le p_t^\lambda(l-j)-K_t^{y,\lambda}(j,l)
 \le CH^{-2}t^{-1/2}+Ce^{-\rho t}.
\]
Choose $H$ using $C/H^2\le b_1/4$, and then take $t$ large enough to
absorb the exponential term. For all sufficiently large $t$, uniformly for $j,l\in\Z$ and
$\lambda\in I$, this proves
\begin{equation}\label{eq:tan-killed-interior}
 K_t^{y,\lambda}(j,l)\ge (b_1/2)t^{-1/2}
 \quad\text{if }y,z\ge H\sqrt t,\quad |z-y|\le\sqrt t.
\end{equation}

Now fix $0<a\le b<\infty$ and suppose
$a\sqrt t\le y,z\le b\sqrt t$. Put $\eta=1/8$ and choose a fixed
integer $N\ge1$ such that
\[
 a\sqrt N-\eta\ge H,\qquad (b-a)/\sqrt N+2\eta\le1.
\]
At the deterministic times $k\tau$, where $\tau=t/N$, restrict the
transverse state to zero and the longitudinal state to
\[
 I_k=\left\{m\in\Z:
 \left|m-\left(j-c_\lambda k\tau+\frac{k}{N}(z-y)\right)\right|
 \le\eta\sqrt\tau\right\},\quad\text{ for }1\le k<N.
\]
For $1\le k<N$ and $m_k\in I_k$, write
\[
 e_k=m_k-j+c_\lambda k\tau-\frac{k}{N}(z-y),
 \qquad e_0=e_N=0.
\]
The height at the $k$th cut is
$(1-k/N)y+(k/N)z+e_k$, with $|e_k|\le\eta\sqrt\tau$.
It is at least $(a\sqrt N-\eta)\sqrt\tau\ge H\sqrt\tau$,
and consecutive heights differ by at most
$((b-a)/\sqrt N+2\eta)\sqrt\tau\le\sqrt\tau$.
By \eqref{eq:tan-killed-interior}, each killed kernel of duration
$\tau$ between consecutive restricted cut states has lower bound
$(b_1/2)\tau^{-1/2}$ once $\tau$ exceeds its uniform time threshold, while each $I_k$ contains at least
$\eta\sqrt\tau$ integers for large $t$. The Markov property at these
deterministic times gives
\[
\begin{gathered}
 K_t^{y,\lambda}(j,l)
 \ge (\eta\sqrt\tau)^{N-1}
       \bigl((b_1/2)\tau^{-1/2}\bigr)^N
 \ge c_{a,b}t^{-1/2}
 \quad\text{ for all sufficiently large }t,
 \\\text{ with }a\sqrt t\le y,z\le b\sqrt t.
\end{gathered}
\]

Choose $y_0$ as in Lemma~\ref{lem:tangential-survival-entrance}.
For $y_0\le y,z\le A\sqrt t$, apply that lemma from heights
$y$ and $z$ in the forward and reversed directions. Denote the resulting entrance events by $E_y,\widehat E_z$ and their
stopping times by $T_{\mathrm{stop}},\widehat{T}_{\mathrm{stop}}\le t/3$, respectively. Define the corresponding subprobability
measures by
\[
 \begin{aligned}
 F_y(ds,m)&=\mathbb P^\lambda_{(j,0)}
       \{T_{\mathrm{stop}}\in ds,N_{T_{\mathrm{stop}}}=m,E_y\},\\
 \widehat F_z(du,n)&=\widehat{\mathbb P}^\lambda_{(l,0)}
       \{\widehat{T}_{\mathrm{stop}}\in du,
         \widehat N_{\widehat{T}_{\mathrm{stop}}}=n,\widehat E_z\}.
 \end{aligned}
\]
For all sufficiently large $t$, their masses are at least
$c_Ay/\sqrt t$ and $c_Az/\sqrt t$. On the supports of $F_y(ds,m)$
and $\widehat F_z(du,n)$, the transverse coordinate is zero and the
heights
\[
 y_s=y+m-j+c_\lambda s,
 \qquad z_u=z+n-l-c_\lambda u
\]
belong to $[\delta\sqrt t,B_A\sqrt t]$. The middle interval has
length $T_{s,u}=t-s-u\in[t/3,t]$, and
\[
 y_s+n-m+c_\lambda T_{s,u}=z_u.
\]
Thus $\delta\sqrt{T_{s,u}}\le y_s,z_u
\le\sqrt3B_A\sqrt{T_{s,u}}$. Apply the lower bound just proved
for $a\sqrt T\le y,z\le b\sqrt T$ with
$a=\delta$, $b=\sqrt3B_A$ and $T=T_{s,u}$. It gives
$K_{T_{s,u}}^{y_s,\lambda}(m,n)\ge c_A t^{-1/2}$ uniformly over these
outer segments. The events $E_y\in\mathcal F_{T_{\mathrm{stop}}}$ and
$\widehat E_z\in\widehat{\mathcal F}_{\widehat{T}_{\mathrm{stop}}}$ are determined
at their respective stopping times; their path segments are disjoint
because $s+u\le2t/3$.
Consequently the
forward prefix determines $(T_{\mathrm{stop}},E_y)$ and the reversed suffix
determines $(\widehat{T}_{\mathrm{stop}},\widehat E_z)$, even after an arbitrary
surviving middle path has been inserted. Every concatenated path
is therefore counted at most once. By
\eqref{eq:tan-finite-path-reversal}, reversal of the terminal
segment has weight one because both transverse endpoints equal
zero. Keeping the survival condition on the interval of length $t-s-u$ gives
{\postdisplaypenalty=10000
\begin{equation}\label{eq:tan-killed-stopping-convolution}
 \begin{aligned}
 K_t^{y,\lambda}(j,l)
 &\ge\sum_{m,n\in\Z}\int_{[0,t/3]}\!\int_{[0,t/3]}
 F_y(ds,m)K_{t-s-u}^{y_s,\lambda}(m,n)\widehat F_z(du,n)\\
 &\ge c_A t^{-1/2}\frac{y}{\sqrt t}\frac{z}{\sqrt t},
 \end{aligned}
\end{equation}
which proves the claimed lower bound.\qedhere
}

\end{proof}

For $t\ge0$, $a\ge a_*$, $y>0$, $j,l\in\Z$, and $\lambda\in I$, put
\[
 b_{t,a}(s)
 =
 \kappa_{\parallel}\log\frac{1+a+t}{1+a+t-s},
 \qquad
 0\le s\le t,
\]
and define
\[
 \mathcal K_{t,a}^{y,\lambda}(j,l)
 =
 \mathbb P^\lambda_{(j,0)}
 \left\{
 (N_t,R_t)=(l,0),\
 y+N_s-j+c_\lambda s>b_{t,a}(s)
 \text{ for }0\le s\le t
 \right\}.
\]
On the event $(N_t,R_t)=(l,0)$, the terminal height
$y+N_t-j+c_\lambda t-b_{t,a}(t)$ equals
\[
 z
 =
 y+l-j+c_\lambda t-b_{t,a}(t).
\]

To estimate survival above $b_{t,a}$, we use the independent
identically distributed cycle pairs of the unrestricted chain in
Lemma~\ref{lem:tangential-regeneration}, together with their uniform
exponential moment.

\begin{lemma}\label{lem:tangential-backward-log}

After shrinking $I$ if necessary, there is $C<\infty$ such that
\[
 \mathcal K_{t,a}^{y,\lambda}(j,l)
 \le
 C\frac{(1+y)(1+z)}{(1+t)^{3/2}}
\]
for every $t\ge0$, $j,l\in\Z$, $\lambda\in I$ and $a\ge a_*$
with $y,z>0$.
There is $Y>0$ such that, for every $A>0$, there are $T,c_A>0$,
independent of $a\ge a_*$ and $\lambda\in I$, such that
\[
 \mathcal K_{t,a}^{y,\lambda}(j,l)
 \ge
 c_A\frac{yz}{t^{3/2}}
\]
whenever $j,l\in\Z$, $\lambda\in I$, $a\ge a_*$ and
\[
 t\ge T,
 \qquad
 Y\le y,z\le A\sqrt t.
\]
Put $Y_s=y+N_s-j+c_\lambda s-b_{t,a}(s)$ for $0\le s\le t$.
If $y,z$ range in a fixed compact subset of $(Y,\infty)$, there is
a constant $C$, depending also on that subset but independent of $t$,
$a\ge a_*$ and $\lambda\in I$, such that, whenever the conditioning
event has positive probability,
\[
 \mathbb E^\lambda_{(j,0)}
 \left[
 \int_0^t
 (1+Y_s)e^{-\lambda_{\parallel}Y_s}\,ds
 \,\middle|\,
 (N_t,R_t)=(l,0),\
 Y_r>0\text{ for }0\le r\le t
 \right]
 \le C.
\]

\end{lemma}

\begin{proof}
We first prove the survival bound for the reversed process above the
decreasing logarithmic boundary.

\medskip
\noindent\textit{Step 1: a deterministic logarithmic boundary for the
cycle walk.}
Fix $\lambda\in I$. In the cycle calculations below, $\mathbb P$
and $\mathbb E$ denote probability and expectation for the reversed
chain started at transverse state zero; horizontal translation makes
the initial longitudinal coordinate irrelevant. Let
\[
 (\widehat T_k,\widehat\xi_k)_{k\ge1}
\]
be the successive cycle lengths and horizontal marks.  The waiting-time
law is the same as in Lemma~\ref{lem:tangential-regeneration}, while
the two horizontal probabilities are interchanged.  Put
\[
 X_k=\widehat\xi_k-c_\lambda\widehat T_k,\qquad k\in\N.
\]
By Lemma~\ref{lem:tangential-regeneration},
$X_k\overset{\mathrm d}=-(\xi_k+c_\lambda T_k)$. Thus the $X_k$ are i.i.d.,
\[
 \mathbb EX_k=0,
\]
and their variances are bounded above and below by positive constants,
uniformly for $\lambda\in I$. There is $\eta_1>0$ such that
$\sup_{\lambda\in I}\mathbb E e^{\eta_1|X_1|}<\infty$.  The waiting time $\widehat T_k$ remains independent
of the mark $\widehat\xi_k$.

Let $K^\lambda$ be the transverse generator killed at the first
horizontal jump, acting on $\ell^2(\pi^\lambda)$. Under the unitary
map $F_r\mapsto\sqrt{\pi_r^\lambda}F_r$ from
$\ell^2(\pi^\lambda)$ to $\ell^2(\N_0)$, it
is unitarily equivalent to
$\mathcal J_{\mathcal G}-2D\langle\delta_0,\cdot\rangle\delta_0
-\Lambda_{\parallel}(\lambda)I$. Thus $\langle F,K^\lambda F\rangle_{\pi^\lambda}
\le-g\|F\|_{\ell^2(\pi^\lambda)}^2$ for
$F\in\ell^2(\pi^\lambda)$, where
$g=\min_{\lambda\in I}
(\Lambda_{\parallel}(\lambda)-\sigma_{\mathcal G})>0$.
Horizontal reversal preserves this killed generator because the total
horizontal rate remains $2D\cosh\lambda$.
Let $\varrho_\lambda$ be the spectral probability measure of
$-K^\lambda$ at $\delta_0/\sqrt{\pi_0^\lambda}$ in
$\ell^2(\pi^\lambda)$; its support lies in $[g,\infty)$.
For $u>0$, a first horizontal jump at time $u$ has density
\[
 2D\cosh\lambda\,(e^{uK^\lambda})(0,0)
 =2D\cosh\lambda\int_{[g,\infty)}e^{-\beta u}\,
                  \dd\varrho_\lambda(\beta).
\]
For $u\ge0$, integrating this density from $u$ to infinity gives
\[
 \widehat{\mathbb P}^\lambda_{(l,0)}\{\widehat T_1>u\}
 =\int_{[g,\infty)}e^{-\beta u}\,\dd\omega_\lambda(\beta),
 \qquad
 \dd\omega_\lambda(\beta)
 =\frac{2D\cosh\lambda}{\beta}\,\dd\varrho_\lambda(\beta).
\]
Since $\widehat T_1<\infty$ almost surely and
$\widehat T_1>0$ almost surely, setting $u=0$ in the tail formula
gives $\omega_\lambda([g,\infty))=1$. Set $\rho_0=g$.
For $u,v\ge0$, use
$e^{-\beta(u+v)}\le e^{-\rho_0v}e^{-\beta u}$
for $\beta\ge g$ and divide by
$\int_{[g,\infty)}e^{-\beta u}\,\dd\omega_\lambda(\beta)>0$.
We obtain
\[
 \widehat{\mathbb P}^\lambda_{(l,0)}\{\widehat T_1>u+v\mid\widehat T_1>u\}
 \le e^{-\rho_0v},
 \quad\text{ for }u,v\ge0.
\]
Integrating in $v$ yields
\[
 \sup_{\lambda\in I}\sup_{u\ge0}
 \widehat{\mathbb E}^\lambda_{(l,0)}[
 \widehat T_1-u\mid\widehat T_1>u]
 \le \rho_0^{-1}.
\]

Fix $\kappa\ge0$ and define
\[
 g_k=\kappa\log(1+k),
 \quad\text{ for }k\ge0.
\]
For $y>0$ and $n\in\N_0$ let
\[
 F_n(y)
 =
 \mathbb P\left\{
 y+\sum_{i=1}^kX_i+g_k>0
 \text{ for }1\le k\le n
 \right\}.
\]
We claim that
\begin{equation}
 \label{eq:cycle-log-survival}
 F_n(y)
 \le
 C_\kappa\frac{1+y}{\sqrt{1+n}},
 \quad\text{ for }n\ge0,
\end{equation}
uniformly for $\lambda\in I$.

Let
\[
 Y_k=y+\sum_{i=1}^kX_i+g_k,\quad k\in\N_0,
 \qquad
 \tau=\inf\{k\ge1:Y_k\le0\}.
\]
Use the cycle filtration
$\mathcal F_k=\sigma(\widehat T_i,\widehat\xi_i:1\le i\le k)$.
On $\{\tau>k-1\}$, condition on $\mathcal F_{k-1}$ and on the horizontal mark
$\widehat\xi_k$.  If
\[
 u=
 \frac{Y_{k-1}+g_k-g_{k-1}+\widehat\xi_k}
 {c_\lambda}
\]
is nonnegative, the exit event is
$\{\widehat T_k\ge u\}$ and
\[
 -Y_k=c_\lambda(\widehat T_k-u).
\]
Since $0<\min_{\lambda\in I}c_\lambda
\le c_\lambda\le\max_{\lambda\in I}c_\lambda<\infty$ and
$\widehat T_k$ is independent of $\mathcal F_{k-1}$ and
$\widehat\xi_k$, the bound
$\widehat{\mathbb E}^\lambda_{(l,0)}[\widehat T_1-u\mid\widehat T_1>u]
\le\rho_0^{-1}$ controls $-Y_k$ on $\{\tau=k\}$.  More precisely, on
the surviving past and when $u\ge0$,
\[
 \mathbb E[-Y_k;\tau=k\mid\mathcal F_{k-1},\widehat\xi_k]
 \le
 \frac{\max_{\lambda\in I}c_\lambda}{\rho_0}
 \mathbb P(\tau=k\mid\mathcal F_{k-1},\widehat\xi_k).
\]
If $u<0$, then necessarily $\widehat\xi_k=-1$ and
$|Y_{k-1}+g_k-g_{k-1}-1|\le1$ at a possible exit; hence
\[
 -Y_k\le c_\lambda\widehat T_k+1,
\]
and the conditional expectation of $c_\lambda\widehat T_k+1$ is
uniformly bounded. When $u<0$, exit at step $k$ is certain conditional on
$\{\tau>k-1\}$, $\mathcal F_{k-1}$ and the mark $\widehat\xi_k$.  In both cases the
conditional mean overshoot is therefore bounded by a constant times
the conditional exit probability.  Summing over the disjoint events
$\{\tau=k\}$ gives a constant $C_1$ such that
\[
 \mathbb E[-Y_\tau;\tau\le m]\le C_1
 \quad\text{ for } m\in\N.
\]
The process $\sum_{i=1}^kX_i$ is an integrable martingale for
$\mathcal F_k$, and $\tau\wedge m\le m$ is a bounded stopping time.
Moreover,
\[
 g_{\tau\wedge m}
 =\sum_{k=1}^m(g_k-g_{k-1})\mathbf1_{\{\tau\ge k\}},
 \qquad \mathbb P(\tau\ge k)=F_{k-1}(y).
\]
Optional stopping \cite[Section~4.8]{Durrett2019} therefore gives
\[
 \begin{aligned}
 \mathbb E[Y_m;\tau>m]
 &=
 y+\mathbb E g_{\tau\wedge m}
   -\mathbb E[Y_\tau;\tau\le m]\\
 &\le
 y+C_1+
 \sum_{k=1}^m(g_k-g_{k-1})F_{k-1}(y).
 \end{aligned}
\]
Now let $n\ge4$ and put $m=\lfloor n/2\rfloor$.
For $1\le r\le n-m$,
\[
 0\le g_{m+r}-g_m
 \le\kappa\log3.
\]
On $\{\tau>m\}$, condition on $\mathcal F_m$. The conditional
probability of survival to step $n$ is at most the probability that
a new centered walk started from $Y_m+\kappa\log3$ stays positive
for $n-m$ steps. Since $X_1\overset{\mathrm d}=-(\xi_1+c_\lambda T_1)$,
the two residual second-moment bounds
\eqref{eq:s9tan-residual-bound} hold with the uniform constant verified
in Lemma~\ref{lem:tangential-survival-entrance}. The straight-line
estimate in Lemma~\ref{lem:centered-cycle-exit} therefore yields
\[
 F_n(y)
 \le
 \frac{C}{\sqrt{1+n}}
 \left[
 1+y+
 \sum_{k=1}^m(g_k-g_{k-1})F_{k-1}(y)
 \right].
\]
Set
\[
 G_n
 =
 1+y+C_1+
 \sum_{k=1}^n(g_k-g_{k-1})F_{k-1}(y),\qquad n\in\N_0.
\]
The increments $g_k-g_{k-1}$ are nonnegative, so $G_k$ is
nondecreasing. The bound obtained by conditioning at
$m=\lfloor(n-1)/2\rfloor$ reads
\[
 F_{n-1}(y)
 \le C\frac{G_{\lfloor(n-1)/2\rfloor}}{\sqrt n}
 \le C\frac{G_{n-1}}{\sqrt n}
\]
for large $n$; increasing $C$ covers the remaining indices. Using
$G_n-G_{n-1}=(g_n-g_{n-1})F_{n-1}(y)$ now gives
\[
 G_n
 \le
 G_{n-1}
 \left(
 1+
 C\frac{g_n-g_{n-1}}{\sqrt n}
 \right),\qquad n\in\N.
\]
Since $0\le g_n-g_{n-1}\le\kappa/n$ and
$\sum_{n\ge1}n^{-3/2}<\infty$, the product $\prod_{n\ge1}(1+C\kappa n^{-3/2})$ is finite.
Iteration of the inequality for $G_n$ therefore gives
\[
 G_n\le C_\kappa(1+y),\qquad n\in\N_0,
\]
which proves \eqref{eq:cycle-log-survival}.

\medskip
\noindent\textit{Step 2: the random renewal-time boundary.}
Let $\widehat\tau_0=0$ and
$\widehat\tau_k=\sum_{i=1}^k\widehat T_i$ for $k\in\N$.
Start the reversed chain at $(l,0)$ and measure its height from
$z>0$ by $z+\widehat N_s-l-c_\lambda s$.
Consider
\[
 \mathcal E_t
 =
 \left\{
 z+\widehat N_s-\widehat N_0-c_\lambda s
 >
 -\kappa_{\parallel}
 \log\left(1+\frac{s}{1+a}\right)
 \text{ for }0\le s\le t
 \right\}.
\]
We prove
\begin{equation}
 \label{eq:random-renewal-log-survival}
 \widehat{\mathbb P}^\lambda_{(l,0)}(\mathcal E_t)
 \le
 C\frac{1+z}{\sqrt{1+t}},
\end{equation}
for $t\ge0$ and $z>0$, uniformly for $\lambda\in I$ and $a\ge a_*$.

On $\mathcal E_t$, for every $k$ with $\widehat\tau_k\le t$,
the post-jump height satisfies
\[
 z+\sum_{i=1}^{k}\widehat\xi_i
 -c_\lambda\widehat\tau_k
 >
 -\kappa_{\parallel}
 \log\left(1+\frac{\widehat\tau_k}{1+a}\right).
\]
Let $c_0=\min_{\lambda\in I}c_\lambda>0$ and choose
$\varepsilon=c_0/(2\kappa_{\parallel})$, recalling that
$\kappa_{\parallel}=3/(2\lambda_{\parallel})>0$.
Since $a\ge a_*$,
\[
 \log\left(1+\frac{s}{1+a}\right)
 \le
 \varepsilon s+C_\varepsilon
 \quad\text{ for } s\ge0,
\]
with $C_\varepsilon$ independent of $a$. Using $\sum_{i\le k}\widehat\xi_i\le k$, the post-jump
survival inequality gives
\[
 \frac{c_0}{2}\widehat\tau_k
 <z+k+\kappa_{\parallel}C_\varepsilon
 \qquad\text{when }\widehat\tau_k\le t.
\]
Thus, on $\mathcal E_t$,
\[
 \widehat\tau_k
 \le
 C_{\mathrm{cl}}(1+z+k)
 \quad\text{ for } \widehat\tau_k\le t,
\]
where $C_{\mathrm{cl}}$ is independent of $a\ge a_*$ and
$\lambda\in I$. Choose $C_{\log}\ge1+C_{\mathrm{cl}}/(1+a_*)$. Then
$1+C_{\mathrm{cl}}(1+z+k)/(1+a)\le C_{\log}(1+z+k)$ for
$a\ge a_*$, $z>0$ and $k\ge0$.

On $\mathcal E_t$, for every integer $k\ge1$ such that
$\widehat\tau_k\le t$,
\[
 z+\sum_{i=1}^kX_i
 >
 -\kappa_{\parallel}
 \log\left(
 1+\frac{\widehat\tau_k}{1+a}
 \right)
 \ge
 -\kappa_{\parallel}
 \log\bigl(C_{\log}(1+z+k)\bigr).
\]
Put, for $z>0$,
\[
 y_z
 =
 z+\kappa_{\parallel}\log\bigl(C_{\log}(1+z)\bigr).
\]
On the same event and for the same indices $k$,
\[
 y_z+\sum_{i=1}^kX_i
 >
 -\kappa_{\parallel}
 \log\left(\frac{1+z+k}{1+z}\right)
 \ge
 -\kappa_{\parallel}\log(1+k).
\]
Thus $\mathcal E_t\cap\{\widehat\tau_n\le t\}$ is contained in
$\{y_z+\sum_{i=1}^kX_i+\kappa_{\parallel}\log(1+k)>0:\ 1\le k\le n\}$.
Applying \eqref{eq:cycle-log-survival} with initial height $y_z$
and using $y_z\le C(1+z)$, we obtain
\[
 \widehat{\mathbb P}^\lambda_{(l,0)}\{
 \mathcal E_t,\ \widehat\tau_n\le t
 \}
 \le
 C\frac{1+z}{\sqrt{1+n}}.
\]

Choose $\alpha>0$ so small that
\[
 \alpha
 \sup_{\lambda\in I}\widehat{\mathbb E}^\lambda_{(l,0)}\widehat T_1<\frac12.
\]
The waiting times $\widehat T_1$ and $T_1$ have the same law, and the
choice of $\alpha$ verifies the mean condition used in
\eqref{eq:s9tan-cycle-concentration}. Exponential Markov inequality
and their common exponential moment therefore give
\[
 \widehat{\mathbb P}^\lambda_{(l,0)}\{
 \widehat\tau_{\lfloor\alpha t\rfloor}>t
 \}
 \le Ce^{-\rho t}.
\]
Taking $n=\lfloor\alpha t\rfloor$ proves
\eqref{eq:random-renewal-log-survival}.

\medskip
\noindent\textit{Step 3: the logarithmically killed kernel.}
Let $T_{\mathrm{stop}}$ be the first horizontal jump after $t/4$ and let
$\widehat{T}_{\mathrm{stop}}$ be the first horizontal jump after $t/4$ for the reversed
path started at $(l,0)$. Write $L_R$ for the transverse generator and
$K^\lambda$ for the killed transverse generator defined in
Lemma~\ref{lem:tangential-regeneration}. On
$\ell^2(\pi^\lambda)$, $L_R\le0$ and $K^\lambda\le-\rho$ for a common
$\rho>0$, by the spectral estimate in
Lemma~\ref{lem:tangential-regeneration}. Since $\pi^\lambda$ is a
probability vector, $\|\mathbf1\|_{\ell^2(\pi^\lambda)}=1$ and
$\|\delta_0/\pi_0^\lambda\|_{\ell^2(\pi^\lambda)}
=(\pi_0^\lambda)^{-1/2}$. The two semigroup norm bounds and
Cauchy--Schwarz therefore give, for $v,u\ge0$,
\[
 \mathbb P^\lambda_{(j,0)}
 \{\text{no horizontal jump on }[v,v+u]\}
 =\left\langle e^{vL_R}\frac{\delta_0}{\pi_0^\lambda},
                    e^{uK^\lambda}\mathbf1\right\rangle_{\pi^\lambda}
 \le(\pi_0^\lambda)^{-1/2}e^{-\rho u}.
\]
The reversed chain has the same transverse generator and total
horizontal rate, so the no-horizontal-jump estimate also holds under
$\widehat{\mathbb P}^\lambda_{(l,0)}$.
Since $\inf_{\lambda\in I}\pi_0^\lambda>0$, the contribution of
$T_{\mathrm{stop}}>t/3$ or $\widehat{T}_{\mathrm{stop}}>t/3$ is at most $Ce^{-\rho t/12}$.

On the remaining paths, both cut durations belong to $[t/4,t/3]$,
both cut states have transverse coordinate zero, and the middle
interval has length at least $t/3$. Define $\alpha_y(ds,m)$ as the
subprobability measure of the first segment, surviving above
$b_{t,a}$ and stopped just after its horizontal jump at $T_{\mathrm{stop}}=s\le t/3$
at $(m,0)$. For the reversed segment put
\[
 \widehat Y_r=z+\widehat N_r-l-c_\lambda r,\qquad 0\le r\le t.
\]
The forward condition at time $t-r$ is
\[
 \widehat Y_r>b_{t,a}(t-r)-b_{t,a}(t)
 =-\kappa_{\parallel}\log\left(1+\frac r{1+a}\right),\qquad 0\le r\le t.
\]
Let $\widehat\alpha_z(dr,n)$ be the subprobability measure of the
reversed segment satisfying this inequality at every time up to and
including its stopping time, stopped just after its jump
at $\widehat{T}_{\mathrm{stop}}=r\le t/3$ at
$(n,0)$. Since $b_{t,a}\ge0$, Lemma~\ref{lem:tangential-survival-entrance}
at time $t/4$ gives
\[
 \sum_{m\in\Z}\int\alpha_y(ds,m)\le C\frac{1+y}{\sqrt{1+t}}.
\]
The reversed boundary just displayed and
\eqref{eq:random-renewal-log-survival}, also at time $t/4$, give
\[
 \sum_{n\in\Z}\int\widehat\alpha_z(dr,n)\le C\frac{1+z}{\sqrt{1+t}}.
\]
Expanding paths at these two cut jumps, and retaining no absorption
condition on the middle segment, bounds the contribution from
$T_{\mathrm{stop}},\widehat{T}_{\mathrm{stop}}\le t/3$ by
\[
 \sum_{m,n\in\Z}\int_{[t/4,t/3]}\!\int_{[t/4,t/3]}
 \alpha_y(ds,m)\,p_{t-s-r}^\lambda(n-m)\,
 \widehat\alpha_z(dr,n).
\]
Here the first cut is a post-jump state in forward time, while the
second is a pre-jump state in forward time; each cut jump is included
once in its outer segment. Reversal uses the invariant measure
$\pi_r^\lambda$ times counting measure in $n$. Since both cut states
have transverse coordinate zero, its weight ratio is
$\pi_0^\lambda/\pi_0^\lambda=1$.
For $t\ge1$, Lemma~\ref{lem:tangential-free-kernel} and
$t-s-r\ge t/3$ give $p_{t-s-r}^\lambda(n-m)\le C(1+t)^{-1/2}$.
Substituting the two mass bounds into this integral, and adding the
exceptional contribution, yields
\[
 \mathcal K_{t,a}^{y,\lambda}(j,l)
 \le C\frac{(1+y)(1+z)}{(1+t)^{3/2}}+Ce^{-\rho t/12}
 \le C\frac{(1+y)(1+z)}{(1+t)^{3/2}}.
\]
Increasing $C$ covers $0\le t\le1$, because a transition probability
is at most one and $y,z>0$.

For the lower estimate, first choose a compact interval
$I_1\Subset(0,\infty)$ whose interior contains $I$ and on which the
uniform estimates of Lemmas~\ref{lem:tangential-killed}
and~\ref{lem:tangential-regeneration} hold, shrinking $I$ if necessary.
Let $Y=y_0$, where $y_0$ is the endpoint threshold in
Lemma~\ref{lem:tangential-killed} for the interval $I_1$. Fix $A>0$
and consider $Y\le y,z\le A\sqrt t$.
For $t>0$, set
$d_{t,a}=\kappa_{\parallel}t^{-1}\log(1+t/(1+a))$.
For sufficiently large $t$, choose $\eta=\eta_{t,a}\in I_1$ near
$\lambda$ so that
\[
 \Lambda_{\parallel}'(\eta)
 =
 c_\lambda-d_{t,a}.
\]
The positive lower bound for $\Lambda_{\parallel}''$ on $I_1$ and
$d_{t,a}\le\kappa_{\parallel}\log(1+t)/t$ give
\[
 |\eta-\lambda|
 \le
 C\frac{\log(1+t)}t.
\]
The function $b_{t,a}$ is convex and vanishes at zero, hence
\[
 b_{t,a}(s)-d_{t,a}s\le0
 \quad\text{ for } 0\le s\le t.
\]
The choice $\eta\in I_1$ is uniform for $\lambda\in I$ and
$a\ge a_*$. Keep the path event fixed:
\[
 E_\lambda=\{y+N_s-j+c_\lambda s>b_{t,a}(s),\ 0\le s\le t\}.
\]
Under $\mathbb P^\eta_{(j,0)}$, the straight-line event
$\{y+N_s-j+c_\eta s>0,\ 0\le s\le t\}$ is contained in
$E_\lambda$, because $c_\eta=c_\lambda-d_{t,a}$.
On $\{(N_t,R_t)=(l,0)\}$, the height centered using $c_\eta$ is
$y+l-j+c_\eta t=z$. Consequently
Lemma~\ref{lem:tangential-killed} gives
\[
 \mathbb P^\eta_{(j,0)}
   (E_\lambda,(N_t,R_t)=(l,0))
 \ge K_t^{y,\eta}(j,l)
 \ge c_A\frac{yz}{t^{3/2}}.
\]
Throughout this change of measure, $E_\lambda$ retains the
centering velocity $c_\lambda$.

We return to $\lambda$ by computing the density on paths, keeping
$E_\lambda$ fixed. Both chains $(N_t,R_t)$ at parameters $\lambda$
and $\eta$ are
Doob transforms of the same linear operator $\mathcal A_{\parallel}$.
Recall that $H_\lambda(n,r)=e^{-\lambda n}h_r^\lambda$ and
$h_0^\lambda=1$. For distinct states
$x,x'\in\Z\times\N_0$ connected by a positive diffusion rate,
the rate ratio is
\[
 \frac{q_\lambda(x,x')}{q_\eta(x,x')}
 =\frac{H_\lambda(x')/H_\eta(x')}
        {H_\lambda(x)/H_\eta(x)}.
\]
For $x=(n,r)$, write
$q_\lambda(x)=\sum_{x'\ne x}q_\lambda(x,x')$.
The eigenfunction equation for $H_\lambda$ gives
$q_\lambda(x)=\Lambda_{\parallel}(\lambda)+D(q+1)-a_r$, and hence
$q_\lambda(x)-q_\eta(x)=\Lambda_{\parallel}(\lambda)
-\Lambda_{\parallel}(\eta)$.
The jump ratios telescope and the holding factors therefore contribute
$e^{(\Lambda_{\parallel}(\eta)-\Lambda_{\parallel}(\lambda))t}$.
Both chains have bounded rates, so finite paths exhaust their laws
on $[0,t]$. Since the initial and terminal transverse states are zero
and $h_0^\eta=h_0^\lambda=1$, for $t\ge0$, $j,l\in\Z$ and every
event $E$ measurable with respect to the path up to time $t$, one obtains
\begin{equation}\label{eq:tan-path-change}
 \mathbb P^\lambda_{(j,0)}(E,(N_t,R_t)=(l,0))
 =
 e^{(\Lambda_{\parallel}(\eta)-\Lambda_{\parallel}(\lambda))t
   +(\eta-\lambda)(l-j)}
 \mathbb P^\eta_{(j,0)}(E,(N_t,R_t)=(l,0)).
\end{equation}
Using $l-j=z-y-c_\eta t$ and
$\Lambda_{\parallel}'(\eta)=c_\lambda-d_{t,a}$, the logarithm of the
exponential factor is
\[
 \bigl[
 \Lambda_{\parallel}(\eta)
 -\Lambda_{\parallel}(\lambda)
 -\Lambda_{\parallel}'(\eta)(\eta-\lambda)
 \bigr]t
 +
 (\eta-\lambda)(z-y).
\]
Its absolute value is bounded by
\[
 Ct|\eta-\lambda|^2
 +
 2A\sqrt t\,|\eta-\lambda|,
\]
which tends to zero as $t\to\infty$, uniformly for $a\ge a_*$ and $\lambda\in I$.
Thus the exponential factor in \eqref{eq:tan-path-change} is bounded
below by a positive constant for all sufficiently large $t$,
uniformly in $a\ge a_*$ and $\lambda\in I$. Applying that identity
to the fixed event $E_\lambda$ and the preceding lower bound under
$\mathbb P^\eta_{(j,0)}$ proves the lower estimate.

\medskip
\noindent\textit{Step 4: the bridge occupation estimate.}
Shrink $I$ around $\lambda_{\parallel}$ so that
$\min_{\lambda\in I}c_\lambda>\kappa_{\parallel}/(1+a_*)$.
For this reduced interval, take $c_0=\min_{\lambda\in I}c_\lambda$.
Then
\[
 c_\lambda-b_{t,a}'(s)
 \ge
 c_0-\frac{\kappa_{\parallel}}{1+a_*}
 =:v_0>0
\]
for $0\le s\le t$, $t\ge0$, $a\ge a_*$ and $\lambda\in I$.
Hence between two horizontal jumps $Y_s$
increases at speed at least $v_0$. On a path that survives throughout
$[0,t]$, let $s_0\in[0,t)$ be zero or a horizontal jump time, put
$h=Y_{s_0}>0$, and let $s_1$ be the smaller of $t$ and the next
horizontal jump time after $s_0$. Then
\[
 \int_{s_0}^{s_1}
 (1+Y_s)e^{-\lambda_{\parallel}Y_s}\,\dd s
 \le
 \frac1{v_0}
 \int_h^\infty
 (1+r)e^{-\lambda_{\parallel}r}\,\dd r
 \le
 C(1+h)e^{-\lambda_{\parallel}h}.
\]
Summing this bound over the initial cycle and the cycles beginning
at horizontal jumps bounds the occupation integral by the corresponding
sum of $C(1+h)e^{-\lambda_{\parallel}h}$.

For $s\in[0,t]$ and $n\in\Z$, define
\[
 h_n=h_n(s)=y+n-j+c_\lambda s-b_{t,a}(s).
\]
At a horizontal jump ending at $(n,0)$ at time $s$,
the post-jump height is $h_n(s)$.
Write $q_\xi=De^{-\lambda\xi}$ for the rate of a horizontal jump
of size $\xi\in\{-1,1\}$.
Put
\[
 E=\{(N_t,R_t)=(l,0),\ Y_v>0\text{ for }0\le v\le t\}.
\]
Apart from the initial cycle, the expectation of the occupation
integral times $\mathbf1_E$ is bounded by
\begin{equation}\label{eq:tan-occupation-jump-convolution}
 \begin{aligned}
 C\int_0^t\sum_{\xi\in\{-1,1\}}
 {\sum_{\substack{n\in\Z:\ h_n>0\\ h_n-\xi>0}}}
 &q_\xi(1+h_n)e^{-\lambda_{\parallel}h_n}\\
 &\times
 \mathcal K_{s,a+t-s}^{y,\lambda}(j,n-\xi)
 \mathcal K_{t-s,a}^{h_n,\lambda}(n,l)
 \,\dd s.
 \end{aligned}
\end{equation}
To justify \eqref{eq:tan-occupation-jump-convolution}, let $J_s=(N_s,R_s)$ and write
$q_\lambda(x,x')$ for its transition rate from $x$ to $x'$, where
$x,x'\in\Z\times\N_0$. On $(0,t]$, predictable means measurable
for the sigma-algebra generated by functions
$\mathbf1_{(s_0,s_1]}(s)H(\omega,x,x')$, where
$0\le s_0<s_1\le t$ are deterministic and $H$ is measurable with
respect to the path filtration at time $s_0$.
For nonnegative predictable
$F(s,x,x')$, the Poisson-clock construction gives
\[
 \mathbb E^\lambda_{(j,0)}\sum_{s\le t:\,J_s\ne J_{s-}}
 F(s,J_{s-},J_s)
 =\mathbb E^\lambda_{(j,0)}\int_0^t\sum_{x'\ne J_s}
 q_\lambda(J_s,x')F(s,J_s,x')\,\dd s.
\]
For step functions depending on the past this follows by conditioning
on each Poisson increment; nonnegative approximation and monotone
convergence give the displayed identity. The pre-jump state is $J_{s-}$. Each path has finitely many jumps
on $[0,t]$, so $J_{s-}=J_s$ for Lebesgue-almost every $s$.
The time integral consequently has the displayed form. The jump-and-hold construction used here is described in
\cite[Section~2.6]{Norris}.
For a horizontal transition $(n-\xi,0)\to(n,0)$, use the predictable
integrand
\[
 \begin{aligned}
 F(s,(n-\xi,0),(n,0))
 ={}&\mathbf1_{\{Y_v>0:\,0\le v<s\}}
      \mathbf1_{\{h_n>0,\ h_n-\xi>0\}}\\
 &\times(1+h_n)e^{-\lambda_{\parallel}h_n}
       \mathcal K_{t-s,a}^{h_n,\lambda}(n,l),
 \end{aligned}
\]
and set $F$ to zero on every other transition. The path indicator
$\mathbf1_{\{Y_v>0:\,0\le v<s\}}$ equals one up to and including
the first exit time and zero afterwards. In the present setting exits
can occur only at horizontal jumps, because $Y$ increases between
those jumps. The indicator is therefore left-continuous and adapted
as a function of $s$, and hence predictable. The final kernel is a deterministic function
of $s$ and the post-jump state. Conditioning the future on $(n,0)$
therefore gives that final kernel. In the compensator integral, the
law of the surviving pre-jump state $(n-\xi,0)$ is
$\mathcal K_{s,a+t-s}^{y,\lambda}(j,n-\xi)$, and the transition
rate is $q_\xi$. This gives the displayed bound for the expectation
with $\mathbf1_E$.
The offset in the first kernel follows from
\[
 b_{t,a}(r)=b_{s,a+t-s}(r)\quad\text{ for } 0\le r\le s,
 \qquad
 b_{t,a}(s+r)-b_{t,a}(s)=b_{t-s,a}(r)
 \quad\text{ for }0\le r\le t-s.
\]
The terminal height of the first kernel is $h_n-\xi$, and that of the
second is $z$.  For $s\in[0,t]$ and $\xi\in\{-1,1\}$, since $q_\xi$ is uniformly bounded and $h_n$ ranges over
a translate of the integers,
\[
 {\sum_{\substack{n\in\Z:\ h_n>0\\h_n-\xi>0}}}
 (1+h_n)^2(1+h_n-\xi)e^{-\lambda_{\parallel}h_n}
 \le C.
\]
The upper bound for $\mathcal K$ proved in Step~3 applies
to $\mathcal K_{s,a+t-s}^{y,\lambda}(j,n-\xi)$ with endpoint
heights $(y,h_n-\xi)$ and to
$\mathcal K_{t-s,a}^{h_n,\lambda}(n,l)$ with endpoint heights
$(h_n,z)$. Since $a+t-s\ge a_*$, both offsets are admissible. Thus, for
$0\le s\le t$, $n\in\Z$, $\xi\in\{-1,1\}$ with
$h_n>0$ and $h_n-\xi>0$, we obtain
\[
 \mathcal K_{s,a+t-s}^{y,\lambda}(j,n-\xi)
 \le C\frac{(1+y)(1+h_n-\xi)}{(1+s)^{3/2}},
 \qquad
 \mathcal K_{t-s,a}^{h_n,\lambda}(n,l)
 \le C\frac{(1+h_n)(1+z)}{(1+t-s)^{3/2}}.
\]
Using these kernel bounds and the sum over $h_n$ in
\eqref{eq:tan-occupation-jump-convolution}, we obtain
\[
 C(1+y)(1+z)
 \int_0^t\frac{\dd s}{(1+s)^{3/2}(1+t-s)^{3/2}}.
\]
For $s\le t/2$, one has $1+t-s\ge(1+t)/2$. Applying this inequality
on $[0,t/2]$ and the substitution $s\mapsto t-s$ on the other half gives
\begin{equation}\label{eq:tan-time-convolution}
 \int_0^t\frac{\dd s}{(1+s)^{3/2}(1+t-s)^{3/2}}
 \le\frac{2^{5/2}}{(1+t)^{3/2}}
      \int_0^\infty\frac{\dd s}{(1+s)^{3/2}}
 =\frac{2^{7/2}}{(1+t)^{3/2}}.
\end{equation}
The integral over the initial cycle is at most
$C(1+y)e^{-\lambda_{\parallel}y}$ on $E$, so its contribution is
bounded by $C(1+y)e^{-\lambda_{\parallel}y}\mathbb P^\lambda_{(j,0)}(E)$.
Every later incomplete cycle is already included in the marked-jump
sum: if it starts at height $h>0$, its integral up to time $t$ is at
most $C(1+h)e^{-\lambda_{\parallel}h}$ by the cycle estimate at the
start of Step~4. For $y,z$ in the fixed compact endpoint set,
Step~3 gives, for all sufficiently large $t$,
\[
 c(1+t)^{-3/2}
 \le\mathbb P^\lambda_{(j,0)}(E)
 =\mathcal K_{t,a}^{y,\lambda}(j,l)
 \le C(1+t)^{-3/2}.
\]
The marked-jump contribution and the initial-cycle contribution are
therefore both at most $C(1+t)^{-3/2}$. Dividing by the probability
in this display proves the conditional occupation bound. For bounded
$t$, the bound follows directly from
\[
 \int_0^t(1+Y_s)e^{-\lambda_{\parallel}Y_s}\,\dd s
 \le t\sup_{r\ge0}(1+r)e^{-\lambda_{\parallel}r}
 \quad\text{on }E,
\]
whenever $\mathbb P^\lambda_{(j,0)}(E)>0$.\qedhere

\end{proof}

We next derive a path representation for the linear equation on a
moving subset of $\Z\times\N_0$, then use it for comparison with
the nonlinear equation.

\paragraph{Two-coordinate moving-domain formula.}
Let
\[
 \Omega_t
 =
 \{(n,r)\in\Z\times\N_0:n>m(t)\},\qquad t\ge0,
\]
where $m\in C^1([0,\infty))$ and $m'(t)>0$ for every $t\ge0$.
Thus each integer is crossed at most once, and the active set changes
at only finitely many times on each bounded interval. {The functions $m_b$ defined at the start of
Subsection~\ref{sec:tangential-bramson}} satisfy $m_b'\!>0$ and belong to
$C^1([0,\infty))$. The diffusion
part of $\mathcal A_{\parallel}$ is the generator of the coordinate
chain with rates, for $n\in\Z$,
\[
 \begin{aligned}
 &(n,0)\to(n\pm1,0)&&\text{at rate }D,\\
 &(n,0)\to(n,1)&&\text{at rate }D(q-1),\\
 &(n,r)\to(n,r-1)&&\text{at rate }D\text{ for }r\ge1,\\
 &(n,r)\to(n,r+1)&&\text{at rate }Dq\text{ for }r\ge1.
 \end{aligned}
\]
In this moving-domain argument, $(N_s,R_s)_{s\ge0}$ denotes this
diffusion chain rather than the tilted chain with generator
$\mathcal Q_{\parallel}^\lambda$. Write
$\mathbb P^{\mathrm d}_{(n,r)}$ and $\mathbb E^{\mathrm d}_{(n,r)}$
for its law and expectation when $(N_0,R_0)=(n,r)$.
Its total jump rate is $D(q+1)$ at every state. If
$V:[0,\infty)\times\Z\times\N_0\to[0,\infty)$ is Borel measurable
and bounded uniformly in $(n,r)$ on each bounded time interval,
consider the linear equation
\[
 v_t=\mathcal A_{\parallel}v-Vv
 \quad\text{for }(n,r)\in\Omega_t\text{ and almost every }t>0
 \text{ between removal times}.
\]
Prescribe zero values outside $\Omega_t$ and a bounded nonnegative
initial datum $v_0:\Z\times\N_0\to[0,\infty)$ vanishing outside
$\Omega_0$. Between removal times take $v$ to satisfy the coordinatewise
integral equation, and at a removal time set the newly removed
coordinates to zero, retaining the other coordinates. With this
convention, the solution has, for $t\ge0$ and
$(n,r)\in\Z\times\N_0$, the path representation
\[
 \begin{aligned}
 v(t,n,r)
 =
 \mathbb E^{\mathrm d}_{(n,r)}
 \Biggl[
 &\exp\left(
 \int_0^t
 [a_{R_s}-V(t-s,N_s,R_s)]\,\dd s
 \right)
 v_0(N_t,R_t);\\
 &N_s>m(t-s)\quad\text{ for } 0\le s\le t
 \Biggr].
 \end{aligned}
\]
To verify the formula, let $q^{\mathrm d}(x,x')$, for distinct
$x,x'\in\Z\times\N_0$, equal the listed diffusion rate when
$x\to x'$ is one of the listed transitions and zero otherwise.
Put $\varrho=D(q+1)$, the common total rate. Fix $t\ge0$ and the
starting state $x_0=(n,r)$ under $\mathbb P^{\mathrm d}_{(n,r)}$. Expand the expectation over paths with $k$ jumps,
states $x_0=(n,r),x_1,\ldots,x_k$, and jump times
$0<s_1<\cdots<s_k<t$. The factor for such a path is
\[
 e^{-\varrho t}
 \prod_{i=1}^k q^{\mathrm d}(x_{i-1},x_i)
 \exp\!\left(\int_0^t[a_{R_s}-V(t-s,N_s,R_s)]\,\dd s\right)
 v_0(x_k)
 \mathbf1_{\{N_s>m(t-s):\,0\le s\le t\}}.
\]
After summing over the states and integrating over the ordered jump
times, the term with $k$ jumps is at most
\[
 \|v_0\|_\infty e^{a_{\max}t}
 e^{-\varrho t}\frac{(\varrho t)^k}{k!}.
\]
Thus the series converges absolutely, uniformly on bounded time
intervals. Conditioning on a first short time interval gives the
coordinatewise integral equation with generator
$\mathcal A_{\parallel}-V$ on each fixed active set, and hence the
stated differential equation almost everywhere in time. If two bounded solutions with this reset convention coincide at
$t_1$, write $d(t,n,r)$ for their difference. On an interval
$[t_1,t_2]$ containing no removal in its interior, compare their
evolutions up to the left limit at a removal endpoint. The
coordinatewise integral equations give, for $t_1\le t<t_2$,
\[
 \|d(t,\cdot)\|_\infty
 \le\left(2\varrho+a_{\max}
       +\sup_{\substack{0\le s\le t_2\\(n,r)\in\Z\times\N_0}}V(s,n,r)\right)
      \int_{t_1}^t\|d(s,\cdot)\|_\infty\,\dd s.
\]
Gronwall's inequality gives $d=0$ there.
At a removal time $t$, the condition at path time $s=0$ in the
survival indicator sets $v(t,n,r)=0$ whenever $n\le m(t)$. Since $m'>0$, only removals occur, and there are finitely many
such times on each bounded interval. Concatenating these integral
equations proves the stated formula, including the resets.

Let $u:[0,\infty)\times\Z\times\N_0\to[0,1]$ solve
$u_t=\mathcal A_{\parallel}u+f_r(u)-a_ru$ and satisfy
$v_0(n,r)\le u(0,n,r)$ for every $(n,r)\in\Z\times\N_0$.
Assume $0\le v(t,n,r)\le1$ and $V(t,n,r)\ge C_2v(t,n,r)$ for
$t\ge0$ and $(n,r)\in\Z\times\N_0$. Then for almost every $t>0$ between removal times
and every $(n,r)\in\Omega_t$
\[
 v_t
 =\mathcal A_{\parallel}v-Vv
 \le\mathcal A_{\parallel}v-C_2v^2
 \le\mathcal A_{\parallel}v+f_r(v)-a_rv
\]
by \eqref{eq:taylor}. At a site $(n,r)\notin\Omega_t$, one has $v(t,n,r)=0$; away from
removal times its time derivative is zero and its diffusion term
is a sum of nonnegative neighboring values with nonnegative rates. The zero value therefore also satisfies the subsolution
inequality. Fix $0\le t_0<t_1<\infty$ with no removal in
$(t_0,t_1)$, and suppose $v(t_0,n,r)\le u(t_0,n,r)$ for every
$(n,r)\in\Z\times\N_0$, using the reset values at $t_0$ when
$t_0$ is a removal time. The coordinate integral equations, the
common total diffusion rate $\varrho=D(q+1)$, and the bound on $V$
give local Lipschitz continuity of $u$ and $v$ in
$\ell^\infty(\Z\times\N_0)$ before the next removal.
Put
\[
 M(s)=\sup_{(n,r)\in\Z\times\N_0}
       (v(s,n,r)-u(s,n,r))_+,
 \qquad \beta=\varrho+L,
 \quad\text{for }t_0\le s<t_1.
\]
Since $z\mapsto Lz+f_r(z)$ is nondecreasing on $[0,1]$ for every
$r\in\N_0$, integration of the coordinate inequalities, as in
\eqref{eq:kern-coordinate-comparison}, gives
\[
 M(t)\le e^{-\beta(t-t_0)}M(t_0)
       +(\beta+L)\int_{t_0}^t e^{-\beta(t-s)}M(s)\,\dd s
 \quad\text{for }t_0\le t<t_1.
\]
The initial ordering gives $M(t_0)=0$, so Gronwall's inequality
yields $M(t)=0$ for $t_0\le t<t_1$. At a removal time the retained
coordinates keep their left limits, and resetting the other
coordinates to zero preserves their order with $u\ge0$.
Induction over the finitely many removal times on each bounded
interval therefore gives $v(t,n,r)\le u(t,n,r)$ for every $t\ge0$
and $(n,r)\in\Z\times\N_0$.

Likewise, let $u:[0,\infty)\times\Z\times\N_0\to[0,1]$ satisfy
$u_t=\mathcal A_{\parallel}u+(f_r(u)-a_ru)$ on $\Z\times\N_0$,
with $u(0,n,r)=0$ for $(n,r)\in\Omega_0$. Fix $t\ge0$ and
$(n,r)\in\Omega_t$, and start the diffusion chain at $(n,r)$. Put
\[
 \tau
 =
 \inf\{s\in[0,t]:N_s\le m(t-s)\},
\]
with $\inf\varnothing=\infty$. Then $\tau$ is a stopping time
for the diffusion chain. Starting at $(n,r)\in\Omega_t$, the height
$N_s-m(t-s)$ strictly increases whenever $N_s$ is constant,
because $m'>0$. Hence an exit occurs through a horizontal jump.
We evaluate the exit state immediately after that jump. For $t\ge0$ and $(n,r)\in\Omega_t$
\begin{equation}
 \label{eq:tangential-stopped-comparison}
 u(t,n,r)
 \le
 \mathbb E^{\mathrm d}_{(n,r)}
 \left[
 \exp\left(\int_0^\tau a_{R_s}\,\dd s\right);
 \tau\le t
 \right].
\end{equation}
For this stopped estimate, put
$\mathcal R_r(z)=a_rz-f_r(z)\ge0$ for $r\in\N_0$ and $0\le z\le1$.
Write $\mathcal L^{\mathrm d}$ for the diffusion generator, so that
$(\mathcal L^{\mathrm d}F)_{n,r}
=(\mathcal A_{\parallel}F)_{n,r}-a_rF_{n,r}$. For $0\le s<t$ and $(n,r)\in\Z\times\N_0$, the nonlinear equation gives the local identity
\[
 -u_t(t-s,n,r)+\mathcal L^{\mathrm d}u(t-s,n,r)
       +a_ru(t-s,n,r)
 =\mathcal R_r\bigl(u(t-s,n,r)\bigr).
\]
The solution and its time derivative are uniformly bounded on each
bounded time interval, by the bounded jump rates and
$|f_r'|\le L$. Differentiation between jumps and compensation of
the bounded-rate jumps therefore show that
\[
 \begin{aligned}
 M_s={}&e^{\int_0^s a_{R_v}\,\dd v}u(t-s,N_s,R_s)-u(t,n,r)\\
 &-\int_0^s e^{\int_0^v a_{R_w}\,\dd w}
       \mathcal R_{R_v}\bigl(u(t-v,N_v,R_v)\bigr)\,\dd v
 \end{aligned}
\]
is a martingale under $\mathbb P^{\mathrm d}_{(n,r)}$, for the
completed, right-continuous natural filtration of the diffusion
chain and $0\le s\le t$. Indeed, the drift of its first
term is $e^{\int_0^s a_{R_v}\,\dd v}\mathcal R_{R_s}
(u(t-s,N_s,R_s))$, which is canceled by the integral. The bounds
$0\le u\le1$ and $0\le\mathcal R_r(u)\le a_{\max}$ bound $M_s$
on the fixed interval $[0,t]$.
For $T_{\mathrm{stop}}=t\wedge\tau$ and each integer $k\ge1$, apply discrete-time
optional stopping to $M$ on the finite grid
$\{m2^{-k}:m\in\N_0,\ m2^{-k}<t\}\cup\{t\}$, at
$T_{\mathrm{stop},k}=(2^{-k}\lceil2^kT_{\mathrm{stop}}\rceil)\wedge t$. Right-continuity and bounded
convergence as $k\to\infty$ give $\mathbb E^{\mathrm d}_{(n,r)}M_{T_{\mathrm{stop}}}=0$,
and hence
\[
 \begin{aligned}
 u(t,n,r)
 ={}&\mathbb E^{\mathrm d}_{(n,r)}
 \left[e^{\int_0^{T_{\mathrm{stop}}} a_{R_v}\,\dd v}u(t-T_{\mathrm{stop}},N_{T_{\mathrm{stop}}},R_{T_{\mathrm{stop}}})\right]\\
 &-\mathbb E^{\mathrm d}_{(n,r)}
 \left[\int_0^{T_{\mathrm{stop}}} e^{\int_0^s a_{R_v}\,\dd v}
       \mathcal R_{R_s}\bigl(u(t-s,N_s,R_s)\bigr)\,\dd s\right].
 \end{aligned}
\]
On $\{\tau>t\}$ the first expectation has zero terminal value,
because $(N_t,R_t)\in\Omega_0$ and $u(0,\cdot)=0$ there. On
$\{\tau\le t\}$, the factor $u(t-\tau,N_\tau,R_\tau)$ is at most one.
Dropping the nonnegative second expectation therefore proves
\eqref{eq:tangential-stopped-comparison}.

\begin{lemma}\label{lem:tangential-linear-seed}

Let $Y_*$ be a choice of the constant $Y$ in
Lemma~\ref{lem:tangential-backward-log}. Fix $b\in\Z$ and $z_0>Y_*$,
and put
\[
 m(t)=b-z_0+c_{\parallel}t
 -\kappa_{\parallel}\log\frac{1+a_*+t}{1+a_*}
 \quad\text{for }t\ge0.
\]
Let $W$ solve $W_t=\mathcal A_{\parallel}W$ on
$\{(n,r)\in\Z\times\N_0:n>m(t)\}$, with
$W(0,n,r)=\mathbf1_{\{(n,r)=(b,0)\}}$ and
$W(t,n,r)=0$ when $n\le m(t)$; newly removed coordinates are set to zero.
There is $C<\infty$ such that
\[
 0\le W(t,n,r)\le C(1+y)e^{-\lambda_{\parallel}y},
 \qquad y=n-m(t)>0,
\]
for every $t\ge0$, $n\in\Z$ and $r\in\N_0$ with $n>m(t)$.
Moreover, for every fixed interval $[Y,B]$ with $Y_*<Y\le B<\infty$,
\[
 W(t,n,0)\ge c_{Y,B}>0
\]
for all sufficiently large $t$ whenever
$Y\le n-m(t)\le B$.

\end{lemma}

\begin{proof}

Write $\Lambda_*=\Lambda_{\parallel}(\lambda_{\parallel})$ in this proof.
Recall that $H_{\lambda_{\parallel}}(n,r)
=e^{-\lambda_{\parallel}n}h_r^{\lambda_{\parallel}}$ and
$h_0^{\lambda_{\parallel}}=1$. Since
$n-b=y-z_0+c_{\parallel}t
-\kappa_{\parallel}\log((1+a_*+t)/(1+a_*))$ and
$\Lambda_*=c_{\parallel}\lambda_{\parallel}$, the moving-domain
representation at $r=0$ gives the exact formula
\begin{equation}\label{eq:tan-linear-seed-kernel}
 W(t,n,0)
 =
 e^{\lambda_{\parallel}z_0}
 \left(\frac{1+a_*+t}{1+a_*}\right)^{3/2}
 e^{-\lambda_{\parallel}y}
 \mathcal K_{t,a_*}^{y,\lambda_{\parallel}}(n,b),
 \qquad t\ge0,\quad n\in\Z,\quad y=n-m(t)>0.
\end{equation}
The terminal height in \eqref{eq:tan-linear-seed-kernel} is $z_0$. The upper bound in
Lemma~\ref{lem:tangential-backward-log} therefore gives
\[
 W(t,n,0)
 \le C e^{-\lambda_{\parallel}y}
 \left(\frac{1+a_*+t}{1+a_*}\right)^{3/2}
 \frac{(1+y)(1+z_0)}{(1+t)^{3/2}}
 \le C(1+y)e^{-\lambda_{\parallel}y}.
\]
Choose $Y$ above the lower endpoint threshold in
Lemma~\ref{lem:tangential-backward-log}; the fixed $z_0$ in the
statement already exceeds that threshold. For $y\in[Y,B]$ and
$t\ge\max\{B,z_0\}^2$, we have $y,z_0\le\sqrt t$. The lower bound
of that lemma with $A=1$ therefore applies once $t$ also exceeds its
time threshold, and yields
\[
 W(t,n,0)
 \ge c y z_0e^{-\lambda_{\parallel}y}
       \left(\frac{1+a_*+t}{1+a_*}\right)^{3/2}t^{-3/2}
 \ge c_{Y,B}>0.
\]
For $r\in\N$, let $\mathbb P_r^{\lambda_{\parallel}}$ and
$\mathbb E_r^{\lambda_{\parallel}}$ denote the law and expectation
of the transverse chain $R$ in Lemma~\ref{lem:tangential-regeneration}
at $\lambda=\lambda_{\parallel}$, started from $r$. Put
$\sigma_0=\inf\{s\ge0:R_s=0\}$,
with $\inf\varnothing=\infty$.
Before $\sigma_0$ the
longitudinal coordinate is constant, while $n-m(t-s)$ increases
because $m'>0$, so no boundary killing occurs.  The initial datum is supported at $(b,0)$, so only paths with
$\sigma_0\le t$ contribute to the moving-domain representation.
The identity
$\mathcal A_{\parallel}(H_{\lambda_{\parallel}}F)
=H_{\lambda_{\parallel}}(\mathcal Q_{\parallel}^{\lambda_{\parallel}}F
+\Lambda_*F)$ from Lemma~\ref{lem:tangential-regeneration}, followed
by the strong Markov property at $\sigma_0$, gives
\begin{equation}\label{eq:tan-first-hit-decomposition}
 W(t,n,r)
 =h_r^{\lambda_{\parallel}}
  \mathbb E_r^{\lambda_{\parallel}}
  \left[e^{\Lambda_*\sigma_0}
        W(t-\sigma_0,n,0);\ \sigma_0\le t\right],
 \qquad t\ge0,\quad n\in\Z,\quad n>m(t),\quad r\in\N.
\end{equation}
To integrate the upper bound for $W(t-\sigma_0,n,0)$ in
\eqref{eq:tan-first-hit-decomposition}, we estimate the law of
$\sigma_0$ after multiplication by $h_r^{\lambda_{\parallel}}$,
uniformly for $r\in\N$.

Let $\mathcal B_1$ be the bounded self-adjoint Jacobi operator on
$\ell^2(\N)$ with diagonal $a_r-D(q+1)$, adjacent off-diagonal entries
$D\sqrt q$, and zero exterior value at $r=0$.
Lemma~\ref{lem:geodesic-resolvent} gives
$\sup\spec(\mathcal B_1)=\mu$. Hence, for every $s\ge0$ and $r\ge1$,
$|(e^{s\mathcal B_1})(r,1)|\le\|e^{s\mathcal B_1}\|\le e^{\mu s}$
by the spectral theorem. In the symmetric transverse coordinates,
the killed Doob semigroup satisfies
\[
 \mathbb P_r^{\lambda_{\parallel}}
 \{R_s=1,\ \sigma_0>s\}
 =
 e^{-\Lambda_*s}
 \frac{\psi_1^{\lambda_{\parallel}}}
      {\psi_r^{\lambda_{\parallel}}}
 (e^{s\mathcal B_1})(r,1).
\]
The jump rate from transverse state $1$ to $0$ is
\[
 \tau_0
 \frac{\psi_0^{\lambda_{\parallel}}}
      {\psi_1^{\lambda_{\parallel}}},
 \qquad
 \tau_0=D\sqrt{q-1}.
\]
For $\beta<\Lambda_*-\mu$, this gives
\[
 \begin{aligned}
 h_r^{\lambda_{\parallel}}
 \mathbb E_r^{\lambda_{\parallel}}
 \left[
 (1+\sigma_0)e^{\beta\sigma_0};
 \sigma_0<\infty
 \right]
 &\le
 C\int_0^\infty
 (1+s)
 e^{-(\Lambda_*-\beta)s}
 |(e^{s\mathcal B_1})(r,1)|\,\dd s\\
 &\le
 C\int_0^\infty
 (1+s)
 e^{-(\Lambda_*-\beta-\mu)s}\,\dd s.
 \end{aligned}
\]
Indeed, multiplying the killed transition probability by the rate
from $1$ to $0$ gives the first-hit density
\begin{equation}\label{eq:tan-first-hit-density}
 \mathbb P_r^{\lambda_{\parallel}}(\sigma_0\in\dd s)
 =\frac{\tau_0}{\psi_r^{\lambda_{\parallel}}}
   e^{-\Lambda_*s}(e^{s\mathcal B_1})(r,1)\,\dd s,
 \qquad s>0,\quad r\in\N.
\end{equation}
Here $\psi_0^{\lambda_{\parallel}}=1$. Multiplying \eqref{eq:tan-first-hit-density} by
$h_r^{\lambda_{\parallel}}
=\psi_r^{\lambda_{\parallel}}/\sqrt{\nu_r}$ leaves the factor
$\tau_0/\sqrt{\nu_r}\le\tau_0$, uniformly for $r\ge1$.
Choose
\[
 \frac{\lambda_{\parallel}\kappa_{\parallel}}{1+a_*}
 <\beta<
 \Lambda_{\parallel}(\lambda_{\parallel})-\mu,
\]
which is possible by the choice of $a_*$. Since
$\Lambda_*-\beta-\mu>0$, the integral
$\int_0^\infty(1+s)e^{-(\Lambda_*-\beta-\mu)s}\,\dd s$
is finite and independent of $r$.

For $0\le s\le t$,
\[
 y+\left(c_{\parallel}
          -\frac{\kappa_{\parallel}}{1+a_*}\right)s
 \le n-m(t-s)\le y+c_{\parallel}s.
\]
Since $n>m(t)\ge m(t-s)$, the upper bound at transverse state zero gives
$W(t-s,n,0)\le C[1+n-m(t-s)]
 e^{-\lambda_{\parallel}[n-m(t-s)]}$ for $0\le s\le t$.
The inequality $n-m(t-s)\le y+c_{\parallel}s$ implies
$1+n-m(t-s)\le C(1+y)(1+s)$. After factoring $e^{-\lambda_{\parallel}y}$ out of the
$r=0$ bound at time $t-s$, the exponential multiplier in
\eqref{eq:tan-first-hit-decomposition} is bounded by
\[
 e^{\Lambda_*s}
 e^{-\lambda_{\parallel}
      (m(t)-m(t-s))}
 \le
 e^{\lambda_{\parallel}\kappa_{\parallel}s/(1+a_*)}.
\]
Since $\beta>\lambda_{\parallel}\kappa_{\parallel}/(1+a_*)$,
the uniform bound on
$h_r^{\lambda_{\parallel}}\mathbb E_r^{\lambda_{\parallel}}
[(1+\sigma_0)e^{\beta\sigma_0};\sigma_0<\infty]$
controls the resulting integral. Substituting the estimate at $r=0$ into
\eqref{eq:tan-first-hit-decomposition}
therefore proves
\[
 W(t,n,r)\le C(1+y)e^{-\lambda_{\parallel}y}
\]
with $C$ independent of $r$.\qedhere

\end{proof}

\begin{lemma}\label{lem:tangential-positive-seed}

Let $b\in\Z$. There is $\varepsilon_0\in(0,1]$ such that the
following holds for every $\varepsilon\in(0,\varepsilon_0]$.
Let $U=(U_{n,r})_{(n,r)\in\Z\times\N_0}$ solve
$\dot U_{n,r}=(\mathcal A_{\parallel}U)_{n,r}+f_r(U_{n,r})-a_rU_{n,r}$
for $t\ge0$, with initial datum $\varepsilon\mathbf1_{\{(b,0)\}}$.
For each such $\varepsilon$, there are $\delta>0$, $C_0<\infty$,
and $T_0<\infty$ such that
\[
 U_{n,0}(t)\ge\delta
\]
whenever
\[
 t\ge T_0,
 \qquad
 |n-b|
 \le
 c_{\parallel}t
 -\kappa_{\parallel}\log(1+t)
 -C_0.
\]

\end{lemma}

\Needspace{7\baselineskip}
\begin{proof}

Fix the integer $b\in\Z$ specifying the single-vertex initial datum.
Choose $z_0$ above the endpoint threshold in
Lemma~\ref{lem:tangential-backward-log}, and put
\[
 m(t)=b-z_0+c_{\parallel}t
 -\kappa_{\parallel}\log\frac{1+a_*+t}{1+a_*}.
\]
Let $W$ solve $W_t=\mathcal A_{\parallel}W$ on $n>m(t)$, with zero
exterior value and initial datum $\mathbf1_{\{(b,0)\}}$.
For $y\ge0$, define
\[
 \omega(y)=(1+y)e^{-\lambda_{\parallel}y}.
\]
Lemma~\ref{lem:tangential-linear-seed} supplies $C_W$ such that
$W(t,n,r)\le C_W\omega(n-m(t))$ at every active site.
On $n>m(t)$, let $v$ solve $v_t=\mathcal A_{\parallel}v-Vv$
with zero exterior value, initial datum
$\varepsilon\mathbf1_{\{(b,0)\}}$, and potential
\[
 V(t,n,r)
 =
 C_2\varepsilon C_W\omega(n-m(t)),
 \qquad t\ge0,\quad(n,r)\in\Z\times\N_0,\quad n>m(t).
\]
{Choose
$0<\varepsilon\le\min\{1,(C_W\sup_{y\ge0}\omega(y))^{-1}\}$.} By linear comparison
with $\varepsilon W$ and the definition of $V$,
\[
 0\le v\le\varepsilon W\le1,
 \qquad
 V\ge C_2v.
\]
Set $v(t,n,r)=V(t,n,r)=0$ when $n\le m(t)$.
For almost every $t>0$ between removal times and every
$(n,r)\in\Z\times\N_0$ with $n>m(t)$, the equation for $v$
and \eqref{eq:taylor} give
\[
 \begin{aligned}
 &\partial_tv(t,n,r)
 -(\mathcal A_{\parallel}v(t,\cdot,\cdot))_{n,r}
 -f_r(v(t,n,r))+a_rv(t,n,r)\\
 &\quad=a_rv(t,n,r)-f_r(v(t,n,r))
        -V(t,n,r)v(t,n,r)\\
 &\quad\le C_2v(t,n,r)^2-V(t,n,r)v(t,n,r)\le0.
 \end{aligned}
\]
At an inactive coordinate $(n,r)$, $v(t,n,r)=0$ and the diffusion
term is a sum of nonnegative neighboring values.
In coordinate variables, $U$ satisfies
$\partial_tU_{n,r}=(\mathcal A_{\parallel}U)_{n,r}
+f_r(U_{n,r})-a_rU_{n,r}$ and
$U_{n,r}(0)=v(0,n,r)$ for every $(n,r)\in\Z\times\N_0$.
The moving-domain comparison applies because $0\le v\le1$,
$V\ge C_2v$, and $v(0,n,r)=U_{n,r}(0)$; the zero assignments at
each removal time preserve the order. It therefore gives
$v(t,n,r)\le U_{n,r}(t)$ for $t\ge0$ and
$(n,r)\in\Z\times\N_0$.

Choose $Y$ above the endpoint threshold in
Lemma~\ref{lem:tangential-backward-log}. For $Y\le n-m(t)\le Y+1$, put $y=n-m(t)$.
Under $\mathbb P^{\lambda_{\parallel}}_{(n,0)}$, write
$Y_s=N_s-m(t-s)
=y+N_s-n+c_{\parallel}s-b_{t,a_*}(s)$ for $0\le s\le t$, and define
\[
 E_t=\{(N_t,R_t)=(b,0),\ N_s>m(t-s)\text{ for }0\le s\le t\}.
\]
Lemma~\ref{lem:tangential-backward-log} gives
$\mathbb P^{\lambda_{\parallel}}_{(n,0)}(E_t)
=\mathcal K_{t,a_*}^{y,\lambda_{\parallel}}(n,b)>0$
for all sufficiently large $t$. Apply the moving-domain path formula
to $v$ and $W$, and use the Doob conjugacy with
$H_{\lambda_{\parallel}}$. On the common terminal event
$(N_t,R_t)=(b,0)$, both expressions have the factor
$\exp\{\Lambda_{\parallel}(\lambda_{\parallel})t
-\lambda_{\parallel}(n-b)\}$, which cancels in their ratio.
Thus $v(t,n,0)/(\varepsilon W(t,n,0))$ is the expectation of
\[
 \exp\left(
 -C_2\varepsilon C_W
 \int_0^t\omega(Y_s)\,ds
 \right)
\]
under $\mathbb P^{\lambda_{\parallel}}_{(n,0)}$ conditional on $E_t$. The two endpoint heights
are $y\in[Y,Y+1]$ and $z_0$, so they lie in a fixed compact interval
above the threshold of that lemma. In particular,
Lemma~\ref{lem:tangential-backward-log} gives
$\mathbb E^{\lambda_{\parallel}}_{(n,0)}
[\int_0^t\omega(Y_s)\,\dd s\mid E_t]\le C$, with $C$ independent
of $t$ and of every integer $n$ satisfying $Y\le n-m(t)\le Y+1$. Jensen's inequality
for $x\mapsto e^{-x}$ therefore gives
\begin{equation}\label{eq:tan-nonlinear-seed-ratio}
 \begin{aligned}
 \frac{v(t,n,0)}{\varepsilon W(t,n,0)}
 &\ge \exp\left(-C_2\varepsilon C_W
       \mathbb E^{\lambda_{\parallel}}_{(n,0)}
       \left[\int_0^t\omega(Y_s)\,\dd s\,\middle|\,E_t\right]
       \right)
 \\&\ge e^{-C\varepsilon}>0.
 \end{aligned}
\end{equation}
Lemma~\ref{lem:tangential-linear-seed} gives
$W(t,n,0)\ge c_{Y,Y+1}>0$ for all sufficiently large $t$ in
$Y\le n-m(t)\le Y+1$. Set
$\delta=\varepsilon e^{-C\varepsilon}c_{Y,Y+1}>0$. Then
\eqref{eq:tan-nonlinear-seed-ratio} gives
\[
 U_{n,0}(t)\ge v(t,n,0)\ge\delta
\]
for every integer $n$ satisfying $Y\le n-m(t)\le Y+1$;
at least one such integer exists because the interval has length one.
Define $\widetilde U_{n,r}(t)=U_{2b-n,r}(t)$. The map
$n\mapsto2b-n$ interchanges the two horizontal neighbors at $r=0$
and leaves every transverse term and $f_r$ unchanged. Thus
$\widetilde U$ solves the same nonlinear equation as $U$ and has
initial datum $\varepsilon\mathbf1_{\{(b,0)\}}$. Uniqueness gives
\[
 U_{2b-n,r}(t)=U_{n,r}(t),
 \qquad t\ge0,\quad(n,r)\in\Z\times\N_0.
\]
In particular the lower bound holds at $2b-n$. To extend it between
these two points, set
\[
 d_{n,r}(t)=U_{b+n,r}(t)-U_{b+n+1,r}(t),
 \quad\text{ for }n,r\in\N_0,\quad t\ge0.
\]
For $(n,r)\in\N_0\times\N_0$ and $t\ge0$, define
\[
 a_{n,r}(t)=\int_0^1 f_r'\bigl(
 U_{b+n+1,r}(t)+\tau(U_{b+n,r}(t)-U_{b+n+1,r}(t))
 \bigr)\,\dd\tau.
\]
The fundamental theorem of calculus gives
$f_r(U_{b+n,r}(t))-f_r(U_{b+n+1,r}(t))
=a_{n,r}(t)d_{n,r}(t)$.
By \eqref{eq:taylor}, $|a_{n,r}(t)|\le L$ for
$n,r\in\N_0$ and $t\ge0$. For $n\ge1$
the horizontal part at $r=0$ is
$D(d_{n-1,0}+d_{n+1,0}-2d_{n,0})$.  At $n=0$, reflection gives
$U_{b-1,r}=U_{b+1,r}$ and hence
\begin{equation}\label{eq:tan-longitudinal-difference}
 \dot d_{0,0}
 =D(d_{1,0}-3d_{0,0})
   +D(q-1)(d_{0,1}-d_{0,0})+a_{0,0}(t)d_{0,0},\qquad t>0.
\end{equation}
For $r\ge1$ the transverse difference equations retain the positive
off-diagonal coefficients $D$ and $Dq$. All diagonal coefficients
are at least $-D(q+2)-L$. Multiplication by
$e^{[D(q+2)+L]t}$ therefore gives a bounded linear system on $\ell^\infty(\N_0\times\N_0)$
with nonnegative coefficients. Its Volterra series preserves
nonnegativity. Here $d_{0,0}(0)=\varepsilon$ and all other
$d_{n,r}(0)$ vanish, so $d_{n,r}(t)\ge0$ for
$n,r\in\N_0$ and $t\ge0$.  Choose an integer $n_t$ with $Y\le n_t-m(t)\le Y+1$.
For all sufficiently large $t$, $n_t\ge b$ and
\[
 n_t-b=c_{\parallel}t-\kappa_{\parallel}\log(1+t)+O(1),
\]
where the error is bounded independently of $t$. The inequalities
$d_{n,r}\ge0$ extend $U_{n_t,0}(t)\ge\delta$ to every integer
between $b$ and $n_t$. Reflection extends the same bound to
$[2b-n_t,b]\cap\Z$. Increasing $C_0$ absorbs the bounded error in
the displayed position and proves the stated interval bound.\qedhere

\end{proof}

\begin{proof}[Proof of Theorem~\ref{thm:tangential-bramson}]

Fix $0<\theta<p_*(\gamma(0))$. For the outer bound, choose
integers $B,R\ge0$ so that every vertex
in $\operatorname{supp}u_0$ has $|\ell_{\mathcal G}(x)|\le B$ and
$s(x)\le R$. Let $\overline U$ solve the nonlinear coordinate equation
in Lemma~\ref{lem:tangential-positive-seed}, with initial datum
\[
 \overline U_{n,r}(0)=\mathbf1_{\{|n|\le B,\ 0\le r\le R\}},
 \qquad(n,r)\in\Z\times\N_0.
\]
The lifted datum $x\mapsto\overline U_{\ell_{\mathcal G}(x),s(x)}(0)$
dominates $u_0$, so comparison gives
$u(t,x)\le\overline U_{\ell_{\mathcal G}(x),s(x)}(t)$ for
$t\ge0$ and $x\in\Tq$. Moreover,
$\overline U_{n,r}(0)=0$ for $n>m_B(0)=B$, as required by the
stopped moving-domain estimate. We first estimate sites with $n>m_B(t)$,
where $m_B$ is defined at the start of
Subsection~\ref{sec:tangential-bramson}.  Since $m_B'(t)>0$,
a backward path can leave the active region only through a horizontal
jump from height in $(0,1]$; there are no continuous boundary contacts
during a transverse excursion.

Fix $t\ge0$ and $n\in\Z$ with $n>m_B(t)$, and apply
\eqref{eq:tangential-stopped-comparison} to $\overline U$ at $(n,0)$.
For a bounded function $F$ on $\Z\times\N_0$,
Lemma~\ref{lem:tangential-regeneration} gives the pointwise identity
$(\mathcal A_{\parallel}-\Lambda_{\parallel}(\lambda))(H_\lambda F)
=H_\lambda\mathcal Q_{\parallel}^\lambda F$, where
$H_\lambda(n,r)=e^{-\lambda n}h_r^\lambda$.
We use it at $\lambda=\lambda_{\parallel}$ to express the stopped
expectation under $\mathbb P^{\lambda_{\parallel}}_{(n,0)}$.
At an exit the transverse state is zero, so
$h_0^{\lambda_{\parallel}}/h_0^{\lambda_{\parallel}}=1$.
Let $\tau=\inf\{s\in[0,t]:N_s\le m_B(t-s)\}$, with
$\inf\varnothing=\infty$, and put
\[
 y=n-m_B(t)>0.
\]
On $\{\tau\le t\}$, write
$\zeta_\tau=N_\tau-m_B(t-\tau)\in(-1,0]$.
By the eigenfunction change of measure in
\eqref{eq:tangential-stopped-comparison}, the exponential factor is
\[
 e^{\Lambda_{\parallel}(\lambda_{\parallel})\tau
       -\lambda_{\parallel}(n-N_\tau)}
 =e^{-\lambda_{\parallel}y+\lambda_{\parallel}\zeta_\tau}
   \left(\frac{1+a_*+t}{1+a_*+t-\tau}\right)^{3/2}.
\]
Since $e^{\lambda_{\parallel}\zeta_\tau}\le1$, the stopped estimate
for $\overline U$ and the order $u(t,\gamma(n))\le\overline U_{n,0}(t)$ give
\[
 u(t,\gamma(n))
 \le\overline U_{n,0}(t)
 \le
 Ce^{-\lambda_{\parallel}y}
 \mathbb E^{\lambda_{\parallel}}_{(n,0)}
 \left[
 \left(
 \frac{1+a_*+t}
 {1+a_*+t-\tau}
 \right)^{3/2};
 \tau\le t
 \right].
\]
An exit by a left horizontal jump at time $s$ has pre-jump height in
$(0,1]$. Put $\xi_s(k)=k-m_B(t-s)$. The left-jump rate is
$De^{\lambda_{\parallel}}$, so marking this jump bounds its density by
\[
 De^{\lambda_{\parallel}}
 \sum_{\substack{k\in\Z\\0<\xi_s(k)\le1}}
 \mathcal K_{s,a_*+t-s}^{y,\lambda_{\parallel}}(n,k)\,\dd s.
\]
There are at most two terms. The kernel's terminal height is
$\xi_s(k)$, and its boundary is
$b_{s,a_*+t-s}(v)
=\kappa_{\parallel}\log((1+a_*+t)/(1+a_*+t-v))$ for
$0\le v\le s$. Thus its offset satisfies $a_*+t-s\ge a_*$,
and Lemma~\ref{lem:tangential-backward-log} gives
\begin{equation}\label{eq:tan-log-exit-density}
 \mathbb P^{\lambda_{\parallel}}_{(n,0)}(\tau\in ds)
 \le
 C(1+y)(1+s)^{-3/2}\,ds,\qquad 0<s\le t.
\end{equation}
Hence
\[
 u(t,\gamma(n))
 \le
 C(1+y)e^{-\lambda_{\parallel}y}
 (1+a_*+t)^{3/2}
 \int_0^t
 \frac{ds}
 {(1+s)^{3/2}(1+a_*+t-s)^{3/2}}.
\]
Write $A=1+a_*>0$. On $[0,t/2]$, one has
$A+t-s\ge(A+t)/2$, while on $[t/2,t]$ one has
$1+s\ge(1+t)/2$. Therefore
\[
 \begin{aligned}
 \int_0^t\frac{\dd s}{(1+s)^{3/2}(A+t-s)^{3/2}}
 &\le 2^{3/2}(A+t)^{-3/2}
           \int_0^\infty(1+s)^{-3/2}\,\dd s\\
 &\quad+2^{3/2}(1+t)^{-3/2}
           \int_0^\infty(A+r)^{-3/2}\,\dd r\\
 &\le C_{a_*}(A+t)^{-3/2}.
 \end{aligned}
\]
Consequently,
\[
 u(t,\gamma(n))\le\overline U_{n,0}(t)
 \le
 C(1+y)e^{-\lambda_{\parallel}y}.
\]
Choose $C_\theta$ so that
$C(1+y)e^{-\lambda_{\parallel}y}<\theta$ for every $y\ge C_\theta$.
Then $u(t,\gamma(n))<\theta$ whenever $n\ge m_B(t)+C_\theta$.
The solution $\overline U$ has initial datum
$\mathbf1_{\{|n|\le B,\ 0\le r\le R\}}$, invariant under
$n\mapsto-n$. This map interchanges the two horizontal diffusion
terms and leaves the transverse and reaction terms unchanged.
Uniqueness makes $\overline U$ even in $n$. The upper bound for $\overline U$ on
$n\ge m_B(t)+C_\theta$, together with comparison with $u$, thus
gives $u(t,\gamma(n))<\theta$ for
$n\le-m_B(t)-C_\theta$. Moreover,
\[
 m_B(t)-m_{\parallel}(t)
 =B+\kappa_{\parallel}\log
   \frac{(1+a_*)(1+t)}{1+a_*+t}
 \quad\text{for }t\ge0.
\]
This difference lies in $[B,B+\kappa_{\parallel}\log(1+a_*)]$.
Increasing $C_\theta$ therefore gives the outer inclusion centered
at $m_{\parallel}(t)$.

For the inner bound, strict positivity at a fixed time $t_0>0$ gives
an integer $b$ and an $\varepsilon>0$ such that
\[
 u(t_0,\cdot)
 \ge
 \varepsilon\mathbf1_{\{\gamma(b)\}}.
\]
Decrease $\varepsilon$, if necessary, so that
Lemma~\ref{lem:tangential-positive-seed} applies.
Let $U$ solve the nonlinear coordinate equation in
Lemma~\ref{lem:tangential-positive-seed} with initial datum
$\varepsilon\mathbf1_{\{(b,0)\}}$, whose tree lift is
$\varepsilon\mathbf1_{\{\gamma(b)\}}$. Comparison gives
$u(t_0+t,\gamma(n))\ge U_{n,0}(t)$ for every $n$ and $t\ge0$.
Lemma~\ref{lem:tangential-positive-seed} gives
\[
 U_{n,0}(t)\ge\delta
\]
throughout
\[
 |n-b|
 \le
 m_{\parallel}(t)-C_0
\]
for all large $t$.

We now raise the fixed positive level $\delta$ to the prescribed
$\theta$. Let $Z(t,r)$ solve the interior equation in
\eqref{eq:radial} and the boundary equation
\eqref{eq:geodesic-boundary}, with initial datum
\[
 Z(0,r)=\delta\mathbf1_{\{r=0\}},\qquad r\in\N_0.
\]
The lift $x\mapsto Z(t,s(x))$ has nonzero nonnegative initial datum.
By Theorem~\ref{thm:geodesic-threshold}, the solution $Z$ satisfies
\[
 \liminf_{t\to\infty}Z(t,0)
 \ge
 p_*(\gamma(0)).
\]
Since $\theta<p_*(\gamma(0))$, choose a fixed $\tau>0$ such that
$Z(\tau,0)>\theta$.
For an integer $J\ge0$, let $Z^J(t,n,r)$ solve
$\partial_tZ^J_{n,r}=(\mathcal A_{\parallel}Z^J)_{n,r}
+f_r(Z^J_{n,r})-a_rZ^J_{n,r}$, with initial datum
\[
 Z^J(0,n,r)=\delta\mathbf1_{\{|n|\le J,\ r=0\}},
 \qquad(n,r)\in\Z\times\N_0.
\]
By comparison, $Z(\tau,r)-Z^J(\tau,n,r)\ge0$.
The Lipschitz bound $|f_r'|\le L$ in \eqref{eq:taylor} bounds this
difference by $e^{L\tau}$ times the coordinate diffusion semigroup
applied to $\delta\mathbf1_{\{|n|>J,\ r=0\}}$.  For $n\in\Z$ with $|n|\le J-R$, a path from $(n,0)$ to
$\{(k,0):|k|>J\}$ requires at least $R$ jumps.  Since the diffusion chain has total rate $D(q+1)$ and the
initial difference is bounded by $\delta$, we obtain
\[
 0\le Z(\tau,0)-Z^J(\tau,n,0)
 \le\delta e^{L\tau}
 \mathbb P\{\operatorname{Poisson}(D(q+1)\tau)\ge R\}
 \quad\text{for }J\ge R,\ |n|\le J-R.
\]
Choose an integer $R\ge1$, independent of $J$, so that the displayed
Poisson-tail bound is below $Z(\tau,0)-\theta$. Then
$Z^J(\tau,n,0)>\theta$ whenever $J\ge R$ and $|n|\le J-R$.

At a large time $t$, compare $U(t,\cdot,\cdot)$ with the datum
$\delta\mathbf1_{\{|n-b|\le J,\ r=0\}}$, where
\[
 J=
 \left\lfloor
 m_{\parallel}(t)-C_0
 \right\rfloor.
\]
The lower bound for $U_{n,0}(t)$ gives
$U_{n,r}(t)\ge Z^J(0,n-b,r)$ for every $(n,r)\in\Z\times\N_0$.
Translation invariance in $n$ and comparison over $[t,t+\tau]$
therefore yield
\[
 U_{n,0}(t+\tau)\ge Z^J(\tau,n-b,0)>\theta
 \quad\text{for }n\in\Z,\ |n-b|\le J-R,
\]
for all sufficiently large $t$ that $J\ge R$.
For $t,h\ge0$,
\[
 m_{\parallel}(t+h)-m_{\parallel}(t)
 =c_{\parallel}h-\kappa_{\parallel}
       \log\frac{1+t+h}{1+t}\le c_{\parallel}h.
\]
Increase $C_\theta$ so that
$C_\theta\ge C_0+R+1+|b|+c_{\parallel}(t_0+\tau)$.
At the original time $s=t_0+t+\tau$, the condition
$|n|\le m_{\parallel}(s)-C_\theta$ then gives
\[
 |n-b|\le m_{\parallel}(t)-C_0-R-1\le J-R.
\]
Comparison from time $t_0$ therefore gives the inner inclusion
for the original solution.

For all sufficiently large $t$, the inner inclusion gives
$L_\theta(t)\ge\lfloor m_{\parallel}(t)-C_\theta\rfloor$, and the
outer inclusion gives
$L_\theta(t)\le\lfloor m_{\parallel}(t)+C_\theta\rfloor$.
Thus $|L_\theta(t)-m_{\parallel}(t)|\le C_\theta+1$.
Since $\log(1+t)-\log t\to0$ as $t\to\infty$, the bounded-error
formula also holds with $\log t$.
\end{proof}
\begin{corollary}\label{thm:geodesic}

Assume \eqref{eq:hyp} and
\[
 f_i'(0)=a>\rho_q
 \quad\text{ for } 0\le i<P.
\]
Let the reactions be assigned on the tree by an arbitrary map
\[
 \iota:\Tq\to\{0,\ldots,P-1\},
 \qquad
 f(x,s)=f_{\iota(x)}(s),\qquad x\in\Tq,\quad s\in[0,1].
\]
Let $u$ solve \eqref{eq:tree-cauchy} with initial datum $u_0$.
Fix a bi-infinite geodesic $\mathcal G$ and suppose
\[
 \eta\ind_{\mathcal N_{R_1}(\mathcal G)}
 \le u_0
 \le
 \ind_{\mathcal N_{R_2}(\mathcal G)},
 \qquad
 0<\eta\le1,\quad
 0\le R_1\le R_2<\infty.
\]
Then, for every $\theta\in(0,1)$, there are
$C_\theta,T_\theta<\infty$ such that
\begin{equation}
\label{eq:geodesic-bramson}
 \mathcal N_{m_a(t)-C_\theta}(\mathcal G)
 \subseteq
 \{x:u(t,x)\ge\theta\}
 \subseteq
 \mathcal N_{m_a(t)+C_\theta}(\mathcal G)
 \quad\text{ for } t\ge T_\theta,
\end{equation}
where
\[
 m_a(t)
 =
 c_*(a)t-\frac{3}{2\lambda_*(a)}\log t,
 \quad\text{ for }t>0.
\]

\end{corollary}

\begin{remark}
For the assignment $f(x,u)=f_{s(x)}(u)$, the inclusions also
follow from Theorem~\ref{thm:periodic-bramson-distance}.
\end{remark}

\begin{proof}

Let $h_-,h_+$ be the scalar reactions from
Lemma~\ref{lem:envelopes}, and let $u^-,u^+$ solve the homogeneous
tree equations with reactions $h_-,h_+$ and initial data
\[
 \eta\ind_{\mathcal N_{R_1}(\mathcal G)}
 \quad\text{and}\quad
 \ind_{\mathcal N_{R_2}(\mathcal G)},
\]
respectively.  Comparison gives
\[
 u^-(t,x)\le u(t,x)\le u^+(t,x),\qquad t\ge0,\quad x\in\Tq.
\]
Because the comparison equations are homogeneous and their initial
data depend only on $s(x)$, write
$u^\pm(t,x)=U^\pm_{s(x)}(t)$. For $t>0$, their coordinate equations are
\[
 \begin{aligned}
 \dot U_r^\pm
 &=D(U_{r-1}^\pm+qU_{r+1}^\pm-(q+1)U_r^\pm)
       +h_\pm(U_r^\pm),&&\text{ for }r\ge1,\\
 \dot U_0^\pm&=D(q-1)(U_1^\pm-U_0^\pm)+h_\pm(U_0^\pm).
 \end{aligned}
\]
The initial sequences are
$U_r^-(0)=\eta\mathbf1_{\{r\le R_1\}}$ and
$U_r^+(0)=\mathbf1_{\{r\le R_2\}}$ for $r\in\N_0$.
Since $0<\eta\le1$ and $0\le R_1\le R_2<\infty$, both sequences
are nonzero, take values in $[0,1]$, and have finite support. Both reactions have derivative
$a$ at zero. Since $a>\rho_q$, Lemma~\ref{lem:geodesic-local}
gives $U_0^\pm(t)\to1$ as $t\to\infty$. These facts verify all
the hypotheses of Lemma~\ref{lem:boundary-log}, with
$b_0=D(q-1)>0$, for both comparison solutions,
and both have the same $c_*(a)$, $\lambda_*(a)$, and logarithmic
coefficient. Choose $C_\theta$ as the maximum of the two spatial
constants and $T_\theta$ as the maximum of the two time thresholds.
The lower inclusion for $U^-$ and the upper inclusion for $U^+$,
with $r=s(x)$, combine with $u^-\le u\le u^+$ to give
\eqref{eq:geodesic-bramson}.\qedhere

\end{proof}

\begin{proof}[Proof of Theorem~\ref{thm:geodesic-persistence}]
For the rates in the statement, Lemma~\ref{lem:geodesic-certificate}
gives
\begin{equation}\label{eq:example-spectral-signs}
 \mu\le-\frac{97}{9800}<0,
 \qquad \sigma_{\mathcal G}>\frac1{25}>0.
\end{equation}
The same lemma provides a connected set $\Omega$ with thirty
vertices and a positive eigenpair $(\eta,\varphi)$ for the
zero-exterior linearized operator on $\Omega$. Normalize
$\max_\Omega\varphi=1$ and extend $\varphi$ by zero. The diagonal
of this operator remains $a_{s(x)}-(q+1)$, and
$\eta>1/25$.

Fix $0<\varepsilon<\min\{1,\eta/a_{\max}\}$ and put
$z=\varepsilon\varphi$. Since $D=1$, the eigenvalue equation gives
\begin{equation}\label{eq:example-strict-subsolution}
 \Delta_{\Tq}z(x)+a_{s(x)}z(x)(1-z(x))
 =z(x)(\eta-a_{s(x)}z(x))>0
 \quad\text{ for }x\in\Omega.
\end{equation}
At $x\notin\Omega$, one has $z(x)=0$ and
$\Delta_{\Tq}z(x)=\sum_{y\sim x}z(y)\ge0$.
Thus $z$ is a stationary subsolution of the differential equation
in \eqref{eq:tree-cauchy}. For every $0<\vartheta\le1$,
\[
 \Delta_{\Tq}(\vartheta z)(x)
 +a_{s(x)}\vartheta z(x)(1-\vartheta z(x))
 =\vartheta z(x)(\eta-a_{s(x)}\vartheta z(x))>0
 \quad\text{for }x\in\Omega.
\]
For $x\notin\Omega$, the expression on the left equals
$\vartheta\sum_{y\sim x}z(y)\ge0$.

Let $u^z$ solve \eqref{eq:tree-cauchy} from $z$.
By \eqref{eq:example-strict-subsolution} and
Lemma~\ref{lem:comparison},
$z\le u^z(t,\cdot)\le1$ for $t\ge0$.
Autonomy then gives
$u^z(t+s,x)\ge u^z(t,x)$ for all $s,t\ge0$.
Consequently, $p(x)=\lim_{t\to\infty}u^z(t,x)$ exists at every
vertex. For any fixed $T>0$, the integral equation on $[t,t+T]$
and bounded convergence give
\begin{equation}\label{eq:example-stationary-limit}
 0=T\bigl[\Delta_{\Tq}p(x)+a_{s(x)}p(x)(1-p(x))\bigr]
 \quad\text{ for }x\in\Tq.
\end{equation}
There are finitely many neighbors at each vertex, and the logistic
reaction is continuous on $[0,1]$, so every term passes to the limit.
Strict positivity in Lemma~\ref{lem:comparison} gives $p(x)>0$.

{To identify $p$, let $v:\Tq\to(0,1]$ satisfy
$\Delta_{\Tq}v(x)+a_{s(x)}v(x)(1-v(x))=0$ for every $x\in\Tq$.}
The finite set $\Omega$ and the positivity of $z$ on $\Omega$ give
$\vartheta_*:=\min_{x\in\Omega}v(x)/z(x)>0$.
If $\vartheta_*<1$, choose $x_*\in\Omega$ where this minimum is
attained. Then $v-\vartheta_*z\ge0$ on the tree and vanishes at
$x_*$, so its graph Laplacian is nonnegative there. At $x_*$, the identity $v(x_*)=\vartheta_*z(x_*)$ gives
$a_{s(x_*)}v(x_*)(1-v(x_*))
=a_{s(x_*)}\vartheta_*z(x_*)(1-\vartheta_*z(x_*))$. Subtracting the strict
subsolution inequality for $\vartheta_*z$ from the stationary
equation for $v$ therefore gives
$\Delta_{\Tq}(v-\vartheta_*z)(x_*)<0$, a contradiction.
It follows that $v\ge z$. By comparison, $v\ge u^z(t,\cdot)$
for all $t\ge0$, and hence $v\ge p$. In view of
\eqref{eq:example-stationary-limit}, $p$ is itself positive and
stationary. Thus $p=p_*$, the smallest solution in
Theorem~\ref{thm:geodesic-threshold}. This argument applies to
every stated $\varepsilon$, which supplies arbitrarily small
finitely supported initial data converging increasingly to $p_*$.

For the transverse decay, put
$h_r=2^{-r/2}\phi_r(\lambda_0)$ for $r\in\N_0$.
By \eqref{eq:example-spectral-signs},
Theorem~\ref{thm:geodesic-threshold}(iii) gives
\begin{equation}\label{eq:example-stationary-bound}
 0<p_*(x)<1,\qquad
 p_*(x)\le\frac{h_{s(x)}}{h_0}
 \quad\text{ for }x\in\Tq.
\end{equation}
The periodic vector $\phi(\lambda_0)$ is bounded, and hence
\eqref{eq:example-stationary-bound} proves
\eqref{eq:localized-steady}.

Finally, let $u_0$ be any nonzero finitely supported datum in
$[0,1]$. By strict positivity at time one and finiteness of
$\Omega$, there is $0<\vartheta\le1$ with
$u(1,x)\ge\vartheta z(x)$ at every vertex. Since $\vartheta z=(\vartheta\varepsilon)\varphi$ and
$0<\vartheta\varepsilon<\min\{1,\eta/a_{\max}\}$, the argument
above applied with $\vartheta\varepsilon$ in place of $\varepsilon$
shows that the solution from $\vartheta z$ increases to $p_*$. Comparison therefore yields
\begin{equation}\label{eq:example-persistence-limit}
 \liminf_{t\to\infty}u(t,x)\ge p_*(x)>0
 \quad\text{for each fixed }x\in\Tq.
\end{equation}
For $i\in\{0,\ldots,P-1\}$ and $0<s\le1$,
$f_i(s)/s=a_i(1-s)$ is strictly decreasing because $a_i>0$; hence
\eqref{eq:strict-quotient} holds. Let $\lambda_{\mathrm R}$ be the
larger positive root of $k(\lambda)=0$. For each fixed
$\lambda\in(\lambda_0,\lambda_{\mathrm R})$,
Theorem~\ref{thm:geodesic-convergence} gives
$\sup_{\{x:|\ell_{\mathcal G}(x)|\le M\}}|u(t,x)-p_*(x)|\to0$
as $t\to\infty$ for every $0\le M<\infty$, and uniqueness of
$p_*$ among stationary solutions $0<p\le1$ satisfying
$\sup_{x\in\Tq}e^{\lambda s(x)}p(x)<\infty$.\qedhere

\end{proof}

\subsection{Data with infinite support}

\begin{proposition}\label{prop:horo-slices}

Assume \eqref{eq:hyp} and $f(x,u)=f_{b(x)}(u)$. Let $B\in\R$
and $0\le u_0\le\ind_{\{b\le B\}}$. Let $u$ solve
\eqref{eq:tree-cauchy} with initial datum $u_0$. Then, for every
$\nu\in\R$ and $t\ge0$,
\begin{equation}
\label{eq:horo-slice-decay}
 \sup_{\{x\in\Tq:\,b(x)\ge\nu\}}u(t,x)
 \le C(1+t)^{-1/2}e^{\mu t-\lambda_0(\nu-B)}.
\end{equation}
{The constant $C$ in \eqref{eq:horo-slice-decay} is independent of
$t$, $\nu$, and $B$. For each fixed $\nu$, the supremum}
$\sup_{b(x)\ge\nu}u(t,x)$ tends to zero as $t\to\infty$ when
$\mu\le0$. If $\mu>0$ and
$u_0(x)\ge v_0(b(x))$ for every $x\in\Tq$, where
$v_0:\Z\to[0,1]$ and $v_0\not\equiv0$, then
\[
 \inf_{\{x\in\Tq:\,|b(x)|\le M\}}u(t,x)\to1\quad\text{ as } t\to\infty
 \quad\text{ for every }M\ge0.
\]
The sets $\{x:|b(x)|\le M\}$ and $\{x:b(x)\ge\nu\}$
are unbounded in graph distance whenever they are nonempty.

\end{proposition}
\begin{proof}

Let $w$ solve \eqref{eq:lattice} with
$w_j(0)=\mathbf1_{\{j\le B\}}$ for $j\in\Z$.
The function $(t,x)\mapsto w_{b(x)}(t)$ solves the differential
equation in \eqref{eq:tree-cauchy} and initially dominates $u_0$.
Lemma~\ref{lem:comparison} therefore gives
$u(t,x)\le w_{b(x)}(t)$ for $t\ge0$ and $x\in\Tq$.
Since $\lambda_0=\tfrac12\log q>0$,
\[
 \sum_{j\le B}e^{\lambda_0j}
 =\frac{e^{\lambda_0\lfloor B\rfloor}}{1-e^{-\lambda_0}}
 \le\frac{e^{\lambda_0B}}{1-e^{-\lambda_0}}.
\]
Lemma~\ref{lem:linear-bound}, with $k(\lambda_0)=\mu$, gives
$w_j(t)\le C(1+t)^{-1/2}e^{\mu t-\lambda_0(j-B)}$.
Taking the supremum over $j\ge\nu$ proves
\eqref{eq:horo-slice-decay}, with a constant independent of $B,\nu,t$.
{For fixed $\nu$ and $B$, the factor $(1+t)^{-1/2}e^{\mu t}$
in \eqref{eq:horo-slice-decay} tends to zero as $t\to\infty$
when $\mu\le0$, including $\mu=0$.} For the lower assertion, Lemma~\ref{lem:matrix} gives
$c_-<0<c_+$ from $\mu>0$. Thus zero belongs to the speed interval
in Lemma~\ref{lem:bump}; apply that lemma with $c=0$ and $R=P$.
Write $v_j=v(0,j)$ for the time-independent subsolution supplied
by that lemma. It is compactly supported, with
$v_j\ge\delta>0$ for $|j|\le P$, such that every
$\vartheta v$, $0<\vartheta\le1$, is a strict subsolution on its
positive set.

Let $W:[0,\infty)\times\Z\to[0,1]$ solve
\eqref{eq:lattice} with $W(0,\cdot)\not\equiv0$. By
Lemma~\ref{lem:comparison}, $W(1,j)>0$ for every $j\in\Z$.
The set $\{j\in\Z:v_j>0\}$ is finite and nonempty, so
\[
 \vartheta=\min\left\{1,
       \min_{\substack{j\in\Z\\v_j>0}}\frac{W(1,j)}{v_j}\right\}>0.
\]
Let $U$ solve \eqref{eq:lattice} with $U_j(0)=\vartheta v_j$.
Comparison and autonomy give $U_j(t+s)\ge U_j(t)$ for
$j\in\Z$ and $s,t\ge0$, while $0\le U_j(t)\le1$.
Thus $p_j=\lim_{t\to\infty}U_j(t)$ exists for every $j\in\Z$.
Passing to the limit in the coordinate integral equation on $[t,t+1]$
by bounded convergence gives
\[
 D(p_{j-1}+qp_{j+1}-(q+1)p_j)+f_j(p_j)=0
 \quad\text{for }j\in\Z.
\]
Moreover, $p_j\ge U_j(1)>0$ by strict positivity.
For every integer $k$, put
$v^{(k)}_j=v_{j-kP}$ for $j\in\Z$. Since $f_{j-kP}=f_j$,
every $\vartheta v^{(k)}$, $0<\vartheta\le1$, is a strict
stationary subsolution on $\{j:v^{(k)}_j>0\}$. If
\[
 s_k:=\min\left\{1,\min_{j:v^{(k)}_j>0}
                         \frac{p_j}{v^{(k)}_j}\right\}<1,
\]
the positive support of $v^{(k)}$ is finite and $p_j>0$ there,
so $0<s_k<1$ and the minimum is attained at an index $j_k$ with
$v^{(k)}_{j_k}>0$ and $p_{j_k}=s_kv^{(k)}_{j_k}$.
{Put $\zeta_j=p_j-s_kv^{(k)}_j$. Then $\zeta_j\ge0$ for every $j\in\Z$
and $\zeta_{j_k}=0$, so
\[
 D\bigl(\zeta_{j_k-1}+q\zeta_{j_k+1}-(q+1)\zeta_{j_k}\bigr)\ge0.
\]
Moreover, $f_{j_k}(p_{j_k})=f_{j_k}(s_kv^{(k)}_{j_k})$.
The stationary equation for $p$ and the strict subsolution inequality
for $s_kv^{(k)}$ would therefore give
\[
 \begin{aligned}
 0&=D\bigl(p_{j_k-1}+qp_{j_k+1}-(q+1)p_{j_k}\bigr)
       +f_{j_k}(p_{j_k})\\
 &\ge D\bigl(s_kv^{(k)}_{j_k-1}+qs_kv^{(k)}_{j_k+1}
                 -(q+1)s_kv^{(k)}_{j_k}\bigr)
       +f_{j_k}(s_kv^{(k)}_{j_k})>0,
 \end{aligned}
\]
which is a contradiction.} Thus $p\ge v^{(k)}$ for every $k$.
For every $j\in\Z$, choose $k\in\Z$ with $|j-kP|\le P$.
Then $p_j\ge v^{(k)}_j=v_{j-kP}\ge\delta$. The stationary solution
$(t,j)\mapsto p_j$ is therefore an entire solution of \eqref{eq:lattice}
with $\delta\le p_j\le1$. Lemma~\ref{lem:liouville} gives $p\equiv1$.
Since $W(1,j)\ge U_j(0)$, Lemma~\ref{lem:comparison} gives
$U_j(t)\le W(1+t,j)\le1$ for $j\in\Z$ and $t\ge0$.
The limit $U_j(t)\to p_j=1$ therefore gives convergence of $W$
to one on every finite set of indices. In particular, take $W(0,j)=v_0(j)$ for $j\in\Z$.
The initial inequality $u_0(x)\ge v_0(b(x))$ and comparison give
$u(t,x)\ge W(t,b(x))$ for $t\ge0$ and $x\in\Tq$.
For each $M\ge0$, the set $\{j\in\Z:|j|\le M\}$ is finite,
so $W(t,j)\to1$ uniformly on that set. This proves the asserted
uniform convergence on $\{x\in\Tq:|b(x)|\le M\}$.

For an integer $j$, take large integers $k$ with $k+j\ge1$.
Starting at the ray vertex $o_k$, follow a path of length $k+j$
along which $b$ increases by one at each step, choosing the first
edge outside the distinguished ray; such a path exists since
$q-1\ge1$. Its endpoint $x_k$ has $b(x_k)=j$ and
$d(o,x_k)=2k+j\to\infty$ as $k\to\infty$.
The sets $\{|b|\le M\}$ and $\{b\ge\nu\}$ contain the horocycles
$\{b=0\}$ and $\{b=\lceil\nu\rceil\}$, respectively, and are
therefore unbounded in graph distance.
\end{proof}

\begin{proposition}\label{prop:lp-restored}
Let $\iota:\Tq\to\{0,\ldots,P-1\}$ and
$f(x,s)=f_{\iota(x)}(s)$, where the family $(f_i)$ satisfies
\eqref{eq:hyp}. Let $u$ solve \eqref{eq:tree-cauchy} with
$0\le u_0\le1$ and $u_0\in\ell^p(\Tq)$, where $1\le p<\infty$.
Then, for every $t\ge0$,
\begin{equation}
\label{eq:lp-restored}
 \|u(t,\cdot)\|_\infty\le
 C_p(1+t)^{-3/(2p)}e^{(a_{\max}-\gamma_p)t}\|u_0\|_p.
\end{equation}
{Here $C_p$ depends only on $p,q,D$, and
$\gamma_p=\rho_q$ for $1\le p\le2$ and}
$\gamma_p=2\rho_q/p$ for $2\le p<\infty$, as in
Lemma~\ref{lem:tree-lp}. For a common derivative
$a_i=a$, replace $a_{\max}$ by $a$. In particular, at $a=\rho_q$ the
bounds for $\ell^1$ and $\ell^2$ data are respectively
$O(t^{-3/2})$ and $O(t^{-3/4})$ as $t\to\infty$. 
\end{proposition}
\begin{proof}
By the KPP bound in \eqref{eq:hyp},
$f_{\iota(x)}(s)\le a_{\max}s$ for $0\le s\le1$.
Lemma~\ref{lem:comparison} therefore gives
$0\le u(t,\cdot)\le e^{a_{\max}t}P_tu_0$, where
$P_t=e^{tD\Delta_{\Tq}}$ is the tree heat semigroup.
Applying Lemma~\ref{lem:tree-lp} to $u_0\in\ell^p(\Tq)$ proves
\eqref{eq:lp-restored}.
\end{proof}

\subsection{Two-sided kernel bounds at a straight line}

Theorem~\ref{thm:killed-two-sided} combines
Lemmas~\ref{lem:killed-upper} and~\ref{lem:killed-straight-lower}.
Throughout,
$\mathbb P_j^\lambda$ is the continuous-time path law defined in
Subsection~\ref{sec:periodic-chain}, $c=k'(\lambda)$, and a path condition imposed for
$0\le s\le t$ includes the times between successive jumps.

\begin{theorem}\label{thm:killed-two-sided}
Fix $\lambda\in\R$ and use the killed kernel in
\eqref{eq:killing-definition}, with $c=k'(\lambda)$. There exists
$y_0\ge1$ such that, for each $A>0$, there are $C_A\ge1$ and $t_A<\infty$
for which
\begin{equation}
 \label{eq:killed-two-sided}
 C_A^{-1}\frac{yz}{t^{3/2}}
 \le K_t^y(j,l)\le C_A\frac{yz}{t^{3/2}}
\end{equation}
whenever $t\ge t_A$, $j,l\in\Z$, and
\[
 y_0\le y\le A\sqrt t,\qquad
 y_0\le z:=y+l-j+ct\le A\sqrt t.
\]
Constants are uniform in the initial and terminal spatial phases.
The boundary is the straight line in
\eqref{eq:killing-definition}.
\end{theorem}
\begin{proof}
Choose $y_0$ from Lemma~\ref{lem:killed-straight-lower} and increase
$t_A$ to at least one. For $y,z\ge y_0\ge1$,
\[
 \frac{(1+y)(1+z)}{(1+t)^{3/2}}
 \le\frac{4yz}{t^{3/2}}.
\]
Thus Lemma~\ref{lem:killed-upper} gives the required upper bound,
and Lemma~\ref{lem:killed-straight-lower} gives the lower bound on
the same endpoint range. Taking $C_A$ at least as large as the
upper-bound constant and the reciprocal of the lower-bound constant
proves the assertion uniformly in the spatial phases.
\end{proof}

\begin{remark}

For each compact set $K\subseteq[y_0,\infty)$, there are
$c_K,C_K>0$ and $T_K<\infty$ such that
\[
 c_Kt^{-3/2}\le K_t^y(j,l)\le C_Kt^{-3/2}
\]
for $t\ge T_K$ and $j,l\in\Z$ whenever
$y,z:=y+l-j+ct\in K$. This follows from
\eqref{eq:killed-two-sided}, since $y,z\ge y_0>0$ and are bounded above
on $K$.

\end{remark}

\subsection{A logarithmic boundary}

For a time horizon $t>0$, we use two deterministic boundaries on $[0,t]$:
$s\mapsto\kappa\log(1+\min\{s,t-s\})$, which is zero at both
endpoints, and $s\mapsto\kappa\log(1+s)$. The parameter satisfies
$\kappa\ge0$, and the endpoint heights are measured from the
particular boundary in each statement. For random walks with
independent identically distributed increments, first-passage times
above deterministic moving boundaries are studied in
\cite{DenisovSakhanenkoWachtel}. For the continuous-time chain with
law $\mathbb P_j^\lambda$, the jump rates are periodic in the spatial
index, and absorption includes crossings between jumps. The proofs
use the straight-line kernel bounds in
Theorem~\ref{thm:killed-two-sided} and the parameter-change identity
\eqref{eq:path-tilt-identity}.

\begin{theorem}\label{thm:logarithmic-arch}

Fix $\lambda\in\R$, $\kappa\ge0$ and $A>0$, and put $c=k'(\lambda)$.
For $t>0$, $j,l\in\Z$ and $y>0$, define
$g_t(s)=\kappa\log(1+\min\{s,t-s\})$ for $0\le s\le t$ and
\[
 K_{t,\kappa}^{y,\mathrm{arch}}(j,l)
 =\mathbb P_j^\lambda\{X_t=l,\ y+X_s-j+cs>g_t(s)
                                  \text{ for every }0\le s\le t\}.
\]
There are $y_{A,\kappa}\ge1$, $t_{A,\kappa}<\infty$ and
$C_{A,\kappa}\ge1$ such that
\begin{equation}
 \label{eq:log-arch-bounds}
 C_{A,\kappa}^{-1}\frac{yz}{t^{3/2}}
 \le K_{t,\kappa}^{y,\mathrm{arch}}(j,l)
 \le C_{A,\kappa}\frac{yz}{t^{3/2}}
\end{equation}
whenever $t\ge t_{A,\kappa}$ and
$y_{A,\kappa}\le y,z:=y+l-j+ct\le A\sqrt t$.
The constants can be chosen uniformly for $\lambda$ in a compact interval.

\end{theorem}
\begin{proof}

The case $\kappa=0$ is Theorem~\ref{thm:killed-two-sided}. Suppose
$\kappa>0$. Since $g_t(s)\ge0$ for $0\le s\le t$, every path counted
by $K_{t,\kappa}^{y,\mathrm{arch}}(j,l)$ also survives above
the straight line. Lemma~\ref{lem:killed-upper} therefore gives
the upper bound. For the lower bound we first restrict paths to stay
above a fixed height $H>0$, which will be chosen large. Let $y_0$ and
$t_A$ be the height and time thresholds in
Theorem~\ref{thm:killed-two-sided}, and put $s_H=e^{H/\kappa}-1$.
In the following estimates assume $t\ge\max\{1,t_A,2s_H\}$ and
$2H+y_0\le y,z\le A\sqrt t$.
Theorem~\ref{thm:killed-two-sided} gives
\begin{equation}
 \label{eq:elevated-line-lower}
 \mathbb P_j^\lambda\{X_t=l,\ y+X_s-j+cs>H\text{ for }0\le s\le t\}
 =K_t^{y-H}(j,l)\ge c_A\frac{yz}{t^{3/2}}.
\end{equation}
The lower bound in Theorem~\ref{thm:killed-two-sided} applies to
$y-H,z-H\in[y_0,A\sqrt t]$. Since $y,z\ge2H+y_0$,
$(y-H)(z-H)\ge yz/4$, and $c_A$ can be chosen independently of $H$.

We estimate the paths in \eqref{eq:elevated-line-lower} that meet or
cross $g_t$. Such a crossing can occur only when
$\min\{s,t-s\}\ge s_H$. Write $Y_s=y+X_s-j+cs$ for the distance from the straight line.
A first time at which $Y_s\le g_t(s)$ is either a downward jump
of $X$ or a time of equality $Y_s=g_t(s)$ between jumps. At a downward jump
$m\to m-1$ at time $s$, the pre-jump distance must
lie in $(g_t(s),g_t(s)+1]$. There are at most two possible integers
$m$ with this property. Put $Y_s(m)=y+m-j+cs$. By the jump-rate
formula and the Markov property, the contribution of all these jumps
is at most
\begin{equation}
 \label{eq:arch-jump-loss}
 \int_{s_H}^{t-s_H}\sum_{\substack{m\in\Z\\
                  g_t(s)<Y_s(m)\le g_t(s)+1}}
 K_s^{y-H}(j,m)\,\alpha_m
 K_{t-s}^{Y_s(m)-1-H}(m-1,l)\,\dd s.
\end{equation}
Here $\alpha_m=D e^\lambda\phi_{m-1}(\lambda)/\phi_m(\lambda)$
is the downward jump rate under $\mathbb P_j^\lambda$.
Every killed-kernel entry with nonpositive initial distance is
interpreted as zero, since its survival event already fails at time
zero. To obtain \eqref{eq:arch-jump-loss}, mark one downward jump
in the finite-path expansion: its factor is $\alpha_m\,\dd s$, and the
portions before and after it have the displayed killed kernels.
Summing over the marked jumps gives an upper bound, since every path
with a downward crossing contributes at least once. Bounded rates and
nonnegativity justify the integrals; the jump-time construction is
given in \cite[Section~2.6]{Norris}.

Between jumps, the distance is affine. A contact with $g_t$ then has
a deterministic time in
\[
 \mathcal C_t=\{s\in[s_H,t-s_H]:
                      j-y+g_t(s)-cs\in\Z\}.
\]
For $\kappa>0$, $g_t(s)-cs$ is strictly concave on each half-interval,
and its derivative has magnitude at most $\kappa+|c|$. Splitting
at $t/2$ and at zeros of the derivative gives at most four strictly
monotone pieces. On each piece within a unit time interval, its range
has length at most $\kappa+|c|$, so it meets at most
$\lceil\kappa+|c|\rceil+1$ values of any translate of $\Z$. Including the finitely many splitting
points yields a bound on the number of contacts per unit interval
independent of $t,j,y$. At a contact time $s\in\mathcal C_t$,
the contact site is the integer $m=j-y+g_t(s)-cs$.
The Markov property
at this deterministic time gives the following probability for paths
that visit $m$ at time $s$, remain above $H$, and end at $l$:
\[
 K_s^{y-H}(j,m)K_{t-s}^{g_t(s)-H}(m,l).
\]
We include a possible contact at $t/2$ once. Exact jump times have no
atom at a prescribed positive time; thus the jump and continuous
contacts together cover the event being estimated.

Apply Lemma~\ref{lem:killed-upper} to both kernels in each product.
Every endpoint height in these kernels is obtained by subtracting
$H$ from its height above zero. Omitting these subtractions
increases each factor in the upper bound of
Lemma~\ref{lem:killed-upper}.
For every summand in \eqref{eq:arch-jump-loss} and every product
at a time in $\mathcal C_t$, the resulting bound is
\[
 C\frac{yz(2+g_t(s))^2}
              {(1+s)^{3/2}(1+t-s)^{3/2}}.
\]
The constant is independent of $H$. If $t\ge2s_H$, integrate this
bound and sum it over $\mathcal C_t$. The uniform number of contacts
per unit interval and symmetry under $s\mapsto t-s$ give total loss
at most
\begin{equation}
 \label{eq:arch-total-loss}
 C\frac{yz}{t^{3/2}}\,\omega_\kappa(s_H),\qquad
 \omega_\kappa(R)=\sum_{n=\max\{0,\lfloor R\rfloor-2\}}^\infty
                 \frac{(2+\kappa\log(3+n))^2}{(1+n)^{3/2}},\quad R\ge0.
\end{equation}
For $0\le s\le t/2$, the summand before
\eqref{eq:arch-total-loss} satisfies
\[
 \frac{t^{3/2}(2+g_t(s))^2}
 {(1+s)^{3/2}(1+t-s)^{3/2}}
 \le 2^{3/2}
       \frac{(2+\kappa\log(1+s))^2}{(1+s)^{3/2}}.
\]
On each unit interval $[n,n+1]$ this is at most
$2^{3/2}(2+\kappa\log(3+n))^2/(1+n)^{3/2}$.
The integral over that interval has length at most one, and the
number of points of $\mathcal C_t$ there is uniformly bounded by
the contact count proved above. Since a contributing interval has
$n\ge\max\{0,\lfloor s_H\rfloor-2\}$, summing gives at most
$C\omega_\kappa(s_H)$. On $[t/2,t]$, put $r=t-s$;
then $g_t(t-r)=\kappa\log(1+r)$ and the two denominator factors
are interchanged, so the identical unit-interval bound applies.
Multiplication by $yz/t^{3/2}$ gives
\eqref{eq:arch-total-loss}. The series defining $\omega_\kappa$
converges, and hence $\omega_\kappa(R)\to0$ as $R\to\infty$.

Choose $H$ large enough that $C\omega_\kappa(s_H)<c_A/2$,
where $C$ is the constant in \eqref{eq:arch-total-loss}. Subtract
\eqref{eq:arch-total-loss} from \eqref{eq:elevated-line-lower} and
{choose $y_{A,\kappa}=2H+y_0$. The lower bound then holds for
$t\ge\max\{1,t_A,2s_H\}$, with $t_A$ as in
Theorem~\ref{thm:killed-two-sided}.} On a compact parameter interval,
Lemma~\ref{lem:killed-straight-lower} supplies common $y_0,c_A,t_A$,
and Lemma~\ref{lem:uniform-tilt} supplies a common killed upper-bound
constant. The contact count is uniform because $k'(\lambda)$ is
bounded there. Hence $H$, $y_{A,\kappa}$, and the final time and
kernel constants can all be chosen uniformly on that interval.
\end{proof}

\begin{theorem}\label{thm:positive-log-boundary}

Fix $\lambda\in\R$, $\kappa\ge0$ and $A>0$, and let $c=k'(\lambda)$.
For $t>0$, $j,l\in\Z$ and $y>0$, define
\[
 L_{t,\kappa}^y(j,l)
 =\mathbb P_j^\lambda\{X_t=l,\ y+X_s-j+cs>\kappa\log(1+s)
                                  \text{ for every }0\le s\le t\}.
\]
There are constants $Y\ge1$, $T<\infty$ and $C\ge1$, depending on
$\lambda,\kappa,A$ and the fixed equation, such that
\begin{equation}
 \label{eq:positive-log-boundary}
 C^{-1}\frac{yz}{t^{3/2}}
 \le L_{t,\kappa}^y(j,l)\le C\frac{yz}{t^{3/2}}
\end{equation}
for $t\ge T$, $j,l\in\Z$, and
\[
 Y\le y\le A\sqrt t,\qquad
 Y\le z:=y+l-j+ct-\kappa\log(1+t)\le A\sqrt t.
\]
The path inequality is imposed at every real time $s\in[0,t]$.

\end{theorem}
\begin{proof}

The assertion at $\kappa=0$ is Theorem~\ref{thm:killed-two-sided}. For fixed
$\kappa>0$ and sufficiently large $t$, put
\[
 d_t=\frac{\kappa\log(1+t)}t,\qquad c_t=c-d_t,
 \qquad k'(\lambda_t)=c_t.
\]
The inverse function theorem and $k''(\lambda)>0$ give a unique
$\lambda_t$ near $\lambda$, with
$|\lambda_t-\lambda|\le C\log(1+t)/t$. Put $Z_s=y+X_s-j+c_t s$. In these coordinates, the boundary
in the path inequality is
\[
 r_t(s)=\kappa\left(\log(1+s)-\frac{s}{t}\log(1+t)\right),
 \qquad 0\le s\le t.
\]
The logarithm lies above its chord on $[0,t]$. Therefore,
for $0\le s\le t$,
\begin{equation}
 \label{eq:log-chord-sandwich}
 0\le r_t(s)\le\kappa\log(1+\min\{s,t-s\}).
\end{equation}
For $s\le t/2$, one has
$r_t(s)\le\kappa\log(1+s)$.
For $s\ge t/2$, put $v=t-s$. The bounds
$\log(1+s)\le\log(1+t)$ and
$\log(1+v)\ge(v/t)\log(1+t)$ give
\[
 r_t(s)\le\kappa\frac vt\log(1+t)
       \le\kappa\log(1+v),
\]
which proves the other half of \eqref{eq:log-chord-sandwich}.

Under $\mathbb P_j^{\lambda_t}$ the chain has velocity $-c_t$,
and $Z_0=y$, $Z_t=z$ on $\{X_t=l\}$. The straight-line and arch
estimates apply for either sign of $c_t$, since they use the centered
coordinate $Z_s=y+X_s-j+c_t s$. By
\eqref{eq:log-chord-sandwich}, the path events satisfy
\[
 \begin{aligned}
 &\{Z_s>\kappa\log(1+\min\{s,t-s\})\text{ for }0\le s\le t\}\\
 &\qquad\subseteq\{Z_s>r_t(s)\text{ for }0\le s\le t\}\\
 &\qquad\subseteq\{Z_s>0\text{ for }0\le s\le t\}.
 \end{aligned}
\]
All three events are intersected with $\{X_t=l\}$ when applying
the kernel estimates. Theorem~\ref{thm:logarithmic-arch} and
Theorem~\ref{thm:killed-two-sided}, with the parameter uniformity from
Lemma~\ref{lem:uniform-tilt}, bound
$\mathbb P_j^{\lambda_t}\{X_t=l,\ Z_s>r_t(s)\text{ for }0\le s\le t\}$
above and below by positive constant multiples of $yz/t^{3/2}$.
The constants and endpoint thresholds are common for $\lambda_t$ in
a fixed compact neighborhood of $\lambda$, on the endpoint window in
\eqref{eq:positive-log-boundary}.

It remains to return to $\mathbb P_j^\lambda$. For each fixed $t$, use
the same path event under both laws, namely
\[
 E=\{y+X_s-j+cs>\kappa\log(1+s)\text{ for }0\le s\le t\}
   =\{Z_s>r_t(s)\text{ for }0\le s\le t\}.
\]
By \eqref{eq:path-tilt-identity}, with $\eta=\lambda_t$,
$\mathbb P_j^\lambda(E,X_t=l)$ equals
$\mathbb P_j^{\lambda_t}(E,X_t=l)$ multiplied by
\[
 e^{(k(\lambda_t)-k(\lambda))t+(\lambda_t-\lambda)(l-j)}
 \frac{\phi_j(\lambda_t)\phi_l(\lambda)}
      {\phi_j(\lambda)\phi_l(\lambda_t)}.
\]
Set $\delta_t=\lambda_t-\lambda$. The logarithm of
its exponential factor is
\[
 [k(\lambda_t)-k(\lambda)-c_t\delta_t]t
                      +\delta_t(z-y),
\]
because $l-j=-c_tt+z-y$. Taylor's formula bounds its absolute value by
$C t\delta_t^2+2A\sqrt t\,|\delta_t|$, which tends to zero as $t\to\infty$. The
ratio
$\phi_j(\lambda_t)\phi_l(\lambda)/
(\phi_j(\lambda)\phi_l(\lambda_t))$ in
\eqref{eq:path-tilt-identity} is bounded above and below by positive
constants for all sufficiently large $t$, uniformly in $j,l\in\Z$.
Indeed, $\lambda_t\to\lambda$, each $\phi_i$ is positive and
continuous, and there are only $P$ phases. Combining this bound with
the exponential-factor estimate transfers
\eqref{eq:positive-log-boundary} from $\mathbb P_j^{\lambda_t}$
to $\mathbb P_j^\lambda$.\qedhere

\end{proof}
\section{Examples and further questions}\label{sec:examples}

\subsection{A two-periodic reaction}

Let $P=2$, $\beta\in\R$ and $a>|\beta|$, and use the
normalization $b(o)=0$. For $x\in\Tq$ and $s\in[0,1]$, set
\[
 f(x,s)=\bigl(a+\beta(-1)^{b(x)}\bigr)s(1-s).
\]
Thus $f_i(s)=(a+\beta(-1)^i)s(1-s)$ for $i\in\{0,1\}$,
with $a_0=f_0'(0)=a+\beta$ and $a_1=f_1'(0)=a-\beta$.
Since $a>|\beta|$, both $f_0$ and $f_1$ satisfy \eqref{eq:hyp}.
Each edge changes both $b(x)$ and $r(x)$ by one in absolute value.
Together with $b(o)=r(o)=0$, this gives $b(x)\equiv r(x)\pmod2$.
Hence $f(x,s)=f_{b(x)}(s)=f_{r(x)}(s)$, with
$f_{j+2}=f_j$ for every $j\in\Z$.

For $P=2$, the indices $i-1$ and $i+1$ coincide modulo $P$.
For $\lambda\in\R$, the corresponding terms in \eqref{eq:matrix} therefore add to give
\[
 M_\lambda=
 \begin{pmatrix}
 a+\beta-D(q+1)&D(e^\lambda+qe^{-\lambda})\\
 D(e^\lambda+qe^{-\lambda})&a-\beta-D(q+1)
 \end{pmatrix}.
\]
The off-diagonal entries are positive, so the larger eigenvalue of
$M_\lambda$ is its principal eigenvalue $k(\lambda)$. Computing this
eigenvalue and using $e^{\lambda_0}=\sqrt q$ gives
\begin{equation}
\label{eq:alternating}
 \begin{aligned}
 k(\lambda)
 &=a-D(q+1)
   +\sqrt{\beta^2+D^2(e^\lambda+qe^{-\lambda})^2},\\
 \mu&=a-D(q+1)+\sqrt{\beta^2+4D^2q}.
 \end{aligned}
\end{equation}
Let $u$ solve \eqref{eq:tree-cauchy} with $u_0\not\equiv0$,
$0\le u_0\le1$, and finite $\supp u_0$.
By Theorems~\ref{thm:threshold} and~\ref{thm:radial-threshold},
\[
 \begin{aligned}
 \|u(t,\cdot)\|_{\ell^\infty(\Tq)}&\to0
 &&\text{ as }t\to\infty,\quad\text{if }\mu\le0,\\
 u(t,x)&\to1
 &&\text{ as }t\to\infty,\quad\text{if }\mu>0,
 \quad x\in\Tq\text{ fixed}.
 \end{aligned}
\]
Since $f(x,u)=f_{r(x)}(u)$ for this assignment,
Theorem~\ref{thm:periodic-bramson-distance} gives the ball inclusions
\eqref{eq:periodic-bramson-balls}, where
$m_*(t)=c_+t-3\log(1+t)/(2\lambda_+)$, for every
$\theta\in(0,1)$ when $\mu>0$.

Formula~\eqref{eq:alternating} shows that $\mu$ depends on both
$a=(a_0+a_1)/2$ and $|\beta|=|a_0-a_1|/2$.
Take $q=2$, $D=1$, $a=17/100$ and $\beta=1/10$.
Then $a_0=27/100$, $a_1=7/100$ and $\bar a=17/100$, while
\[
 \frac{17}{100}<3-2\sqrt2=\rho_2,
 \qquad
 \mu=-\frac{283}{100}+\sqrt{\frac{801}{100}}>0.
\]
The inequalities $a-\rho_2<0<\mu$ follow from
\[
 8<\frac{80089}{10000}
   =\left(\frac{283}{100}\right)^2
   <\frac{801}{100}.
\]
Taking positive square roots gives $2\sqrt2<283/100=3-17/100$ and
$\sqrt{801/100}>283/100$.
Replacing $a+\beta$ and $a-\beta$ by $a=17/100$ therefore changes the
upper spectral endpoint of the linearized tree operator from
$\mu>0$ in \eqref{eq:alternating} to $a-\rho_2<0$.

For a general period $P$, suppose that the numbers
$a_i=f_i'(0)$, $0\le i<P$, are positive and not all equal.
Recall $\bar a=P^{-1}\sum_{i=0}^{P-1}a_i$ and
$d_q(\lambda)=D(e^\lambda+qe^{-\lambda}-(q+1))$.
The strict bound in Lemma~\ref{lem:matrix}, together with
$d_q(\lambda_0)=-\rho_q$ and $k(\lambda_0)=\mu$, gives
\[
 k(\lambda)>\bar a+d_q(\lambda)\quad\text{ for } \lambda\in\R,
 \qquad \mu>\bar a-\rho_q.
\]
Evaluating $k(\lambda)>\bar a+d_q(\lambda)$ at the positive
minimizer $\lambda_+$ gives
\begin{equation}
\label{eq:heterogeneity-speed}
 c_+=\frac{k(\lambda_+)}{\lambda_+}
 >\frac{\bar a+d_q(\lambda_+)}{\lambda_+}
 \ge c_*(\bar a).
\end{equation}
The last inequality uses
$c_*(\bar a)=\min_{\lambda>0}(\bar a+d_q(\lambda))/\lambda$
from \eqref{eq:homogeneous-parameters}.

\subsection{Dependence on the branching parameter}

Fix $a,D>0$. For this calculation, allow real $q>1$ in the
homogeneous dispersion formula \eqref{eq:homogeneous-parameters};
the tree equations use integer values $q\ge2$.
For $\lambda\in\R$ and $q>1$, write
\[
 \Psi_a(\lambda,q)
 =a+D(e^\lambda+qe^{-\lambda}-(q+1)),
 \qquad
 c_*(a,q)=\min_{\lambda>0}\frac{\Psi_a(\lambda,q)}{\lambda}.
\]
To locate the minimizer, define
\[
 G(\lambda,q)
 :=\lambda\partial_\lambda\Psi_a(\lambda,q)-\Psi_a(\lambda,q).
\]
For every fixed $q>1$, one has $G(0,q)=-a<0$,
$G(\lambda,q)\to\infty$ as $\lambda\to\infty$, and
\[
 G_\lambda(\lambda,q)=D\lambda(e^\lambda+qe^{-\lambda})>0
 \quad\text{ for } \lambda>0.
\]
Thus $G(\cdot,q)$ has exactly one positive zero. Since
$\partial_\lambda(\Psi_a(\lambda,q)/\lambda)=G(\lambda,q)/\lambda^2$,
the quotient decreases before this zero and increases after it.
Consequently, the unique positive solution
$\lambda_*=\lambda_*(a,q)$ of $G(\lambda,q)=0$ minimizes
$\lambda\mapsto\Psi_a(\lambda,q)/\lambda$ over $(0,\infty)$.
For integer $q\ge2$ and constant phase derivatives $a_i=a$,
\eqref{eq:homogeneous-parameters} gives
$\lambda_*(a,q)=\lambda_+$ and $c_*(a,q)=c_*(a)=c_+$.
Also,
\[
 G_q(\lambda,q)=D\bigl(1-(1+\lambda)e^{-\lambda}\bigr)>0
 \quad\text{ for } \lambda>0,
\]
because $e^\lambda>1+\lambda$. Since
$G_\lambda(\lambda_*(a,q),q)>0$, the implicit function theorem shows
that $q\mapsto\lambda_*(a,q)$ is continuously differentiable on
$(1,\infty)$, with
\begin{equation}
\label{eq:q-derivatives}
 \frac{\partial\lambda_*}{\partial q}
 =-\frac{1-(1+\lambda_*)e^{-\lambda_*}}
 {\lambda_*(e^{\lambda_*}+qe^{-\lambda_*})}<0.
\end{equation}
At $\lambda=\lambda_*(a,q)$, one has
$\partial_\lambda(\Psi_a(\lambda,q)/\lambda)=0$. Differentiating
$c_*(a,q)=\Psi_a(\lambda_*(a,q),q)/\lambda_*(a,q)$ therefore yields
\[
 \frac{\partial c_*}{\partial q}
 =\frac{D(e^{-\lambda_*}-1)}{\lambda_*}<0,
 \qquad
 \frac{\partial}{\partial q}\left(\frac3{2\lambda_*}\right)
 =-\frac3{2\lambda_*^2}\frac{\partial\lambda_*}{\partial q}>0.
\]
Consequently, $\lambda_*(a,q)$ and $c_*(a,q)$ strictly decrease
with real $q>1$, while $3/(2\lambda_*(a,q))$ strictly increases.

Now let $q\ge2$ be an integer, assume \eqref{eq:hyp}, and suppose
$f_i'(0)=a$ for $0\le i<P$. If $f(x,u)=f_{b(x)}(u)$ and,
for integers $A\le B$,
\[
 \mathbf1_{\{b(x)\le A\}}\le u_0(x)\le
 \mathbf1_{\{b(x)\le B\}}
 \quad\text{ for }x\in\Tq,
\]
then Corollary~\ref{cor:common-step} applies whether
$c_*(a,q)<0$, $c_*(a,q)=0$ or $c_*(a,q)>0$.

The assertions for $r(x)$ and $s(x)$ in
Theorem~\ref{thm:periodic-bramson-distance}\textup{(i)--(ii)}
require the following initial data.
For $f(x,u)=f_{r(x)}(u)$, take $0\le u_0\le1$ with
$u_0\not\equiv0$ and finite $\supp u_0$.
For $f(x,u)=f_{s(x)}(u)$, assume instead that
\[
 \eta\mathbf1_{\mathcal N_{R_1}(\mathcal G)}\le u_0
 \le\mathbf1_{\mathcal N_{R_2}(\mathcal G)},
 \qquad 0<\eta\le1,\quad0\le R_1\le R_2<\infty.
\]
Since $\mu=a-\rho_q$ when the phase derivatives coincide, both
distance assertions require
\[
 a>\rho_q=D(\sqrt q-1)^2
 \quad\Longleftrightarrow\quad
 1<q<\left(1+\sqrt{a/D}\right)^2.
\]
For $\mathbf1_{\{b\le A\}}\le u_0\le\mathbf1_{\{b\le B\}}$,
Corollary~\ref{cor:common-step} locates the levels in $b(x)$.
For the stated finitely supported data and geodesic-neighborhood
data, respectively, \eqref{eq:periodic-bramson-balls} and
\eqref{eq:periodic-bramson-tubes} locate the levels in $r(x)$
and $s(x)$ when $a>\rho_q$. In all three conclusions, the level
position differs by a bounded amount from
\[
 m_*(t)=c_*(a,q)t-\frac{3}{2\lambda_*(a,q)}\log(1+t),\qquad t\ge0.
\]
By \eqref{eq:q-derivatives}, $3/(2\lambda_*(a,q))$ increases
with real $q>1$.

\subsection{Level positions and profile convergence}

Assume \eqref{eq:hyp}. Let $w$ solve \eqref{eq:lattice} with
$\mathbf1_{\{j\le A\}}\le w_j(0)\le\mathbf1_{\{j\le B\}}$
for every $j\in\Z$, where $A\le B$ are integers. Put
$H_\theta(t)=\max\{j\in\Z:w_j(t)\ge\theta\}$ whenever this
maximum exists, and recall
$m_*(t)=c_+t-3\log(1+t)/(2\lambda_+)$ for $t\ge0$.
For every fixed $\theta\in(0,1)$,
Theorem~\ref{thm:periodic-bramson-step} gives
\begin{equation}
\label{eq:periodic-location-summary}
 H_\theta(t)
 =c_+t-\frac3{2\lambda_+}\log t+O_{\theta,A,B}(1)
 \quad\text{ as } t\to\infty.
\end{equation}
Replacing $\log(1+t)$ by $\log t$ changes $m_*(t)$ by $o(1)$
as $t\to\infty$. Formula \eqref{eq:periodic-location-summary} holds
for all three signs of $c_+$.
For a solution $u$ of \eqref{eq:tree-cauchy}, assume $\mu>0$ and
fix $\theta\in(0,1)$.

If $f(x,u)=f_{b(x)}(u)$ and $u_0:\Tq\to[0,1]$ is nonzero
and finitely supported, Theorem~\ref{thm:horo-compact-bramson} gives
$H_\theta^b(t)=m_*(t)+O_{\theta,u_0}(1)$ as $t\to\infty$,
where $H_\theta^b(t)=\max\{b(x):x\in\Tq,\ u(t,x)\ge\theta\}$.
The maximum exists for all sufficiently large $t$.

If $f(x,u)=f_{r(x)}(u)$ and $u_0$ is
nonzero, takes values in $[0,1]$, and has finite support, then
Theorem~\ref{thm:periodic-bramson-distance} gives
$R_\theta^o(t)=m_*(t)+O(1)$ as $t\to\infty$, where
$R_\theta^o(t)=\max\{r(x):x\in\Tq,\ u(t,x)\ge\theta\}$.
If $f(x,u)=f_{s(x)}(u)$ and
$\eta\mathbf1_{\mathcal N_{R_1}(\mathcal G)}\le u_0
\le\mathbf1_{\mathcal N_{R_2}(\mathcal G)}$, with
$0<\eta\le1$ and $0\le R_1\le R_2<\infty$, then the same
theorem gives $R_\theta^{\mathcal G}(t)=m_*(t)+O(1)$ as
$t\to\infty$, where
$R_\theta^{\mathcal G}(t)=\max\{s(x):x\in\Tq,\ u(t,x)\ge\theta\}$.
Both maxima exist for all sufficiently large $t$ under their
respective hypotheses.
The earlier estimate \eqref{eq:step-position} reads
$H_\theta(t)\le c_+t-\log(1+t)/(2\lambda_+)+C_\theta$
for all sufficiently large $t$.
The lower and upper bounds giving the bounded remainder in
\eqref{eq:periodic-location-summary} come from Lemmas~\ref{lem:sharp-step-upper}
and~\ref{lem:sharp-step-lower}, as used in the proof of
Theorem~\ref{thm:periodic-bramson-step}.

For propagation along $\mathcal G$, assume
$\mu\le0<\sigma_{\mathcal G}$ and $f(x,u)=f_{s(x)}(u)$.
Let $u$ solve \eqref{eq:tree-cauchy} with $u_0\not\equiv0$,
$0\le u_0\le1$, and finite $\supp u_0$. For the level position on $\mathcal G$,
put $L_\theta(t)=\max\{|n|:n\in\Z,\ u(t,\gamma(n))\ge\theta\}$
whenever the maximum exists. Theorem~\ref{thm:tangential-bramson}
gives, as $t\to\infty$,
\[
 L_\theta(t)
 =c_{\parallel}t-\frac3{2\lambda_{\parallel}}\log t+O_\theta(1),
 \quad\text{ for }0<\theta<p_*(\gamma(0)).
\]
Here $\lambda_{\parallel}$ is the unique minimizer of
$\Lambda_{\parallel}(\lambda)/\lambda$ over $\lambda>0$, as in
\eqref{eq:tangential-speed-intro}. The proof uses regeneration at jumps that change
$\ell_{\mathcal G}(x)$ in Lemma~\ref{lem:tangential-regeneration},
the two-sided killed-kernel bounds in
Lemma~\ref{lem:tangential-killed}, and the conditional occupation
estimate in Lemma~\ref{lem:tangential-backward-log}.

For the homogeneous equation in \cite[(1.1), Theorem~2]{BFRZ},
let $u^0=(u_j^0)_{j\in\Z}$ be nonzero, with $0\le u_j^0\le1$
for $j\in\Z$ and $u_j^0=0$ for $j\ge J$, for some $J\in\Z$.
For these data, Besse, Faye, Roquejoffre and Zhang obtain approximation
by a logarithmically shifted minimal front with a bounded
time-dependent phase correction, uniformly on $j\ge0$.
They also prove locally uniform convergence in space and time after
recentering at level sets. Guo and Wu \cite[(1.1), Theorems~1.1--1.2]{GuoWu} prove
stability for periodic fronts of speed strictly above the minimum
under prescribed exponential asymptotics. For a front in their equation with speed $c$ strictly above its
minimal speed, normalized to have leading term
$e^{-\lambda(j-ct)}v_j$, where $\lambda>0$ is the smaller positive
exponent corresponding to its speed and $v=(v_j)_{j\in\Z}$ is a
positive periodic vector normalized by $\max_jv_j=1$, their stability
theorem assumes data
$u^0=(u_j^0)_{j\in\Z}$ satisfying
$0\le u_j^0\le\min\{1,e^{-\lambda j}v_j\}$,
$\liminf_{j\to-\infty}u_j^0>0$, and
$u_j^0/(e^{-\lambda j}v_j)\to1$ as $j\to+\infty$.
Their uniqueness theorem identifies profiles of the same speed up to
time translation when their leading exponential terms differ by a positive
constant factor.

Assume \eqref{eq:hyp} and \eqref{eq:derivative-kpp}.
For \eqref{eq:lattice}, assume $\mu>0$ and fix
$W_j(t)=\Phi_j(j-ct)$ from Theorem~\ref{thm:all-speeds},
with $c\ge c_+$. Let $V$ solve \eqref{eq:lattice} with
\[
 \Phi_j(j+A)\le V_j(0)\le\Phi_j(j-A)
 \quad\text{ for }j\in\Z
\]
for some $A\ge0$. If
$(e^{\lambda_+j}(V_j(0)-W_j(0)))_{j\in\Z}\in\ell^p(\Z)$
for some $1\le p<\infty$, then Theorem~\ref{thm:wave-stability}
gives $\sup_{j\in\Z}|V_j(t)-W_j(t)|\to0$ as $t\to\infty$.

For \eqref{eq:tree-cauchy} with $f(x,u)=f_{b(x)}(u)$, let
$(\widehat W_j(t))_{j\in\Z}$ be a front from
Theorem~\ref{thm:all-speeds} and set
$\widehat W(t,x)=\widehat W_{b(x)}(t)$.
Let $u$ solve \eqref{eq:tree-cauchy}. Under the pointwise bounds
\[
 \widehat W_{b(x)+mP}(0)\le u_0(x)\le\widehat W_{b(x)-mP}(0)
 \quad\text{ for }x\in\Tq
\]
for some integer $m\ge0$, as in \eqref{eq:tree-front-bracket}, and
\[
 \bigl(e^{\lambda_+b(x)}(u_0(x)-\widehat W(0,x))\bigr)_{x\in\Tq}
 \in\ell^p(\Tq)
 \quad\text{ for some }1\le p<\infty,
\]
Theorem~\ref{thm:tree-wave-stability} gives
$\|u(t,\cdot)-\widehat W(t,\cdot)\|_{\ell^\infty(\Tq)}\to0$
as $t\to\infty$.

Let $u$ solve \eqref{eq:tree-cauchy} with $u_0\not\equiv0$,
$0\le u_0\le1$, and finite $\supp u_0$, for $f(x,u)=f_{s(x)}(u)$
and $\mu\le0<\sigma_{\mathcal G}$. By
\eqref{eq:finite-mass-level-count}, $u(t,\cdot)\in\ell^1(\Tq)$ for
every finite $t\ge0$, so $u(t,\gamma(n))\to0$ as $|n|\to\infty$.
Since $p_*(\gamma(n))=p_*(\gamma(0))>0$, taking the supremum along
$\mathcal G$ gives
$\|u(t,\cdot)-p_*(\cdot)\|_{\ell^\infty(\Tq)}
\ge p_*(\gamma(0))$ for every finite $t\ge0$, as in
\eqref{eq:no-global-geodesic-convergence}.
Under \eqref{eq:strict-quotient},
Theorem~\ref{thm:tangential-spreading} gives
\[
 \sup_{\{x\in\Tq:\,|\ell_{\mathcal G}(x)|\le ct\}}
 |u(t,x)-p_*(x)|\to0
 \quad\text{ as }t\to\infty,\quad\text{for }0<c<c_{\parallel}.
\]
The supremum imposes no restriction on $s(x)$. Each fixed
$x\in\Tq$ belongs to $\{|\ell_{\mathcal G}|\le ct\}$ for all
sufficiently large $t$, so $u(t,x)\to p_*(x)$ as $t\to\infty$.
The final subsection asks whether traveling profiles approximate
$u(t,x)$ on the sets where
$|\ell_{\mathcal G}(x)-m_{\parallel}(t)|\le R$ or
$|-\ell_{\mathcal G}(x)-m_{\parallel}(t)|\le R$ for fixed $R>0$,
with $m_{\parallel}(t)=c_{\parallel}t
-3\log(1+t)/(2\lambda_{\parallel})$ for $t\ge0$.

\subsection{Stationary solutions with transverse decay}

\begin{proof}[Proof of Theorem~\ref{thm:stationary-family}]

Assume \eqref{eq:hyp}, \eqref{eq:strict-quotient} and
$\mu\le0<\sigma_{\mathcal G}$, with $f(x,u)=f_{s(x)}(u)$.
Write $p=(p_r)_{r\in\N_0}$, where $p_r=p_*(x)$ whenever
$s(x)=r$; in particular, $p_0=p_*(\gamma(0))$.
By Theorem~\ref{thm:geodesic-threshold}, $p$ satisfies the stationary
interior equation in \eqref{eq:radial} and the stationary boundary
equation \eqref{eq:geodesic-boundary}. Thus $x\mapsto p_{s(x)}$ is
a stationary solution of the differential equation in
\eqref{eq:tree-cauchy}.
Fix $\alpha\in(p_0,1)$ and write $U_r=U_r^\alpha$ throughout
this proof.
By \eqref{eq:geodesic-minimal-tail},
$0<p_r<1$ and $p_r\le h_r/h_0$ for every $r\ge0$, where
$h_r=q^{-r/2}\phi_r(\lambda_0)$. In particular,
$p_r=O(h_r)$ as $r\to\infty$.

For the fixed $\alpha\in(p_0,1)$, determine the entries $U_r$ by
\begin{equation}\label{eq:stationary-family-recursion}
 \begin{aligned}
 U_0&=\alpha,\qquad
 U_1=\alpha-\frac{f_0(\alpha)}{D(q-1)},\\
 qU_{r+1}&=(q+1)U_r-U_{r-1}-D^{-1}f_r(U_r)
 \quad\text{ for }r\ge1.
 \end{aligned}
\end{equation}
We prove by induction that each value belongs to $(0,1)$.
Thus every evaluation of $f_r$ remains in its prescribed domain.
We will prove
\begin{equation}\label{eq:stationary-family-order}
 \begin{gathered}
 p_r<U_r<1,\qquad U_{r+1}<U_r
 \quad\text{ for }r\in\N_0,\\
 \lim_{r\to\infty}U_r=0.
 \end{gathered}
\end{equation}
For $r\in\N_0$ and $0<s\le1$, put $g_r(s)=f_r(s)/s$.
Recall $\nu_0=1$ and $\nu_r=(q-1)q^{r-1}$ for $r\ge1$.
Define $B_r=D\nu_{r+1}p_rp_{r+1}>0$ for every $r\in\N_0$.
For each entry $U_j$ already constructed by
\eqref{eq:stationary-family-recursion}, put $V_j=U_j/p_j$.
For $r\in\N_0$, whenever both $U_r$ and $U_{r+1}$ have been
constructed, set
\[
 J_r=B_r(V_{r+1}-V_r)
     =D\nu_{r+1}(p_rU_{r+1}-p_{r+1}U_r).
\]
At $r=0$, \eqref{eq:stationary-family-recursion} and stationarity
of $p$ give
\[
 \begin{aligned}
 D(q-1)(U_1-U_0)+f_0(U_0)&=0,\\
 D(q-1)(p_1-p_0)+f_0(p_0)&=0.
 \end{aligned}
\]
Since $J_0=D(q-1)(p_0U_1-p_1U_0)$ and $U_0=\alpha$, the two
boundary equations give
$J_0=\alpha f_0(p_0)-p_0f_0(\alpha)$. Thus
\eqref{eq:strict-quotient} with $p_0<\alpha$ gives
\begin{equation}\label{eq:stationary-family-boundary-flux}
 J_0=p_0\alpha\bigl(g_0(p_0)-g_0(\alpha)\bigr)>0.
\end{equation}
Since $B_0>0$, it follows that $V_1>V_0=\alpha/p_0>1$.
By \eqref{eq:stationary-family-recursion} and $f_0(\alpha)>0$,
we also have $U_1<U_0=\alpha<1$. Thus $p_1<U_1<\alpha$.

Fix an integer $r\ge1$. Suppose that $p_j<U_j<1$ has been proved
for $0\le j\le r$, that $U_j<U_{j-1}$ for $1\le j\le r$,
and that $J_{r-1}>0$. Since $U_r\in(0,1)$, the value $f_r(U_r)$ is
defined. The recursion therefore determines a finite real number
$U_{r+1}$. At index $r$, the recursion and the stationary equation
for $p$ read
\[
 \begin{aligned}
 D(U_{r-1}+qU_{r+1}-(q+1)U_r)+f_r(U_r)&=0,\\
 D(p_{r-1}+qp_{r+1}-(q+1)p_r)+f_r(p_r)&=0.
 \end{aligned}
\]
By $\nu_{r+1}=q\nu_r$ and the definition of $J_r$,
\[
 J_r-J_{r-1}
 =D\nu_r\bigl[p_r(U_{r-1}+qU_{r+1})
                 -U_r(p_{r-1}+qp_{r+1})\bigr].
\]
Subtracting the equation for $p$ multiplied by $\nu_rU_r$
from the equation for $U$ multiplied by $\nu_rp_r$ gives, by
\eqref{eq:strict-quotient} with $p_r<U_r$,
\begin{equation}\label{eq:stationary-family-flux}
 J_r-J_{r-1}
 =\nu_rp_rU_r\bigl(g_r(p_r)-g_r(U_r)\bigr)>0.
\end{equation}
In view of $J_{r-1}>0$ and $B_r>0$, we obtain $J_r>0$ and
$V_{r+1}>V_r>1$, so $U_{r+1}>p_{r+1}$. Put
$\delta_j=U_j-U_{j-1}$ for $1\le j\le r+1$. The induction
hypothesis gives $\delta_r<0$, while \eqref{eq:hyp} gives
$f_r(U_r)>0$. Hence \eqref{eq:stationary-family-recursion} yields
\[
 q\delta_{r+1}=\delta_r-D^{-1}f_r(U_r)<0.
\]
Combining $q\delta_{r+1}<0$ with $V_{r+1}>1$ yields
$p_{r+1}<U_{r+1}<U_r<1$. The induction therefore defines
$U_r$ for every $r\ge0$ and proves both order inequalities in
\eqref{eq:stationary-family-order}.
The positive decreasing sequence $U_r$ therefore has a limit
$L\in[0,\alpha]$. Fix $i\in\{0,\ldots,P-1\}$ and take
$r=nP+i$ with $n\to\infty$. Then $U_{r-1}$, $U_r$, and
$U_{r+1}$ all tend to $L$. Periodicity gives $f_r=f_i$, so
continuity of $f_i$ allows passage to
\eqref{eq:stationary-family-recursion} and yields $f_i(L)=0$.
The equality $f_i(L)=0$ holds for $0\le i<P$. By
\eqref{eq:hyp}, $f_i$ is positive
on $(0,1)$; since $L\le\alpha<1$, we have $L=0$.
This proves \eqref{eq:stationary-family-order}.

For each $\alpha\in(p_0,1)$, the sequence $U^\alpha=U$ defined
by \eqref{eq:stationary-family-recursion} therefore makes
$x\mapsto U_{s(x)}^\alpha$ a positive stationary solution of the
differential equation in \eqref{eq:tree-cauchy}.
If $\alpha,\widetilde\alpha\in(p_0,1)$ and $\alpha\ne\widetilde\alpha$, then
$U_0^\alpha=\alpha\ne\widetilde\alpha=U_0^{\widetilde\alpha}$.
Moreover, by \eqref{eq:stationary-family-order},
\[
 \sup_{\{x\in\Tq:\,s(x)\ge R\}}U_{s(x)}^\alpha
 =U_R^\alpha\to0
 \quad\text{ as }R\to\infty,\quad R\in\N_0.
\]
If $\sup_{r\ge0}U_r/h_r<\infty$, set
$C=\max\{h_0^{-1},\sup_{r\ge0}U_r/h_r\}$. Then
$0<U_r,p_r\le1$ and $U_r,p_r\le Ch_r$ for every $r\ge0$.
For each of $v=U$ and $v=p$, the stationary equations are
$D(q-1)(v_1-v_0)+f_0(v_0)=0$ and
$D(v_{r-1}+qv_{r+1}-(q+1)v_r)+f_r(v_r)=0$ for $r\ge1$.
Together with $0<U_r,p_r\le1$ and $U_r,p_r\le Ch_r$, these equations
verify the hypotheses of Lemma~\ref{lem:cutoff-steady-uniqueness}, which gives
$U_r=p_r$ for every $r\ge0$, contradicting $U_0=\alpha>p_0$.
Hence $\sup_{r\ge0}U_r/h_r=\infty$, as asserted in
\eqref{eq:stationary-family-decay-intro}.

Conversely, let $W=(W_r)_{r\in\N_0}$ satisfy $0<W_r\le1$ for
every $r\in\N_0$, and suppose that $x\mapsto W_{s(x)}$ is a
stationary solution of the differential equation in
\eqref{eq:tree-cauchy}. By the minimality of $p_*$ in
Theorem~\ref{thm:geodesic-threshold}, $W_r\ge p_r$ for every
$r\in\N_0$, and hence $W_0\in[p_0,1]$.
The stationary boundary equation determines $W_1$ from $W_0$, and
the interior equations successively determine $W_{r+1}$ from
$W_{r-1}$ and $W_r$, exactly as in
\eqref{eq:stationary-family-recursion}. Thus there is at most one
such sequence for each value of $W_0$. Every value in $(p_0,1)$
is realized by the construction above. The endpoint values $W_0=p_0$ and $W_0=1$ give $W_r=p_r$ and
$W_r=1$ for every $r\ge0$, respectively. Indeed, each subsequent
entry is fixed by \eqref{eq:stationary-family-recursion}, while
stationarity of $p$ and $f_i(1)=0$ in \eqref{eq:hyp} verify the
two endpoint sequences. Thus the possible values of $W_0$ are
exactly $[p_0,1]$, and each value determines exactly one stationary
sequence $W$ with $0<W_r\le1$ for $r\in\N_0$.\qedhere

\end{proof}

\subsection{Further questions}

We ask whether the estimate \eqref{eq:periodic-location-summary}
can be accompanied by convergence to a pulsating front of
\eqref{eq:lattice}. Assume \eqref{eq:hyp}, let $w$ solve
\eqref{eq:lattice} with
$\mathbf1_{\{j\le A\}}\le w_j(0)\le\mathbf1_{\{j\le B\}}$
for every $j\in\Z$, where $A\le B$ are integers, and suppose
$c_+\ne0$.
Does there exist a pulsating front profile
$\Phi=(\Phi_i)_{i\in\Z}$ of speed $c_+$, as defined in
Subsection~\ref{sec:waves}, and a function
$a:[1,\infty)\to\R$ with $\sup_{t\ge1}|a(t)|<\infty$ such that
\[
 \sup_{j\in\Z}\left|
 w_j(t)-\Phi_j\left(j-c_+t+
        \frac{3}{2\lambda_+}\log(1+t)+a(t)\right)
 \right|\to0\quad\text{ as } t\to\infty?
\]
One may further ask whether $a(t)$ can be chosen so that
$a(t)\to a_\infty\in\R$ as $t\to\infty$.
The periodic continuum result \cite[Theorem~1.2]{HNRR}
and the homogeneous lattice result \cite[Theorem~2]{BFRZ} give
approximation by minimal fronts with bounded time-dependent phase
corrections on $x>0$ and $j\ge0$, respectively.
For the homogeneous continuum logistic equation,
\cite{NRRsingle} proves convergence with a constant phase correction
for compact perturbations of step data.

Assume \eqref{eq:hyp} and $c_+=0$, and let $w$ solve
\eqref{eq:lattice} with
$\mathbf1_{\{j\le A\}}\le w_j(0)\le\mathbf1_{\{j\le B\}}$
for integers $A\le B$. A stationary front $p=(p_j)_{j\in\Z}$,
as defined in Subsection~\ref{sec:waves}, satisfies
\[
 D(p_{j-1}+qp_{j+1}-(q+1)p_j)+f_j(p_j)=0
 \quad\text{ for }j\in\Z.
\]
The shifts
$p_j\mapsto p_{j-mP}$, $m\in\Z$, preserve the coefficients.
For every $t_n\to\infty$ and $m_n\in\Z$ satisfying
$Pm_n=-3\log(1+t_n)/(2\lambda_+)+O(1)$ as $n\to\infty$,
do there exist a subsequence $n_k\to\infty$ and a stationary front
$p=(p_j)_{j\in\Z}$ of \eqref{eq:lattice} such that
\[
 \max_{\substack{j\in\Z\\|j|\le M}}|w_{j+m_{n_k}P}(t_{n_k})-p_j|\to0
 \quad\text{ as }k\to\infty
 \quad\text{for every }M\in\N_0?
\]
For $c_+\ne0$, what is the leading asymptotic expansion of
$\Phi_i(z)$ as $z\to+\infty$, for $0\le i<P$, when
$\Phi=(\Phi_i)_{i\in\Z}$ is a pulsating front profile of
\eqref{eq:lattice} at speed $c_+$?
For moving fronts at speeds $c>c_+$ with $c\ne0$,
\eqref{eq:front-tail} specifies the exponential tail.

For the assignment $f(x,u)=f_{s(x)}(u)$, assume \eqref{eq:hyp}
and $\mu\le0<\sigma_{\mathcal G}$. Define
$h_r=q^{-r/2}\phi_r(\lambda_0)$ for $r\in\N_0$.
Consider functions $p:\Tq\to(0,1]$ satisfying
\[
 D\Delta_{\Tq}p(x)+f_{s(x)}(p(x))=0
 \quad\text{ for }x\in\Tq,
 \qquad
 \sup_{x\in\Tq}\frac{p(x)}{h_{s(x)}}<\infty.
\]
Theorem~\ref{thm:geodesic-critical-convergence} proves uniqueness
in this class under \eqref{eq:strict-quotient}.
Can uniqueness in this class be proved under \eqref{eq:hyp} without
\eqref{eq:strict-quotient}, or under a condition weaker than
\eqref{eq:strict-quotient}?
Under \eqref{eq:strict-quotient}, put $p_{*,0}=p_*(\gamma(0))$.
For each $\alpha\in(p_{*,0},1)$,
Theorem~\ref{thm:stationary-family} gives
$U^\alpha=(U_r^\alpha)_{r\in\N_0}$ with $U_0^\alpha=\alpha$.
The sequence $U^\alpha$ satisfies
\eqref{eq:stationary-family-recursion} and
\eqref{eq:stationary-family-order}, and
$x\mapsto U_{s(x)}^\alpha$ is stationary for the differential equation
in \eqref{eq:tree-cauchy}.
What is the asymptotic decay of $U_r^\alpha$ as $r\to\infty$?
For a positive stationary solution $p:\Tq\to(0,1]$ of the differential
equation in \eqref{eq:tree-cauchy}, which data
$u_0:\Tq\to[0,1]$ yield $u(t,x)\to p(x)$ as $t\to\infty$ for
every $x\in\Tq$? Which of these data also yield
$\|u(t,\cdot)-p(\cdot)\|_{\ell^\infty(\Tq)}\to0$
as $t\to\infty$?
Here $u$ solves \eqref{eq:tree-cauchy} with initial datum $u_0$.
{Can one classify stationary solutions $p$ of the differential equation
in \eqref{eq:tree-cauchy} for which} $p(x)\ne p(y)$ for some $x,y\in\Tq$ with
$s(x)=s(y)$?

Finally, assume \eqref{eq:hyp}, $f(x,u)=f_{s(x)}(u)$, and
$\mu\le0<\sigma_{\mathcal G}$. Fix $c>0$ and seek a solution
$u:\R\times\Tq\to(0,1)$ of the differential equation in
\eqref{eq:tree-cauchy} of the form
\[
 u(t,x)=\Phi_{s(x)}(\ell_{\mathcal G}(x)-ct)
 \quad\text{ for }(t,x)\in\R\times\Tq.
\]
Write $p_{*,r}=p_*(x)$ when $s(x)=r$. For every $r\in\N_0$, require
\begin{equation}\label{eq:geodesic-profile-conditions}
 \begin{gathered}
 \Phi_r\in C^1(\R),\qquad
 0<\Phi_r(z)<p_{*,r}\quad\text{ for }z\in\R,\\
 \lim_{z\to-\infty}\Phi_r(z)=p_{*,r},\qquad
 \lim_{z\to+\infty}\Phi_r(z)=0.
 \end{gathered}
\end{equation}
Substitution into \eqref{eq:tree-cauchy} gives, for $z\in\R$,
\begin{equation}\label{eq:geodesic-profile-equations}
 \begin{aligned}
 -c\Phi_0'(z)
 &=D\bigl(\Phi_0(z-1)+\Phi_0(z+1)+(q-1)\Phi_1(z)
                -(q+1)\Phi_0(z)\bigr)+f_0(\Phi_0(z)),\\
 -c\Phi_r'(z)
 &=D\bigl(\Phi_{r-1}(z)+q\Phi_{r+1}(z)-(q+1)\Phi_r(z)\bigr)
       +f_r(\Phi_r(z))\quad\text{ for }r\ge1.
 \end{aligned}
\end{equation}
By \eqref{eq:geodesic-minimal-tail}, for each $\varepsilon>0$
there is an integer $R\ge1$ such that $p_{*,r}<\varepsilon$ for
$r\ge R$. On these indices, both $|\Phi_r(z)-p_{*,r}|$ and
$|\Phi_r(z)|$ are at most $p_{*,r}<\varepsilon$ for every $z\in\R$.
For the finitely many indices $0\le r<R$,
\eqref{eq:geodesic-profile-conditions} supplies $Z_\varepsilon>0$
such that $|\Phi_r(z)-p_{*,r}|<\varepsilon$ for
$z<-Z_\varepsilon$ and $|\Phi_r(z)|<\varepsilon$ for
$z>Z_\varepsilon$. Combining these inequalities with the bound
$p_{*,r}<\varepsilon$ for $r\ge R$ gives
\[
 \lim_{z\to-\infty}\sup_{r\in\N_0}|\Phi_r(z)-p_{*,r}|=0,
 \qquad
 \lim_{z\to+\infty}\sup_{r\in\N_0}|\Phi_r(z)|=0.
\]
For which $c>0$ does \eqref{eq:geodesic-profile-equations} admit
a family $\Phi=(\Phi_r)_{r\in\N_0}$ satisfying
\eqref{eq:geodesic-profile-conditions}?
If $\Phi$ and $\widetilde\Phi$ satisfy
\eqref{eq:geodesic-profile-equations} and
\eqref{eq:geodesic-profile-conditions} at the same speed, must there exist $\zeta\in\R$ such that
$\widetilde\Phi_r(z)=\Phi_r(z+\zeta)$ for all
$(r,z)\in\N_0\times\R$?
Assume that $\Phi$ satisfies \eqref{eq:geodesic-profile-equations}
at $c=c_{\parallel}$ and \eqref{eq:geodesic-profile-conditions}.
Let $u$
solve \eqref{eq:tree-cauchy} with
$u_0\not\equiv0$, $0\le u_0\le1$, and finite $\supp u_0$.
Recall from Theorem~\ref{thm:tangential-bramson} that
$m_{\parallel}(t)=c_{\parallel}t
-3\log(1+t)/(2\lambda_{\parallel})$.
Do there exist $a_+,a_-:[1,\infty)\to\R$ with
$\sup_{t\ge1}(|a_+(t)|+|a_-(t)|)<\infty$ such that,
for every $R>0$ and both choices of the sign,
\[
 \lim_{t\to\infty}
 \sup_{\{x\in\Tq:\,|\pm\ell_{\mathcal G}(x)-m_{\parallel}(t)|\le R\}}
 \left|u(t,x)-\Phi_{s(x)}\bigl(\pm\ell_{\mathcal G}(x)
              -m_{\parallel}(t)+a_\pm(t)\bigr)\right|=0?
\]
Here the supremum is taken to be zero if its indexing set is empty.
The minus sign evaluates the same family $\Phi$ at
$-\ell_{\mathcal G}(x)-m_{\parallel}(t)+a_-(t)$.

\end{document}